\documentclass{memo-l}

\DeclareFontFamily{OT1}{pzc}{}
\DeclareFontShape{OT1}{pzc}{m}{it}{<-> s * [1.10] pzcmi7t}{}
\DeclareMathAlphabet{\mathpzc}{OT1}{pzc}{m}{it}

\usepackage{amsmath,amssymb,amscd,latexsym,mathrsfs,epic,wasysym,latexsym,tikz,mathrsfs,cite,hyperref}
\usepackage{stmaryrd}
\usepackage{pb-diagram}
\usepackage[matrix,arrow]{xy}
\usepackage{commath}
\usepackage{tikz-cd}
\usepackage{bbm}
\usepackage{ytableau}
\usepackage{mathdots}
\usepackage{imakeidx} 
\usepackage{mathtools}
\usepackage{sansmath}
\usepackage{upgreek}

\numberwithin{equation}{section}

\newtheorem{Proposition}[equation]{Proposition}
\newtheorem{Lemma}[equation]{Lemma}
\newtheorem{Theorem}[equation]{Theorem}
\newtheorem{Corollary}[equation]{Corollary}
\theoremstyle{definition}  
\newtheorem{Definition}[equation]{Definition}
\newtheorem{Remark}[equation]{Remark}

\newtheorem{Example}[equation]{Example}

\makeatletter \@addtoreset{equation}{chapter}

\makeatother

\numberwithin{section}{chapter}
\numberwithin{equation}{section}

\let\<\langle
\let\>\rangle

\newcommand\Comment[2][\relax]{\space\par\medskip\noindent%
   \fbox{\begin{minipage}{\textwidth}\textbf{Comment\ifx\relax#1\else---#1\fi}\newline%
        #2\end{minipage}}\medskip
}

\DeclarePairedDelimiterX{\infdivx}[2]{(}{)}{%
  #1\;\delimsize\|\;#2%
}

\def\bi{\text{\boldmath$i$}}
\def\bj{\text{\boldmath$j$}}
\def\bm{\text{\boldmath$m$}}
\def\bk{\text{\boldmath$k$}}
\def\bl{\text{\boldmath$l$}}

\def\b1{\text{\boldmath$1$}}

\def\bg{\text{\boldmath$g$}}

\def\hatbg{\text{\boldmath$\hat g$}}

\def\ba{\text{\boldmath$a$}}

\def\bP{\text{\boldmath$P$}}
\def\bQ{\text{\boldmath$Q$}}

\newcommand{\nth}{\mathop{\rm th}\nolimits}
\newcommand{\gsM}{\mathop{\rm gsMor}\nolimits}
\newcommand{\sM}{\mathop{\rm sMor}\nolimits}
\newcommand{\cC} {\mathcal{C}}
\newcommand{\cY} {\mathcal{TC}}
\newcommand{\cB} {\mathcal{B}}
\newcommand{\cX} {\mathcal{X}}

\newcommand{\cc} {\mathpzc{c}}
\newcommand{\cb} {\mathfrak{b}}

\newcommand{\cM} {\mathcal{M}}

\newcommand{\cT} {\mathcal{T}}
\newcommand{\ct} {\mathpzc{t}}

\newcommand{\cm} {\mathpzc{m}}

\newcommand{\super}{\mathop{\tt s}\nolimits}

\newcommand{\der}{\mathop{\rm der}\nolimits}

\newcommand{\swr}{\wr_{\super}}

\def\a{\mathfrak{a}}

\def\pmod#1{\text{ }(\text{\rm mod } #1)\,}

\newcommand{\Hom}{\operatorname{Hom}}

\newcommand{\End}{\operatorname{End}}
\newcommand{\Br}{\operatorname{Br}}

\newcommand{\im}{\operatorname{im}}

\newcommand{\res}{{\operatorname{res}}}
\newcommand{\soc}{\operatorname{soc}\,}
\newcommand{\head}{\operatorname{head}}

\newcommand{\cha}{\operatorname{char}}

\newcommand{\bideg}{\operatorname{bideg}}

\def\cont{{\operatorname{cont}}}
\def\core{{\operatorname{core}}}
\newcommand{\cus}{{\operatorname{cus}}}

\newcommand{\Z}{\mathbb{Z}}

\newcommand{\F}{\mathbb{F}}

\newcommand{\N}{\mathbb{Z}_{\geq 0}}
\newcommand{\0}{{\bar 0}}
\renewcommand{\1}{{\bar 1}}
\def\eps{{\varepsilon}}
\def\phi{{\varphi}}

\newcommand{\zc}{{\mathsf{c}}}
\newcommand{\zv}{{\mathsf{v}}}
\newcommand{\zz}{{\mathsf{z}}}
\newcommand{\ze}{{\mathsf{e}}}

\newcommand{\za}{{\mathsf{a}}}
\newcommand{\zf}{{\mathsf{f}}}
\newcommand{\zb}{{\mathsf{b}}}
\newcommand{\zu}{{\mathsf{u}}}

\newcommand{\zB}{{\mathscr{B}}}

\newcommand{\funA}{{\mathcal A}}
\newcommand{\funF}{{\mathcal F}}
\newcommand{\funG}{{\mathcal G}}

\newcommand{\funE}{{\mathcal E}}
\newcommand{\funK}{{\mathcal K}}

\newcommand{\cR}{{\mathcal R}}
\newcommand{\cRC}{{\mathcal{RC}}}

\newcommand{\di}{{\operatorname{div}}}

\newcommand{\Fock}{{\mathscr F}}

\newcommand{\quot}{\operatorname{quot}}

\newcommand{\ga}{\gamma}
\newcommand{\Ga}{\Gamma}
\newcommand{\la}{\lambda}
\newcommand{\La}{\Lambda}
\newcommand{\al}{\alpha}
\newcommand{\be}{\beta}

\def\Si{\mathfrak{S}}
\def\Ai{\mathfrak{A}}
\newcommand{\si}{\sigma}
\newcommand{\om}{\omega}
\newcommand{\Om}{\Omega}

\newcommand{\de}{\delta}
\newcommand{\De}{\Delta}
\newcommand{\ka}{\kappa}

\newcommand{\Irr}{{\mathrm {Irr}}}

\renewcommand{\Im}{{\mathrm {Im}}}
\newcommand{\Ind}{{\mathrm {Ind}}}

\newcommand{\Mor}{{\mathrm {Mor}}}

\newcommand{\supp}{\operatorname{supp}}

\def\Ab{{\tt Ab}}

\newcommand{\tr}{{\mathrm {tr}}}

\newcommand{\pr}{{\mathrm {pr}}}
\newcommand{\rad}{{\mathrm {rad}\,}}

\newcommand{\Res}{{\mathrm {Res}}}

\newcommand{\C}{{\mathbb C}}
\newcommand{\Q}{{\mathbb Q}}
\newcommand{\R}{{\mathbb R}}

\newcommand{\Cent}{Z}

\newcommand{\ZZ}{{\mathbb Z}}

\newcommand{\D}{{\mathscr D}}

\renewcommand{\mod}{\bmod \,}

\newcommand{\Zig}{{\mathsf{A}}}
\newcommand{\Zag}{{\mathsf{B}}}

\newcommand{\ttA}{{\tt A}}
\newcommand{\ttB}{{\tt B}}
\newcommand{\ttC}{{\tt C}}
\newcommand{\ttf}{{\tt f}}

\newcommand{\tty}{{\tt y}}

\def\bM{\text{\boldmath$M$}}

\def\ttB{{\mathtt B}}

\def\col{{\operatorname{col}}}

\newcommand{\Std}{\operatorname{Std}}

\newcommand{\Add}{\operatorname{Ad}}
\newcommand{\Rem}{\operatorname{Re}}
\newcommand{\PA}{\operatorname{PAd}}
\newcommand{\PR}{\operatorname{PRe}}

\def\g{{\mathfrak g}}

\def\Par{{\mathscr P}}
\def\Cores{{\mathscr{C}}}

\def\umu{{\underline{\mu}}}

\def\ua{{\underline{a}}}
\def\ub{{\underline{b}}}
\def\um{{\underline{m}}}
\def\un{{\underline{n}}}

\def\b{\mathfrak{b}}
\def\k{\Bbbk}

\def\spa{\operatorname{span}}
\def\height{{\operatorname{ht}}}

\def\wt{{\operatorname{wt}}}

\def\op{{\mathrm{op}}}
\def\re{{\mathrm{re}}}
\def\im{{\mathrm{im}\,}}
\def\onto{{\twoheadrightarrow}}
\def\into{{\hookrightarrow}}
\def\Mod#1{#1\!\operatorname{-Mod}}
\def\mod#1{#1\!\operatorname{-mod}}

\def\proj#1{#1\!\operatorname{-proj}}

\def\underlineproj#1{#1\!\operatorname{-\underline{proj}}}

\def\iso{\stackrel{\sim}{\longrightarrow}}

\def\X{{\mathcal X}}
\def\Y{{\mathcal Y}}

\def\Blo{\cB}
\def\blo{\cb}
\def\BC{\mathcal{BC}}

\def\T{{\sf T}}
\def\Stab{{\sf S}}
\def\U{{\sf U}}
\def\V{{\sf V}}

\def\lan{\langle}
\def\ran{\rangle}

\def\Stand{\Delta}

\def\DIM{{\operatorname{dim}_q\,}}

\def\words{I}

\def\Seq{\operatorname{Seq}}

\def\I{{\mathcal I}}

\def\ttB{{\mathtt B}}

\def\ttc{{\mathtt c}}

\def\ttK{{\mathtt K}}

\def\bla{\text{\boldmath$\lambda$}}
\def\bmu{\text{\boldmath$\mu$}}
\def\bnu{\text{\boldmath$\nu$}}

\makeatletter
\newcommand*{\rom}[1]{\expandafter\@slowromancap\romannumeral #1@}
\makeatother

{\catcode`\|=\active
  \gdef\set#1{\mathinner{\lbrace\,{\mathcode`\|"8000%
  \let|\midvert #1}\,\rbrace}}
}
\def\midvert{\egroup\mid\bgroup}

{\catcode`\|=\active
  \gdef\set#1{\mathinner{\lbrace\,{\mathcode`\|"8000%
  \let|\midvert #1}\,\rbrace}}
}
\def\midvert{\egroup\mid\bgroup}

\colorlet{darkgreen}{green!50!black}
\tikzset{dots/.style={very thick,loosely dotted},
         greendot/.style={fill,circle,color=darkgreen,inner sep=1.5pt,outer sep=0},
         blackdot/.style={fill,circle,color=black,inner sep=1.1pt,outer sep=0},
         smallgreendot/.style={fill,circle,color=darkgreen,inner sep=0.9pt,outer sep=0},
         smallreddot/.style={fill,circle,color=red,inner sep=0.9pt,outer sep=0},
         smallbluedot/.style={fill,circle,color=blue,inner sep=0.9pt,outer sep=0},
         graydot/.style={fill,circle,color=gray,inner sep=1.1pt,outer sep=0},
         reddot/.style={fill,circle,color=red,inner sep=1.1pt,outer sep=0},
         bluedot/.style={fill,circle,color=blue,inner sep=1.1pt,outer sep=0}
}
\def\greendot(#1,#2){\node[greendot] at(#1,#2){}}
\def\blackdot(#1,#2){\node[blackdot] at(#1,#2){}}
\def\smallgreendot(#1,#2){\node[smallgreendot] at(#1,#2){}}
\def\smallreddot(#1,#2){\node[smallreddot] at(#1,#2){}}
\def\smallbluedot(#1,#2){\node[smallbluedot] at(#1,#2){}}
\def\graydot(#1,#2){\node[graydot] at(#1,#2){}}
\def\reddot(#1,#2){\node[reddot] at(#1,#2){}}
\def\bluedot(#1,#2){\node[bluedot] at(#1,#2){}}

\newenvironment{braid}{
  \begin{tikzpicture}[baseline=6mm,black,line width=.7pt, scale=0.32,
                      draw/.append style={rounded corners},
                      every node/.append style={font=\fontsize{5}{5}\selectfont}]%
  }{\end{tikzpicture}
}

\def\Grid(#1,#2){
  \draw[very thin,gray,step=2mm] (0,0)grid(#1,#2);
  \draw[very thin,darkgreen,step=10mm] (0,0)grid(#1,#2);
}
\newcommand{\braidbox}[6][0.2]{
	\draw[rounded corners, color=black] (#2-#1,#4) -- (#3+#1,#4) -- (#3+#1,#5) -- (#2-#1,#5) -- cycle;
	\draw ({(#2+#3)/2},{(#4+#5)/2}) node{\footnotesize #6};
}
\newcommand{\redbraidbox}[6][0.2]{
	\draw[rounded corners, color=red] (#2-#1,#4) -- (#3+#1,#4) -- (#3+#1,#5) -- (#2-#1,#5) -- cycle;
	\draw ({(#2+#3)/2},{(#4+#5)/2}) node{\footnotesize #6};
}

\newcommand{\darkgreenbraidbox}[6][0.2]{
	\draw[rounded corners, color=darkgreen] (#2-#1,#4) -- (#3+#1,#4) -- (#3+#1,#5) -- (#2-#1,#5) -- cycle;
	\draw ({(#2+#3)/2},{(#4+#5)/2}) node{\footnotesize #6};
}

\newcommand\Tableau[2][\relax]{
  \begin{tikzpicture}[scale=0.5,draw/.append style={thick,black}]
    \ifx\relax#1\relax%
    \else 
      \foreach\box in {#1} { \filldraw[blue!30]\box+(-.5,-.5)rectangle++(.5,.5); }
    \fi
    \newcount\row\newcount\col
    \row=0
    \foreach \Row in {#2} {
       \col=1
       \foreach\k in \Row {
          \draw(\the\col,\the\row)+(-.5,-.5)rectangle++(.5,.5);
          \draw(\the\col,\the\row)node{\k};
          \global\advance\col by 1
       }
       \global\advance\row by -1
    }
  \end{tikzpicture}
}

\newcommand\YoungDiagram[2][\relax]{
  \begin{tikzpicture}[scale=0.5,draw/.append style={thick,black}]
    \ifx\relax#1\relax%
    \else 
    \foreach\box in {#1} {
      \filldraw[blue!30]\box rectangle ++(1,1);
    }
    \fi
    \newcount\row
    \row=0
    \foreach \col in {#2} {
       \draw(1,\the\row)grid ++(\col,1);
       \global\advance\row by -1
    }
  \end{tikzpicture}
}

\newenvironment{Young}{\begingroup
       \def\vr{\vrule height0.89\hoogte width\dikte depth 0.2\hoogte}
       \def\fbox##1{\vbox{\offinterlineskip
                    \hrule height\dikte
                    \hbox to \breedte{\vr\hfill##1\hfill\vr}
                    \hrule height\dikte}}
       \vbox\bgroup \offinterlineskip \tabskip=-\dikte \lineskip=-\dikte
            \halign\bgroup &\fbox{##\unskip}\unskip  \crcr }
       {\egroup\egroup\endgroup}
\def\diagram#1{\relax\ifmmode\vcenter{\,\begin{Young}#1\end{Young}\,}\else%
              $\vcenter{\,\begin{Young}#1\end{Young}\,}$\fi}

\makeindex

\begin{document}

\frontmatter

\title[RoCK blocks for double covers of symmetric groups and quiver Hecke superalgebras]{{\bf RoCK blocks for double covers of symmetric groups and quiver Hecke superalgebras}}

\author{\sc Alexander Kleshchev}
\address{Department of Mathematics\\ University of Oregon\\ Eugene\\ OR 97403, USA}
\email{klesh@uoregon.edu}

\author{\sc Michael Livesey}
\address{School of Mathematics, University of Manchester, Manchester, M139PL, U.K.}
\email{michael.livesey@manchester.ac.uk}


\thanks{
The first author was supported by the NSF grant DMS-2101791 and Charles Simonyi Endowment at the Institute for Advanced Study. The second author was supported by the EPSRC (grant no EP/T004606/1). }

\begin{abstract}
We define and study RoCK blocks for double covers of symmetric groups. We prove that RoCK blocks of double covers are Morita equivalent to standard `local' blocks. The analogous result for blocks of symmetric groups, a theorem of Chuang and Kessar, was an important step in Chuang and Rouquier ultimately proving Brou\'e's abelian defect group conjecture for symmetric groups. Indeed we prove Brou\'e's conjecture for the RoCK blocks defined in this article. Our methods involve translation into the quiver Hecke superalgebras setting using the Kang-Kashiwara-Tsuchioka isomorphism and categorification methods of Kang-Kashiwara-Oh. As a consequence we construct Morita equivalences between more general objects than blocks of finite groups. In particular, our results extend to certain blocks of level one cyclotomic Hecke-Clifford superalgebras. We also study imaginary cuspidal categories of quiver Hecke superalgebras and classify  irreducible representations of quiver Hecke superalgebras in terms of cuspidal systems. 
\end{abstract}

\maketitle

\setcounter{page}{4}

\tableofcontents

\chapter{Introduction}
\section{Brou\'e's conjecture for symmetric groups and beyond}
\label{SI1}

A fundamental conjecture in representation theory of finite groups is:

\vspace{2mm}
\noindent
{\bf Brou\'e's abelian defect group conjecture. }\index{Brou\'e's abelian defect group conjecture}
{\em 
For any finite group $G$, a block $B$ of $G$ with abelian defect group $D$ is derived equivalent to its Brauer correspondent~$b$.}

\vspace{2mm}
\noindent

\subsection{Brou\'e's conjecture for symmetric groups}
\label{SSIRock}
For symmetric groups this conjecture is a theorem. The proof has two major steps. First, Chuang and Kessar \cite{CK} proved Brou\'e's conjecture for the remarkable special blocks of symmetric groups---the so-called RoCK blocks (so named after Rouquier, Chuang and Kessar); then Chuang and Rouquier \cite{CR} extended the result to all blocks by constructing derived equivalences between the blocks of symmetric groups with isomorphic defect groups.

To be more precise, let $\F$ be a field of characteristic $p>0$ and $\Si_n$ denote the symmetric group on $n$ letters. The blocks of the symmetric group algebra $\F\Si_n$ are labeled by pairs $(\rho,d)$, where $\rho$ is a $p$-core partition and $d$ is a non-negative integer such that $|\rho|+dp=n$ (for a partition $\la$ of $r$, we always write $|\la|$ for $r$). We denote by  $B^{\rho,d}$\index{b@$B^{\rho,d}$} the block of $\F\Si_n$ corresponding to such a pair $(\rho,d)$. The integer $d$ is referred to as the {\em weight} of the block $B^{\rho,d}$.  

Let $\rho$ be a $p$-core, $d\in\Z_{\geq 0}$ and $n=|\rho|+pd$. The defect group $D^{\rho,d}$\index{d@$D^{\rho,d}$} of $B^{\rho,d}$ is abelian if and only if $d<p$. Moreover, in this case $D^{\rho,d}\cong C_p^{\times d}$ and 
$$
N_{\Si_n}(D^{\rho,d})\cong \Si_{|\rho|}\times\big((C_p\rtimes C_{p-1})\wr \Si_d\big).
$$ 
The Brauer correspondent $b^{\rho,d}$\index{b@$b^{\rho,d}$} of $B^{\rho,d}$ is then the block of 
$N_{\Si_n}(D^{\rho,d})$ 
isomorphic to 
$B^{\rho,0}\otimes \big(\F(C_p\rtimes C_{p-1})\wr \Si_d\big)$, which is Morita equivalent to $\F(C_p\rtimes C_{p-1})\wr \Si_d$
since $B^{\rho,0}$ is a matrix algebra. 

In fact, $\F(C_p\rtimes C_{p-1})= b^{\varnothing,1}$, and a special case of Brou\'e's conjecture is that $b^{\varnothing,1}$ is derived equivalent to $B^{\varnothing,1}$,  
where $B^{\varnothing,1}$ is the principal block of $\F\Si_p$. 
This special case follows from a result of Rickard \cite[Theorem 4.2]{Rickard_1}. Now using a theorem of  Marcus \cite{Mar}, we have that $b^{\rho,d}$ derived equivalent to the wreath product $B^{\varnothing,1}\wr \Si_d$. Thus, Brou\'e's conjecture for symmetric groups is equivalent to the assertion that, under the assumption $d<p$, the block $B^{\rho,d}$ is derived equivalent to $B^{\varnothing,1}\wr \Si_d$. 

For every $d\in\Z_{\geq 0}$, Rouquier \cite{RoTh} defined special $p$-cores, referred to as {\em $d$-Rouquier cores},\index{r@$d$-Rouquier $p$-core} which are generic in a certain (combinatorial or Lie theoretic) sense, and conjectured that these should play a special role in representation theory of symmetric groups. The blocks $B^{\rho,d}$, with $\rho$ being $d$-Rouquier, are called {\em RoCK blocks}. \index{RoCK block}

\vspace{2mm}
\noindent
{\bf Theorem (Chuang-Kessar).} 
{\em 
Let $d<p$ and the block $B^{\rho,d}$ be RoCK. Then $B^{\rho,d}$ is {\em Morita}\, equivalent to $B^{\varnothing,1}\wr \Si_d$. 
}

\vspace{2mm}
\noindent

In particular, this proves Brou\'e's conjecture for RoCK blocks (with abelian defect). 
Brou\'e's conjecture for general blocks (with abelian defect) now comes from: 

\vspace{1mm}
\noindent
{\bf Theorem (Chuang-Rouquier).} 
{\em 
Arbitrary blocks $B^{\rho,d}$ and $B^{\rho',d'}$ of symmetric groups are derived equivalent if and only if $d=d'$.
}

\vspace{2mm}
\noindent

To summarize, the scheme of the proof of Brou\'e's conjecture is as follows. Let $d<p$ and $B^{\rho',d}$ be any block of weight $d$. Choose $\rho$ to be a $d$-Rouquier core. Then

\vspace{.1mm}
\begin{align}
b^{\rho',d}\,
&\cong \,\,B^{\rho',0}\otimes \big(b^{\varnothing,1}\wr \Si_d\big) 
\label{E111}
\\
&\sim_\Mor\,\,b^{\varnothing,1}\wr \Si_d&(\text{since $B^{\rho',0}$ is a matrix algebra})
\label{E112}
\\
&\sim_{\der}\,\,B^{\varnothing,1}\wr \Si_d
&(\text{by Rickard and Marcus})
\label{E113}
\\
&\sim_\Mor\,\,B^{\rho,d}&(\text{by Chuang-Kessar})
\label{E114}
\\
&\sim_{\der}\,\,B^{\rho',d}&(\text{by Chuang-Rouquier}).
\label{E115}
\end{align}

\vspace{4mm}

\subsection{Beyond abelian defect for blocks of symmetric groups}
\label{SSBeyondS}
Let us now drop the assumption $d<p$, in particular we do not suppose anymore that the defect group of a block is abelian. 

By Chuang-Rouquier, any block of weight $d$ is still derived equivalent to a RoCK block $B^{\rho,d}$. For any RoCK block $B^{\rho,d}$, 
Evseev \cite[Theorem 1.1]{Evseev} constructed an explicit idempotent $\mathsf{f}_{\rho,d}$ such that  
$$\mathsf{f}_{\rho,d}B^{\rho,d}\mathsf{f}_{\rho,d}\cong B^{\varnothing,1}\wr \Si_d.$$ 
Moreover, the basic algebra of $B^{\varnothing,1}$ is the Brauer tree algebra $\mathsf{Br}(p)$ corresponding to a linear tree with $p$ vertices. So the idempotent $\mathsf{f}_{\rho,d}$ can be refined to an idempotent $\mathsf{f}'_{\rho,d}$ such that $\mathsf{f}'_{\rho,d}B^{\rho,d}\mathsf{f}'_{\rho,d}\cong \mathsf{Br}(p)\wr \Si_d$. Thus we get (without assuming $d<p$):

\vspace{2mm}
\noindent
{\bf Theorem (Evseev).} 
{\em 
Suppose that $B^{\rho,d}$ is a RoCK block. Then there exists an idempotent $\mathsf{f}$ such that $\mathsf{f}B^{\rho,d}\mathsf{f}\cong \mathsf{Br}(p)\wr \Si_d$.
}

\vspace{2mm}
\noindent

However, for $d\geq p$, the algebra $B^{\rho,d}$ has more irreducible modules than the algebra $\mathsf{Br}(p)\wr \Si_d$, so they are certainly neither Morita nor derived equivalent. Thus we need to `upgrade' the local object $\mathsf{Br}(p)\wr \Si_d$ to something more subtle. 
Such an upgrade was suggested by Turner in \cite{Turner}. 

Turner defined new algebras $T^{\mathsf{Br}(p)}(d,d)$\index{t@$T^{\mathsf{Br}(p)}(d,d)$} 
such that $eT^{\mathsf{Br}(p)}(d,d)e\cong \mathsf{Br}(p)\wr \Si_d$
for some idempotent $e\in T^{\mathsf{Br}(p)}(d,d)$, but $T^{\mathsf{Br}(p)}(d,d)$ is Morita equivalent to $\mathsf{Br}(p)\wr \Si_d$ if and only if $d<p$. Turner conjectured that {\em any}\, RoCK block $B^{\rho,d}$ is Morita equivalent to $T^{\mathsf{Br}(p)}(d,d)$, hence any weight $d$ block of symmetric groups is derived equivalent to $T^{\mathsf{Br}(p)}(d,d)$, see Turner \cite[Conjecture 165]{Turner}. Thus, $T^{\mathsf{Br}(p)}(d,d)$ was conjectured to be the sought after `local object' which should replace wreath products of Brauer tree algebras in the context of `Brou\'e's conjecture' for blocks of symmetric groups with non-abelian defect groups.

The definition of Turner's algebra $T^{\mathsf{Br}(p)}(d,d)$ is rather subtle. In the first approximation, it can be thought of as Schur's  algebra of $\mathsf{Br}(p)\wr \Si_d$. However, a rigorous definition of $T^{\mathsf{Br}(p)}(d,d)$ involves working over integers first and choosing some explicit maximal rank $\Z$-subalgebra in the endomorphism algebra of `permutation modules' over $\mathsf{Br}(p)\wr \Si_d$. Alternatively, to define $T^{\mathsf{Br}(p)}(d,d)$ one can use Turner's double construction (although it is rather difficult to work with explicitly), see \cite{EK1} for more details on this. 

Turner's conjecture was proved in \cite{EK2}. The main result of \cite{EK2} can be interpreted as a construction of an explicit idempotent $\ga_{\rho,d}$ with the following properties:

\vspace{2mm}
\noindent
{\bf Theorem (Evseev-Kleshchev).} 
{\em 
Suppose that $B^{\rho,d}$ is a RoCK block. Then there exists an idempotent $\ga_{\rho,d}\in B^{\rho,d}$ such that 
$$
B^{\rho,d}\sim_\Mor \ga_{\rho,d}B^{\rho,d}\ga_{\rho,d}\cong T^{\mathsf{Br}(p)}(d,d).$$
}


\subsection{Khovanov-Lauda-Rouquier algebras}
The techniques of \cite{EK2} rely crucially on the theory of {\em Khovanov-Lauda-Rouquier (KLR) algebras}, also known as quiver Hecke algebras, which categorify the basic module over the quantum group $U_q(\g)$ for $\g$ of affine type $A_\ell^{(1)}$. Here $\ell=p-1$, where $p=\cha \F$, but in fact we do not need to suppose that $\ell+1$ is prime. 

To be more precise, by \cite{BKyoung} and \cite{RoKLR}, there is an algebra isomorphism $$B^{\rho,d}\cong R^{\La_0}_{\cont(\rho)+d\de},$$ 
where $R^{\La_0}_{\cont(\rho)+d\de}$ is the cyclotomic KLR algebra\index{cyclotomic KLR algebra} of Lie type $A_{p-1}^{(1)}$ corresponding to the fundamental dominant weight $\La_0$ and the element $\cont(\rho)+d\de$ of the non-negative part of the root lattice $Q_+$, where $\cont(\rho)\in Q_+$ is the residue content of the core partition $\rho$, and $\de$ is the null-root. 

Let $\ell\in\Z_{\geq 1}$ and $R^{\La_0}_{\cont(\rho)+d\de}$ be the cyclotomic KLR algebra of Lie type $A_{\ell}^{(1)}$. If $\rho$ is a $d$-Rouquier $(\ell+1)$-core we say that  $R^{\La_0}_{\cont(\rho)+d\de}$ is a RoCK block. Let $R^{\La_0}_{\cont(\rho)+d\de}$ be a RoCK block. The main result of  \cite{EK2} is then a construction of an idempotent $\ga_{\rho,d}\in R^{\La_0}_{\cont(\rho)+d\de}$ such that 
$$
R^{\La_0}_{\cont(\rho)+d\de}\sim_\Mor\ga_{\rho,d}R^{\La_0}_{\cont(\rho)+d\de}\ga_{\rho,d}\cong T^{\mathsf{Br}(\ell+1)}(d,d).
$$

The result of \cite{EK2} described in the previous paragraph is proved over any field, in fact even over the integers. So it covers not only blocks of symmetric groups but also blocks of Iwahori-Hecke algebras over arbitrary fields with parameter a root of unity. Note also that the KLR algebras $R^{\La_0}_{\cont(\rho)+d\de}$ come with $\Z$-grading, and Turner's algebra $T^{\mathsf{Br}(\ell+1)}(d,d)$ also inherits a natural $\Z$-grading from $\mathsf{Br}(\ell+1)$. The isomorphism $\ga_{\rho,d}R^{\La_0}_{\cont(\rho)+d\de}\ga_{\rho,d}\cong T^{\mathsf{Br}(\ell+1)}(d,d)$ and the Morita equivalence $R^{\La_0}_{\cont(\rho)+d\de}\sim_\Mor\ga_{\rho,d}R^{\La_0}_{\cont(\rho)+d\de}\ga_{\rho,d}$ established in \cite{EK2} are actually {\em graded}. Moreover, gradings play a crucial role in the proofs. 

The cyclotomic KLR algebra $R^{\La_0}_\theta$ is a quotient of the (infinite dimensional) KLR algebra $R_\theta$. The category of finitely generated graded $R_\theta$-modules is stratified with strata built out of the  categories of {\em cuspidal} $R_\theta$-modules, see \cite{KMJAlg}. The cuspidal categories correspond to the positive roots of affine Kac-Moody Lie algebra $\g$. The real cuspidal categories (i.e. the ones corresponding to the real positive roots) are well-understood, while the {\em imaginary cuspidal categories} (the ones corresponding to the imaginary roots $d\de$) are rather mysterious, see for example \cite{KM, MM}. However, it is these imaginary cuspidal categories that are directly related to RoCK blocks. 

The {\em imaginary cuspidal algebra $\bar R_{d\de}$} is the quotient of $R_{d\de}$ by the ideal generated by the standard idempotents $e(\bi)$ corresponding to the non-cuspidal words $\bi$. The imaginary cuspidal category can be described as the category of $R_{d\de}$-modules factoring through the surjection $R_{d\de}\,\onto\,\bar R_{d\de}$. 

The homomorphism $\Om:R_{d\de}\to R^{\La_0}_{\cont(\rho)+d\de}$ defined as the composition of natural maps
$$
R_{d\de}\hookrightarrow R_{\cont(\rho)}\otimes R_{d\de}\hookrightarrow R_{\cont(\rho)+d\de}\,\onto\, R^{\La_0}_{\cont(\rho)+d\de},
$$
factors through the surjection $R_{d\de}\,\onto\,\bar R_{d\de}$ if and only if $R^{\La_0}_{\cont(\rho)+d\de}$ is RoCK. Moreover, in that case  $\Om(R_{d\de})$ is Morita equivalent to $R^{\La_0}_{\cont(\rho)+d\de}$ and can be further truncated, preserving the Morita equivalence class, to get the Turner's algebra $T^{\mathsf{Br}(\ell+1)}(d,d)$. In this way we get a connection between RoCK blocks, imaginary cuspidal categories for KLR algebras, and Turner's algebras, which is crucial for understanding all three, see \cite{EK2}.

\section{RoCK blocks for double covers of symmetric groups}
Let $\F$ be an algebraically closed field of odd characteristic $p>0$. 
\subsection{Spin blocks of symmetric and alternating groups}
\label{SSSpinBlocks}

We denote by $\tilde\Si_n$ one of the two double covers of the symmetric group $\Si_n$, and by $\tilde\Ai_n$ the double cover of the alternating group $\Ai_n$, see \S\ref{SDoubleCovers}. 

We have a canonical central element $z\in\tilde\Si_n$ and a central idempotent $e_z:=(1-z)/2\in\F\tilde\Si_n.$ \index{e@$e_{z}$} 
This yields an ideal decomposition 
$$\F\tilde\Si_n=\F\tilde\Si_ne_z\oplus \F\tilde\Si_n(1-e_z).$$ 
So the blocks of $\F\tilde\Si_n$ split into blocks of $\F\tilde\Si_ne_z$ and the blocks of $\F\tilde\Si_n(1-e_z)$. Note that $\F\tilde\Si_n(1-e_z)\cong\F\Si_n$, so the latter blocks are relatively well understood, see \S\ref{SI1}. From now on we concentrate on the blocks of $\F\tilde\Si_n e_z$ sometimes referred to as the {\em spin blocks}\index{spin block} of the symmetric group. Little is known about these spin blocks. In particular, Brou\'e's conjecture for them is wide open.

In \S\ref{sec:twist_group} we introduce the algebra $\cT_n$, a twisted group algebra of the symmetric group with basis $\{\ct_w\mid w\in \Si_n\}$. In fact, $\cT_n \cong \F\tilde\Si_n e_z$ and it is often important to consider $\cT_n$ as a superalgebra with 
$$(\cT_n)_\0\cong\F\Ai_ne_z$$ 
spanned by the $\ct_w$ for even $w$, and $(\cT_n)_\1$ spanned by the $\ct_w$ for odd $w$. We can also speak of {\em superblocks} of $\cT_n$, i.e. indecomposable superideals of $\cT_n$. It is well-known (see \S\ref{SSB}) that superblocks are indecomposable as usual ideals (i.e. superblocks are blocks) unless we deal with trivial defect blocks. Trivial defect blocks are simple algebras so they are completely elementary, in particular Brou\'e's conjecture for them is trivial. So working with superblocks vs. blocks does not create any serious problems.

Superblocks of $\cT_n$ are labeled by pairs $(\rho,d)$, where $\rho$ is a $\bar p$-core partition and $d$ is a non-negative integer such that $|\rho|+dp=n$, see \S\ref{SSB}. We denote by  $\Blo^{\rho,d}$\index{b@$\Blo^{\rho,d}$} the block of $\cT_n$ corresponding to such a pair $(\rho,d)$. The integer $d$ is referred to as the {\em weight} of the block $\Blo^{\rho,d}$. We note that the case $d=0$ corresponds to the defect zero situations, so we will assume that $d>0$ throughout this section. The even part $\Blo^{\rho,d}_\0$ is then a spin block of the alternating group algebra $\F\tilde\Ai_n$.

We denote by $D^{\rho,d}$ a defect group of $\Blo^{\rho,d}$; it is known that $D^{\rho,d}$ is also a defect group of $\Blo^{\rho,d}_\0$, see \S\ref{SSB}. The defect group $D^{\rho,d}$ is abelian if and only if $d<p$. Let $\blo^{\rho,d}$\index{b@$\blo^{\rho,d}$} be the Brauer correspondent block of $\Blo^{\rho,d}$ and $\blo^{\rho,d}_\0$\index{b@$\blo^{\rho,d}_\0$} be the Brauer correspondent block of $\Blo^{\rho,d}_\0$ (note that $\blo^{\rho,d}_\0$ is indeed the even part of the superalgebra $\blo^{\rho,d}$), see \S\ref{SSB} for more details. Thus, $\blo^{\rho,d}$ is a block of $\F N_{\tilde{\Si}_n}(D^{\rho,d})$ and $\blo^{\rho,d}_\0$ is a block of $\F N_{\tilde{\Ai}_n}(D^{\rho,d})$. 

Here and everywhere we denote by `$\otimes$' a {\em tensor product of superalgebras}. For a superalgebra $A$ a {\em twisted wreath superproduct}\, $A \swr \cT_d$ is defined in \S\ref{sec:twist_group}. Now the analogue of (\ref{E111}) is given by the following theorem, which comes from  Proposition~\ref{prop:Brauer_corr}.

\vspace{2mm}
\noindent
{\bf Theorem A.} 
{\em 
Let $0<d<p$ and $\rho$ be any $\bar p$-core. Then 
$$
\blo^{\rho,d}\cong  \Blo^{\rho,0}\otimes (\blo^{\varnothing,1}\swr \cT_d)
$$
and
$$
\blo^{\rho,d}_\0\cong  \big(\Blo^{\rho,0}\otimes (\blo^{\varnothing,1}\swr \cT_d)\big)_\0.
$$
}
\vspace{2mm}

For a partition $\la$ we denote by $h(\la)$ the number of non-zero parts of $\la$. 
The analogue of (\ref{E112}) is given by the following theorem, which comes easily from Theorem A (see (\ref{mor_even_brauer}) and (\ref{mor_odd_brauer})): 

\vspace{2mm}
\noindent
{\bf Theorem B.} 
{\em 
Let $0<d<p$ and $\rho$ be any $\bar p$-core. Then 
$$
\blo^{\rho,d}
\sim_{\Mor}
\begin{cases}
\blo^{\varnothing,1} \swr \cT_d&\text{if }|\rho|-h(\rho)\text{ is even,}\\
(\blo^{\varnothing,1} \swr \cT_d)_{\0}&\text{if }|\rho|-h(\rho)\text{ is odd,}
\end{cases}
$$
and
$$
\blo^{\rho,d}_\0
\sim_{\Mor}
\begin{cases}
(\blo^{\varnothing,1} \swr \cT_d)_{\0}&\text{if }|\rho|-h(\rho)\text{ is even,}\\
\blo^{\varnothing,1} \swr \cT_d&\text{if }|\rho|-h(\rho)\text{ is odd.}
\end{cases}
$$
}
\vspace{2mm}

The analogue of the Rickard and Marcus step (\ref{E113}) is now given by the following theorem, which is Proposition~\ref{lem:Brauer_der}:

\vspace{2mm}
\noindent
{\bf Theorem C.} 
{\em 
Let $0<d<p$ and $\rho$ be any $\bar p$-core. Then 
$$
\blo^{\rho,d}
\sim_{\der}
\begin{cases}
\Blo^{\varnothing,1} \swr \cT_d&\text{if }|\rho|-h(\rho)\text{ is even,}\\
(\Blo^{\varnothing,1} \swr \cT_d)_{\0}&\text{if }|\rho|-h(\rho)\text{ is odd,}
\end{cases}
$$
and
$$
\blo^{\rho,d}_\0
\sim_{\der}
\begin{cases}
(\Blo^{\varnothing,1} \swr \cT_d)_{\0}&\text{if }|\rho|-h(\rho)\text{ is even,}\\
\Blo^{\varnothing,1} \swr \cT_d&\text{if }|\rho|-h(\rho)\text{ is odd.}
\end{cases}
$$
}
\vspace{2mm}

While Theorems A and B are established using standard techniques in modular representation theory and can be considered as a folklore, Theorem C requires more work since it involves extension of the results of Marcus from wreath products to twisted wreath products. 

\subsection{RoCK spin blocks}
We now discuss the spin analogue of the Chuang-Kessar step (\ref{E114}). 
In Section~\ref{SRock}, we define the notion of a {\em RoCK spin block}. As in the classical case described in \S\ref{SSIRock}, for every $d\in\Z_{\geq 0}$, we define {\em $d$-Rouquier $\bar p$-cores}\, using $\bar p$-abaci (see \S\ref{SSRoCore}). Then, if $\rho$ is a $d$-Rouquier $\bar p$-core, the spin block $\Blo^{\rho,d}$ is called {\em RoCK}. The corresponding spin block $\Blo^{\rho,d}_\0$ of $\F\tilde\Ai_n$ is also called {\em RoCK}. 

One of the main goals of this paper is to establish nice properties of RoCK spin blocks analogous to the properties of classical RoCK blocks established in \cite{CK} and \cite{Evseev}, and to prove Brou\'e's conjecture for RoCK spin blocks of symmetric and alternating groups. In particular, we have the following analogue of the Chuang-Kessar step (\ref{E114}):


\vspace{2mm}
\noindent
{\bf Theorem D.} 
{\em 
Let $0<d<p$ and $\Blo^{\rho,d}$ be a RoCK block. 
Then 
\begin{align*}
\Blo^{\rho,d} \sim_{\Mor}
\begin{cases}
\Blo^{\varnothing,1} \swr \cT_d&\text{if }|\rho|-h(\rho)\text{ is even,}\\
(\Blo^{\varnothing,1} \swr \cT_d)_{\0}&\text{if }|\rho|-h(\rho)\text{ is odd,}
\end{cases}
\end{align*}
and
\begin{align*}
\Blo^{\rho,d}_\0 \sim_{\Mor}
\begin{cases}
(\Blo^{\varnothing,1} \swr \cT_d)_{\0}&\text{if }|\rho|-h(\rho)\text{ is even,}\\
\Blo^{\varnothing,1} \swr \cT_d&\text{if }|\rho|-h(\rho)\text{ is odd.}
\end{cases}
\end{align*}
}
\vspace{2mm}

Theorems D and C have the following immediate corollary:

\vspace{2mm}
\noindent
{\bf Theorem E.} 
{\em 
Brou\'e's abelian defect group conjecture holds for RoCK spin blocks of symmetric and alternating groups.
}
\vspace{2mm}

Now, to deduce Brou\'e's conjecture for all spin blocks of symmetric and alternating groups it would be sufficient to prove the following spin analogue of Chuang-Rouquier, the so-called Kessar-Schaps conjecture, cf. \cite{AriSch}:

\vspace{2mm}
\noindent
{\bf Conjecture (Kessar-Schaps).} 
{\em 
Arbitrary spin blocks $\Blo^{\rho,d}$ and $\Blo^{\rho',d'}$ of symmetric groups (or the spin blocks $\Blo^{\rho,d}_{\0}$ and $\Blo^{\rho',d'}_{\0}$ of alternating groups) are derived equivalent if and only if $d=d'$ and $|\rho|-h(\rho)\equiv |\rho'|-h(\rho')\pmod{2}$. Moreover, $\Blo^{\rho,d}$ and $\Blo^{\rho',d'}_{\0}$ are derived equivalent if and only if $d=d'$ and $|\rho|-h(\rho)\not\equiv |\rho'|-h(\rho')\pmod{2}$.
}
\vspace{2mm}

The conjecture is known to hold if one replaces `derived equivalent' with the weaker `perfectly isometric', cf. Brunat-Gramain \cite{BruGra}. Indeed, the papers \cite{BruGra} and \cite{Liv} prove that Brou\'e's perfect isometries conjecture holds for spin blocks of symmetric and alternating groups. There have also been a number of results proving the Kessar-Schaps conjecture for specific pairs of blocks by showing they are even Morita equivalent, cf. Kessar \cite{Ke}, Kessar-Schaps \cite{KesSch}, Arisha-Schaps \cite{AriSch} and Leabovich-Schaps\cite{LeaSch}. Recently, proofs of the Kessar-Schaps Conjecture were announced in \cite{ELV,BKBroue}.

As demonstrated in \cite{FKM}, Theorem D can be used to provide complete information on decomposition numbers for RoCK spin blocks of abelian defect, analogous to the well-known results for classical RoCK blocks \cite{CT,LM,CT2,JLM}. We point out that little is known about decomposition numbers of double covers of symmetric groups in general, see for example \cite{Wales} for basic and second basic spin modules, \cite{MY2,BMO,Maas} for some small rank cases, \cite{LT} for a related general conjecture, \cite[\S10]{BK2} for some general facts on spin decomposition numbers and connections to the supergroup $Q(n)$, \cite{BKReg} for some leading term decomposition numbers, \cite{Mu} for weight $1$ spin blocks, and \cite{Fayers} for weight $2$ spin blocks.

Theorem D is Corollary~\ref{C021221_4}, which is just a  `desuperization' of the stronger Proposition~\ref{P021221_3} describing the Morita {\em superequivalence} class of the RoCK block $\Blo^{\rho,d}$ in the case $0<d<p$ as follows: 
\begin{align*}
\Blo^{\rho,d} \sim_{\sM}
\begin{cases}
\Blo^{\varnothing,1} \swr \cT_d&\text{if }|\rho|-h(\rho)\text{ is even,}\\
(\Blo^{\varnothing,1} \swr \cT_d)\otimes \cC_1&\text{if }|\rho|-h(\rho)\text{ is odd.}
\end{cases}
\end{align*}
Here and below, we denote by $\cC_m$ the rank $m$ {\em Clifford superalgebra}, and Morita superequivalence is a stronger notion than Morita equivalence reviewed in  Section~\ref{sec:grdd_supalg}, see especially \S\ref{sec:more_Morita}. 

\subsection{Sergeev and related superalgebras}
\label{SSSergI}
It is well-known that Clifford superalgebras play an important role in the theory of spin representations of symmetric groups. For a superalgebra $A$, we have that $A\otimes \cC_m$ is Morita superequivalent to $A$ if $m$ is even and $A\otimes \cC_m$ is Morita superequivalent to $A\otimes \cC_1$ if $m$ is odd, see Lemma~\ref{lem:clif_mat1}. We also have an algebra isomorphism $(A \otimes \cC_1)_{\0}\cong A$, and, under a reasonable additional assumption on $A$, a Morita equivalence $A\otimes \cC_1 \sim_{\Mor} A_{\0}$, see Lemma~\ref{lem:clif_mat2}. Finally, $\cC_m\otimes\cC_l\cong \cC_{m+l}$. 
So it is harmless for us to (super)tensor any superalgebra of interest with a Clifford superalgebra of arbitrary rank. 

In particular, it is often convenient to consider the {\em Sergeev superalgebra}\index{Sergeev superalgebra}
$$
\cY_n:=\cT_n\otimes \cC_n.\index{t@$\cY_n$}
$$
and its `superblocks' 
$$
\BC^{\rho,d}:=\Blo^{\rho,d}\otimes \cC_n,\index{b@$\BC^{\rho,d}$}
$$
so that
\begin{equation}\label{E201221}
\cY_n=\bigoplus \BC^{\rho,d},
\end{equation}
where the sum is over all $\bar p$-cores $\rho$ and non-negative integers $d$ such that $|\rho|+dp=n$. 

The Sergeev superalgebra $\cY_n$ is often more convenient to work with than the twisted group algebra $\cT_n$; in fact, it has the following three important advantages: 

\vspace{.8mm}
\noindent
{\sf (1)} $\cY_n$ has a useful affinization $\cX_n$ called the {\em affine Sergeev superalgebra}. 
Moreover, $\cX_n$ has a family of finite dimensional quotients $\cX_n^\La$ labeled by the integral dominant weights $\La$ of the affine root system of type $A_{p-1}^{(2)}$. Taking $\La$ to be the fundamental dominant weight $\La_0$ recovers $\cY_n$ as one of these quotients:
$
\cY_n\cong\cX_n^{\La_0}.
$
We refer the reader to \S\ref{SSS} for more details on this.

\vspace{.8mm}
\noindent
{\sf (2)} 
The superalgebras $\cX_n$ and $\cX_n^\La$ can be quantized to give a family of superalgebras $\cX_n(q)$ and $\cX_n^\La(q)$ for any parameter $q\in\F^\times$. The case $q=1$ is interpreted as $\cX_n(1)=\cX_n$ and $\cX_n^\La(1)=\cX_n^\La$. In particular, $\cX_n^{\La_0}(q)$ is isomorphic to {\em Olshanski's superalgebra}\, $\cY_n(q)$ which is a well-known quantization of $\cY_n\cong\cX_n^{\La_0}(1)$. This  quantization is based on Sergeev's  isomorphism $\cY_n\cong \cC_1\swr\Si_n$ and the standard Iwahori-Hecke quantization of the group algebra of $\Si_n$. We again refer the reader to \S\ref{SSS} for more details on this.

\vspace{.8mm}
\noindent
{\sf (3)} The algebras $\cX_n^\La(q)$ and their superblocks are related to {\em quiver Hecke superalgebras} via the Kang-Kashiwara-Tsuchioka isomorphism. 
We refer the reader to \S\ref{SSKKT} for more details on this.

\vspace{.8mm}
Point {\sf (3)} is a fundamental connection which places the study of spin block of symmetric groups into a much broader context of quiver Hecke superalgebras and categorification. Powerful methods which come from this broader context are key to our work. This will be reviewed in detail in the next section. 

We now give some additional comments on point {\sf (2)} referring the reader to \S\ref{SSS} for more details. 
To avoid trivial `semisimple situations' and stay in the modular representation theory context, we always suppose that we are in one of the following two interesting cases:
\begin{enumerate}
\item[{\sf (G)}] $\F$ is a field of odd characteristic $p$ written as $p=2\ell+1$ and $q=1$ (the `group case' considered above);
\item[{\sf (M)}] $\F$ is a field of any characteristic not equal to $2$ and $q\in \F$ is a primitive $(2\ell+1)$st root of unity for some $\ell\in\Z_{>0}$ (the `mixed case'). 
\end{enumerate}
Note that in the case {\sf (G)}, we have that $2\ell+1$ is the prime number $p=\cha\F$. In the case {\sf (M)}, we still set $p:=2\ell+1$, but now $p$ can be any odd number $\geq 3$, which is prime to $\cha\F$ if $\cha \F>0$. 

The combinatorics of $\bar p$-cores and $\bar p$-abaci does not depend on $p$ being prime, and we have all the combinatorial notions defined above, including $d$-Rouquier cores and RoCK blocks for this more general case. We have the `superblock decomposition' 
$$
\cY_n(q)=\bigoplus \BC^{\rho,d}(q)
$$
where the sum is over all $\bar p$-cores $\rho$ and non-negative integers $d$ such that $|\rho|+dp=n$. When $q=1$, i.e. we are in case {\sf (G)}, this recovers (\ref{E201221}). We call $\BC^{\rho,d}(q)$ a {\em RoCK block}\, if $\rho$ is a $d$-Rouquier $\bar p$-core, i.e. $\BC^{\rho,d}(q)$ is RoCK if and only if the corresponding block $\Blo^{\rho,d}$ of $\cT_n$ is RoCK. 

\subsection{Beyond abelian defect for spin blocks}
In \S\ref{SSZig}, we introduce a basic algebra $\Zig_\ell$\index{a@$\Zig_\ell$} which is defined as the path algebra of the quiver 
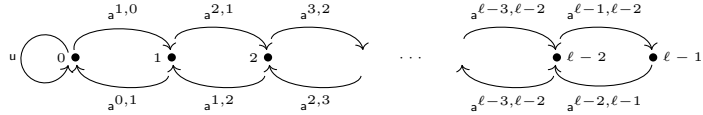
\begin{align*}
\begin{braid}\tikzset{baseline=3mm}
\coordinate (1) at (0,0);
\coordinate (2) at (4,0);
\coordinate (3) at (8,0);
\coordinate (4) at (12,0);
\coordinate (6) at (16,0);
\coordinate (L1) at (20,0);
\coordinate (L) at (24,0);
\draw [thin, black, ->] (-0.3,0.2) arc (15:345:1cm);
\draw [thin, black,->,shorten <= 0.1cm, shorten >= 0.1cm]   (1) to[distance=1.5cm,out=100, in=100] (2);
\draw [thin,black,->,shorten <= 0.25cm, shorten >= 0.1cm]   (2) to[distance=1.5cm,out=-100, in=-80] (1);
\draw [thin,black,->,shorten <= 0.25cm, shorten >= 0.1cm]   (2) to[distance=1.5cm,out=80, in=100] (3);
\draw [thin,black,->,shorten <= 0.25cm, shorten >= 0.1cm]   (3) to[distance=1.5cm,out=-100, in=-80] (2);
\draw [thin,black,->,shorten <= 0.25cm, shorten >= 0.1cm]   (3) to[distance=1.5cm,out=80, in=100] (4);
\draw [thin,black,->,shorten <= 0.25cm, shorten >= 0.1cm]   (4) to[distance=1.5cm,out=-100, in=-80] (3);
\draw [thin,black,->,shorten <= 0.25cm, shorten >= 0.1cm]   (6) to[distance=1.5cm,out=80, in=100] (L1);
\draw [thin,black,->,shorten <= 0.25cm, shorten >= 0.1cm]   (L1) to[distance=1.5cm,out=-100, in=-80] (6);
\draw [thin,black,->,shorten <= 0.25cm, shorten >= 0.1cm]   (L1) to[distance=1.5cm,out=80, in=100] (L);
\draw [thin,black,->,shorten <= 0.1cm, shorten >= 0.1cm]   (L) to[distance=1.5cm,out=-100, in=-100] (L1);
\blackdot(0,0);
\blackdot(4,0);
\blackdot(8,0);
\blackdot(20,0);
\blackdot(24,0);
\draw(0,0) node[left]{$0$};
\draw(4,0) node[left]{$1$};
\draw(8,0) node[left]{$2$};
\draw(14,0) node {$\cdots$};
\draw(20,0) node[right]{$\ell-2$};
\draw(24,0) node[right]{$\ell-1$};
 \draw(-2.6,0) node{$\zu$};
\draw(2,1.2) node[above]{$\za^{1,0}$};
\draw(6,1.2) node[above]{$\za^{2,1}$};
\draw(10,1.2) node[above]{$\za^{3,2}$};
\draw(18,1.2) node[above]{$\za^{\ell-3,\ell-2}$};
\draw(22,1.2) node[above]{$\za^{\ell-1,\ell-2}$};
\draw(2,-1.2) node[below]{$\za^{0,1}$};
\draw(6,-1.2) node[below]{$\za^{1,2}$};
\draw(10,-1.2) node[below]{$\za^{2,3}$};
\draw(18,-1.2) node[below]{$\za^{\ell-3,\ell-2}$};
\draw(22,-1.2) node[below]{$\za^{\ell-2,\ell-1}$};
\end{braid}
\end{align*}
modulo the following relations:
\begin{enumerate}
\item all paths of length three or greater are zero;
\item all paths of length two that are not cycles are zero;
\item the length-two cycles based at the vertex $i\in\{1,\dots,\ell-1\}$ are equal;
\item $\zu^2=\za^{0,1}\za^{1,0}$. 
\end{enumerate}
As an algebra, $\Zig_\ell$ is a special Brauer tree algebra, but it is important to consider $\Zig_\ell$ as a {\em superalgebra} with $|\ze^i|=|\za^{i,j}|=\0$ and $|\zu|=\1$ (where $\ze^j$ denotes the length zero path from $j$ to $j$). 

The role of the algebra $\Zig_\ell$ is explained by the fact that in the `group case' {\sf (G)}  
we have 
$$
\Blo^{\varnothing,1}\sim_{\sM}\Zig_{\ell}\otimes\cC_1\qquad\text{and}\qquad
\BC^{\varnothing,1}=\Blo^{\varnothing,1}\otimes\cC_p\sim_{\sM}\Zig_{\ell},
$$
see Lemmas~~\ref{lem:zigzag} and \ref{lem:sym_S}. But as is clear from Theorem F and Conjecture 1 below, the algebra $\Zig_\ell$ plays the same role for both the `group case' {\sf (G)} and the `mixed case' {\sf (M)}. 

Let $r_0$ denote the number of nodes of $\bar p$-residue $0$ in the Young diagram of $\rho$, see \S\ref{SSGen}. We note that $r_0\equiv h(\rho)\pmod{2}$. 
For a superalgebra $A$, a {\em wreath superproduct}\, $A \swr \Si_d$ is defined in \S\ref{SSBasicRep}.  
Now the spin analogue of Evseev's result \cite[Theorem 1.1]{Evseev} cited above is (see Theorem~\ref{T231221}):

\vspace{2mm}
\noindent
{\bf Theorem F.} 
{\em 
Suppose that we are in one of the cases {\sf (G)} or {\sf (M)}, and $\BC^{\rho,d}(q)$ is a RoCK block. Then there exists an idempotent $f\in\BC^{\rho,d}(q)$ such that 
$$
f\BC^{\rho,d}(q)f\cong (\Zig_\ell\swr\Si_d)\otimes \cC_{r_0+2d}.
$$ 
}
\vspace{2mm}


As in the classical situation of \S\ref{SSBeyondS}, we have that $\BC^{\rho,d}(q)\sim_{\sM}f\BC^{\rho,d}(q)f$ if and only if $d<p$, and so to give a local description of a RoCK block $\BC^{\rho,d}(q)$ for $d\geq p$ one needs an `upgrade' of the wreath superproduct $\Zig_\ell\swr\Si_d$. 
Let $T^{\Zig_\ell}_{\a_\ell}(d,d)$ be the generalized Schur algebra of \cite[\S3B]{KM3} corresponding to the choice of the subalgebra
$$
\a_\ell:= (\Zig_\ell)_\0=
\spa(\ze_0,\ze_0,\dots,\ze_{\ell-1},\za^{0,1},\za^{1,2},\dots,\za^{\ell-2,\ell-1})\subseteq \Zig_\ell.
$$
The following is an analogue of Turner's conjecture for symmetric groups discussed in \S\ref{SSBeyondS}:

\vspace{2mm}
\noindent
{\bf Conjecture 1.} 
{\em 
Suppose that we are in one of the cases {\sf (G)} or {\sf (M)}, and let  $\BC^{\rho,d}(q)$ be a RoCK block. 
Then there exists an idempotent $f\in\BC^{\rho,d}(q)$ such that 
$$
\BC^{\rho,d}(q)\sim_{\sM} f\BC^{\rho,d}(q)f\cong T^{\Zig_\ell}_{\a_\ell}(d,d)\otimes \cC_{r_0+2d}.
$$ 
}
\vspace{2mm}

\section{Quiver Hecke superalgebras}
\label{SIntro3}

\subsection{Kang-Kashiwara-Tsuchioka isomorphism}

To prove Theorems~D and F, we use the Kang-Kashiwara-Tsuchioka isomorphism to translate the problem into the language of {\em quiver Hecke superalgebras}. Quiver Hecke superalgebras are non-trivial superanalogues of Khovanov-Lauda-Rouquier algebras defined in \cite{KKT} and used to categorify quantum Kac-Moody algebras and superalgebras, cf. \cite{KKO,KKOII} and \cite{BE2}.

The Kang-Kashiwara-Tsuchioka isomorphism relates spin blocks and certain cyclotomic quiver Hecke superalgebras, which allows us to deduce Theorems~D and F from the corresponding results  on these cyclotomic quiver Hecke superalgebras. The proofs use categorification techniques of \cite{KKO} and the theory of imaginary cuspidal representations for quiver Hecke superalgebras, which is a superanalogue of \cite{Kcusp}. 

To be more precise, let $\g$ be the Kac-Moody Lie algebra of type $A_{2\ell}^{(2)}$. We refer the reader to \S\ref{ChBasicNotLie} for the corresponding standard Lie theoretic notation. In particular, the vertices of the Dynkin diagram of $\g$ are labeled by the elements of the set $I=\{0,1,\dots,\ell\}$, $\de$ is the null-root,  
$\{\al_i\mid i\in I\}$ are the simple roots, $\{\La_i\mid i\in I\}$ are the corresponding fundamental dominant weights, $Q_+$ denotes the non-negative part of the root lattice, and $P_+$ denotes the set of integral dominant weights. 

To every element $\theta=\sum_{i\in I}m_i\al_i\in Q_+$, following \cite{KKT}, we associate a (graded) superalgebra $\cR_\theta$ called a {\em quiver Hecke superalgebra} and an auxiliary superalgebra $\cRC_\theta$ called a {\em quiver Hecke-Clifford superalgebra}. These algebras also have {\em cyclotomic quotients}  $\cR_\theta^\La$ and $\cRC_\theta^\La$ corresponding to each $\La\in P_+$. 

We now review the connections between the quiver Hecke superalgebras and blocks of the superalgebras $\cX_n^\La(q)$ established in \cite{KKT}, for more details see \S\ref{SSKKT}. 
First, it is pointed out in \cite{KKT} that there is an idempotent $e_I\in \cRC_\theta$ such that 
$$
\cRC_\theta\sim_{\sM} e_I \cRC_\theta e_I\cong \cR_\theta\otimes \cC_{m_0}
$$
and
$$
\cRC_\theta^\La\sim_{\sM} e_I \cRC_\theta^\La e_I\cong \cR_\theta^\La\otimes \cC_{m_0}.
$$
In particular, the superalgebras $\cRC_\theta^\La$ and $\cR_\theta^\La$ are Morita superequivalent up to tensoring with a Clifford superalgebra. 

On the other hand, Kang, Kashiwara and Tsuchioka \cite{KKT} establish a non-trivial isomorphism of superalgebras 
$$
\cRC^\La_\theta\cong \X_n^\La(q)e_\theta,
$$
where 
$e_\theta\in \cX_n^\La(q)$ is a certain explicit central idempotent. Note that the choice of the parameter $q$ is always made as in {\sf (G)} or {\sf (M)} in \S\ref{SSSergI}.

In the most important special case where $q=1$ and $\La=\La_0$, this specializes to an isomorphism 
$$
\cRC^{\La_0}_{\cont(\rho)+d\de}\cong \Blo^{\rho,d}\otimes \cC_{|\rho|+dp},
$$
where 
$\cont(\rho)\in Q_+$ is the residue content of $\rho$, see \S\ref{ChBasicNotLie} and (\ref{SEResCont}). 
In particular, up to tensoring with Clifford superalgebras, 
the spin block $\Blo^{\rho,d}$ of a symmetric group is Morita superequivalent to 
the quiver Hecke superalgebra $\cR_{\cont(\rho)+d\de}^{\La_0}$. This motivates our study of the representation theory of $\cR_\theta^{\La_0}$. 

In view of the Kang-Kashiwara-Oh categorification theorem \cite{KKO}, we have $\cR^{\La_0}_\theta\neq 0$ if and only if $\theta$ is of the form $\cont(\rho)+d\de$ for some $\bar p$-core $\rho$ and $d\in\Z_{\geq 0}$. Now the Kessar-Schaps conjecture follows from the following:

\vspace{2mm}
\noindent
{\bf Conjecture 2.}
{\em Let $\rho$ and $\rho'$ be $\bar p$-cores, and $d,d'\in\Z_{\geq 0}$. Then, if $d=d'$, the algebras $\cR^{\La_0}_{\cont(\rho)+d\de}$ and $\cR^{\La_0}_{\cont(\rho')+d'\de}$ are derived equivalent and the algebras $\cR^{\La_0}_{\cont(\rho)+d\de} \otimes \cC_1$ and $\cR^{\La_0}_{\cont(\rho')+d'\de} \otimes \cC_1$ are also derived equivalent.
}
\vspace{2mm}


\subsection{RoCK quiver Hecke superalgebras}
We call $\cR^{\La_0}_{\cont(\rho)+d\de}$ {\em  RoCK} if $\rho$ is a $d$-Rouquier $\bar p$-core, i.e. $\cR^{\La_0}_{\cont(\rho)+d\de}$ is RoCK if and only if the {\em spin}\, block $\Blo^{\rho,d}$ of the symmetric group $\Si_{|\rho|+dp}$ is RoCK.  Now Conjecture 1 follows easily from the following: 

\vspace{2mm}
\noindent
{\bf Conjecture 3.} 
{\em 
If $\cR^{\La_0}_{\cont(\rho)+d\de}$ is RoCK, then there exists an idempotent $\ga\in \cR^{\La_0}_{\cont(\rho)+d\de}$ such that 
$$
\cR^{\La_0}_{\cont(\rho)+d\de}\sim_{\sM} \ga \cR^{\La_0}_{\cont(\rho)+d\de}\ga\cong T^{\Zig_\ell}_{\a_\ell}(d,d).
$$ 
}
\vspace{2mm}

In Chapter~\ref{ChRockqHs}, which is a central part of this work, we undertake a detailed study of RoCK quiver Hecke superalgebras $\cR^{\La_0}_{\cont(\rho)+d\de}$ and establish the following weak version of Conjecture 3: 

\vspace{2mm}
\noindent
{\bf Theorem G.} 
{\em 
Let $\cR^{\La_0}_{\cont(\rho)+d\de}$ be RoCK. Then there exists an idempotent $f\in \cR^{\La_0}_{\cont(\rho)+d\de}$ 
such that 
$$
f \cR^{\La_0}_{\cont(\rho)+d\de}f\cong \Zig_\ell\swr\Si_d.
$$ 
Moreover, $\cR^{\La_0}_{\cont(\rho)+d\de}\sim_{\sM}\Zig_\ell\swr\Si_d$ if and only if $d<p$. 
}
\vspace{2mm}

Theorems~D and F are eventually deduced from Theorem~G. 

The construction of the idempotent $f$ and the isomorphism 
$f \cR^{\La_0}_{\cont(\rho)+d\de}f\cong \Zig_\ell\swr\Si_d$ in Theorem G are  explicit. The idempotent $f$  is a version (and a superanalogue) of a so-called Gelfand-Graev idempotent considered in \cite[\S7.2]{KMimag} and \cite[\S6.1]{EK2}. The elements of $f \cR^{\La_0}_{\cont(\rho)+d\de}f$ corresponding to the standard generators $\ze^i$, $\zu$ and $\za^{i,j}$ of $\Zig_\ell$ are described explicitly by formulas (\ref{EEIFormula}), (\ref{EUFormula}), (\ref{EAII-1Formula}) and (\ref{EAI-1IFormula}). On the other hand, the elements of $f \cR^{\La_0}_{\cont(\rho)+d\de}f$ corresponding to the standard generators $s_r\in\Si_d$ are described explicitly using Kang-Kashiwara-Oh intertwiners, cf. \S\ref{SSKKO} and Section~\ref{SIntertwiners}. 

We now describe some remarkable algebraic properties of RoCK quiver Hecke algebras $\cR^{\La_0}_{\cont(\rho)+d\de}$ which in some sense explain their nice behavior, referring the reader to \S\ref{SSFurtherRoCK} for more details and generalizations. First of all, in general, there is an idempotent $e_{\cont(\rho),d\de}\in \cR^{\La_0}_{\cont(\rho)+d\de}$ and a superalgebra homomorphism 
$$
\Om:\cR_{d\de}\to e_{\cont(\rho),d\de}\cR^{\La_0}_{\cont(\rho)+d\de}e_{\cont(\rho),d\de},
$$
which arises from the composition of the following natural maps
$$
\cR_{d\de}\hookrightarrow \cR_{\cont(\rho)}\otimes \cR_{d\de}\hookrightarrow \cR_{\cont(\rho)+d\de}\,\onto\, \cR^{\La_0}_{\cont(\rho)+d\de},
$$
see (\ref{EOmega}). Denote $Z_{\rho,d}^{\La_0}:=\Im(\Om)$. 
If $\cR^{\La_0}_{\cont(\rho)+d\de}$ is RoCK, then
$$
e_{\cont(\rho),d\de}\cR^{\La_0}_{\cont(\rho)+d\de}e_{\cont(\rho),d\de}\cong \cR^{\La_0}_{\cont(\rho)}\otimes Z_{\rho,d}^{\La_0},
$$
$\cR^{\La_0}_{\cont(\rho)+d\de}$ is Morita superequivalent to $Z_{\rho,d}^{\La_0}$, and, most importantly, the homomorphism $\Om$ factors through the surjection $\cR_{d\de}\onto \bar\cR_{d\de}$, where $\bar\cR_{d\de}$ is the {\em imaginary cuspidal algebra}. 

The imaginary cuspidal algebra $\bar \cR_{d\de}$ is the quotient responsible for the category of imaginary cuspidal representations of $\cR_{d\de}$. 
We discuss cuspidality in the next subsection. 


\subsection{Cuspidal algebras for quiver Hecke superalgebras $\cR_\theta$}
We now discuss cuspidal algebras and a classification of irreducible $\cR_\theta$-supermodules via cuspidal systems. The results are analogous to \cite{Kcusp,KM}, see also \cite{McN1,TW}. We have (non-unital) parabolic subalgebras 
$$
\cR_{\theta_1}\otimes\dots\otimes \cR_{\theta_k}\cong \cR_{\theta_1,\dots,\theta_k}\subseteq \cR_{\theta_1+\dots+\theta_k}
$$ 
and the corresponding induction and restriction functors $\Ind_{\theta_1,\dots,\theta_k}$ and $\Res_{\theta_1,\dots,\theta_k}$. A family of standard modules induced from the cuspidal modules is defined, and irreducible $\cR_\theta$-supermodules arise as simple heads of the standard modules. 

To be more precise, we fix a {\em convex preorder} $\preceq$ on the set of positive roots $\Phi_+$ of $\g$, see \S\ref{SSConv} for more details. We point out that only a special convex order is important in connection to RoCK blocks, cf. \S\ref{SSRoCKCusp}, but at this stage we can work with an arbitrary convex preorder. 

We denote by $\Psi:=\Phi_+^\re\cup\{\de\}\subseteq \Phi_+$ the set of all indivisible positive roots. Let $\be\in\Psi$ and $m\in\Z_{>0}$. A finitely generated $\cR_{m\be}$-supermodule $M$ is called {\em cuspidal} if $\Res_{\theta,\eta}M\neq 0$ for $\theta,\eta\in Q_+$ implies that $\theta$ is a sum of positive roots $\preceq \be$ and $\eta$ is a sum of positive roots $\succeq\be$. 

There is a largest quotient $\bar \cR_{m\be}$ of $\cR_{m\be}$ such that an $\cR_{m\be}$-supermodule is cuspidal if and only if it factors through the surjection $\cR_{m\be}\,\onto\,\bar \cR_{m\be}$, and the category of cuspidal $\cR_{m\be}$-supermodules is nothing but the category of $\bar \cR_{m\be}$-supermodules, see \S\ref{SSCuspidalModules}. 
The algebra $\bar \cR_{m\be}$ is called a {\em cuspidal algebra}. 
For the case where $\be\in \Phi_+^\re$ we speak of real cuspidal algebras and modules, while for $\be=\de$ we speak of {\em imaginary cuspidal algebras} and modules. 

Denote by $\Par^\ell(m)$ the set of all $\ell$-multipartitions of $m$. 

\vspace{2mm}
\noindent
{\bf Theorem H.} 
{\em 
Let $\be\in \Psi$ and $m\in \Z_{>0}$.
\begin{enumerate}
\item[{\rm (i)}] If $\be\in \Phi_+^\re$, then there exists a unique (up to isomorphism) irreducible cuspidal $\cR_{m\be}$-supermodule $L_{m\be}$. 
\item[{\rm (ii)}] For every multipartition $\bmu\in\Par^\ell(m)$ there exists an  irreducible cuspidal $\cR_{m\de}$-supermodule $L_\bmu$, and $\{L_\bmu\mid \bmu\in\Par^\ell(m)\}$ is a complete irredundant set of irreducible cuspidal $\cR_{m\de}$-supermodules up to isomorphism.
\end{enumerate}
}

\vspace{2mm}

Let $\theta\in Q_+$. A {\em root partition of $\theta$} is a pair $(\um,\bmu)$, where $\um=(m_\be)_{\be\in\Psi}$ is a tuple of non-negative integers such that $\sum_{\be\in\Psi}m_\be\be=\theta$  and $\bmu\in\Par^\ell(m_\de)$. Denote by $\Par(\theta)$ the set of all root partitions of $\theta$. Given $(\um,\bmu)\in \Par(\theta)$, note that $m_\be=0$ for almost all $\be$, so we can choose a finite subset
$$
\be_1\succ\dots\succ\be_s\succ\de\succ\be_{-t}\succ\dots\succ\be_{-1}
$$
of $\Psi$ such that $m_\be=0$ for $\be$'s outside of this subset. Then, denoting $m_u:=m_{\be_u}$, we can write $(\um,\bmu)$ in the form
\begin{equation*}
(\um,\bmu)=(\be_1^{m_1},\dots,\be_s^{m_s},\bmu,\be_{-t}^{m_{-t}},\dots,\be_{-1}^{m_{-1}}).
\end{equation*}
Denote also
$$
\lan\um\ran:=(m_1\be_1,\dots,m_s\be_s,m_\de\de,m_{-t}\be_{-t},\dots,m_{-1}\be_{-1})\in Q_+^{s+t+1},
$$
so we have a parabolic subalgebra $\cR_{\lan\um\ran}\subseteq \cR_\theta$. 

\vspace{2mm}
\noindent
{\bf Theorem I.} 
{\em 
Let $\theta\in Q_+$. For every $(\um,\bmu)\in\Par(\theta)$, the induced module  
$$
\Ind_{\lan\um\ran}L_{m_1\be_1} \boxtimes \dots\boxtimes L_{m_s\be_s}\boxtimes L(\bmu)\boxtimes L_{m_{-t}\be_{-t}}\boxtimes\dots\boxtimes  L_{m_{-1}\be_{-1}}
$$
has an irreducible head denoted $L(\um,\bmu)$. Moreover, 
$$\{L(\um,\bmu)\mid (\um,\bmu)\in \Par(\theta)\}$$ is a complete irredundant set of irreducible $\cR_\theta$-supermodules up to isomorphism. 
}
\vspace{2mm}

Theorems~H and I have to be proved together. 

The real cuspidal algebra $\bar \cR_{d\be}$ has an easy structure and it is Morita equivalent to the algebra of symmetric polynomials in $d$ variables. On the other hand, the imaginary cuspidal algebras $\bar \cR_{d\de}$ is very interesting. In Theorems~\ref{TAffIso} and \ref{T221221}, we show that it is  related to the {\em affine Brauer tree algebra}\, $H_d(\Zig_\ell)$ defined in \S\ref{SSAff}:

\vspace{2mm}
\noindent
{\bf Theorem J.} 
{\em 
There exists an idempotent $\ga\in \bar \cR_{d\de}$ 
such that 
$
\ga \bar \cR_{d\de}\ga \cong H_d(\Zig_\ell).
$ 
Moreover, $\bar \cR_{d\de}\sim_{\sM}H_d(\Zig_\ell)$ if and only if $d<p$. 
}
\vspace{2mm}

Theorem J is closely connected to Theorem G and the two have to be proved together. The proof depends on Theorem~H. A crucial role in the proofs is played by the dimension formulas discussed in the next subsection.

\subsection{Graded dimension formula for $\cR^{N\La_0}_\theta$}
We complete this introduction with a discussion of another key tool, which is a combinatorial formula for the graded dimension of a general cyclotomic quotient of the form $\cR^{N\La_0}_\theta$. 

The {\em graded} superalgebra $\cR^{N\La_0}_\theta$ comes with a family of standard orthogonal idempotents $\{e(\bi)\mid \bi\in I^\theta\}$, see \S\ref{SQHDef}. In Theorem~\ref{TDim}, we prove:

\vspace{2mm}
\noindent
{\bf Theorem K.} 
{\em 
Let $N\in\Z_{> 0}$. For $\theta\in Q_+$ set $n=\height(\theta)$, and let $\bi,\bj\in I^\theta$. Then 
$$
\dim_qe(\bi)\cR^{N\La_0}_\theta e(\bj)=\sum_{\bla\in\Par_p^N(n)}\sum_{\substack{\Stab\in \Std_p(\bla,\bi)\\ \T\in\Std_p(\bla,\bj)}} \deg(\Stab)\deg(\T)\norm{\bla}.
$$
}
\vspace{2mm}

Here $\bla$ runs over the set of $p$-strict $N$-multipartitions of $n$, $\Std_p(\bla,\bi)$ denotes the set of $p$-standard tableaux of shape $\bla$ and type $\bi$, and $\deg(\T)$ is the degree of a tableau $\T$ defined in \S\ref{SRemAdd}. The dimension formula of Theorem~K is an analogue of the dimension formula from \cite[Theorem 4.20]{BKdec}. The proof follows the same idea but uses the combinatorics of the type $A_{2\ell}^{(2)}$ Fock space developed in \cite{KMPY} and the main categorification result of \cite{KKO}. 

Putting $q$ to $1$ and taking $N=1$ in Theorem~K yields the formula for the ungraded dimension $\dim e(\bi)\cR^{\La_0}_\theta e(\bj)$ obtained in \cite[Theorem 3.4]{AP}, see \S\ref{subsec:dim_form}. However, in this paper we do need to use the formula of Theorem K in full generality. 

As we were preparing this work for publication, the preprint \cite{HS} has appeared, which contains a proof of a {\em different}\, graded dimension formula for arbitrary cyclotomic quotients $\cR_\theta^\La$. The dimension formula of Theorem K is related to a basis of cellular type in our algebra, while the  Hu-Shi dimension formula is related to a basis of monomial type. 
We do not know if the Hu-Shi dimension formula can be used instead of Theorem K to perform the dimension calculations of \S\S\ref{SSGrDimYRho1},\,\ref{SSDimYRhoD}. 

\subsection*{Acknowledgement} This work is strongly influenced by the ideas of Anton Evseev, who sketched the $d=1$ case of Theorem \ref{TMorWr}. In particular, he suggested an explicit form of the isomorphism in Theorem~\ref{TMaind=1}.

We are also very grateful to the referee who clearly took a large amount of time and care to read the finer details of the article. We feel this significantly helped improve the quality of the paper.

\mainmatter

\chapter{Background material}

\section{Basic notation}
\label{ChBasicNot}

\subsection{Generalities}
\label{ChBasicNotGen}

\begin{enumerate}
\item[$\bullet$] $\F$ is an arbitrary ground field.\index{f@$\F$}

\item[$\bullet$] $q$ is an indeterminate. \index{q@$q$}

\item[$\bullet$] $\ell$\index{l@$\ell$} is a fixed positive integer and $p:=2\ell+1$.\index{p@$p$}

\item[$\bullet$] $I := \{0,1,\dots,\ell\}$\index{i@$I$}, $J:=\{0,1,\dots,\ell-1\}$\index{j@$J$} and $K:=\{0,1,\dots,\ell-2\}$\index{k@$K$}.
\end{enumerate}

\subsection{Partitions, compositions, symmetric group}
\label{ChBasicNotPar}
\begin{enumerate}
\item[$\bullet$] A {\em composition}\index{composition} of $d\in\Z_{\geq 0}$ is a tuple 
$\la=(\la_1,\dots,\la_n)$ of non-negative integers summing to $d$.

\item[$\bullet$]  A {\em partition}\index{partition} of $d\in\Z_{\geq 0}$ is an infinite tuple $(\la_1,\la_2,\dots)$ of non-negative integers such that $\la_1\geq \la_2\geq\dots$ and $\la_1+\la_2+\dots=d$. 

\item[$\bullet$] The trivial partition\index{trivial partition} $(0,0,\dots)$ of $0$ is denoted $\varnothing$. 

\item[$\bullet$] The {\em length} of a partition $\la$\index{length of a partition} is defined as 
$h(\la):=\max\{k\mid \la_k>0\}$ \index{h@$h(\la)$}
(by convention, $h(\varnothing)=0$). If $k\geq h(\la)$ we also write $\la=(\la_1,\dots,\la_k)$.  

\item[$\bullet$] We denote by $\Par(d)$\index{p@$\Par(d)$} the set of all partitions of $d$, and set $\Par:=\bigsqcup_{d\geq 0}\Par(d).$\index{p@$\Par$} For $\la\in\Par(d)$ we write $d=|\la|$.\index{$\lvert\la\rvert$} 

\item[$\bullet$] For $m\in\Z_{\geq 1}$, we have the set $\Par^m$\index{p@$\Par^m$} of $m$-tuples of partitions, referred to as {\em $m$-multipartitions}.\index{multipartition} Typically, we use the notation $\bla=(\la^{(1)},\dots,\la^{(m)})$ for elements of $\Par^m$. For $1\leq k\leq m$, we refer to $\la^{(k)}$ as the {\em $k$th component}\index{component}\, of the multipartition $\bla$. 

\item[$\bullet$] For  $\bla=(\la^{(1)},\dots,\la^{(m)})\in \Par^m$, we write $|\bla|:=|\la^{(1)}|+\dots+|\la^{(m)}|$,
\index{$\lvert\bla\rvert$}
 and for $d\in \Z_{\geq 0}$, denote $\Par^m(d):=\{\bla\in\Par^m\mid|\bla|=d\}$. \index{p@$\Par^m(d)$}

\item[$\bullet$] ${\Si}_d$\index{s@${\Si}_d$} is the {\em symmetric group}\index{symmetric group} on $\{1,\dots,d\}$ (by convention, $\Si_0$ is the trivial group), and $s_r:=(r,r+1) \in {\Si}_d$\index{s@$s_r$} for $r=1,\dots,d-1$ are the {\em elementary transpositions}. 

\item[$\bullet$] We denote by $\leq$ the Bruhat order\index{Bruhat order} on $\Si_d$, i.e. for $u,w\in \Si_d$, we have $u\leq w$ of and only if $u=s_{r_{a_1}}\cdots s_{r_{a_b}}$ for some $1\leq a_1<\cdots<a_b\leq l$, where $w=s_{r_1}\cdots s_{r_l}$  is a reduced decomposition for $w$. 

\item[$\bullet$] For a composition $\la=(\la_1,\dots,\la_n)$ of $d$, we have the standard {\em parabolic subgroup}\index{parabolic subgroup} ${\Si}_\la:= {\Si}_{\la_1}\times\dots\times {\Si}_{\la_n}\leq {\Si}_d$.\index{s@${\Si}_\la$}

\item[$\bullet$] We denote by $\D^\la$\index{d@$\D^\la$} (resp. ${}^\la\D$)\index{d@${}^\la\D$} the set of the minimal length coset representatives for ${\Si}_d/{\Si}_\la$ (resp. ${\Si}_\la\backslash {\Si}_d$).
\end{enumerate}

\subsection{Words and shuffles}
\label{ChBasicNotWords}
\begin{enumerate}
\item[$\bullet$]
For $d\in \Z_{\geq 0}$, we write elements of $I^d$ as $\bi=(i_1,\dots,i_d)=i_1\cdots i_d$ and refer to them as {\em words}.\index{words} The group $\Si_d$ acts on $I^d$ by place permutations. 

\item[$\bullet$]
 Given words $\bi\in I^d$ and $\bj\in I^e$, we have the concatenation $\bi\bj\in I^{d+e}$.  
 
\item[$\bullet$] If $\la=(\la_1,\dots,\la_n)$ is a composition of $d$ and $\bi^r\in I^{\la_r}$ for $r=1,\dots,n$, the words of the form $w\cdot(\bi^1\cdots \bi^n)$ with $w\in \D^\la$ are called {\em shuffles} of $\bi^1,\dots, \bi^n$. \index{shuffles}

\end{enumerate}

\begin{Example} \label{ExExtraWord}
{\rm 
Let $\bi=(\ell-1)(\ell-2)\cdots100$, $\bj=(\ell-1)(\ell-2)\cdots1$ and $\bl=(\ell-1)(\ell-2)\cdots10$. Note that any shuffle $\bm=m_1\cdots m_{2\ell}$ of $\bl$ with itself is a word 
which satisfies the following condition: for any $i=1,\dots,\ell-1$ and any $k=1,\dots, 2\ell$, the number of $i$'s among $m_1,\dots,m_k$ is not smaller than the number of $(i-1)$'s among $m_1,\dots,m_k$.
It follows that $\bm$ is also a shuffle of $\bi$ and $\bj$; indeed proceed along $\bm$ from left to right and take into $\bi$ the first $\ell-1$, then the first $\ell-2$, etc.  until the first $1$, and both $0$'s. The rest will be the word $\bj$. This manifests the fact that $\bm$ is a shuffle of $\bi$ and $\bj$. 
}
\end{Example}

\subsection{Lie theory}
\label{ChBasicNotLie}
\begin{enumerate}
\item[$\bullet$] $\g$ is the Kac-Moody Lie algebra of type $A_{2\ell}^{(2)}$,  
see \cite[Ch.4]{Kac}. \index{g@$\g$}\index{a@$A_{2\ell}^{(2)}$}

\item[$\bullet$] The Dynkin diagram of $\g$ has vertices labeled by $I= \{0,1,\dots,\ell\}$: 
$$
{\begin{picture}(330, 15)%
\put(6,5){\circle{4}}%
\put(101,2.45){$<$}%
\put(12,2.45){$<$}%
\put(236,2.42){$<$}%
\put(25, 5){\circle{4}}%
\put(44, 5){\circle{4}}%
\put(8, 4){\line(1, 0){15.5}}%
\put(8, 6){\line(1, 0){15.5}}%
\put(27, 5){\line(1, 0){15}}%
\put(46, 5){\line(1, 0){1}}%
\put(49, 5){\line(1, 0){1}}%
\put(52, 5){\line(1, 0){1}}%
\put(55, 5){\line(1, 0){1}}%
\put(58, 5){\line(1, 0){1}}%
\put(61, 5){\line(1, 0){1}}%
\put(64, 5){\line(1, 0){1}}%
\put(67, 5){\line(1, 0){1}}%
\put(70, 5){\line(1, 0){1}}%
\put(73, 5){\line(1, 0){1}}%
\put(76, 5){\circle{4}}%
\put(78, 5){\line(1, 0){15}}%
\put(95, 5){\circle{4}}%
\put(114,5){\circle{4}}%
\put(97, 4){\line(1, 0){15.5}}%
\put(97, 6){\line(1, 0){15.5}}%
\put(6, 11){\makebox(0, 0)[b]{$_{0}$}}%
\put(25, 11){\makebox(0, 0)[b]{$_{1}$}}%
\put(44, 11){\makebox(0, 0)[b]{$_{{2}}$}}%
\put(75, 11){\makebox(0, 0)[b]{$_{{\ell-2}}$}}%
\put(96, 11){\makebox(0, 0)[b]{$_{{\ell-1}}$}}%
\put(114, 11){\makebox(0, 0)[b]{$_{{\ell}}$}}%
\put(230,5){\circle{4}}%
\put(249,5){\circle{4}}%
\put(231.3,3.2){\line(1,0){16.6}}%
\put(232,4.4){\line(1,0){15.2}}%
\put(232,5.6){\line(1,0){15.2}}%
\put(231.3,6.8){\line(1,0){16.6}}%
\put(230, 11){\makebox(0, 0)[b]{$_{{0}}$}}%
\put(249, 11){\makebox(0, 0)[b]{$_{{1}}$}}%
\put(290,2){\makebox(0,0)[b]{\quad if $\ell = 1$.}}%
\put(175,2){\makebox(0,0)[b]{if $\ell \geq 2$, \qquad and}}%
\end{picture}}
$$

\item[$\bullet$] $P$ is the {\em weight lattice} of $\g$. \index{weight lattice}\index{p@$P$}

\item[$\bullet$] $\{\alpha_i \:|\:i \in I\} \subset P$ are the {\em simple roots} of $\g$.\index{simple roots}\index{a@$\alpha_i$}

\item[$\bullet$] $\{h_i\:|\:i \in I\} \subset P^*$ and the 
{\em simple coroots} of $\g$.\index{simple coroots}\index{h@$h_i$}

\item[$\bullet$] 
The {\em Cartan matrix} \index{Cartan matrix}
$\left(\langle h_i, \al_j\rangle\right)_{0\leq i,j\leq \ell}$ is:  
$$
\left(
\begin{matrix}
2 & -2 & 0 & \cdots & 0 & 0 & 0 \\
-1 & 2 & -1 & \cdots & 0 & 0 & 0 \\
0 & -1 & 2 & \cdots & 0 & 0 & 0 \\
 & & & \ddots & & & \\
0 & 0 & 0 & \dots & 2 & -1& 0 \\
0 & 0 & 0 & \dots & -1 & 2& -2 \\
0 & 0 & 0 & \dots & 0 & -1& 2 \\
\end{matrix}
\right)
\quad
\text{if\ \  $\ell\geq 2$,\ \  and}
\quad
\left(
\begin{matrix}
2 & -4 \\
-1 & 2
\end{matrix}
\right)
\quad \text{if\ \  $\ell=1$.}
$$

\item[$\bullet$]  $\{\Lambda_i\:|\:i \in I\} \subset P$ 
are the {\em fundamental dominant weights}, so 
$\langle h_i, \Lambda_j \rangle = \delta_{i,j}$.
\index{fundamental dominant weight}
\index{l@$\Lambda_i$}

\item[$\bullet$] $P_+\subset P$ is the set of all {\em dominant integral weights}.
\index{dominant integral weight}
\index{p@$P_+$}

\item[$\bullet$] $Q$ is the sublattice of $P$ generated by the simple roots. 
\index{q@$Q$}

\item[$\bullet$]
$Q_+:=\big\{\sum_{i\in I}m_i\al_i\mid m_i\in\Z_{\geq 0}\ \text{for all $i\in I$}\big\}\subset Q.
$ 
\index{q@$Q_+$}

\item[$\bullet$]
For $\theta=\sum_{i\in I}m_i\al_i\in Q_+$, its {\em height} is  
$
\height(\theta):=\sum_{i\in I}m_i
$.
\index{height}
\index{h@$\height(\theta)$}

\item[$\bullet$]  $(.|.)$\index{$(.\vert.)$} is a normalized invariant form on $P$ whose Gram matrix with respect to the linearly independent set $\al_0,\al_1,\dots,\al_\ell$ is: 
$$
\left(
\begin{matrix}
2 & -2 & 0 & \cdots & 0 & 0 & 0 \\
-2 & 4 & -2 & \cdots & 0 & 0 & 0 \\
0 & -2 & 4 & \cdots & 0 & 0 & 0 \\
 & & & \ddots & & & \\
0 & 0 & 0 & \dots & 4 & -2& 0 \\
0 & 0 & 0 & \dots & -2 & 4& -4 \\
0 & 0 & 0 & \dots & 0 & -4& 8 \\
\end{matrix}
\right)
\quad
\text{if $\ell\geq 2$, and}
\quad
\left(
\begin{matrix}
2 & -4 \\
-4 & 8
\end{matrix}
\right)
\quad \text{if $\ell=1$.}
$$
In particular, we have 
$(\al_0|\al_0)=2$, $(\al_\ell|\al_\ell)=8$, and $(\al_i|\al_i)=4$ for all other $i\in I$. We also have $(\al_i|\La_j)=\de_{i,j}(\al_i|\al_i)/2$. 


\item[$\bullet$] For all $i\in I$, we denote $\al_i^\vee:=2\al_i/(\al_i|\al_i)$. In particular, $(\al_i^\vee|\La_j)=\de_{i,j}$. \index{a@$\al_i^\vee$}

\item[$\bullet$] We denote 
$
\delta = \sum_{i=0}^{\ell - 1} 2 \alpha_i + \alpha_\ell.
$ \index{d@$\de$}
Then $\{\Lambda_0,\dots,\Lambda_\ell, \delta\}$ is a $\Z$-basis of $P$, cf. \cite[\S6.2]{Kac}. Note that $\height(\de)=p$. The coefficient of $\al_i$ in $\de$ is denoted $a_i$; thus $a_i=1$ if $i=\ell$ and $a_i=2$ otherwise.
\index{a@$a_i$}

\item[$\bullet$] We denote $\ttK:=h_0+2\sum_{i=1}^\ell h_i$. Then $\lan \ttK,\La_0\ran=1$ and $\lan \ttK,\al_i\ran=0$ for all $i\in I$, cf. \cite[\S6.2]{Kac}. The coefficient of $h_i$ in $\de$ is denoted $a_i^\vee$; thus $a_i^\vee=1$ if $i=0$ and $a_i^\vee=2$ otherwise.
\index{k@$\ttK$}\index{a@$a_i^\vee$}

\item[$\bullet$]
We denote for $i\in I$ and $n\in\Z_{\geq 0}$:
$$
q_i:=q^{(\al_i|\al_i)/2},\quad [n]_i:=\frac{q_i^{n}-q_i^{-n}}{q_i-q_i^{-1}},\quad [n]_i^!:=[1]_i[2]_i\cdots[n]_i. 
$$
\index{q@$q_i$}
\index{n@$[n]_i$}
\index{n@$[n]_i^{^^21}$}

\item[$\bullet$] For $\bi=i_1\cdots i_n\in I^n$, we denote 
$\|\bi\|:=\al_{i_1}+\dots+\al_{i_n}\in Q_+$.
\index{$\norm{\bi}$}

\item[$\bullet$]  For $\theta\in Q_+$, we set
$I^\theta:=\{\bi=i_1\cdots i_n\in I^n\mid \|\bi\|=\theta\}.$ 
\index{i@$I^\theta$}

\item[$\bullet$] We define $I^{\theta}_{\di}$ to be the set of all expressions of the form
$i_1^{(m_1)} \cdots i_r^{(m_r)}$ with
$m_1,\dots,m_r\in \Z_{\ge 1}$, $i_1,\dots,i_r\in I$
 and $m_1 \al_{i_1} + \cdots + m_r \al_{i_r} = \theta$ (note we do not insists that $i_k\neq i_{k+1}$). We refer to such expressions as {\em divided power words}. We identify $I^\theta$ with the subset of $I^\theta_\di$ which consists of all expressions as above with all $m_k=1$. We use the same notation for concatenation of divided power words as for concatenation of words.
 \index{i@$I^{\theta}_{\di}$}
 \index{divided power word}
 
 \item[$\bullet$] For $\bi=i_1^{(m_1)} \cdots i_r^{(m_r)}\in I^{\theta}_{\di}$, we denote 
$\hat\bi := i_1^{m_r}\cdots i_r^{m_r}\in I^{\theta}.$ 
\index{i@$\hat\bi$}\index{$\hat\bi$}

\end{enumerate}

\section{Graded superalgebras and supermodules}\label{sec:grdd_supalg}

Throughout the paper, unless otherwise stated, we work over the ground field~$\F$.


\subsection{Graded superalgebras}\label{SSBasicRep}
By a {\em graded superspace}\index{graded superspace} we understand a (vector) space over $\F$ with decomposition $$V=\bigoplus_{n\in\Z,\,\eps\in\Z/2}V_{n,\eps}.$$ We refer to the elements of $V_{n,\eps}$ as the (homogeneous) elements of {\em bidegree} $(n,\eps)$.\index{bidegree} We also refer to $n$ as {\em degree} and $\eps$ as {\em parity} and write for $v\in V_{n,\eps}$:
$$
\bideg(v)=(n,\eps),\quad \deg(v)=n,\quad |v|=\eps.\index{$\lvert v\rvert$}\index{d@$\deg(v)$}
$$
\index{b@$\bideg(\cdot)$}
A (homogeneous) {\em subsuperspace}\, $W\subseteq V$ is a subspace such that $$W=\sum_{n,\eps}(W\cap V_{n,\eps}).$$ 

For $n\in\Z$ we set $V_n:=V_{n,\0}\oplus V_{n,\1}$, thus considering $V$ as a graded  space $V=\bigoplus_{n\in\Z} V_n$. We also set $V_{\0}:=\bigoplus_{n\in\Z}V_{n,\0}$ and $V_{\1}:=\bigoplus_{n\in\Z}V_{n,\1},$ thus considering $V$ as a superspace $V=V_\0\oplus V_\1$. 
Conversely, a graded vector space is considered as a graded superspace concentrated in bidegrees $(*,\0)$, a superspace is considered as a graded superspace concentrated in bidegrees $(0,*)$, and a vector space is considered as a graded superspace concentrated in bidegree $(0,\0)$. 
For a graded vector space $V=\bigoplus_{n\in \Z} V_n$ with finite dimensional graded components, its {\em graded dimension}\index{graded dimension} is 
$$\dim_q V:=\sum_{n\in\Z}(\dim V_n)q^n\in\Z((q)).$$ 
\index{d@$\dim_q$}
For graded vector spaces $V,W$ with finite dimensional graded components, we write 
\begin{equation}\label{EGDinLeq}
\dim_q V\leq \dim_q W
\end{equation}
to indicate that $\dim V_n\leq \dim W_n$ for all $n$.

Let $V,W$ be graded superspaces. We define $V^*:= \bigoplus_{n\in\Z,\,\eps\in\Z/2}(V_{n,\eps})^*$, considered as a graded superspace via $(V^*)_{n,\eps}:=(V_{-n,\eps})^*$. 
The tensor product $V\otimes W$ is considered as a graded superspace via $\bideg(v\otimes w)=(\deg(v)+\deg(w),|v|+|w|)$ (here and below in places like this it is assumed that $v$ and $w$ are homogeneous). For $m\in\Z$ and $\de\in\Z/2$, a bidegree $(m,\de)$ (homogeneous) linear map $f:V\to W$ is a linear map satisfying $f(V_{n,\eps})\subseteq W_{n+m,\eps+\de}$ for all $n,\eps$.

Let $V$ be a superspace and $d\in\Z_{\geq 1}$.
The symmetric group $\Si_d$ acts on $V^{\otimes d}$ via
\begin{equation}\label{EWSignAction}
{}^{w}(v_1\otimes\dots\otimes v_d):=(-1)^{[w;v_1,\dots,v_d]}v_{w^{-1}(1)}\otimes\dots\otimes v_{w^{-1}(d)},
\end{equation}
where
$$
[w;v_1,\dots,v_d]:=\sum_{1\leq a<c\leq d, w(a)>w(c)}|v_a||v_c|.
$$
A {\em graded superalgebra}\index{graded superalgebra} is a graded superspace $A$ which is a unital algebra such that $A_{n,\eps}A_{m,\de}\subseteq A_{n+m,\eps+\de}$ for all $n,m,\eps,\de$. For example, given a finite dimensional graded superspace $V$, we have the graded superalgebra 
$$
\End_\F(V)=\bigoplus_{n\in\Z,\eps\in\Z/2}\End_\F(V)_{n,\eps},
$$
where $\End_\F(V)_{n,\eps}$ consists of all homogeneous linear maps of bidegree $(n,\eps)$. We refer to $\End_\F(V)$ as a {\em graded matrix superalgebra} or simply as a {\em matrix algebra}.

Let $A,B$ be graded superalgebras. An {\em isomorphism} $f:A \to B$ is an algebra isomorphism which is homogeneous of bidegree $(0,\0)$. 
The tensor product $A\otimes B$ of graded superspaces is considered as a graded superalgebra via 
\begin{align*}
(a\otimes b)(a'\otimes b') = (-1)^{|b| |a'|}aa'\otimes bb'.
\end{align*}
For $1\leq r<s\leq d$, we have maps 
\begin{align*}
\iota^d_{r}\index{i@$\iota^d_{r}$}&:A\to A^{\otimes d},\ a\mapsto 1^{\otimes (r-1)}\otimes a\otimes 1^{\otimes(d-r)},\\
\iota^d_{r,s}\index{i@$\iota^d_{r,s}$}&:A\otimes A\to A^{\otimes d},\ a\otimes b\mapsto 1^{\otimes (r-1)}\otimes a\otimes 1^{\otimes(s-r-1)}\otimes b\otimes 1^{\otimes(d-s)}.
\end{align*}
For $x\in A$ and $y\in A\otimes A$, we also denote 
\begin{equation}\label{EInsertion}
x_{r}:= \iota^d_r(x),\quad y_{r,s}:=\iota^d_{r,s}(y).
\end{equation}

We define $A^{\op}$, the {\em opposite graded superalgebra} of $A$, to be equal to $A$ as a graded superspace but with multiplication given by $a.b=(-1)^{|a||b|}ba$, for all $a,b\in A$.

An {\em antiautomorphism} of a graded superalgebra $A$ is a bidegree $(0,\0)$ linear map such that $\tau(ab)=\tau(b)\tau(a)$ (note there is no sign). 
We have an algebra involution 
\begin{equation}\label{ESi}
\si=\si_A:A\to A,\ a\mapsto (-1)^{|a|}a.\index{s@$\si_A$}
\end{equation}

Let $A$ be a graded superalgebra and $B\subseteq A$ be a (unital) graded subsuperalgebra. The 
{\em supercentralizer}\index{supercentralizer} $\Cent_A(B)$ of $B$ in $A$ is the 
graded subsuperalgebra of $A$ defined as the linear span of all (homogeneous) $a\in A$ such that $ba=(-1)^{|b||a|}ab$ for all (homogeneous) $b\in B$

\begin{Lemma}\label{lem:sup_cent}
Let $A$ be a graded superalgebra and $B$ a unital graded subsuperalgebra isomorphic to a graded matrix superalgebra. Then we have an isomorphism of graded superalgebras 
$
B \otimes \Cent_{A}(B) \iso A,\ b\otimes z \mapsto bz.
$
\end{Lemma}

\begin{proof}
This is an analogue of \cite[Proposition 4.10]{Evseev}, from which it easily follows. 
\end{proof}

Given a graded superalgebra $A$, we consider the {\em wreath superproduct} $A\swr \Si_d$,\index{wreath superproduct}\index{$\swr$} where $A\swr \Si_d=A^{\otimes d}\otimes \F\Si_d$ as graded superspaces, with $\F\Si_d$ concentrated in bidegree $(0,\0)$. To define the algebra structure,  recall (\ref{EWSignAction}). We identify $A^{\otimes d}$ and $\F\Si_d$ with subspaces of $A^{\otimes d}\otimes \F\Si_d$ in the obvious way, and postulate that $A^{\otimes d}$, $\F\Si_d$ are subalgebras of  $A\swr \Si_d$, and 
\begin{align*}
w\,(a_1\otimes\dots\otimes a_d)={}^w(a_1\otimes\dots\otimes a_d)\,w\qquad(w\in\Si_d,\ a_1,\dots, a_d\in A).
\end{align*}

Suppose now that $A$ is a finite dimensional graded superalgebra, and $m\in \Z$. An element $\tr\in A^*$ is called a {\em trace map of degree $-m$}\index{trace map} if $\tr(A_{n,\eps})=0$ for $(n,\eps)\neq (m,\0)$, $\tr(ab)=\tr(ba)$ for all $a,b \in A$ and the bilinear form $(a,b)_\tr:=\tr(ab)
$ is non-degenerate. If $A$ is equipped with a degree $-m$ trace map, we say \(A\) is {\em ($m$-)graded symmetric.}\index{graded symmetric superalgebra} In that case, pick a pair $\{a^1,\dots,a^k\}$, $\{b^1,\dots,b^k\}$ of dual bases with respect to $(\cdot,\cdot)_\tr$, and define the {\em distinguished element} of $A \otimes A$:
\begin{align*}
\nabla:= \sum_{l=1}^k a^l\otimes b^l\in A\otimes A.
\end{align*}
Note that $\bideg(\nabla)=(m,\0)$. 


If, in addition, $B$ is an $n$-graded symmetric superalgebra, then it is a standard check to verify that $A \otimes B$ is $(m+n)$-graded symmetric via
\begin{align}\label{tensor_symm_algs}
a \otimes b \mapsto \tr_A(a)\tr_B(b),
\end{align}
where $\tr_A$ and $\tr_B$ are the trace maps of $A$ and $B$ respectively.

\subsection{Graded supermodules}

Let $A$ be a graded superalgebra. 
A {\em graded $A$-supermodule}\index{graded supermodule} $M$ is an $A$-module which is also a graded superspace such that $A_{n,\eps}M_{m,\de}\subseteq M_{n+m,\eps+\de}$ for all $n,m,\eps,\de$. 
A graded {\em $(A,B)$-bisupermodule}\index{graded bisupermodule}  $M$ is an $(A,B)$-bimodule which is also a graded superspace such that $A_{n,\eps}M_{m,\de}\subseteq M_{n+m,\eps+\de}$ and $M_{m,\de}B_{n,\eps}\subseteq M_{n+m,\eps+\de}$ for all $n,m,\eps,\de$. We note that we can also view the graded $(A,B)$-bisupermodule $M$ as a graded $A \otimes B^{\op}$-supermodule via
\begin{equation}\label{bisupmodAB}
(a\otimes b).m=(-1)^{|b||m|}amb,
\end{equation}
for all $a\in A$, $b\in B$ and $m\in M$. In this way, we identify the notions of an $(A,B)$-bisupermodule and an $A \otimes B^{\op}$-supermodule.

A {\em homomorphism} of graded $(A,B)$-bisupermodules $f:V\to W$ is a linear map $f:V\to W$ satisfying $$f(avb)=(-1)^{|f||a|}af(v)b$$ for all (homogeneous) $a\in A,v\in V,b\in B$. In particular, a homomorphism of graded $A$-supermodules has the property $f(av)=(-1)^{|f||a|}af(v)$.

For a graded $(A,C)$-bisupermodule $M$ and a graded $(B,D)$-bisupermodule $N$, we define the graded $(A\otimes B,C\otimes D)$-bisupermodule $M\boxtimes N$ to be the graded superspace $M\otimes N$ with the action

\begin{align*}
(a\otimes b)(m\otimes n)&=(-1)^{|b||m|}(am\otimes bn),\\
(m\otimes n)(c\otimes d)&=(-1)^{|c||n|}(mc\otimes nd).
\end{align*}

In particular, given a graded $A$-supermodule $M$ and a graded $B$-supermodule $N$, we have the graded $(A\otimes B)$-supermodule $M\boxtimes N$ with the action 
$$
(a\otimes b)(m\otimes n)=(-1)^{|b||m|}(am\otimes bn).
$$

\begin{Remark}\label{rem:(bi)supmod}
Let $M$ be an $(A,C)$-bisupermodule and let $N$ be a $(B,D)$-bisupermodule. We can view $M$ as an $A \otimes C^{\op}$-supermodule and $N$ as a $B \otimes D^{\op}$-supermodule via (\ref{bisupmodAB}). The $(A \otimes C^{\op}) \otimes (B \otimes D^{\op})$-supermodule $M\boxtimes N$ can then be identified with the $(A\otimes B, C\otimes D)$-bisupermodule $M\boxtimes N$ using (\ref{bisupmodAB}) and the superalgebra isomorphism
\begin{align*}
(A \otimes C^{\op}) \otimes (B \otimes D^{\op}) &\to (A \otimes B) \otimes (C \otimes D)^{\op}\\
(a \otimes c) \otimes (b \otimes d) &\mapsto (-1)^{|b||c|}(a \otimes b) \otimes (c \otimes d).
\end{align*}
\end{Remark}

\vspace{1mm}
A graded $A$-supermodule $M$ is {\em irreducible}\index{irreducible graded supermodule}  if it has exactly two homogeneous submodules: $0$ and $M$. In general, an irreducible graded $A$-supermodule $M$ may or may not be irreducible when considered as a usual $A$-module, although for the graded superalgebras considered in Part~\ref{Part2}, the irreducible graded supermodules will have this property, see Lemma~\ref{LTypeM}.
We denote by $\Irr(A)$\index{i@$\Irr(A)$} a complete and irredundant set of irreducible graded $A$-supermodules, i.e. every irreducible graded $A$-supermodule is isomorphic to some member of $\Irr(A)$ and any two distinct members of $\Irr(A)$ are not isomorphic to each other.

Suppose in addition that $A$ is left Noetherian (as a graded superalgebra). We denote by 
$\mod{A}$\index{m@$\mod{A}$} the category of all finitely generated graded $A$-supermodules and homomorphisms defined as above. Then for $M,N\in \mod{A}$, we have 
$$\Hom_A(M,N)=\bigoplus_{n\in\Z,\,\eps\in\Z/2}\Hom_A(M,N)_{n,\eps}.$$
We denote by $\proj{A}$\index{p@$\proj{A}$} the full subcategory of finitely generated projective graded $A$-supermodules. 

In fact, $\mod{A}$ and $\proj{A}$ 
are {\em graded $(Q,\Pi)$-supercategories}\index{graded $(Q,\Pi)$-supercategory} in the sense of \cite[Definition 6.4]{BE}, with $\Pi$ the parity change functor\index{parity change functor}\index{p@$\Pi$} and $Q$ the degree shift functor\index{degree shift functor}\index{q@$Q$}:   
$$(\Pi M)_{n,\eps}=M_{n,\eps+\1}\quad \text{with the new action $a\cdot v=(-1)^{|a|}av$},$$ and $$(QM)_{n,\eps}=M_{n-1,\eps}.$$ Note that for any $m\in\Z$ and $\eps\in\Z/2$, the identity map considered as the map $M\to Q^m\Pi^\eps M$ has bidegree $(m,\eps)$ and is an isomorphism in $\mod{A}$. 
For an irreducible graded $A$-supermodule $L$ and a graded supermodule $M$ with a finite composition series (for example if $M$ is finite dimensional), we denote by $[M:L]$ the composition multiplicity of $L$ in $M$. As $L\cong Q^n\Pi^\eps L$, we of course have $[M:L]=[M:Q^n\Pi^\eps L]$ for all $n\in\Z$ and $\eps\in\Z/2$.


Let $\underlineproj{A}$\index{p@$\underlineproj{A}$} be the category whose objects are the same as those of $\proj{A}$ but morphisms are  homogeneous $A$-module homomorphisms of bidegree $(0,\0)$. This is a $(Q,\Pi)$-category in the sense of \cite[Definition 6.12(i)]{BE}. As in \cite[p.1083]{BE}, its Grothendieck group $[\underlineproj{A}]_{q}^\pi$\index{Grothendieck group}\index{$[\underlineproj{A}]_{q}^\pi$} is a $\Z^\pi[q,q^{-1}]$-module, where  
$$\Z^\pi:=\Z[\pi]/(\pi^2-1)\index{z@$\Z^\pi$},$$ 
with $\pi$ acting as $[\Pi]$ and $q$ acting as $[Q]$. Forgetting about $\Pi$, we get the {\em Grothendieck group} 
\begin{equation}\label{EGrothGroup}
[\underlineproj{A}]_q:=[\underlineproj{A}]^{\pi:=1}_q, 
\index{$[\underlineproj{A}]_q$}
\end{equation}
which is just a $\Z[q,q^{-1}]$-module. 
We extend scalars to $\Q(q)$ to get the $\Q(q)$-vector space  
 $$[\underlineproj{A}]_{\Q(q)}:=[\underlineproj{A}]_q\otimes_{\Z[q,q^{-1}]} \Q(q).\index{$[\underlineproj{A}]_{\Q(q)}$}
$$

\begin{Remark}\label{rem:non-sup_non-graded0}
Let $A$ be a graded superalgebra. We will often use the notation $|A|$\index{$\lvert A\rvert$} to denote the underlying $\Z$-graded algebra, where we forget about the superstructure. For example, if $B$ is another graded superalgebra and $M$ a graded $A$-supermodule, then $|A|\otimes |B|$ denotes the usual tensor product of the graded algebras $|A|$ and $|B|$ (i.e. no sign), $\End_{|A|}(M)$ is the usual endomorphism algebra of the $|A|$-module $M$ and $|A|^{\op}$ is the usual opposite algebra of $|A|$. However, it is important to note that, even though we define these algebraic objects without reference to their super structure, $|A|\otimes |B|$, $|A|^{\op}$ and $\End_{|A|}(M)$ are still graded superalgebras in the obvious way.
\end{Remark}

\begin{Remark}\label{rem:non-sup_non-graded}
We will also often have objects with superstructure but no $\Z$-grading such as superspaces, superalgebras and supermodules. We will always think of such objects as graded superspaces with  the $\Z$-grading concentrated in degree zero. In that way we can utilize all the notation and results from this section.
\end{Remark}



\subsection{Graded Morita superequivalence}
\label{SMoritaSuper}
Let $A$ and $B$ be graded superalgebras. A {\em graded Morita superequivalence}\index{graded Morita superequivalence} between $A$ and $B$ is a Morita equivalence between $A$ and $B$ induced by a graded $(A,B)$-bisupermodule $M$ and a graded $(B,A)$-bisupermodule $N$, i.e. $M\otimes_B N\cong A$ and $N\otimes_A M\cong B$ both bimodule isomorphisms of bidegree $(0,\0)$. Forgetting the $\Z$-gradings, we get a {\em Morita superequivalence}\index{Morita superequivalence} induced by bisupermodules $M,N$ as above with $M\otimes_B N\cong A$ and $N\otimes_A M\cong B$, both bimodule isomorphisms being even. We write $$A\sim_{\gsM}B\index{$\sim_{\gsM}$}$$ for graded Morita superequivalence $$A\sim_{\sM}B\index{$\sim_{\sM}$}$$ for Morita superequivalence. If $A$ and $B$ are arbitrary algebras, then we write $$A \sim_{\Mor} B\index{$\sim_{\Mor}$}$$ if $A$ and $B$ are Morita equivalent simply as algebras.

Suppose we have a Morita superequivalence $A\sim_{\sM}B$ induced by an $(A,B)$-bisupermodule $M$ and a $(B,A)$-bisupermodule $N$. For $b\in B$ and $m\in M$, we define $f_b(m):=(-1)^{|b||m|}mb$. Then, as in the classical situation, we have an even isomorphism of superalgebras $\phi:B^{\op}\iso \End_A(M),\  b\mapsto f_b$. 
Moreover, 
$N\cong\Hom_A(M,A)$, so 
$N$ can be recovered from $M$, and 
we sometimes just say that $A\sim_{\sM}B$ is induced by an $(A,B)$-bisupermodule $M$.


A standard argument shows  that that $A$ and $B$ are graded Morita superequivalent if and only if the graded supercategories $\mod{A}$ and $\mod{B}$ are graded superequivalent, i.e. there are functors $\funF: \mod{A}\to\mod{B}$ and $\funG:\mod{B}\to \mod{A}$ such that $\funF\circ\funG$ and $\funG\circ \funF$ are isomorphic to identities via graded supernatural transformations of degrees $(0,\0)$, cf.~\cite[Deinition 1.1(iv)]{BE}. 
Forgetting the gradings, $A$ and $B$ are Morita superequivalent if and only if the supercategories $\mod{A}$ and $\mod{B}$ are superequivalent, i.e. there are functors $\funF: \mod{A}\to\mod{B}$ and $\funG:\mod{B}\to \mod{A}$ such that $\funF\circ\funG$ and $\funG\circ \funF$ are isomorphic to identities via even  supernatural transformations, as in \cite[Deinition 1.1(iv)]{BE}.

\begin{Lemma}\label{lem:idmpt_Mor}
Let $A$ be a graded superalgebra and $e\in A_{0,\0}$ an idempotent such that $eL\neq 0$ for any irreducible graded $A$-supermodule $L$ (equivalently $AeA=A$). Then $eA \otimes_A ?$ and $Ae \otimes_{eAe} ?$ induce a graded Morita superequivalence between $A$ and $eAe$.
\end{Lemma}

\begin{proof}
$Ae$ (resp. $eA$) is a graded $(A,eAe)$-bisupermodule (resp.  graded $(eAe,A)$-bisupermodule) and the usual bimodule isomorphisms $Ae\otimes_{eAe}eA \cong A$ and $eA\otimes_{A}Ae \cong eAe$ are homogeneous isomorphisms of graded bisupermodules of bidegree $(0,\0)$.
\end{proof}

\begin{Corollary}\label{cor:idmpt_Mor}
Let $A$ be a graded superalgebra and $e\in A_{0,\0}$ an idempotent. 
If\, $|\Irr(eAe)|\geq |\Irr(A)|<\infty$ then $A\sim_{\gsM}eAe$.
\end{Corollary}
\begin{proof}
If $L$ is an irreducible graded $A$-supermodule then $eL$ is either $0$ or an irreducible graded $eAe$-supermodule, and all irreducible graded $eAe$-supermodules arise this way, cf. \cite[Theorem 6.2g]{Green}. So the result follows from Lemma~\ref{lem:idmpt_Mor}. 
\end{proof}

\begin{Lemma}\label{Lem_Mor_to_F}
Let $A$ be a graded superalgebra. The following are equivalent:
\begin{enumerate}
\item[{\rm (i)}] $A\sim_{\gsM}\F$.
\item[{\rm (ii)}] $A$ is finite dimensional and there exists an idempotent $e\in A_{0,\0}$ such that the homomorphism $A\to  \End_\F(Ae)$, which maps $a\in A$ to the left multiplication by $a$, is an isomorphism of graded superalgebras. 

\item[{\rm (iii)}] $A$ has a unique irreducible graded supermodule $L$ (up to isomorphism), and there exists an idempotent $e\in A_{0,\0}$ such that $eAe\cong \F$ and $eL\neq 0$.
\end{enumerate}
\end{Lemma}
\begin{proof}
(i)$\Rightarrow$(ii) 
Let $M$ (resp. $N$) be a graded $(A,\F)$-bisupermodule (resp.  graded $(\F,A)$-bisupermodule) such that $M \otimes_{\F} N\cong A$ as graded $A$-$A$-bisupermodules (resp. $N \otimes_A M \cong \F$ as graded $\F$-$\F$-bisupermodules) via an isomorphism of bidegree $(0,\0)$. From the classical (non-graded non-super) case, $A$ is isomorphic to a matrix algebra over $\F$ and $M$ is isomorphic to the column vector module over $A$. Since $A$ acts as the full $\F$-algebra of $\F$-linear transformations on $M$, we have $A\cong \End_\F(M)$ as graded superalgebras. Now it is clear that $A\cong Ae$ for an idempotent $e$ being the projection onto a $1$-dimensional graded subsuperspace of $M$. 

(ii)$\Rightarrow$(iii) is clear. 

(iii)$\Rightarrow$(i) Follows from Lemma~\ref{lem:idmpt_Mor}.
\end{proof}


The following lemma is a standard check:

\begin{Lemma}\label{lem:MSE}
Let $A_1,A_2,B_1,B_2$ be graded superalgebras. If 
$A_i$ and $B_i$ are graded Morita superequivalent via the graded $(A_i,B_i)$-bisupermodule $M_i$ and the graded $(B_i,A_i)$-bisupermodule $N_i$, for $i=1,2$, then $A_1\otimes A_2$ and $B_1\otimes B_2$ are graded Morita superequivalent via the graded $(A_1\otimes A_2,B_1\otimes B_2)$-bisupermodule $M_1\boxtimes M_2$ and the graded $(B_1\otimes B_2,A_1\otimes A_2)$-bisupermodule $N_1\boxtimes N_2$. 
\end{Lemma}


\subsection{More on Morita superequivalence}\label{sec:more_Morita}
In this subsection we will be considering only superalgebras, so there is no $\Z$-grading or $\Z$-grading is trivial, cf. Remark \ref{rem:non-sup_non-graded}. 

\begin{Remark}\label{rem:sup_Mor_eq}
Let $A,B$ be superalgebras and $M$ be an $(A,B)$-bisupermodule. 
To show that $M$ induces a Morita superequivalence between the superalgebras $A$ and $B$, it is enough to show that $M$ induces an ordinary Morita equivalence between $|A|$ and $|B|$ and the natural algebra homomorphism $|B|^{\op}\to \End_{|A|}(M)$ is an isomorphism of superalgebras, cf. 
Remark~\ref{rem:non-sup_non-graded0}.
Indeed, set $N$ to be the $(B,A)$-bisupermodule $\Hom_{|A|}(M,A)$. It is now easily observed that the standard $(|A|,|A|)$-bimodule isomorphism 
$$M\otimes_B N \to A,\ m\otimes n \mapsto n(m)$$ 
and the standard $(|B|,|B|)$-bimodule isomorphism 
$$N\otimes_A M \to \Hom_{|A|}(M,M)\cong B,\ n\otimes m \mapsto (m' \mapsto n(m')m)$$ are both even isomorphisms of bisupermodules.
\end{Remark}

Let $A$ be a superalgebra. 
We denote by $A^\times$ the set of units in $A$. 
If $A_{\1}\cap A^\times \neq \varnothing$, we call $A$ a {\em superalgebra with superunit}, and any $s\in A_{\1}\cap A^\times$ is called a {\em superunit}.\index{superunit}

\begin{Lemma}\label{lem:Mor_A0}
Let $A$ and $B$ be superalgebras with superunit. If the $(A,B)$-bisupermodule $M$ induces a Morita superequiavlence between $B$ and $A$, then $M_{\0}$ induces a Morita equiavlence between $B_{\0}$ and $A_{\0}$.
\end{Lemma}

\begin{proof}
Let $M$ be an $(A,B)$-bisupermodule and $N$ be a $(B,A)$-bisupermodule which induce a Morita superequivalence between $A$ and $B$. We first claim that
\begin{equation}\label{bimod_even_M}
\varphi:M_{\0}\otimes_{B_{\0}} N \to M\otimes_B N,\ 
m\otimes n\mapsto m\otimes n
\end{equation}
is an even isomorphism of $(A_{\0},A)$-bisupermodules. Well-definedness and the fact that it is a homomorphism of $A_{\0}$-$A$-supermodules are immediate. To show $\varphi$ is an isomorphism we fix some superunit $s\in B$ and construct the inverse homomorphism
\begin{align*}
\vartheta: M\otimes_B N\to M_{\0}\otimes_{B_{\0}} N,\ 
(m_0+m_1)\otimes n\mapsto (m_0\otimes n)+(m_1 s^{-1}\otimes s n)
\end{align*}
for all $m_0\in M_{\0}$, $m_1\in M_{\1}$ and $n\in N$. The only non-trivial thing to check is that $\vartheta$ is well-defined. Let $b_0\in B_{\0}$, $b_1\in B_{\1}$ and set $b=b_0+b_1$. Then
\begin{align*}
&\vartheta((m_0+m_1)b\otimes n) \\
=\,& \vartheta((m_0+m_1)(b_0+b_1)\otimes n)\\
=\,& \vartheta(m_0b_0\otimes n + m_0b_1\otimes n + m_1b_0\otimes n + m_1b_1\otimes n)\\
=\,& (m_0b_0\otimes n) + (m_0b_1s^{-1}\otimes sn) + (m_1b_0s^{-1}\otimes sn) + (m_1b_1\otimes n)\\
=\,& (m_0\otimes b_0n) + (m_0\otimes b_1s^{-1}sn) + (m_1s^{-1}\otimes sb_0s^{-1}sn) + (m_1s^{-1}\otimes sb_1n)\\
=\,& \vartheta(m_0\otimes b_0n + m_0\otimes b_1n + m_1\otimes b_0n + m_1\otimes b_1n)\\
=\,& \vartheta((m_0+m_1)\otimes (b_0+b_1)n) 
\\
=\, &\vartheta((m_0+m_1)\otimes bn),
\end{align*}
for all $m_0\in M_{\0}$, $m_1\in M_{\1}$ and $n\in N$, as desired.

Taking the even part of both sides of (\ref{bimod_even_M}) now gives $$M_{\0}\otimes_{B_{\0}}N_{\0}\cong A_{\0}$$ as $A_{\0}$-$A_{\0}$-bimodules. Similarly $$N_{\0}\otimes_{A_{\0}}M_{\0}\cong B_{\0}$$ as $B_{\0}$-$B_{\0}$-bimodules. The claim follows.
\end{proof}

\vspace{2mm}
The following superalgebras will play a central role in Chapter~\ref{P4} :
\begin{itemize}
\item[$\bullet$] The {\em matrix superalgebra\, $\cM_{m,n}(\F)$},
\index{matrix superalgebra}\index{m@$\cM_{m,n}(\F)$}
 which is just the 
matrix algebra $\cM_{m+n}(\F)$ of $(m+n)\times(m+n)$ matrices over $\F$ with parity given on the matrix units by requiring that 
$|E_{r,s}|=\0$ if and only if $r,s\leq m$ or $r,s>m$.
\item[$\bullet$] The rank $n$ {\em Clifford superalgebra\, $\cC_n$} 
\index{Clifford superalgebra}\index{c@$\cC_n$}
which is the superalgebra
given by odd generators $\cc_1,\dots,\cc_n$ subject to the relations $\cc_r^2=1$ for $r=1,\dots,n$ and $\cc_r\cc_s=-\cc_s\cc_r$ for all $1\leq r\neq s\leq n$. Recalling the tensor product of superalgebras from  \S\ref{SSBasicRep}, it is easy to see that $$\cC_n\otimes \cC_m\cong \cC_{n+m}.$$ In particular, $\cC_n\cong \cC_1^{\otimes n}$. If we work with $\cC_1$, we sometimes write $\cc=\cc_1$. 
\end{itemize}

\vspace{2mm}

\begin{Lemma}\label{lem:clif_mat1}
Let $A$ be a superalgebra. Then, for all $m,n$, we have that \, $A\otimes \cM_{m,n}(\F) \sim_{\sM}\hspace{-.3mm}A$. If, in addition, $\F$ has a primitive $4^{\nth}$ root of unity, then $\cC_{2n}\cong \cM_{2^{n-1},2^{n-1}}(\F)$. In particular, $A\otimes \cC_{2n} \sim_{\sM} A$. 
\end{Lemma}
\begin{proof}
For the matrix unit $E_{1,1}\in \cM_{m,n}(\F)_{\0}$, we have $E_{11} \cM_{m,n}(\F) E_{11} \cong \F$. So 
$
A\otimes \cM_{m,n}(\F) \sim_{\sM} A\otimes \F \cong A,
$ 
where the Morita superequivalence follows from Lemmas~\ref{lem:idmpt_Mor} and~\ref{lem:MSE}. 
For the isomorphism $\cC_{2n}\cong \cM_{2^{n-1},2^{n-1}}(\F)$, see for example \cite[Example 12.1.3]{Kbook}. 
\end{proof}


\begin{Lemma}\label{lem:clif_mat2}
Let $A$ be a superalgebra. We have an isomorphism of algebras $(A \otimes \cC_1)_{\0}\cong |A|$, and if $A$ has a superunit, then $A\otimes \cC_1 \sim_{\Mor} A_{\0}$.
\end{Lemma}

\begin{proof}
The isomorphism $(A \otimes \cC_1)_{\0}\iso |A|$ is given by 
$a\otimes 1\mapsto a$ if $a\in A_\0$, and $a\otimes \cc\mapsto \sqrt{-1}a$ if $a\in A_\1$. 
Using this isomorphism for the superalgebra $A\otimes \cC_1$ instead of $A$, 
 we get  $$(A\otimes \cC_2)_\0\cong (A\otimes \cC_1\otimes \cC_1)_\0\cong |A\otimes \cC_1|.$$ On the other hand, by Lemma~\ref{lem:clif_mat1}, $A\sim_{\sM}A\otimes \cC_2$. 
 So by Lemma~\ref{lem:Mor_A0}, we have $(A \otimes \cC_2)_{\0} \sim_{\Mor} A_{\0}$. 
\end{proof}

Let $G$ be a finite group. A {\em $G$-graded crossed product}\index{graded crossed product} will refer to an algebra $A$ with a decomposition $\bigoplus_{g\in G} A_g$ into subspaces such that $A_gA_h \subseteq A_{gh}$, for all $g,h\in G$, and such that, for all $g\in G$, we have $A_g\cap A^{\times} \neq \varnothing$. 
If, in addition, $A$ is a superalgebra and each $A_g$ is a subsuperspace then we refer to $A$ as a {\em $G$-graded crossed superproduct}.\index{graded crossed superproduct}
If $A$ is a $G$-graded crossed (super)product, then so is $A^{\op}$ by defining $A^{\op}_g=A_{g^{-1}}$, for all $g\in G$.
Note that $A_{1_G}$ is always a sub(super)algebra of $A$, and $(A_{1_G})^\op=A^\op_{1_G}$.

A $G$-graded crossed (super)product $A$ is called {\em symmetric} if it has a symmetrizing form $\tr: A \to \F$ that turns $A$ into a symmetric superalgebra, as in $\S$\ref{SSBasicRep} and restricts to a symmetrizing form on $A_{1_G}$.

If $A$ and $B$ are $G$-graded crossed (super)products, we define
\begin{equation}\label{EDiag}
(A,B)_{G}:= \sum_{g\in G}A_g\otimes B_{g^{-1}}\index{$(A,B)_{G}$}
=
\sum_{g\in G}A_g\otimes B^{\op}_g  \subseteq A \otimes B^{\op}.
\end{equation}
The definition of $G$-graded crossed (super)product ensures that $$A_{1_G}\otimes B_{1_G}^\op\subseteq (A,B)_{G}\subseteq A \otimes B^{\op}$$ 
are sub(super)algebras.

In part \ref{P4} we make repeated use of the following result which is essentially just a super version of \cite[Theorem 3.4(a)]{Mar}. Recall that we have identified 
$(A,B)$-bisupermodules and $(A\otimes B^{\op})$-supermodules via (\ref{bisupmodAB}).

\begin{Proposition}\label{prop:ext_marcus}
Let $G$ be a finite group, and $A,\, B$ be symmetric $G$-graded crossed superproducts. 
Suppose $M$ is an $(A_{1_G}\otimes B_{1_G}^\op)$-supermodule inducing a Morita superequivalence between $A_{1_G}$ and $B_{1_G}$. If $M$ extends to an $(A,B)_{G}$-supermodule, then\, $\Ind_{(A,B)_{G}}^{A\otimes B^{\op}}M$ induces a Morita superequivalence between $A$ and~$B$.
\end{Proposition}

\begin{proof}
Note that $|A|$ and $|B|$ are symmetric $G$-graded crossed products with $|A|_g=A_g$ and $|B|_g=B_g$ for all $g\in G$. 

View $M$ as an $(A_{1_G},B_{1_G})$-bisupermodule, then forget the super structures and view $M$ as an $(|A_{1_G}|,|B_{1_G}|)$-bimodule, which is the same as an  
$(|A_{1_G}|\otimes |B_{1_G}|^{\op})$-module. 
Since by assumption the $(A_{1_G}\otimes B_{1_G}^\op)$-supermodule structure on $M$ extends to $(A,B)_{G}$, the $(|A_{1_G}|\otimes |B_{1_G}|^\op)$-module structure on $M$ extends to $(|A|,|B|)_{G}$ via:
\begin{equation}\label{sup_nonsup}
(a_g\otimes b_{g^{-1}})*m := (-1)^{|b_{g^{-1}}||m|}(a_g\otimes b_{g^{-1}})\cdot m,
\end{equation}
for all $g\in G$, $a_g\in A_g$, $b_{g^{-1}}\in B_{g^{-1}}$ and $m\in M$, where `$*$' denotes the new action of $(|A|,|B|)_{G}$ on $M$ and `$\cdot$' the given action of $(A,B)_{G}$ on $M$. 

We now have the $A\otimes B^{\op}$-supermodule
$$
\bM:=\Ind_{(A,B)_{G}}^{A\otimes B^{\op}}M
$$
and the 
$|A|\otimes |B|^{\op}$-module
$$
\bM':=\Ind_{(|A|,|B|)_{G}}^{|A|\otimes |B|^{\op}}M.
$$ 
By \cite[Theorem 3.4(a)]{Mar}, the module $\bM'$ induces a Morita equivalence between $|A|$ and $|B|$. (Note that \cite{Mar} primarily deals with blocks of finite groups. However, \cite[Remarks 3.2(e)]{Mar} asserts that one can also apply the theorem in this more general, symmetric algebra setting. In particular, there is no condition on the field/ring the algebras are defined over.)

Recall that $\bM'$ is naturally a superspace with parity given by $|(a \otimes b)\otimes m| = |a|+|b|+|m|$ for $a\in A$, $b\in B$ and $m\in M$. So we can view $\bM'$ as an $(A, B)$-bisupermodule or as an 
$(A\otimes B^\op)$-supermodule using  (\ref{bisupmodAB}) so that:  
\begin{align*}
(a'\otimes b').[(a \otimes b) \otimes m] = (-1)^{|b'|(|a|+|b|+|m|)}(a'a \otimes bb') \otimes m,
\end{align*}
for all $a,a'\in A$, $b,b'\in B$ and $m \in M$. (Note above that $a\otimes b \in |A|\otimes |B|^{\op}$ whereas $a'\otimes b' \in A\otimes B^{\op}$.) 
We now deduce using Remark~\ref{rem:sup_Mor_eq} that $\bM'$ induces a Morita superequivalence between $A$ and $B$ since 
the natural isomorphism
$
|B|^{\op} \cong \End_{|A|}\left(\bM'\right)
$ 
is certainly an isomorphism of superalgebras.

The action of $A \otimes B^{\op}$ on $\bM$ is given by
\begin{align*}
(a'\otimes b').[(a\otimes b)\otimes m] = (-1)^{|b'|(|a|+|b|)}(a'a\otimes bb')\otimes m,
\end{align*}
for all $a,a'\in A$, $b,b'\in B$ and $m \in M$. (Here $a\otimes b$ and $a'\otimes b'$ both live in $A\otimes B^{\op}$.) Therefore, using (\ref{sup_nonsup}), we obtain an $(A\otimes B^{\op})$-supermodule isomorphism 
\begin{align*}
\bM' \iso \bM,\ 
(a \otimes b) \otimes m \mapsto (-1)^{|b||m|} (a\otimes b)\otimes m.
\end{align*}
Since $\bM'$ induces a Morita superequivalence between $A$ and $B$ by the previous paragraph, we now deduce that so does $\bM$. 
\end{proof}

\subsection{The Brauer tree algebra $\Zig_\ell$}\label{SSZig}
Recall the notation $J,K$ from \S\ref{ChBasicNotGen}. 
We consider a special {\em Brauer tree algebra}\, $\Zig_\ell$\index{a@$\Zig_\ell$} which is defined as the path algebra of the quiver 
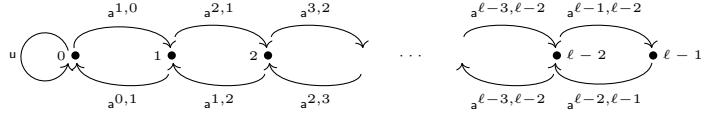
\begin{align*}
\begin{braid}\tikzset{baseline=3mm}
\coordinate (1) at (0,0);
\coordinate (2) at (4,0);
\coordinate (3) at (8,0);
\coordinate (4) at (12,0);
\coordinate (6) at (16,0);
\coordinate (L1) at (20,0);
\coordinate (L) at (24,0);
\draw [thin, black, ->] (-0.3,0.2) arc (15:345:1cm);
\draw [thin, black,->,shorten <= 0.1cm, shorten >= 0.1cm]   (1) to[distance=1.5cm,out=100, in=100] (2);
\draw [thin,black,->,shorten <= 0.25cm, shorten >= 0.1cm]   (2) to[distance=1.5cm,out=-100, in=-80] (1);
\draw [thin,black,->,shorten <= 0.25cm, shorten >= 0.1cm]   (2) to[distance=1.5cm,out=80, in=100] (3);
\draw [thin,black,->,shorten <= 0.25cm, shorten >= 0.1cm]   (3) to[distance=1.5cm,out=-100, in=-80] (2);
\draw [thin,black,->,shorten <= 0.25cm, shorten >= 0.1cm]   (3) to[distance=1.5cm,out=80, in=100] (4);
\draw [thin,black,->,shorten <= 0.25cm, shorten >= 0.1cm]   (4) to[distance=1.5cm,out=-100, in=-80] (3);
\draw [thin,black,->,shorten <= 0.25cm, shorten >= 0.1cm]   (6) to[distance=1.5cm,out=80, in=100] (L1);
\draw [thin,black,->,shorten <= 0.25cm, shorten >= 0.1cm]   (L1) to[distance=1.5cm,out=-100, in=-80] (6);
\draw [thin,black,->,shorten <= 0.25cm, shorten >= 0.1cm]   (L1) to[distance=1.5cm,out=80, in=100] (L);
\draw [thin,black,->,shorten <= 0.1cm, shorten >= 0.1cm]   (L) to[distance=1.5cm,out=-100, in=-100] (L1);
\blackdot(0,0);
\blackdot(4,0);
\blackdot(8,0);
\blackdot(20,0);
\blackdot(24,0);
\draw(0,0) node[left]{$0$};
\draw(4,0) node[left]{$1$};
\draw(8,0) node[left]{$2$};
\draw(14,0) node {$\cdots$};
\draw(20,0) node[right]{$\ell-2$};
\draw(24,0) node[right]{$\ell-1$};
 \draw(-2.6,0) node{$\zu$};
\draw(2,1.2) node[above]{$\za^{1,0}$};
\draw(6,1.2) node[above]{$\za^{2,1}$};
\draw(10,1.2) node[above]{$\za^{3,2}$};
\draw(18,1.2) node[above]{$\za^{\ell-3,\ell-2}$};
\draw(22,1.2) node[above]{$\za^{\ell-1,\ell-2}$};
\draw(2,-1.2) node[below]{$\za^{0,1}$};
\draw(6,-1.2) node[below]{$\za^{1,2}$};
\draw(10,-1.2) node[below]{$\za^{2,3}$};
\draw(18,-1.2) node[below]{$\za^{\ell-3,\ell-2}$};
\draw(22,-1.2) node[below]{$\za^{\ell-2,\ell-1}$};
\end{braid}
\end{align*}
generated by length $0$ paths $\{\ze^j\mid j\in J\}$\index{e@$\ze^j$}, and length $1$ paths $\zu$\index{u@$\zu$} and $\{\za^{k,k+1},\za^{k+1,k}\mid k\in K\}$\index{a@$\za^{j,k}$}, modulo the following relations:
\begin{enumerate}
\item all paths of length three or greater are zero;
\item all paths of length two that are not cycles are zero;
\item the length-two cycles based at the vertex $i\in\{1,\dots,\ell-1\}$ are equal;
\item $\zu^2=\za^{0,1}\za^{1,0}$. 
\end{enumerate}
For example, if $\ell=1$ the algebra $\Zig_\ell$ is the truncated polynomial algebra $\F[\zu]/(\zu^3)$. 
The algebra $\Zig_\ell$ is considered as a graded superalgebra by setting
$$
\bideg(\ze^j)=(0,\0),\ \bideg(\zu)=(2,\1),\ \bideg(\za^{k+1,k})=(4,\0),\ \bideg(\za^{k,k+1})=(0,\0).
$$
We have the bidegree $(4,\0)$ elements 
$$
\zc^0:=\zu^2\qquad \text{and}\qquad \zc^i:=\za^{i,i-1}\za^{i-1,i}\quad \text{for $i=1,\dots,\ell-1$}. \index{c@$\zc^i$}
$$
We set 
$$\zc:=\zc^0+\zc^1+\dots+\zc^{\ell-1},$$ 
so that $\zc^j=\zc\ze^j=\ze^j\zc$ for all $j\in J$. 
Note that
\begin{align*}
\zB_\ell:=\{\ze^j,\zc^j\mid j\in J\}\index{b@$\zB_\ell$}
\cup\{\za^{k,k+1},\za^{k+1,k}\mid k\in K\}
\cup\{\zu\}
\end{align*}
is a basis of $\Zig_\ell$. We refer to it as the {\em standard basis of $\Zig_\ell$}. 
It follows that  for any $i,j\in J$, we have:
\begin{equation}\label{EZigBiWt}
\dim_q\ze^i\Zig_\ell\ze^j=
\left\{
\begin{array}{ll}
1+q^2+q^4 &\hbox{if $i=j=0$,}\\
1+q^4 &\hbox{if $i=j\neq 0$,}\\
1&\hbox{if $j=i+1$,}\\
q^4&\hbox{if $j=i-1$,}\\ 
0&\hbox{if $|i-j|>1$,}\\
\end{array}
\right.
\end{equation}
and
\begin{equation}\label{EDimAL}
\dim_q\Zig_\ell=2\ell-1+q^2+(2\ell-1)q^4.
\end{equation}

Alternatively, $\Zig_\ell$ can be defined as the algebra generated by
\begin{equation}\label{ESomegeneratorsd=1}
\{\ze^j,\,  \za^{k,k+1},\, \za^{k+1,k},\, \zc,\, \zu\mid j\in J,\, k\in K\}
\end{equation}
subject only to the following relations:
\begin{align}
\label{RZig1}
\sum_{j\in J}\ze^j=1,\quad\ze^i\ze^j=\de_{i,j}\ze^i,
\\
\ze^{j}\zu=\zu\ze^{j}=\de_{j,0}\zu,\quad\ze^{j}\zc=\zc\ze^{j},
\\
\za^{i,j}\ze^k=\de_{j,k}\za^{i,j},\quad
\ze^k\za^{i,j}=\de_{i,k}\za^{i,j},
\\
\za^{i,j}\za^{k,l}=\de_{j,k}\de_{i,l}\zc\ze^i,
\quad
\zu^2=\zc\ze^0,
\\
\label{RZig5}
\zu^3=0,\quad\zc\za^{i,j}=\za^{i,j}\zc=0,\quad
\zu\za^{i,j}= \za^{i,j}\zu=0.
\end{align}

Recall the graded symmetricity notion from \S\ref{SSBasicRep}. 
The algebra $\Zig_\ell$ is $4$-graded symmetric with trace map
$$
\tr:\Zig_\ell\to\F,\ \zc^0\mapsto 1, \dots,\ \zc^{\ell-1}\mapsto 1,\  \zb\mapsto 0\quad\text{for all}\quad \zb\in\zB_\ell\setminus\{\zc^0,\dots,\zc^{\ell-1}\}.
$$
The corresponding distinguished element is 
$$
\nabla=\zu\otimes \zu +\sum_{j\in J}(\zc^j\otimes \ze^j+\ze^j\otimes \zc^j)+\sum_{k\in K}(\za^{k+1,k}\otimes \za^{k,k+1}+\za^{k,k+1}\otimes \za^{k+1,k}).
$$
Using the notation (\ref{EInsertion}), we can also write it as 
\begin{equation}\label{ENablaZig}
\nabla=\zu_1\zu_2 + \sum_{j\in J}(\zc^j_1\ze^j_2+\ze^j_1 \zc^j_2)+\sum_{k\in K}(\za^{k+1,k}_1 \za^{k,k+1}_2+\za^{k,k+1}_1 \za^{k+1,k}_2).
\end{equation}

\subsection{The affine Brauer tree algebra \mbox{$H_d(\Zig_\ell)$}}\label{SSAff}

In this subsection we introduce the rank $d$ affine Brauer tree Hecke superalgebra $H_d(\Zig_\ell)$. Its construction and properties fit into the general construction of \cite[\S3]{KM}, with the only difference that we need to take into account superstructures. 

Let $\zz$\index{z@$\zz$} be a variable of bidegree $(4,\0)$ and consider the polynomial algebra $\F[z]$. We denote by $\Zig_\ell[\zz]$ the free product $\F[z]\star\Zig_\ell$ subject to the relations 
\begin{equation}\label{ETwistedRelations}
\zu\zz=-\zz\zu\qquad \text{and}\qquad \zb\zz=\zz\zb\quad \text{for all $\zb\in\zB_\ell\setminus\{\zu\}$}.
\end{equation} 
It is easy to check that as graded superspaces $\Zig_\ell[\zz]\cong  \F[\zz]\otimes\Zig_\ell$. Equivalently, 
$\Zig_\ell[\zz]$ has basis 
$
\{\zz^k\zb\mid k\in\Z_{\geq 0},\,\zb\in\zB_\ell\}.
$
Note that $\bideg(\zz^k\zb)=(4k+\deg(\zb),|\zb|)$. 
(If $\zz$ was odd we could just define $\Zig_\ell[\zz]$ to be the tensor product $\F[\zz]\otimes\Zig_\ell$ of superalgebras, but we want $\zz$ to be even). 

Now fix $d\in\Z_{\geq 1}$. Consider the graded superalgebra 
$\Zig_\ell[\zz]^{\otimes d}$. Using the notation (\ref{EInsertion}), this algebra has basis
$$
\{\zz_1^{a_1}\cdots\zz_d^{a_d}\zb^1_1\cdots \zb^d_d\mid 
a_1,\dots,a_d\in\Z_{\geq 0},\ \zb^1,\dots,\zb^d\in\zB_\ell\}.
$$
We always consider the group algebra $\F\Si_d$ as a graded superalgebra concentrated in bidegree $(0,\0)$. Recall the notation (\ref{EWSignAction}), (\ref{ENablaZig}) and (\ref{EInsertion}); in particular we have elements $\nabla_{r,r+1}\in \Zig_\ell^{\otimes d}\subseteq \Zig_\ell[\zz]^{\otimes d}$. Define the {\em affine Brauer tree algebra $H_d(\Zig_\ell)$}\index{affine Brauer tree algebra}\index{h@$H_d(\Zig_\ell)$} to be the free product $\Zig_\ell[\zz]^{\otimes d}\star\F\Si_d$ subject to the following relations:
\begin{align}
w\,(\zb^1\otimes\dots\otimes \zb^d)={}^w(\zb^1\otimes\dots\otimes \zb^d)\,w\qquad(w\in\Si_d,\ \zb^1,\dots\zb^d\in\Zig_\ell)
\label{ERAff1}
\\
s_r\zz_t-\zz_{s_r(t)}s_r=\de_{r,t}\nabla_{r,r+1}-\de_{r+1,t}({}^{s_r}\nabla_{r,r+1})\qquad(1\leq r<d,\, 1\leq t\leq d). 
\label{ERAff2}
\end{align}

There are natural graded superalgebra homomorphisms (of bidegree  $(0,\0)$): 
\begin{align*}
\iota^{(1)}:\Zig_\ell[z]^{\otimes d} \to H_d(\Zig_\ell)
\quad\text{and}
\quad \iota^{(2)}:\F\Si_d\to H_d(\Zig_\ell).
\end{align*}
By the following theorem, these maps are in fact embeddings. 
We will use the same symbols for elements of the domain of these maps as for their images in $H_d(\Zig_\ell)$. 

\begin{Theorem}\label{TAffBasis} \cite[Theorem 3.8]{KM}
The map 
$$\Zig_\ell[z]^{\otimes d} \otimes \F\Si_d\to H_d(\Zig_\ell),\ x \otimes y \mapsto \iota^{(1)}(x)\iota^{(2)}(y)
$$
is an isomorphism of graded superspaces. In particular
$$
\{\zz_1^{a_1}\cdots\zz_d^{a_d}\zb^1_1\cdots \zb^d_dw\mid 
a_1,\dots,a_d\in\Z_{\geq 0},\ \zb^1,\dots,\zb^d\in\zB_\ell,\,w\in\Si_d\}.
$$
is a basis of $H_d(\Zig_\ell)$. 
\end{Theorem}

By the theorem and the relations, $\Zig_\ell^{\otimes d}\otimes\F\Si_d$ is a subalgebra of $H_d(\Zig_\ell)$ isomorphic to the wreath superproduct $\Zig_\ell\swr \Si_d$. 

For $\bi=i_1\cdots i_d\in  J^d$, define
\begin{equation}\label{EZEBI}
\ze^\bi:=\ze^{i_1}\otimes\dots\otimes \ze^{i_d}\in \Zig_\ell^{\otimes d}\subseteq H_d(\Zig_\ell).
\end{equation}
Then the relation (\ref{ERAff2}) is equivalent to the following relations for all $\bi\in J^d$:
\begin{equation}\label{ESZId}
(s_r \zz_t - \zz_{s_r(t)} s_r)\ze^\bi
=
\begin{cases}
\big((\delta_{r,t}- \delta_{r+1,t})(\zc_r + \zc_{r+1})
&
\\
\qquad\qquad\qquad+\de_{i_r,0}\zu_r\zu_{r+1}\big)\ze^\bi
&
\text{if}\ i_r = i_{r+1},\\
(\delta_{r,t}- \delta_{r+1,t})\za_r^{i_{r+1},i_r} \za_{r+1}^{i_r,i_{r+1}} \ze^\bi
&
\text{if}\ |i_r-i_{r+1}|=1,\\
0
& 
\textup{otherwise}.
\end{cases}
\end{equation}

\section{Combinatorics}
\label{ChComb}
\subsection{$p$-strict partitions}
\label{SSGen}


Let $\la$ be a partition. Collecting equal parts of $\la$, we  can write it in the form 
\begin{equation}\label{EStrictForm}
\la=(l_1^{m_1},\dots, l_k^{m_k})\ \text{with $l_1>\dots>l_k>0$ and $m_1,\dots,m_k\geq 1$}.
\end{equation} 
We denote 
\begin{equation}\label{ELaNorm}
\norm{\la}\,:=\prod_{r\,\, \text{with\,\,$p| l_r$}}\,\prod_{s=1}^{m_r}(1-(-q^2)^{s}).\index{$\norm{\la}$}
\end{equation}
If $m_r>1$ implies $p\mid l_r$ for all $1\leq r\leq k$ then $\la$ is called {\em $p$-strict}.\index{p@$p$-strict partition} Note that $0$-strict means simply {\em strict}\index{strict partition}, i.e. all parts are distinct. 
If, in addition, we have 
$$
\left\{
\begin{array}{ll}
\la_r-\la_{r+1}<p &\hbox{if $p\mid\la_r$,}\\
\la_r-\la_{r+1}\leq p &\hbox{if $p\nmid\la_r$.}
\end{array}
\right.
$$
then $\la$ is called {\em $p$-restricted}.\index{p@$p$-restricted  partition} We denote by $\Par_p(n)$\index{p@$\Par_p(n)$} the set of all $p$-strict partitions of $n$, and let 
$
\Par_p:=\bigsqcup_{n\geq 0}\Par_p(n).\index{p@$\Par_p$}
$
We denote by $\Par^\res_p(n)$\index{p@$\Par^\res_p(n)$} the set of all $p$-restricted $p$-strict partitions of $n$, and let  
$
\Par^\res_p:=\bigsqcup_{n\geq 0}\Par^\res_p(n).
\index{p@$\Par^\res_p$}
$

Let $\la$ be a $p$-strict partition. As usual, we identify $\la$ with its {\em Young diagram}\index{Young diagram} 
$$\la=\{(r,s)\in\Z_{>0}\times \Z_{>0}\mid s\leq \la_r\}.$$ 
We refer to the element $(r,s)\in\Z_{>0}\times \Z_{>0}$ as the {\em node}\index{node} in row $r$ and column $s$. 
We define a preorder `$\leq$' on the nodes via 
$
(r,s)\leq(r',s')
$
if and only if $s\leq s'$. 

We label the nodes with the elements of the set $I = \{0,1,\dots,\ell\}$ as follows: the labeling follows the repeating pattern
$$
0,1,\dots,\ell-1,\ell,\ell-1,\dots,1,0,
$$
starting from the first column and going to the right, 
see Example~\ref{Ex140821} below. If a node $\ttA\in \la$ is labeled with $i$, we say that $\ttA$ has {\em residue $i$}\index{residue} and write $\Res\, \ttA=i$. \index{r@$\Res\, \ttA$}

Following \cite{Morris,LT}, we can associate to every $\la\in\Par_p$ its {\em $\bar p$-core}\index{p@$\bar p$-core} $\core(\la)\in\Par_p^\res$.\index{c@$\core(\la)$} It is obtained from $\la$ by first removing all parts of $\la$ divisible by $p$ to get a strict partition $\mu$ (after reordering so that there are no zero parts); next, we remove {\em $p$-bars} as usual to get the usual $\bar p$-core of $\mu$. 
This means that we either remove two rows $i$ and $j$ such that $\mu_i+\mu_j=p$ or we remove $p$ boxes from row $i$ providing no part of $\mu$ equals $\mu_i-p$. After that we reorder the remaining parts so that we have a strict partition again and repeat the process until we arrive to a partition $\core(\lambda)$ such that no further removals are possible. It follows from \cite[Theorem~1]{MorYas} that $\core(\lambda)$ is well defined, see also \S\ref{SSAb}. It is clear from the definition that the number of nodes removed to go from $\la$ to $\core(\la)$ is divisible by $p$, we we can define a non-negative integer called the {\em $\bar p$-weight}\index{p@$\bar p$-weight} of $\la$:
\begin{equation}\label{EWi}
\wt(\la):=(|\la|-|\core(\la)|)/p\in\Z_{\geq 0}.\index{w@$\wt$}
\end{equation}

A partition $\rho\in\Par_p$ is called a {\em $\bar p$-core} if $\core(\rho)=\rho$. For such $\rho$, we denote 
\begin{equation}\label{EParRhoD}
\Par_p(\rho,d):=\{\la\in\Par_p\mid \core(\la)=\rho,\ \wt(\la)=d\}.
\index{p@$\Par_p(\rho,d)$}
\end{equation}
We also denote by $\Cores_p$\index{c@$\Cores_p$} the set of all $\bar p$-core partitions:
\begin{equation}\label{ECores}
\Cores_p:=\{\la\in\Par_p\mid \text{$\la$ is a $\bar p$-core}\}.
\end{equation}

\vspace{2mm}
\begin{Example}\label{Ex140821}
Let $\ell=2$, so $p=5$. The partition $\la=(16, 11,10,10,9,4,1)$ is $5$-strict. The residues of the nodes are as follows:
\vspace{3mm}
$$
\begin{ytableau}
$0$ & $1$ & $2$ & $1$ & $0$ & $0$& $1$ & $2$ & $1$ & $0$ & $0$ & $1$ & $2$ & $1$ & $0$ & $0$\cr 
$0$ & $1$ & $2$ & $1$ & $0$ & $0$ & $1$ & $2$ & $1$ & $0$ & $0$ \cr
$0$ & $1$ & $2$ & $1$ & $0$ & $0$ & $1$ & $2$ & $1$ & $0$  \cr
$0$ & $1$ & $2$ & $1$ & $0$ & $0$ & $1$ & $2$ & $1$ & $0$  \cr
$0$ & $1$ & $2$ & $1$ & $0$ & $0$ & $1$ & $2$ & $1$ \cr
$0$ & $1$ & $2$ & $1$ \cr
$0$ \cr
\end{ytableau}\vspace{3 mm}
$$
The $\bar 5$-core of $\la$ is $(1)$.
\end{Example}

\vspace{2mm}
Recalling $\al_i$'s and $Q_+$ from \S\ref{ChBasicNotLie}, define the {\em residue content}\index{residue content} of $\la\in\Par_p$ 
\begin{equation}\label{SEResCont}
\cont(\la):=\sum_{\ttA\in \la}\al_{\Res\, \ttA}\in Q_+.
\index{c@$\cont(\la)$}
\end{equation}

\begin{Lemma}\label{LLT}
\cite[Theorem 5]{MorYas} 
\label{LCore}
Let $\la,\mu\in \Par_p(n)$. Then $\core(\la)=\core(\mu)$ if and only if 
$\cont(\la)=\cont(\mu)$.
\end{Lemma}

Moves taking $\la$ to $\core(\la)$ reduce the content by a multiple of $\de$, so from Lemma~\ref{LLT} we get:

\begin{Lemma}  \label{LCoreCont} 
Let $\rho\in \Cores_p$, $d\in\Z_{\geq 0}$, and $\la\in\Par_p$. Then  $\la\in \Par_p(\rho,d)$ if and only if 
$\cont(\la)=\cont(\rho)+d\de$. 
\end{Lemma}

\begin{Lemma}  \label{LCoreNoDeSubtr}
Let $\la\in \Par_p$. Then $\la$ is a $\bar p$-core if and only is there is not $\mu\in\Par_p$ with $\cont(\mu)+\de=\cont(\la)$.  
\end{Lemma}
\begin{proof}
Suppose $\la$ is not a $p$-bar core. If $\la$ has a non-zero part divisible $p$, we can reduce this part by $p$ (and reorder the parts of necessary) to get a partition $\mu\in\Par_p$ with $\cont(\mu)=\cont(\la)-\de$. If $\la$ has no parts divisible by $p$, a removal of a $p$-bar also reduces $\cont(\la)$ by $\de$. 

Conversely, suppose there is $\mu\in\Par_p$ with $\cont(\mu)+\de=\cont(\la)$. Add a part equal to $p$ to $\mu$. This will create a partition $\nu\in \Par_p$ with $\cont(\nu)= \cont(\mu)+\de=\cont(\la)$, which is not a $\bar p$-core. By  Lemma~\ref{LLT}, $\la$ is also not a $p$-bar core. 
\end{proof}

\subsection{Abaci}\label{SSAb}
We will use the $\bar p$-abacus notation for partitions introduced in~\cite{MorYas}. As we work slightly more generally allowing $p$-strict (not just strict) partitions, we give all the necessary definitions.  
 
We define the {\em abacus}\index{abacus} $\Ab:=\Z_{\geq 0} \times \Z/p\Z$. \index{a@$\Ab$}
When convenient we identify $\Z/p\Z$ with the subset $\{0,1,\dots,p-1\}\subset \Z$. 
The element $(n,j)\in\Ab$ is referred to as the {\em position in row $n$ and on runner $j$}.\index{position}\index{runner} The position $(n,j)$ is also referred to as the $(pn+j)^{\nth}$ position (or the position $pn+j$), so that the  positions are labeled with non-negative integers starting with the position in row $0$ on runner $0$ and going along the rows. Thus the positions are labeled as follows (note rows increase from top to bottom): 
$$
\begin{picture}(120, 65)%
\put(0,50){\makebox(0,0)[b]{${}_0$}}
\put(30,50){\makebox(0,0)[b]{${}_1$}}
\put(60,50){\makebox(0,0)[b]{$\cdots$}}
\put(90,50){\makebox(0,0)[b]{${}_{p-2}$}}
\put(120,50){\makebox(0,0)[b]{${}_{p-1}$}}
\put(0,30){\makebox(0,0)[b]{${}_p$}}
\put(30,30){\makebox(0,0)[b]{${}_{p+1}$}}
\put(60,30){\makebox(0,0)[b]{$\cdots$}}
\put(90,30){\makebox(0,0)[b]{${}_{2p-2}$}}
\put(120,30){\makebox(0,0)[b]{${}_{2p-1}$}}
\put(0,10){\makebox(0,0)[b]{$\vdots$}}
\put(30,10){\makebox(0,0)[b]{$\vdots$}}
\put(90,10){\makebox(0,0)[b]{$\vdots$}}
\put(120,10){\makebox(0,0)[b]{$\vdots$}}
\end{picture}
$$

Let $\la\in\Par_p$. Fix an integer $N\geq h(\la)$ and write $\la$ using $N$ parts: $\la=(\la_1,\dots,\la_N)$ (of course $\la_k=0$ for all $h(\la)\leq k\leq N$). The corresponding {\em abacus display of $\la$}\index{abacus display}, denoted 
$\Ab_\la$\index{a@$\Ab_\la$}, consists of $N$ {\em beads}\index{bead} occupying positions $\la_1,\la_2,\dots,\la_N$. Note that positions on runner $0$ can be occupied with several beads, since $\la$ can have repeated parts as long as they are divisible by $p$; in particular the $0^{\nth}$ position is occupied with $N-h(\la)$ beads. Usually, when working with $\la$ of weight $d$ we choose $N\geq h(\la)+d$. 

\begin{Example} \label{EAb} 
Continue with $\ell$ and $\la$ from Example~\ref{Ex140821} and take  $N=10$. Then, denoting beads by $\times$'s and empty positions with $\cdot$'s, we have
\vspace{1mm}
$$
\begin{picture}(60, 85)%
\put(-40,35){\makebox(0,0)[b]{$\Ab_\la=$}}
\put(0,70){\makebox(0,0)[b]{$\times^3$}}
\put(30,70){\makebox(0,0)[b]{$\times$}}
\put(60,70){\makebox(0,0)[b]{$\cdot$}}
\put(90,70){\makebox(0,0)[b]{$\cdot$}}
\put(120,70){\makebox(0,0)[b]{$\times$}}
\put(0,50){\makebox(0,0)[b]{$\cdot$}}
\put(30,50){\makebox(0,0)[b]{$\cdot$}}
\put(60,50){\makebox(0,0)[b]{$\cdot$}}
\put(90,50){\makebox(0,0)[b]{$\cdot$}}
\put(120,50){\makebox(0,0)[b]{$\times$}}
\put(0,30){\makebox(0,0)[b]{$\times^2$}}
\put(30,30){\makebox(0,0)[b]{$\times$}}
\put(60,30){$\cdot$}
\put(90,30){\makebox(0,0)[b]{$\cdot$}}
\put(120,30){\makebox(0,0)[b]{$\cdot$}}
\put(0,10){\makebox(0,0)[b]{$\cdot$}}
\put(30,10){\makebox(0,0)[b]{$\times$}}
\put(60,10){$\cdot$}
\put(90,10){\makebox(0,0)[b]{$\cdot$}}
\put(120,10){\makebox(0,0)[b]{$\cdot$}}
\end{picture}
$$
Note that we only drew the rows $0,1,2,3$ of $\Ab_\la$ since there are no beads in larger rows. 
\end{Example}

Let $\la\in\Par_p$. For $j\in \{0,1,\dots,p-1\}$, we denote
\begin{equation}\label{EBJ}
b_j^\la:=\{k\mid \la_k>0\ \text{and}\ \la_k\equiv j\pmod{p}\}.
\end{equation}
Note that
$$
b_j^\la:=\index{b@$b_j^\la$}
\left\{
\begin{array}{ll}
\text{number of beads on the $j^{\nth}$ runner of $\Ab_\la$} &\hbox{if $j\neq 0$,}\\
\text{number of beads in rows $>0$ on the $0^{\nth}$ runner of $\Ab_\la$} &\hbox{if $j= 0$.}
\end{array}
\right.
$$

Let $1\leq r\leq h(\la)$ be such that $\la_r\geq p$. Let position $\la_r$ be on runner $j$. 
Note that the removal of $p$ boxes from row $r$ of $\la$ corresponds to moving a bead of $\Ab_\la$ from position $\la_r$ to position $\la_r-p$ which is exactly above it on the same runner $j$. Moreover, if $j\neq 0$, i.e. $p\nmid \la_r$, removing $p$ boxes from row $r$ is an allowed $p$-bar removal if and only if the position $\la_r-p$ is not occupied in $\Ab_\la$. Moving a bead on runner $j$ of $\Ab_\la$ one position up with the requirement that the new position is not occupied if $j\neq 0$ will be called an {\em elementary slide up (on runner $j$)}\index{elementary slide up/down}. If, in addition, we demand that the new position is not occupied even when $j=0$, then we call this a {\em strict elementary slide up}.\index{strict elementary slide up/down}
 {\em Elementary slides down} and {\em strict elementary slides down} are defined similarly: for example a strict elementary slide down is moving a bead from position $a$ to position $a+p$ under the condition that position $a+p$ is not occupied.

Let $1\leq r,s\leq h(\la)$ be such that $\la_r+\la_s=p$. Let  position $\la_r$ be on runner $j$. Note that $j\neq 0$ and position $\la_s$ is on runner $p-j$. Moreover, positions $\la_r,\la_s$ are in row $0$. The removal of rows $r$ and $s$ of $\la$ corresponds to the removal of beads in row $0$ on runners $j$ and $p-j$. Removing two such beads will be called an {\em elementary removal  (on runners $j,p-j$)}. \index{elementary removal}

It follows from the previous two paragraphs that  
$$b_j^{\core(\la)}=
\left\{
\begin{array}{ll}
b_j^\la-\min(b_j^\la,b_{p-j}^\la) &\hbox{if $j\neq 0$}\\
0 &\hbox{if $j=0$}
\end{array}
\right.
$$
From this it is clear that $\core(\la)$ is well defined.

\begin{Lemma} \label{LComputation} 
Let $\la\in\Par_p(n)$ and $c_\ell(\la)$ be the number of nodes of residue $\ell$ in $\la$. Then
$$
pc_\ell(\la)-n=\sum_{i=\ell+1}^{p-1}(p-i)b_i^\la-\sum_{i=1}^\ell ib_i^\la.
$$
\end{Lemma}
\begin{proof}
Note that the left hand side and the right hand side do not change 
under elementary removals and slides up, so we may assume that $\la$ is a $\bar p$-core. In this case it is easy to see that
\begin{align*}
c_\ell(\la)=\sum_{i=1}^\ell\sum_{k=1}^{b_i^\la-1}k+\sum_{i=\ell+1}^{p-1}\sum_{k=1}^{b_i^\la}k
\end{align*}
and
\begin{align*}
n=\sum_{i=1}^{p-1}\sum_{k=0}^{b_i^\la-1}(i+kp),
\end{align*}
which easily implies the result.
\end{proof}

Let $\la\in\Par_p$ with $\wt(\la)=d$. Following \cite[p.27]{MorYas}, we define the multipartition $\quot(\la)=(\la^{(0)},\dots,\la^{(\ell)})\in\Par^{\ell+1}(d)$\index{q@$\quot(\la)$} referred to as the the {\em $\bar p$-quotient}\index{p@$\bar p$-quotient} of $\la$. The partition $\la^{(0)}$ records the configuration of beads on the $0$ runner of $\Ab_\la$ as follows: to every bead in row $r>0$ and runner $0$ we associate a part equal to $r$ of $\la^{(0)}$. Let $1\leq j\leq \ell$. The partition $\la^{(j)}$ is determined by the bead configuration on runners $j$ and $p-j$. Suppose that one of the runners $j,p-j$ has beads in rows $b_1>\dots>b_r$ while the other in positions $a_1>\dots>a_{s+r}$. 
Then we put 
$$
\la^{(j)}:=(a_1-s+1,a_2-s+2,\dots,a_s;(a_{s+1},\dots,a_{s+r}\mid b_1,\dots,b_r)),
$$
where $(a_{s+1},\dots,a_{s+r}\mid b_1,\dots, b_r)$ is a partition in Frobenius notation and for partitions $\mu,\nu$ we denote by $(\mu;\nu)$ the partition obtained by concatenation of $\mu$ and $\nu$ (provided the last non-zero part of $\mu$ is not less than the first part of $\nu$). 

\begin{Example}
Continuing with $\la$ and $\ell$ as in Examples~\ref{Ex140821}, \ref{EAb}, we get $\la^{(0)}=(2,2)$, 
$
\la^{(1)}=(3;(2,0\mid 1,0))=(3,3,2)
$ 
and 
$\la^{(2)}=\varnothing$. 
\end{Example}

It is easy to see that one can recover $\la$ from $\core(\la)$ and $\quot(\la)$. In fact,

\begin{Lemma} \label{TMY} {\rm \cite[Theorem 2]{MorYas}}
Let $\rho\in \Cores_p$. Then assignment $\la\mapsto\quot(\la)$ is a bijection between $\Par_p(\rho,d)$ and $\Par^{\ell+1}(d)$. 
\end{Lemma}

\subsection{Removable and addable nodes and tableaux}
\label{SRemAdd}
Let $\la$ be a $p$-strict partition and $i \in I$.
A node $\ttA
\in \la$ is called {\em $i$-removable}\index{removable node} (for $\la$) if one of the following
holds:
\begin{enumerate}
\item[(R1)] $\Res\, \ttA = i$ and
$\la_\ttA:=\la\setminus \{\ttA\}$\index{$\la_\ttA$} is again a $p$-strict partition; such $\ttA$'s are also called {\em properly $i$-removable};
\index{properly removable node}
\item[(R2)] the node $\ttB$ immediately to the right of $\ttA$
belongs to $\la$,
$\Res\, \ttA = \Res\, \ttB = i$,
and both $\la_\ttB = \la \setminus  \{\ttB\}$ and
$\la_{\ttA,\ttB} := \la \setminus  \{\ttA,\ttB\}$ are $p$-strict partitions.
\end{enumerate}
A node $\ttB\notin\la$ is called 
{\em $i$-addable}\index{addable node} (for $\la$) if one of the following holds:
\begin{enumerate}
\item[(A1)] $\Res\, \ttB = i$ and
$\la^\ttB:=\la\cup\{\ttB\}$\index{$\la^\ttB$} is again a $p$-strict partition; such $\ttB$'s are also called {\em properly $i$-addable};\index{properly addable node}
\item[(A2)] 
the node $\ttA$
immediately to the left of $\ttB$ does not belong to $\la$,
$\Res\, \ttA = \Res\, \ttB = i$, and both 
$\la^\ttA = \la \cup \{\ttA\}$ and 
$\la^{\ttA, \ttB} := \la \cup\{\ttA,\ttB\}$ are $p$-strict partitions.
\end{enumerate}
We note that (R2) and (A2) above are only possible if $i = 0$.
For $i \in I$, we denote by $\Add_i(\la)$\index{a@$\Add_i(\la)$} (resp. $\Rem_i(\la)$)\index{r@$\Rem_i(\la)$} the set of all $i$-removable (resp. $i$-addable) nodes for $\la$. 
We also denote by $\PA_i(\la)$\index{p@$\PA_i(\la)$} (resp. $\PR_i(\la)$)\index{p@$\PA_i(\la)$}\index{p@$\PR_i(\la)$} the set of all properly $i$-removable (resp. properly $i$-addable) nodes for $\la$.

Recall the element $\ttK=\sum_{i=0}^\ell a_i^\vee h_i$ from \S\ref{ChBasicNotLie}.

\begin{Lemma} \label{LemmaK} 
Let $\la\in\Par_p$. Then $\lan\La_0-\cont(\la),h_i\ran=|\Add_i(\la)|-|\Rem_i(\la)|$ for every $i\in I$, and  
$$\sum_{i\in I}a_i^\vee(|\Add_i(\la)|-|\Rem_i(\la)|)=1.$$
\end{Lemma}
\begin{proof}
The first statement is checked by induction on $|\la|$. The second statement follows from the first:
\begin{align*}
\sum_{i\in I}a_i^\vee(|\Add_i(\la)|-|\Rem_i(\la)|)&=\sum_{i\in I}a_i^\vee\lan h_i,\La_0-\cont(\la)\ran
\\
&=\lan\sum_{i\in I}a_i^\vee h_i, \La_0-\cont(\la)\ran
\\
&=\lan \ttK, \La_0-\cont(\la)\ran=1
\end{align*} 
since $\lan \ttK,\La_0\ran=1$ and $\lan \ttK,\al_i\ran=0$ for all $i$. 
\end{proof}

Let $\la\in\Par_p$ be written in the form (\ref{EStrictForm}). Suppose $\ttA\in \PR_i(\la)$. Then there is $1\leq r\leq k$ such that $\ttA=(m_1+\dots+m_r,l_r)$.  
Recalling the preorder `$\leq$' on the nodes from \S\ref{SSGen}, we define 
\begin{align*}
\eta_\ttA(\la)&:=\sharp\{\ttC\in\Rem_i(\la)\mid \text{$\ttC>\ttA$}\}
-\sharp\{\ttC\in\Add_i(\la)\mid \text{$\ttC>\ttA$}\},
\index{e@$\eta_\ttA(\la)$}
\\
\zeta_\ttA(\la)&:=
\left\{
\begin{array}{ll}
(1-(-q^2)^{m_r}) &\hbox{if $p\mid l_r$,}\\
1 &\hbox{otherwise.}
\end{array}
\right.
\index{z@$\zeta_\ttA(\la)$}
\\
d_\ttA(\la)&:= q_i^{\eta_\ttA(\la)} \zeta_\ttA(\la).
\index{d@$d_\ttA(\la)$}
\end{align*}
Suppose $\ttB\in\PA_i(\la)$. Then there is $r$ such that $1\leq r\leq k+1$ and $\ttB=(m_1+\dots+m_{r-1}+1,l_r+1)$, where we interpret $l_{k+1}$ as $0$.   We define 
\begin{align}
\label{EEtaDef}
\eta^\ttB(\la)&:=\sharp\{\ttC\in\Add_i(\la)\mid \text{$\ttC<\ttB$}\}
-\sharp\{\ttC\in\Rem_i(\la)\mid \text{$\ttC<\ttB$}\},
\index{e@$\eta^\ttB(\la)$}
\\
\zeta^\ttB(\la)&:=
\left\{
\begin{array}{ll}
(1-(-q^2)^{m_r}) &\hbox{if $r\leq k$ and $p\mid l_r$,}\\
1 &\hbox{otherwise.}
\end{array}
\right.
\index{z@$\zeta^\ttB(\la)$}
\\
d^\ttB(\la)&:= q_i^{\eta^\ttB(\la)} \zeta^\ttB(\la).
\index{d@$d^\ttB(\la)$}
\end{align}

\begin{Example} 
{\rm 
Let $\ell=2$ so $p=5$. The partition $\la=(16, 11,10,10,9,5,1)$ is $5$-strict, and the residues of its boxes are labeled on the diagram below:
\vspace{3mm}
$$
\begin{ytableau}
$0$ & $1$ & $2$ & $1$ & $0$ & $0$& $1$ & $2$ & $1$ & $0$ & $0$ & $1$ & $2$ & $1$ & $0$ & $0$\cr 
$0$ & $1$ & $2$ & $1$ & $0$ & $0$ & $1$ & $2$ & $1$ & $0$ & $0$ \cr
$0$ & $1$ & $2$ & $1$ & $0$ & $0$ & $1$ & $2$ & $1$ & $0$  \cr
$0$ & $1$ & $2$ & $1$ & $0$ & $0$ & $1$ & $2$ & $1$ & $0$  \cr
$0$ & $1$ & $2$ & $1$ & $0$ & $0$ & $1$ & $2$ & $1$ \cr
$0$ & $1$ & $2$ & $1$ & $0$\cr
$0$ \cr
\end{ytableau}
$$

\vspace{4mm} 
\noindent We mark the removable nodes as $\ttA_r$ and addable nodes as $\ttB_s$:
\vspace{4mm} 
$$
\ytableausetup{mathmode}
\begin{ytableau}
\, &  &  &  &  &  &  &  &  &  &  &  &  &  & \ttA_5 & \ttA_6 & \none[\ttB_5] \\ 
 &  &  &  &  &  &  &  &  &  & \ttA_4 & \none[\ttB_4] \\
 &  &  &  &  &  &  &  &  &   \\
 &  &  &  &  &  &  &  &  &  \\
 &  &  &  &  &  &  &  & \ttA_3 & \none[\ttB_3] \\ 
 &  &  & & \ttA_2 & \none[\ttB_2] \\
\ttA_1 & \none[\ttB_1] \\
\end{ytableau}
$$
\vspace{2mm} 

\noindent Note that $\ttA_1,\ttA_2,\ttA_4,\ttA_5,\ttA_6$ are $0$-removable (with $\ttA_5$ not properly $0$-removable), 
$\ttB_2,\ttB_3$ are $0$-addable, 
$\ttA_3$ is $1$-removable, and $\ttB_1,\ttB_4,\ttB_5$ are $1$-addable. We have
\begin{align*}
&d_{\ttA_1}(\la)=q^2,\ d_{\ttA_2}(\la)=q(1+q^{2}),\ d_{\ttA_3}(\la)=q^{-4},\ d_{\ttA_4}(\la)=q^{2},\ d_{\ttA_6}(\la)=1,
\\
&d^{\ttB_1}(\la)=1,\ d^{\ttB_2}(\la)=q^{-2}(1+q^{2}),\ d^{\ttB_3}(\la)=q^{-1},\  d^{\ttB_4}(\la)=1,\ d^{\ttB_5}(\la)=q^{2}.
\end{align*}
}
\end{Example}

\begin{Example} \label{Ex1} 
{\rm 
Let $\ell=2$ so $p=5$.
The partition $\la=(5, 5,2)$ is $5$-strict,
and the residues of its boxes are labeled on the diagram below: 

\vspace{2mm} 
$$
\begin{ytableau}
$0$ & $1$ & $2$ & $1$ & $0$ \cr 
$0$ & $1$ & $2$ & $1$ & $0$\cr
$0$ & $1$\cr
\end{ytableau}\vspace{2 mm}
$$

\vspace{3 mm}
\noindent
The only $0$-removable node is marked as $\ttA_1$, and the $0$-addable nodes are marked as $\ttB_1,\ttB_2$:

\vspace{2mm} 
$$
\ytableausetup{mathmode}
\begin{ytableau}
\, &  &  &  &  & \none[\ttB_2] \cr 
 &  &  &  & \ttA_1 \cr
 & \cr
 \none[\ttB_1]
\end{ytableau}\vspace{2 mm}
$$

\vspace{1mm} 
\noindent
We have $d^{\ttB_1}(\la)=1$ and $d^{\ttB_2}(\la)=(1-q^4)$. On the other hand for the partition $\mu=(5)$ and the node $\ttB=(1,6)$, we have $d^\ttB(\mu)=(1+q^2)$. 
}
\end{Example}

\vspace{2mm}
Let $\la\in\Par_p(n)$. 
A {\em $\la$-tableau} \index{tableau} is a bijection $\T: \{1,\dots,n\}\to \la$, where we have identified $\la$ with the set of nodes of its Young diagram. In this case we denote
$$
\bi^\T:=(\Res \T(1))\,\cdots\, (\Res \T(n))\in I^n.\index{i@$\bi^\T$}
$$
A $\la$-tableau $\T$ is called {\em $p$-standard}\index{p@$p$-standard tableau} if 
$\T(\{1,\dots,k\})$ is a $p$-strict partition for all $k=1,\dots,n$. 
In this case, the restriction $\T_{\leq k}:=\T|_{\{1,\dots,k\}}$ can be considered as a $p$-standard $\T(\{1,\dots,k\})$-tableau. We denote by $\Std_p(\la)$\index{s@$\Std_p(\la)$} the set of all $p$-standard $\la$-tableaux. 
For $\bj\in I^n$, we denote
$$
\Std_p(\la,\bj)=\{\T\in\Std_p(\la)\mid \bi^\T=\bj\}.
\index{s@$\Std_p(\la,\bj)$}
$$ 

Let $\T\in\Std_p(\la)$. 
The {\em degree}\index{degree} of $\T$ is defined as 
\begin{equation}\label{EDeg}
\deg(\T):=\prod_{k=1}^n d^{\T(k)}(\T(\{1,\dots,k-1\})).
\index{d@$\deg(\T)$}
\end{equation}
Note that 
$$\deg(\T)=d^{\T(n)}(\T(\{1,\dots,n-1\}))\deg(\T_{\leq n-1}).$$

\subsection{$p$-strict multipartitions}
\label{SRemAddMulti}
We fix $N\in\Z_{\geq 1}$ and consider the set $\Par^N_p$ of  all {\em $p$-strict $N$-multipartitions}\index{p@$p$-strict multipartition} of $n$, i.e. the $N$-tuples $\bla=(\la^{(1)},\dots,\la^{(N)})$ such that $\la^{(1)},\dots,\la^{(N)}\in\Par_p$. 
For $\bla=(\la^{(1)},\dots,\la^{(N)})\in\Par_p^N$ we write $|\bla|:=|\la^{(1)}|+\dots+|\la^{(N)}|$. For $n\in\Z_{\geq 0}$, let 
$$\Par_p^N(n)=\{\bla\in\Par_p^N\mid|\bla|=n\}. \index{p@$\Par_p^N(n)$}
$$
Recalling (\ref{ELaNorm}), we denote 
\begin{equation}\label{ELaNormM}
\norm{\bla}\,:=\,\prod_{t=1}^{N}\Vert\la^{(t)}\Vert.
\index{$\norm{\bla}$}
\end{equation}


Generalizing (\ref{EParRhoD}), for a $\bar p$-core $\rho$, we set
\begin{equation}\label{EParNRhoD}
\Par_p^N(\rho^N;d):=\{\bla\in\Par^N_p\mid \core(\la^{(t)})=\rho
\ \text{for all $t$ and}\ \sum_{t=1}^N\wt(\la^{(t)})=d\}.
\index{p@$\Par_p^N(\rho^N;d)$}
\end{equation}

Let $\bla\in\Par_p^N$. We identify $\bla$ with its {\em Young diagram}\index{Young diagram} $$\la=\{(r,s,t)\in\Z_{>0}\times \Z_{>0}\times\{1,\dots,N\}]\mid s\leq \la_r^{(t)}\},$$ 
which can be thought of as a disjoint union of the Young diagrams 
$\la^{(1)}\sqcup\dots\sqcup \la^{(N)}$. 
We refer to the element 
$$(r,s,t)\in\Z_{>0}\times \Z_{>0}\times\{1,\dots,N\}$$ as the node\index{node} in row $r$, column $s$ and component $t$\index{component}. 
We define a preorder `$\leq$' on the nodes via 
$
(r,s,t)\leq(r',s',t')
$
if and only if either $t>t'$, or $t=t'$ and $s\leq s'$. 

The residue $\Res\, \ttA$\index{residue} of a node $\ttA=(r,s,t)$ does not depend on $t$ and is set to be the same as the residue of the node $(r,s)$ of the Young diagram $\la^{(t)}$.  
The {\em residue content}\index{residue content} of $\bla$ is  
\begin{equation}\label{SEResContM}
\cont(\bla)=\sum_{\ttA\in \bla}\al_{\Res\, \ttA}=\sum_{t=1}^N\cont(\la^{(t)})\in Q_+.
\index{c@$\cont(\bla)$}
\end{equation}

Let $\bla\in\Par_p^N(n)$ and $i \in I$.
A node $\ttA=(r,s,t)\in \bla$ is called {\em $i$-removable}\index{removable node} (resp. properly $i$-removable)\index{properly removable node} for $\bla$ if it is so for the component $\la^{(t)}$. If $\ttA$ is properly removable for $\bla$, we have $\bla_\ttA:=\bla\setminus\{\ttA\}\in\Par_p^N(n-1)$. \index{$\bla_\ttA$}
A node $\ttB=(r,s,t)\notin\bla$ is called 
{\em $i$-addable}\index{addable node} (resp. properly $i$-addable)\index{properly addable node} for $\bla$ if it is so for the component $\la^{(t)}$. If $\ttB$ is properly addable for $\bla$, we have $\bla^\ttB:=\bla\cup\{\ttB\}\in\Par_p^N(n+1)$. \index{$\bla^\ttB$}
We denote by $\Add_i(\bla)$\index{a@$\Add_i(\bla)$} (resp. $\Rem_i(\bla)$\index{r@$\Rem_i(\bla)$}) the set of all $i$-removable (resp. $i$-addable) nodes for $\bla$. 
We also denote by $\PA_i(\bla)$\index{p@$\PA_i(\bla)$} (resp. $\PR_i(\bla)$\index{p@$\PR_i(\bla)$}) the set of all properly $i$-removable (resp. properly $i$-addable) nodes for $\bla$.

Let $\ttA
\in \PR_i(\bla)$ and $\ttB\in\PA_i(\la)$ be nodes in a component $t$. 
We define 
\begin{align}
\label{EEtaMDef}
\eta_\ttA(\bla)&:=\sharp\{\ttC\in\Rem_i(\bla)\mid \text{$\ttC>\ttA$}\}
-\sharp\{\ttC\in\Add_i(\bla)\mid \text{$\ttC>\ttA$}\},
\index{e@$\eta_\ttA(\bla)$}
\\
\zeta_\ttA(\bla)&:=\zeta_\ttA(\la^{(t)}),
\index{z@$\zeta_\ttA(\bla)$}
\quad
d_\ttA(\bla):= q_i^{\eta_\ttA(\bla)} \zeta_\ttA(\bla).
\index{d@$d_\ttA(\bla)$}
\\
\eta^\ttB(\bla)&:=\sharp\{\ttC\in\Add_i(\bla)\mid \text{$\ttC<\ttB$}\}
-\sharp\{\ttC\in\Rem_i(\bla)\mid \text{$\ttC<\ttB$}\},
\index{e@$\eta^\ttB(\bla)$}
\\
 \label{EZetaEta}
\zeta^\ttB(\bla)&:=\zeta^\ttB(\la^{(t)}),
\index{z@$\zeta^\ttB(\bla)$}
\quad d^\ttB(\la):= q_i^{\eta^\ttB(\la)} \zeta^\ttB(\la).
\index{d@$d^\ttB(\la)$}
\end{align}

\begin{Example} 
{\rm 
Let $\ell=N=2$.
The $2$-multipartition $\bla=((5,5);(6))$ is $5$-strict,
and the residues of its boxes are labeled on the diagram below: 

\vspace{2mm}
\begin{align*}
&\begin{ytableau}
$0$ & $1$ & $2$ & $1$ & $0$ \cr 
$0$ & $1$ & $2$ & $1$ & $0$ \cr
\end{ytableau}\\
&\\
&\begin{ytableau}
$0$ & $1$ & $2$ & $1$ &$0$ & $0$\cr
\end{ytableau}
\end{align*}

\vspace{5mm}
\noindent
The $0$-removable nodes are marked as $\ttA_r$ ($\ttA_1$ is not properly removable), and $0$-addable nodes are marked as $\ttB_s$:

\vspace{1mm}
\begin{align*}
&\ytableausetup{mathmode}
\begin{ytableau}
\, &  &  &  &  & \none[\ttB_3] \cr 
 &  &  &  &  \ttA_3 \cr
\none[\ttB_2] \cr
\end{ytableau}\\
&\\
&\ytableausetup{mathmode}
\begin{ytableau}
\, &  &  &  & \ttA_1 & \ttA_2 \cr
\none[\ttB_1] \cr
\end{ytableau}
\end{align*}

\vspace{1mm}
\noindent
We have:
\begin{align*}&d_{\ttA_2}(\bla)=q^{-1},\ d_{\ttA_3}(\bla)=q^{-1}(1-q^{4}),\\ 
&d^{\ttB_1}(\bla)=1,\ d^{\ttB_2}(\bla)=q^{-1},\ d^{\ttB_3}(\bla)=q^{-1}(1-q^4).
\end{align*}
}
\end{Example}

\vspace{2mm}
Let $\bla\in\Par_p^N(n)$. 
A {\em $\bla$-tableau}\index{tableau} is a bijection $\T: \{1,\dots,n\}\to \bla$, where we have identified $\bla$ with the set of nodes of its Young diagram. In this case we denote
$$
\bi^\T:=(\Res \T(1))\,\cdots\, (\Res \T(n))\in I^n.\index{i@$\bi^\T$}
$$
A $\bla$-tableau $\T$ is called {\em $p$-standard}\index{p@$p$-standard tableau} if 
$\T(\{1,\dots,k\})$ is a $p$-strict multipartition for all $k=1,\dots,n$. 
In this case, the restriction $\T_{\leq k}:=\T|_{\{1,\dots,k\}}$ can be considered as a $p$-standard $\T(\{1,\dots,k\})$-tableau.
We denote by $\Std_p(\bla)$\index{s@$\Std_p(\bla)$} the set of all $p$-standard $\bla$-tableaux. 
For $\bj\in I^n$, we denote
$$
\Std_p(\bla,\bj)=\{\T\in\Std_p(\bla)\mid \bi^\T=\bj\}.\index{s@$\Std_p(\bla,\bj)$}
$$ 
The {\em degree}\index{degree} of a $p$-standard $\bla$-tableau $\T$ is defined as 
\begin{equation}\label{EDegM}
\deg(\T):=\prod_{k=1}^n d^{\T(k)}(\T(\{1,\dots,k-1\})).
\index{d@$\deg(\T)$}
\end{equation}
Note that 
$$\deg(\T)=d^{\T(n)}(\T(\{1,\dots,n-1\}))\deg(\T_{\leq n-1}).$$

\section{Lie theory}
\label{ChLie}

\subsection{Weyl group and positive roots}
\label{SSLT}
We continue with the Lie theoretic set up of \S\ref{ChBasicNotLie}. In addition, recall from \cite[\S\S3.7,3.13]{Kac} the (affine) {\em Weyl group}\index{Weyl group} $W$\index{w@$W$} generated by the fundamental reflections $\{r_i\mid i\in I\}$\index{r@$r_i$} as a Coxeter group. The form $(.|.)$ is $W$-invariant, see \cite[Proposition 3.9]{Kac}.

\begin{Lemma} \label{LMultiple} 
Let $N\in\Z_{\geq 1}$ and $\la,\la_1,\dots,\la_N\in W\La_0$. If $N\la=\la_1+\dots+\la_N$ then $\la_1=\dots=\la_N=\la$. 
\end{Lemma}
\begin{proof}
By assumption, $\la=w\La_0$, $\la_1=w_1\La_0,\dots,\la_N=w_N\La_0$ for some $w,w_1,\dots,w_N\in W$. Using \cite[(6.2.2)]{Kac}, we have 
$$(\la|\la)=(w\La_0|w\La_0)=(\La_0|\La_0)=0.$$ Similarly 
$$
(\la_1|\la_1)=\dots=(\la_N|\la_N)=0.
$$ 
Moreover, if $u\La_0\neq \La_0$ for some $u\in W$ then $u\La_0=\la_0-\theta$ for some $\theta=\sum_{i\in I}m_i\al_i\in Q_+$ with $m_0>0$, and so $(\La_0|u\La_0)<0$. We conclude that if $w_t\La_0\neq w_s\La_0$ then $(w_t\La_0|w_s\La_0)<0$. In other words, if $\la_t\neq \la_s$ then $(\la_t|\la_s)<0$. The equality 
$$
0=(N\la|N\la)=(\sum_{t=1}^N\la_t|\sum_{t=1}^N\la_t)
$$
now implies that $\la_1=\dots=\la_N$ and the lemma.
\end{proof}

As in \cite[\S5]{Kac}, the set $\Phi$\index{f@$\Phi$} of {\em roots}\index{root} of $\g$ is a disjoint union of the set $\Phi^\im=\{n\de\mid n\in \Z\}$\index{f@$\Phi^\im$} of {\em imaginary roots} \index{imaginary root}and the set $\Phi^\re$\index{f@$\Phi^\re$} of {\em real roots}\index{real root}. The real roots are exactly the roots in $\Phi$ which are $W$-conjugate to simple roots. 

Let $\Phi'$\index{f@$\Phi'$} be the root system of type $C_\ell$ whose Dynkin diagram is obtained by dropping the simple root $\al_0$ from our type $A_{2\ell}^{(2)}$ Dynkin diagram. Then $\Phi'=\Phi'_{\text{s}}\sqcup \Phi'_{\text{l}}$ where $\Phi'_{\text{s}}=\{\al\in\Phi'\mid (\al|\al)=4\}$ and $\Phi'_{\text{l}}=\{\al\in\Phi'\mid (\al|\al)=8\}$. 
By \cite[\S6]{Kac}, we have $\Phi^\re=\Phi^\re_{\text{s}}\sqcup \Phi^\re_{\text{m}}\sqcup \Phi^\re_{\text{l}}$ for 
\begin{align}
\Phi^\re_{\text{s}}&= \{(\al+(2n-1)\de)/2 \mid \al\in\Phi'_l,\, n\in\Z\},
\label{EPhiS}
\\
\Phi^\re_{\text{m}}&=\{\al+n\de \mid \al\in\Phi'_s,\, n\in\Z\},
\label{EPhiM}
\\ 
\Phi^\re_{\text{l}}&=\{\al+2n\de \mid \al\in\Phi'_l,\, n\in\Z\}.
\label{EPhiL}
\end{align}
The set of positive roots is then 
$$\Phi_+\index{f@$\Phi_+$}=\Phi_+^\im\sqcup\Phi_+^\re,$$ 
where $$\Phi_+^\im=\{n\de\mid n\in\Z_{>0}\},$$ while $\Phi_+^\re$ consists of the roots in $\Phi'_+$ together with the roots in (\ref{EPhiS})-(\ref{EPhiL}) with $n\in\Z_{>0}$, cf. \cite[Proposition 6.3]{Kac}.

We consider the set of {\em indivisible positive roots}:
\begin{equation}\label{EPsi}
\Psi:=\Phi_+^\re\cup\{\de\}.\index{$\Psi$}\index{y@$\Psi$}
\end{equation}

Let $\theta\in Q_+$. A {\em root partition of $\theta$}\index{root partition} is a pair $(\um,\bmu)$, where $\um=(m_\be)_{\be\in\Psi}$ is a tuple of non-negative integers such that $\sum_{\be\in\Psi}m_\be\be=\theta$, and $\bmu\in\Par^\ell(m_\de)$ is an $\ell$-multipartition of $m_\de$, see \S\ref{ChBasicNotPar}. Denote by $\Par(\theta)$\index{p@$\Par(\theta)$} the set of all root partitions of~$\theta$.

\subsection{Quantized enveloping algebra}

Let $U_q(\g)$\index{u@$U_q(\g)$} be the {\em quantized enveloping algebra of type $A_{2\ell}^{(2)}$},\index{quantized enveloping algebra} i.e the the associative unital $\Q(q)$-algebra with generators $\{E_i,F_i,K_i^{\pm 1}\mid i\in I\}$\index{e@$E_i$}\index{f@$F_i$}\index{k@$K_i$} subject only to the quantum Serre relations:
\begin{eqnarray*}
T_iE_jT_i^{-1}&=&q^{(\al_i|\al_j)}E_j,
\\
T_iF_jT_i^{-1}&=&q^{-(\al_i|\al_j)}F_j, 
\\
E_iF_j-F_jE_i&=&\de_{i,j}\frac{T_i-T_i^{-1}}{q_i-q_i^{-1}},
\\
\sum_{r=0}^{1-a_{i,j}}(-1)^rE_i^{(r)}E_jE_i^{(1-a_{i,j}-r)}&=&0\qquad(i\neq j),
\\
\sum_{r=0}^{1-a_{i,j}}(-1)^rF_i^{(r)}F_jF_i^{(1-a_{i,j}-r)}&=&0\qquad(i\neq j),
\end{eqnarray*}
where we have set
$$
T_i\index{t@$T_i$}:=K_i^{(\al_i|\al_i)/2},\quad E_i^{(r)}:=\frac{E_i}{[r]_i^!},\quad F_i^{(r)}:=\frac{F_i}{[r]_i^!}.
\index{e@$E_i^{(r)}$}\index{f@$F_i^{(r)}$}
$$

We denote by $U_q^-(\g)$ the subalgebra of $U_q(\g)$ generated by the $F_i$. We have 
$$U_q^-(\g)=\bigoplus_{\theta\in Q_+}  U_q^-(\g)_\theta,$$ 
where $U_q^-(\g)_\theta$ is the span of the monomials $F_{i_1}\cdots F_{i_n}$ such that $\al_{i_1}+\dots+\al_{i_n}=\theta$. Root partitions label the elements of a PBW basis of $U_q^-(\g)$ so for any $\theta\in Q_+$ we have 
\begin{equation}\label{EDimRP}
\dim U_q^-(\g)_\theta=|\Par(\theta)|.
\end{equation}
We will also consider the $\Z[q,q^{-1}]$-subalgebra $U^-_{\Z[q,q^{-1}]}(\g)\subseteq U_q^-(\g)$ generated by all $F_i^{(r)}$'s.

We consider $U_q(\g)$ as a Hopf algebra with respect to the coproduct given for all $i\in I$ as follows (cf. \cite[(2.2.3)]{KMPY}):
\begin{equation}\label{ECoproduct}
\De\,:\,K_i\mapsto K_i\otimes K_i,\quad E_i\mapsto E_i\otimes 1+T_i^{-1}\otimes E_i,\quad
F_i\mapsto F_i\otimes T_i+1\otimes F_i.
\index{d@$\De$}
\end{equation}

For $\La\in P_+$, we denote by $V(\La)$\index{v@$V(\La)$} the {\em irreducible integrable module} for $U_q(\g)$ of highest weight $\La$. We fix a non-zero highest weight vector $v_+\in V(\La)_\La$\index{v@$v_+$}, so $E_i v_+=0$ and $T_iv_+=q^{(\al_i|\La)}$ for all $i\in I$.

The formal character of the irreducible module $V(\La_0)$ is well understood. The reference to \cite{Kac} in the following lemma applies because of \cite[Theorem 4.12 and \S4.14]{Lusztig}.

\begin{Lemma} {\rm \cite[(12.6.1),(12.6.2),(12.13.5)]{Kac}}\label{LKac}
The weights of $V(\La_0)$ are of the form $
w\La_0-d\de$ with $w\in W$ and $d\in\Z_{\geq 0}$. Moreover, 
$$\dim V(\La_0)_{w\La_0-d\de}=|\Par^\ell(d)|.$$ 
\end{Lemma}

There is a ($\Q(q)$-linear) anti-involution $\boldsymbol\sigma:U_q(\g)\to U_q(\g)$ with 
\begin{align*}
\boldsymbol\sigma\,:\, E_i&\mapsto q_iF_iT_i^{-1}=q_i^{-1}T_i^{-1}F_i,\quad 
F_i\mapsto q_i^{-1}T_iE_i=q_iE_iT_i,
\quad 
T_i\mapsto T_i.
\index{s@$\boldsymbol\sigma$}
\end{align*}

There is a unique symmetric bilinear form $(\cdot,\cdot)$ \index{$(\cdot,\cdot)$}on $V(\La)$ such that  
\begin{equation}\label{EF()}
(v_+,v_+)=1\ \text{and $(xv,w)=(v,\boldsymbol\sigma(x)w)$ for all $x\in U_q(\g)$ and $v,w\in V(\La)$},
\end{equation} 
cf. \cite[Appendix D]{KMPY}. 
We will need to evaluate an arbitrary $$
(F_{i_n}\cdots F_{i_1}v_+,F_{j_n}\cdots F_{j_1}v_+)
$$
for the case $\La=N\La_0$ with $N\in\Z_{\geq 1}$ in terms of combinatorics of ($p$-strict) $N$-multipartitions, see Proposition~\ref{PForm}.

The following standard fact is easy to check using the commutation formula from \cite[Lemma 1.7]{Jantzen}:

\begin{Lemma} \label{LRefl} 
Let $\La\in P_+$, $\mu\in P$ and $v$ be a non-zero vector of the weight space $V(\La)_\mu$. If $i\in I$ is such that  $E_iv=0$, then $a:=(\mu|\al_i^\vee)\geq 0$, $F_i^{(a)}v\neq 0$, and $(F_i^{(a)}v,F_i^{(a)}v)=(v,v)$.
\end{Lemma}

\begin{Lemma} \label{LWExt} 
Let $\La\in P_+$ and $w\in W$ with a reduced decomposition $w=r_{i_l}\cdots r_{i_1}$. Then  $$a_k:=(r_{i_{k-1}}\cdots r_{i_1}\La|\al_{i_k}^\vee)\geq 0 \qquad(k=1,\dots,l),
$$
$F_{i_l}^{(a_l)}\cdots F_{i_1}^{(a_1)}v_+$ is a non-zero vector of the weight space $V(\La)_{w\La}$, and $$(F_{i_l}^{(a_l)}\cdots F_{i_1}^{(a_1)}v_+,F_{i_l}^{(a_l)}\cdots F_{i_1}^{(a_1)}v_+)=1.$$ 
\end{Lemma}
\begin{proof}
Induction on $l=0,1\dots$, the base being trivial. Let $l>0$ and set $u:=r_{i_{l-1}}\cdots r_{i_1}$, $v:=F_{i_{l-1}}^{(a_{l-1})}\cdots F_{i_1}^{(a_1)}v_+$. 
By \cite[Lemma 3.11]{Kac}, $u^{-1}\al_{i_l}$ is a positive root. So $(u\La|\al_{i_l})=(\La|u^{-1}\al_{i_l})\geq 0$, and $u\La+\al_{i_l}=u(\La+u^{-1}\al_{i_l})$ is not a weight of $V(\La)$. Hence  $E_{i_l}v=0$. We can now apply Lemma~\ref{LRefl} and the inductive assumption. 
\end{proof}

\subsection{Fock space}

Recall the combinatorial notions defined in \S\ref{ChComb}.  

The  {\em ($q$-deformed) level $1$ Fock space}\index{Fock space} $\Fock$\index{f@$\Fock$}, as defined in \cite{KMPY}, see also \cite{LT}, is the $\Q(q)$ vector space with basis $\{u_\la\mid \la\in \Par_p\}$\index{u@$u_\la$} labeled by the $p$-strict partitions:
$$
\Fock:=\bigoplus_{\la\in\Par_p} \Q(q)\cdot  u_\la.
$$
There is a structure of a $U_q(\g)$-module on $\Fock$ such that $u_\varnothing$ is a highest weight vector, with the submodule $U_q(\g)\cdot u_\varnothing \subseteq \Fock$ isomorphic to $V(\La_0)$, and 
\begin{eqnarray}\label{EActionFockE}
E_iu_\la&=&\sum_{\ttA\in\PR_i(\la)}d_\ttA(\la)u_{\la_\ttA},\\
\label{EActionFockF}
F_iu_\la&=&\sum_{\ttB\in\PA_i(\la)}d^\ttB(\la)u_{\la^\ttB},
\\
\label{EActionFockT}
T_iu_\la&=&
q^{(\al_i|\La_0-\cont(\la))}u_\la.
\end{eqnarray}
For instance, in the set up of Example~\ref{Ex1} we have $$F_0u_{(5,5,2)}=(1-q^4)u_{(6,5,2)}+u_{(5,5,2,1)}.$$

Moreover, as established in \cite[Appendix D]{KMPY}, there is a bilinear form $(\cdot,\cdot)$\index{$(\cdot,\cdot)$} on $\Fock$ which satisfies
\begin{equation}\label{EOrth}
(u_\la,u_\mu)=\de_{\la,\mu}\norm{\la}
\end{equation}
and $(xv,w)=(v,\boldsymbol\sigma(x)w)$ for all $x\in U_q(\g)$ and $v,w\in\Fock$. 
 
For $N\in\Z_{\geq 1}$, we now define $\Fock^N:=\Fock^{\otimes N}$. \index{f@$\Fock^N$}This is a $U_q(\g)$-module via the coproduct (\ref{ECoproduct}). For $\bla=(\la^{(1)},\dots,\la^{(N)})$ we denote 
$$
u_\bla:=u_{\la^{(1)}}\otimes\dots\otimes u_{\la^{(1)}}\in\Fock^N.
\index{u@$u_\bla$}
$$
The formulas  (\ref{EActionFockE})-(\ref{EActionFockT}) and (\ref{ECoproduct}) imply:
\begin{eqnarray}\label{EActionFockEM}
E_iu_\bla&=&\sum_{\ttA\in\PR_i(\bla)}d_\ttA(\bla)u_{\bla_\ttA},\\
\label{EActionFockFM}
F_iu_\bla&=&\sum_{\ttB\in\PA_i(\bla)}d^\ttB(\bla)u_{\bla^\ttB},
\\
\label{EActionFockTM}
T_iu_\bla&=&
q^{(\al_i|N\La_0-\cont(\bla))}u_\bla.
\end{eqnarray}
The bilinear form (\ref{EOrth}) is now extended to $\Fock^N$ as follows:
\begin{equation}\label{EOrthM}
(u_\bla,u_\bmu)=\de_{\bla,\bmu}\norm{\bla},\index{$(\cdot,\cdot)$}
\end{equation}
and it satisfies $(xv,w)=(v,\boldsymbol\sigma(x)w)$ for all $x\in U_q(\g)$ and $v,w\in\Fock^N$. Comparing with (\ref{EF()}), we deduce:

\begin{Lemma} \label{LForm} 
The form $(\cdot,\cdot)$ on $V(N\La_0)$ is the restriction of the form $(\cdot,\cdot)$ on $\Fock^N$ to $V(N\La_0)=U_q(\g)\cdot u_{(\varnothing,\dots,\varnothing)} \subseteq \Fock^N$.
\end{Lemma}

\begin{Proposition} \label{PForm} 
Let $\bi=i_1\cdots i_d,\,\bj=j_1\cdots j_n\in I^n$ and $v_+\in V(N\La_0)$ be a highest weight vector with $(v_+,v_+)=1$. Then 
$$
(F_{i_n}\cdots F_{i_1}v_+,F_{j_n}\cdots F_{j_1}v_+)=\sum_{\bla\in\Par_p^N(n)}\sum_{\substack{\Stab\in\Std_p(\bla,\bi)\\ \T\in\Std_p(\bla,\bj)}} \deg(\Stab)\deg(\T)\norm{\bla}.
$$
\end{Proposition}
\begin{proof}
By Lemma~\ref{LForm} and (\ref{EOrthM}), it suffices to prove that 
$$
F_{i_n}\cdots F_{i_1}u_{(\varnothing,\dots,\varnothing)}=\sum_{\bla\in\Par_p^N(n)}\bigg(\sum_{\substack{\T\in\Std_p(\bla,\bi)}} \deg(\Stab)\bigg)u_\bla,
$$
which follows by induction on $n$, using (\ref{EActionFockFM}) and (\ref{EDegM}). 
\end{proof}

\chapter{Quiver Hecke superalgebras}
\label{Part2}

\section{Quiver Hecke superalgebras and a dimension formula}

\subsection{Definition of quiver Hecke superalgebras}
\label{SQHDef}
In this subsection we introduce the quiver Hecke (KLR) superalgebra of type $A_{2\ell}^{(2)}$ defined in \cite{KKT}. 
Recall the Lie theoretic set up of \S\ref{ChBasicNotLie}. 
We assign a parity (i.e. an element of $\Z/2=\{\0,\1\}$) to the set $I$ of vertices of the Dynkin diagram as follows:
$$
|i|:=
\left\{
\begin{array}{ll}
\1 &\hbox{if $i=0$,}\\
\0 &\hbox{if $i=1,\dots,\ell$.}
\end{array}
\right.
\index{$\lvert i\rvert$}
$$

For $i,j\in I$, we define polynomials $Q_{i,j}(u,v)\in\F[u,v]$ as follows:
$$
Q_{i,j}(u,v):=\index{q@$Q_{i,j}(u,v)$}
\left\{
\begin{array}{ll}
0 &\hbox{if $i=j$,}\\
1 &\hbox{if $|i-j|>1$,}\\
u-v &\hbox{if $1\leq i=j-1<\ell-1$,}\\
v-u &\hbox{if $1\leq j=i-1<\ell-1$,}\\
u^2-v &\hbox{if $\ell>1$, and $(i,j)=(0,1)$ or $(\ell-1,\ell)$,}\\
v^2-u &\hbox{if $\ell>1$, and $(i,j)=(1,0)$ or $(\ell,\ell-1)$,}\\
u^4-v &\hbox{if $\ell=1$, and $(i,j)=(0,1)$,}\\
v^4-u &\hbox{if $\ell=1$, and $(i,j)=(1,0)$,}\\
\end{array}
\right.
$$
For $i,j,k\in I$, we also define polynomials $B_{i,j,k}(u,v)\in\F[u,v]$ as follows:
\begin{equation*}
B_{i,j,k}(u,v)\index{b@$B_{i,j,k}(u,v)$}
:=\left\{
\begin{array}{ll}
-1 &\hbox{if $i=k=j+1$,}\\
1 &\hbox{if $i=k= j-1\not\in \{0,\ell-1\}$,}\\
(u+v) &\hbox{if $i=k= j-1=\ell-1>0$,}\\
(v-u) &\hbox{if $i=k=j-1=0$ and $\ell>1$,}\\
(u^2+v^2)(v-u) &\hbox{if $i=k= j-1=0$ and $\ell=1$,}\\
0 &\hbox{otherwise,}
\end{array}
\right.
\label{EB}
\end{equation*}

Let $\theta\in Q_+$ with $\height(\theta)=n$. 
The {\em quiver Hecke superalgebra}\index{quiver Hecke superalgebra} $R_\theta$\index{r@$R_\theta$} is the unital graded 
superalgebra generated by the elements 
$$
\{e(\bi)\mid\bi\in I^\theta\}\cup \{y_1,\dots,y_n\}\cup \{\psi_1,\dots,\psi_{n-1}\}
\index{e@$e(\bi)$}\index{y@$y_r$}\index{y@$\psi_r$}\index{$\psi_r$}
$$
and the following defining relations (for all admissible $r,s,\bi$, etc.)
\begin{equation}
\label{R1}
e(\bi)e(\bj)=\de_{\bi,\bj}e(\bi),
\end{equation}
\begin{equation}
\sum_{\bi\in I^\theta}e(\bi)=1,
\label{R2}
\end{equation}
\begin{equation}
\label{R2.5}
y_re(\bi)=e(\bi)y_r,
\end{equation}
\begin{equation}
\label{R2.75}
\psi_re(\bi)=e(s_r\cdot \bi)\psi_r,
\end{equation}
\begin{equation}
y_ry_se(\bi)=
\left\{
\begin{array}{ll}
-y_sy_re(\bi) &\hbox{if $r\neq s$ and $|i_r|=|i_s|=\1$,}\\
 y_sy_re(\bi)&\hbox{otherwise,}
\end{array}
\right.
\label{R3}
\end{equation}
\begin{equation}
\psi_r y_se(\bi)=(-1)^{|i_r||i_{r+1}||i_s|}y_s \psi_r e(\bi)\qquad(s\neq r,r+1),
\label{R4}
\end{equation}
\vspace{2mm}
\begin{equation}
\begin{split}
(\psi_ry_{r+1}-(-1)^{|i_r||i_{r+1}|}y_r \psi_r) e(\bi)&=
(y_{r+1} \psi_r-(-1)^{|i_r||i_{r+1}|} \psi_ry_r) e(\bi)
\\&=
\left\{
\begin{array}{ll}
e(\bi) &\hbox{if $i_r=i_{r+1}$,}\\
0 &\hbox{otherwise,}
\end{array}
\right.
\end{split}
\label{R5}
\end{equation}
\vspace{2mm}
\begin{equation}\label{R6}
\psi_r^2e(\bi)=Q_{i_r,r_{r+1}}(y_r,y_{r+1}),
\end{equation}
\begin{equation}\label{R65}
\psi_r\psi_s e(\bi)=(-1)^{|i_r||i_{r+1}||i_s||i_{s+1}|}\psi_s\psi_r e(\bi)\qquad(|r-s|>1),
\end{equation}
\begin{equation}
(\psi_{r+1}\psi_r\psi_{r+1}-\psi_{r}\psi_{r+1} \psi_r) e(\bi)=B_{i_r,i_{r+1},i_{r+2}}(y_r,y_{r+2})e(\bi)
\label{R7}
\end{equation}
The structure of a graded superalgebra on $R_\theta$ is defined by setting 
\begin{align}
\bideg(\bi)&:=(0,\0),
\label{EGrading1}
\\ \bideg(y_se(\bi))&:=((\al_{i_s}|\al_{i_s}),|i_s|),
\label{EGrading2}
\\ \bideg(\psi_r e(\bi))&:=-((\al_{i_r}|\al_{i_{r+1}}),|i_r||i_{r+1}|).
\label{EGrading3}
\end{align}

Sometimes we denote the identity in $R_\theta$ by $e_\theta$\index{e@$e_\theta$} rather than $1$.

\begin{Remark} 
{\rm 
In the Introduction, the algebra $R_\theta$ was denoted $\cR_\theta$ to distinguish it from its purely even analogue---KLR algebra. Since KLR algebras do not arise in the main body of this paper, we will always use the notation $R_\theta$ for quiver Hecke superalgebras from now on.
}
\end{Remark}

We will use the usual Khovanov-Lauda diagrams \cite{KL1} to represent the elements of $R_\theta$, so if $\height(\theta)=n$, $\bi=i_1\cdots i_n\in I^\theta$, $1\leq r<n$ and $1\leq s\leq n$, we have 
\vspace{2mm}
\begin{align*}
e(\bi)&=
\begin{braid}\tikzset{baseline=3mm}
  \draw (0,0)node[below]{$i_1$}--(0,3);
  \draw (1,0)node[below]{$i_2$}--(1,3);
  \draw[dots] (1.5,2.9)--(4.7,2.9);
  \draw (5,0)node[below]{$i_{n}$}--(5,3);
  \draw[dots] (1.5,0)--(4.7,0);
 \end{braid},
 \\ 
y_s&=
 \begin{braid}\tikzset{baseline=3mm}
  \draw (0,0)node[below]{$i_1$}--(0,3);
  \draw (3,0)node[below]{$i_{s-1}$}--(3,3);
  \draw[dots] (0.5,2.9)--(2.7,2.9);
  \draw (4.4,0)node[below]{$i_{s}$}--(4.4,3);
  \draw (5.7,0)node[below]{$i_{s+1}$}--(5.7,3);
  \draw[dots] (0.5,0)--(2.7,0);
  \draw[dots] (6.1,2.9)--(8.3,2.9);
  \draw[dots] (6.1,0)--(8.3,0);
  \draw (8.5,0)node[below]{$i_{n}$}--(8.5,3);
\blackdot (4.4,1.5);
 \end{braid},
\\ 
\psi_r e(\bi)&=
 \begin{braid}\tikzset{baseline=3mm}
  \draw (0,0)node[below]{$i_1$}--(0,3);
  \draw (3,0)node[below]{$i_{r-1}$}--(3,3);
  \draw[dots] (0.5,2.9)--(2.7,2.9);
  \draw (4.6,0)node[below]{$i_{r}$}--(5.8,3);
  \draw (5.8,0)node[below]{$i_{r+1}$}--(4.6,3);
  \draw (7.4,0)node[below]{$i_{r+2}$}--(7.4,3);
  \draw[dots] (0.5,0)--(2.7,0);
  \draw[dots] (7.8,2.9)--(10,2.9);
  \draw[dots] (7.8,0)--(10,0);
  \draw (10.2,0)node[below]{$i_{n}$}--(10.2,3);
\end{braid}.
\end{align*}

\vspace{3mm}
For every $w\in\Si_n$, we choose a reduced decomposition $w=s_{r_1}\dots s_{r_l}$ and define $\psi_{w}:=\psi_{r_1}\cdots\psi_{r_l}$.\index{$\psi_{w}$}\index{y@$\psi_{w}$} In general $\psi_w$ depends on the choice of a reduced decomposition for $w$, but:

\begin{Lemma} \label{LBKW} {\rm \cite[Proposition 2.5]{BKW}}
Let $\theta\in Q_+$ with $\height(\theta)=n$, $\bi\in I^\theta$, and $w$ be an element of $\Si_n$ written as a product of simple transpositions:  $w=s_{t_1}\dots s_{t_m}$ for some $1\leq t_1\dots,t_m<n$. 
\begin{enumerate}
\item[{\rm (i)}] If the decomposition $w=s_{t_1}\dots s_{t_m}$ is reduced, and $w=s_{r_1}\dots s_{r_m}$ is another reduced decomposition of $w$, then in $R_\theta$ we have 
$$\psi_{t_1}\dots\psi_{t_m}e(\bi)=\psi_{r_1}\dots\psi_{r_m}e(\bi)+(*),$$
where $(*)$ is a linear combination of elements of the form $\psi_uf(y)e(\bi)$ such that $u< w$ and $f(y)$ is a polynomial in $y_1,\dots,y_n$.



\item[{\rm (ii)}] If the decomposition $w=s_{t_1}\dots s_{t_m}$ is not reduced, then   $\psi_{t_1}\dots\psi_{t_m}e(\bi)$ can be written as a linear combination of elements of the form
$$\psi_{t_{a_1}}\dots\psi_{t_{a_b}}f(y)e(\bi)$$ 
such that $1\leq a_1<\dots<a_b\leq m$, $b<m$, $s_{t_{a_1}}\dots s_{t_{a_b}}$ is a reduced word and $f(y)$ is a polynomial in $y_1,\dots,y_n$.
\end{enumerate}
\end{Lemma}

\begin{Theorem}\label{TBasis}{\cite[Corollary 3.15]{KKT}} 
Let $\theta\in Q_+$ and $n=\height(\theta)$. Then  
\begin{align*}
&\{\psi_w y_1^{k_1}\dots y_n^{k_n}e(\bi)\mid w\in {\Si}_n,\ k_1,\dots,k_n\in\Z_{\geq 0}, \ \bi\in I^\theta\},
\\
&\{ y_1^{k_1}\dots y_n^{k_n}\psi_we(\bi)\mid w\in {\Si}_n,\ k_1,\dots,k_n\in\Z_{\geq 0}, \ \bi\in I^\theta\}
\end{align*}
are bases of  $R_\theta$. 
\end{Theorem}

By \cite[Proposition 6.15]{HW}, \cite[\S4.2]{KKO}, \cite[Theorem 8.6]{KKOII}, we have:

\begin{Lemma} \label{LTypeM} 
Every irreducible graded $R_\theta$-supermodule is finite dimensional and irreducible as an $R_\theta$-module. 
\end{Lemma}

Moreover, by (\ref{EDimRP}) and \cite[Corollary 10.3]{KKO}, we have:

\begin{Lemma} \label{LAmount}
We have $|\Irr(R_\theta)|=|\Par(\theta)|$.  
\end{Lemma}



Let $\La=\sum_{i\in I}a_i\La_i\in P_+.$ 
The {\em cyclotomic quiver Hecke superalgebra $R_\theta^\La$} \index{cyclotomic quiver Hecke superalgebra}\index{r@$R_\theta^\La$}is defined as $R_\theta$ modulo the relations
\begin{equation}
y_1^{a_{i_1}}e(\bi)=0\qquad(\text{for all}\ \bi=i_1\cdots i_n\in I^\theta).
\label{RCyc}
\end{equation}
For example, in the important special case $\La=\La_0$, the cyclotomic relations (\ref{RCyc}) are equivalent to $y_1=0$ and  $e(\bi)=0$ if $i_1\neq 0$. We have the natural projection maps
\begin{equation}\label{EPi}
\pi_\theta^\La \colon R_\theta\onto R_\theta^{\La}.\index{p@$\pi_\theta^\La$}
\end{equation} 
Inflating along $\pi_\theta^\La$, every graded $R_\theta^{\La}$-supermodule can be considered as a graded $R_\theta$-supermodule. From Lemma~\ref{LTypeM}, we now get:

\begin{Lemma} \label{LTypeMCyc} 
Every irreducible graded $R_\theta^\La$-supermodule is  irreducible as an $R_\theta^\La$-module. 
\end{Lemma}

We have antiautomorphisms 
\begin{equation}\label{ETauAntiAuto}
\tau:R_\theta\to R_\theta\quad\text{and}\quad\tau:R_\theta^\La\to R_\theta^\La\index{t@$\tau$}
\end{equation}
which are identity on the generators.


\subsection{Induction and Restriction}
\label{SSIndRes}
Let $\theta_1,\dots,\theta_n\in Q_+$ and $\theta=\theta_1+\dots+\theta_n$. Denote 
$$\underline{\theta}:=(\theta_1,\dots,\theta_n)\in Q_+^n$$ 
and 
$$R_{\underline{\theta}}\index{r@$R_{\underline{\theta}}$}
=R_{\theta_1,\dots,\theta_n}:=R_{\theta_1}\otimes\dots\otimes R_{\theta_n}$$ (tensor product of superalgebras). 
Let
$$
e_{\underline{\theta}}
\index{e@$e_{\underline{\theta}}$}
=e_{\theta_1,\dots,\theta_n}:=\sum_{\bi^1\in I^{\theta_1},\dots,\bi^n\in I^{\theta_n}}e(\bi^1\cdots\bi^n)\in R_{\theta}.
$$
By \cite[\S4.1]{KKO}, we have the natural embedding 
\begin{equation}\label{EIota}
\iota_{\underline{\theta}}:R_{\underline{\theta}}\to e_{\underline{\theta}}R_{\theta}e_{\underline{\theta}},
\index{i@$\iota_{\underline{\theta}}$}
\end{equation}
and we identify $R_{\underline{\theta}}$ with a subalgebra of $e_{\underline{\theta}}R_{\theta}e_{\underline{\theta}}$ via this embedding. We refer to this subalgebra as a {\em parabolic subalgebra}\index{parabolic subalgebra}. Note that $\iota_{\underline{\theta}}(e_{\theta_1}\otimes\dots\otimes e_{\theta_n})=e_{\underline{\theta}}$.

There are the exact functors preserving finite dimensional modules:  
\begin{align*}
\Ind_{\underline{\theta}} 
\index{i@$\Ind_{\underline{\theta}}$}
&:= R_{\theta} e_{\underline{\theta}}
\otimes_{R_{\underline{\theta}}} ?:\mod{R_{\underline{\theta}}} \rightarrow \mod{R_{\theta}},\\
\Res_{\underline{\theta}} 
\index{r@$\Res_{\underline{\theta}}$}
&:= e_{\underline{\theta}} R_{\theta}
\otimes_{R_{\theta}} ?:\mod{R_{\theta}}\rightarrow \mod{R_{\underline{\theta}}}.
\end{align*}
The functor $\Ind_{\underline{\theta}}
$ is left adjoint to $\Res_{\underline{\theta}}
$. 
If $M_1\in\Mod{R_{\theta_1}},\dots,M_n\in\Mod{R_{\theta_n}}$, we define 
\begin{equation*}\label{ECircProd}
M_1\circ\dots\circ M_n:=\Ind_{\underline{\theta}}
M_1\boxtimes\dots\boxtimes M_n. 
\index{$\circ$}
\end{equation*}

The functors of induction and restriction have obvious parabolic analogues. Given a family $(\theta^a_b)_{1\leq a\leq n,\ 1\leq b\leq m}$ of elements of $Q_+$, set 
$\eta_b:=\sum_{a=1}^n\theta^{a}_b$ for all $1\leq b\leq m$. Then we have  functors
\begin{align*}
\Ind_{\theta^{1}_1,\dots,\theta^{n}_{1}\,;\,\dots\,;\,\theta^{1}_{m},\dots,\theta^{n}_{m}}^{\,\eta_1;\,\dots\,;\,\eta_m}
&:
\mod{R_{\theta^{1}_1,\dots,\theta^{n}_{1}\,,\,\dots\,,\,\theta^{1}_{m},\dots,\theta^{n}_{m}}}\to\mod{R_{\,\eta_1,\,\dots\,,\,\eta_m}}
\\
\Res_{\theta^{1}_1,\dots,\theta^{n}_{1}\,;\,\dots\,;\,\theta^{1}_{m},\dots,\theta^{n}_{m}}^{\,\eta_1;\,\dots\,;\,\eta_m}
&:\mod{R_{\,\eta_1,\,\dots\,,\,\eta_m}}\to \mod{R_{\theta^{1}_1,\dots,\theta^{n}_{1}\,,\,\dots\,,\,\theta^{1}_{m},\dots,\theta^{n}_{m}}}
\end{align*}

For $\si\in \Si_n$ and $\underline{\theta}=(\theta_1,\dots,\theta_n)\in Q_+^n$, let 
$$
\si\underline{\theta}:=(\theta_{\si^{-1}(1)},\dots,\theta_{\si^{-1}(n)}), 
$$
and 
\begin{align*}
s(\si,\underline{\theta})&:=-\sum_{1\leq m<k\leq n,\ \si(m)>\si(k)}(\theta_m \mid \theta_k)\in\Z,
\\
t(\si,\underline{\theta})&:=-\sum_{1\leq m<k\leq n,\ \si(m)>\si(k)}|\theta_m|\,|\theta_k|\in\Z/2.
\end{align*}
where for $\theta=\sum_{i\in I}n_i\al_i$ we have set $|\theta|:=n_0$. 
There is a superalgebra isomorphism 
$$
\phi^\si:R_{\si\underline{\theta}}\to R_{\underline{\theta}},\ 
x_1\otimes\dots\otimes x_n\mapsto 
(-1)^{\sum_{1\leq a<c\leq n, w(a)>w(c)}|x_a||x_c|}
x_{\si(1)}\otimes\dots\otimes x_{\si(n)}.
$$
Composing with this isomorphism, we get a functor
$$
\mod{R_{\underline{\theta}}}\to \mod{R_{\si\underline{\theta}}},\  M\mapsto M^{\phi^\si}.
$$
Making in addition degree and parity shifts, we get a functor 
\begin{equation}\label{ETwist}
\mod{R_{\underline{\theta}}}\to \mod{R_{\si\underline{\theta}}},\  M\mapsto {}^\si M:=\Pi^{t(\si,\underline{\theta})}Q^{s(\si,\underline{\theta})}M^{\phi^\si}.
\end{equation}

The Mackey Theorem below for $m=n=2$ follows from \cite[Proposition 4.5]{KKO}. The general case can be deduced from the case $m=n=2$ by induction. See also \cite[Proposition 3.7]{Evseev} or the proof of \cite[Proposition 2.18]{KL1}.

\begin{Theorem} \label{TMackeyKL} 
{\rm (Mackey Theorem)}
Let $\underline{\theta}=(\theta_1,\dots,\theta_n)\in Q_+^n$ and $\underline{\eta}=(\eta_1,\dots,\eta_m)\in Q_+^m$ with 
$$\theta:=\theta_1+\dots+\theta_n=\eta_1+\dots+\eta_m.$$ 
\begin{enumerate}
\item[{\rm (i)}] Let $k:=\height(\theta)$, and consider the compositions 
$\la:=(\height(\theta_1),\dots,\height(\theta_n))$ and  $\mu:=(\height(\eta_1),\dots,\height(\eta_m))$ of $k$. Then
$$e_{\underline{\eta}}R_\theta e_{\underline{\theta}}=\sum_{w\in{}^\mu\D^{\la}}R_{\underline{\eta}}\psi_w R_{\underline{\theta}}.
$$ 
\item[{\rm (ii)}] For any $M\in\mod{R_{\underline{\theta}}}$ we have that $\Res_{\underline{\eta}}\,\Ind_{\underline{\theta}} M$ has  filtration with factors of the form
$$
\Ind_{\ga^{1}_1,\dots,\ga^{n}_{1}\,;\,\dots\,;\,\ga^{1}_{m},\dots,\ga^{n}_{m}}^{\,\eta_1;\,\dots\,;\,\eta_m}
{}^{\si(\underline{\ga})}\big(\Res_{\ga^{1}_1,\dots,\ga^{1}_{m}\,;\,\dots\,;\,\ga^{n}_{1},\dots,\ga^{n}_{m}}^{\,\theta_1;\,\dots\,;\,\theta_n}
\,M \big)
$$
with $\underline{\ga}=(\ga^a_b)_{1\leq a\leq n,\ 1\leq b\leq m}$ running over all tuples of elements of $Q_+$ such that  $\sum_{b=1}^m\ga^{a}_b=\theta_a$ for all $1\leq a\leq n$ and $\sum_{a=1}^n\ga^{a}_b=\eta_b$ for all $1\leq b\leq m$, and $\si(\underline{\ga})$ is the permutation of $mn$ which maps 
$$
(\ga^{1}_1,\dots,\ga^{1}_{m};\ga^{2}_1,\dots,\ga^{2}_{m};\dots;\ga^{n}_{1},\dots,\ga^{n}_{m})
\mapsto
(\ga^{1}_1,\dots,\ga^{n}_{1};\ga^{1}_2,\dots,\ga^{n}_{2};\dots;\ga^{1}_{m},\dots,\ga^{n}_{m}).
$$
\end{enumerate} 
\end{Theorem}

\vspace{2mm}
A special case where $n=2$ and $\height(\eta_k)=1$ for all $k$ yields: 

\vspace{1mm}
\begin{Corollary} \label{CShuffle}
Let $M\in\mod{R_{\theta}}$, $N\in\mod{R_{\eta}}$, and $\bi\in I^{\theta+\eta}$. Then $e(\bi)(M\circ N)\neq 0$ if and only if $\bi$ is a shuffle of words $\bj\in I^\theta$ and $\bk\in I^\eta$ such that $e(\bj)M\neq 0$ and $e(\bk)N\neq 0$. 
\end{Corollary}


Let $\theta,\eta\in Q_+$ and $\La\in P_+$. 
Recalling (\ref{EPi}), we have the {\em cyclotomic parabolic subalgebra}\index{cyclotomic parabolic subalgebra}
\begin{equation}\label{ECycPar}
R^{\La}_{\theta,\eta}:=\pi_\theta^\La(R_{\theta,\eta})\subseteq e_{\theta,\eta}R^{\La}_\theta e_{\theta,\eta}.
\index{r@$R^{\La}_{\theta,\eta}$}
\end{equation}
We have a natural embedding 
$$\zeta_{\theta,\eta}\colon R_{\theta+\eta}\to R_{\theta,\eta}, \ x\mapsto 
x\otimes e_\eta.
\index{z@$\zeta_{\theta,\eta}$}
$$ 
The map $\pi_{\theta+\eta}^\La\circ \zeta_{\theta,\eta}$ factors through the quotient $R_\theta^{\La}$ to give the natural {\em unital} algebra homomorphism
\begin{equation}\label{EZetaHom}
\zeta_{\theta,\eta}^\La\colon R_\theta^{\La}\to R_{\theta,\eta}^{\La}.
\index{z@$\zeta_{\theta,\eta}^\La$}
\end{equation}


\subsection{KKO-supercategorification theorem}
We now review the categorification theory of Kang-Kashiwara-Oh \cite{KKO,KKOII}. 

Let $\La\in P_+$ and $\theta,\eta\in Q_+$. 
Composing the homomorphism from (\ref{EZetaHom}) with the embedding (\ref{ECycPar}), we get the algebra homomorphism 
$$
R_\theta^\La\to e_{\theta,\eta}R_{\theta+\eta}^\La e_{\theta,\eta}.
$$ 
This homomorphism makes $e_{\theta,\eta}R^\La_{\theta+\eta}$ into an $(R^\La_\theta, R^\La_{\theta+\eta})$-bimodule, and $R^\La_{\theta+\eta}e_{\theta,\eta}$ into an $(R^\La_{\theta+\eta},R^\La_\theta)$-bimodule. 
Specializing to $\eta=\al_i$ for some $i\in I$, we consider the functors
\begin{align*}
\funE_i^{\La}
\index{e@$\funE_i^{\La}$}
&:=q_i^{1-\langle h_i,\La-\theta\rangle}e_{\theta,\al_i}R^\La_{\theta+\al_i}\otimes_{R^\La_{\theta+\al_i}}-:\Mod{R^\La_{\theta+\al_i}}\to\Mod{R^\La_\theta},
\\
\funF_i^{\La}
\index{f@$\funF_i^{\La}$}
&:=R^\La_{\theta+\al_i}e_{\theta,\al_i}\otimes_{R^\La_{\theta}}-:\Mod{R^\La_{\theta}}\to\Mod{R^\La_{\theta+\al_i}},
\\
\funK_i^{\La}&:=q_i^{\langle h_i,\La-\theta\rangle}:\Mod{R^\La_{\theta}}\to\Mod{R^\La_{\theta}}.
\index{k@$\funK_i^{\La}$}
\end{align*}
Note the degree shift in the definition of the functor $\funE_i^\La$. 
These functors are exact, preserve finite dimensionality, and map projective modules to projective modules. In particular, recalling  (\ref{EGrothGroup}), they induce $\Z[q,q^{-1}]$-linear operators 
$$E_i^{\La}:=[\funE_i^{\La}],\index{e@$E_i^{\La}$} \quad F_i^{\La}:=[\funF_i^{\La}]\index{f@$F_i^{\La}$}\quad \text{and}\quad K_i^{\La}:=[\funK_i^{\La}]\index{k@$K_i^{\La}$}
$$ on the Grothendieck group 
$$
[\underlineproj{R^\La}]_q:=\bigoplus_{\theta\in Q_+}[\underlineproj{R^\La_\theta}]_q,
$$
see \cite[Theorem 8.9]{KKO}. 

Given a finitely generated projective graded $R_\theta^\La$-supermodule $P$, we denote by $P^\tau$\index{p@$P^\tau$} the finitely generated right graded  $R_\theta^\La$-supermodule which is $P$ as a graded superspace and the action is given by $v\cdot x=\tau(x)v$ for all $x\in R_\theta^\La$ and $v\in P$. For example, $(R_\theta^\La e(\bi))^\tau\cong e(\bi)R_\theta^\La$.  
There is a bilinear form on $[\underlineproj{R^\La_\theta}]_q$, defined from 
\begin{equation}\label{EFG}
([P],[P'])=\dim_q(P^\tau\otimes_{R^\La_\theta} P')\index{$(\cdot,\cdot)$}
\end{equation}
for $P,P'\in \proj{R^\La_\theta}$.

\begin{Theorem} \label{TCat} 
Let $\La\in P_+$. There is an isomorphism
$$\iota:V(\La)\iso [\underlineproj{R^\La}]_{\Q(q)}
\index{i@$\iota$}
$$
with the following properties:
\begin{enumerate}
\item[{\rm (1)}] $\iota(E_i v)=E_i^\La(\iota(v))$, $\iota(F_i v)=F_i^\La(\iota(v))$, $\iota(K_i v)=K_i^\La(\iota(v))$ for all $i\in I$ and all $v\in V(\La)$; in particular we have 
$$\iota(V(\La)_{\La-\theta})=[\proj{R^\La_\theta}]_{\Q(q)},
$$ 
and $v_+:=\iota^{-1}([R^\La_0])$ is a highest weight vector of $V(\La)$.
\item[{\rm (2)}] $\iota^{-1}( [\proj{R^\La}]_q)=
U^-_{\Z[q,q^{-1}]}(\g)\cdot v_+$. 
\item[{\rm (3)}] $(v,w)=(\iota(v),\iota(w))$, where on the left we have the form defined in (\ref{EF()}) and on the right the form defined in (\ref{EFG}).
\end{enumerate}
\end{Theorem}
\begin{proof}
Everything but (3) is contained in \cite[Theorem 10.2]{KKO}. For (3), we clearly have $([R^\La_0],[R^\La_0])=1$, so it suffices to check the property (\ref{EF()}) for the form (\ref{EFG}). The property (\ref{EF()}) can be checked just for the generators $E_i,F_i,K_i$. The property is in fact clear for the $K_i$, and by the symmetricity of the form, the property follows for $E_i$'s if we can check it for $F_i$'s. To check the property for $F_i$, for $P\in \proj{R^\La_\theta}$ and $P'\in\proj{R^\La_{\theta+\al_i}}$ we  compute:
\begin{align*}
(F_i^\La[P],[P'])&=([\funF_i^\La P],P')
\\
&=\dim_q\big((\funF_i^\La P)^\tau\otimes_{R_{\theta+\al_i}^\La}P'\big)
\\
&=\dim_q\big((R_{\theta+\al_i}^\La e_{\theta,\al_i}\otimes_{R_{\theta}^\La} P)^\tau\otimes_{R_{\theta+\al_i}^\La}P'\big)
\\
&=\dim_q\big((P^\tau\otimes_{R_{\theta}^\La} e_{\theta,\al_i}R_{\theta+\al_i}^\La)\otimes_{R_{\theta+\al_i}^\La}P'\big)
\\
&=\dim_q\big(P^\tau\otimes_{R_{\theta}^\La} (e_{\theta,\al_i}R_{\theta+\al_i}^\La\otimes_{R_{\theta+\al_i}^\La}P')\big)
\\
&=\dim_q\big(P^\tau\otimes_{R_{\theta}^\La} (q_i^{-1+\lan h_i,\la_\theta\ran}q_i^{1-\lan h_i,\la_\theta\ran}e_{\theta,\al_i}R_{\theta+\al_i}^\La\otimes_{R_{\theta+\al_i}^\La}P')\big)
\\
&=\dim_q\big(P^\tau\otimes_{R_{\theta}^\La} (q_i^{-1+\lan h_i,\la_\theta\ran}\funE_i^\La P')\big)
\\
&=([P], q_i^{-1}T_i^\La E_i^\La [P']\big),
\end{align*}
which is what we need, since $\si(F_i)=q_i^{-1}T_iE_i$.
\end{proof}

\subsection{Graded dimension formula for $R_\theta^{N\La_0}$}\label{subsec:dim_form}
Throughout the subsection we fix $N\in\Z_{> 0}$. 

\begin{Theorem} \label{TDim} 
For $\theta\in Q_+$ set $n=\height(\theta)$, and let $\bi,\bj\in I^\theta$. Then 
$$
\dim_qe(\bi)R^{N\La_0}_\theta e(\bj)=\sum_{\bla\in\Par_p^N(n)}\sum_{\substack{\Stab\in \Std_p(\bla,\bi)\\ \T\in\Std_p(\bla,\bj)}} \deg(\Stab)\deg(\T)\norm{\bla}.
$$
\end{Theorem}
\begin{proof}
Writing $\bi=i_1\cdots i_n$, note for any $\La\in P_+$ that 
$$R_\theta^\La e(\bi)\cong \funF^\La_{i_n}\cdots \funF^\La_{i_1} R_0^\La.$$  
On the other hand, 
$$
e(\bi)R^\La_\theta e(\bj)\cong e(\bi)R^\La_\theta\otimes_{R^\La_\theta} R^\La_\theta e(\bj)
\cong (R^\La_\theta e(\bi))^\tau\otimes_{R^\La_\theta} R^\La_\theta e(\bj).
$$
So, by Theorem~\ref{TCat}, 
\begin{equation}\label{EDim}
\dim_q e(\bi)R^\La_\theta e(\bj)=([R^\La_\theta e(\bi)],[R^\La_\theta e(\bj)])=(F_{i_n}\cdots F_{i_1}v_+,F_{j_n}\cdots F_{j_1}v_+).
\end{equation}
The theorem now follows from Proposition~\ref{PForm}. 
\end{proof}

\begin{Corollary} \label{CIdNZ} 
Let $\theta\in Q_+$ and $\bi\in I^\theta$. If $e(\bi)\neq 0$ in $R^{N\La_0}_\theta$ then 
$\Std_p(\bla,\bi)\neq \varnothing$ for some $\bla\in\Par_p^N$. 
\end{Corollary}
\begin{proof}
If $e(\bi)\neq 0$ in $R^{N\La_0}_\theta$ then $e(\bi)R^{N\La_0}_\theta e(\bi)\neq 0$, so the result follows immediately from Theorem~\ref{TDim}. 
\end{proof}

\begin{Example} 
{\rm 
Suppose $p=3$. Then $\Par_3(3)=\{(3),(2,1)\}$. We have  
$$\norm{(3)}=1+q^2,\quad \norm{(2,1)}=1,$$ and 
\begin{align*}
\Std_p((3))=\Big\{\Stab:=
\vcenter{\hbox{\begin{ytableau}
$1$ & $2$ & $3$\cr 
\end{ytableau}}}\,\Big\},\quad  
\Std_p((2,1))=\Big\{\T:=
\vcenter{\hbox{\begin{ytableau}
$1$ & $2$ \cr
 $3$ \cr
\end{ytableau}}}\,\Big\}.
\end{align*}
Note that 
\begin{align*}
\bi^\Stab=010=\bi^\T, 
\quad
\deg(\Stab)=q,
\quad 
\deg(\T)=1.
\end{align*}
From the theorem, it follows that $e(010)$ is the only non-zero standard idempotent in $R^{\La_0}_{2\al_0+\al_1}$, and 
\begin{align*}
\dim_q R^{\La_0}_{2\al_0+\al_1}&=\dim_q e(010)R^{\La_0}_{2\al_0+\al_1}e(010)
\\
&=q^2(1+q^2)+1
=1+q^2+q^4.
\end{align*}
From this it is not hard to deduce using defining relations that $\{1,y_3,y_3^2\}$ is a basis of $R^{\La_0}_{2\al_0+\al_1}$. 
}
\end{Example}

\vspace{1mm}
\begin{Example} 
{\rm 
Suppose $p=3$. Then $\Par_3(4)=\{(4),(3,1)\}$. We have   
$$\norm{(4)}=1,\quad \norm{(3,1)}=1+q^2,$$ 
and 
\begin{align*}
\Std_p((4))&=\Big\{\Stab:=
\vcenter{\hbox{\begin{ytableau}
$1$ & $2$ & $3$& $4$\cr 
\end{ytableau}}}\,\Big\},\  
\\
\Std_p((3,1))&=\Big\{\T_1:=
\vcenter{\hbox{\begin{ytableau}
$1$ & $2$ & $4$\cr
 $3$ \cr
\end{ytableau}}},\,  
\T_2:=
\vcenter{\hbox{\begin{ytableau}
$1$ & $2$ & $3$\cr
 $4$ \cr
\end{ytableau}}}\,\Big\}.
\end{align*}
Note that 
\begin{align*}
&\bi^\Stab=\bi^{\T_1}=\bi^{\T_2}=0100, 
\\
&\deg(\Stab)=q(1+q^2),  \quad
\deg(\T_1)=q^{-1},
\quad 
\deg(\T_2)=q.
\end{align*}
From the theorem, it follows that $e(0100)$ is the only non-zero standard idempotent in $R^{\La_0}_{3\al_0+\al_1}$, and 
\begin{align*}
\dim_q R^{\La_0}_{3\al_0+\al_1}&=\dim_q e(0100)R^{\La_0}_{3\al_0+\al_1}e(0100)
\\
&=(q(1+q^2))^2+(1+q^2)(q+q^{-1})^2
\\
&=(1+q^2+q^4)(q+q^{-1})^2.
\end{align*}
}
\end{Example}

\vspace{2mm}
\begin{Example}
Suppose $p=3$. We compute $\dim_q R_{4\al_0+2\al_1}$. Let 
$$\bi:=010010\quad \text{and}\quad \bj:=010001.$$ 
In Table I below, we list the partitions $\la\in\Par_3(6)$ together with $\norm{\la}$ as well as the corresponding standard tableaux $\T\in\Std_p(\la)$ with $\bi^\T$ and $\deg(\T)$ (recall that $q_0=q$ and $q_1=q^4$).
So
\vspace{1mm}
\begin{align*}
\dim_q e(\bj)R_{4\al_0+2\al_1}e(\bj)=\,&((q^5+q^3+q)(1+q^2))^2+
((q+q^{-1}+q^{-3})(1+q^2))^2,
\end{align*}
\begin{align*}
\dim_q e(\bi)R_{4\al_0+2\al_1}e(\bi)=\,&(1+q^2)(q^2(1+q^2))^2+
(q(1+q^2))^2+((q^3+q)(1+q^2))^2
\\
&+(1+q^2)(1-q^4)(1+q^2)^2+(1+q^2)(q+q^{-1})^2,
\\
\dim_q e(\bi)R_{4\al_0+2\al_1}e(\bj)=
\,&\dim_q e(\bj)R_{4\al_0+2\al_1}e(\bi)
\\
=\,&q(1+q^2)(q^5+q^3+q)(1+q^2)
\\
&+(q^3+q)(1+q^2)(q+q^{-1}+q^3)(1+q^2).
\end{align*}

\vspace{5mm}
\begin{center}
\begin{tabular}{c | c | c | c | c  }
$\la$ & $\norm{\la}$ & $\T\in \Std_p(\la)$ & $\bi^\T$ & $\deg(\T)$
\\ \hline\hline
$(6)$ & $1+q^2$ & $\vcenter{\vspace{0.5mm}\hbox{\begin{ytableau}
$1$ & $2$ & $3$ &$4$ &$5$ &$6$\cr
\end{ytableau}\vspace{0.5mm}}}$
& $\bi$ &$q^2(1+q^2)$
\\ \hline
$(5,1)$ & $1$ & $\vcenter{\vspace{0.5mm}\hbox{\begin{ytableau}
$1$ & $2$ & $3$ & $4$ & $5$\cr  $6$\cr
\end{ytableau}\vspace{0.5mm}}}$
& $\bi$ &$q(1+q^2)$
\\ 
 & & $\vcenter{\vspace{0.5mm}\hbox{\begin{ytableau}
$1$ & $2$ & $3$ & $4$ & $6$\cr  $5$\cr
\end{ytableau}\vspace{0.5mm}}}$
& $\bj$ &$q^5(1+q^2)$
\\ 
 & & $\vcenter{\vspace{0.5mm}\hbox{\begin{ytableau}
$1$ & $2$ & $3$ & $5$ & $6$\cr  $4$\cr
\end{ytableau}\vspace{0.5mm}}}$
& $\bj$ &$q^3(1+q^2)$
\\ 
 & & $\vcenter{\vspace{0.5mm}\hbox{\begin{ytableau}
$1$ & $2$ & $4$ & $5$ & $6$\cr  $3$\cr
\end{ytableau}\vspace{0.5mm}}}$
& $\bj$ &$q(1+q^2)$
\\ \hline
$(4,2)$ & $1$ & $\vcenter{\vspace{0.5mm}\hbox{\begin{ytableau}
$1$ & $2$ & $3$ & $4$ \cr $5$ &  $6$\cr
\end{ytableau}\vspace{0.5mm}}}$ & $\bj$ &$q(1+q^2)$
\\ 
 &  & $\vcenter{\vspace{0.5mm}\hbox{\begin{ytableau}
$1$ & $2$ & $3$ & $5$ \cr $4$ &  $6$\cr
\end{ytableau}\vspace{0.5mm}}}$
& $\bj$ &$q^{-1}(1+q^2)$
\\ 
 &  & $\vcenter{\vspace{0.5mm}\hbox{\begin{ytableau}
$1$ & $2$ & $3$ & $6$ \cr $4$ &  $5$\cr
\end{ytableau}\vspace{0.5mm}}}$
& $\bi$ &$q^3(1+q^2)$
\\ 
 &  & $\vcenter{\vspace{0.5mm}\hbox{\begin{ytableau}
$1$ & $2$ & $4$ & $5$ \cr $3$ &  $6$\cr
\end{ytableau}\vspace{0.5mm}}}$
& $\bj$ &$q^{-3}(1+q^2)$
\\ 
 &  & $\vcenter{\vspace{0.5mm}\hbox{\begin{ytableau}
$1$ & $2$ & $4$ & $6$ \cr $3$ &  $5$\cr
\end{ytableau}\vspace{0.5mm}}}$
& $\bi$ &$q(1+q^2)$
\\ \hline
$(3,3)$ & $(1+q^2)(1-q^4)$ & $\vcenter{\vspace{0.5mm}\hbox{\begin{ytableau}
$1$ & $2$ & $3$ \cr $4$ & $5$ &  $6$\cr
\end{ytableau}\vspace{0.5mm}}}$
& $\bi$ &$q^2$
\\ 
 &  & $\vcenter{\vspace{0.5mm}\hbox{\begin{ytableau}
$1$ & $2$ & $4$ \cr $3$ & $5$ &  $6$\cr
\end{ytableau}\vspace{0.5mm}}}$
& $\bi$ &$1$
\\ \hline
$(3,2,1)$ & $(1+q^2)$ & $\vcenter{\vspace{0.5mm}\hbox{\begin{ytableau}
$1$ & $2$ & $3$ \cr $4$ & $5$ \cr  $6$\cr
\end{ytableau}\vspace{0.5mm}}}$
& $\bi$ &$q$
\\ 
 &  & $\vcenter{\vspace{0.5mm}\hbox{\begin{ytableau}
$1$ & $2$ & $4$ \cr $3$ & $5$ \cr  $6$\cr
\end{ytableau}\vspace{0.5mm}}}$
& $\bi$ &$q^{-1}$
\\
\hline
\end{tabular}
\end{center}
\vspace{3mm}
\centerline{\sc Table I.}
\end{Example}
 
\newpage
Putting $q$ to $1$ in the right hand side of the dimension formula of Theorem~\ref{TDim} yields the formula for the ungraded dimension $\dim e(\bi)R^{N\La_0}_\theta e(\bj)$. For $N=1$, this ungraded dimension has been essentially computed in \cite[Theorem 3.4]{AP}, but in slightly different combinatorial terms, cf. Corollary~\ref{CDimAP} below. We now explain how to deduce the Ariki-Park dimension formula from our graded dimension formula. 

First, observe from (\ref{ELaNorm}) that\, $\norm{\la}|_{q=1}=0$ unless $\la\in\Par_0$, i.e. $\la$ is {\em strict}, which means that all parts of $\la$ are distinct. Hence $\norm{\bla}|_{q=1}=0$ unless $\bla=(\la^{(1)},\dots,\la^{(N)})$ is {\em strict}, which means that every component $\la^{(t)}$ is strict. 
Thus we are left with summing only over strict multipartitions. Let $\T$ be a standard $\bla$-tableau for some strict multipartition $\bla$ of $n$. If for some $m<n$ the multipartition $\T(\{1,\dots,m\})$ is not strict, we pick a minimal such $m$ and note that then $d^{\T(m)}(\T(\{1,\dots,m\}))$ is divisible by $q^2-1$ and so $\deg(\T)|_{q=1}=0$. Thus we are left with the subset $\Std_0(\bla)$ of $0$-standard $\bla$-tableaux, i.e. such that each $\T(\{1,\dots,m\})$ is a strict multipartition. 

To complete the proof of the fact that our graded dimension formula specializes to the Ariki-Park dimension formula, it remains to observe 
that for a strict partition $\bla$ of residue content $\theta$ and strictly standard $\la$-tableaux $\Stab$ and $\T$ we have 
$$
(\deg(\Stab)\deg(\T)\norm{\bla})|_{q=1}=2^{m_0-h(\bla)},
$$
where $m_0$ is defined from $\theta=\sum_{i\in I}m_i\al_i$, and $h(\bla):=\sum_{t=1}^Nh(\la^{(t)})$ for $\bla=(\la^{(1)},\dots,\la^{(N)})$.  
Thus:

\begin{Corollary} \label{CDimAP} 
For $\theta=\sum_{i\in I}m_i(\theta)\al_i\in Q_+$ set $\height(\theta)=n$,  and let $\bi,\bj\in I^\theta$. Then 
$$
\dim e(\bi)R^{N\La_0}_\theta e(\bj)=\sum_{\bla\in\Par_0^N(n)}\sum_{\substack{\Stab\in\Std_0(\bla,\bi)\\ \T\in\Std_0(\bla,\bj)}} 2^{m_0(\theta)-h(\bla)}.
$$
\end{Corollary}

\begin{Corollary} \label{CWeights} 
Let $\theta\in Q_+$. Then $V(N\La_0)_{N\La_0-\theta}\neq 0$ if and only $R^\La_\theta\neq 0$ if and only if there exists $\bla\in\Par_p^N$ with $\cont(\bla)=\theta$. 
\end{Corollary}
\begin{proof}
The first equivalence comes from Theorem~\ref{TCat}(1). Moreover, by Corollary~\ref{CDimAP}, $R^\La_\theta\neq 0$ if and only if there exists $\bla\in\Par_0^N$ with $\cont(\bla)=\theta$. It remains to notice that this is equivalent to the fact that there exists $\bla\in\Par_p^N$ with $\cont(\bla)=\theta$ (if $\bla\in\Par_p^N$ is not $0$-strict, in each component remove all parts divisible by $p$, add a new part equal to the sum of the removed parts and reorder to get a $0$-strict multipartition of the same content). 
\end{proof}

Recall the notation (\ref{EParRhoD}) and (\ref{ECores}).

\begin{Lemma} \label{LCombLie} 
Consider the map 
$$\ka:\Par_p\to P,\ \la\mapsto \La_0-\cont(\la)
\index{k@$\ka$}
$$
from the set of $p$-strict partitions to the set of weights. Then:
\begin{enumerate}
\item[{\rm (i)}] $\Im(\ka)=\{w\La_0-d\de\mid w\in W,\ d\in\Z_{\geq 0}\}$.
\item[{\rm (ii)}] $\ka$ restricts to a bijection between $\Cores_p$ and $W\La_0$.
\item[{\rm (iii)}] For any $\rho\in \Cores_p$ and $d\in\Z_{\geq 0}$, we have $\Par_p(\rho,d)=\ka^{-1}(\ka(\rho)-d\de)$.
\end{enumerate}
\end{Lemma}
\begin{proof}
(i) By Corollary~\ref{CWeights}, for any $\theta\in Q_+$, we have that $\La_0-\theta$ is a weight of $V(\La_0)$ if and only if there exists a partition $\la\in\Par_p$ with $\cont(\la)=\theta$. By Lemma~\ref{LKac}, the set of weights of $V(\La_0)$ is exactly $\{w\La_0-d\de\mid w\in W,\ d\in\Z_{\geq 0}\}$. 
Part  (i) follows. 

(ii) Observe that on the one hand, $\{w\La_0\mid w\in W\}$ is the set of all weights $\mu$ of $V(\la)$ such that $\mu+\de$ is not a weight of $V(\La_0)$. On the other hand, by Lemma~\ref{LCoreNoDeSubtr}, the set $\Cores_p$ of all $\bar p$-cores is the set of all $p$-strict partitions $\rho$ such that $\cont(\rho)-\de$ is not a content of any $p$-strict partition. Now part (ii) follows from part (i). 

(iii) Part (iii) follows from parts (i) and (ii).
\end{proof}

\begin{Lemma} \label{LNumberOfIrrCyc} 
Let $\rho\in\Cores_p$ and $d\in\Z_{\geq 0}$. Then $|\Irr(R^{\La_0}_{\cont(\rho)+d\de})|=|\Par^\ell(d)|$ . 
\end{Lemma}
\begin{proof}
By Theorem~\ref{TCat}(1), we have $$|\Irr(R^{\La_0}_{\cont(\rho)+d\de})|=\dim V(\La_0)_{\La_0-\cont(\rho)-d\de}.$$ 
By Lemma~\ref{LCombLie}(ii), we have that $\La_0-\cont(\rho)=w\La_0$ for some $w\in W$. It remains to apply Lemma~\ref{LKac}. 
\end{proof}

\begin{Lemma} \label{LRhoNUnique} 
Let $\rho\in\Cores_p$, $N\in\Z_{\geq 1}$ and $\rho^N:=(\rho,\dots,\rho)\in\Par_p^N$. If  $\bmu
\in\Par_p^N$ satisfies $\cont(\bmu)=N\cont(\rho)$ then 
$\bmu=\rho^N$. 
\end{Lemma}
\begin{proof}
By Lemma~\ref{LCombLie}(ii), $\La_0-\cont(\rho)= w\La_0$ for some $w\in W$. So $N\La_0-N\cont(\rho)=w(N\La_0)$. In particular, $N\La_0-N\cont(\rho)+\de$ is not a weight of $V(N\La_0)$. 
Note that the components $\mu^{(1)},\dots,\mu^{(N)}$ of $\bmu$ must be $\bar p$-cores, for otherwise, using Corollary~\ref{CWeights}, we conclude that $N\La_0-N\cont(\rho)+\de=N\La_0-\cont(\bmu)+\de$ is a weight of $V(N\La_0)$ giving a contradiction. So, by Lemma~\ref{LCombLie}(ii), for all $t=1,\dots,N$, we have $\La_0-\cont(\mu^{(t)})=w_t\La_0$ for some $w_t\in W$. Hence 
\begin{align*}
Nw\La_0&=N\La_0-N\cont(\rho)=
N\La_0-\cont(\bmu)
\\
&=(\La_0-\cont(\mu^{(1)}))+\dots+(\La_0-\cont(\mu^{(N)}))=
w_1\La_0+\dots+w_N\La_0.
\end{align*} 
By Lemma~\ref{LMultiple}, $w\La_0=w_1\La_0=\dots=w_N\La_0$, hence $\cont(\rho)=\cont(\mu^{(1)})=\dots=\cont(\mu^{(N)})$, therefore $\mu^{(1)}=\dots=\mu^{(N)}=\rho$ by Lemma~\ref{LCombLie}(ii) again. 
\end{proof}

\begin{Lemma} \label{LNonZeroSpecialCut} 
If $\rho\in \Cores_p$, $N\in\Z_{>0}$ and $d\in\Z_{\geq 0}$, then $$e_{N\cont(\rho),d\de}R_{N\cont(\rho)+d\de}^{N\La_0}e_{N\cont(\rho),d\de}\neq 0.$$ 
\end{Lemma}
\begin{proof}
We apply Corollary~\ref{CDimAP}. 
Let $r:=|\rho|$ and $\la\in\Par_0(r+dp)$ be a partition obtained from $\rho$ by adding a part equal to $dp$. It is easy to see that  there exists $\T\in\Std_0(\la)$ with $\T(\{1,\dots,r\})=\rho$.  Now setting $\bla:=(\la,\rho,\dots,\rho)\in\Par_0^N$ it is clear that there exists $\Stab\in\Std_0(\bla)$ with $\Stab(\{1,\dots,Nr\})=(\rho,\dots,\rho)$. It remains to note that $e_{N\cont(\rho),d\de}e(\bi^\Stab) e_{N\cont(\rho),d\de}=e(\bi^\Stab)$, and $e(\bi^\Stab)R_{N\cont(\rho)+d\de}^{N\La_0}e(\bi^\Stab)\neq 0$ by Corollary~\ref{CDimAP}. 
\end{proof}

\section{Further properties of quiver Hecke superalgebras}
\subsection{Divided power idempotents}\label{SSDPI}
Let $i\in I$ and $m\in \Z_{\ge 1}$. We denote by $w_0$ the longest element of $\Si_m$. If $i\neq 0$, the algebra $R_{m\al_i}$ is known to be the {\em nil-Hecke algebra}\index{nil-Hecke algebra}  
and has an idempotent 
$
e(i^{(m)}):=\psi_{w_0}\prod_{s=1}^m y_s^{s-1},\index{e@$e(i^{(m)})$}
$ 
cf. \cite{KL1}. If $i=0$, the algebra $R_{m\al_0}$ is known to be the {\em odd nil-Hecke algebra}\index{odd nil-Hecke algebra} and has an idempotent 
$e(0^{(m)})$\index{e@$e(0^{(m)})$} of the form $\pm\psi_{w_0}\prod_{s=1}^m y_s^{s-1}$, cf. \cite{EKL}. 
Diagrammatically we will denote 
$$e(i^{(m)}):=
\begin{braid}\tikzset{baseline=0.25em}
	\braidbox{0}{1.3}{0}{1}{$i^m$}
\end{braid}
=
\begin{braid}\tikzset{baseline=0.25em}
	\braidbox{0}{2.2}{0}{1}{$i\dots i$}
\end{braid}.
$$
For example, for $m=2$, we have for any $i\in I$:
\begin{equation}\label{EDivDiag}
e(i^{(2)}):=
\begin{braid}\tikzset{baseline=0.25em}
	\braidbox{0}{1.3}{0}{1}{$i^2$}
\end{braid}
=
\begin{braid}\tikzset{baseline=0.25em}
	\braidbox{0}{1.6}{0}{1}{$i\,\,i$}
\end{braid}
=
\begin{braid}\tikzset{baseline=3mm}
  \draw (0,0)node[below]{$i$}--(1,2);
  \draw (1,0)node[below]{$i$}--(0,2);
\blackdot (0.8,0.45);
 \end{braid}.
\end{equation}

\begin{Lemma} \label{L150921_2}
For any $i,j\in I$, we have 
$$
\begin{braid}\tikzset{baseline=1em}
	\braidbox{0}{1.3}{-0.5}{0.5}{\normalsize$i\,\, i$}
	\draw (2,-0.1) node{\normalsize$j$};
	\draw (1.7,3) node{\normalsize$i\,\,i$};
	\draw (0,3) node{\normalsize$j$};
	\draw(0.2,0.5)--(1.2,2.5);
	\draw(1.1,0.5)--(2.1,2.5);
	\draw(2,0.5)--(0.2,2.5);
\end{braid}
=
\begin{braid}\tikzset{baseline=1em}
	\braidbox{0}{1.3}{-0.5}{0.5}{\normalsize$i\,\, i$}
	\draw (2,-0.1) node{\normalsize$j$};
	\braidbox{1}{2.3}{2.5}{3.5}{\normalsize$i\,\, i$}
	\draw (0,3) node{\normalsize$j$};
	\draw(0.2,0.5)--(1.2,2.5);
	\draw(1.1,0.5)--(2.1,2.5);
	\draw(2,0.5)--(0.2,2.5);
\end{braid}
 \qquad\text{and}
 \qquad
 \begin{braid}\tikzset{baseline=1em}
	\braidbox{1}{2.3}{-0.5}{0.5}{\normalsize$i\,\, i$};
	\draw (0,-0.1) node{\normalsize$j$};
	\draw (0.7,3) node{\normalsize$i\,\,i$};
	\draw (2,3) node{\normalsize$j$};
	\draw(0.2,0.5)--(2,2.5);
	\draw(1.1,0.5)--(0.2,2.5);
	\draw(2,0.5)--(1.2,2.5);
\end{braid}
=
\begin{braid}\tikzset{baseline=1em}
	\braidbox{0}{1.3}{2.5}{3.5}{\normalsize$i\,\, i$}
	\draw (0,-0.1) node{\normalsize$j$};
	\draw(1.6,0) node{\normalsize$i\,\, i$};
	\draw (2,3) node{\normalsize$j$};
	\draw(0.2,0.5)--(2,2.5);
	\draw(1.1,0.5)--(0.2,2.5);
	\draw(2,0.5)--(1.2,2.5);
\end{braid}
=
\begin{braid}\tikzset{baseline=1em}
	\braidbox{1}{2.3}{-0.5}{0.5}{\normalsize$i\,\, i$};
	\draw (0,-0.1) node{\normalsize$j$};
	\braidbox{0}{1.3}{2.5}{3.5}{\normalsize$i\,\, i$}
	\draw (2,3) node{\normalsize$j$};
	\draw(0.2,0.5)--(2,2.5);
	\draw(1.1,0.5)--(0.2,2.5);
	\draw(2,0.5)--(1.2,2.5);
\end{braid}
$$
\end{Lemma}
\begin{proof}
We prove the first equality, the second and third ones being similar. We have using (\ref{EDivDiag}), the fact that $e(i^{(2)})$ is an idempotent, and the relations (\ref{R7}), (\ref{R5}):
\begin{align*}
\begin{braid}\tikzset{baseline=1em}
	\braidbox{0}{1.3}{-0.5}{0.5}{\normalsize$i\,\, i$}
	\draw (2,-0.1) node{\normalsize$j$};
	\draw (1.7,3) node{\normalsize$i\,\,i$};
	\draw (0,3) node{\normalsize$j$};
	\draw(0.2,0.5)--(1.2,2.5);
	\draw(1.1,0.5)--(2.1,2.5);
	\draw(2,0.5)--(0.2,2.5);
\end{braid}
&=
\begin{braid}\tikzset{baseline=1em}
	\braidbox{0}{1.3}{-0.5}{0.5}{\normalsize$i\,\, i$}
	\draw (2,-0.1) node{\normalsize$j$};
	\draw (1.7,3) node{\normalsize$i\,\,i$};
	\draw (0,3) node{\normalsize$j$};
	\draw(0.2,0.5)--(2.1,2.5);
	\draw(1.1,0.5)--(0.3,1.5)--(1.2,2.5);
	\draw(2,0.5)--(0.2,2.5);
	\blackdot(0.9,0.76);
\end{braid}
=
\begin{braid}\tikzset{baseline=1em}
	\braidbox{0}{1.3}{-0.5}{0.5}{\normalsize$i\,\, i$}
	\draw (2,-0.1) node{\normalsize$j$};
	\draw (1.7,3) node{\normalsize$i\,\,i$};
	\draw (0,3) node{\normalsize$j$};
	\draw(0.2,0.5)--(2.1,2.5);
	\draw(1.1,0.5)--(2,1.5)--(1.2,2.5);
	\draw(2,0.5)--(0.2,2.5);
	\blackdot(1.3,0.76);
\end{braid}
=
\begin{braid}\tikzset{baseline=1em}
	\braidbox{0}{1.3}{-0.5}{0.5}{\normalsize$i\,\, i$}
	\draw (2,-0.1) node{\normalsize$j$};
	\draw (1.7,3) node{\normalsize$i\,\,i$};
	\draw (0,3) node{\normalsize$j$};
	\draw(0.2,0.5)--(2.1,2.5);
	\draw(1.1,0.5)--(2,1.5)--(1.2,2.5);
	\draw(2,0.5)--(0.2,2.5);
	\blackdot(1.84,1.5);
\end{braid}
\pm
\begin{braid}\tikzset{baseline=1em}
	\braidbox{0}{1.3}{-0.5}{0.5}{\normalsize$i\,\, i$}
	\draw (2,-0.1) node{\normalsize$j$};
	\draw (1.7,3) node{\normalsize$i\,\,i$};
	\draw (0,3) node{\normalsize$j$};
	\draw(0.2,0.5)--(2.1,2.5);
	\draw(1.2,0.6)--(1.2,0.9)--(0.2,2.5);
	\draw(2,0.5)--(1.7,0.7)--(1.2,2.5);
\end{braid}
=\begin{braid}\tikzset{baseline=1em}
	\braidbox{0}{1.3}{-0.5}{0.5}{\normalsize$i\,\, i$}
	\draw (2,-0.1) node{\normalsize$j$};
	\braidbox{1}{2.3}{2.5}{3.5}{\normalsize$i\,\, i$}
	\draw (0,3) node{\normalsize$j$};
	\draw(0.2,0.5)--(1.2,2.5);
	\draw(1.1,0.5)--(2.1,2.5);
	\draw(2,0.5)--(0.2,2.5);
\end{braid}
\end{align*}
since 
$
\begin{braid}\tikzset{baseline=0.25em}
	\braidbox{0}{1.3}{-0.5}{0.5}{\normalsize$i\,\, i$};
	\draw(0.2,0.5)--(1.2,1.5);
	\draw(1.2,0.5)--(0.2,1.5);
\end{braid}=0
$
by (\ref{R6}). 
\end{proof}

Let $\theta\in Q_+$. Recalling the notation $I^{\theta}_{\di}$ from \S\ref{ChBasicNotLie}, fix a divided power word $\bi = (i_1^{(m_1)},\ldots, i_r^{(m_r)}) \in I^{\theta}_{\di}$. We have the {\em divided power idempotent}\index{divided power idempotent} 
$$
e(\bi):=\iota_{m_1\al_{i_1},\dots,m_r\al_{i_r}}(e(i_1^{(m_1)})\otimes\dots\otimes e(i_r^{(m_r)})) \in R_{\theta}.
\index{e@$e(\bi)$}
$$
Setting
\begin{equation}\label{EHatBi}
\hat\bi := (i_1,\ldots,i_1,\ldots,i_r,\ldots,i_r)\in I^{\theta},
\index{i@$\hat\bi$}\index{$\hat\bi$}
\end{equation}
with $i_k$ repeated $m_k$ times, note that $e(\bi) e(\hat\bi) = e(\bi) =e(\hat \bi)e(\bi)$. 
Define 
$$\bi! := [m_1]_{i_1}^! \cdots [m_r]_{i_r}^!
\qquad\text{and}\qquad
\langle \bi \rangle := \sum_{k=1}^r (\al_{i_k}\mid\al_{i_k})m_k (m_k-1)/4.
\index{$\bi^^21$}\index{$\langle \bi \rangle$}\index{i@$\bi^^21$}\index{i@$\langle \bi \rangle$}
$$ 

As in \cite[\S2.5]{KL1}, one proves

\begin{Lemma}\label{LFact}  
 Let $\bi\in I^{\theta}_{\di}$, and\, $U$ (resp.~$W$) be a left (resp.~right)
graded $R_{\theta}$-supermodule with finite dimensional $e(\bi)U$ (resp. $We(\bi)$). 
 Then 
 \[
  \DIM (e(\hat\bi) U) = \bi! q^{\langle \bi\rangle} \DIM(e(\bi) U)
  \quad \text{and} \quad 
  \DIM (W e(\hat\bi))= \bi! q^{-\langle \bi\rangle} \DIM(W e(\bi)). 
 \]
\end{Lemma}

\begin{Lemma} \label{LMatrix} 
Let $\La\in P_+$ and $\theta\in Q_+$ be such that $\La-\theta\in W\La$. Then there exists $\bi\in I^{\theta}_\di$ such that $\dim_q e(\bi)R^{\La}_{\theta}e(\bi)=1$ and there is an isomorphism of graded $\F$-superalgebras $R_{\theta}^{\La}\iso \End_\F(R_{\theta}^{\La}e(\bi))$, where $x\in R_{\theta}^{\La}$ gets mapped to the left multiplication by $x$. 
In particular, $R_{\theta}^{\La}$ is a matrix algebra over $\F$. 
\end{Lemma}
\begin{proof}
Let $v_+=\iota^{-1}([R^{\La}_0])\in V(\La)$ be the highest weight vector of $V(\La)$ corresponding to the trivial module of $R^{\La}_0$ under the categorification as in Theorem~\ref{TCat}(1). By assumption, we can write  $\La-\theta =w\La$ for some $w\in W$. Applying Lemma~\ref{LWExt}, we find $i_1,\dots,i_l\in I$ and $a_1,\dots,a_l\in\Z_{\geq 0}$ such that 
$F_{i_l}^{(a_l)}\cdots F_{i_1}^{(a_1)}v_+\in V(\La)_{\La-\theta}$ and $$
(F_{i_l}^{(a_l)}\cdots F_{i_1}^{(a_1)}v_+, F_{i_l}^{(a_l)}\cdots F_{i_1}^{(a_1)}v_+)=1.
$$
By (\ref{EDim}), we have 
\begin{align*}
\dim_q e(i_1^{a_1}\cdots i_l^{a_l})R^{\La}_\theta e(i_1^{a_1}\cdots i_l^{a_l})
&=(F_{i_l}^{a_l}\cdots F_{i_1}^{a_1}v_+, F_{i_l}^{a_l}\cdots F_{i_1}^{a_1}v_+)
\\
&=([a_1]^!_{i_1})^2\cdots ([a_l]^!_{i_l})^2.
\end{align*}
Let $\bi=i_1^{(a_1)}\cdots  i_1^{(a_1)}\in I^\theta_\di$. Now by  Lemma~\ref{LFact}, we have  $\dim_q e(\bi) R^{\La}_\theta e(\bi)=1$. Moreover, 
$$\dim V(\La)_{\La_0-\theta}=\dim V(\La)_{w\La}=\dim V(\La)_{\La}=1.$$ 
So by Theorem~\ref{TCat}(1), the algebra $R^{\La}_{\theta}$ has only one graded irreducible supermodule $L$. Considering the composition series of the left regular module over $R^{\La}_{\theta}$, we see that $e(\bi)L\neq 0$ (in fact, $\dim e(\bi) L=1$). 
It remains to apply Lemma~\ref{Lem_Mor_to_F}.
\end{proof}

\begin{Corollary} \label{CMatrix} 
Let $\rho\in \Cores_p$ and $N\in\Z_{>0}$. Then there exists $\bi\in I^{N\cont(\rho)}_\di$ such that $\dim_q(e(\bi)R^{N\La_0}_{N\cont(\rho)}e(\bi))=1$ and 
$R^{N\La_0}_{N\cont(\rho)}$ is a matrix algebra over $\F$.
\end{Corollary}
\begin{proof}
We apply Lemma~\ref{LMatrix}, observing by Lemma~\ref{LCombLie}(ii), that $\La_0-\cont(\rho)=w\La_0$ for some $w\in W$, hence $N\La_0-N\cont(\rho)=w(N\La_0)$. 
\end{proof}

Recall the unital algebra homomorphism $\zeta_{\theta,\eta}^\La$ from (\ref{EZetaHom}). 

\begin{Corollary} \label{CZeta}
If $\rho\in \Cores_p$, $N\in\Z_{>0}$ and $d\in\Z_{\geq 0}$, then the homomorphism 
$$\zeta_{N\cont(\rho),d\de}^{N\La_0}\colon R_{N\cont(\rho)}^{N\La_0}\to  e_{N\cont(\rho),d\de}R_{N\cont(\rho)+d\de}^{N\La_0}e_{N\cont(\rho),d\de}$$ is injective. 
\end{Corollary}
\begin{proof}
By Corollary~\ref{CMatrix}, $R_{N\cont(\rho)}^{N\La_0}$ is a simple algebra, so it is enough to note by Lemma~\ref{LNonZeroSpecialCut} that $e_{N\cont(\rho),d\de}R_{N\cont(\rho)+d\de}^{N\La_0}e_{N\cont(\rho),d\de}\neq 0$. \end{proof}

By Corollary~\ref{CZeta}, we can identify 
\begin{equation}\label{ECoreSubalg}
R_{N\cont(\rho)}^{N\La_0}=
\zeta_{N\cont(\rho),d\de}^{N\La_0}(R_{N\cont(\rho)}^{N\La_0})\subseteq 
e_{N\cont(\rho),d\de}R_{N\cont(\rho)+d\de}^{N\La_0}e_{N\cont(\rho),d\de}.
\end{equation}
We consider the supercentralizer of $R_{N\cont(\rho)}^{N\La_0}$ in $e_{N\cont(\rho),d\de}R_{N\cont(\rho)+d\de}^{N\La_0}e_{N\cont(\rho),d\de}$:
\begin{equation}\label{ECentralizer}
\Cent_{\rho,d}^{N\La_0} := \Cent_{e_{N\cont(\rho),d\de}R_{N\cont(\rho)+d\de}^{N\La_0}e_{N\cont(\rho),d\de}} (R_{N\cont(\rho)}^{N\La_0}).
\index{z@$\Cent_{\rho,d}^{N\La_0}$}
\end{equation}

By Corollary~\ref{CMatrix} and Lemma~\ref{lem:sup_cent}, we get:

\begin{Lemma} \label{L270116New} 
We have an isomorphism of graded superalgebras  
$$R_{N\cont(\rho)}^{N\La_0} \otimes \Cent_{\rho,d}^{N\La_0} \iso e_{N\cont(\rho),d\de}R_{N\cont(\rho)+d\de}^{N\La_0}e_{N\cont(\rho),d\de}$$ 
given by $a\otimes b \mapsto ab$. 
\end{Lemma}


\begin{Lemma} \label{LIdSymNew} 
Let $e(\bi)\in R_{N\cont(\rho)}^{N\La_0}$ be the idempotent from 
Corollary~\ref{CMatrix}, and set $\eps:=\zeta_{N\cont(\rho),d\de}^{N\La_0}(e(\bi))$. Then $\eps$ is an idempotent in the algebra $R_{N\cont(\rho)+d\de}^{N\La_0}$ such that $\eps e_{N\cont(\rho),d\de}=\eps=e_{N\cont(\rho),d\de}\eps$ and there is an isomorphism of graded algebras $$\eps R_{N\cont(\rho)+d\de}^{N\La_0}\eps\cong \Cent_{\rho,d}^{N\La_0}.$$
In particular,  $\Cent_{\rho,d}^{N\La_0}$ is a symmetric algebra whenever $R_{N\cont(\rho)+d\de}^{N\La_0}$ is so. 
\end{Lemma}
\begin{proof}
By Corollary~\ref{CMatrix}, $\dim_q e(\bi)R_{N\cont(\rho)}^{N\La_0}e(\bi)=1$, so the isomorphism comes from Lemma~\ref{L270116New}. 
By \cite[Theorem IV.4.1]{SY}, an idempotent truncation of a symmetric algebra is symmetric, with a symmetrizing form obtained by restriction. This gives the second statement. \end{proof}

\subsection{Rank $2$ nil-Hecke (super)algebras}
We have pointed out in \S\ref{SSDPI} that $R_{n\al_i}$ is the rank $n$ nil-Hecke algebra if $i\neq 0$, while $R_{n\al_0}$ is the rank $n$ odd nil-Hecke (super)algebra. In this subsection we will consider the special case $n=2$. By Theorem~\ref{TBasis}, the subalgebra $\F[y_1,y_2]$ of $R_{2\al_i}$ generated by $y_1,y_2$ has basis $\{y_1^{a_1}y_2^{a_2}\mid a_1,a_2\in\Z_{\geq 0}\}$. Moreover, for $i\neq 0$, $y_1y_2=y_2y_1$, so the subalgebra $\F[y_1,y_2]$ generated by $y_1,y_2$ is the usual polynomial algebra in $y_1,y_2$. On the other hand, for $i=0$ we have $y_1y_2=-y_2y_1$.

Recall from (\ref{EDivDiag}) the idempotent $e(i^{(2)})=\psi_1y_2\in\R_{2\al_i}$. Define
$$
\la_i:=e(i^{(2)})(y_1+y_2)e(i^{(2)})\quad\text{and}\quad\mu_i:=e(i^{(2)})y_1y_2e(i^{(2)}).
\index{l@$\la_i$}\index{m@$\mu_i$}
$$
Note that $\la_i\mu_i=\mu_i\la_i$ if $i\neq 0$ and $\la_0\mu_0=-\mu_0\la_0$ (see also Lemma~\ref{LNH2}(iv) below). Note also that for all $i$, we have 
$$\la_i=e(i^{(2)})(y_1+y_2)\quad \text{and} \quad \mu_i=e(i^{(2)})y_1y_2.
$$ 
So the subalgebra $\F[\la_i,\mu_i]\subseteq e(i^{(2)})R_{2\al_i}e(i^{(2)})$ generated by $\la_i,\mu_i$ equals the idempotent truncation 
$e(i^{(2)})\F[y_1+y_2,y_1y_2]e(i^{(2)})$, where $\F[y_1+y_2,y_1y_2]$ is the subalgebra of $R_{2\al_i}$ generated by $y_1+y_2$ and $y_1y_2$. It is known from \cite{KL1} for $i\neq 0$ and \cite{EKL} for $i=0$  that $e(i^{(2)})$ is a primitive idempotent in $\R_{2\al_i}$, and in fact:

\begin{Lemma} \label{LNH1}
We have 
$$
e(i^{(2)})R_{2\al_i}e(i^{(2)})=e(i^{(2)})\F[y_1+y_2,y_1y_2]e(i^{(2)})=\F[\la_i,\mu_i].
$$ 
Moreover, $\{\la_i^{a}\mu_i^b\mid a,b\in\Z_{\geq 0}\}$ is a basis of $\F[\la_i,\mu_i]$, $\la_i\mu_i=\mu_i\la_i$ for $i\neq 0$, and $\la_0\mu_0=-\mu_0\la_0$. 
\end{Lemma}

We need to make some additional computations in $e(i^{(2)})R_{2\al_i}e(i^{(2)})$.

\begin{Lemma} \label{LNHEven}
Let $i\neq 0$. Then:
\begin{enumerate}
\item[{\rm (i)}] $e(i^{(2)})y_1e(i^{(2)})=0$ and $e(i^{(2)})y_1^2e(i^{(2)})=-\mu_i$;
\item[{\rm (ii)}] $e(i^{(2)})\psi y_1^2e(i^{(2)})=-\la_i$ and $e(i^{(2)})\psi y_2^2e(i^{(2)})=\la_i$.
\end{enumerate} 
\end{Lemma}
\begin{proof}
(i) Using the fact that $y_2y_1$ commutes with $\psi_1$, we get  
$$e(i^{(2)})y_1e(i^{(2)}=\psi_1y_2y_1e(i^{(2)})=y_2y_1\psi_1e(i^{(2)})=y_2y_1\psi_1^2y_2=0.$$ 
Moreover, 
\begin{align*}
e(i^{(2)})y_1^2e(i^{(2)})&=\psi_1y_2y_1^2e(i^{(2)})=y_1\psi_1y_1^2e(i^{(2)})+y_1^2e(i^{(2)})
\\
&=y_1y_2\psi_1y_1e(i^{(2)})-y_1y_1e(i^{(2)})+y_1^2e(i^{(2)})
=-y_1y_2e(i^{(2)})=-\mu_i,
\end{align*}
where for the last equality we used the fact that $y_1y_2$ commutes with $e(i^{(2)})$. 

Part (ii) is proved similarly. 
\end{proof}

From now on, until the end of the subsection, we consider the case $i=0$. In this case we denote:
\begin{align*}
&\tau:=\psi_1,\ x:=y_1,\ y=y_2,\ e:=e(0^{(2)})=\tau y,\\ &u:=\la_0=e(x+y)e,\ z:=\mu_0=exye.
\index{t@$\tau$}
\index{x@$x$}
\index{y@$y$}\index{e@$e$}
\index{u@$u$}
\index{z@$z$}
\end{align*}
With this notation we now have:

\begin{Lemma} \label{LNH2}
In $R_{2\al_0}$, we have:
\begin{enumerate}
\item[{\rm (i)}] $xe=y-ey=ye-eye$;
\item[{\rm (ii)}] $ex=ye-y=eye-ey$;
\item[{\rm (iii)}] $xye=exy$;
\item[{\rm (iv)}] $uz=-zu$.
\end{enumerate}
\end{Lemma}
\begin{proof}
To prove (i), we use relations in $R_{2\al_0}$ to get
$$
xe=x\tau y=-\tau y y+y=-ey+y,
$$
giving the first equality. The second equality is obtained by multiplying the first equality with $e$ on the right. 

Part (ii) is proved similarly. 

For (iii), note using (ii) and (i) that 
$$xye=x(ex+y)=xex+xy
=(y-ey)x+xy=exy.
$$ 

For (iv), note using (iii) that   
$$uz=e(x+y)exye=e(x+y)xye=-exy(x+y)e=-exye(x+y)e=-zu,$$
completing the proof. 
\end{proof}

\begin{Lemma} \label{LNH3}
In $eR_{2\al_0}e$, we have:
\begin{enumerate}
\item[{\rm (i)}] $eye=u$ and $ey^ne=uey^{n-1}e-zey^{n-2}e$ for all $n\geq 2$; in particular, $ey^2e=u^2-z$, $ey^3e=u^3$ and $ey^4e=u^4-u^2z+z^2$.
\item[{\rm (ii)}] $exe=0$ and $ex^ne=uex^{n-1}e+zex^{n-2}e$ for all $n\geq 2$; in particular, $ex^2e=z$, $ex^3e=uz$ and $ex^4e=u^2z+z^2$.
\end{enumerate}
\end{Lemma}
\begin{proof}
Note using relations that 
$$exe=\tau yx\tau y
=-x\tau x\tau y+x\tau y
=xy\tau \tau y-x\tau y+x\tau y=0,$$
so $u=e(x+y)e=eye$. Now, using Lemma~\ref{LNH2}(ii)(iii), we have
\begin{align*}
ey^ne&=(ey)y^{n-1}e=(eye-ex)y^{n-1}e=(eye)ey^{n-1}e-exyy^{n-2}e
\\
&=
uey^{n-1}e-exyeey^{n-2}e=
uey^{n-1}e-zey^{n-2}e.
\end{align*}
This gives the recurrent relation in (i). The remaining equalities in (i) follow by applying this recurrent relation. The proof of (ii) is similar. 
\end{proof}

\begin{Lemma} \label{LNH4}
In $eR_{2\al_0}e$, we have:
\begin{enumerate}
\item[{\rm (i)}] $e\tau xe=e$, $e\tau x^2e=-u$, $e\tau x^3e=u^2+z$ and $e\tau x^4e=-u^3$;
\item[{\rm (ii)}] $e\tau ye=e$, $e\tau y^2e=u$, $e\tau y^3e=u^2-z$ and $e\tau y^4e=u^3$.
\end{enumerate}
\end{Lemma}
\begin{proof}
(i) We give details for the most difficult last equality. Using  the fact that $eyx^2e=yxexe=0$, Lemma~\ref{LNH3} and Lemma~\ref{LNH2}(iv), we get
\begin{align*}
e\tau x^4e&=ex^3e-ey\tau x^3e=uz-eyx^2e+ey^2\tau x^2e
=uz+ey^2 xe-ey^3\tau xe
\\
&=uz+exy^2 e-ey^3e=uz+zu-ey^3e=-u^3.
\end{align*}

(ii) Use $e\tau y^ne=ey^{n-1}e$ and  Lemma~\ref{LNH3}(i). 
\end{proof}

\subsection{Crystal operators}
The theory of crystal operators for $R_\theta$ has been developed in \cite[\S6]{KKO}. We review necessary facts for reader's convenience. 

Let $i\in I$. Then $R_{n\al_i}$ is the rank $n$ nil-Hecke algebra if $i\neq 0$ and $R_{n\al_0}$ is the rank $n$ odd nil-Hecke algebra. In any case, $R_{n\al_i}$ has a unique (up to isomorphism) irreducible module of dimension $n!$, which we denote by $L(i^n)$\index{l@$L(i^n)$}, see \cite[(4.22)]{KKO}. 
We have functors 
\begin{align*}
&e_i: \mod{R_\theta}\to\mod{R_{\theta-\al_i}},\ M\mapsto \Res^{R_{\theta-\al_i,\al_i}}_{R_{\theta-\al_i}}\circ \Res_{\theta-\al_i,\al_i}M,
\index{e@$e_i$}
\\
&f_i: \mod{R_\theta}\to\mod{R_{\theta+\al_i}},\ M\mapsto \Ind_{\theta,\al_i}M\boxtimes L(i).
\index{f@$f_i$}
\end{align*}
If $L\in\mod{R_\theta}$ is irreducible, we define {\em crystal operators}\index{crystal operators}
$$
\tilde f_i L:=\head (f_i L),\quad \tilde e_i L:=\soc (e_i L).
\index{f@$\tilde f_i$}\index{e@$\tilde e_i$}
$$
A fundamental fact is that for an irreducible $L$ we have $\tilde f_i L$ is again irreducible and $\tilde e_i L$ is irreducible or zero (cf. \cite[\S6, I1]{KKO}). Moreover,  $\tilde e_i\tilde f_iL\cong L$,  and if $\tilde e_i L\neq 0$ then $\tilde f_i\tilde e_iL\cong L$; (cf. 
\cite[\S6, I4]{KKO}). 
For any $M\in\mod{R_\theta}$, we define
$$
\eps_i(M):=\max\{k\geq 0\mid e_i^k(M)\neq 0\}.
\index{e@$\eps_i$}
$$
Then for an irreducible $L\in\mod{R_\theta}$, we have $\eps_i(L)=\max\{k\geq 0\mid \tilde e_i^k(L)\neq 0\}$. Moreover,  setting $\eps:=\eps_i(L)$, we have  
$$\Res_{\theta-\eps_i(L)\al_i,\eps\al_i}L\cong \tilde e_i^\eps L\boxtimes L(i^\eps)$$ (cf. \cite[\S6, I3]{KKO}).

A word $i_1^{a_1}\dots i_b^{a_b}\in I^\theta$, with  $a_1,\dots,a_b\in\Z_{\geq 0}$, is called {\em extremal}\index{extremal word} for $M\in\mod{R_\theta}$ if 
$$a_b=\eps_{i_b}(M),\ a_{b-1}=\eps_{i_{b-1}}(e_{i_b}^{a_b}M)\ ,\ \dots\ ,\  a_1=\eps_{i_1}(e_{i_2}^{a_2}\dots e_{i_b}^{a_b}M).
$$   


We denote by $R_0$ the trivial module $\F$ over the trivial algebra $R_0\simeq \F$. The following useful results are versions of \cite[Theorem 2.16, Corollary 2.17]{BKdurham}, and their proofs are the same, using the standard properties of the crystal operators listed above, cf. \cite[\S2.8]{Kcusp} .

\begin{Lemma} \label{LExtrMult}
Let $L\in \mod{R_\theta}$ be irreducible, and $\bi=i_1^{a_1}\cdots i_b^{a_b}\in\words^\theta$ be an extremal word for $L$. Then we have  
$\dim_q e(\bi)L=q^n[a_1]^!_{i_1}\cdots [a_b]^!_{i_b}$ for some $n\in \Z$, and 
$L\cong \tilde f_{i_b}^{a_b} \tilde f_{i_{b-1}}^{a_{b-1}}\dots\tilde f_{i_1}^{a_1}R_0.$  Moreover, $\bi$ is not an extremal word for any irreducible module $L'\in\mod{R_\theta}$ not isomorphic to  $L$. 
\end{Lemma}

\begin{Corollary} \label{CExtrNew}
Let $M\in\mod{R_\theta}$, and $\bi=i_1^{a_1}\cdots i_b^{a_b}\in\words^\theta$ be an extremal word for $M$. Then we can write $\dim e(\bi)M=ma_1!\cdots a_b!$ for some $m\in\Z_{\geq 0}$. Moreover, if $L:= \tilde f_{i_b}^{a_b} \tilde f_{i_{b-1}}^{a_{b-1}}\dots\tilde f_{i_1}^{a_1}R_0$ then we have $[M:L]=m$. 
\end{Corollary}

The following result shows that any induction product of irreducible modules always has a multiplicity one composition factor, cf. \cite[Corollary 2.12]{Kcusp}.

\begin{Corollary} \label{CPowerIrr}
Let  $\bi=i_1^{a_1}\cdots i_k^{a_k}\in I^\theta$ be an extremal word for an irreducible $R_\theta$-module $L$. Then for any $n\in\Z_{>0}$, we have that $\bj:=i_1^{na_1}\cdots i_k^{na_k}$ is an extremal word for $L^{\circ n}$, and 
$$[L^{\circ n} :\tilde f_{i_k}^{na_k} \tilde f_{i_{k-1}}^{na_{k-1}}\dots\tilde f_{i_1}^{na_1}R_0]=1.$$ 
\end{Corollary}
\begin{proof}
It is easy to see using Corollary~\ref{CShuffle} that $\bj$  
is an extremal word for $L^{\circ n}$. By Lemma~\ref{LExtrMult}   
$\dim e(\bi)L=a_1!\cdots a_k!$ . From this we deduce 
that $\dim e(\bj)(L^{\circ n})=(na_1)!\cdots (na_k)!$ and apply Corollary~\ref{CExtrNew}. We refer the reader to the proof of  \cite[Proposition 2.11, Corollary 2.12]{Kcusp} for more details. 
\end{proof}

\subsection{Kang-Kashiwara-Oh intertwiners}\label{SSKKO}
Let $\theta\in Q_+$ and $n:=\height(\theta)$. Following \cite[\S8]{KKO}, we define the {\em KKO intertwiners}\index{KKO intertwiners} $\phi_1,\dots,\phi_{n-1}\in R_\theta$\index{f@$\phi_r$} by setting
\begin{equation}\label{EPhi}
\phi_re(\bi):=
\left\{
\begin{array}{ll}
\psi_re(\bi) &\hbox{if $i_r\neq i_{r+1}$,}\\
(1+(y_r-y_{r+1})\psi_r)e(\bi) &\hbox{if $i_r=i_{r+1}\neq 0$,}
\\
(y_{r+1}-y_r+(y_r^2-y_{r+1}^2)\psi_r)e(\bi) &\hbox{if $i_r=i_{r+1}= 0$.}
\end{array}
\right.
\end{equation}

We extend Khovanov-Lauda diagrams to denote 
\begin{equation}\label{EPhiDiagram}
\phi_r e(\bi)=
 \begin{braid}\tikzset{baseline=3mm}
  \draw (0,0)node[below]{$i_1$}--(0,3);
  \draw (3,0)node[below]{$i_{r-1}$}--(3,3);
  \draw[dots] (0.5,2.9)--(2.7,2.9);
  \draw (4.6,0)node[below]{$i_{r}$}--(5.8,3);
  \draw (5.8,0)node[below]{$i_{r+1}$}--(4.6,3);
  \draw (7.4,0)node[below]{$i_{r+2}$}--(7.4,3);
  \draw[dots] (0.5,0)--(2.7,0);
  \draw[dots] (7.8,2.9)--(10,2.9);
  \draw[dots] (7.8,0)--(10,0);
  \draw (10.2,0)node[below]{$i_{n}$}--(10.2,3);
 \filldraw[color=black!60, fill=black!5, thick] (5.2,1.5) circle (0.2);
\end{braid}.
\end{equation}
Thus we have
\begin{eqnarray}
\label{EPhiDifCol}
\begin{braid}\tikzset{baseline=2mm}
  \draw (0,0)node[below]{$i$}--(1,2);
  \draw (1,0)node[below]{$j$}--(0,2);
 \filldraw[color=black!60, fill=black!5, thick] (0.5,1) circle (0.2);
\end{braid}&=&\begin{braid}\tikzset{baseline=3mm}
  \draw (0,0)node[below]{$i$}--(1,2);
  \draw (1,0)node[below]{$j$}--(0,2);
\end{braid}
\qquad(\text{if $i\neq j$}),
\\
\label{EOpening}
\begin{braid}\tikzset{baseline=2mm}
  \draw (0,0)node[below]{$i$}--(1,2);
  \draw (1,0)node[below]{$i$}--(0,2);
 \filldraw[color=black!60, fill=black!5, thick] (0.5,1) circle (0.2);
\end{braid}&=&
\begin{braid}\tikzset{baseline=2mm}
  \draw (0,0)node[below]{$i$}--(0,2);
  \draw (1,0)node[below]{$i$}--(1,2);
\end{braid}
+\begin{braid}\tikzset{baseline=2mm}
  \draw (0,0)node[below]{$i$}--(1,2);
  \draw (1,0)node[below]{$i$}--(0,2);
  \blackdot (0.1,1.8);
\end{braid}
-\begin{braid}\tikzset{baseline=2mm}
  \draw (0,0)node[below]{$i$}--(1,2);
  \draw (1,0)node[below]{$i$}--(0,2);
  \blackdot (0.9,1.8);
\end{braid}\qquad(\text{if $i\neq 0$}),
\\
\label{EOpening0}
\begin{braid}\tikzset{baseline=2mm}
  \draw[red] (0,0)node[below]{$0$}--(1,2);
  \draw[red] (1,0)node[below]{$0$}--(0,2);
 \filldraw[color=red!60, fill=red!5, thick] (0.5,1) circle (0.2);
\end{braid}&=&
\begin{braid}\tikzset{baseline=2mm}
  \draw[red] (0,0)node[below]{$0$}--(1,2);
  \draw[red] (1,0)node[below]{$0$}--(0,2);
  \reddot (0.1,1.8);
  \reddot (0.3,1.4);
\end{braid}
-\begin{braid}\tikzset{baseline=2mm}
  \draw[red] (0,0)node[below]{$0$}--(1,2);
  \draw[red] (1,0)node[below]{$0$}--(0,2);
  \reddot (0.9,1.8);
  \reddot (0.7,1.4);
\end{braid}
-\begin{braid}\tikzset{baseline=2mm}
  \draw[red] (0,0)node[below]{$0$}--(0,2);
  \draw[red] (1,0)node[below]{$0$}--(1,2);
  \reddot (0,1);
\end{braid}
+\begin{braid}\tikzset{baseline=2mm}
  \draw[red] (0,0)node[below]{$0$}--(0,2);
  \draw[red] (1,0)node[below]{$0$}--(1,2);
  \reddot (1,1);
\end{braid}.
\end{eqnarray}

The following key properties of the intertwiners are  established in \cite{KKO}: 

\begin{Lemma} \label{LKKO} {\rm \cite[Lemma 8.3]{KKO}}
For $1\leq r,s<n$, $1\leq t\leq n$ and $\bi\in I^\theta$, we have:
\begin{align}
\label{KKO1}
\phi_re(\bi)=e(s_r\cdot\bi)\phi_r,\\
\label{KKO2}
\phi_ry_te(\bi)=(-1)^{|i_r||i_{r+1}||i_t|}y_{s_r(t)}\phi_re(\bi),
\\
\label{KKO3}
\phi_r\psi_se(\bi)=(-1)^{|i_r||i_{r+1}||i_s||i_{s+1}|}\psi_{s}\phi_re(\bi)\quad\text{if $|r-s|>1$},
\\
\label{KKO4}
\psi_r\phi_{r+1}\phi_{r}=\phi_{r+1}\phi_{r}\psi_{r+1}\quad\text{if $r<n-1$}.
\end{align}
\end{Lemma}

In particular, in terms of diagrams, (\ref{KKO1}) implies $
\begin{braid}\tikzset{baseline=2mm}
  \draw (0,0)node[below]{$i$}--(1,2);
  \draw (1,0)node[below]{$j$}--(0,2);
 \filldraw[color=black!60, fill=black!5, thick] (0.5,1) circle (0.2);
\end{braid}
=\begin{braid}\tikzset{baseline=2mm}
  \draw (0,2)node[above]{$j$}--(1,0);
  \draw (1,2)node[above]{$i$}--(0,0);
   \filldraw[color=black!60, fill=black!5, thick] (0.5,1) circle (0.2);
\end{braid}
$, 
(\ref{KKO2}) implies
\begin{equation}\label{EDotPastCircle}
\begin{braid}\tikzset{baseline=2mm}
  \draw (0,0)node[below]{$i$}--(1,2);
  \draw (1,0)node[below]{$j$}--(0,2);
 \filldraw[color=black!60, fill=black!5, thick] (0.5,1) circle (0.2);
 \blackdot (0.9,0.2);
\end{braid}
=(-1)^{|i||j|}\begin{braid}\tikzset{baseline=2mm}
   \draw (0,0)node[below]{$i$}--(1,2);
  \draw (1,0)node[below]{$j$}--(0,2);
   \filldraw[color=black!60, fill=black!5, thick] (0.5,1) circle (0.2);
   \blackdot (0.1,1.8);
\end{braid}
\qquad\text{and}\qquad
\begin{braid}\tikzset{baseline=2mm}
  \draw (0,0)node[below]{$i$}--(1,2);
  \draw (1,0)node[below]{$j$}--(0,2);
 \filldraw[color=black!60, fill=black!5, thick] (0.5,1) circle (0.2);
 \blackdot (0.1,0.2);
\end{braid}
=(-1)^{|i||j|}\begin{braid}\tikzset{baseline=2mm}
   \draw (0,0)node[below]{$i$}--(1,2);
  \draw (1,0)node[below]{$j$}--(0,2);
   \filldraw[color=black!60, fill=black!5, thick] (0.5,1) circle (0.2);
   \blackdot (0.9,1.8);
\end{braid},
\end{equation}
and 
(\ref{KKO4}) implies
\begin{equation}\label{EPartialCircleBraid}
\begin{braid}\tikzset{baseline=2mm}
  \draw (0,0)--(2,2);
  \draw (2,0)--(0,2);
  \draw (1,0)--(0,1)--(1,2);
 \filldraw[color=black!60, fill=black!5, thick] (1,1) circle (0.2);
 \filldraw[color=black!60, fill=black!5, thick] (.5,.5) circle (0.2);
\end{braid}
=
\begin{braid}\tikzset{baseline=2mm}
  \draw (0,0)--(2,2);
  \draw (2,0)--(0,2);
  \draw (1,0)--(2,1)--(1,2);
 \filldraw[color=black!60, fill=black!5, thick] (1,1) circle (0.2);
 \filldraw[color=black!60, fill=black!5, thick] (1.5,1.5) circle (0.2);
\end{braid}.
\end{equation}

\begin{Corollary} \label{C180621} 
Let $\bi\in I^\theta$ and $1\leq r,s,t< n$. 
\begin{enumerate}
\item[{\rm (i)}] If $|s-t|>1$ then 
$\phi_s\phi_te(\bi)=(-1)^{|i_s||i_{s+1}||i_t||i_{t+1}|}\phi_t\phi_se(\bi)$.
\item[{\rm (ii)}] If $r<n-1$ then $\phi_r\phi_{r+1}\phi_r=\phi_{r+1}\phi_r\phi_{r+1}$.
\end{enumerate} 
\end{Corollary}
\begin{proof} (i) comes from (\ref{R3}),(\ref{R4}),(\ref{R65}). For  (ii), it suffices to show that 
$\phi_r\phi_{r+1}\phi_re(\bi)=\phi_{r+1}\phi_r\phi_{r+1}e(\bi)$ for any $\bi\in I^\theta$. If $i_{r+1}\neq i_{r+2}$ then we have $\phi_{r+1}e(\bi)=\psi_{r+1}e(\bi)$, and $\phi_{r}e(s_rs_{r+1}\bi)=\psi_{r}e(s_rs_{r+1}\bi)$, so 
\begin{align*}
\phi_{r+1}\phi_r\phi_{r+1}e(\bi)
&=
\phi_{r+1}\phi_r\psi_{r+1}e(\bi)
\\&
\stackrel{(\ref{KKO4})}{=}\psi_{r}\phi_{r+1}\phi_{r}e(\bi)
\\
&
\stackrel{(\ref{KKO1})}{=}
\psi_{r}e(s_{r+1}s_r\cdot\bi)\phi_{r+1}\phi_{r}e(\bi)
\\
&
=\phi_r\phi_{r+1}\phi_re(\bi).
\end{align*}
Suppose $i_{r+1}= i_{r+2}=0$. Then 
\begin{align*}
\phi_{r+1}\phi_r\phi_{r+1}e(\bi)
&=
\phi_{r+1}\phi_r\big(y_{r+2}-y_{r+1}+(y_{r+1}^2-y_{r+2}^2)\psi_{r+1}\big)e(\bi)
\\
&=
\phi_{r+1}\phi_r(y_{r+2}-y_{r+1})e(\bi)
+\phi_{r+1}\phi_r(y_{r+1}^2-y_{r+2}^2)\psi_{r+1}e(\bi)
\\
&
\stackrel{(\ref{KKO2})}{=}
(y_{r+1}-y_r)\phi_{r+1}\phi_re(\bi)+(y_r^2-y_{r+1}^2)\phi_{r+1}\phi_r\psi_{r+1}e(\bi)
\\
&\stackrel{(\ref{KKO4})}{=}(y_{r+1}-y_r)\phi_{r+1}\phi_re(\bi)+(y_r^2-y_{r+1}^2)\psi_{r}\phi_{r+1}\phi_{r}e(\bi)
\\
&\stackrel{(\ref{KKO1})}{=}\big(y_{r+1}-y_r+(y_r^2-y_{r+1}^2)\psi_{r}\big)e(s_{r+1}s_r\cdot\bi)\phi_{r+1}\phi_{r}e(\bi)\\
&=\phi_r\phi_{r+1}\phi_re(\bi).
\end{align*}
The case $i_{r+1}= i_{r+2}=0$ is similar.
\end{proof}

In particular, in terms of diagrams, Corollary~\ref{C180621}(ii) implies
\begin{equation}\label{ECircleBraid}
\begin{braid}\tikzset{baseline=2mm}
  \draw (0,0)--(2,2);
  \draw (2,0)--(0,2);
  \draw (1,0)--(0,1)--(1,2);
 \filldraw[color=black!60, fill=black!5, thick] (1,1) circle (0.2);
 \filldraw[color=black!60, fill=black!5, thick] (.5,.5) circle (0.2);
 \filldraw[color=black!60, fill=black!5, thick] (.5,1.5) circle (0.2);
\end{braid}
=
\begin{braid}\tikzset{baseline=2mm}
  \draw (0,0)--(2,2);
  \draw (2,0)--(0,2);
  \draw (1,0)--(2,1)--(1,2);
 \filldraw[color=black!60, fill=black!5, thick] (1,1) circle (0.2);
 \filldraw[color=black!60, fill=black!5, thick] (1.5,1.5) circle (0.2);
 \filldraw[color=black!60, fill=black!5, thick] (1.5,.5) circle (0.2);
\end{braid}.
\end{equation}

Let $w\in \Si_n$ with a reduced decomposition $w=s_{r_1}\cdots s_{r_k}$. Define 
\begin{equation}\label{EPhiW}
\phi_w:=\phi_{r_1}\cdots\phi_{r_k}\in R_\theta.
\end{equation}
In view of the previous corollary, choosing a different reduced decomposition for $w$ might only lead to a sign change for $\phi_w$. 

Recalling the polynomials $Q_{i,j}(u,v)$ from \S\ref{SQHDef}, we have:

\begin{Lemma} \label{LphiQuadr} 
For $1\leq r<n$ and $\bi\in I^\theta$, we have
$$
\phi_r^2e(\bi)=
\left\{
\begin{array}{ll}
Q_{i_r,i_{r+1}}(y_r,y_{r+1})e(\bi) &\hbox{if $i_r\neq i_{r+1}$,}\\
e(\bi) &\hbox{if $i_r= i_{r+1}\neq 0$,}\\
(y_r^2+y_{r+1}^2)e(\bi) &\hbox{if $i_r= i_{r+1}=0$.}
\end{array}
\right.
$$
\end{Lemma}
\begin{proof}
If we are not in the case $i_r= i_{r+1}=0$, the result follows from \cite[Lemma 1.5(i)]{KKK}. Let $i_r= i_{r+1}=0$. By \cite[(8.10)]{KKO}, we have 
$$\phi_re(\bi)=(y_r-y_{r+1}-\psi_r(y_r^2-y_{r+1}^2))e(\bi),$$ 
so $\phi_r^2e(\bi)$ equals
\begin{align*}
\phi_r^2e(\bi)\,=\ &\big(y_{r+1}-y_r+(y_r^2-y_{r+1}^2)\psi_{r}\big)
\big(y_r-y_{r+1}-\psi_r(y_r^2-y_{r+1}^2)\big)e(\bi)
\\
=\ 
&\Big((y_{r+1}-y_r)(y_r-y_{r+1})-(y_{r+1}-y_r)\psi_r(y_r^2-y_{r+1}^2)
\\&
+(y_r^2-y_{r+1}^2)\psi_r(y_r-y_{r+1})\Big)e(\bi)
\\
\stackrel{(\ref{R5})}{=}
&\Big((y_{r+1}-y_r)(y_r-y_{r+1})-(y_{r+1}-y_r)(y_{r+1}^2-y_{r}^2)\psi_r
\\
&
-2(y_{r+1}-y_r)(y_r-y_{r+1})
-(y_r^2-y_{r+1}^2)(y_{r+1}-y_{r})\psi_r\Big)e(\bi)
\\
=\ 
&-(y_{r+1}-y_r)(y_r-y_{r+1})e(\bi)
\stackrel{(\ref{R3})}{=}
(y_r^2+y_{r+1}^2)e(\bi),
\end{align*}
as required. 
\end{proof}

\section{Cuspidality}
In this section we develop the theory of cuspidal systems and (proper) standard modules for quiver Hecke superalgebras by analogy with a similar theory for KLR algebras \cite{KR, McN1, Kcusp, TW, McN2}.  
Recall from \S\ref{SSLT} the set $\Phi_+$ of positive roots.

\subsection{Convex preorders}
\label{SSConv}
 Following \cite{BKT}, we define a {\em convex preorder}\index{convex preorder} on $\Phi_+$ as a preorder $\preceq$\index{$\preceq$} such that the following three conditions hold for all $\be,\ga\in\Phi_+$:
\begin{eqnarray}
\label{EPO1}
&\be\preceq\ga \ \text{or}\ \ga\preceq \be;
\\
\label{EPO2}
&\text{if $\be\preceq \ga$ and $\be+\ga\in\Phi_+$, then $\be\preceq\be+\ga\preceq\ga$};
\\
&\label{EPO3}
\text{$\be\preceq\ga$ and $\ga\preceq\be$ if and only if $\be$ and $\ga$ are proportional}.
\end{eqnarray}
We write $\be\prec\ga$ if $\be\preceq\ga$ but $\ga\not\preceq\be$.

By (\ref{EPO3}), $\be\preceq\ga$ and $\ga\preceq\be$ happens for $\be\neq\ga$ if and only if both $\be$ and $\ga$ are imaginary. 
In particular, the set $\Psi$ from (\ref{EPsi}) is totally ordered with respect to $\preceq$. 
The set of real roots splits into two disjoint infinite sets
$$
\Phi^\re_{\succ\de}:=\{\be\in \Phi_+^\re\mid \be\succ\de\}\ \text{and}\ 
\Phi^\re_{\prec\de}:=\{\be\in \Phi_+^\re\mid \be\prec\de\}. 
\index{f@$\Phi^\re_{\succ\de}$}\index{f@$\Phi^\re_{\prec\de}$}
$$

\begin{Example} \label{ExConPr} 
{\rm 
Following \cite[Example 3.5]{McN2}, consider $\R$ as a vector space over $\Q$ and let $\chi:\Q\Phi\to \R$ be an injective linear map. Then it is easy to see that the following defines a convex preorder on $\Phi_+$: 
$$\be\preceq\ga \quad \iff\quad \chi(\be)/\height(\be)\leq \chi(\ga)/\height(\ga).
$$ 
Recall the description of the real roots from (\ref{EPhiS})-(\ref{EPhiL}), and the finite root system $\Phi'$ from \S\ref{SSLT}. Let  
let $\tilde\al=2\al_1+\dots+2\al_{l-1}+\al_\ell\in\Phi'$ be the highest root in $\Phi'$. Denote by $\Phi'_{\sharp}$ the set of all roots in 
$\Phi'$ which are non-negative linear combinations of $\al_1,\dots,\al_{\ell-1},-\tilde\al$ (this is a non-standard choice of a system of positive roots in $\Phi'$). The linear function $\chi:\Q\Phi\to \R$ is determined by a choice of $\Q$-linearly independent real numbers $\chi(\al_0),\dots,\chi(\al_{\ell-1}),\chi(\al_\ell)$. These can be chosen so that $\chi(\al_1),\dots,\chi(\al_{\ell-1}),\chi(-\tilde\al)$ are positive and $\chi(\de)$ is very close to zero. 
For the corresponding convex order $\preceq$ we have 
\begin{equation}\label{ESharp}
\begin{split}
\Phi^\re_{\succ\de}=\{\be\in\Phi_+\mid\text{$\be$ is of the form $r\al+s\de$ with $r\in\{1/2,1\}$, $\al\in\Phi'_{\sharp}$}\},
\\
\Phi^\re_{\prec\de}=\{\be\in\Phi_+\mid\text{$\be$ is of the form $-r\al+s\de$ with $r\in\{1/2,1\}$, $\al\in\Phi'_{\sharp}$}\}.
\end{split}
\end{equation}
We will see that such convex orders are related to RoCK blocks of spin covers of symmetric groups.
}
\end{Example}

\begin{Lemma}\label{LCones} \cite[Theorem 3.2]{McN2}
Suppose $A$ and $B$ are disjoint subsets of $\Phi_+$ such that $\al\prec\be$ for any $\al\in A$ and $\be\in B$. Then the cones formed by the $\R_{\geq 0}$ spans of $A$ and $B$ meet only at the origin.
\end{Lemma}

The following properties of convex preorders are well-known and can be easily deduced from Lemma~\ref{LCones}, see for example \cite[\S3.1]{Kcusp}:

\begin{enumerate}
\item[{\rm (Con1)}] Let $\be\in\Phi_+^\re$, $m\in\Z_{>0}$, and $m\be=\sum_{a=1}^b \ga_a$ for some positive roots $\ga_a$. Assume that either  $\ga_a\preceq \be$ for all $a=1,\dots,b$ or $\ga_a\succeq \be$ for all $a=1,\dots,b$. Then  $b=m$ and $\ga_a=\be$ for all $a=1,\dots,b$. 
\item[{\rm (Con2)}] Let $\al,\be$ be two positive roots, not both imaginary. If $\al+\be=\sum_{a=1}^b \ga_a$ for some positive roots $\ga_a\preceq \be$, then $\be\succeq \al$. 
\item[{\rm (Con3)}] Let $\be\in\Phi_+^\im$, and $\be=\sum_{a=1}^b \ga_a$ for some positive roots $\ga_a$. If either $\ga_a\preceq \be$ for all $a=1,\dots,b$ or $\ga_a\succeq \be$ for all $a=1,\dots,b$, then all $\ga_a$ are imaginary. 
\end{enumerate}

\subsection{Cuspidal modules}
\label{SSCuspidalModules}
Fix a convex preorder $\preceq$ on $\Phi_+$. 
Let $\be\in\Psi$, $m\in\Z_{\geq 1}$ and $M\in\mod{R_{m\be}}$. Then $M$ is called {\em cuspidal}\index{cuspidal module} if $\Res_{\theta,\eta}M\neq 0$ for $\theta,\eta\in Q_+$ implies that $\theta$ is a sum of positive roots $\preceq \be$ and $\eta$ is a sum of positive roots $\succeq\be$. Recalling the notation $\|\bi\|$ from \S\ref{ChBasicNotLie}, a word $\bi\in I^{m\be}$ is called {\em cuspidal}\index{cuspidal word} if $\bi=\bj\bk$ for words $\bj,\bk$ implies that $\|\bj\|$ is a sum of positive roots $\preceq \be$ and $\|\bk\|$ is a sum of positive roots $\succeq\be$. (Note that we eschew the terminology `semicuspidal' used in some literature). 

We denote by 
$I^{m\be}_\cus$\index{i@$I^{m\be}_\cus$} the set of all cuspidal words in $I^{m\be}$. 
Then $M\in\mod{R_{m\be}}$ is cuspidal if and only if every word of $M$ is cuspidal if and only if $e(\bi)M=0$ for all $\bi\in I^{m\be}\setminus I^{m\be}_\cus$ if and only if $M$ 
 factors through a module over the {\em cuspidal algebra} \index{cuspidal algebra}
\begin{equation}\label{ECuspidalAlgebra}
\bar R_{m\be}:=R_{m\be}/(e(\bi)\mid \bi\in I^{m\be}\setminus I^{m\be}_\cus),
\index{r@$\bar R_{m\be}$}
\end{equation}
which is the quotient of $R_{m\be}$ by the two-sided ideal generated by all  idempotents $e(\bi)$ corresponding to non-cuspidal words. Thus $\mod{\bar R_{m\be}}$ can be identified with the full subcategory of $\mod{R_{m\be}}$ which consists of  the cuspidal (graded super)modules. 

\begin{Lemma} \label{LTensImagIsImag}
Let $\be\in\Psi$ and $m,n\in\Z_{>0}$. 
\begin{enumerate}
\item[{\rm (i)}] If $\bi\in I^{m\be}_\cus$ and $\bj\in I^{n\be}_\cus$ then any shuffle of $\bi$ and $\bj$ is in $ I^{(m+n)\be}_\cus$. 
\item[{\rm (ii)}] If $M\in\mod{R_{m\be}}$ and $N\in\mod{R_{n\be}}$ are cuspidal then so is $M\circ N$.
\end{enumerate}
\end{Lemma}
\begin{proof}
(i) Suppose that $\bk\bl$. Then we can write $\bi=\bi'\bi''$ and  $\bj=\bj'\bj''$ so that $\bk$ is a shuffle of $\bi'$ and $\bj'$, and $\bl$ is a shuffle of $\bi''$ and $\bj''$. In particular,  
$\|\bk\|=\|\bi'\|+\|\bj'\|$ and $\|\bl\|=\|\bi''\|+\|\bj''\|$. 
Since $\bi$ and $\bj$ are cuspidal,  $\|\bi'\|$ and $\|\bj'\|$ are sums of roots $\preceq \be$, and $\|\bi''\|$ and $\|\bi''\|$ are sums of roots $\succeq \be$, proving (i). 

(ii) follows from (i) by Corollary~\ref{CShuffle}. 
\end{proof}

A {\em cuspidal system}\index{cuspidal system} (for the fixed convex preorder $\preceq$) is the following data:
\begin{enumerate}
\item[{\rm (Cus1)}] An irreducible cuspidal $L_\be\in\mod{R_\be}$ for every $\be\in \Phi_+^\re$;
\item[{\rm (Cus2)}] An irreducible cuspidal $L_\bmu\in\mod{R_{n\de}}$ for every $n\in\Z_{>0}$ and every  $\bmu\in\Par^\ell(n)$, such that $L_\bmu\not\cong L_\bnu$ if $\bmu\neq\bnu$. 
\end{enumerate}

We call the (graded super)modules $L_\be$ from (Cus1) {\em real irreducible cuspidal modules}, and 
the modules $L(\bmu)$ from (Cus2) {\em  imaginary irreducible cuspidal modules}. 
It will be proved that for every convex preorder, there exists a cuspidal systems, unique up to isomorphism of modules and  permutation of the imaginary cuspidal modules $L_\bmu$ for $\bmu\in \Par^\ell(n)$ for every $n$.

\subsection{Root partitions and proper standard modules}
We continue working with a fixed convex preorder $\preceq$ on $\Phi_+$. Recall the notation $\Par(\theta)$ from \S\ref{SSLT}. 
Let $(\um,\bmu)\in \Par(\theta)$ for some $\theta\in Q+$. Since almost all $m_\be$ are zero, we can choose a finite subset
$
\be_1\succ\dots\succ\be_s\succ\de\succ\be_{-t}\succ\dots\succ\be_{-1}
$
of $\Psi$ such that $m_\be=0$ for $\be$'s outside of this subset. Then, denoting $m_u:=m_{\be_u}$, we can write $(\um,\bmu)$ in the form
\begin{equation}\label{EStandForm}
(\um,\bmu)=(\be_1^{m_1},\dots,\be_s^{m_s},\bmu,\be_{-t}^{m_{-t}},\dots,\be_{-1}^{m_{-1}}).
\end{equation}
Denote also 
$$
\lan\um\ran:=(m_1\be_1,\dots,m_s\be_s,m_\de\de,m_{-t}\be_{-t},\dots,m_{-1}\be_{-1})\in Q_+^{s+t+1},
$$
so we have a parabolic subalgebra $R_{\lan\um\ran}\subseteq e_{\lan\um\ran}R_\theta e_{\lan\um\ran}$, cf. \S\ref{SSIndRes}.

We define a partial order on $\Par(\theta)$. 
First, there are left and write lexicographic orders $\leq_l$ and $\leq_r$ on the finitary tuples $\um=(m_{\be})_{\be\in\Psi}$, and for two such tuples $\um$ and $\un$ we set 
$\um\leq\un$ if $\um\leq_l \un$ and $\um\geq_r\un$. 
As in \cite{McN1,McN2}, for $(\um,\bmu),(\un,\bnu)\in\Par(\theta)$, we now set \begin{equation}\label{EBilex}
(\um,\bmu)\leq (\un,\bnu)\ \iff  \um< \un,\ \text{or $\um=\un$ and}\ \bmu=\bnu.
\end{equation}

Now suppose that $\{L_\be,L_\bmu\mid\be\in\Phi_+^\re,\bmu\in\Par^\ell\}$ is a cuspidal system for our fixed convex preorder. For a root partition 
$(\um,\bmu)\in\Pi(\theta)$ 
written in the form (\ref{EStandForm}), we have 
$$
L_{\um,\bmu}:=L_{\be_1}^{\circ m_1} \boxtimes \dots\boxtimes L_{\be_s}^{\circ m_s}\boxtimes L(\bmu)\boxtimes L_{\be_{-t}}^{\circ m_{-t}}\boxtimes\dots\boxtimes  L_{\be_{-1}}^{\circ m_{-1}}\in\mod{R_{\lan\um\ran}}
\index{l@$L_{\um,\bmu}$}
$$
and the {\em standard (graded super)module}\index{standard module} over $R_\theta$:
\begin{equation}\label{EStand}
\Stand(\um,\bmu):=\Ind_{\lan\um\ran}L_{\um,\bmu}=L_{\be_1}^{\circ m_1} \circ \dots\circ L_{\be_s}^{\circ m_s}\circ L(\bmu)\circ L_{\be_{-t}}^{\circ m_{-t}}\circ\dots\circ  L_{\be_{-1}}^{\circ m_{-1}}.
\index{d@$\Stand(\um,\bmu)$}
\end{equation}
(In light of the theory of stratified algebras, the `correct'  terminology and notation are `proper standard' and $\bar\De(\um,\bmu)$, but since no other kinds of standard modules will appear the terminology and notation we chose will cause no confusion). 

The proof of the following proposition follows that of \cite[Proposition 3.5]{Kcusp} closely, but we provide it for reader's convenience.

\begin{Proposition}\label{P1}
Let $(\um,\bmu),(\un,\bnu)\in\Pi(\theta)$. Then:
\begin{enumerate}
\item[{\rm (i)}] $\Res_{\lan\un\ran} \Stand(\un,\bnu)\cong L_{\un,\bnu}$.
\item[{\rm (ii)}] $\Res_{\lan\um\ran} \Stand(\un,\bnu)\neq 0$ implies $\um\leq \un$. 
\end{enumerate}
\end{Proposition}
\begin{proof}
We write the root partitions $(\um,\bmu)$ and $(\un,\bnu)$ in the form (\ref{EStandForm}):
\begin{align*}
(\um,\bmu)&=(\be_1^{m_1},\dots,\be_s^{m_s},\bmu,\be_{-t}^{m_{-t}},\dots,\be_{-1}^{m_{-1}}),
\\
(\un,\bnu)&=(\be_1^{n_1},\dots,\be_s^{n_s},\bnu,\be_{-t}^{n_{-t}},\dots,\be_{-1}^{n_{-1}}).
\end{align*}
Let $\Res_{\lan\um\ran} \Stand(\un,\bnu)\neq 0$. It suffices to prove that $\um\geq_l \un$ or $\um\leq_r \un$ implies that $\um=\un$ and  $\Res_{\lan\um\ran} \Stand(\un,\bnu)\cong L_{\un,\bnu}$. We assume that $\um\geq_l \un$, the case $\um\leq _r \un$ being similar. 
We apply induction on $\height(\theta)$ and consider three cases. 

Case 1: $m_\be>0$ for some $\be>\de$. 
Pick the maximal such $\be$. Since $\un\leq_l \um$, we have that $n_\ga=0$ for all $\ga>\be$ and $n_\be\leq m_\be$. 
Let $(\um',\bmu')\in\Par(\theta-m_\be\be)$ and $(\un',\bnu')\in\Par(\theta-n_\be\be)$ be the root partitions obtained by eliminating $\be^{m_\be}$ and $\be^{n_\be}$ from $(\um,\bmu)$ and $(\un,\bnu)$, respectively. 
By the Mackey Theorem~\ref{TMackeyKL}(ii), $\Res_{\lan\um\ran} \Stand(\un,\bnu)$ has  filtration with factors of the form 
$
\Ind^{m_\be\be;|\um'|}_{\eta_1,\dots,\eta_c;\underline{\ga}}V,
$
where $\eta_1,\dots,\eta_c\in Q_+\setminus\{0\}$ satisfy $\eta_1+\dots+\eta_c=m_\be\be$, and $\underline{\ga}$ is a refinement of $\lan\um'\ran$. Moreover, the module $V$ is obtained by twisting as in (\ref{ETwist}) of a module obtained by restriction of 
$$
L_{\be_1}^{\boxtimes n_1}\boxtimes\dots\boxtimes  L_{\be_s}^{\boxtimes n_s}\boxtimes L(\bnu)\boxtimes L_{\be_{-t}}^{\boxtimes n_{-t}}\boxtimes\dots\boxtimes L_{\be_{-1}}^{\boxtimes n_{-1}}
$$
to a parabolic which has $\eta_1,\dots,\eta_c$ in the beginnings of the corresponding blocks. In particular, if $V\neq 0$, then for each $b=1,\dots,c$ we have that $\Res_{\eta_b,\be_k-\eta_b}L_{\be_k}\neq 0$ for some $k=k(b)$ with $n_k\neq 0$ or $\Res_{\eta_b,n_\de\de-\eta_b}L(\bnu)\neq 0$. 

If $\Res_{\eta_b,\be_k-\eta_b}L_{\be_k}\neq 0$, then by (Cus1), $\eta_b$ is a sum of roots $\preceq \be_k$. Moreover, since $\um\geq_l \un$ and $n_k\neq0$, we have that $\be_k\preceq\be$. So $\eta_b$ is a sum of roots $\preceq \be$. On the other hand, if $\Res_{\eta_b,n_\de\de-\eta_b}L(\bnu)\neq 0$, then by (Cus2), either $\eta_b$ is an imaginary root or it is a sum of real roots less than $\de$. Thus in either case, $\eta_b$ is a sum of roots $\preceq \be$. Using (Con1), we conclude that $c=m_\be$, and $\eta_b=\be=\be_{k(b)}$ for all $b=1,\dots,c$. Hence $n_\be\geq m_\be$. So $n_\be=m_\be$, and 
$$
\Res_{\lan\um\ran} \Stand(\un,\bnu)\cong L_{\be}^{\circ m_\be}\boxtimes \Res_{\lan\um'\ran}^{\theta-m_\be\be} \Stand(\un',\bnu'). 
$$
Since $\height(\theta-m_\be\be)<\height(\theta)$, we can now apply the inductive hypothesis. 

Case 2: $m_\be=0$ for all $\be>\de$, but $m_\de\neq 0$. 
Since $\un\leq_l \um$, we also have that $n_\be=0$ for all $\be>\de$. Let $(\um',\bmu')\in\Par(\theta-m_\de\de)$ and $(\un',\bnu')\in\Par(\theta-n_\be\de)$ be the root partitions obtained by eliminating $\bmu$ and $\bnu$ from $(\um,\bmu)$ and $(\un,\bnu)$, respectively. By the Mackey Theorem~\ref{TMackeyKL}(ii), $\Res_{\lan\um\ran} \Stand(\un,\bnu)$ has filtration with factors of the form 
$
\Ind^{m_\de\de;\lan\um'\ran}_{\eta_1,\dots,\eta_c;\underline{\ga}}V,
$
where $m_\de\de=\eta_1+\dots+\eta_c$, with $\eta_1,\dots,\eta_c\in Q_+\setminus\{0\}$, and $\underline{\ga}$ is a refinement of $\lan\um'\ran$. Moreover, the module $V$ is obtained by twisting of a module obtained by parabolic restriction of the module 
$
L(\bnu)\boxtimes L_{\be_{-t}}^{\boxtimes n_{-t}}\boxtimes\dots\boxtimes L_{\be_{-1}}^{\boxtimes n_{-1}}
$
to a parabolic which has $\eta_1,\dots,\eta_c$ in the beginnings of the corresponding blocks. In particular, if $V\neq 0$, then either 

(1) $\Res_{\eta_1,n_\de\de-\eta_1}L(\bnu)\neq 0$ and for  $b=2,\dots,c$, there is $k=k(b)<0$ such that $\Res_{\eta_b,\be_k-\eta_b}L_{\be_k}\neq 0$,  or 

(2) for $b=1,\dots,c$ there is $k=k(b)<0$ such that $\Res_{\eta_b,\be_k-\eta_b}L_{\be_k}\neq 0.$

\noindent
By (Cus1) and (Con3), only (1) is possible, and in that case, using also (Cus2), we must have $c=1$ and $\eta_1=m_\de\de$. Since $\um\geq _l \un$, we conclude that $n_\de=m_\de$, and 
$$
\Res_{\lan\um\ran} \Stand(\un,\bnu)\cong L(\bnu)\boxtimes \Res_{\lan\um'\ran}^{\theta-m_\de\de} \Stand(\un',\bnu). 
$$
Since $\height(\theta-m_\de\de)<\height(\theta)$, we can now apply the inductive hypothesis. 

Case 3: $m_\be=0$ for all $\be\geq \de$. This case is similar to Case 1. 
\end{proof}

\subsection{Classification of irreducible modules}\label{SRough}
We continue to work with a fixed convex preorder $\preceq$ on $\Phi_+$. We will prove:

\begin{Theorem} \label{THeadIrr}
For a given convex preorder there exists a cuspidal system $\{L_\be,L_\bmu\mid \be\in \Phi_+^\re,\,\bmu\in\Par^\ell\}$. Moreover: 
\begin{enumerate}
\item[{\rm (i)}] For every root partition $(\um,\bmu)$, the standard module  
$
\Stand(\um,\bmu)
$ has an irreducible head; denote this irreducible module $L(\um,\bmu)$.\index{l@$L(\um,\bmu)$} 

\item[{\rm (ii)}] $\{L(\um,\bmu)\mid (\um,\bmu)\in \Par(\theta)\}=\Irr(R_\theta)$. 

\item[{\rm (iii)}] $[\Stand(\um,\bmu):L(\um,\bmu)]=1$, and $[\Stand(\um,\bmu):L(\un,\bnu)]\neq 0$ implies $(\un,\bnu)\leq (\um,\bmu)$. 

\item[{\rm (iv)}] $\Res_{\lan\um\ran}L(\um,\bmu)\cong L_{\um,\bmu}$ and $\Res_{\lan\un\ran}L(\um,\bmu)\neq 0$ implies $\un\leq \um$.  

\item[{\rm (v)}] Let $\be\in \Phi_+^\re$ and $n\in\Z_{>0}$. Then $\De(\be^n)=L_\be^{\circ n}$ is irreducible, i.e. $L(\be^n)=L_\be^{\circ n}$. 

\item[{\rm (vi)}] Let $n\in\Z_{>0}$. Then $\{L(\bmu)\mid\bmu\in \Par^\ell(n)\}=\Irr(\bar R_{n\de})$\index{l@$L(\bmu)$} and $\{L(\be^n)\}=\Irr(\bar R_{n\be})$ for any $\be\in\Phi_+^\re$. 
\end{enumerate}
\end{Theorem}

The rest of \S\ref{SRough} is devoted to the proof of Theorem~\ref{THeadIrr}. Since part (vi) of Theorem~\ref{THeadIrr} follows from parts (ii) and (iv), it suffices to prove the following statements for all $\theta\in Q_+$ by induction on $\height(\theta)$: 
\begin{enumerate}
\item[{\rm (1)}] For each $\be\in\Phi_+^\re$ with $\height(\be)\leq\height(\theta)$ there exists a unique up to isomorphism irreducible $L_\be\in\mod{R_\be}$ which satisfies the property (Cus1). Moreover, $L_\be^{\circ n}$ is irreducible 
if $\height(n\be)\leq\height(\theta)$.

\item[{\rm (2)}] For each $n\in\Z_{\geq 0}$ with $\height(n\de)\leq \height(\theta)$ there exist irreducible $\{L(\bmu)\in\mod{R_{n\de}}\mid \bmu\in\Par_n\}$ which satisfy the property (Cus2).

\item[{\rm (3)}] The standard modules $\Stand(\um,\bmu)$ for all $(\um,\bmu)\in\Par(\theta)$, defined as in (\ref{EStand}) using the modules from (1) and (2), satisfy the properties (i)--(iv) of Theorem~\ref{THeadIrr}.
\end{enumerate}

The induction starts with $\height(\theta)=0$, and for $\height(\theta)=1$ the theorem is also clear since $R_{\al_i}$ is a polynomial algebra, which has only the trivial irreducible denoted $L_{\al_i}$. The inductive assumption will stay valid throughout \S\ref{SRough}. 

Let $\Par_{\operatorname{nc}}(\theta)\subseteq\Par(\theta)$ be the set of all root partitions $(\um,\bmu)\in \Par(\theta)$, such that   that $m_\ga\neq 0\neq m_\beta$ for some 
$\ga\neq \beta$ in $\Psi$. Note that for $(\um,\bmu)\in \Par_{\operatorname{nc}}(\theta)$, the modules $\De(\um,\bmu)$ and $L_{\um,\bmu}$ are already defined by induction. 

\begin{Proposition} \label{PHeadIrr} 
Let $(\um,\bmu),(\un,\bnu)\in \Par_{\operatorname{nc}}(\theta)$. 
\begin{enumerate}
\item[{\rm (i)}] 
$
\Stand(\um,\bmu)
$ has an irreducible head; denote this irreducible module $L(\um,\bmu)$. 

\item[{\rm (ii)}] If $(\un,\bnu)\neq (\um,\bmu)$, then $L(\um,\bmu)\not\cong L(\un,\bnu)$. 

\item[{\rm (iii)}] $[\Stand(\um,\bmu):L(\um,\bmu)]=1$, and $[\Stand(\um,\bmu):L(\un,\bnu)]\neq 0$ implies $(\un,\bnu)\leq (\um,\bmu)$. 

\item[{\rm (iv)}] $\Res_{\lan\um\ran}L(\um,\bmu)\cong L_{\um,\bmu}$ and $\Res_{\lan\un\ran}L(\um,\bmu)\neq 0$ implies $\un\leq \um$.  
\end{enumerate}
\end{Proposition}

\begin{proof}
(i) If $L$ is an irreducible quotient of $\Stand(\um,\bmu)=\Ind_{\lan\um\ran}L_{\um,\bmu}$, then by adjointness of $\Ind_{\lan\um\ran}$ and $\Res_{\lan\um\ran}$ and the irreducibility of the $R_{\lan\um\ran}$-module $L_{\um,\bmu}$, which holds by the inductive assumption, $L_{\um,\bmu}$ is a submodule of $\Res_{\lan\um\ran} L$. 
On the other hand, by Proposition~\ref{P1}(i), $[\Res_{\lan\um\ran} \Stand(\um,\bmu):L_{\um,\bmu}]=1$, so (i) follows. 

(iv) Note that we have also proved the first statement in (iv), while the second statement in (iv) follows from Proposition~\ref{P1}(ii) and the exactness of $\Res_{\lan\um\ran}$. 

(iii) Suppose that $[\Stand(\um,\bmu):L(\un,\bnu)]\neq 0$. By (iv), $\Res_{\lan\un\ran}L(\un,\bnu)\cong L_{\un,\bnu}\neq 0$. Therefore $\Res_{\lan\un\ran}\Stand(\um,\bmu)\neq 0$ by exactness of $\Res_{\lan\un\ran}$. By Proposition~\ref{P1}, we then have $\un\leq \um$ and the first equality in (iii). If $\un=\um$, and $\bnu\neq \bmu$, then $[\Res_{\lan\um\ran}\Stand(\um,\bmu):L_{\un,\bnu}]\neq 0$  contradicts (iv). 

(ii) If $L(\um,\bmu)\cong L(\un,\bnu)$, then we deduce from (iii) that $(\um,\bmu)\leq (\un,\bnu)$ and $(\un,\bnu)\leq (\um,\bmu)$, whence $(\um,\bmu)=(\un,\bnu)$. 
\end{proof}

Suppose that $\theta=n\de$ for some $n\in\Z_{\geq 0}$. Then Proposition~\ref{PHeadIrr}, yields $|\Par(\theta)|-|\Par^\ell(n)|$ pairwise non-isomorphic irreducible graded $R_\theta$-supermodules, namely the modules $L(\um,\bmu)$ corresponding to the root partitions $(\um,\bmu)$ such that $m_\rho\neq 0$ for some $\rho\in \Phi_+^\re$. 
By Lemma~\ref{LAmount}, there are exactly $|\Par^\ell(n)|$ irreducible graded $R_{n\de}$-supermodules left, so we can label the remaining  irreducible graded $R_{n\de}$-supermodules by the elements of $\Par^\ell(n)$ in {\em some} way. We get the irreducible graded $R_{n\de}$-supermodules 
$\{L(\bmu)\mid\bmu\in\Par^\ell(n)\}$, and then $\{L(\um,\bmu)\mid (\um,\bmu)\in\Par(n\de)\}=\Irr(R_{n\de})$. Our next goal is Lemma~\ref{LMcNamaraImag} which proves that the modules $\{L(\bmu)\mid\bmu\in\Par^\ell(n)\}$ are cuspidal. 

We need some terminology. Let $(\um,\bmu)$ be a root partition.
We define the {\em support}\index{support} $\supp(\um):=\{\be\in \Psi\mid m_\be\neq 0\}$.\index{s@$\supp(\um)$} We denote by $\max(\um)$ the largest root appearing in $\supp(\um)$. Suppose that  $\beta\in\Phi_+$ satisfies $\beta\succeq \max(\um)$. If $\be$ is real then $L_\be\circ \Stand(\um,\bmu)$ is again a standard module. If $\be=n\de$ is imaginary, $\bnu\in\Par^\ell(n)$, and $\max(\um)$ is real, then $L(\bnu)\circ \Stand(\um,\bmu)$ is again a standard module (in this case $m_\de=0$ and $\bmu=\varnothing^\ell$). 

\begin{Lemma} \label{LMcNamaraImag}
Let $\theta=n\de$ and  $\bla\in\Par^\ell(n)$. Then $L(\bla)$ is cuspidal. 
\end{Lemma}
\begin{proof}
To prove the lemma it suffices to show that for non-zero $\zeta,\eta\in Q_+\setminus\Phi_+^\im$ such that $n\de=\zeta+\eta$ and $\Res_{\zeta,\eta}L(\bla)\neq 0$, we must have that $\zeta$ is a sum of real roots $\prec\de$ and $\eta$ is a sum of real roots  $\succ\de$. 
We prove that $\zeta$ is a sum of real roots $\prec\de$, the proof that $\eta$ is a sum of real roots $\succ\de$ being similar. 

Let $L(\um,\bmu)\boxtimes L(\un,\bnu)$, with $(\um,\bmu)\in\Par(\zeta)$ and $(\un,\bnu)\in\Par(\eta)$, be an irreducible submodule of $\Res_{\zeta,\eta} L(\bla)$. Note that $\height(\zeta),\height(\eta)<\height(\theta)$, so the modules $L(\um,\bmu), L(\un,\bnu)$ are defined by induction.

Let $\be:=\max(\um)$. If $\be\preceq \de$, then, since $\zeta$ is not an imaginary root, $\zeta$  is a sum of real roots less than $\de$. So we may assume that $\be\succ\de$. 
 Moreover, $\Res_{\be,\zeta-\be}L(\um,\bmu)\neq 0$, and hence $\Res_{\be,n\de-\be}L(\bla)\neq0$. So we may assume from the beginning that $\zeta\in\Phi^+$, $\zeta\succ\de$ and $L(\um,\bmu)\cong L_\zeta$. Moreover, we may assume that $\zeta$ is the largest possible real root for which $\Res_{\zeta,\eta} L(\bla)\neq 0$. 

Now, let $\ga:=\max(\un)$. If $\ga$ is real, we have the cuspidal module $L_\ga$. If $\ga$ is imaginary, we interpret  $L_\ga$ as $L(\bnu)$. 
We have a non-zero map $L_\zeta\boxtimes L_\ga\boxtimes V\to \Res_{\zeta,\ga,\eta-\ga}L(\bla)$, for some non-zero $V\in\mod{R_{\eta-\ga}}$. By adjunction, this yields a non-zero map
$$
f: (\Ind_{\zeta,\ga} L_\zeta\boxtimes L_\ga)\boxtimes V\to \Res_{\zeta+\ga,\eta-\ga}L(\bla)
$$

Suppose $\ga=\eta$. Then $\zeta\neq \eta$, since $\zeta,\eta\not\in\Phi_+^\im$. If  $\zeta\succ\eta$ then $L(\bla)$ is a quotient of the standard module $L_\zeta\circ L_\eta$, which contradicts the definition of $L(\bla)$. So $\zeta\prec\eta$,  and since $n\de=\zeta+\ga$, we have by (Con3) that $\zeta\prec \de\prec\eta$, giving a contradiction. 

Suppose $\ga\neq\eta$. Pick a composition factor  $L(\um',\bmu')$ of $\Ind_{\zeta,\ga} L_\zeta\boxtimes L_\ga$, which is not in the kernel of $f$. As $\zeta$ is maximal with $\Res_{\zeta,\eta} L(\bla)\neq 0$, every root $\ga'$ in $\supp(\um')$ satisfies $\ga'\preceq \zeta$. 
So $\zeta+\ga$ is a sum of roots $\preceq \zeta$. Now 
(Con2) implies that $\ga\preceq \zeta$, and so by adjointness, $L(\bla)$ is a quotient of the standard module $L_\zeta\circ\Stand(\un,\bnu)$, which is a contradiction. 
\end{proof}

Suppose now that $\theta=\be\in\Phi_+^\re$. Then 
Proposition~\ref{PHeadIrr} yields $|\Par(\be)|-1$ irreducible graded $R_\be$-supermodules, corresponding to the root partitions $\Par_{\operatorname{nc}}(\be)=\Par(\be)\setminus\{(\be)\}$. We define $L_\be$\index{l@$L_\be$} to be the missing irreducible graded $R_\be$-supermodule, cf. Lemma~\ref{LAmount}. Then $\{L(\um,\bmu)\mid(\um,\bmu)\in\Par(\be)\}$ is a complete irredundant system of irreducible $R_\be$-modules up to isomorphism. 


\begin{Lemma} \label{LMcNamara}
Let $\theta=\be\in\Phi_+^\re$. Then $L_\be$ is cuspidal. 
\end{Lemma}
\begin{proof}
To prove the lemma it suffices to show that for $\zeta,\eta\in Q_+\setminus\{0\}$ such that $\be=\zeta+\eta$ and $\Res_{\zeta,\eta}L_\be\neq 0$, we must have that  $\zeta$ is a sum of roots $\prec\be$ and $\eta$ is a sum of roots $\succ \be$. 
We prove that $\zeta$ is a sum of roots $\prec\be$, the proof that $\eta$ is a sum of roots $\succ\be$ being similar. 

Let $L(\um,\bmu)\boxtimes L(\un,\bnu)$, with $(\um,\bmu)\in\Par(\zeta)$ and $(\un,\bnu)\in\Par(\eta)$, be an irreducible submodule of $\Res_{\zeta,\eta} L_\be$. Let $\ga:=\max(\um)$.  Then $\Res_{\ga,\zeta-\ga}L(\um,\bmu)\neq 0$, and hence $\Res_{\ga,\beta-\ga}L_\be\neq0$. If we can prove that $\ga$ is a sum of roots less than $\be$, then by (Con1), (Con3), $\ga$ is a root less than $\be$,  whence, by the maximality of $\ga$, we have that $\zeta$ is a sum of roots less than $\be$. 
So we may assume from the beginning that $\zeta$ is a root and $L(\um,\bmu)=L_\zeta$ (if $\zeta$ is imaginary, $L_\zeta$ is interpreted as $L(\bmu)$). Moreover, we may assume that $\zeta$ is the largest possible root for which $\Res_{\zeta,\eta} L_\be\neq 0$. 

Let $\al:=\max(\un)$. If $\al$ is real we have the cuspidal module $L_\al$. If $\al$ is imaginary, we interpret $L_\al$ as $L(\bnu)$. We have a non-zero map 
$L_\zeta\boxtimes L_\al\boxtimes V\to \Res_{\zeta,\al,\eta-\al}L_\be,$ 
for some $V\in\mod{R_{\eta-\al}}$. By adjunction, this yields a non-zero map
$$
f: (\Ind_{\zeta,\al} L_\zeta\boxtimes L_\al)\boxtimes V\to \Res_{\zeta+\al,\eta-\al}L_\be.
$$

If $\al=\eta$, then we must have $\zeta\prec\eta$, for otherwise $L_\be$ is a quotient of the standard module $L_\zeta\circ L_\eta$, which contradicts the definition of the cuspidal module $L_\be$. Now, since $\be=\zeta+\al$, we have by (Con1) that $\zeta\prec\be\prec\eta$, in particular $\zeta\prec\be$. 

Next, let $\al\neq\eta$, and pick a composition factor  $L(\um',\bmu')$ of $\Ind_{\zeta,\al} L_\zeta\boxtimes L_\al$, which is not in the kernel of $f$. 
Since $\zeta$ is maximal with $\Res_{\zeta,\eta} L_\be\neq 0$,  every root $\al'$ in the support of $\um'$ satisfies $\al'\preceq \zeta$. 
Thus $\zeta+\al$ is a sum of roots $\preceq \zeta$. If $\zeta$ and $\al$ are not both imaginary, then (Con2) implies that $\al\preceq \zeta$, and so by adjointness, $L_\be$ is a quotient of the standard module $L_\zeta\circ\De(\un,\bnu)$, which is a contradiction. 

If $\zeta$ and $\al$ are both imaginary, then $\De(\un,\bnu)=L(\bnu)\circ\De(\un',\varnothing^\ell)$ for $\un'$ such that $\max(\un')\prec\de$. In this case, we have by adjunction that $L_\be$ is a quotient of $L(\bmu)\circ L(\bnu)\circ L(\un',\varnothing^\ell)$. It now follows from Lemma~\ref{LTensImagIsImag} that $L_\be$ is a quotient of the standard module of the form $L(\bla)\circ L(\un',\varnothing^\ell)$ for some composition factor $L(\bla)$ of $L(\bmu)\circ L(\bnu)$, so we get a contradiction again. 
\end{proof} 

\begin{Lemma} \label{LCuspPower} 
Let $\be\in \Phi_+^\re$ and $n\in\Z_{>0}$. Then $L_\be^{\circ n}$ is irreducible. 
\end{Lemma}
\begin{proof}
In view of Proposition~\ref{PHeadIrr}, we have the irreducible graded $R_{n\be}$-supermodules $L(\um,\bmu)$ for all root partitions $(\um,\bmu)\in\Par(n\be)$, except for $(\um,\bmu)=(\be^n)$ for which $\Stand(\be^n)=L_\be^{\circ n}$. By (Con1), we have that $(\be^n)$ is the smallest element of $\Par(n\be)$. By Proposition~\ref{PHeadIrr}(v), we conclude that 
$L_\be^{\circ n}$ has only one composition factor $L$ appearing with certain multiplicity $c$, and $L\not\cong L(\um,\bmu)$ for all $(\um,\bmu)\in\Par_{\operatorname{nc}}(n\be)$. By Corollary~\ref{CPowerIrr}, $c=1$. 
\end{proof}

The proof of Theorem~\ref{THeadIrr} is now complete. 

\subsection{Induction and restriction for cuspidal modules}
We continue to work with a fixed convex order $\preceq$. 
The material of this subsection is parallel to \cite[\S4.6]{EK2}.

Let $d\in\Z_{>0}$ and $\la=(\la_1,\dots,\la_n)$ be a composition of $d$. 
For $\be\in\Psi$, we denote $\la\be:=(\la_1\be,\dots,\la_n\be)\in Q_+^n$. In particular, we have the idempotent $e_{\la\be}\in R_{d\be}$ and the parabolic subalgebra 
$R_{\la\be}\subseteq e_{\la\be}R_{d\be}e_{\la\be}.$ Recall that $R_{\la\be}$ has been identified with $R_{\la_1\be}\otimes\dots\otimes R_{\la_n\be}$ via the embedding $\iota_{\la\be}$ of (\ref{EIota}). So we can consider $\bar R_{\la_1\be}\otimes\dots\otimes \bar R_{\la_n\be}$ as the quotient of $R_{\la\be}$, and we call a module $M\in\mod{R_{\la\be}}$ cuspidal if it factors through $\bar R_{\la_1\be}\otimes\dots\otimes \bar R_{\la_n\be}$; in other words $M$ is cuspidal  if and only if $e(\bi^{(1)}\cdots\bi^{(n)})M\neq 0$ for $\bi^{(1)}\in I^{\la_1\be},\dots, \bi^{(n)}\in I^{\la_n\be}$ implies $\bi^{(1)}\in I^{\la_1\be}_\cus,\dots, \bi^{(n)}\in I^{\la_n\be}_\cus$. 
The following lemma shows that cuspidality is preserved under parabolic restriction. The analogous result for parabolic induction is Lemma~\ref{LTensImagIsImag}. 

\begin{Lemma} \label{LResCusp}
If $M\in\mod{R_{d\be}}$ is cuspidal then so is 
$\Res_{\la\be}M\in\mod{R_{\la\be}}$. 
\end{Lemma}
\begin{proof}
Suppose $L_1\boxtimes \dots\boxtimes L_n$ is a composition factor of $\Res_{\la\be}M$ with $L_r\in\mod{R_{\la_r\be}}$ for $r=1,\dots,n$. We need to prove that $L_1,\dots,L_n$ are cuspidal. If not, let $s$ be minimal such that $L_s$ is non-cuspidal. Then $L_s$ corresponds to a root partition $(\um,\bmu)$ with $\ga:=\max\um\succ \be$. By Theorem~\ref{THeadIrr}(iv), we deduce that $\Res_{m_\ga\ga,\la_s\be-m_\ga\ga}L_s\neq 0$. 
Let $\theta=\sum_{t=1}^{s-1}\la_t\be+m_\ga\ga$. 
It follows that $\Res_{\theta,d\be-\theta}M\neq 0$. Since $M$ is cuspidal, we must have that $\sum_{t=1}^{s-1}\la_t\be+m_\ga\ga$ is a sum of the roots $\preceq \be$, which contradicts Lemma~\ref{LCones}. 
\end{proof}

\begin{Corollary} \label{C130921} 
If $\bi^{(1)}\in I^{\la_1\be},\dots,\bi^{(n)}\in I^{\la_n\be}$ and the concatenation word $\bi^{(1)}\cdots\bi^{(n)}$ is cuspidal then so are $\bi^{(1)},\dots,\bi^{(n)}$.
\end{Corollary}
\begin{proof}
Apply Lemma~\ref{LResCusp} to the cuspdial $\bar R_{d\be}$-module $\bar R_{d\be}e(\bi^{(1)}\cdots\bi^{(n)})$. 
\end{proof}

\begin{Lemma} \label{LIndCuspProj} 
If $\bi^{1}\in I^{\la_1\be}_\cus,\dots,\bi^{n)}\in I^{\la_n\be}_\cus$, then there is an isomorphism of $R_{d\be}$-modules
\begin{align*}
\bar R_{d\be}e({\bi^{(1)}\cdots\bi^{(n)}})&\iso \bar R_{\la_1\be}e({\bi^{(1)}})\circ\dots\circ \bar R_{\la_n\be}e({\bi^{(n)}}),
\\
e({\bi^{(1)}\cdots\bi^{(n)}})
&\mapsto
e_{\la\be}\otimes e({\bi^{(1)}})\otimes \dots\otimes e({\bi^{(n)}}).\end{align*}
\end{Lemma}
\begin{proof}
Since $\bar R_{\la_1\be}e({\bi^{(1)}})\circ\dots\circ \bar R_{\la_n\be}e({\bi^{(n)}})$ is cuspidal by Lemma~\ref{LTensImagIsImag}, we can consider it as a $\bar R_{d\be}$-module. So there exists a homomorphism as in the lemma. To construct the inverse homomorphism, use adjointness of induction and restriction together with Lemma~\ref{LResCusp}.
\end{proof}

Define the {\em cuspidal parabolic subalgebra}\index{cuspidal parabolic subalgebra}
$\bar R_{\la\be}\subseteq e_{\la\be} \bar R_{d\be}e_{\la\be}$\index{r@$\bar R_{\la\be}$} 
to be the image of $R_{\la\be}$ under the natural projection $e_{\la\be}R_{d\be}e_{\la\be}\onto e_{\la\be}\bar R_{d\be}e_{\la\be}$. Whereas the parabolic subalgebra $R_{\la\be}$ has been identified with $R_{\la_1\be}\otimes\dots\otimes R_{\la_n\be}$ via the embedding $\iota_{\la_1\be,\dots,\la_n\be}$ of (\ref{EIota}), it is not obvious that $\bar R_{\la\be}\cong \bar R_{\la_1\be}\otimes\dots\otimes \bar R_{\la_n\be}$. This is proved in the following lemma:

\begin{Lemma}\label{L030216}
Let $\be\in\Psi$, $h=\height(\be)$ and $h\la:=(h\la_1,\dots,h\la_n)$.
\begin{enumerate}
\item[{\rm (i)}] The natural map
$$
R_{\la_1\be}\otimes\dots\otimes R_{\la_n\be}\stackrel{\iota_{\la\be}}{\longrightarrow} e_{\la\be}R_{d\be}e_{\la\be}\onto e_{\la\be}\bar R_{d\be}e_{\la\be}
$$
factors through $\bar R_{\la_1\be}\otimes\dots\otimes \bar R_{\la_n\be}$ and induces an isomorphism
$$\bar R_{\la_1\be}\otimes\dots\otimes \bar R_{\la_n\be}\iso \bar R_{\la\be}.$$ 
\item[{\rm (ii)}] $\bar R_{d\be} e_{\la\be}$ is a free right $\bar R_{\la\be}$-module with basis
$\{ \psi_we_{\la\be} \mid w\in \D^{h\la}\}$.
\item[{\rm (iii)}] $e_{\la\be}\bar R_{d\be}$ is a free left $\bar R_{\la\be}$-module with basis
$\{ e_{\la\be}\psi_w \mid w\in {}^{h\la}\D\}$.
\end{enumerate}
\end{Lemma}
\begin{proof}
The map factors through $\bar R_{\la_1\be}\otimes\dots\otimes \bar R_{\la_n\be}$ thanks to Lemma~\ref{LResCusp}.
Consider the $R_{d\be}$-module
$W:=\bar R_{\la_1\be}\circ\dots\circ \bar R_{\la_n\be}$. By Lemma~\ref{LTensImagIsImag}, the module $W$ factors through $\bar R_{d\be}$. On the other hand, by Theorem~\ref{TBasis}, we can decompose 
$$W=\bigoplus_{w\in \D^{h\la}} \psi_w e_{\la\be}\otimes \bar R_{\la_1\be}\otimes\dots\otimes \bar R_{\la_n\be}$$ as a vector space, with each summand being naturally isomorphic to $\bar R_{\la_1\be}\otimes\dots\otimes \bar R_{\la_n\be}$ as a vector space. Parts (i) and (ii) follow. Part (iii)  follows from part (ii) using the antiautomorphism (\ref{ETauAntiAuto}).
\end{proof}

In view of the lemma we identify $\bar R_{\la_1\be}\otimes\dots\otimes \bar R_{\la_n\be}$ with $\bar R_{\la\be}$. Then:

\begin{Corollary} 
Suppose that $W_r\in\mod{\bar R_{\la_r\be}}$ for $r=1,\dots,n$. Then there is a natural isomorphism of cuspidal $R_{d\be}$-modules
\begin{align*}
W_1\circ\dots\circ W_n&\iso \bar R_{d\be}e_{\la\be}\otimes_{\bar R_{\la\be}} (W_1\boxtimes\dots\boxtimes W_n), \\
u e_{\la\be}\otimes w_1\otimes\dots\otimes w_n&\mapsto
\bar u e_{\la\be}\otimes w_1\otimes\dots\otimes w_n,
\end{align*}
where $\bar u\in \bar R_{d\be}$ is the image of $u\in R_{d\be}$ under the natural projection $R_{d\be}\onto \bar R_{d\be}$.
\end{Corollary}

From now on we identify the induced modules as in the corollary.

\chapter{RoCK blocks of quiver Hecke superalgebras}
\label{ChRockqHs}
\section{RoCK blocks}
\label{SRock}
Let $d\in\Z_{\geq 0}$. In this section we always work with abaci which have at least $d$ beads on the $0^{\nth}$ runner, cf. \S\ref{SSAb}.

\subsection{Rouquier cores}
\label{SSRoCore}
In this subsection we define Rouquier cores and RoCK blocks (after Rouquier, Chuang and Kessar). 
These notions are inspired by the corresponding notions for symmetric groups \cite{RoTh,CK}.

Recall the notation (\ref{EBJ}) and let $\rho$ be $\bar p$-core and $d\in Z_{\geq 0}$. Then $\rho$ is called a {\em $d$-Rouquier $\bar p$-core}\index{d@$d$-Rouquier $\bar p$-core}\index{r@$d$-Rouquier $\bar p$-core} if $b_1^\rho\geq d$ and $b_j^\rho-b_{j-1}^\rho\geq d-1$ for all $j$ with $1<j\leq\ell$. It is then automatic that $b_j^\rho=0$ for $j$ with $\ell<j< p$. If $\rho$ is a $d$-Rouquier $\bar p$-core, we refer to the cyclotomic quiver Hecke superalgebra $R^{\La_0}_{\cont(\rho)+d\de}$ as a {\em RoCK block}\index{RoCK block} of {\em weight $d$}. We will also consider the higher level cyclotomic algebras $R^{N\La_0}_{N\cont(\rho)+d\de}$ (but do not refer to them as RoCK blocks.)

Recall the notation (\ref{EParRhoD}) and the notion of an elementary slide down from \S\ref{SSAb}. The following lemma follows easily from the definitions. 
Part (iii) is an analogue of \cite[Lemma 4(1)]{CK}, while part (ii) is an analogue of \cite[Lemma 1.1(1)]{Paget}.  

\begin{Lemma} \label{LCK1} 
Let $\rho$ be a $d$-Rouquier $\bar p$-core and $\la\in\Par_p(\rho,c)$ for $c\leq d$. Then:
\begin{enumerate}
\item[{\rm (i)}] $\Ab_\la$ is obtained from $\Ab_\rho$ by $c$ consecutive elementary slides down on runners $0,1,\dots,\ell$. In particular, all positions on runners $\ell+1,\dots,p-1$ are not occupied. 

\item[{\rm (ii)}] $\la$ is $p$-restricted if and only if 
$\Ab_\la$ is obtained from $\Ab_\rho$ by $c$ consecutive elementary slides down on runners $0,1,\dots,\ell-1$. 

\item[{\rm (iii)}] Suppose the position in row $m$ on runner $j$ of $\Ab_\la$ is occupied. Then the positions in rows $<m$ on runners $j+1,\dots,\ell$  of $\Ab_\la$ are occupied. If $c<d$ then the positions in rows $\leq m$ on runners $j+1,\dots,\ell$  of $\Ab_\la$ are occupied.
\end{enumerate} 
\end{Lemma}

The following lemma is an analogue of \cite[Lemma 4(2)]{CK}.

\begin{Lemma}\label{ChuKesLem} 
Let $\rho$ be a $d$-Rouquier $\bar p$-core, $\la\in\Par_p(\rho,d)$ and $\mu\in\Par_p(\rho,d-1)$ is such that $\mu\subseteq \la$. Then there exists $i\in I$ such that:
\begin{enumerate}
\item[{\rm (i)}] $\Ab_\la$ is obtained from $\Ab_\mu$ by an elementary slide down on runner $i$.
\item[{\rm (ii)}] There exists $s> \ell-i$ such that $\res(s-\ell+i,\mu_{s-\ell+i}+1)=\ell$, and 
$$
\la_t-\mu_t=
\left\{
\begin{array}{ll}
\ell+1+i &\hbox{if $t=s-\ell+i$,}\\
1 &\hbox{if $s-\ell+i< t\leq s$,}\\
0 &\hbox{otherwise.}
\end{array}
\right.
$$
\end{enumerate}  
\end{Lemma}
\begin{proof}
Let $N>h(\la)+d$, and write 
$$\la=(\la_1,\dots,\la_N),\quad \mu=(\mu_1,\dots,\mu_N).$$ 
We use the abaci $\Ab_\la$ and $\Ab_\mu$ with $N$ beads. 

By assumption we have $\la_t\geq \mu_t$ for $t=1,2\dots,N$. 
Let $s$ be such that $\la_s>\mu_s$ and $\la_t=\mu_t$ for all $t>s$. Let the position $a:=\mu_s$ be in row $m$ on runner $i$, so $a=i+mp$. Note that the positions $b$ with $0\leq b<a$ on $\Ab_\la$ and $\Ab_\mu$ are identical, the position $a$ in $\Ab_\mu$ is occupied with $k>0$ beads ($k=1$ unless $i=0$), while the position $a$ in $\Ab_\la$ is occupied with $l<k$ beads. 

By Lemma~\ref{LCK1}(i), the amounts of beads on runner $i$ in $\Ab_\la$ and $\Ab_\mu$ are the same (if $i=0$, we do count beads in position $0$). 
So in $\Ab_\la$, there must be a bead on runner $i$  in row $>m$. By Lemma~\ref{LCK1}(iii), in both $\Ab_\la$ and $\Ab_\mu$ the positions in row $m$ on runners $j$ satisfying $i<j\leq\ell$ are occupied. Moreover, again by Lemma~\ref{LCK1}(iii), in $\Ab_\la$ the positions $m+1$ on runners $j$ satisfying $0\leq j<i$ are {\em not} occupied. Hence $\mu_t\leq\la_t-1$ for all $t$ satisfying $s-(\ell-i-1)\leq t\leq s$, and $\mu_{s-\ell+i}\leq \la_{s-\ell+i}-(p+i-\ell)$. 
As $\mu\subset\la$ and $|\mu|=|\la)-p$, the inequalities must be equalities, and we also must have $\mu_t=\la_t$ for $t<s-\ell+i$. Thus $\Ab_\la$ is obtained from $\Ab_\mu$ by moving a bead on runner $i$ from row $m$ to row $m+1$.
\end{proof}

\begin{Remark}\label{R010921}
The condition (ii) in the lemma can be interpreted as the fact that for some $0\leq i\leq \ell$, the boxes with residues of the skew shape $\la\setminus\mu$ look as follows:
\vspace{2mm}
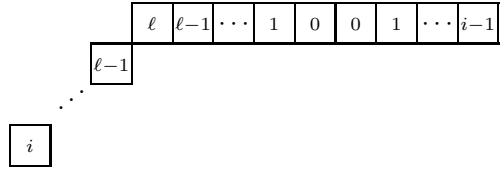
$$
\ytableausetup{mathmode}
\begin{ytableau}
\none[\,] & \none & \none & {\scriptstyle\ell} & {\scriptstyle\ell-1} & \cdots &{\scriptstyle 1} & {\scriptstyle 0} & {\scriptstyle 0} & {\scriptstyle 1} & \cdots & {\scriptstyle i-1} \cr
\none & \none & {\scriptstyle\ell-1} \cr
\none & \none[\iddots] \cr
{\scriptstyle i} \cr
\end{ytableau}
$$

\vspace{2mm}
\noindent
If $i=0$ this is interpreted as  
\vspace{2mm}
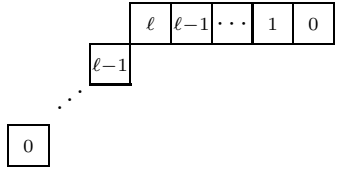
$$
\ytableausetup{mathmode}
\begin{ytableau}
\none[\,] & \none & \none & {\scriptstyle\ell} & {\scriptstyle\ell-1} & \cdots &{\scriptstyle 1} & {\scriptstyle 0} \cr
\none & \none & {\scriptstyle\ell-1} \cr
\none & \none[\iddots] \cr
{\scriptstyle 0} \cr
\end{ytableau}
$$
\end{Remark}

\vspace{4mm}
For $i=0,\dots,\ell$ we denote
\begin{equation}\label{EJIKI}
\begin{split}
\bj^{i}&:=(\ell-1)\,(\ell-2)\,\cdots \,1\,0\,0\,1\,\cdots\, (i-1)\in I^{\ell+i},
\index{j@$\bj^{i}$}
\\
\quad\bk^i&:=(\ell-1)\,(\ell-2)\,\cdots\, i\in I^{\ell-i}.
\index{k@$\bk^{i}$}
\end{split}
\end{equation}
If $i=0$, we interpret $\bj^0=\bk^0=(\ell-1)\,(\ell-2)\,\cdots\, 1\,0$. 

\begin{Corollary} \label{CSkewShuffle} 
Let $\rho$ be a $d$-Rouquier $\bar p$-core, $\la\in\Par_p(\rho,d)$ and $\mu\in\Par_p(\rho,d-1)$ such that $\mu\subseteq \la$.
Suppose that $\T\in\Std_p(\mu)$ and $\Stab\in\Std_p(\la)$ satisfy $\Stab|_{\{1,\dots,n\}}=\T$ for $n:=|\mu|$. 
Then there exists $i\in I$ such that $\Ab_\la$ is obtained from $\Ab_\mu$ by an elementary slide down on runner $i$ and $\bi^{\Stab}=\bi^{\T}\bi$ for $\bi\in I^\de$ such that $i_1=\ell$ and $i_2\cdots i_p$ is a shuffle of the words $\bj^{i},\,\bk^i$. 
\end{Corollary}

\begin{Lemma}\label{Lcoreab}
Let $\rho$ be a $d$-Rouquier $\bar p$-core, $r:=|\rho|$, $\T\in\Std_p(\rho)$ and  $\bi^\T=i_1\cdots i_{r}$.
\begin{enumerate}
\item[{\rm (i)}] If $1\leq m\leq r$, then then $|\{k\mid r-m+1\leq k\leq r\ \text{and}\ i_{k}=\ell\}|<m/p$.
\item[{\rm (ii)}] Let $\la\in\Par_p(r+s)$ for $s\leq dp$, $\Stab\in\Std_p(\la)$ and suppose that $\bi^\Stab=\bi^\T\bj$ for some $\bj=j_1\cdots j_{s}\in I^{s}$. Then $\{k\mid j_k=\ell\}\geq s/p$, and the equality holds if and only if $\core(\la)=\rho$. 
\end{enumerate} 
\end{Lemma}

\begin{proof}
(i) Let $\mu:=\T(\{1,\dots,r-m\})\in\Par_p(r-m)$. Recall the notation of Lemma~\ref{LComputation}. 
We need to prove $c_\ell(\rho)-c_\ell(\mu)<(|\rho|-|\mu|)/p$. By 
Lemma~\ref{LComputation}, this is equivalent to 
$$
\sum_{i=\ell+1}^{p-1}(p-i)b_i^\rho-\sum_{i=1}^\ell ib_i^\rho
<
\sum_{i=\ell+1}^{p-1}(p-i)b_i^\mu-\sum_{i=1}^\ell ib_i^\mu,
$$
By definition of a $d$-Rouquier core, we have $b_j^\rho=0$ for $\ell<j<p$. So the last inequality is equivalent to 
$$
-\sum_{i=1}^\ell ib_i^\rho
<
\sum_{i=\ell+1}^{p-1}(p-i)b_i^\mu-\sum_{i=1}^\ell ib_i^\mu,
$$
which is true because $\Ab_\mu$ is obtained from $\Ab_\rho$ by moving beads to smaller positions. 

(ii) Note by Lemma~\ref{LLT} that $\rho=\Stab(\{1,\dots,r\})\subseteq \la$. We need to prove $c_\ell(\la)-c_\ell(\rho)\geq (|\la|-|\rho|)/p$ and that the equality holds if and only if $\core(\la)=\rho$. By 
Lemma~\ref{LComputation} and the definition of a $d$-Rouquier core, the inequality $c_\ell(\la)-c_\ell(\rho)\geq (|\la|-|\rho|)/p$  is equivalent to 
\begin{equation}\label{ERequired}
\sum_{i=\ell+1}^{p-1}(p-i)b_i^\la-\sum_{i=1}^\ell ib_i^\la
\geq -\sum_{i=1}^\ell ib_i^\rho.
\end{equation}
Note that $\Ab_\la$ is obtained from $\Ab_\rho$ by making $s\leq dp$ consecutive moves of the following kind: move a bead from position $k$ to position $k+1$ if position $k+1$ is not occupied or if position $k+1$ is on runner $0$. It follows from the definition of a $d$-Rouquier core that 
\begin{equation}\label{EPartialSum}
\sum_{j=i}^\ell b_j^\la\leq \sum_{j=i}^\ell b_{j}^\rho\qquad(i=\ell, \ell-1,\dots,1). 
\end{equation}
Hence 
$$-\sum_{i=1}^\ell ib_i^\la
\geq -\sum_{i=1}^\ell ib_i^\rho,
$$ 
which implies (\ref{ERequired}).  

If we have an equality in (\ref{ERequired}) then we must have  equalities in (\ref{EPartialSum}) and $b_i^\la$ must be $0$ for all $i=\ell+1,\dots,p-1$. It follows that $\Ab_\la$ is obtained from $\Ab_\rho$ by consecutive elementary slides down on runners $0,1,\dots,\ell$. Conversely, if $\core(\la)=\rho$ then $\la\in\Par(\rho,c)$ and $\bj\in I^{c\de}$ for some $c\leq d$, which implies $s=cp$ and $\{k\mid j_k=\ell\}=c= s/p$. 
\end{proof}

\begin{Corollary} \label{CIneq} 
Let $\rho$ be a $d$-Rouquier $\bar p$-core, $r=|\rho|$, $N\in\Z_{\geq 1}$, $\T\in\Std(\rho^N)$ and  $\bi^\T=i_1\cdots i_{Nr}$. If $1\leq m\leq Nr$, then 
$$|\{k\mid r-m+1\leq k\leq r\ \text{and}\ i_{k}=\ell\}|<m/p.$$
\end{Corollary}
\begin{proof}
This follows immediately from Lemma~\ref{Lcoreab}(i).
\end{proof}

Recall the notation (\ref{EParNRhoD}). 

\begin{Lemma} \label{LMichael} 
Let $\rho$ be a $d$-Rouquier $\bar p$-core and $\bla=(\la^{(1)},\dots,\la^{(N)})\in\Par_p^N$ for some $N\in\Z_{\geq 1}$. 
Suppose there exists $\T\in\Std_p(\bla)$ such that $\bi^\T=\bj\bk$ for $\bj\in I^{N\cont(\rho)}$ and $\bk\in I^{d\de}$. Then 
$\bla\in \Par_p^N(\rho^N,d)$ and $\T(\{1,\dots,N|\rho|\})=\rho^N$. 
\end{Lemma}
\begin{proof}
By Lemma~\ref{LRhoNUnique}, we have 
$$\T(\{1,\dots,N|\rho|\})=\rho^N\subseteq \bla.$$ 
By Lemma~\ref{Lcoreab}(ii) applied to each $\la^{(t)}$ with $t=1,\dots,N$, we have that 
$$c_\ell(\la^{(t)})-c_\ell(\rho)\geq (|\la^{(t)}|-|\rho|)/p,$$ with equality if and only if $\core(\la^{(t)})=\rho$. By assumption, we have 
$$
\sum_{t=1}^N(c_\ell(\la^{(t)})-c_\ell(\rho))=d=\sum_{t=1}^N(|\la^{(t)}|-|\rho|)/p.
$$
So $$c_\ell(\la^{(t)})-c_\ell(\rho)= (|\la^{(t)}|-|\rho|)/p$$ and $\core(\la^{(t)})=\rho$ for all $t$, as required. 
\end{proof}

\subsection{RoCK blocks and cuspidal algebras}
\label{SSRoCKCusp}
{\em From now on, we fix a convex preorder $\preceq$ on $\Phi_+$ as in Example~\ref{ExConPr}}, in particular (\ref{ESharp}) holds. When we speak of cuspidality from now on, we mean cuspidality with respect to this specially chosen $\preceq$. Recall the words $\bj^i$ and $\bk^i$ from (\ref{EJIKI}).

\begin{Lemma} \label{LCuspExpl}
Let $\bi=i_1\cdots i_p\in I^\de$. Then $\bi\in I^\de_\cus$ if and only if $i_1=\ell$ and $i_2\cdots i_p$ is a shuffle of the words $\bj^{i},\,\bk^i$ for some $1\leq i\leq \ell$.
\end{Lemma}
\begin{proof}
Using \cite[Table III]{Bourb} and the notation of Example~\ref{ExConPr}, we get 
\begin{align*}\Phi'_\sharp=&\Big\{\sum_{r\leq t<s}\al_t\mid 1\leq r<s\leq \ell\Big\}
\\
&\cup
\Big\{-\sum_{r\leq t<s}\al_t-2\sum_{s\leq t<\ell}\al_t-\al_\ell\mid1\leq r<s\leq \ell\Big\}
\\
&\cup
\Big\{-2\sum_{r\leq t<\ell}\al_t-\al_\ell\mid1\leq r\leq \ell\Big\}.
\end{align*}

For $1\leq n\leq p$ denote 
$$
\theta_n:=\al_{i_1}+\dots+\al_{i_n}.
$$ 
Suppose $\bi\in I^\de_\cus$ and $1\leq n<p$. Then  
$\theta_n$ is a sum of roots from $\Phi^\re_{\prec\de}$  and 
$\de-\theta_n$ is a sum of roots from $\Phi^\re_{\succ\de}$. Consider the roots 
\begin{align*}
\al(r,s)&:=\sum_{r\leq t<s}\al_t+2\sum_{s\leq t<\ell}\al_t+\al_\ell\in\Phi^\re_+
\qquad (1\leq r\leq s\leq \ell),
\\
\be(r,s)&:= 2\sum_{0\leq t<r}\al_t+\sum_{r\leq t<s}\al_t+2\sum_{s\leq t<\ell}\al_t+\al_\ell\in\Phi^\re_+
\qquad (0\leq r< s\leq \ell).
\end{align*}
By (\ref{ESharp}), (\ref{EPhiS})-(\ref{EPhiL}), and the description of $\Phi'_\sharp$ above, 
$\theta_n$ is a sum of roots from 
\begin{align*}
R:=&\{\al(r,s)
\mid 1\leq r\leq s\leq \ell \}
\cup \{\be(r,s)
\mid 0\leq r< s\leq \ell \},
\end{align*}
while $\de-\theta_n$ is a sum of roots from 
$$
\{2\sum_{0\leq t<r}\al_t+\sum_{r\leq t<s}\al_t\mid 0\leq r<s\leq \ell \}\cup\{\sum_{r\leq t<s}\al_t\mid 1\leq r<s\leq \ell \}
$$
Since $\al_\ell$ appears in $\de$ with coefficient $1$, we deduce that in fact $\theta_n\in R$. Thus $\bi\in I^\de_\cus$ implies $\theta_n\in R$ for all $1\leq n<p$. The converse is also true since $\al_0,\dots,\al_{\ell-1}\succ\de$ and so $\theta_n\in R$ implies that $\de-\theta$ is a sum of roots $\succ\de$. We have proved that $\bi$ is cuspidal if and only if $\theta_1,\dots,\theta_{p-1}\in R$. 

Now, if the word $\bi\in I^\de$ is such that $i_1=\ell$ and $i_2\cdots i_p$ is a shuffle of $\bj^i$ and $\bk^i$,we clearly have $\theta_1,\dots,\theta_{p-1}\in R$. Conversely, suppose that $\theta_1,\dots,\theta_{p-1}\in R$. Then $i_1=\ell$. Moreover, there exists $2\leq m<p$ such that $\theta_1,\dots,\theta_{m-1}$ are of the form $\al(r,s)$ and $\theta_m,\dots,\theta_{p-1}$ are of the form $\be(r,s)$. We have the following options:

\begin{enumerate}
\item[{\rm (A1)}] $\theta_n=\al(r,s)$, $i_{n+1}=r-1$ and 
$$\theta_{n+1}=
\left\{
\begin{array}{ll}
\al(r-1,s) &\hbox{if $r>1$,}\\
\be(0,s) &\hbox{if $r=1$.}
\end{array}
\right.
$$
\item[{\rm (A2)}] $\theta_n=\al(r,s)$, $i_{n+1}=s-1$ and 
$$\theta_{n+1}=
\left\{
\begin{array}{ll}
\al(r,s-1) &\hbox{if $r<s$,}\\
\al(r-1,r) &\hbox{if $r=s> 1$,}\\
\be(0,1) &\hbox{if $r=s= 1$.}
\end{array}
\right.
$$
\item[{\rm (B1)}] $\theta_n=\be(r,s)$, $i_{n+1}=r$ and 
$\theta_{n+1}=\be(r+1,s)$ (interpret $\be(s,s)$ as $\de$). 

\item[{\rm (B2)}] $\theta_n=\be(r,s)$, $i_{n+1}=s-1$ and 
$\theta_{n+1}=\be(r,s-1)$ (interpret $\be(r,r)$ as $\de$). 
\end{enumerate}

Now define two words $\bj(n)$ and $\bk(n)$ for $n=0,1,\dots,p-1$ and  as follows: for $n=0$ we set $\bj(0)=\bk(0)=\varnothing$;  let $n>0$ and suppose we have already constructed $\bj(n-1)$ and $\bk(n-1)$. If we are in cases (A1) or (B1) we add $i_{n+1}$ to the right end of $\bj(n-1)$ and do not change $\bk(n-1)$; otherwise we add $i_{n+1}$ to the right end of $\bk(n-1)$ and do not change $\bj(n-1)$. It is clear that $i_2\cdots i_p$ is a shuffle of $\bj(p-1)$ and $\bk(p-1)$. It remains to  notice that $\bj(p-1)=\bj^i$ and $\bk(p-1)=\bk^i$ for some $1\leq i\leq \ell$.
\end{proof}

\begin{Example}\label{ExNonCusp}
For any $i\in J$, the word 
$$
\ell(\ell-1)^{2}\cdots (i+1)^{2}i\cdots 1 0^{2}1\cdots i\in I^\de
$$
is $\ell$ followed by a shuffle of the words $\bj^{i+1}$ and $\bk^{i+1}$, and so it is cuspidal by Lemma~\ref{LCuspExpl}. 
On the other hand, using by Lemma~\ref{LCuspExpl}, one checks: 
\begin{enumerate}
\item[{\rm (a)}] 
$\ell(\ell-1)^2\cdots(i+1)^2i\cdots (j+1)(j-1)j(j-2)\cdots10^21\cdots i$ is non-cuspidal for $1\leq j\leq i<\ell$.  
\item[{\rm (b)}] 
$\ell(\ell-1)^2\cdots(i+1)^2i\cdots 10^21\cdots(j-1)(j+1)j(j+2)\cdots i$ is non-cuspidal for $1\leq j\leq i-2$.  
\item[{\rm (c)}] any word in $I^\de$ starting with 
$\ell(\ell-1)^2\cdots (j+1)^2j(j-1)^2$, for $2\leq j<\ell$, is non-cuspidal. 
\end{enumerate}
\end{Example}

By Example~\ref{ExExtraWord}, any shuffle of $\bj^0$ and $\bk^0$ can be obtained as a shuffle of $\bj^1$ and $\bk^1$. So from Corollary~\ref{CSkewShuffle} and Lemma~\ref{LCuspExpl}, we immediately conclude:

\begin{Lemma} \label{LSkewCusp} 
Let $\rho$ be a $d$-Rouquier $\bar p$-core, $\la\in\Par_p(\rho,d)$ and $\mu\in\Par_p(\rho,d-1)$ is such that $\mu\subseteq \la$.
Suppose that $\T\in\Std_p(\mu)$ and $\Stab\in\Std_p(\la)$ satisfy $\Stab|_{\{1,\dots,n\}}=\T$ for $n:=|\mu|$. Then $\bi^{\Stab}=\bi^{\T}\bj$ for  $\bj\in I^\de_\cus$. 
\end{Lemma}

Recall that by Lemma~\ref{LTensImagIsImag} a shuffle of cuspidal words is cuspidal. 

\begin{Proposition} \label{LCuspSecWord} 
Let $\rho$ be a $d$-Rouquier $\bar p$-core and set $r:=|\rho|$. 
Let $\bla
\in\Par_p^N$ and $\T\in\Std_p(\bla)$. Suppose that 
$\bi^\T=\bk\bj$ for 
$\bk\in I^{N\cont(\rho)}$ and  $\bj\in I^{d\de}$. Then:
\begin{enumerate}
\item[{\rm (i)}] $\T(\{1,\dots,Nr\})=\rho^N$.
\item[{\rm (ii)}]  $\bla\in \Par_p^N(\rho^N,d)$.
\item[{\rm (iii)}] $\bj$ is a shuffle of $d$ words in $I^\de_\cus$; in particular, $\bj\in I^{d\de}_\cus$.
\end{enumerate} 
\end{Proposition}
\begin{proof}
Part (i) follows from Lemma~\ref{LRhoNUnique}, and part (ii) from 
 Lemma~\ref{LMichael}. 

(iii) Write $\bi^\T=i_1\cdots,i_{Nr+dp}$.
Since $\rho=\core(\la^{(t)})$ for every $t=1,\dots,N$, there is a sequence 
$$\rho^N=\bla^0\subset \bla^1\subset\dots\subset\bla^d=\la$$ such that $\bla^c\in \Par_p^N(\rho^N,c)$ for $c=0,\dots,d$. For $c=1,\dots,d$ let 
$$
\T^{-1}(\bla^c\setminus \bla^{c-1})=\{s_{c,1}<\dots<s_{c,p}\}.
$$ 
By Lemma~\ref{LSkewCusp}, $\bj^c:=i_{s_{c,1}}\cdots i_{s_{c,p}}\in I^\de_\cus$. Since $\T(\{1,\dots,Nr\})=\rho^N$, we get the partition of the set $\{1,\dots,dp\}$ into the subsets 
$$\{s_{b,1}-Nr,\dots,s_{b,p}-Nr\},\quad (b = 1,\dots,d).
$$ 
This partition witnesses the fact that $\bj$ is a shuffle of the cuspidal words $\bj^1,\dots,\bj^d$.
\end{proof}

For $d\in\Z_{\geq 0}$ and a 
$\bar p$-core $\rho$,  recalling (\ref{ECycPar}), (\ref{EIota}) and (\ref{EPi}), we consider the algebra homomorphism  
\begin{equation}\label{EOmega}
\Omega^{N\La_0}_{\rho,d}\colon R_{d\de}\to R_{N\cont(\rho),d\de}^{N\La_0},\ x\mapsto \pi_\theta^{N\La_0}\big(\iota_{N\cont(\rho),d\de}(e_{N\cont(\rho)}\otimes x)\big).
\index{o@$\Omega^{N\La_0}_{\rho,d}$}
\index{w@$\Omega^{N\La_0}_{\rho,d}$}
\index{$\Omega^{N\La_0}_{\rho,d}$}
\end{equation}

\begin{Lemma}\label{LOmegaFact} 
Let $d\in\Z_{\geq 0}$ and $\rho$ be a $d$-Rouquier $\bar p$-core
If $\bi\in I^{d\de}\setminus I^{d\de}_\cus$ then $\Om^{N\La_0}_{\rho,d}(e(\bi))=0$. In particular, $\Omega^{N\La_0}_{\rho,d}$ factors through $\bar R_{d\de}$. 
\end{Lemma}
\begin{proof}
Let $\bj\in I^{d\de}$. Then $\Om^{N\La_0}_{\rho,d}(e(\bj))\neq 0$ implies that $e(\bk\bj)\neq 0$ in $R_{N\cont(\rho)+d\de}^{N\La_0}$ for some $\bk\in I^{N\cont(\rho)}$. By Corollary~\ref{CIdNZ} there exists $\bla\in\Par_p^N$ and $\T\in\Std_p(\bla)$ such that $\bi^\T=\bk\bj$.
Now $\bj$ is cuspidal by Proposition~\ref{LCuspSecWord}. 
\end{proof}

In view of the lemma, whenever $\rho$ is a $d$-Rouquier $\bar p$-core, we will consider $\Om^{N\La_0}_{\rho,d}$ as a homomorphism 
\begin{equation}\label{EOm}
\Om^{N\La_0}_{\rho,d}: \bar R_{d\de}\to R^{N\La_0}_{N\cont(\rho),d\de}.
\end{equation}

\subsection{Further properties of RoCK blocks}
\label{SSFurtherRoCK}
Throughout the subsection we fix $N\in \Z_{\geq 1}$, $d\in\Z_{\geq 0}$ and a $d$-Rouquier $\bar p$-core $\rho$.  Set 
$\theta := N\cont(\rho)+d\de\in Q_+.$ 

\begin{Lemma}\label{Ptens1}
We have $$e_{N\cont(\rho),d\de}R^{N\La_0}_{\theta}e_{N\cont(\rho),d\de}=R^{N\La_0}_{N\cont(\rho),d\de}.$$
\end{Lemma}

\begin{proof} 
By Theorem~\ref{TMackeyKL}(i), the algebra $e_{N\cont(\rho),d\de}R_{\theta}^{N\La_0}e_{N\cont(\rho),d\de}$ is generated by $R_{N\cont(\rho),d\de}^{N\La_0}$ together with the elements $\psi_w$ for $w\in {}^{(|\rho|,dp)}\D^{(|\rho|,dp)} \setminus \{1\}$. Hence it suffices to prove  that $$e_{N\cont(\rho),d\de}\psi_we_{N\cont(\rho),d\de}=0$$ in $R^{N\La_0}_{\theta}$ for all such $w$. 
 If not, then there exist $\bi,\bl\in I^{N\cont(\rho)}$ and $\bj,\bk\in I^{d\de}$ such that $e(\bi\bk) \psi_w e(\bl \bj)\ne 0$. By Corollary~\ref{CIdNZ}, $\bi\bk=\bi^\T$ and $\bl\bj=\bi^\Stab$ with $\T,\Stab\in\Std_p(\bla)$ for some $\bla\in\Par_p$. By Proposition~\ref{LCuspSecWord}, 
  $\bj$ and $\bk$ are shuffles of $d$ words in $I^\de_\cus$, and 
  $$\T(\{1,\dots,Nr\})=\Stab(\{1,\dots,Nr\})=\rho^N,$$
 where $r:=|\rho|$. 
 In particular, $\bi=\bi^{\T'}$ for $\T'\in\Std_p(\rho^N)$. 
Moreover, 
$$w=\prod_{t=1}^m (Nr-m+t, Nr+t)$$ for some $m>0$, and therefore
 the last $m$ entries of $\bi$ are 
 $j_1,\ldots,j_m$. Since 
 $\bj$ is a shuffle of $d$ words in $I^\de_\cus$, in view of Lemma~\ref{LCuspExpl}, we have $|\{k\mid k\leq m\ \text{and}\ j_k=\ell\}|\geq m/p$. This contradicts Corollary~\ref{CIneq}. 
\end{proof}

Recall from (\ref{ECoreSubalg}) (and Corollary~\ref{CZeta}) that we have identified 
$
R_{N\cont(\rho)}^{N\La_0}$ with the subalgebra 
$$\zeta_{N\cont(\rho),d\de}^{N\La_0}(R_{N\cont(\rho)}^{N\La_0})\subseteq e_{N\cont(\rho),d\de}R_{N\cont(\rho)+d\de}^{N\La_0}e_{N\cont(\rho),d\de}.
$$ 
By Lemma~\ref{Ptens1}, the latter algebra is $R^{N\La_0}_{N\cont(\rho),d\de}$, so the supercentralizer $\Cent_{\rho,d}^{N\La_0}$ from (\ref{ECentralizer}) becomes 
\begin{equation}\label{ECent}
\Cent_{\rho,d}^{N\La_0}=
\Cent_{R_{N\cont(\rho),d\de}^{N\La_0}} (R^{N\La_0}_{N\cont(\rho)}),
\end{equation}
and Lemmas~\ref{L270116New},~\ref{LIdSymNew} become:

\begin{Lemma} \label{L270116} 
We have an isomorphism of graded superalgebras  
$$R^{N\La_0}_{N\cont(\rho)} \otimes \Cent_{\rho,d}^{N\La_0} \iso R_{N\cont(\rho),d\de}^{N\La_0},\ a\otimes b \mapsto ab.
$$ 
\end{Lemma}

\begin{Lemma} \label{LIdSym} 
There is an idempotent $\eps_d\in R_{N\cont(\rho)+d\de}^{N\La_0}$\index{e@$\eps_d$} such that 
$$\eps_d e_{N\cont(\rho),d\de}=\eps_d=e_{N\cont(\rho),d\de}\eps_d$$ 
and 
$$\eps_d R_{N\cont(\rho)+d\de}^{N\La_0}\eps_d\cong \Cent_{\rho,d}^{N\La_0}.$$
In particular,  $\Cent_{\rho,d}^{N\La_0}$ is a symmetric algebra whenever $R_{N\cont(\rho)+d\de}^{N\La_0}$ is so. 
\end{Lemma}

Recall the map $\Om_{\rho,d}^{N\La_0}$ from (\ref{EOmega}).

\begin{Lemma}\label{Ltens2}
We have $\Cent_{\rho,d}^{N\La_0}=\Omega_{\rho,d}^{N\La_0}(\bar R_{d\de})$. 
\end{Lemma}

\begin{proof}
It is clear from the definitions that $\Omega_{\rho,d}^{N\La_0}(\bar R_{d\de}) \subseteq \Cent_{\rho,d}^{N\La_0}$.  

Conversely, let $x\in \Cent_{\rho,d}^{N\La_0}$. 
Since $x$ is an element of $R_{N\cont(\rho),d\de}^{N\La_0}$, it  can be written as 
$x= \sum_{i=1}^m a_i b_i$ for some $a_1,\ldots,a_m \in R^{N\La_0}_{N\cont(\rho)}$
and $b_1,\ldots,b_m \in \Omega_{\rho,d}^{N\La_0}(\bar R_{d\de})$, and we may assume that $a_1,\ldots,a_m$ are linearly independent and  $a_1=1$. 
By Lemma~\ref{L270116},  
$x=b_1$, so $x\in \Omega_{\rho,d}^{N\La_0}(\bar R_{d\de})$. 
\end{proof}

\section{Gelfand-Graev truncation of \texorpdfstring{$\bar R_\de$}{} and \texorpdfstring{$\Cent_{\rho,1}$}{}}
Recall that, starting with Section~\ref{SRock}, we have fixed a convex preorder $\preceq$ on $\Phi_+$ as in Example~\ref{ExConPr}, in particular (\ref{ESharp}) holds. When we speak of cuspidality from now on, we mean cuspidality with respect to this specially chosen preorder. 

Throughout the section we fix a $1$-Rouquier $\bar p$-core $\rho$ and $N\in\Z_{\geq 1}$.

In this subsection, we define an idempotent $\ga\in \bar R_\de$ such that $\ga \bar R_\de\ga$ is the basic algebra Morita equivalent to $\bar R_\de$ and construct an explicit isomorphism of graded superalgebras $f: \Zig_\ell[\zz]\iso \ga \bar R_\de\ga$, where $\Zig_\ell$ is the Brauer tree superalgebra  introduced in \S\ref{SSZig} and $\Zig_\ell[\zz]$ is a twisted polynomial algebra with coefficients in $\Zig_\ell$ and degree $4$ variable $\zz$, see \S\ref{SSAff}. 

Importantly, the isomorphism $f$ will induce an isomorphism of graded superalgebras $\Zig_\ell\iso \ga \Cent_{\rho,1}^{\La_0}\ga$, and $\ga \Cent_{\rho,1}^{\La_0}\ga$ will turn out to be the basic algebra that is graded Morita superequivalent to $ \Cent_{\rho,1}^{\La_0}$ and to the weight $1$ Rock block $R^{\La_0}_{\cont(\rho)+\de}$.

\subsection{Gelfand-Graev idempotents}\label{SSGG}
Recall the notation $\hat\bi$ from \S\ref{ChBasicNotLie}. 
For all $i\in J$ we will consider {\em Gelfand-Graev words} \index{Gelfand-Graev word}
\begin{equation}\label{EGGW}
\begin{split}
\bg^i&:=\ell\,(\ell-1)^{(2)}\,\cdots\, (i+1)^{(2)}\,i\,\cdots\, 1\, 0^{(2)}\,1\,\cdots\, i\,\in\, I^\de_\di,
\index{g@$\bg^i$}
\\
\hat\bg^i&=\widehat{\bg^i}=\ell\,(\ell-1)^{2}\,\cdots\, (i+1)^{2}\,i\,\cdots \,1\, 0^{2}\,1\,\cdots\, i\,\in \,I^\de.
\index{g@$\hatbg^i$}
\end{split}
\end{equation}
We denote the entries of the word $\hat\bg^i$ by $g^i_1,\dots, g^i_p$, so that:
\begin{equation}\label{EHatGI}
\hat\bg^i=g^i_1\cdots g^i_p\qquad(i\in J).
\index{g@$g^i_r$}
\end{equation}

Recalling (\ref{EJIKI}), note that $\hat\bg^i\in I^{\de}$ is $\ell$ followed by a shuffle of the words $\bj^{i+1}$ and $\bk^{i+1}$, so by Lemma~\ref{LCuspExpl}, $\hat\bg^i$ is cuspidal. 
So it makes sense to consider the {\em Gelfand-Graev idempotents}\index{Gelfand-Graev idempotent}
$$
\ga^i:=e(\bg^i)\in \bar R_\de\quad\text{and}\qquad
\ga:=\sum_{i\in J}\ga^i.
\index{g@$\ga^i$}\index{c@$\ga^i$}\index{g@$\ga$}\index{c@$\ga$}
$$

\begin{Example} 
{\rm 
Let $\ell=2$. Then $\ga^0=\psi_2\psi_4y_3y_5e(21100)$, $\ga^1=\psi_3y_4e(21001)$. 
}
\end{Example}

Recalling Lemma~\ref{Ltens2}, we use the same symbols $\ga^i$ and $\ga$ to denote the Gelfand-Graev idempotents in $\Cent_{\rho,1}^{N\La_0}$: 
$$\ga^i:=\Om^{N\La_0}_{\rho,1}(\ga^i)\in \Cent_{\rho,1}^{N\La_0}\quad\text{and}\quad\ga:=\Om^{N\La_0}_{\rho,1}(\ga)\in \Cent_{\rho,1}^{N\La_0}.$$
We consider the idempotent truncations 
\begin{equation}\label{ECTr}
B:=\ga \bar R_{\de}\ga \qquad \text{and}\qquad
Y_{\rho,1}^{N\La_0}:=\ga\Cent_{\rho,1}^{N\La_0}\ga. 
\index{b@$B$}\index{y@$Y_{\rho,1}^{N\La_0}$}
\end{equation}
The surjective homomorphism $\Om^{N\La_0}_{\rho,1}:\bar R_\de\to \Cent_{\rho,1}^{N\La_0}$ restricts to the surjective homomorphism denoted by the same symbol: 
\begin{equation}\label{ECTrRes}
\Om_{\rho,1}^{N\La_0}:B\to 
Y_{\rho,1}^{N\La_0}. 
\end{equation}

By Lemma~\ref{LIdSym}, we have that
\begin{equation}\label{LIdTrCentd=1}
Y_{\rho,1}^{N\La_0}\cong \ga\eps_1 R_{N\cont(\rho)+\de}^{N\La_0} \eps_1\ga
\end{equation}
is an idempotent truncation of the algebra $R_{N\cont(\rho)+\de}^{N\La_0}$, so by \cite[Theorem IV.4.1]{SY}:

\begin{Lemma} 
The algebra $Y_{\rho,1}^{N\La_0}$ is symmetric whenever $R_{N\cont(\rho)+\de}^{N\La_0}$ is so. 
\end{Lemma}

\subsection{Graded dimension formula for $Y_{\rho,1}$}
\label{SSGrDimYRho1}

Throughout this subsection we set $r:=|\rho|$. 

Recall the numbers $a_i$ and $a_i^\vee$ from \ref{ChBasicNotLie}. 

\begin{Lemma} \label{LRhoRemAdd}
We have 
$$
\prod_{i\in I}q_i^{a_i(|\Add_i(\rho)|-|\Rem_i(\rho)|)}=q^2.
$$
\end{Lemma}
\begin{proof}
We have 
\begin{align*}
\prod_{i\in I}q_i^{a_i(|\Add_i(\rho)|-|\Rem_i(\rho)|)}
&=
q^{2(|\Add_0(\rho)|-|\Rem_0(\rho)|)+4\sum_{i=1}^\ell(|\Add_i(\rho)|-|\Rem_i(\rho)|)}
\\
&=q^{2\sum_{i\in I}a_i^\vee(|\Add_i(\rho)|-|\Rem_i(\rho)|)},
\end{align*}
and the result follows from Lemma~\ref{LemmaK}.
\end{proof}

For $i\in J$ and $k\in I$, we now define
\begin{equation}\label{EChi}
\chi^i_k:=
\begin{cases}
q(q^2 + q^{-2})^{\ell-i-1}(q^2+1)&\text{if }i=k-1,\\
q^{-1}(q^2 + q^{-2})^{\ell-i-1}(q^2+1)&\text{if }i=k,\\
0&\text{otherwise.}
\end{cases}
\end{equation}

\begin{Lemma}\label{lem:deg_calc1}
Let $k\in I$, $i\in J$, and $\U\in\Std_p(\rho)$. Let $\la \in \Par_p(r+p)$ be such that the abacus $\Ab_{\la}$ is obtained from $\Ab_{\rho}$ by performing an elementary slide down on runner $k$. 
Denote 
$$
\Std^{\U,i}_p(\la):=\{\Stab \in \Std_p(\la)\mid 
\Stab_{\leq r}=\U\ \text{and}\ \bi^{\Stab} = \bi^\U\hat\bg^i\}.
\index{s@$\Std^{\U,i}_p(\la)$}
$$ 
Then 
\begin{align*}
\sum_{\Stab \in \Std^{\U,i}_p(\la)}\frac{\deg(\Stab)}{\deg(\U)}=
\chi^i_k.
\end{align*}
\end{Lemma}
\begin{proof}
Throughout the proof we abbreviate $\sum_\Stab:=\sum_{\Stab \in \Std^{\U,i}_p(\la)}$. 

By assumption, $\la\in \Par_p(\rho,1)$. So 
by Remark~\ref{R010921}, we have 
\vspace{2mm}
$$
\la\setminus\rho =
\vcenter{\hbox{
\ytableausetup{mathmode}
\begin{ytableau}
\none[\,] & \none & \none & {\scriptstyle\ell} & {\scriptstyle\ell-1} & \cdots &{\scriptstyle 1} & {\scriptstyle 0} & {\scriptstyle 0} & {\scriptstyle 1} & \cdots & {\scriptstyle k-1} \cr
\none & \none & {\scriptstyle\ell-1} \cr
\none & \none[\iddots] \cr
{\scriptstyle k} \cr
\end{ytableau}}}
$$

\vspace{2mm}
\noindent
interpreted as 
\vspace{2mm}
$$
\ytableausetup{mathmode}
\begin{ytableau}
\none[\,] & \none & \none & {\scriptstyle\ell} & {\scriptstyle\ell-1} & \cdots &{\scriptstyle 1} & {\scriptstyle 0} \cr
\none & \none & {\scriptstyle\ell-1} \cr
\none & \none[\iddots] \cr
{\scriptstyle 0} \cr
\end{ytableau}
$$

\vspace{2mm}
\noindent
if $i=0$.  
Since the last entry of $\hat\bg^i$ is $i$, we conclude that $\sum_\Stab \deg(\Stab)/\deg(\U)=0$ unless $i=k-1$ or $k$. We therefore assume from now on that $i=k-1$ or $k$.

Recall that $p=2\ell+1$. Using the notation (\ref{EHatGI}), for  $u=1,\dots,p$, we define 
\begin{equation}\label{EQU}
q(u):=q_{g_u^i}=
\left\{
\begin{array}{ll}
q &\hbox{if }g_u^i=0,\\
q^4 &\hbox{if }g_u^i=\ell,\\
q^2 &\hbox{otherwise}.
\end{array}
\right.
\end{equation}

For $\Stab \in \Std^{\U,i}_p(\la)$ and $1\leq u\leq p$, we denote  
$A_u(\Stab):=\Stab(r+u)$. Let   
\begin{align*}
\la(u,\Stab)&:=\la\setminus\{A_{u+1}(\Stab),\dots,A_p(\Stab)\}\in\Par_p(r+u)\qquad(0\leq u\leq p),
\end{align*} 
and set 
$$
\eta_u(\Stab):=q(u)^{\eta^{A(u,\Stab)}(\la(u-1,\Stab))}\quad\text{and}\quad \zeta_u(\Stab):=\zeta^{A(u,\Stab)}(\la(u-1,\Stab))\qquad(1\leq u\leq p).
$$
By (\ref{EDeg}), we have 
$$\frac{\deg(\Stab)}{\deg(\U)}=\prod_{u=1}^{p}\eta_u(\Stab)\zeta_u(\Stab).
$$

Note that 
\begin{equation}\label{E030921_1}
\prod_{u=1}^{p}\zeta_u(\Stab)=
\left\{\begin{array}{ll} 1+q^2  &\hbox{if $k\neq 0$,}\\ 1 &\hbox{if $k=0$,}\end{array}\right.
\end{equation}
since when $k\neq 0$ there is exactly one $u$ such that $A_u(\Stab)$ is added to a row of $\la(u-1,\Stab)$ of length divisible by $p$, while when $k=0$ there is no such $u$. In particular, the product $\prod_{u=1}^{p}\zeta_u(\Stab)$ does not depend on $\Stab$, and we denote it simply by $\zeta$. Now 
\begin{equation}\label{E030921}
\sum_{\Stab}\frac{\deg(\Stab)}{\deg(\U)}=
\zeta\cdot\sum_{\Stab} \prod_{u=1}^{p}\eta_u(\Stab). 
\end{equation}

To complete the proof, we now calculate $\sum_\Stab \prod_{u=1}^p\eta_u(\Stab).$ We abbreviate $b_j:=b_j^\rho$ and consider three cases. 

\vspace{2mm}
{\sf Case 1:} $k=\ell$. In this case we must have $i=\ell-1$, $\la-\rho$ is one row, so there is only one tableau $\Stab$ in $\Std_p^{(\U,i)}(\la)$. For this $\Stab$, we have  
$$
\eta^{A(u,\Stab)}(\la(u-1,\Stab))=
\left\{
\begin{array}{ll}
b_\ell-1 &\hbox{if $u=1$,}\\
1+b_{\ell-1}-b_\ell &\hbox{if $u=2$,}\\
b_{u-1}-b_u &\hbox{if $3\leq u\leq \ell$,}\\
1-2b_1 &\hbox{if $u= \ell+1$,}\\
-2b_1 &\hbox{if $u= \ell+2$,}\\
b_{u-\ell-2}-b_{u-\ell-1} &\hbox{if $\ell+3\leq u\leq 2\ell$,}
\\
1+b_{\ell-1}-b_{\ell} &\hbox{if $u= 2\ell+1$.}
\end{array}
\right.
$$ 
Since $q(1)=q^4$, $q(\ell+1)=q(\ell+2)=q$, and $q(u)=q^2$ for all other $u$, we deduce that 
$$\sum_\Stab \prod_{u=1}^p\eta_u(\Stab)=q,$$ 
and 
$$\zeta\cdot\sum_\Stab \prod_{u=1}^p\eta_u(\Stab)=q(1+q^2),$$ 
as desired.  

\vspace{2mm}
\noindent
{\sf Case 2:} $1\leq k\leq \ell-1$ and $i=k$. In this case there are $2^{\ell-k-1}$ tableaux in $\Std_p^{(U,i)}$. Moreover, for $S\in \Std_p^{(U,i)}$, we have:
\begin{align*}
\eta_1(\Stab)&=q_\ell^{b_k-1}=q^{4(b_k-1)},
\\
(\eta_2(\Stab),\eta_3(\Stab))&=(1,q_{\ell-1}^{-1})=(1,q^{-2})\quad \text{or} \quad(q_{\ell-1},1)=(q^2,1),
\\
&\vdots
\\
(\eta_{2\ell-2k-2}(\Stab),\eta_{2\ell-2k-1}(\Stab))&=(1,q_{k+1}^{-1})=(1,q^{-2})\quad \text{or} \quad(q_{k+1},1)=(q^2,1),
\\
\eta_{2\ell-2k}(\Stab)&=q_k=q^2,\\
\eta_{2\ell-2k+1}(\Stab)&=q_{k-1}^{b_{k-1}-b_k}=q^{2(b_{k-1}-b_k)},
\\
&\vdots
\\
\eta_{2\ell-k-2}(\Stab)&=q_1^{b_{1}-b_2}=q^{2(b_{1}-b_2)},
\\
\eta_{2\ell-k-1}(\Stab)&=q_0^{1-2b_{1}}=q^{1-2b_{1}},
\\
\eta_{2\ell-k}(\Stab)&=q_0^{-2b_{1}}=q^{-2b_{1}},
\\
\eta_{2\ell-k+1}(\Stab)&=q_1^{b_{1}-b_2}=q^{2(b_{1}-b_2)},
\\
&\vdots
\\
\eta_{2\ell}(\Stab)&=q_{k-1}^{b_{k-1}-b_k}=q^{2(b_{k-1}-b_k)},
\\
\eta_{2\ell+1}(\Stab)&=q_{k}^{0}=1.
\end{align*}
So 
\begin{align*}
\sum_\Stab \prod_{u=1}^p\eta_u(\Stab)=\,
&q^{4(b_k-1)}(q^2+q^{-2})^{\ell-k-1} q^2 q^{2(b_{k-1}-b_k)}\cdots q^{2(b_1-b_2)} 
\\&\cdot q^{1-2b_1}q^{-2b_1}q^{2(b_1-b_2)} \cdots q^{2(b_{k-1}-b_k)}
\\
=\,&q^{-1}(q^2+q^{-2})^{\ell-k-1},
\end{align*}
and 
$$\zeta\cdot\sum_\Stab \prod_{u=1}^p\eta_u(\Stab)=(1+q^2)q^{-1}(q^2+q^{-2})^{\ell-i-1},$$ 
as desired.

\vspace{2mm}
\noindent
{\sf Case 3:} $1\leq k\leq \ell-1$ and $i=k-1$. In this case, arguing as in Case 2, we get 
\begin{align*}
\sum_\Stab \prod_{u=1}^p\eta_u(\Stab)=\,
&q^{4(b_k-1)}(q^2+q^{-2})^{\ell-k} q^{2(1+b_{k-1}-b_k)}
q^{2(b_{k-2}-b_{k-1})}\cdots q^{2(b_1-b_2)}
\\&\cdot q^{1-2b_1}q^{-2b_1}
q^{2(b_1-b_2)} \cdots q^{2(b_{k-2}-b_{k-1})}q^{2(b_{k-1}-b_k-1)}
\\
=\,&q(q^2+q^{-2})^{\ell-k},
\end{align*}
and 
$$\zeta\cdot\sum_\Stab \prod_{u=1}^p\eta_u(\Stab)=(1+q^2)q(q^2+q^{-2})^{\ell-i-1},$$ 
as desired.

\vspace{2mm}
{\sf Case 4:} $k=0$. In this case we must have $i=0$. 
Arguing as in Case 2, we get 
\begin{align*}
\sum_\Stab \prod_{u=1}^p\eta_u(\Stab)=\,(q^2 + q^{-2})^{l-1}(q+q^{-1}),
\end{align*}
and 
$$\zeta\cdot\sum_\Stab \prod_{u=1}^p\eta_u(\Stab)=(q^2+q^{-2})^{\ell-1}(q+q^{-1})=q^{-1}(q^2+q^{-2})^{\ell-1}(q^2+1),$$ 
as desired.
\end{proof}

\begin{Lemma}\label{lem:deg_calc2}
Let $1\leq s\leq N$, $k\in I$, $i\in J$, and $\U\in\Std_p(\rho^N)$. Let $\bla=(\la^{(1)},\dots,\la^{(N)}) \in \Par_p^N(Nr+p)$ be such that $\la^{(t)}=\rho$ for all $t\neq s$ and the abacus $\Ab_{\la^{(s)}}$ is obtained from $\Ab_{\rho}$ by performing an elementary slide down on runner $k$. 
Denote 
$$
\Std^{\U,i}_p(\bla):=\{\Stab \in \Std_p(\bla)\mid 
\Stab_{\leq Nr}=\U\ \text{and}\ \bi^{\Stab} = \bi^\U\hat\bg^i\}.
\index{s@$\Std^{\U,i}_p(\bla)$}
$$ 
Then 
\begin{align*}
\sum_{\Stab \in \Std^{\U,i}_p(\bla)}\frac{\deg(\Stab)}{\deg(\U)}=
q^{2(N-s)}\chi^i_k.
\end{align*}
\end{Lemma}
\begin{proof}
We explain how to reduce the result to Lemma~\ref{lem:deg_calc1}. 

If $\bmu\in\Par_p^N$ and $\T\in\Std_p(\bmu)$, let $\{i_1< \dots<i_l\}=\T^{-1}(\mu^{(s)})$. Then $\T^{(s)}: k\mapsto \T(i_k)$ for $1\leq k\leq l$ is a $p$-standard $\mu^{(s)}$-tableau. 

Note that 
$$
\Std^{\U,i}_p(\bla)\to \Std^{\U^{(s)},i}_p(\la^{(s)}),\ \Stab\mapsto \Stab^{(s)}
$$
is a bijection. Moreover, by (\ref{EDegM}), using the notation (\ref{EQU}), we have 
\begin{align*}
\frac{\deg(\Stab)}{\deg(\U)}&=\frac{\deg(\Stab^{(s)})}{\deg(\U^{(s)})}\cdot \prod_{u=1}^pq(u)^{(N-s)(|\Rem_{g_u}(\rho)|-|\Add_{g_u}(\rho)|)}
\\
&=\frac{\deg(\Stab^{(s)})}{\deg(\U^{(s)})}\cdot q^{2(N-s)}, 
\end{align*}
where we have used Lemma~\ref{LRhoRemAdd} for the last equality. The result now follows from Lemma~\ref{lem:deg_calc1}.
\end{proof}

For $i,j\in J$, let 
\begin{equation}\label{EMIJ}
m_{i,j}(q):=(1 + q^2)(1+q^{-2})(1 + q^4)^{\ell-i-1}(1 + q^{-4})^{\ell-j-1}.\index{m@$m_{i,j}(q)$}
\end{equation}
Recalling the notation of Lemma~\ref{LFact}, it is easy to check that
\begin{equation}\label{EMIJEq}
m_{i,j}(q)=\bg^i!q^{\lan\bg^i\ran}\bg^j!q^{-\lan\bg^j\ran}.
\end{equation}

Recall the Brauer tree algebra $\Zig_\ell$ introduced in \S\ref{SSZig}, especially (\ref{EZigBiWt}). 

\begin{Lemma} \label{L010921}
Let $i,j\in J$. Then 
$$
\sum_{ \bla \in \Par_p^N(\rho^N,1)}
\| \bla\|
\sum_{\substack{\Stab \in\Std^{\U,i}_p(\bla) \\ \T \in \Std^{\V,j}_p(\bla)}}  \frac{\deg(\Stab)}{\deg(\U)}  \frac{\deg(\T)}{\deg(\V)}
=
m_{i,j}(q)\,(\dim_q \ze^i\Zig_\ell\ze^j)\, \sum_{s=0}^{N-1}q^{4s}.
$$
\end{Lemma}
\begin{proof}
By Lemma~\ref{ChuKesLem}(i), we have 
$$
\Par_p^N(\rho^N,1)=\{\bla^{s,k}\mid 1\leq s\leq N,\  k\in I\},$$
where $\bla^{k,s}=(\la^{(1)},\dots,\la^{(N)})\in\Par_p^N$ is defined from $\la^{(t)}=\rho$ for all $t\neq s$ and the abacus $\Ab_{\la^{(s)}}$ is obtained from $\Ab_{\rho}$ by performing an elementary slide down on runner $k$. Note that $\| \bla^{k,s}\|=(1+q^2)^{\de_{k,0}}$. Now, the left hand side in the lemma equals
\begin{align*}
&\sum_{s=1}^N\sum_{k\in I}
\| \bla^{k,s}\|
\bigg[\sum_{\Stab \in\Std^{\U,i}_p(\bla^{k,s})}  \frac{\deg(\Stab)}{\deg(\U)}\bigg]\bigg[ \sum_{\T \in \Std^{\V,j}_p(\bla^{k,s})} \frac{\deg(\T)}{\deg(\V)}\bigg]
\\
=\,&\sum_{s=1}^N\sum_{k\in I}
(1+q^2)^{\de_{k,0}}
\big[q^{2(N-s)}\chi^i_k\big]\big[ q^{2(N-s)}\chi^j_k\bigg]
\\
=\,&\bigg[\sum_{s=1}^Nq^{4(N-s)}\bigg]\sum_{k\in I}
(1+q^2)^{\de_{k,0}}\,\chi^i_k\, \chi^j_k,
\end{align*}
where we have used Lemma~\ref{lem:deg_calc2} for the first  equality. So it remains to prove that 
\begin{equation}\label{E040921}
\sum_{k\in I}
(1+q^2)^{\de_{k,0}}\,\chi^i_k\, \chi^j_k=m_{i,j}(q)\, \dim_q \ze^i\Zig_\ell\ze^j.
\end{equation}

By definition, $\chi^i_k=0$ unless $k\in\{i,i+1\}$, and $\chi^j_k=0$  unless $k\in \{j,j+1\}$. In particular, comparing with the formula for $\dim_q \ze^i\Zig_\ell\ze^j$ from (\ref{EZigBiWt}), we may assume that $|i-j|\leq 1$. We now proceed in cases.

\vspace{2mm}
\noindent
{\sf Case 1:} $i=j=0$. Then 
\begin{align*}
\sum_{k\in I}
(1+q^2)^{\de_{k,0}}\,\chi^i_k\, \chi^j_k=\,
&(1+q^2)(\chi^0_0)^2+(\chi^0_1)^2,
\\
=\,&(1+q^2) \big(q^{-1}(q^2+q^{-2})^{\ell-1}(q^2+1)\big)^2
\\&+
\big(q(q^2+q^{-2})^{\ell-1}(q^2+1)\big)^2
\\
=\,&m_{0,0}(q)(1+q^2+q^4)
\\
=\,
&m_{0,0}(q)\,\dim_q \ze^0\Zig_\ell\ze^0,
\end{align*}
where we have used (\ref{EZigBiWt}) for the last equality.

\vspace{2mm}
\noindent
{\sf Case 2:} $1\leq i=j\leq \ell-1$. Then 
\begin{align*}
\sum_{k\in I}
(1+q^2)^{\de_{k,0}}\,\chi^i_k\, \chi^j_k=\,
&(\chi^i_i)^2+(\chi^i_{i+1})^2,
\\
=\,& \big(q^{-1}(q^2+q^{-2})^{\ell-i-1}(q^2+1)\big)^2
\\&+
\big(q(q^2+q^{-2})^{\ell-i-1}(q^2+1)\big)^2
\\
=\,&m_{i,i}(q)(1+q^4) 
\\=\,&m_{i,i}(q)\,\dim_q \ze^i\Zig_\ell\ze^i,
\end{align*}
where we have used (\ref{EZigBiWt}) for the last equality.

\vspace{2mm}
\noindent
{\sf Case 3:} $i=j-1$. Then 
\begin{align*}
\sum_{k\in I}
(1+q^2)^{\de_{k,0}}\,\chi^i_k\, \chi^j_k=\,
&
\chi^i_{i+1}\chi^{i+1}_{i+1}
\\
=\,& 
q(q^2+q^{-2})^{\ell-i-2}(q^2+1)\cdot q^{-1}(q^2+q^{-2})^{\ell-i-1}(q^2+1)
\\
=\,&m_{i,i+1}(q) 
\\=\,&m_{i,i+1}(q)\,\dim_q \ze^i\Zig_\ell\ze^{i+1},
\end{align*}
where we have used (\ref{EZigBiWt}) for the last equality.

\vspace{2mm}
\noindent
{\sf Case 4:} $i=j+1$. Then 
\begin{align*}
\sum_{k\in I}
(1+q^2)^{\de_{k,0}}\,\chi^i_k\, \chi^j_k=\,
&
\chi^{j+1}_{j+1}\chi^{j}_{j+1}
\\
=\,& 
q^{-1}(q^2+q^{-2})^{\ell-j-2}(q^2+1)\cdot q(q^2+q^{-2})^{\ell-j-1}(q^2+1)
\\
=\,&m_{j+1,j}(q)\cdot q^4
\\=\,&m_{j+1,j}(q)\,\dim_q \ze^{j+1}\Zig_\ell\ze^{j},
\end{align*}
where we have used (\ref{EZigBiWt}) for the last equality.
\end{proof}

\begin{Proposition}\label{PDimd=1Level N}
For all $i,j\in J$, we have 
$$\dim_q \ga^iY_{\rho,1}^{N\La_0}\ga^j=(1+q^4+\dots+q^{4(N-1)})\dim_q \ze^i\Zig_\ell\ze^j.$$
In particular, $\dim_q Y_{\rho,1}^{N\La_0}=(1+q^4+\dots+q^{4(N-1)})\dim_q \Zig_\ell.$
\end{Proposition}

\begin{proof}
Recall that $r$ denotes $|\rho|$. In addition, we abbreviate 
$$R:=R^{N\Lambda_0}_{N\cont(\rho) + \de},\ Q:=R^{N\Lambda_0}_{N\cont(\rho)},\ \al:=N\cont(\rho).
$$  
For any $i\in I$, we have the following elements of $R_{\al + \de}$, and hence we can also consider them as elements of $R$:
\begin{align*}
\ttf_i:=\iota_{\al,\de}(e_{\al}\otimes e(\bg^i))
\quad\text{and}\quad
\hat \ttf_i:=\iota_{\al,\de}(e_{\al}\otimes e(\hat\bg^i)).
\end{align*}
Note that
$
\ga^iY_{\rho,1}^{N\La_0}\ga^j = \ga^i\Cent_{\rho,1}^{N\La_0}\ga^j
$
and that, by Lemma~\ref{L270116New},
$$
Q \otimes \ga^i\Cent_{\rho,1}^{N\La_0}\ga^j \cong \ttf_i R\ttf_j.
$$ 
It will therefore suffice to prove that for all $i,j\in J$, we have
\begin{align}\label{algn:sub_result}
\dim_q \ttf_iR \ttf_j
=(\dim_q Q)\cdot(\dim_q \ze^i\Zig_\ell\ze^j)\cdot \sum_{s=0}^{N-1}q^{4s},
\end{align}

Recalling the notation of Lemma~\ref{LFact}, we denote 
\begin{align*}
m_{i,j}(q)&:=\bg^i!q^{\lan\bg^i\ran}\bg^j!q^{-\lan\bg^j\ran}=
(1 + q^2)(1 + q^4)^{\ell-i-1}(1 + q^{-2})(1 + q^{-4})^{\ell-j-1}.
\end{align*}
Then by Lemma~\ref{LFact}, we have 
\begin{equation}\label{E300821}
\dim_q \ttf_iR \ttf_j = \sum_{\bi,\bj\in I^{\al}} \dim_q e(\bi\bg^i) R e(\bj\bg^j)
= \sum_{\bi,\bj\in I^{\al}} \frac{\dim_q e(\bi\hat\bg^i) R e(\bj\hat\bg^j)}{m_{i,j}(q)}.
\end{equation}
On the other hand, using Theorem~\ref{TDim}, for fixed $\bi,\bj\in I^{\al}$, we have
\begin{align*}
\dim_q e(\bi\hat\bg^i) R e(\bj\hat\bg^i)
&=  \sum_{ \bla}\sum_{\substack{\Stab\in \Std_p(\bla,\bi\hat\bg^i) \\ \T \in \Std_p(\bla,\bj\hat\bg^j)}}  \deg(\Stab) \deg(\T)\| \bla\|
\\
&=  \sum_{ \bla}\| \bla\|
\Bigg[\sum_{\Stab \in \Std_p(\bla,\bi\hat\bg^i)}  \deg(\Stab)\Bigg] 
\Bigg[\sum_{\T \in \Std_p(\bla,\bj\hat\bg^j) }  \deg(\T)\Bigg],
\end{align*}
where the first summations are over all $\bla \in \Par_p^N(Nr+p)$. 
But if $\Std_p(\bla,\bi\hat\bg^i)\neq\varnothing$, then by Lemma~\ref{LMichael}, we have 
 $\Stab_{\leq Nr} \in \Std_p(\rho^N)$ and 
 $\bla\in \Par_p^N(\rho^N,1)$. 
So $\dim_q e(\bi\hat\bg^i) R e(\bj\hat\bg^i)$ equals
\begin{align*}
&\sum_{\substack{\U\in\Std_p(\rho^N,\bi) \\ \V\in\Std_p(\rho^N,\bj)}}\sum_{ \bla \in \Par_p^N(\rho^N,1)}
\hspace{-1mm}\| \bla\|
\Bigg[\sum_{\Stab \in \Std^{\U,i}_p(\bla)}  \deg(\Stab)\Bigg] \Bigg[\sum_{\T \in \Std^{\V,j}_p(\bla)}  \deg(\T)\Bigg]
\\
=&\hspace{-1mm}\sum_{\substack{\U\in\Std_p(\rho^N,\bi) \\ \V\in\Std_p(\rho^N,\bj)}}
\hspace{-1.5mm}
\deg(\U)\deg(\V)
\sum_{ \bla \in \Par_p^N(\rho^N,1)}\hspace{-1.5mm}
\| \bla\|
\Bigg[\sum_{\Stab \in\Std^{\U,i}_p(\bla)}  \frac{\deg(\Stab)}{\deg(\U)}\Bigg] \Bigg[\sum_{\T \in \Std^{\V,j}_p(\bla)}  \frac{\deg(\T)}{\deg(\V)}\Bigg]
\\
=&\hspace{-1mm}\sum_{\substack{\U\in\Std_p(\rho^N,\bi) \\ \V\in\Std_p(\rho^N,\bj)}}
\hspace{-1mm}\deg(\U)\deg(\V)
\bigg[m_{i,j}(q)\,(\dim_q \ze^i\Zig_\ell\ze^j)\, \sum_{s=0}^{N-1}q^{4s}\bigg],
\end{align*}
where we have used Lemma~\ref{L010921} for the last equality.
In view of (\ref{E300821}), the required equality (\ref{algn:sub_result}) follows from
$$
\sum_{\substack{\U,\V\in\Std_p(\rho^N) \\ \text{with } \bi^\U=\bi,\,\bi^\V=\bj}}\hspace{-5mm}\deg(\U)\deg(\V)=\dim_qQ,
$$
which in turn follows from Theorem~\ref{TDim} upon observing that $\|\rho^N\|=1$ and $\{\rho^N\}=\{\bmu\in\Par_p^N\mid \cont(\bmu)=\al\}$ by Lemma~\ref{LRhoNUnique}. 
\end{proof}

\subsection{Setting up the isomorphism $f$}
\label{SSSettingd=1}
Let $\Zig_\ell$ be the Brauer tree algebra  introduced in \S\ref{SSZig} and $\Zig_\ell[\zz]=H_1(\Zig_\ell)$ be the rank $1$ affine Brauer tree algebra introduced in \S\ref{SSAff}. Our goal is to prove that there is an isomorphisms of graded superalgebras $Y_{\rho,1}^{\La_0}\cong \Zig_\ell$ and $B\cong \Zig_\ell[\zz]$.
Recall that $\bar R_\de$ is a quotient of $R_\de$ which is generated by the elements $\psi_r,y_s,e(\bi)$.  We use the same notation for the corresponding elements of $\bar R_\de$, so we have $y_1,\dots,y_p,\psi_1,\dots,\psi_{p-1}\in \bar R_\de$. 
To construct an isomorphism $f: \Zig_\ell[\zz]\to B$, we define the elements 
$$e^i,u,a^{i,j},c^i,z^i\in B$$ which will represent $f(\ze^i),f(\zu),f(\za^{i,j}),f(\zc^i),f(\ze^i \zz)$, respectively:
\begin{align}
\label{EEIFormula}
e^i&:=\ga^i\qquad (i\in J),
\index{e@$e^i$}
\\
\label{EUFormula}
u&:=e^0y_pe^0
\index{u@$u$}
\\
\label{EZFormula}
z^i&:=e^iy_{p-i-1}y_{p-i}e^i\qquad (i\in J),
\index{z@$z^i$}
\\
\label{ECFormula}
c^i&:=
\left\{
\begin{array}{ll}
(-1)^i(z^i+e^iy_pe^i) &\hbox{if $i=1,\dots,\ell-1$,}\\
e^0(y_{p-1}^2+y_p^2)e^0 &\hbox{if $i=0$,}
\end{array}
\right.
\index{c@$c^i$}
\\
\label{EAII-1Formula}
a^{i,i-1}&:=e^i\psi_{p-1}\cdots\psi_{p-2i}e^{i-1}\qquad (i=1,\dots,\ell-1),
\index{a@$a^{j,k}$}
\\
\label{EAI-1IFormula}
a^{i-1,i}&:=(-1)^ie^{i-1}\psi_{p-2i-1}\cdots\psi_{p-1}e^{i}\qquad (i=1,\dots,\ell-1).
\end{align}
We also set 
$$
z=\sum_{i\in J}z^i\qquad\text{and}\qquad c=\sum_{i\in J}c^i.
\index{z@$z$}\index{c@$c$}
$$
Recalling Khovanov-Lauda diagrams and (\ref{EDivDiag}), we have (colors are for convenience only):
\vspace{2mm}
\begin{align*}
e^i&=\begin{braid}\tikzset{baseline=-0.3em}
	\draw (0,0) node{\color{blue}\normalsize$\ell$\color{black}};
	\braidbox{1}{4.2}{-0.8}{0.7}{$(\ell-1)^2$};
	\draw[dots] (5,0)--(6.5,0);
	\braidbox{7}{10.2}{-0.8}{0.7}{$(i+1)^2$};
	\draw (11,0) node{\normalsize$i$};
	\draw[dots] (11.6,0)--(13.3,0);
	\draw (13.6,0) node{\normalsize$1$};
	\redbraidbox{14.5}{16.6}{-0.8}{0.7}{$\color{red}0\,\,0\color{black}$};
	\draw (17.5,0) node{\normalsize$1$};
	\draw[dots] (18.1,0)--(19.5,0);
	\draw (20,0) node{\normalsize$i$};
	\draw (20,-1) node{};
\end{braid},
\\
u&=\begin{braid}\tikzset{baseline=1.3em}
	\draw (0,0) node{\color{blue}\normalsize$\ell$\color{black}};
	\braidbox{1}{4.2}{-0.8}{0.7}{$(\ell-1)^2$};
	\draw[dots] (5,0)--(6.5,0);
	\braidbox{7}{9.1}{-0.8}{0.7}{$1\,\,1$};
	\redbraidbox{10}{12.1}{-0.8}{0.7}{$\color{red}0\,\,0\color{black}$};
	\draw (0,4) node{\color{blue}\normalsize$\ell$\color{black}};
	\braidbox{1}{4.2}{3.2}{4.7}{$(\ell-1)^2$};
	\draw[dots] (5,4)--(6.5,4);
	\braidbox{7}{9.1}{3.2}{4.7}{$1\,\,1$};
	\redbraidbox{10}{12.1}{3.2}{4.7}{$\color{red}0\,\,0\color{black}$};
	 \draw[blue](0,3)--(0,1);
	 \draw(1.5,3)--(1.5,1);
	 \draw(3.5,3)--(3.5,1);
	 \draw(7.6,3)--(7.6,1);
	 \draw(8.5,3)--(8.5,1);
	  \draw[red](10.6,3)--(10.6,1);
	  \draw[red](11.5,3)--(11.5,1);
	  \reddot (11.5,2);
	  \draw (10,-1) node{};
	  \end{braid},
\\
z^i&=\begin{braid}\tikzset{baseline=1.3em}
	\draw (0,0) node{\color{blue}\normalsize$\ell$\color{black}};
	\braidbox{1}{4.2}{-0.8}{0.7}{$(\ell-1)^2$};
	\draw[dots] (5,0)--(6.5,0);
	\braidbox{7}{10.2}{-0.8}{0.7}{$(i+1)^2$};
	\draw (11,0) node{\normalsize$i$};
	\draw[dots] (11.6,0)--(13.3,0);
	\draw (13.6,0) node{\normalsize$1$};
	\redbraidbox{14.5}{16.6}{-0.8}{0.7}{$\color{red}0\,\,0\color{black}$};
	\draw (17.5,0) node{\normalsize$1$};
	\draw[dots] (18.1,0)--(19.5,0);
	\draw (20,0) node{\normalsize$i$};
	\draw (0,4) node{\color{blue}\normalsize$\ell$\color{black}};
	\braidbox{1}{4.2}{3.2}{4.7}{$(\ell-1)^2$};
	\draw[dots] (5,4)--(6.5,4);
	\braidbox{7}{10.2}{3.2}{4.7}{$(i+1)^2$};
	\draw (11,4) node{\normalsize$i$};
	\draw[dots] (11.6,4)--(13.3,4);
	\draw (13.6,4) node{\normalsize$1$};
	\redbraidbox{14.5}{16.6}{3.2}{4.7}{$\color{red}0\,\,0\color{black}$};
	\draw (17.5,4) node{\normalsize$1$};
	\draw[dots] (18.1,4)--(19.5,4);
	\draw (20,4) node{\normalsize$i$};
	\draw[blue](0,3)--(0,1);
	 \draw(1.5,3)--(1.5,1);
	 \draw(3.5,3)--(3.5,1);
	  \draw(7.5,3)--(7.5,1);
	  \draw(9.5,3)--(9.5,1);
	  \draw(11,3)--(11,1);
	  \draw(13.6,3)--(13.6,1);
	   \draw[red](15.1,3)--(15.1,1);
	   \draw[red](16,3)--(16,1);
	   \draw(17.5,3)--(17.5,1);
\draw(20,3)--(20,1);
 \reddot (15.1,2.2);
 \reddot (16,1.8);
  \draw (10,-1) node{};
\end{braid},
\\
c^i&=(-1)^iz^i+(-1)^i\begin{braid}\tikzset{baseline=1.3em}
	\draw (0,0) node{\color{blue}\normalsize$\ell$\color{black}};
	\braidbox{1}{4.2}{-0.8}{0.7}{$(\ell-1)^2$};
	\draw[dots] (5,0)--(6.5,0);
	\braidbox{7}{10.2}{-0.8}{0.7}{$(i+1)^2$};
	\draw (11,0) node{\normalsize$i$};
	\draw[dots] (11.6,0)--(13.3,0);
	\draw (13.6,0) node{\normalsize$1$};
	\redbraidbox{14.5}{16.6}{-0.8}{0.7}{$\color{red}0\,\,0\color{black}$};
	\draw (17.5,0) node{\normalsize$1$};
	\draw[dots] (18.1,0)--(19.5,0);
	\draw (20,0) node{\normalsize$i$};
	\draw (0,4) node{\color{blue}\normalsize$\ell$\color{black}};
	\braidbox{1}{4.2}{3.2}{4.7}{$(\ell-1)^2$};
	\draw[dots] (5,4)--(6.5,4);
	\braidbox{7}{10.2}{3.2}{4.7}{$(i+1)^2$};
	\draw (11,4) node{\normalsize$i$};
	\draw[dots] (11.6,4)--(13.3,4);
	\draw (13.6,4) node{\normalsize$1$};
	\redbraidbox{14.5}{16.6}{3.2}{4.7}{$\color{red}0\,\,0\color{black}$};
	\draw (17.5,4) node{\normalsize$1$};
	\draw[dots] (18.1,4)--(19.5,4);
	\draw (20,4) node{\normalsize$i$};
	\draw[blue](0,3)--(0,1);
	 \draw(1.5,3)--(1.5,1);
	 \draw(3.5,3)--(3.5,1);
	  \draw(7.5,3)--(7.5,1);
	  \draw(9.5,3)--(9.5,1);
	  \draw(11,3)--(11,1);
	  \draw(13.6,3)--(13.6,1);
	   \draw[red](15.1,3)--(15.1,1);
	   \draw[red](16,3)--(16,1);
	   \draw(17.5,3)--(17.5,1);
\draw(20,3)--(20,1);
 \blackdot (20,2);
 \draw (10,-1) node{};
\end{braid}\qquad(i\neq 0),
\\
c^0&=\begin{braid}\tikzset{baseline=1.3em}
	\draw (0,0) node{\color{blue}\normalsize$\ell$\color{black}};
	\braidbox{1}{4.2}{-0.8}{0.7}{$(\ell-1)^2$};
	\draw[dots] (5,0)--(6.5,0);
	\braidbox{7}{9.1}{-0.8}{0.7}{$1\,\,1$};
	\redbraidbox{10}{12.1}{-0.8}{0.7}{$\color{red}0\,\,0\color{black}$};
	\draw (0,4) node{\color{blue}\normalsize$\ell$\color{black}};
	\braidbox{1}{4.2}{3.2}{4.7}{$(\ell-1)^2$};
	\draw[dots] (5,4)--(6.5,4);
	\braidbox{7}{9.1}{3.2}{4.7}{$1\,\,1$};
	\redbraidbox{10}{12.1}{3.2}{4.7}{$\color{red}0\,\,0\color{black}$};
	 \draw[blue](0,3)--(0,1);
	 \draw(1.5,3)--(1.5,1);
	 \draw(3.5,3)--(3.5,1);
	 \draw(7.6,3)--(7.6,1);
	 \draw(8.5,3)--(8.5,1);
	  \draw[red](10.6,3)--(10.6,1);
	  \draw[red](11.5,3)--(11.5,1);
	  \reddot (10.6,2.2);
	  \reddot (10.6,1.8);
	  \end{braid}
	  +
	  \begin{braid}\tikzset{baseline=1.3em}
	\draw (0,0) node{\color{blue}\normalsize$\ell$\color{black}};
	\braidbox{1}{4.2}{-0.8}{0.7}{$(\ell-1)^2$};
	\draw[dots] (5,0)--(6.5,0);
	\braidbox{7}{9.1}{-0.8}{0.7}{$1\,\,1$};
	\redbraidbox{10}{12.1}{-0.8}{0.7}{$\color{red}0\,\,0\color{black}$};
	\draw (0,4) node{\color{blue}\normalsize$\ell$\color{black}};
	\braidbox{1}{4.2}{3.2}{4.7}{$(\ell-1)^2$};
	\draw[dots] (5,4)--(6.5,4);
	\braidbox{7}{9.1}{3.2}{4.7}{$1\,\,1$};
	\redbraidbox{10}{12.1}{3.2}{4.7}{$\color{red}0\,\,0\color{black}$};
	 \draw[blue](0,3)--(0,1);
	 \draw(1.5,3)--(1.5,1);
	 \draw(3.5,3)--(3.5,1);
	 \draw(7.6,3)--(7.6,1);
	 \draw(8.5,3)--(8.5,1);
	  \draw[red](10.6,3)--(10.6,1);
	  \draw[red](11.5,3)--(11.5,1);
	  \reddot (11.5,2.2);
	  \reddot (11.5,1.8);
	  \draw (10,-1) node{};
	  \end{braid},
	  \\
a^{i,i-1}&=\begin{braid}\tikzset{baseline=1.3em}
	\draw (0,0) node{\color{blue}\normalsize$\ell$\color{black}};
	\braidbox{1}{4.2}{-0.8}{0.7}{$(\ell-1)^2$};
	\draw[dots] (5,0)--(6.5,0);
	\braidbox{7}{10.2}{-0.8}{0.7}{$(i+1)^2$};
	\braidbox{11}{13}{-0.8}{0.7}{$i\,\,i$};
	\draw (14.8,0) node{\normalsize$i-1$};
		\draw[dots] (16.6,0)--(18.3,0);
	\draw (18.6,0) node{\normalsize$1$};
	\redbraidbox{19.5}{21.6}{-0.8}{0.7}{$\color{red}0\,\,0\color{black}$};
	\draw (22.5,0) node{\normalsize$1$};
	\draw[dots] (23.1,0)--(24.5,0);
	\draw (26,0) node{\normalsize$i-1$};
	\draw (0,4) node{\color{blue}\normalsize$\ell$\color{black}};
	\braidbox{1}{4.2}{3.2}{4.7}{$(\ell-1)^2$};
	\draw[dots] (5,4)--(6.5,4);
	\braidbox{7}{10.2}{3.2}{4.7}{$(i+1)^2$};
	\draw (11.4,4) node{\normalsize$i$};
	\draw (13.3,4) node{\normalsize$i-1$};
	\draw[dots] (14.9,4)--(16.6,4);
	\draw (16.9,4) node{\normalsize$1$};
	\redbraidbox{17.8}{19.9}{3.2}{4.7}{$\color{red}0\,\,0\color{black}$};
	\draw (20.8,4) node{\normalsize$1$};
	\draw[dots] (21.6,4)--(23,4);
	\draw (24.7,4) node{\normalsize$i-1$};
	\draw (26.8,4) node{\normalsize$i$};
	\draw[blue](0,3)--(0,1);
	 \draw(1.5,3)--(1.5,1);
	 \draw(3.5,3)--(3.5,1);
	  \draw(7.5,3)--(7.5,1);
	  \draw(9.5,3)--(9.5,1);
	  \draw(11.4,3)--(11.4,1);
	  \draw(13.3,3)--(14.6,1);
	  \draw(16.9,3)--(18.4,1);
	   \draw[red](18.4,3)--(19.9,1);
	   \draw[red](19.3,3)--(20.8,1);
\draw(20.9,3)--(22.4,1);
\draw(24.5,3)--(25.9,1);
\draw(26.7,3)--(12.5,1);
\draw (20,-1) node{};
\end{braid},
\\
a^{i-1,i}&=(-1)^i\begin{braid}\tikzset{baseline=1.3em}
	\draw (0,4) node{\color{blue}\normalsize$\ell$\color{black}};
	\braidbox{1}{4.2}{3.2}{4.7}{$(\ell-1)^2$};
	\draw[dots] (5,4)--(6.5,4);
	\braidbox{7}{10.2}{3.2}{4.7}{$(i+1)^2$};
	\braidbox{11}{13}{3.2}{4.7}{$i\,\,i$};
	\draw (14.8,4) node{\normalsize$i-1$};
		\draw[dots] (16.6,4)--(18.3,4);
	\draw (18.6,4) node{\normalsize$1$};
	\redbraidbox{19.5}{21.6}{3.2}{4.7}{$\color{red}0\,\,0\color{black}$};
	\draw (22.5,4) node{\normalsize$1$};
	\draw[dots] (23.1,4)--(24.5,4);
	\draw (26,4) node{\normalsize$i-1$};
	\draw (0,0) node{\color{blue}\normalsize$\ell$\color{black}};
	\braidbox{1}{4.2}{-0.8}{0.7}{$(\ell-1)^2$};
	\draw[dots] (5,0)--(6.5,0);
	\braidbox{7}{10.2}{-0.8}{0.7}{$(i+1)^2$};
	\draw (11.4,0) node{\normalsize$i$};
	\draw (13.3,0) node{\normalsize$i-1$};
	\draw[dots] (14.9,0)--(16.6,0);
	\draw (16.9,0) node{\normalsize$1$};
	\redbraidbox{17.8}{19.9}{-0.8}{0.7}{$\color{red}0\,\,0\color{black}$};
	\draw (20.8,0) node{\normalsize$1$};
	\draw[dots] (21.6,0)--(23,0);
	\draw (24.7,0) node{\normalsize$i-1$};
	\draw (26.8,0) node{\normalsize$i$};
	\draw[blue](0,1)--(0,3);
	 \draw(1.5,1)--(1.5,3);
	 \draw(3.5,1)--(3.5,3);
	  \draw(7.5,1)--(7.5,3);
	  \draw(9.5,1)--(9.5,3);
	  \draw(11.5,1)--(12.3,3);
	  \draw(13.3,1)--(14.6,3);
	  \draw(16.9,1)--(18.4,3);
	   \draw[red](18.4,1)--(19.9,3);
	   \draw[red](19.3,1)--(20.8,3);
\draw(20.9,1)--(22.4,3);
\draw(24.5,1)--(25.9,3);
\draw(26.7,1)--(11.6,3);
\end{braid}.
\end{align*}

\vspace{2mm}
\begin{Lemma} \label{EAOneSide} 
Let $i\in\{1,\dots,\ell-1\}$. Then 
$$
a^{i,i-1}=\psi_{p-1}\cdots\psi_{p-2i}e^{i-1}\quad\text{and}\quad
a^{i-1,i}=(-1)^ie^{i-1}\psi_{p-2i-1}\cdots\psi_{p-1}.$$
\end{Lemma}
\begin{proof}
Note using (\ref{R5}) and (\ref{R7}) that
$$
\psi_{p-1}\cdots\psi_{p-2i}e^{i-1}=e^i\psi_{p-1}\cdots\psi_{p-2i}e^{i-1}=a^{i,i-1}.
$$
The second equality is proved similarly. 
\end{proof}

\begin{Lemma} \label{LPsiE} 
Let $i\in J$. Then $\psi_re^i=0$ in $\bar R_\de$ unless $r=p-1$ or $r=2t+1$ for $1\leq t\leq \ell-1-i$. 
\end{Lemma}
\begin{proof}
Recall the notation(\ref{EHatGI}). 
Then for $1\leq t\leq\ell-1-i$, we have $g_{2t}^i=g_{2t+1}^i$ , so $\psi_{2t}e^i=0$ by the relation (\ref{R6}) with $i_r=i_{r+1}$. Similarly $\psi_{p-i-1}e^i=0$. 

Suppose now that $r=p-i-1-j$ for some $1\leq j\leq i$. Then 
$\psi_{r}e^i=e(s_{r}\cdot\hat\bg^i)\psi_{p-i-1-j}e^i$, and the word 
$$s_{r}\cdot\hat\bg^i=\ell(\ell-1)^2\cdots(i+1)^2i\cdots (j+1)(j-1)j(j-2)\cdots10^21\cdots i,$$
is not cuspidal by Example~\ref{ExNonCusp}(a), so $\psi_{r}e^i=0$. 

Finally, suppose that $r=p-i+j$ for some $1\leq j\leq i-2$. Then 
$\psi_{r}e^i=e(s_{r}\cdot\hat\bg^i)\psi_{r}e^i$, and the word 
$$s_{r}\cdot\hat\bg^i=\ell(\ell-1)^2\cdots(i+1)^2i\cdots 10^21\cdots(j-1)(j+1)j(j+2)\cdots i,$$
is not cuspidal by Example~\ref{ExNonCusp}(b), so $\psi_{r}e^i=0$. 
\end{proof}

\subsection{Spanning set for $B$}
Recall the notation (\ref{EHatGI}).

\begin{Lemma} \label{LIBI}
Let $i\in J$. Then $e^iBe^i=e^i\F[y_1,\dots y_p]e^i$. 
\end{Lemma}
\begin{proof}
Let $H:=\{w\in {\Si}_p\mid  w\cdot \hat \bg^i=\hat \bg^i\}$. 
By Theorem~\ref{TBasis}, we have 
$$e^iBe^i=\sum_{w\in H} e^i\F[y_1,\dots y_p]\psi_we^i.
$$  
So it suffices to show that for $w\in H$ we have $\psi_we^i=0$ in $\bar R_\de$ unless $w=1$. Note that, while $\psi_w$ depends on the choice of the reduced expression for $w$, we can chose any such reduced expression. 

Let $\psi_we^i\neq 0$ in $\bar R_\de$. Since $g_1^i=\ell$ and $\ell$ appears in $\hat \bg^i$ only once, we must have $w(1)=1$. Suppose  $g_2^i=g_3^i=\ell-1$ (i.e. $i\leq \ell-2$). Since $\ell-1$ appears in $\hat \bg^i$ only  twice, we must have either $w(1)=1,w(2)=2$ or $w(1)=2,w(2)=w(1)$. In the latter case we can choose the reduced expression for $w$ so that $\psi_w=\psi_{w'}\psi_2$, and since $\psi_2e^i=0$ by Lemma~\ref{LPsiE}, we have $\psi_we^i=0$. So $w(2)=2,w(3)=3$. Similarly, we prove that $w(r)=r$ for all $r=2,\dots,p-2i-2$ and for $r=p-i-1,p-i$. 
It remains to prove that 
\begin{equation}\label{070421}
w(p-i-1-k)=p-i-1-k\quad \text{and}\quad w(p-i+k)=p-i+k
\end{equation} 
for all $k=1,\dots,i$. Choose the smallest $k$ for which (\ref{070421}) fails. Since $g_{p-i-1-k}^i=k=g_{p-i+k}^i$, we then have 
$$w(p-i-1-k)=p-i+k\quad\text{and}\quad w(p-i+k)=p-i-1-k.$$ 
Since $w(p-i-k)=p-i-k$ by the choice of $k$, we can choose a reduced expression for $w$ so that $\psi_w=\psi_{w'}\psi_{p-i-1-k}$. But 
$\psi_{p-i-1-k}e^i=0$ by Lemma~\ref{LPsiE}, hence $\psi_we^i=0$.
\end{proof}

\begin{Lemma} \label{LIBI-1}
Let $i\in\{1,\dots,\ell-1\}$. Then $e^iBe^{i-1}=e^i\F[y_1,\dots y_p]a^{i,i-1}.$ 
\end{Lemma}
\begin{proof}
Let $H:=\{w\in {\Si}_p\mid w\hat\bg^{i-1}=\hat\bg^{i}$. 
By Theorem~\ref{TBasis}, 
$$e^iBe^{i-1}\sum_{w\in H}e^i\F[y_1,\dots y_p]\psi_we^{i-1}.$$ So it suffices to show that for any $w\in H$ we have $\psi_we^{i-1}=0$ in $\bar R_\de$ unless $w=s_{p-1}\cdots s_{p-2i}$. Note that, while $\psi_w$ depends on the choice of the reduced expression for $w$, we can chose any such reduced expression. The argument now is similar to that of Lemma~\ref{LIBI}. We use Lemma~\ref{LPsiE} to first show that 
$w(r)=r$ for all $r=1,\dots,p-2i-1$ and $w(p-i)=p-i-1,w(p-i+1)=p-i$, and then using Lemma~\ref{LPsiE} again to deduce $w(r)=r-1$ for the remaining $r$.
\end{proof}

\begin{Lemma} \label{LI-1BI}
Let $i\in\{1,\dots,\ell-1\}$. Then $e^{i-1}Be^{i}=a^{i-1,i}\F[y_1,\dots y_p]e^i.$ 
\end{Lemma}
\begin{proof}
By Theorem~\ref{TBasis}, $e^{i-1}Be^{i}$ is the sum of  $e^{i-1}\F[y_1,\dots y_p]\psi_we^{i}$ with $w\in {\Si}_p$ such that $w\hat\bg^{i}=\hat\bg^{i-1}$. Using Lemma~\ref{LPsiE} as in the proof of Lemma~\ref{LIBI}, we can show that one of the following two cases happens:
\begin{enumerate}
\item[{\rm (1)}] $\psi_w=\psi_{p-2i-1}\cdots\psi_{p-1}$;
\item[{\rm (2)}] $\psi_w=\psi_{p-2i}\cdots\psi_{p-1}$.
\end{enumerate}
Note using the relations (\ref{R4}), (\ref{R5}) that 
\begin{align*}
&e^{i-1}\F[y_1,\dots,y_p]\psi_{p-2i-1}\cdots\psi_{p-1}e^i
+e^{i-1}\F[y_1,\dots,y_p]\psi_{p-2i}\cdots\psi_{p-1}e^i
\\
=\,&e^{i-1}\psi_{p-2i-1}\cdots\psi_{p-1}\F[y_1,\dots,y_p]e^i
+e^{i-1}\psi_{p-2i}\cdots\psi_{p-1}\F[y_1,\dots,y_p]e^i.
\end{align*}
Note that $e^{i-1}=e^{i-1}\psi_{p-2i-1}y_{p-2i}$, so the second summand equals
\begin{align*}
&e^{i-1}\psi_{p-2i-1}y_{p-2i}\psi_{p-2i}\cdots\psi_{p-1}\F[y_1,\dots,y_p]e^i
\\
=
\,
&e^{i-1}\psi_{p-2i-1}\psi_{p-2i}\cdots\psi_{p-1}y_p\F[y_1,\dots,y_p]e^i,
\end{align*}
which lies in the first summand. On the other hand, by Lemma~\ref{EAOneSide}, the first summand equals
$
a^{i-1,i}\F[y_1,\dots,y_p]e^i,
$
which yields the claim. \end{proof}

\begin{Lemma} \label{LIBJ}
Let $i,j\in J$ be such that $|i-j|>1$. Then $e^jBe^{i}=0.$ 
\end{Lemma}
\begin{proof}
We explain the case $i<j-1$, the argument for $i>j+1$ is similar. 
Let $H:=\{w\in {\Si}_p\mid w\hat\bg^{i}=\hat\bg^{j}$. 
By Theorem~\ref{TBasis}, 
$$e^jBe^{i}=\sum_{w\in H}e^j\F[y_1,\dots y_p]\psi_we^{i}.$$ 
It suffices to show that for any $w\in H$ we have $\psi_we^i=0$ in $\bar R_\de$. Note that, while $\psi_w$ depends on the choice of the reduced expression for $w$, we can chose any such reduced expression. Using Lemma~\ref{LPsiE}, we deduce that 
\begin{align*}
&w(2\ell-2j)= 2\ell-2j, 
&w(2\ell-2j+1)=p,
\\ 
&w(2\ell-2j+2)= 2\ell-2j+1,
&w(2\ell-2j+3)=p-1.
\end{align*}
So we may assume that $\psi_w=\psi_u\psi_{2\ell-2j+2}\psi_{2\ell-2j+1}$ for some $u\in {\Si}_p$. We deduce that $\psi_we^j=0$ since the word 
$$
s_{2\ell-2j+2}s_{2\ell-2j+1}\cdot\hat\bg^j=\ell(\ell-1)^2\cdots (j+1)^2j(j-1)^2\cdots
$$
is non-cuspidal by Example~\ref{ExNonCusp}(c). 
\end{proof}

Recalling the notation (\ref{EHatGI}), we set for all $i,j\in J$:
\begin{equation}\label{ERS}
r^i_j=\min\{r\mid g^i_r=j\}\quad\text{and}\quad s^i_j=\max\{r\mid g^i_r=j\}.\index{r@$r^i_j$}\index{s@$s^i_j$}
\end{equation}
For example, $r^i_0=p-i-1$ and $s^i_0=p-i$ so 
$$
z^i=e^iy_{r^i_0}y_{s^i_0}e^i\quad\text{and}\quad 
u=e^0y_{s^0_0}e^0=e^0y_{p}e^0=e^0(y_{p-1}+y_{p})e^0=e^0(y_{r^0_0}+y_{s^0_0})e^0
$$
since $e^0y_{p-1}e^0=0$ thanks to Lemma~\ref{LNH3}(ii). By the same lemma, we have $e^iy_{r^i_j}e^i=0$ if $j>i$. 

Fix $i\in J$. For $j>i$, or $j=0\leq i$, we set 
$$
\la^i_j:=e^i(y_{r^i_j}+y_{s^i_j})e^i=e^iy_{s^i_j}e^i\in e^iBe^i\quad\text{and}\quad
\mu^i_j:=e^i(y_{r^i_j}y_{s^i_j})e^i\in e^iBe^i.
\index{l@$\la^i_j$}\index{m@$\mu^i_j$}
$$ 
Note that $\la^0_0=u$ and $\mu^i_0=z^i$. 
For $1\leq j\leq i$, we set 
$$
\xi^i_j:=e^iy_{r^i_j}e^i\in e^iBe^i\quad\text{and}\quad \eta^i_j:=e^iy_{s^i_j}e^i\in e^iBe^i.
\index{x@$\xi^i_j$}\index{$\xi^i_j$}\index{e@$\eta^i_j$}\index{$\eta^i_j$}
$$
Finally, set 
$$\ka^i:=e^iy_1e^i\in e^iBe^i.\index{k@$\ka^i$}$$ 
Note that all the new elements introduced in this paragraph  commute, with one exception: $$\la^i_0\mu^i_0=-\mu^i_0\la^i_0.$$ 

By Lemma~\ref{LNH1}, we can now slightly improve on Lemma~\ref{LIBI}:

\begin{Lemma} \label{LIBIImprove}
Let $i\in J$. Then 
$e^iBe^i$ is generated by the elements 
$$\ka^i,\la^i_0,z^i,\xi^i_1,\eta^i_1,\dots,\xi^i_i,\eta^i_i,\la^i_{i+1},\mu^i_{i+1},\dots,\la^i_{\ell-1},\mu^i_{\ell-1}.
$$
\end{Lemma}

To further improve on Lemma~\ref{LIBI}, we describe some relations between the elements in Lemma~\ref{LIBIImprove}:

\begin{Lemma} \label{LRelB}
Let $i\in J$. In $e^iBe^i$, we have: 
\begin{eqnarray}
&&-\ka^i=\mu^i_{\ell-1}=\dots=\mu^i_{i+1}=
\left\{
\begin{array}{ll}
-(\xi^i_i)^2 &\hbox{if $i\neq 0$,}\\
-u^2z^0-(z^0)^2 &\hbox{if $i=0$,}
\end{array}
\right.
\label{ERelB1}
\\
&&\la^i_{\ell-1}=\dots=\la^i_{i+1}=0,
\label{ERelB2}
\\
&&
\xi^i_i=\dots=\xi^i_1=z^i=-\eta^i_1=\dots=-\eta^i_{i-1},
\label{ERelB3}
\\
&&
\la^i_0=0\qquad(\text{if $i\neq 0$}),
\label{ERelB4}
\\
&&
u^3=0.
\label{ERelB5}
\end{eqnarray} 
\end{Lemma}
\begin{proof}
We will drop the superscripts $i$ to simplify the notation. 

First, we deal with the exceptional case $\ell=1$. In this case $i=0$ and the lemma claims that $\ka=u^2z+z^2$ and $u^3=0$. Using the fact that the word $010\in I^\de$ is not cuspidal by Lemma~\ref{LCuspExpl}, we deduce that in $eBe$, we have $e\psi_1^2e=0=e\psi_2\psi_1^2e$. So, using (\ref{R6}), we get 
\begin{align*}
0=e\psi_2\psi_1^2e=-e\psi_2y_1e+e\psi_2y_2^4e
=-ey_1\psi_2e+e\psi_2y_2^4e
=-u^3,
\\
0=e\psi_1^2e=-ey_1e+ey_2^4e=-\ka+u^2z+z^2,
\end{align*}
where for the last equality in the first line we have used Lemmas~\ref{LPsiE} and \ref{LNH4}(i), and for the last equality in the second line we have used Lemma~\ref{LNH3}(ii).

Now let $\ell>1$. Since the word $s_1\cdot\hat \bg^i$ is not cuspidal by Lemma~\ref{LCuspExpl}, we have using Lemma~\ref{LNHEven}(i) for the last equality:
$$
0=e\psi_1^2e=e(-y_1+y_2^2)e=-ey_1e+ey_2^2e=-\ka-\mu_{\ell-1},
$$
which gives the first equality in (\ref{ERelB1}). Similarly, using Lemma~\ref{LNHEven}(ii), we have 
$$
0=e\psi_2\psi_1^2e=e\psi_2(-y_1+y_2^2)e=e\psi_2y_2^2e=-\la_{\ell-1},
$$
which gives the equality $\la_{\ell-1}=0$ in (\ref{ERelB2}). 

Since the word $s_4s_3\cdot\hat\bg^i=\ell\,(\ell-1)\,(\ell-2)^2\cdots\cdots$ is non-cuspidal by Lemma~\ref{LCuspExpl}, we have, using (\ref{R6}) and (\ref{R5}): 
\begin{align*}
0&=e\psi_2\psi_3\psi_4^2\psi_3e=e\psi_2(y_3^2-y_3y_4-y_3y_5+y_4y_5)e=
\la_{\ell-1}-\la_{\ell-2},
\end{align*}
where we have used that 
$e\psi_2y_3^2e=\la_{\ell-1}$ by Lemma~\ref{LNHEven}(ii), 
$e\psi_2y_3y_4e=0$ by Lemma~\ref{LNHEven}(i), 
$e\psi_2y_4y_5e=0$ by Lemma~\ref{LPsiE}, and $e\psi_2y_4y_5e=0$ by definition. 
This gives the equality $\la_{\ell-1}=\la_{\ell-2}$ in (\ref{ERelB2}). 

Since the word $s_2s_3\cdot\hat\bg^i=\ell\,(\ell-2)\cdots$ is non-cuspidal by Lemma~\ref{LCuspExpl}, we obtain, similarly to the previous paragraph,
\begin{align*}
0&=e\psi_3\psi_2^2\psi_3\psi_3e=
\mu_{\ell-1}-\mu_{\ell-2}.
\end{align*}
This gives the equality $\mu_{\ell-1}=\mu_{\ell-2}$ in (\ref{ERelB1}). 

Continuing like this, we get $0=\la_{\ell-1}=\la_{\ell-2}=\dots=\la_{i+1}$, which is (\ref{ERelB2}), and $\mu_{\ell-1}=\mu_{\ell-2}=\dots=\mu_{i+1}$. 

Now it will be more convenient to use Khovanov-Lauda diagrams instead of formulas. Suppose first that $i>0$. Since a word starting with $\ell\,(\ell-1)^2\,\cdots\,(i+2)^2\,i$ is non-cuspidal by Lemma~\ref{LCuspExpl}, we get 
$$
0=\begin{braid}\tikzset{baseline=1.3em}
	\draw (0,0) node{\color{blue}\normalsize$\ell$\color{black}};
	\braidbox{1}{4.2}{-0.8}{0.7}{$(\ell-1)^2$};
	\draw[dots] (5,0)--(6.5,0);
	\braidbox{7}{10.2}{-0.8}{0.7}{$(i+1)^2$};
	\draw (11,0) node{\normalsize$i$};
	\draw[dots] (11.6,0)--(13.3,0);
	\draw (13.6,0) node{\normalsize$1$};
	\redbraidbox{14.5}{16.6}{-0.8}{0.7}{$\color{red}0\,\,0\color{black}$};
	\draw (17.5,0) node{\normalsize$1$};
	\draw[dots] (18.1,0)--(19.5,0);
	\draw (20,0) node{\normalsize$i$};
	\draw (0,4) node{\color{blue}\normalsize$\ell$\color{black}};
	\braidbox{1}{4.2}{3.2}{4.7}{$(\ell-1)^2$};
	\draw[dots] (5,4)--(6.5,4);
	\braidbox{7}{10.2}{3.2}{4.7}{$(i+1)^2$};
	\draw (11,4) node{\normalsize$i$};
	\draw[dots] (11.6,4)--(13.3,4);
	\draw (13.6,4) node{\normalsize$1$};
	\redbraidbox{14.5}{16.6}{3.2}{4.7}{$\color{red}0\,\,0\color{black}$};
	\draw (17.5,4) node{\normalsize$1$};
	\draw[dots] (18.1,4)--(19.5,4);
	\draw (20,4) node{\normalsize$i$};
	\draw[blue](0,3)--(0,1);
	 \draw(1.5,3)--(1.5,1);
	 \draw(3.5,3)--(3.5,1);
	  \draw(7.5,3)--(7.5,1);
	  \draw(9.5,3)--(9.5,1);
	  \draw(11,3)--(11,2.7)--(7.2,2.4)--(7.2,1.6)--(11,1.3)--(11,1);
	  \draw(13.6,3)--(13.6,1);
	   \draw[red](15.1,3)--(15.1,1);
	   \draw[red](16,3)--(16,1);
	   \draw(17.5,3)--(17.5,1);
\draw(20,3)--(20,1);
\end{braid}=\mu_{i+1}-\xi_i\la_{i+1}+\xi_i^2=\mu_{i+1}+\xi_i^2,
$$
giving the equality $\mu_{i+1}=-\xi_i^2$ and completing the proof of (\ref{ERelB1}) for $i>0$. Let $i=0$. Since a word $\ell\,(\ell-1)^2\,\cdots\,2^2\,0\,1\,1\,0$ is non-cuspidal by Lemma~\ref{LCuspExpl}, we have
$$
0=\begin{braid}\tikzset{baseline=1.3em}
	\draw (0,0) node{\color{blue}\normalsize$\ell$\color{black}};
	\braidbox{1}{4.2}{-0.8}{0.7}{$(\ell-1)^2$};
	\draw[dots] (5,0)--(6.5,0);
	\braidbox{7}{9.1}{-0.8}{0.7}{$1\,\,1$};
	\redbraidbox{10}{12.1}{-0.8}{0.7}{$\color{red}0\,\,0\color{black}$};
	\draw (0,4) node{\color{blue}\normalsize$\ell$\color{black}};
	\braidbox{1}{4.2}{3.2}{4.7}{$(\ell-1)^2$};
	\draw[dots] (5,4)--(6.5,4);
	\braidbox{7}{9.1}{3.2}{4.7}{$1\,\,1$};
	\redbraidbox{10}{12.1}{3.2}{4.7}{$\color{red}0\,\,0\color{black}$};
	 \draw[blue](0,3)--(0,1);
	 \draw(1.5,3)--(1.5,1);
	 \draw(3.5,3)--(3.5,1);
	 \draw(7.6,3)--(7.6,1);
	 \draw(8.5,3)--(8.5,1);
	  \draw[red](10.6,3)--(10.6,2.7)--(7.2,2.4)--(7.2,1.6)--(10.6,1.3)--(10.6,1);
	  \draw[red](11.5,3)--(11.5,1);
	  \end{braid}
= \mu_1-ey_{p-1}^2\la_{1}+ey_{p-1}^4e=\mu_1+u^2z+z^2,
$$
where we have used $\la_{1}=0$ and Lemma~\ref{LNH3}(ii)  for the last equality. This completes the proof of (\ref{ERelB1}).

Note that  (\ref{ERelB3}) is not vacuous only if $i>0$. Since a word starting with $\ell\,(\ell-1)^2\,\cdots\,(i+1)^2\,(i-1)$ is non-cuspidal by Lemma~\ref{LCuspExpl}, we have  
$$
0=\begin{braid}\tikzset{baseline=1.3em}
	\draw (0,4) node{\color{blue}\normalsize$\ell$\color{black}};
	\braidbox{1}{4.2}{3.2}{4.7}{$(\ell-1)^2$};
	\draw[dots] (5,4)--(6.5,4);
	\braidbox{7}{10.2}{3.2}{4.7}{$(i+1)^2$};
	\draw (11.4,4) node{\normalsize$i$};
	\draw (13.3,4) node{\normalsize$i-1$};
	\draw[dots] (14.9,4)--(16.6,4);
	\draw (16.9,4) node{\normalsize$1$};
	\redbraidbox{17.8}{19.9}{3.2}{4.7}{$\color{red}0\,\,0\color{black}$};
	\draw (20.8,4) node{\normalsize$1$};
	\draw[dots] (21.6,4)--(23,4);
	\draw (24.7,4) node{\normalsize$i-1$};
	\draw (26.8,4) node{\normalsize$i$};
	\draw[blue](0,3)--(0,1);
	 \draw(1.5,3)--(1.5,1);
	 \draw(3.5,3)--(3.5,1);
	  \draw(7.5,3)--(7.5,1);
	  \draw(9.5,3)--(9.5,1);
	  \draw(11.4,3)--(13.3,2)--(11.4,1);
	  \draw(13.3,3)--(11.4,2)--(13.3,1);
	  \draw(16.9,3)--(16.9,1);
	   \draw[red](18.4,3)--(18.4,1);
	   \draw[red](19.3,3)--(19.3,1);
\draw(20.9,3)--(20.9,1);
\draw(24.5,3)--(24.5,1);
\draw(26.7,3)--(26.7,1);
	\draw (0,0) node{\color{blue}\normalsize$\ell$\color{black}};
	\braidbox{1}{4.2}{-0.8}{0.7}{$(\ell-1)^2$};
	\draw[dots] (5,0)--(6.5,0);
	\braidbox{7}{10.2}{-0.8}{0.7}{$(i+1)^2$};
	\draw (11.4,0) node{\normalsize$i$};
	\draw (13.3,0) node{\normalsize$i-1$};
	\draw[dots] (14.9,0)--(16.6,0);
	\draw (16.9,0) node{\normalsize$1$};
	\redbraidbox{17.8}{19.9}{-0.8}{0.7}{$\color{red}0\,\,0\color{black}$};
	\draw (20.8,0) node{\normalsize$1$};
	\draw[dots] (21.6,0)--(23,0);
	\draw (24.7,0) node{\normalsize$i-1$};
	\draw (26.8,0) node{\normalsize$i$};
\end{braid}=-\xi_i+\xi_{i-1},
$$
which gives the equality $\xi_i=\xi_{i-1}$ in (\ref{ERelB3}). Continuing like this, we get the equalities $\xi_i=\xi_{i-1}=\dots=\xi_1$. Since the word starting with $\ell\,(\ell-1)^2\,\cdots\,(i+1)^2\,i\,\cdots\,2\,0$ is non-cuspidal by Lemma~\ref{LCuspExpl}, we have  
$$
0=\begin{braid}\tikzset{baseline=1.3em}
	\draw (0,4) node{\color{blue}\normalsize$\ell$\color{black}};
	\braidbox{1}{4.2}{3.2}{4.7}{$(\ell-1)^2$};
	\draw[dots] (5,4)--(6.5,4);
	\braidbox{7}{10.2}{3.2}{4.7}{$(i+1)^2$};
	\draw (11.4,4) node{\normalsize$i$};
	\draw (13.3,4) node{\normalsize$i-1$};
	\draw[dots] (14.9,4)--(16.6,4);
	\draw (16.9,4) node{\normalsize$1$};
	\redbraidbox{17.8}{19.9}{3.2}{4.7}{$\color{red}0\,\,0\color{black}$};
	\draw (20.8,4) node{\normalsize$1$};
	\draw[dots] (21.6,4)--(23,4);
	\draw (24.7,4) node{\normalsize$i-1$};
	\draw (26.8,4) node{\normalsize$i$};
	\draw[blue](0,3)--(0,1);
	 \draw(1.5,3)--(1.5,1);
	 \draw(3.5,3)--(3.5,1);
	  \draw(7.5,3)--(7.5,1);
	  \draw(9.5,3)--(9.5,1);
	  \draw(11.4,3)--(11.4,1);
	  \draw(13.3,3)--(13.3,1);
	  \draw(16.9,3)--(18.4,2)--(16.9,1);
	   \draw[red](18.4,3)--(16.9,2)--(18.4,1);
	   \draw[red](19.3,3)--(19.3,1);
\draw(20.9,3)--(20.9,1);
\draw(24.5,3)--(24.5,1);
\draw(26.7,3)--(26.7,1);
	\draw (0,0) node{\color{blue}\normalsize$\ell$\color{black}};
	\braidbox{1}{4.2}{-0.8}{0.7}{$(\ell-1)^2$};
	\draw[dots] (5,0)--(6.5,0);
	\braidbox{7}{10.2}{-0.8}{0.7}{$(i+1)^2$};
	\draw (11.4,0) node{\normalsize$i$};
	\draw (13.3,0) node{\normalsize$i-1$};
	\draw[dots] (14.9,0)--(16.6,0);
	\draw (16.9,0) node{\normalsize$1$};
	\redbraidbox{17.8}{19.9}{-0.8}{0.7}{$\color{red}0\,\,0\color{black}$};
	\draw (20.8,0) node{\normalsize$1$};
	\draw[dots] (21.6,0)--(23,0);
	\draw (24.7,0) node{\normalsize$i-1$};
	\draw (26.8,0) node{\normalsize$i$};
\end{braid}=-\xi_1+ey_{p-i-1}^2e,
$$
which gives the equality $\xi_1=z$ in (\ref{ERelB3}), using Lemma~\ref{LNH3}(ii). 
Since the word starting with $\ell\,(\ell-1)^2\,\cdots\,(i+1)^2\,i\,\cdots\,1\,0\,1$ is non-cuspidal by Lemma~\ref{LCuspExpl}, we have  
$$
0=\begin{braid}\tikzset{baseline=1.3em}
	\draw (0,4) node{\color{blue}\normalsize$\ell$\color{black}};
	\braidbox{1}{4.2}{3.2}{4.7}{$(\ell-1)^2$};
	\draw[dots] (5,4)--(6.5,4);
	\braidbox{7}{10.2}{3.2}{4.7}{$(i+1)^2$};
	\draw (11.4,4) node{\normalsize$i$};
	\draw (13.3,4) node{\normalsize$i-1$};
	\draw[dots] (14.9,4)--(16.6,4);
	\draw (16.9,4) node{\normalsize$1$};
	\redbraidbox{17.8}{19.9}{3.2}{4.7}{$\color{red}0\,\,0\color{black}$};
	\draw (20.8,4) node{\normalsize$1$};
	\draw[dots] (21.6,4)--(23,4);
	\draw (24.7,4) node{\normalsize$i-1$};
	\draw (26.8,4) node{\normalsize$i$};
	\draw[blue](0,3)--(0,1);
	 \draw(1.5,3)--(1.5,1);
	 \draw(3.5,3)--(3.5,1);
	  \draw(7.5,3)--(7.5,1);
	  \draw(9.5,3)--(9.5,1);
	  \draw(11.4,3)--(11.4,1);
	  \draw(13.3,3)--(13.3,1);
	  \draw(16.9,3)--(18.4,2)--(16.9,1);
	   \draw[red](18.4,3)--(19.3,1);
	   \draw[red](19.3,3)--(16.9,2)--(18.4,1);
\draw(20.9,3)--(20.9,1);
\draw(24.5,3)--(24.5,1);
\draw(26.7,3)--(26.7,1);
	\draw (0,0) node{\color{blue}\normalsize$\ell$\color{black}};
	\braidbox{1}{4.2}{-0.8}{0.7}{$(\ell-1)^2$};
	\draw[dots] (5,0)--(6.5,0);
	\braidbox{7}{10.2}{-0.8}{0.7}{$(i+1)^2$};
	\draw (11.4,0) node{\normalsize$i$};
	\draw (13.3,0) node{\normalsize$i-1$};
	\draw[dots] (14.9,0)--(16.6,0);
	\draw (16.9,0) node{\normalsize$1$};
	\redbraidbox{17.8}{19.9}{-0.8}{0.7}{$\color{red}0\,\,0\color{black}$};
	\draw (20.8,0) node{\normalsize$1$};
	\draw[dots] (21.6,0)--(23,0);
	\draw (24.7,0) node{\normalsize$i-1$};
	\draw (26.8,0) node{\normalsize$i$};
\end{braid}=e\psi_{p-i-1}y_{p-i-1}^2e,
$$
which equals $-\la_0$ thanks to Lemma~\ref{LNH4}(i). This proves (\ref{ERelB4}). 

Again, since the word starting with $\ell\,(\ell-1)^2\,\cdots\,(i+1)^2\,i\,\cdots\,1\,0\,1$ is non-cuspidal by Lemma~\ref{LCuspExpl}, we have  
$$
0=\begin{braid}\tikzset{baseline=1.3em}
	\draw (0,4) node{\color{blue}\normalsize$\ell$\color{black}};
	\braidbox{1}{4.2}{3.2}{4.7}{$(\ell-1)^2$};
	\draw[dots] (5,4)--(6.5,4);
	\braidbox{7}{10.2}{3.2}{4.7}{$(i+1)^2$};
	\draw (11.4,4) node{\normalsize$i$};
	\draw (13.3,4) node{\normalsize$i-1$};
	\draw[dots] (14.9,4)--(16.6,4);
	\draw (16.9,4) node{\normalsize$1$};
	\redbraidbox{17.8}{19.9}{3.2}{4.7}{$\color{red}0\,\,0\color{black}$};
	\draw (20.8,4) node{\normalsize$1$};
	\draw[dots] (21.6,4)--(23,4);
	\draw (24.7,4) node{\normalsize$i-1$};
	\draw (26.8,4) node{\normalsize$i$};
	\draw[blue](0,3)--(0,1);
	 \draw(1.5,3)--(1.5,1);
	 \draw(3.5,3)--(3.5,1);
	  \draw(7.5,3)--(7.5,1);
	  \draw(9.5,3)--(9.5,1);
	  \draw(11.4,3)--(11.4,1);
	  \draw(13.3,3)--(13.3,1);
	  \draw(16.9,3)--(16.9,1);
	   \draw[red](18.4,3)--(18.4,1);
	   \draw[red](19.3,3)--(20.9,2)--(19.3,1);
\draw(20.9,3)--(19.3,2)--(20.9,1);
\draw(24.5,3)--(24.5,1);
\draw(26.7,3)--(26.7,1);
	\draw (0,0) node{\color{blue}\normalsize$\ell$\color{black}};
	\braidbox{1}{4.2}{-0.8}{0.7}{$(\ell-1)^2$};
	\draw[dots] (5,0)--(6.5,0);
	\braidbox{7}{10.2}{-0.8}{0.7}{$(i+1)^2$};
	\draw (11.4,0) node{\normalsize$i$};
	\draw (13.3,0) node{\normalsize$i-1$};
	\draw[dots] (14.9,0)--(16.6,0);
	\draw (16.9,0) node{\normalsize$1$};
	\redbraidbox{17.8}{19.9}{-0.8}{0.7}{$\color{red}0\,\,0\color{black}$};
	\draw (20.8,0) node{\normalsize$1$};
	\draw[dots] (21.6,0)--(23,0);
	\draw (24.7,0) node{\normalsize$i-1$};
	\draw (26.8,0) node{\normalsize$i$};
\end{braid}=ey_{p-i}^2e-\eta_1,
$$
which equals $\la_0^2-z-\eta_1$ thanks to Lemma~\ref{LNH3}(i). Since $\la_0^2=0$ by (\ref{ERelB4}), we get 
the equality $z=-\eta_1$ in (\ref{ERelB3}).  

Since the word starting with $\ell\,(\ell-1)^2\,\cdots\,(i+1)^2\,i\,\cdots\,1\,0^2\,2$ is non-cuspidal by Lemma~\ref{LCuspExpl}, we have  
$$
0=\begin{braid}\tikzset{baseline=1.3em}
	\draw (0,4) node{\color{blue}\normalsize$\ell$\color{black}};
	\braidbox{1}{4.2}{3.2}{4.7}{$(\ell-1)^2$};
	\draw[dots] (5,4)--(6.5,4);
	\braidbox{7}{10.2}{3.2}{4.7}{$(i+1)^2$};
	\draw (11.4,4) node{\normalsize$i$};
	\draw (13.3,4) node{\normalsize$i-1$};
	\draw[dots] (14.9,4)--(16.6,4);
	\draw (16.9,4) node{\normalsize$1$};
	\redbraidbox{17.8}{19.9}{3.2}{4.7}{$\color{red}0\,\,0\color{black}$};
	\draw (20.8,4) node{\normalsize$1$};
	\draw[dots] (22.6,4)--(24,4);
	\draw (21.8,4) node{\normalsize$2$};
	\draw (24.8,4) node{\normalsize$i$};
	\draw[blue](0,3)--(0,1);
	 \draw(1.5,3)--(1.5,1);
	 \draw(3.5,3)--(3.5,1);
	  \draw(7.5,3)--(7.5,1);
	  \draw(9.5,3)--(9.5,1);
	  \draw(11.4,3)--(11.4,1);
	  \draw(13.3,3)--(13.3,1);
	  \draw(16.9,3)--(16.9,1);
	   \draw[red](18.4,3)--(18.4,1);
	   \draw[red](19.3,3)--(19.3,1);
\draw(20.8,3)--(21.8,2)--(20.8,1);
\draw(21.8,3)--(20.8,2)--(21.8,1);
\draw(24.8,3)--(24.8,1);
	\draw (0,0) node{\color{blue}\normalsize$\ell$\color{black}};
	\braidbox{1}{4.2}{-0.8}{0.7}{$(\ell-1)^2$};
	\draw[dots] (5,0)--(6.5,0);
	\braidbox{7}{10.2}{-0.8}{0.7}{$(i+1)^2$};
	\draw (11.4,0) node{\normalsize$i$};
	\draw (13.3,0) node{\normalsize$i-1$};
	\draw[dots] (14.9,0)--(16.6,0);
	\draw (16.9,0) node{\normalsize$1$};
	\redbraidbox{17.8}{19.9}{-0.8}{0.7}{$\color{red}0\,\,0\color{black}$};
	\draw (20.8,0) node{\normalsize$1$};
	\draw[dots] (22.6,0)--(24,0);
	\draw (21.8,0) node{\normalsize$2$};
	\draw (24.8,0) node{\normalsize$i$};
\end{braid}=\eta_1-\eta_2.
$$
Continuing like this we get the equalities $\eta_1=\eta_2=\dots=\eta_{i-1}$ and complete the proof of (\ref{ERelB3}). 

Finally, to prove (\ref{ERelB5}), note that the word $\ell\,(\ell-1)^2\,\cdots\, 2^2\,0\,1\,1\,0$ is non-cuspidal by Lemma~\ref{LCuspExpl}. So using the quadratic relations (\ref{R6}) to get 
\begin{align*}
0&=\begin{braid}\tikzset{baseline=1.3em}
	\draw (0,0) node{\color{blue}\normalsize$\ell$\color{black}};
	\braidbox{1}{4.2}{-0.8}{0.7}{$(\ell-1)^2$};
	\draw[dots] (5,0)--(6.5,0);
	\braidbox{7}{9.1}{-0.8}{0.7}{$1\,\,1$};
	\redbraidbox{10}{12.1}{-0.8}{0.7}{$\color{red}0\,\,0\color{black}$};
	\draw (0,4) node{\color{blue}\normalsize$\ell$\color{black}};
	\braidbox{1}{4.2}{3.2}{4.7}{$(\ell-1)^2$};
	\draw[dots] (5,4)--(6.5,4);
	\braidbox{7}{9.1}{3.2}{4.7}{$1\,\,1$};
	\redbraidbox{10}{12.1}{3.2}{4.7}{$\color{red}0\,\,0\color{black}$};
	 \draw[blue](0,3)--(0,1);
	 \draw(1.5,3)--(1.5,1);
	 \draw(3.5,3)--(3.5,1);
	 \draw(7.6,3)--(7.6,1);
	 \draw(8.5,3)--(8.5,1);
	  \draw[red](11.5,3)--(10.6,2.7)--(7.2,2.4)--(7.2,1.6)--(10.6,1.3)--(10.6,1);
	  \draw[red](10.6,3)--(11.5,1);
	  \end{braid}
\\&= ey_{p-3}y_{p-2}\psi_{p-1}e-
e(y_{p-3}+y_{p-2})\psi_{p-1}y_{p-1}^2e
+e\psi_{p-1}y_{p-1}^4e.
\end{align*}
Note that $ey_{p-3}y_{p-2}\psi_{p-1}e=0$ by Lemma~\ref{LPsiE} and 
$$e(y_{p-3}+y_{p-2})\psi_{p-1}y_{p-1}^2e=e(y_{p-3}+y_{p-2})e\psi_{p-1}y_{p-1}^2e
=\la_1\psi_{p-1}y_{p-1}^2e=
0$$ since $\la_1=0$ by (\ref{ERelB3}). 
Finally, $e\psi_{p-1}y_{p-1}^4e=-u^3$ by Lemma~\ref{LNH4}(i). 
\end{proof}

\begin{Corollary} \label{CIBI}
For $i\in J$ we have
$$
e^iBe^i=
\left\{
\begin{array}{ll}
\spa\big\{(z^i)^a(c^i)^b\mid a,b\in\Z_{\geq 0}\big\} &\hbox{if $i\neq 0$,}\\
\spa\big\{(z^0)^a u^b\mid a\in\Z_{\geq 0},\, b\in\{0,1,2\}\big\} &\hbox{if $i=0$.}
\end{array}
\right.
$$
\end{Corollary}
\begin{proof}
By Lemma~\ref{LIBIImprove}, 
$e^0Be^0$ is generated by the elements 
$$\ka^0,\la^0_0=u,z^0,\la^0_{1},\mu^0_{1},\dots,\la^0_{\ell-1},\mu^0_{\ell-1}.
$$
Moreover, by Lemma~\ref{LRelB}, we have $\la^0_{1}=\dots=\la^0_{\ell-1}=0$, 
$-\ka^0=\mu^0_{1}=\dots=\mu^0_{\ell-1}=
-u^2z^0-(z^0)^2$, and 
$
u^3=0.$ This implies the corollary for $i=0$. 
For $i\neq 0$ the proof is similar recalling that 
by definition 
$c^i=
(-1)^i(z^i+\eta^i_i).$
\end{proof}

\begin{Corollary} \label{CIBI-1}
For $1\leq i\leq \ell-1$, we have
\begin{align*}
e^{i}Be^{i-1}&=\spa\big\{(z^i)^a(c^i)^b a^{i,i-1}\mid a,b\in\Z_{\geq 0}\big\},
\\
e^{i-1}Be^i&=\spa\big\{a^{i-1,i}(z^i)^a(c^i)^b\mid a,b\in\Z_{\geq 0}\big\}.
\end{align*}
\end{Corollary}
\begin{proof}
Note that $e^i\F[y_1,\dots,y_p]a^{i,i-1}=e^i\F[y_1,\dots,y_p]e^ia^{i,i-1}$. So the first equality follows from Lemma~\ref{LIBI-1} and Corollary~\ref{CIBI}. The proof of the second equality is similar, using Lemmas~\ref{LI-1BI} instead of Lemma~\ref{LIBI-1}. 
\end{proof}

\begin{Corollary} \label{CBGen}
The algebra $B$ is generated by all $e^i,z^i,c^i,a^{i,i-1},a^{i-1,i}$ and $u$. 
\end{Corollary}
\begin{proof}
This follows from Lemma~\ref{LIBJ} and Corollaries~\ref{CIBI} and \ref{CIBI-1}.
\end{proof}

\begin{Corollary} \label{CSumProduct}
Let $i\in J$. Then in $B$ we have:
\begin{enumerate}
\item[{\rm (i)}] $e^i(y_{r^i_j}+y_{s^i_j})=(y_{r^i_j}+y_{s^i_j})e^i=\de_{i,j}(-1)^ic^i$ for any  $j\in J\setminus\{0\}$.
\item[{\rm (ii)}] $e^iy_{r^i_j}y_{s^i_j}=y_{r^i_j}y_{s^i_j}e^i=-(z^i)^2-\de_{i,0}u^2z^0+\de_{i,j}(-1)^iz^ic^i$ for any  $j\in J\setminus\{0\}$.
\item[{\rm (iii)}] $e^i(y_{r^i_0}^2+y_{s^i_0}^2)=(y_{r^i_0}^2+y_{s^i_0}^2)e^i=\de_{i,0}u^2$.
\item[{\rm (iv)}] $e^iy_{r^i_0}^2y_{s^i_0}^2=y_{r^i_0}^2y_{s^i_0}^2e^i=-(z^i)^2$.
\end{enumerate}
\end{Corollary}
\begin{proof}
(i) and (ii). That $(y_{r^i_j}+y_{s^i_j})$ and $y_{r^i_j}y_{s^i_j}$ commute with $e^i$ follows from the fact that the polynomials are symmetric in $y_{r^i_j}$ and $y_{s^i_j}$. 

If $i<j$, then $e^i(y_{r^i_j}+y_{s^i_j})e^i=\la^i_j$ and $e^i(y_{r^i_j}y_{s^i_j})e^i=\mu^i_j$. So in this case (i) follows from (\ref{ERelB2}) and (ii) from (\ref{ERelB1}),(\ref{ERelB3}). 

If $i\geq j$ then $e^i(y_{r^i_j}+y_{s^i_j})e^i=\xi^i_j+\eta^i_j$ and $e^iy_{r^i_j}y_{s^i_j}e^i=\xi^i_j\eta^i_j$. If $i>j$, then by (\ref{ERelB3}) we have $\xi^i_j+\eta^i_j=0$  and $\xi^i_j\eta^i_j=z^i(-z^i)=-(z^i)^2$, as desired. 
If $i=j$ then by (\ref{ERelB3}) and (\ref{EZFormula}), we get 
\begin{align*}
&e^i(y_{r^i_i}+y_{s^i_i})e^i=\xi^i_i+e^iy_pe^i=z^i+((-1)^ic^i-z^i)=(-1)^ic^i,
\\
&e^iy_{r^i_i}y_{s^i_i}e^i=(\xi^i_i)(e^iy_pe^i)=z^i((-1)^ic^i-z^i)=-(z^i)^2+(-1)^ic^iz^i.
\end{align*}

(iii) follows from Lemmas~\ref{LNH3} and (\ref{ERelB4}).

(iv) follows from the definition of $z^i$. 
\end{proof}

\subsection{Relations in $B$}
We continue with the notation of the previous subsection.

\begin{Lemma} \label{LAA0} 
For $1\leq i\leq \ell-2$, we have $a^{i-1,i}a^{i,i+1}=0$ and $a^{i+1,i}a^{i,i-1}=0$. 
\end{Lemma}
\begin{proof}
We prove the first equality, the proof of the second one being similar. 
If $i>1$ then 
\begin{align*}
a^{i-1,i}a^{i,i+1}&=\pm e^{i-1}\psi_{p-2i-1}\cdots\psi_{p-2}\psi_{p-1}e^{i}
\psi_{p-2i-3}\cdots\psi_{p-2}\psi_{p-1}e^{i+1}
\\
&=\pm e^{i-1}\psi_{p-2i-1}\cdots\psi_{p-2}e^{i}
\psi_{p-2i-3}\cdots\psi_{p-1}\psi_{p-2}\psi_{p-1}e^{i+1}
\\
&=\pm e^{i-1}\psi_{p-2i-1}\cdots\psi_{p-2}e^{i}
\psi_{p-2i-3}\cdots\psi_{p-2}\psi_{p-1}\psi_{p-2}e^{i+1}=0,
\end{align*}
since $s_{p-2}\cdot\hat\bg^{i+1}$ is non-cuspidal by Lemma~\ref{LCuspExpl}.
On the other hand, if $i=1$ then, in terms of diagrams, we have that $a^{0,1}a^{1,2}$ equals 
\begin{align*}
\begin{braid}\tikzset{baseline=-0.3em}
\draw (0,-4) node{\color{blue}\normalsize$\ell$\color{black}};
	\braidbox{1}{4.2}{-4.8}{-3.3}{$(\ell-1)^2$};
	\draw[dots] (5,-4)--(6.5,-4);
	\braidbox{7}{9.1}{-4.8}{-3.3}{$3\,\,3$};
	\draw (10.3,-4) node{\normalsize$2$};
	\draw (11.6,-4) node{\normalsize$1$};
	\redbraidbox{12.9}{15}{-4.8}{-3.3}{$\color{red}0\,\,0\color{black}$};
	\draw (16.1,-4) node{\normalsize$1$};
	\draw (17.8,-4) node{\normalsize$2$};
	\draw (0,0) node{\color{blue}\normalsize$\ell$\color{black}};
	\braidbox{1}{4.2}{-0.8}{0.7}{$(\ell-1)^2$};
	\draw[dots] (5,0)--(6.5,0);
	\braidbox{7}{9.1}{-0.8}{0.7}{$3\,\,3$};
	\braidbox{10}{12.1}{-0.8}{0.7}{$2\,\,2$};
	\draw (13.4,0) node{\normalsize$1$};
	\redbraidbox{14.5}{16.6}{-0.8}{0.7}{$\color{red}0\,\,0\color{black}$};
	\draw (17.8,0) node{\normalsize$1$};
	\draw (0,4) node{\color{blue}\normalsize$\ell$\color{black}};
	\braidbox{1}{4.2}{3.2}{4.7}{$(\ell-1)^2$};
	\draw[dots] (5,4)--(6.5,4);
	\braidbox{7}{9.1}{3.2}{4.7}{$3\,\,3$};
	\braidbox{10}{12.1}{3.2}{4.7}{$2\,\,2$};
	\braidbox{13}{15.1}{3.2}{4.7}{$1\,\,1$};
	\redbraidbox{16}{18.1}{3.2}{4.7}{$\color{red}0\,\,0\color{black}$};
	 \draw[blue](0,3)--(0,1);
	 \draw[blue](0,-1)--(0,-3);
	 \draw(1.5,3)--(1.5,1);
	  \draw(1.5,-1)--(1.5,-3);
	 \draw(3.5,3)--(3.5,1);
	  \draw(3.5,-1)--(3.5,-3);
	 \draw(7.6,3)--(7.6,1);
	  \draw(7.6,-1)--(7.6,-3);
	 \draw(8.5,3)--(8.5,1);
	  \draw(8.5,-1)--(8.5,-3);
	 \draw(10.5,3)--(10.5,1);
	  \draw(11.5,3)--(11.5,1);
	  \draw(13.6,3)--(17.7,1);
	  \draw(14.6,3)--(13.6,1);
	  \draw(11.5,-1)--(10.4,-3);
	  \draw[red](16.6,3)--(15.1,1);
	  \draw[red](17.5,3)--(16,1);
	  \draw[red](15.1,-1)--(13.6,-3);
	  \draw[red](16.1,-1)--(14.6,-3);
	  \draw(13.4,-1)--(11.7,-3);
	  \draw(10.5,-1)--(17.7,-3);
	  \draw(17.7,-1)--(16.2,-3);
	  \end{braid}
	  =
\begin{braid}\tikzset{baseline=-0.3em}
\draw (0,-4) node{\color{blue}\normalsize$\ell$\color{black}};
	\braidbox{1}{4.2}{-4.8}{-3.3}{$(\ell-1)^2$};
	\draw[dots] (5,-4)--(6.5,-4);
	\braidbox{7}{9.1}{-4.8}{-3.3}{$3\,\,3$};
	\draw (10.3,-4) node{\normalsize$2$};
	\draw (11.6,-4) node{\normalsize$1$};
	\redbraidbox{12.9}{15}{-4.8}{-3.3}{$\color{red}0\,\,0\color{black}$};
	\draw (16.1,-4) node{\normalsize$1$};
	\draw (17.8,-4) node{\normalsize$2$};
	\draw (0,0) node{\color{blue}\normalsize$\ell$\color{black}};
	\braidbox{1}{4.2}{-0.8}{0.7}{$(\ell-1)^2$};
	\draw[dots] (5,0)--(6.5,0);
	\braidbox{7}{9.1}{-0.8}{0.7}{$3\,\,3$};
	\braidbox{10}{12.1}{-0.8}{0.7}{$2\,\,2$};
	\redbraidbox{14.5}{16.6}{-0.8}{0.7}{$\color{red}0\,\,0\color{black}$};
	\draw (17.8,0) node{\normalsize$1$};
	\draw (0,4) node{\color{blue}\normalsize$\ell$\color{black}};
	\braidbox{1}{4.2}{3.2}{4.7}{$(\ell-1)^2$};
	\draw[dots] (5,4)--(6.5,4);
	\braidbox{7}{9.1}{3.2}{4.7}{$3\,\,3$};
	\braidbox{10}{12.1}{3.2}{4.7}{$2\,\,2$};
	\braidbox{13}{15.1}{3.2}{4.7}{$1\,\,1$};
	\redbraidbox{16}{18.1}{3.2}{4.7}{$\color{red}0\,\,0\color{black}$};
	 \draw[blue](0,3)--(0,1);
	 \draw[blue](0,-1)--(0,-3);
	 \draw(1.5,3)--(1.5,1);
	  \draw(1.5,-1)--(1.5,-3);
	 \draw(3.5,3)--(3.5,1);
	  \draw(3.5,-1)--(3.5,-3);
	 \draw(7.6,3)--(7.6,1);
	  \draw(7.6,-1)--(7.6,-3);
	 \draw(8.5,3)--(8.5,1);
	  \draw(8.5,-1)--(8.5,-3);
	 \draw(10.5,3)--(10.5,1);
	  \draw(11.5,3)--(11.5,1);
	  \draw(13.6,3)--(14.7,2)--(12.8,-2.2)--(16,-3);
	  \draw(14.6,3)--(11.7,-3);
	  \draw(11.5,-1)--(10.4,-3);
	  \draw[red](16.6,3)--(15.1,1);
	  \draw[red](17.5,3)--(16,1);
	  \draw[red](15.1,-1)--(13.6,-3);
	  \draw[red](16.1,-1)--(14.6,-3);
	  \draw(10.5,-1)--(17.7,-3);
	  \end{braid}
\end{align*}
which is $0$, as the word $\ell\,(\ell-1)^2\,\cdots\, 3^2\,2\,1^2\,0^2\,2$ is non-cuspidal by Lemma~\ref{LCuspExpl}.
\end{proof}

\begin{Lemma} \label{LLoopU^2} 
We have $a^{0,1}a^{1,0}=u^2$. 
\end{Lemma}
\begin{proof}
Using braid relations (\ref{R7}) and dot-crossing relations  (\ref{R5}) we get
\begin{align*}
a^{0,1}a^{1,0}
&=-
\begin{braid}\tikzset{baseline=-0.3em}
\draw (0,-4) node{\color{blue}\normalsize$\ell$\color{black}};
	\braidbox{1}{4.2}{-4.8}{-3.3}{$(\ell-1)^2$};
	\draw[dots] (5,-4)--(6.5,-4);
	\braidbox{7}{9.1}{-4.8}{-3.3}{$2\,\,2$};
	\braidbox{10}{12.1}{-4.8}{-3.3}{$1\,\,1$};
	\redbraidbox{13}{15.1}{-4.8}{-3.3}{$\color{red}0\,\,0\color{black}$};
	\draw (0,0) node{\color{blue}\normalsize$\ell$\color{black}};
	\braidbox{1}{4.2}{-0.8}{0.7}{$(\ell-1)^2$};
	\draw[dots] (5,0)--(6.5,0);
	\braidbox{7}{9.1}{-0.8}{0.7}{$2\,\,2$};
	\draw (10.4,0) node{\normalsize$1$};
	\redbraidbox{11.5}{13.6}{-0.8}{0.7}{$\color{red}0\,\,0\color{black}$};
	\draw (14.8,0) node{\normalsize$1$};
	\draw (0,4) node{\color{blue}\normalsize$\ell$\color{black}};
	\braidbox{1}{4.2}{3.2}{4.7}{$(\ell-1)^2$};
	\draw[dots] (5,4)--(6.5,4);
	\braidbox{7}{9.1}{3.2}{4.7}{$2\,\,2$};
	\braidbox{10}{12.1}{3.2}{4.7}{$1\,\,1$};
	\redbraidbox{13}{15.1}{3.2}{4.7}{$\color{red}0\,\,0\color{black}$};
	 \draw[blue](0,3)--(0,1);
	 \draw[blue](0,-1)--(0,-3);
	 \draw(1.5,3)--(1.5,1);
	  \draw(1.5,-1)--(1.5,-3);
	 \draw(3.5,3)--(3.5,1);
	  \draw(3.5,-1)--(3.5,-3);
	 \draw(7.6,3)--(7.6,1);
	  \draw(7.6,-1)--(7.6,-3);
	 \draw(8.5,3)--(8.5,1);
	  \draw(8.5,-1)--(8.5,-3);
	  \draw(11.5,3)--(10.5,1);
	  \draw(10.6,3)--(14.7,1);
	  \draw[red](13.6,3)--(12.1,1);
	  \draw[red](14.5,3)--(13,1);
	  \draw[red](12.1,-1)--(13.6,-3);
	  \draw[red](13.1,-1)--(14.6,-3);
	  \draw(14.7,-1)--(11.7,-3);
	  \draw(10.5,-1)--(10.5,-3);
	  \end{braid}
\\&=-
\begin{braid}\tikzset{baseline=-0.3em}
\draw (0,-4) node{\color{blue}\normalsize$\ell$\color{black}};
	\braidbox{1}{4.2}{-4.8}{-3.3}{$(\ell-1)^2$};
	\draw[dots] (5,-4)--(6.5,-4);
	\braidbox{7}{9.1}{-4.8}{-3.3}{$2\,\,2$};
	\braidbox{10}{12.1}{-4.8}{-3.3}{$1\,\,1$};
	\redbraidbox{13}{15.1}{-4.8}{-3.3}{$\color{red}0\,\,0\color{black}$};
	\draw (0,4) node{\color{blue}\normalsize$\ell$\color{black}};
	\braidbox{1}{4.2}{3.2}{4.7}{$(\ell-1)^2$};
	\draw[dots] (5,4)--(6.5,4);
	\braidbox{7}{9.1}{3.2}{4.7}{$2\,\,2$};
	\braidbox{10}{12.1}{3.2}{4.7}{$1\,\,1$};
	\redbraidbox{13}{15.1}{3.2}{4.7}{$\color{red}0\,\,0\color{black}$};
	 \draw[blue](0,3)--(0,-3);
	 \draw(1.5,3)--(1.5,-3);
	 \draw(3.5,3)--(3.5,-3);
	 \draw(7.6,3)--(7.6,-3);
	 \draw(8.5,3)--(8.5,-3);
	  \draw(11.5,3)--(10.5,-3);
	  \draw(10.6,3)--(14.7,1)--(14.7,-1)--(11.6,-3);
	  \draw[red](13.6,3)--(12.8,1)--(12.8,-1)--(13.6,-3);
	  \draw[red](14.5,3)--(13.7,1)--(13.7,-1)--(14.5,-3);
	  \end{braid},
\end{align*}
which, using quadratic relations (\ref{R6}), equals 
$$
e^0\psi_{p-3}(-y_{p-2}^2+y_{p-2}(y_{p-1}^2+y_p^2)-y_{p-1}^2y_{p}^2)e^0.
$$
Note that 
$$-e^0\psi_{p-3}y_{p-1}^2y_{p}^2e^0=-e^0y_{p-1}^2y_{p}^2\psi_{p-3}e^0=0
$$ using Lemma~\ref{LPsiE} for the last equality.
Moreover, 
$$-e^0\psi_{p-3}y_{p-2}^2=-e^0(y_{p-3}+y_{p-2})e^0=-\la^1_0=0$$ by (\ref{ERelB2}). Finally, 
$\psi_{p-3}(y_{p-1}^2+y_p^2))e^0=(y_{p-1}^2+y_p^2))\psi_{p-3}e^0=0$ by Lemma~\ref{LPsiE}, so we have 
$$
e^0\psi_{p-3}(y_{p-2}(y_{p-1}^2+y_p^2))e^0=
e^0(y_{p-1}^2+y_p^2)e^0=z^0+(u^2-z^0)=u^2
$$
thanks to Lemma~\ref{LNH3}. 
\end{proof}

\begin{Lemma} \label{LAAC}
For $1\leq i\leq \ell-2$ we have $a^{i,i-1}a^{i-1,i}=a^{i,i+1}a^{i+1,i}=c^i$. 
\end{Lemma}
\begin{proof}
We have using quadratic and dot-crossing relations:
\begin{align*}
a^{i,i-1}a^{i-1,i}&=(-1)^i
\begin{braid}\tikzset{baseline=-.3em}
	\draw (0,4) node{\color{blue}\normalsize$\ell$\color{black}};
	\braidbox{1}{4.2}{3.2}{4.7}{$(\ell-1)^2$};
	\draw[dots] (5,4)--(6.5,4);
	\braidbox{7}{10.2}{3.2}{4.7}{$(i+1)^2$};
	\draw (11.4,4) node{\normalsize$i$};
	\draw (13.3,4) node{\normalsize$i-1$};
	\draw[dots] (14.9,4)--(16.6,4);
	\draw (16.9,4) node{\normalsize$1$};
	\redbraidbox{17.8}{19.9}{3.2}{4.7}{$\color{red}0\,\,0\color{black}$};
	\draw (20.8,4) node{\normalsize$1$};
	\draw[dots] (21.6,4)--(23,4);
	\draw (24.7,4) node{\normalsize$i-1$};
	\draw (26.8,4) node{\normalsize$i$};
	\draw[blue](0,3)--(0,1);
	 \draw(1.5,3)--(1.5,1);
	 \draw(3.5,3)--(3.5,1);
	  \draw(7.5,3)--(7.5,1);
	  \draw(9.5,3)--(9.5,1);
	  \draw(11.4,3)--(11.4,1);
	  \draw(13.3,3)--(14.7,1);
	  \draw(16.9,3)--(18.4,1);
	   \draw[red](18.4,3)--(19.9,1);
	   \draw[red](19.3,3)--(20.8,1);
\draw(20.9,3)--(22.3,1);
\draw(24.7,3)--(26.2,1);
\draw(26.7,3)--(12.5,1);
	\draw (0,0) node{\color{blue}\normalsize$\ell$\color{black}};
	\braidbox{1}{4.2}{-0.8}{0.7}{$(\ell-1)^2$};
	\draw[dots] (5,0)--(6.5,0);
	\braidbox{7}{10.2}{-0.8}{0.7}{$(i+1)^2$};
	\braidbox{11}{13.1}{-0.8}{0.7}{$i\,\, i$};
	\draw (15,0) node{\normalsize$i-1$};
	\draw[dots] (16.5,0)--(18.2,0);
	\draw (18.5,0) node{\normalsize$1$};
	\redbraidbox{19.4}{21.5}{-0.8}{0.7}{$\color{red}0\,\,0\color{black}$};
	\draw (22.4,0) node{\normalsize$1$};
	\draw[dots] (23.3,0)--(24.7,0);
	\draw (26.5,0) node{\normalsize$i-1$};
	\draw (0,-4) node{\color{blue}\normalsize$\ell$\color{black}};
	\braidbox{1}{4.2}{-4.8}{-3.3}{$(\ell-1)^2$};
	\draw[dots] (5,-4)--(6.5,-4);
	\braidbox{7}{10.2}{-4.8}{-3.3}{$(i+1)^2$};
	\draw (11.4,-4) node{\normalsize$i$};
	\draw (13.3,-4) node{\normalsize$i-1$};
	\draw[dots] (14.9,-4)--(16.6,-4);
	\draw (16.9,-4) node{\normalsize$1$};
	\redbraidbox{17.8}{19.9}{-4.8}{-3.3}{$\color{red}0\,\,0\color{black}$};
	\draw (20.8,-4) node{\normalsize$1$};
	\draw[dots] (21.6,-4)--(23,-4);
	\draw (24.7,-4) node{\normalsize$i-1$};
	\draw (26.8,-4) node{\normalsize$i$};
	\draw[blue](0,-3)--(0,-1);
	 \draw(1.5,-3)--(1.5,-1);
	 \draw(3.5,-3)--(3.5,-1);
	  \draw(7.5,-3)--(7.5,-1);
	  \draw(9.5,-3)--(9.5,-1);
	  \draw(11.4,-3)--(12.4,-1);
	  \draw(11.5,-1)--(26.8,-3);
	  \draw(13.4,-3)--(14.9,-1);
	  \draw(16.9,-3)--(18.4,-1);
	   \draw[red](18.4,-3)--(19.9,-1);
	   \draw[red](19.3,-3)--(20.8,-1);
\draw(20.9,-3)--(22.3,-1);
\draw(24.7,-3)--(26.2,-1);
\end{braid}
\\
&=(-1)^i
\begin{braid}\tikzset{baseline=-.3em}
	\draw (0,4) node{\color{blue}\normalsize$\ell$\color{black}};
	\braidbox{1}{4.2}{3.2}{4.7}{$(\ell-1)^2$};
	\draw[dots] (5,4)--(6.5,4);
	\braidbox{7}{10.2}{3.2}{4.7}{$(i+1)^2$};
	\draw (11.4,4) node{\normalsize$i$};
	\draw (13.3,4) node{\normalsize$i-1$};
	\draw[dots] (14.9,4)--(16.6,4);
	\draw (16.9,4) node{\normalsize$1$};
	\redbraidbox{17.8}{19.9}{3.2}{4.7}{$\color{red}0\,\,0\color{black}$};
	\draw (20.8,4) node{\normalsize$1$};
	\draw[dots] (21.6,4)--(23,4);
	\draw (24.7,4) node{\normalsize$i-1$};
	\draw (26.8,4) node{\normalsize$i$};
	\draw[blue](0,3)--(0,-3);
	 \draw(1.5,3)--(1.5,-3);
	 \draw(3.5,3)--(3.5,-3);
	  \draw(7.5,3)--(7.5,-3);
	  \draw(9.5,3)--(9.5,-3);
	  \draw(11.4,3)--(12,0)--(26.8,-3);
	  \draw(13.3,3)--(13.3,-3);
	  \draw(16.9,3)--(16.9,-3);
	   \draw[red](18.4,3)--(18.4,-3);
	   \draw[red](19.3,3)--(19.3,-3);
\draw(20.9,3)--(20.9,-3);
\draw(24.7,3)--(24.7,-3);
\draw(26.8,3)--(12,0)--(11.4,-3);
		\draw (0,-4) node{\color{blue}\normalsize$\ell$\color{black}};
	\braidbox{1}{4.2}{-4.8}{-3.3}{$(\ell-1)^2$};
	\draw[dots] (5,-4)--(6.5,-4);
	\braidbox{7}{10.2}{-4.8}{-3.3}{$(i+1)^2$};
	\draw (11.4,-4) node{\normalsize$i$};
	\draw (13.3,-4) node{\normalsize$i-1$};
	\draw[dots] (14.9,-4)--(16.6,-4);
	\draw (16.9,-4) node{\normalsize$1$};
	\redbraidbox{17.8}{19.9}{-4.8}{-3.3}{$\color{red}0\,\,0\color{black}$};
	\draw (20.8,-4) node{\normalsize$1$};
	\draw[dots] (21.6,-4)--(23,-4);
	\draw (24.7,-4) node{\normalsize$i-1$};
	\draw (26.8,-4) node{\normalsize$i$};
	\end{braid}.
\end{align*}
Using the braid relation (\ref{R7}), this equals 
\begin{align*}
&=(-1)^i
\begin{braid}\tikzset{baseline=-.3em}
	\draw (0,4) node{\color{blue}\normalsize$\ell$\color{black}};
	\braidbox{1}{4.2}{3.2}{4.7}{$(\ell-1)^2$};
	\draw[dots] (5,4)--(6.5,4);
	\braidbox{7}{10.2}{3.2}{4.7}{$(i+1)^2$};
	\draw (11.4,4) node{\normalsize$i$};
	\draw (13.3,4) node{\normalsize$i-1$};
	\draw[dots] (14.9,4)--(16.6,4);
	\draw (16.9,4) node{\normalsize$1$};
	\redbraidbox{17.8}{19.9}{3.2}{4.7}{$\color{red}0\,\,0\color{black}$};
	\draw (20.8,4) node{\normalsize$1$};
	\draw[dots] (21.6,4)--(23,4);
	\draw (24.7,4) node{\normalsize$i-1$};
	\draw (26.8,4) node{\normalsize$i$};
	\draw[blue](0,3)--(0,-3);
	 \draw(1.5,3)--(1.5,-3);
	 \draw(3.5,3)--(3.5,-3);
	  \draw(7.5,3)--(7.5,-3);
	  \draw(9.5,3)--(9.5,-3);
	  \draw(11.4,3)--(12,0)--(26.8,-3);
	  \draw(13.3,3)--(11.3,0)--(13.3,-3);
	  \draw(16.9,3)--(16.9,-3);
	   \draw[red](18.4,3)--(18.4,-3);
	   \draw[red](19.3,3)--(19.3,-3);
\draw(20.9,3)--(20.9,-3);
\draw(24.7,3)--(24.7,-3);
\draw(26.8,3)--(12,0)--(11.4,-3);
		\draw (0,-4) node{\color{blue}\normalsize$\ell$\color{black}};
	\braidbox{1}{4.2}{-4.8}{-3.3}{$(\ell-1)^2$};
	\draw[dots] (5,-4)--(6.5,-4);
	\braidbox{7}{10.2}{-4.8}{-3.3}{$(i+1)^2$};
	\draw (11.4,-4) node{\normalsize$i$};
	\draw (13.3,-4) node{\normalsize$i-1$};
	\draw[dots] (14.9,-4)--(16.6,-4);
	\draw (16.9,-4) node{\normalsize$1$};
	\redbraidbox{17.8}{19.9}{-4.8}{-3.3}{$\color{red}0\,\,0\color{black}$};
	\draw (20.8,-4) node{\normalsize$1$};
	\draw[dots] (21.6,-4)--(23,-4);
	\draw (24.7,-4) node{\normalsize$i-1$};
	\draw (26.8,-4) node{\normalsize$i$};
	\end{braid}
\\
&\hspace{5mm}-(-1)^i
\begin{braid}\tikzset{baseline=-.3em}
	\draw (0,4) node{\color{blue}\normalsize$\ell$\color{black}};
	\braidbox{1}{4.2}{3.2}{4.7}{$(\ell-1)^2$};
	\draw[dots] (5,4)--(6.5,4);
	\braidbox{7}{10.2}{3.2}{4.7}{$(i+1)^2$};
	\draw (11.4,4) node{\normalsize$i$};
	\draw (13.3,4) node{\normalsize$i-1$};
	\draw[dots] (14.9,4)--(16.6,4);
	\draw (16.9,4) node{\normalsize$1$};
	\redbraidbox{17.8}{19.9}{3.2}{4.7}{$\color{red}0\,\,0\color{black}$};
	\draw (20.8,4) node{\normalsize$1$};
	\draw[dots] (21.6,4)--(23,4);
	\draw (24.7,4) node{\normalsize$i-1$};
	\draw (26.8,4) node{\normalsize$i$};
	\draw[blue](0,3)--(0,-3);
	 \draw(1.5,3)--(1.5,-3);
	 \draw(3.5,3)--(3.5,-3);
	  \draw(7.5,3)--(7.5,-3);
	  \draw(9.5,3)--(9.5,-3);
	  \draw(11.4,3)--(11.4,-3);
	  \draw(13.3,3)--(13.3,-3);
	  \draw(16.9,3)--(16.9,-3);
	   \draw[red](18.4,3)--(18.4,-3);
	   \draw[red](19.3,3)--(19.3,-3);
\draw(20.9,3)--(20.9,-3);
\draw(24.7,3)--(24.7,-3);
\draw(26.8,3)--(13.8,0)--(26.8,-3);
		\draw (0,-4) node{\color{blue}\normalsize$\ell$\color{black}};
	\braidbox{1}{4.2}{-4.8}{-3.3}{$(\ell-1)^2$};
	\draw[dots] (5,-4)--(6.5,-4);
	\braidbox{7}{10.2}{-4.8}{-3.3}{$(i+1)^2$};
	\draw (11.4,-4) node{\normalsize$i$};
	\draw (13.3,-4) node{\normalsize$i-1$};
	\draw[dots] (14.9,-4)--(16.6,-4);
	\draw (16.9,-4) node{\normalsize$1$};
	\redbraidbox{17.8}{19.9}{-4.8}{-3.3}{$\color{red}0\,\,0\color{black}$};
	\draw (20.8,-4) node{\normalsize$1$};
	\draw[dots] (21.6,-4)--(23,-4);
	\draw (24.7,-4) node{\normalsize$i-1$};
	\draw (26.8,-4) node{\normalsize$i$};
	\end{braid}.
\end{align*}
The first summand above is $0$ since the word $\ell\,(\ell-1)^2\,\cdots\,(i+1)^2\,(i-1)\,\cdots$ is non-cuspidal by Lemma~\ref{LCuspExpl}. So we get: 
\begin{align*}
&(-1)^{i+1}
\begin{braid}\tikzset{baseline=-.3em}
	\draw (0,4) node{\color{blue}\normalsize$\ell$\color{black}};
	\braidbox{1}{4.2}{3.2}{4.7}{$(\ell-1)^2$};
	\draw[dots] (5,4)--(6.5,4);
	\braidbox{7}{10.2}{3.2}{4.7}{$(i+1)^2$};
	\draw (11.4,4) node{\normalsize$i$};
	\draw (13.3,4) node{\normalsize$i-1$};
	\draw[dots] (14.9,4)--(16.6,4);
	\draw (16.9,4) node{\normalsize$1$};
	\redbraidbox{17.8}{19.9}{3.2}{4.7}{$\color{red}0\,\,0\color{black}$};
	\draw (20.8,4) node{\normalsize$1$};
	\draw[dots] (21.6,4)--(23,4);
	\draw (24.7,4) node{\normalsize$i-1$};
	\draw (26.8,4) node{\normalsize$i$};
	\draw[blue](0,3)--(0,-3);
	 \draw(1.5,3)--(1.5,-3);
	 \draw(3.5,3)--(3.5,-3);
	  \draw(7.5,3)--(7.5,-3);
	  \draw(9.5,3)--(9.5,-3);
	  \draw(11.4,3)--(11.4,-3);
	  \draw(13.3,3)--(13.3,-3);
	  \draw(16.9,3)--(16.9,-3);
	   \draw[red](18.4,3)--(18.4,-3);
	   \draw[red](19.3,3)--(19.3,-3);
\draw(20.9,3)--(20.9,-3);
\draw(24.7,3)--(24.7,-3);
\draw(26.8,3)--(24,0)--(26.8,-3);
		\draw (0,-4) node{\color{blue}\normalsize$\ell$\color{black}};
	\braidbox{1}{4.2}{-4.8}{-3.3}{$(\ell-1)^2$};
	\draw[dots] (5,-4)--(6.5,-4);
	\braidbox{7}{10.2}{-4.8}{-3.3}{$(i+1)^2$};
	\draw (11.4,-4) node{\normalsize$i$};
	\draw (13.3,-4) node{\normalsize$i-1$};
	\draw[dots] (14.9,-4)--(16.6,-4);
	\draw (16.9,-4) node{\normalsize$1$};
	\redbraidbox{17.8}{19.9}{-4.8}{-3.3}{$\color{red}0\,\,0\color{black}$};
	\draw (20.8,-4) node{\normalsize$1$};
	\draw[dots] (21.6,-4)--(23,-4);
	\draw (24.7,-4) node{\normalsize$i-1$};
	\draw (26.8,-4) node{\normalsize$i$};
	\end{braid},
\end{align*}
which equals 
$$
(-1)^{i+1}e^i(y_{p-1}-y_p)e^i=(-1)^i(-\eta^i_{i-1}+e^iy_pe^i)=(-1)^i(z^i+e^iy_pe^i)=c^i,
$$
where we have used (\ref{ERelB3}) for the penultimate equality. 
On the other hand, using the defining relations (\ref{R1})-(\ref{R7}), we get 
\begin{align*}
a^{i,i+1}a^{i+1,i}&=(-1)^{i+1}e^i\psi_{p-2i-3}\cdots\psi_{p-2}\psi_{p-1}e^{i+1}\psi_{p-1}\psi_{p-2}\cdots\psi_{p-2i-2}e_i
\\
&=(-1)^{i+1}e^i\psi_{p-2i-3}\cdots\psi_{p-2}\psi_{p-1}\psi_{p-1}\psi_{p-2}\cdots\psi_{p-2i-2}e_i
\\
&=(-1)^{i+1}e^i\psi_{p-2i-3}\cdots\psi_{p-2}(-y_{p-1}+y_p)\psi_{p-2}\cdots\psi_{p-2i-2}e_i
\\
&=(-1)^{i}e^i\psi_{p-2i-3}(y_{p-2i-2}-y_p)\psi_{p-2i-2}^2e_i
\\
&=(-1)^{i}e^i\psi_{p-2i-3}(y_{p-2i-2}-y_p)(-y_{p-2i-2}+y_{p-2i-1})e_i
\\
&=(-1)^{i+1}e^i(y_{p-2i-3}+y_{p-2i-2})e^i+(-1)^ie^iy_{p-2i-1}e_i+(-1)^ie^iy_pe^i
\\
&=(-1)^{i+1}\la^i_{i+1}+(-1)^i\xi^i_i+(-1)^ie^iy_pe^i
\\
&=(-1)^iz^i+(-1)^ie^iy_pe^i
\\
&=c^i,
\end{align*}
where we have used  (\ref{ERelB2}) and (\ref{ERelB3}) for the penultimate equality. 
\end{proof}

\begin{Lemma} \label{LAU0} 
We have $a^{1,0}u=0=ua^{0,1}$.
\end{Lemma}
\begin{proof}
Using a dot-crossing relation we get
$$
a^{1,0}u=e^1\psi_{p-1}\psi_{p-2}e^0y_pe^0
=e^1y_pe^1\psi_{p-1}\psi_{p-2}e^0-e^1\psi_{p-1}\psi_{p-2}y_{p-1}e^0e^0.
$$
The first summand is $\la^1_0a^{1,0}$ which is $0$ using (\ref{ERelB4}). To see that the second summand is also $0$, note that we can write $e^1=X\psi_{p-2}y_{p-1}$ for some $X$, so
\begin{align*}
e^1\psi_{p-1}\psi_{p-2}y_{p-1}e^0e^0
&=X\psi_{p-2}y_{p-1}\psi_{p-1}\psi_{p-2}y_{p-1}e^0
\\
&=X\psi_{p-2}\psi_{p-1}\psi_{p-2}y_{p}y_{p-1}e^0
\\
&=X\psi_{p-1}\psi_{p-2}\psi_{p-1}y_{p}y_{p-1}e^0
\\
&=X\psi_{p-1}\psi_{p-2}y_{p-1}y_{p}\psi_{p-1}e^0
\\
&=0,
\end{align*}
since $\psi_{p-1}e^0=0$ by Lemma~\ref{LPsiE}.
\end{proof}

\begin{Lemma} \label{LAAA} 
We have  $a^{i,i-1}a^{i-1,i}a^{i,i-1}=0$ for all $1\leq i\leq \ell-1$ and $a^{i,i+1}a^{i+1,i}a^{i,i+1}=0$ for all $0\leq i\leq \ell-2$.
\end{Lemma}
\begin{proof}
We prove the first equality, the proof of the second one is similar. For $i=1$, by Lemmas~\ref{LLoopU^2} and \ref{LAU0}, we have $a^{1,0}a^{0,1}a^{1,0}=a^{1,0}u^2=0$. Let $i>1$. By Lemma~\ref{LAAC}
\begin{align*}
a^{i,i-1}a^{i-1,i}a^{i,i-1}&=a^{i,i-1}c^{i-1}
=\pm a^{i,i-1}(z^{i-1}+y_pe^{i-1})
=\pm(z^{i}+\eta^i_{i-1})a^{i,i-1}=0,
\end{align*}
thanks to (\ref{ERelB3}). 
\end{proof}

\begin{Lemma} \label{LZCentral} 
Set $z:=\sum_{i\in J}z^i\in B$. Then $zu=-uz$ and $z$ commutes with all $e^i,c^i$ and $a^{i,j}$. 
\end{Lemma}
\begin{proof}
The equality $zu=-uz$ follows from Lemma~\ref{LNH2}(iv). 
By the dot-crossing relations (\ref{R5}), for all $1\leq i\leq \ell-1$, we have $z^ia^{i,i-1}=a^{i,i-1}z^{i-1}$ and 
$z^{i-1}a^{i-1,i}=a^{i-1,i}z^{i}$. In particular $z$ commutes with $a^{i,i-1}$ and $a^{i-1,i}$. Now, $z$ commutes with all $c^i$ by Lemma~\ref{LAAC} and with all $e^i$ by definition. 
\end{proof}

\subsection{Main results for $d=1$}
\label{SSMainD=1}
In this subsection we prove that $B\cong \Zig_\ell[\zz]$ and $Y_{\rho,1}^{\La_0}\cong \Zig_\ell$. 

\begin{Proposition} \label{PSurjCuspd=1} 
There exists a surjective homomorphism of graded superalgebras 
$$f: \Zig_\ell[\zz]\to B,\ \ze^i\mapsto e^i,\ \zu\mapsto u,\ \za^{i,j}\mapsto a^{i,j},\ \zc^i\mapsto c^i,\ \ze^i \zz\mapsto z^i.
\index{f@$f$}$$
\end{Proposition}
\begin{proof}
In view of the definition of $\Zig_\ell[\zz]$ and the defining relations (\ref{RZig1})-(\ref{RZig5}) for $\Zig_\ell$, the desired homomorphism $f$ exists in view of the relations checked in Lemmas~\ref{LAA0}, \ref{LLoopU^2}, \ref{LAAC}, \ref{LAU0}, \ref{LAAA} and \ref{LZCentral}. The homomorphism is surjective by Corollary~\ref{CBGen}. 
\end{proof}

We now use graded dimensions to prove that the homomorphism $f$ from Proposition~\ref{PSurjCuspd=1} is an isomorphism: 

\begin{Theorem} \label{TCuspIsoRank1} 
The homomorphism  $f: \Zig_\ell[\zz]\to B$ is an isomorphism. 
\end{Theorem}

\begin{proof}
Recall from (\ref{ECTrRes}) that we have surjective homomorphisms $\Om^{N\La_0}_{\rho,1}:B\to Y^{N\La_0}_{\rho,1}$ for all $N\in\Z_{\geq 1}$. So, using Proposition~\ref{PDimd=1Level N}, we get   
$$\dim_q B\geq \dim_qY^{N\La_0}_{\rho,1}=(1+q^4+\dots+q^{4(N-1)})\dim_q\Zig_\ell$$ 
for all $N$, where we have use the notation (\ref{EGDinLeq}). 
Hence 
$$\dim_q B\geq \frac{1}{1-q^4}
\dim_q\Zig_\ell=\dim_q \Zig_\ell[\zz].$$
So the theorem follows from Proposition~\ref{PSurjCuspd=1}. 
\end{proof}

\begin{Theorem} 
The graded superalgebra $B$ is graded Morita superequivalent to the graded superalgebra $\bar R_\de$.
\end{Theorem}
\begin{proof}
Recall that $B=\ga \bar R_\de\ga$ is an idempotent truncation of $\bar R_\de$. So, by Corollary~\ref{cor:idmpt_Mor}, it suffices to show that $|\Irr(B)|\geq |\Irr(\bar R_\de)|<\infty$. By  Theorem~\ref{THeadIrr}(vi), we have $|\Irr(\bar R_\de)|=|\Par^\ell(1)|=\ell$.  On the other hand,  every irreducible $\Zig_\ell$-module can be considered as an irreducible $\Zig_\ell[\zz]$-module with $\zz$ acting as $0$. Moreover,  $B\cong\Zig_\ell[\zz]$ by Theorem~\ref{TCuspIsoRank1}. Since $|\Irr(\Zig_\ell)|=\ell$, we deduce that  $|\Irr(B)|\geq \ell$ completing the proof. 
\end{proof}

\begin{Theorem} 
The graded superalgebra $Y_{\rho,1}^{\La_0}$ is graded Morita superequivalent to the RoCK block $R^{\La_0}_{\cont(\rho)+\de}$.
\end{Theorem}
\begin{proof}
Recall from (\ref{LIdTrCentd=1}) that $Y_{\rho,1}^{\La_0}$ is an idempotent truncation of $R^{\La_0}_{\cont(\rho)+\de}$. 
So, 
by Corollary~\ref{cor:idmpt_Mor}, it suffices to show that $|\Irr(Y_{\rho,1}^{\La_0})|\geq |\Irr(R^{\La_0}_{\cont(\rho)+\de})|<\infty$. 

By  Lemma~\ref{LNumberOfIrrCyc}, the amount of irreducible $|\Irr(R^{\La_0}_{\cont(\rho)+\de})|=|\Par^\ell(1)|=\ell$.  

On the other hand, recalling (\ref{ECTrRes}) and using Theorem~\ref{TCuspIsoRank1}, we have a surjective homomorphism $\Om^{\La_0}_{\rho,1}\circ f:\Zig_\ell[\zz]\to Y_{\rho,1}^{\La_0}$. Moreover, by  Proposition~\ref{PDimd=1Level N}, $\dim_q \ga^iY_{\rho,1}^{\La_0}\ga^j=\dim_q e^i\Zig_\ell e^j$ for all $i,j\in J$. We deduce that $Y_{\rho,1}^{\La_0}$ is non-negatively graded, and for every $i\in J$, the degree zero component of  $\ga^iY_{\rho,1}^{\La_0}\ga^i$ is $1$-dimensional and is spanned by $\ga^i$. It follows that the commutative algebra $\bigoplus_{i\in J}\F$ is a quotient of the algebra $Y_{\rho,1}^{\La_0}$, so  $|\Irr(Y_{\rho,1}^{\La_0})|\geq |J|=\ell$.
\end{proof}

Recall that a finite dimensional algebra $A$ is called a {\em $QF2$-algebra} if its projective, indecomposable modules have simple socles, see \cite[\S31]{AF}. This condition implies that $A$ is quasi-Frobenius, see \cite[Corollary 31.8]{AF}. 

\begin{Theorem} \label{TMaind=1} 
Let $\rho$ be a $1$-Rouquier $\bar p$-core and 
suppose that $R^{\La_0}_{\cont(\rho)+\de}$ is a $QF2$-algebra. Then the restriction of $\Om^{\La_0}_{\rho,1}\circ f$ to $\Zig_\ell$ is an isomorphism of graded superalgebras $\Zig_\ell\iso Y^{\La_0}_{\rho,1}$. 
\end{Theorem}
\begin{proof}
Recalling (\ref{ECTrRes}), by Theorem~\ref{TCuspIsoRank1}, we have a surjective homomorphism 
\begin{equation*}\label{E170421}
\psi:=\Om^{\La_0}_{\rho,1}\circ f:\Zig_\ell[\zz]\to Y_{\rho,1}^{\La_0}. 
\end{equation*}
By Proposition~\ref{PDimd=1Level N}, $\dim_q \ga^iY_{\rho,1}^{\La_0}\ga^j=\dim_q e^i\Zig_\ell e^j$ for all $i,j\in J$. So by degrees, $\psi(z^n)=0$ for all $n\geq 2$ and $\psi$ is injective upon restriction to $\spa(u,e^i,z^{i,j})$. Moreover, for every $i\in J$, we have that the degree $4$ component of $\ga^iY_{\rho,1}^{\La_0}\ga^i$ is $1$-dimensional and is spanned by $\psi(c^i)$ and $\psi(\ze^i \zz)$. It follows that $Y_{\rho,1}^{\La_0}$ is a basic algebra with projective, indecomposable modules $Y_{\rho,1}^{\La_0}\ga^i$ for $i\in J$. Let $L(i)$ be the corresponding irreducible $Y_{\rho,1}^{\La_0}$-module. 

If $\psi(c^i)\neq 0$ for all $i\in J$, the theorem follows. Otherwise there is $i\in J$ such that $\psi(c^i)=0$ and $\ga^iY_{\rho,1}^{\La_0}\ga^i$ is spanned by $\psi(\ze^i \zz)\neq 0$. It follows that the projective, indecomposable module $Y_{\rho,1}^{\La_0}\ga^i$ has basis 
$$
\left\{
\begin{array}{ll}
\{\psi(\ze^i),\psi(\za^{i-1,i}),\psi(\za^{i+1,i}),\psi(\ze^i \zz)\} &\hbox{if $1\leq i\leq \ell-2$},\\
\{\psi(\ze^{\ell-1}),\psi(\za^{\ell-2,i}),\psi(\ze^{\ell-1} \zz)\} &\hbox{if $i=\ell-1>0$},\\
\{\psi(\ze^{0}),\psi(\zu),\psi(\za^{1,0}),\psi(\ze^0 \zz)\} &\hbox{if $i=0<\ell-1$},\\
\{\psi(\ze^{0}),\psi(\zu),\psi(\ze^0 \zz)\} &\hbox{if $i=0=\ell-1$}.
\end{array}
\right.
$$
From this we deduce that the projective, indecomposable module $Y_{\rho,1}^{\La_0}\ga^i$ contains a submodule isomorphic to
$$
\left\{
\begin{array}{ll}
L(i-1)\oplus L(i+1)\oplus L(i)\} &\hbox{if $1\leq i\leq \ell-2$},\\
\{L(\ell-2)\oplus L(\ell-1) &\hbox{if $i=\ell-1>0$},\\
L(0)\oplus L(1)\oplus L(0) &\hbox{if $i=0<\ell-1$},\\
L(0)\oplus L(0) &\hbox{if $i=0=\ell-1$}.
\end{array}
\right.
$$
This contradicts the assumption that  $R^{\La_0}_{\cont(\rho)+\de}$ has the property $QF2$. 
\end{proof}

\begin{Remark}
We comment on the assumption in the theorem above that $R^{\La_0}_{\cont(\rho)+\de}$ is a $QF2$-algebra. We conjecture that, much more generally, any cyclotomic quiver Hecke superalgebra $R_\theta^\La$ is {\em symmetric}. For quiver Hecke algebras the analogous result is known \cite[Proposition 3.10]{SVV}, see also \cite[Remark 3.19]{Web}.
For quiver Hecke superalgebras, we also know this in many cases. Namely, let us say that the pair $(\cha \F,\ell)$ is bad if 
 $0<\cha \F<2\ell+1$ and $\cha \F$ divides $2\ell+1$. If $(\cha \F,\ell)$ is not bad, then $R^{\La}_{\theta}$ is symmetric, see Corollary~\ref{Cor:QF2}. 
\end{Remark}

\section{Gelfand-Graev truncation of \texorpdfstring{$\bar R_{d\de}$}{}}

Recall that we have fixed a convex preorder $\preceq$ on $\Phi_+$ as in Example~\ref{ExConPr}, and cuspidality always means cuspidality with respect to this preorder. Throughout the section we also fix $d\in\Z_{\geq 0}$ and a $d$-Rouquier $\bar p$-core $\rho$. From now on we only consider the case $N=1$, so we abbreviate $\Cent_{\rho,d}:=\Cent_{\rho,d}^{\La_0}$, $\Om_{\rho,d}:=\Om_{\rho,d}^{\La_0}$, etc. 
\index{z@$\Cent_{\rho,d}$}\index{o@$\Om_{\rho,d}$}\index{w@$\Om_{\rho,d}$}\index{$\Om_{\rho,d}$}
In this section, we will define an idempotent $\ga_{1^d}\in \bar R_{d\de}$ and 
begin to study the idempotent truncations $B_d:=\ga_{1^d}\bar R_{d\de}\ga_{1^d}$ and $Y_{\rho,d}:=\ga_{1^d}\Cent_{\rho,d}\ga_{1^d}$. These truncations are important since if $\cha \F> d$, we will eventually prove that $\bar R_{d\de}$ is Morita equivalent to $B_d$, and $\Cent_{\rho,d}$ is Morita equivalent to $Y_{\rho,d}$ and to $R^{\La_0}_{\cont(\rho)+d\de}$.

Recall that $\bar R_{d\de}$ is a quotient of $R_{d\de}$. The algebra  $R_{d\de}$ is generated by the elements $\psi_r,y_s,e(\bi)$, and we denote by the same symbols the corresponding elements in the quotient $\bar R_{d\de}$.  In particular, we have the elements $y_1,\dots,y_{dp},\psi_1,\dots,\psi_{dp-1}\in \bar R_{d\de}$.

\subsection{Higher Gelfand-Graev idempotents}
We recall the Gelfand-Graev words $\bg^i\in I^\de_\di$ and $\hat\bg^i\in I^\de$ defined in \S\ref{SSGG}. For a $d$-tuple $\bi=i_1\cdots i_d\in J^d$, we define the (higher) {\em Gelfand-Graev words} \index{Gelfand-Graev word (higher)}
to be concatenations 
\begin{equation*}\label{EGGWConcatBold}
\begin{split}
\bg^\bi&:=\bg^{i_1}\cdots \bg^{i_d}
\in I^{d\de}_\di,
\index{g@$\bg^\bi$}
\\
\hat\bg^\bi&:=\widehat{\bg^\bi}=\hat\bg^{i_1}\cdots \hat\bg^{i_d}
\in I^{d\de}.
\index{g@$\hatbg^\bi$}
\end{split}
\end{equation*}
The word $\hat\bg^\bi$, being a concatenation (hence a shuffle) of cuspidal words, is also cuspidal, see Lemma~\ref{LTensImagIsImag}(i). 
So it makes sense to consider the {\em Gelfand-Graev idempotents}
$$
\ga^\bi:=e(\bg^\bi)\in \bar R_{d\de} \quad \index{g@$\ga^\bi$}\index{c@$\ga^\bi$}
\text{and}\qquad
\ga_{1^d}:=\sum_{\bi\in J^d}\ga^{\bi}.
\index{g@$\ga_{1^d}$}\index{c@$\ga_{1^d}$}
$$
(There are more general Gelfand-Graev idempotents $\ga_\mu$ corresponding to all compositions $\mu$ of $d$, but only the case  $\mu=(1^d)$ will be needed in this paper). For $d=1$, $\ga_1$ is what we denoted $\ga$ in \S\ref{SSGG}. For $d=2$, we often write $\ga^{i,j}$ instead of $\ga^{ij}$. 
Recalling Lemma~\ref{Ltens2}, we use the same symbols
to denote the Gelfand-Graev idempotents in $\Cent_{\rho,d}$:
$$
\ga^\bi:=\Om_{\rho,d}(\ga^\bi)\in\Cent_{\rho,d}
\quad\text{and}\quad
\ga_{1^d}:=\Om_{\rho,d}(\ga_{1^d})\in\Cent_{\rho,d}.
$$

We consider the idempotent truncations 
\begin{equation}\label{ECTrBold}
B_{d}:=\ga_{1^d}\bar R_{d\de}\ga_{1^d}\qquad \text{and}\qquad
Y_{\rho,d}:=\ga_{1^d}\Cent_{\rho,d}\ga_{1^d}.
\index{b@$B_{d}$}\index{y@$Y_{\rho,d}$}
\end{equation}
For $d=1$ we get the algebras $B_1=B$ and $Y_{\rho,1}=Y_{\rho,1}^{\La_0}$ from (\ref{ECTr}). Restricting $\Om_{\rho,d}$ to the idempotent truncations, we get a surjective algebra homomorphism
\begin{equation}\label{ETrOmd}
\Om_{\rho,d}:B_d\to Y_{\rho,d}.
\end{equation}

By Lemma~\ref{LIdSym}, we have that $Y_{\rho,d}$ is an idempotent truncation of $R_{\cont(\rho)+d\de}^{\La_0}$:
\begin{equation}\label{EIdTrd}
Y_{\rho,d}\cong \ga_{1^d}\eps_d R_{\cont(\rho)+d\de}^{\La_0} \eps_d\ga_{1^d}.  
\end{equation}

Recall Lemma~\ref{L030216} and abbreviate $\de^d:=(1^d)\de$. By the lemma, we have the parabolic subalgebra 
$
\bar R_\de^{\otimes d}\cong \bar R_{\de^d}\subseteq e_{\de^d}\bar R_{d\de}e_{\de^d}.
$
Note that here and below all tensor products are {\em tensor products of superalgebras}. 
Upon truncation with $\ga_{1^d}$ we get the parabolic subalgebra 
\begin{equation}\label{EBPar}
B^{\otimes d} \cong B_{1^d}:=\ga_{1^d}\bar R_{\de^d}\ga_{1^d}\subseteq B_d.
\index{b@$B_{1^d}$}
\end{equation}
We identify the subalgebra $B^{\otimes d}$ with $B_{1^d}$. Given an element $b\in B$ and $1\leq r\leq d$, we now denote
\begin{equation}\label{EAR}
b_r:=\ga^{\otimes (r-1)}\otimes b\otimes \ga^{\otimes(d-r)}\in B^{\otimes d}=B_{1^d}.
\end{equation}
This agrees with (\ref{EInsertion}). We 
have an injective algebra homomorphism 
\begin{equation}\label{EIotaR}
\iota_r:B\to B^{\otimes d}=B_{1^d}\subseteq B_d, \ b \mapsto b_r \qquad(1\leq r\leq d).
\end{equation}
Recalling the notation (\ref{EEIFormula}) - (\ref{EAI-1IFormula}), we now have the elements 
$$e^i_r,a^{i,j}_r,c^i_r,c_r,z^i_r,z_r\in B_{1^d}.
\index{e@$e^i_r$}\index{a@$a^{i,j}_r$}\index{c@$c^i_r$}\index{c@$c_r$}\index{z@$z^i_r$}\index{z@$z_r$}
$$

Recalling (\ref{EZetaHom}), we have a unital algebra homomorphism
$$
\zeta_{\cont(\rho)+\de,(d-1)\de}^{\La_0}:R^{\La_0}_{\cont(\rho)+\de}\to R^{\La_0}_{\cont(\rho)+d\de}.
$$
It is easy to see that upon restriction, this homomorphism yields a unital algebra homomorphism
\begin{equation}\label{EZeta}
\zeta:Y_{\rho,1}=\ga_1\eps_1R^{\La_0}_{\cont(\rho)+\de}\eps_1\ga_1\to \ga_{1^d}\eps_dR^{\La_0}_{\cont(\rho)+d\de}\eps_d\ga_{1^d}=Y_{\rho,d}
\index{z@$\zeta$}
\end{equation}
so that:

\begin{Lemma} \label{LY_1Y_d} 
The following diagram commutes:
$$
\begin{tikzcd}
  B_1 \arrow[r, "\iota_1"] \arrow[d, "\Om_{\rho,1}"]
    & B_d \arrow[d, "\Om_{\rho,d}"] \\
  Y_{\rho,1} \arrow[r, "\zeta" ]
&  Y_{\rho,d}  \end{tikzcd}
$$
\end{Lemma}

By Theorem~\ref{TCuspIsoRank1} and (\ref{EBPar}), we have $B_{1^d}\cong \Zig_\ell[\zz]^{\otimes d}$, so as vector spaces:
\begin{equation}\label{EParDecB}
B_{1^d}=\F[z_1,\dots,z_d]\otimes \Zig_\ell^{\otimes d}\ =\,\bigoplus_{a_1,\dots,a_d\in\Z_{\geq 0}}z_1^{a_1}\cdots z_d^{a_d} \,\Zig_\ell^{\otimes d}.
\end{equation}
For $n\in\Z_{\geq 0}$, set
\begin{equation}\label{EBNPar}
B_{1^d}^{(n)}:=\bigoplus_{a_1+\dots+a_d\leq n}z_1^{a_1}\cdots z_d^{a_d} \,\Zig_\ell^{\otimes d}\subseteq B_{1^d}.
\end{equation}
By (\ref{ETwistedRelations}), we have for all $n,m\in\Z_{\geq 0}$:
\begin{equation}\label{EFiltrPar}
B_{1^d}^{(n)}B_{1^d}^{(m)}\subseteq B_{1^d}^{(n+m)}.
\end{equation}

Let now $1\leq r<d$. By Lemma~\ref{L030216} again, denoting 
$(\de^{r-1},2\de,\de^{d-r-1}):=(1^{r-1},2,1^{d-r-s})\de,$  
we have a parabolic subalgebra 
$$
\bar R_\de^{\otimes (r-1)}\otimes \bar R_{2\de}\otimes \bar R_\de^{(d-r-1)}\cong \bar R_{\de^r,2\de,\de^{d-r-1}}\subseteq e_{\de^r,2\de,\de^{d-r-1}}\bar R_{d\de}e_{\de^r,2\de,\de^{d-r-1}}.
$$
Upon truncation with $\ga_{1^d}$ we get the parabolic subalgebra
$$
B^{\otimes (r-1)}\otimes B_2\otimes B^{\otimes (d-r-1)}\cong B_{1^r,2,1^{d-r-1}}:=\ga_{1^d}\bar R_{\de^r,2\de,\de^{d-r-1}}\ga_{1^d}\subseteq B_d.
$$
For $b\in B_2$ we thus have an algebra embedding
\begin{equation}\label{EIotaRR+1}
\iota_{r,r+1}:B_2\to B_d,\ b\mapsto \ga^{\otimes (r-1)}\otimes b\otimes \ga^{\otimes (d-r-1)}.
\end{equation}

\subsection{Dimension formula for $Y_{\rho,d}$}
\label{SSDimYRhoD}

We need a few preliminary results before stating and proving our main dimension formula of this subsection. If $\lambda\in\Par_0$ is a strict partition, we denote $K_\lambda':=|\Std_0(\la)|$\index{k@$K_\lambda'$}, the number of strict $\la$-tableaux, see \S\ref{SRemAdd}. On the other hand, if  $\lambda\in\Par$ is any (not necessarily strict) partition, we denote by $K_\lambda$ the number of the usual standard $\la$-tableaux, so that $K_\la$\index{k@$K_\lambda$} equals the dimension of the Specht module corresponding to the partition $\la$, see \cite[3.1.13]{JK}. This immediately implies the first equality in the following lemma:

\begin{Lemma}\label{lem:dim_sqd}
Let $n\in\Z_{\geq 0}$. Then
$$
\sum_{\lambda\in\Par(n)}K_\lambda^2=n!=\sum_{\lambda\in\Par_0(n)}2^{n-h(\lambda)}(K_\lambda')^2.
$$
\end{Lemma}

\begin{proof}
In~\cite[(7.2)]{Stem} we have that
\begin{align*}
\varphi^\lambda(1)=\frac{1}{\varepsilon_\lambda}2^{-(n+h(\lambda))/2}2^n K_\lambda',
\end{align*}
for each $\lambda\in\Par_0(n)$, where $\varphi^\lambda$ is one of the (possibly two) spin characters of $\tilde{\Si}_n$ labeled by $\lambda$ and $\varepsilon_\lambda$ takes the value $1$ or $\sqrt{2}$ depending on whether $n-h(\lambda)$ is even or odd. (Note that in~\cite{Stem}, $K_{\lambda,(1^n)}'$ counts the number of ways of filling in a shifted $\lambda$-tableau but each box can be filled in with either an $i$ or $i'$, hence the factor of $2^n$ when compared to our $K_\la'$.) 

Now $\varphi^\lambda$ labels one or two spin characters of $\tilde{\Si}_n$ if $n-h(\lambda)$ is even or odd respectively. Therefore,
\begin{align*}
n!=\sum_{\lambda\in\Par_0(n)}\varepsilon_\lambda^2\varphi^\lambda(1)^2=\sum_{\lambda\in\Par_0(n)}2^{n-h(\lambda)}(K_\lambda')^2,
\end{align*}
proving the second equality. 
\end{proof}

Recall the notion of the $\bar p$-quotient from \S\ref{SSAb}.
For $\la\in\Par_p(\rho,d)$, we always denote $\quot(\la)=(\la^{(0)},\dots,\la^{(\ell)})$. Recall that $\Par_0(\rho,d)$ denoted the set of all strict partitions $\la$ with $\core(\la)=\rho$ and $\wt(\la)=d$. 

\begin{Lemma} \label{L110921_2} 
We have 
$$\sum_{\la\in\Par_0(\rho,d)}
2^{2d-h(\la^{(0)})-2|\la^{(\ell)}|}\binom{d}{|\la^{(0)}|,\dots,|\la^{(\ell)}|}^2\,(K_{\la^{(0)}}')^2\,\prod_{j=1}^\ell K_{\la^{(j)}}^2
=d!(2p-3)^d.
$$
\end{Lemma}
\begin{proof}
Note that $\la\in\Par_p(\rho,d)$ is strict if and only if $\la^{(0)}$ is. So, in view of Theorem~\ref{TMY}, $\la\mapsto \quot(\la)$ is a bijection between $\Par_0(\rho,d)$ and the multipartitions $(\la^{(0)},\dots,\la^{(\ell)})$ of $d$ such that $\la^{(0)}\in\Par_0$. So,  
\begin{align*}
&\sum_{\la\in\Par_0(\rho,d)}
2^{2d-h(\la^{(0)})-2|\la^{(\ell)}|}\binom{d}{|\la^{(0)}|,\dots,|\la^{(\ell)}|}^2\,(K_{\la^{(0)}}')^2\,\prod_{j=1}^\ell K_{\la^{(j)}}^2
\\
=\,
&\sum_{d_0+d_1+\dots+d_\ell=d}\sum_{\substack{\la^{(0)}\in\Par_0(d_0)\\ \la^{(1)}\in\Par(d_1)\\ \vdots\\ \la^{(\ell)}\in\Par(d_\ell)}}
2^{2d-h(\la^{(0)})-2d_\ell}\binom{d}{d_0,\dots,d_\ell}^2\,(K_{\la^{(0)}}')^2\,\prod_{j=1}^\ell K_{\la^{(j)}}^2
\end{align*}
\begin{align*}
=\,
&\sum_{d_0+d_1+\dots+d_\ell=d}2^{2d-d_0-2d_\ell}\binom{d}{d_0,\dots,d_\ell}^2\,d_0!\,\prod_{j=1}^\ell d_j!
\\
=\,
&d!\sum_{d_0+d_1+\dots+d_\ell=d}2^{2d-d_0-2d_\ell}\binom{d}{d_0,\dots,d_\ell}
\\
=\,
&d!\sum_{d_0+d_1+\dots+d_\ell=d}2^{d_0}\cdot 4^{d_1}\,\cdots \,4^{d_{\ell-1}}\cdot 1^{d_\ell}\binom{d}{d_0,\dots,d_\ell}
\\
=\,
&d!(2+4(\ell-1)+1)^d
\\
=\,
&d!(2p-3)^d,
\end{align*}
where we have used Lemma~\ref{lem:dim_sqd} for the second equality.
\end{proof}

\begin{Definition}\label{DRhoLa}
{\rm 
Let $\la\in\Par_0(\rho,d)$. 
A {\em strict $(\rho,\la)$-sequence}\index{strict $(\rho,\la)$-sequence} is a sequence $(\mu^0,\mu^1,\dots,\mu^d)$ of partitions such that 
$$
\mu^0=\rho\subset \mu^1\subset \dots\subset \mu^{d-1}\subset \mu^d=\la
$$
and $\mu^c\in\Par_0(\rho,c)$ for all $c=0,\dots,d$. 
}
\end{Definition}

We denote by $\Seq(\rho,\la)$\index{s@$\Seq(\rho,\la)$} the set of all strict $(\rho,\la)$-sequences. 

Let $(\mu^0,\mu^1,\dots,\mu^d)\in \Seq(\rho,\la)$. 
By Lemma~\ref{ChuKesLem}(i), for $c=1,\dots,d$, we have that the abacus $\Ab_{\mu^c}$ is be obtained from $\Ab_{\mu^{c-1}}$ by performing a strict elementary slide down on a runner $i_c\in I$. In this case we say that the $(\rho,\la)$-sequence $(\mu^0,\mu^1,\dots,\mu^d)$ has color sequence $\bi=i_1\dots i_d$. We denote by $\Seq_\bi(\rho,\la)$\index{s@$\Seq_\bi(\rho,\la)$} the set of all $(\rho,\la)$-sequences with color sequence $\bi=i_1\dots i_d$. Then
$$
\Seq(\rho,\la)=\bigsqcup_{\bi\in I^d}\Seq_\bi(\rho,\la).
$$

\begin{Lemma}\label{lem:bead_slide_count}
Let $\la\in\Par_0(\rho,d)$ with $\quot(\la)=(\la^{(0)},\dots\la^{(\ell)})$. Denote $d_i:=|\la^{(i)}|$ for all $i\in I$ and set $\eta:=\sum_{i\in I}d_i\al_i\in Q_+$. 
Then 
$
\Seq_\bi(\rho,\la)\neq \varnothing$ if and only if $\bi\in I^\eta$, in which case 
$$
|\Seq_\bi(\rho,\la)|=K_{\la^{(0)}}'\prod_{j=1}^\ell K_{\la^{(j)}}.
$$
In particular, 
$$
|\Seq(\rho,\la)|=
\binom{d}{d_0,\dots,d_\ell}\,K_{\la^{(0)}}'\,\prod_{j=1}^\ell K_{\la^{(j)}}.
$$
\end{Lemma}
\begin{proof}
If $(\mu^0,\mu^1,\dots,\mu^d)\in \Seq_\bi(\rho,\la)$, then $\Ab_\bla$ is obtained from $\Ab_\rho$ by $d$ (strict) consecutive slides down on runners $i_1,\dots,i_d$, so for all $i\in I$ we must have $\sharp\{k\mid 1\leq k\leq d\ \text{and}\ i_k=i\}=d_i$, i.e. $\bi\in I^\eta$. Moreover, $|I^\eta|=\binom{d}{d_0,\dots,d_\ell}$, so it remains to prove that for any $\bi\in I^\eta$, we have $
|\Seq_\bi(\rho,\la)|=K_{\la^{(0)}}'\prod_{j=1}^\ell K_{\la^{(j)}}.
$

Let $i\in I$. A sequence of $d_i$ consecutive strict elementary slides down the $i^{\nth}$ runner of $\Ab_{\rho}$ that leads to the arrangement of beads on the $i^{\nth}$ runner of $\Ab_{\la}$ will be called an {\em $i$-good sequence}. Denote by $S_i$ the set of all $i$-good sequences. Since slides down are performed on different runners independently of each other, 
it suffices to prove that 
$$
|S_i|=\left\{
\begin{array}{ll}
K_{\la^{(0)}}' &\hbox{if $i=0$},\\
K_{\la^{(i)}} &\hbox{if $1\leq i\leq \ell$.}
\end{array}
\right.
$$

If $1\leq i\leq \ell$ there is an obvious bijection between the elements of $S_i$ and the set of the standard $\la^{(i)}$-tableaux (see~\cite[p. 252]{Tu}), so $|S_i|=K_{\la^{(i)}}$. Similarly, there is an obvious bijection between $S_0$ and the set $\Std_0(\la)$ of the strict standard $\la^{(0)}$-tableaux, so $|S_0|=K_{\la^{(0)}}'$. 
\end{proof}

Let $\la\in\Par_0(\rho,d)$ and denote $r:=|\rho|$. For $\U\in\Std_p(\rho)$ and $\bi\in J^d$, we denote:
$$
\Std^{\U,\bi}_0(\la):=\{\Stab\in\Std_0(\la)\mid \Stab_{\leq r}=\U\ \text{and}\ \bi^{\Stab}=\bi^\U\hat\bg^\bi\}.
\index{s@$\Std^{\U,\bi}_0(\la)$}
$$

\begin{Lemma} \label{L110921} 
Let $\la\in\Par_0(\rho,d)$ with $\quot(\la)=(\la^{(0)},\dots,\la^{(\ell)})$ and $\U\in\Std_p(\rho)$. Then
$$
\sum_{\bi\in J^d}\sum_{\Stab\in \Std^{\U,\bi}_0(\la)}2^{\sum_{s=1}^d (i_s-\ell)}=2^{-d_\ell}\binom{d}{d_0,\dots,d_\ell}\,K_{\la^{(0)}}'\,\prod_{j=1}^\ell K_{\la^{(j)}},
$$
where $d_i:=|\la^{(i)}|$ for all $i\in I$.
\end{Lemma}
\begin{proof}
Set $r:=|\rho|$. 

If $\Stab\in\Std^{\U,\bi}_0(\la)$ for some $\bi\in J^d$, then, setting  
$\mu^c:=\Stab(\{1,\dots,r+cp\})$ for $c=0,\dots, d$, we note by Lemma~\ref{LMichael}, that $\umu:=(\mu^0,\dots,\mu^d)$ is a strict $(\rho,\la)$-sequence as in  Definition~\ref{DRhoLa}. In this case we write $\Stab\in\Std^{\U,\bi}_0(\umu)$.  Therefore  
$$
\sum_{\bi\in J^d}\sum_{\Stab\in \Std^{\U,\bi}_0(\la)}2^{\sum_{s=1}^d (i_s-\ell)}=
\sum_{\umu\in\Seq(\rho,\la)}\sum_{\bi\in J^d}\sum_{\Stab\in \Std^{\U,\bi}_0(\umu)}2^{\sum_{s=1}^d (i_s-\ell)},
$$
so, in view of Lemma~\ref{lem:bead_slide_count}, it suffices to prove that for every $\umu\in\Seq(\rho,\la)$, we have 
$$
\sum_{\bi\in J^d}\sum_{\Stab\in \Std^{\U,\bi}_0(\umu)}2^{\sum_{s=1}^d (i_s-\ell)}=2^{-d_\ell}.
$$

Fix $\umu=(\mu^0,\dots,\mu^d)\in\Seq(\rho,\la)$, and let $j_1\cdots j_d$ be the color sequence of $\umu$. By Corollary~\ref{CSkewShuffle}, for every $c=1,\dots,d$, we have 
$$
i_c=
\left\{
\begin{array}{ll}
0 &\hbox{if $j_c=0$,}\\
\ell-1&\hbox{if $j_c=\ell$}\\
j_c\ \text{or}\ j_{c-1}&\hbox{if $1\leq j_c\leq\ell-1$.}
\end{array}
\right.
$$ 

Let $1\leq c\leq d$. 

If $j_c=0$ then $i_c=0$ and, by Remark~\ref{R010921}, 
$$
\mu^c\setminus\mu^{c-1} =
\vcenter{\hbox{
\ytableausetup{mathmode}
\begin{ytableau}
\none[\,] & \none & \none & {\scriptstyle\ell} & {\scriptstyle\ell-1} & \cdots &{\scriptstyle 1} & {\scriptstyle 0} \cr
\none & \none & {\scriptstyle\ell-1} \cr
\none & \none[\iddots] \cr
{\scriptstyle 0} \cr
\end{ytableau}}}
$$
From this it is easy to see that there are $2^\ell$ options for filling these boxes when constructing a tableau $\Stab\in \Std^{\U,\bi}_0(\umu)$. 

If $j_c=\ell$ then $i_c=\ell-1$ and, by Remark~\ref{R010921}, 
$$
\mu^c\setminus\mu^{c-1} = 
\vcenter{\hbox{
\ytableausetup{mathmode}
\begin{ytableau}
{\scriptstyle\ell} & {\scriptstyle\ell-1} & \cdots &{\scriptstyle 1} & {\scriptstyle 0} & {\scriptstyle 0} & {\scriptstyle 1} & \dots & {\scriptstyle\ell-1}\cr
\end{ytableau}}}.
$$
From this it is easy to see that there is only one option for filling these boxes when constructing a tableau $\Stab\in \Std^{\U,\bi}_0(\umu)$. 

If $j:= j_c\in\{ 1,\dots,\ell-1\}$ then $i_c=j$ or $i_c=j-1$, and, by Remark~\ref{R010921}, 
$$
\mu^c\setminus\mu^{c-1} =
\vcenter{\hbox{
\ytableausetup{mathmode}
\begin{ytableau}
\none[\,] & \none & \none & {\scriptstyle\ell} & {\scriptstyle\ell-1} & \cdots &{\scriptstyle 1} & {\scriptstyle 0} & {\scriptstyle 0} & {\scriptstyle 1} & \cdots & {\scriptstyle j-1} \cr
\none & \none & {\scriptstyle\ell-1} \cr
\none & \none[\iddots] \cr
{\scriptstyle j} \cr
\end{ytableau}}}
$$
If $i_c=j$ then there are $2^{\ell-1-j}$ options for filling these boxes when constructing a tableau $\Stab\in \Std^{\U,\bi}_0(\umu)$, and if $i_c=j-1$ then there are $2^{\ell-j}$ options.

Thus,
\begin{align*}
\sum_{\bi\in J^d}\sum_{\Stab\in \Std^{\U,\bi}_0(\umu)}2^{\sum_{s=1}^d (i_s-\ell)}&=(2^\ell\, 2^{-\ell})^{d_0}(1\cdot 2^{-1})^{d_\ell}\prod_{j=1}^{\ell-1}(2^{\ell-j}\, 2^{j-1-\ell}+2^{\ell-1-j}\, 2^{j-\ell})^{d_j},
\end{align*}
which is easily seen to equal to $2^{-d_\ell}$, as desired. 
\end{proof}

\begin{Proposition}\label{PDimYd}
We have $\dim Y_{\rho,d}=d!(2p-3)^d = \dim (\Zig_\ell\swr \Si_d).$
\end{Proposition}

\begin{proof}
The second equality is clear since $\dim \Zig_\ell=2p-3$ by  (\ref{EDimAL}). Denote 
$$
\al:=\cont(\rho)=\sum_{i\in I}m_i^\al\al_i\quad\text{and}\quad
\theta:=\cont(\rho)+d\de=\sum_{i\in I}m_i^\theta\al_i.
$$
Note that 
\begin{equation}\label{E110921}
m_0^\theta-m_0^\al=2d.
\end{equation} 

Let 
$$\ttf:=\sum_{\bi\in J^d}\iota_{\al,d\de}(e_\al\otimes e(\bg^\bi))
=\sum_{\bk\in I^{\al}}
\sum_{\bi\in J^d } e(\bk\bg^\bi)
\in R^{\La_0}_{\theta}.$$
By Lemma~\ref{L270116New}, we have
\begin{align*}
R_{\al}^{\La_0} \otimes Y_{\rho,d} = R_{\al}^{\La_0} \otimes \ga_{1^d}\Cent_{\rho,d}\ga_{1^d} \cong \ttf R_{\theta}^{\La_0}\ttf.
\end{align*}
So it suffices to prove that
\begin{align*}
\dim \ttf R^{\Lambda_0}_{\theta}\ttf  = d!(2p-3)^d\,\dim R^{\Lambda_0}_{\al}.
\end{align*}


Let $\bk,\bl\in I^{\al}$ and $\bi=i_1\cdots i_d,\bj=j_1\cdots j_d\in J^d$. By Lemma~\ref{LFact}, we have 
\begin{align*}
\dim e(\bk\bg^\bi) R^{\Lambda_0}_{\theta} e(\bl\bg^\bj) 
&=
2^{\sum_{s=1}^d -(\ell-i_s)}2^{\sum_{t=1}^d -(\ell-j_t)}
\dim e(\bk\hat\bg^\bi) R^{\Lambda_0}_{\theta} e(\bl\hat\bg^\bj).
\end{align*}
So 
\begin{align*}
\dim \ttf R^{\Lambda_0}_{\theta}\ttf
&=\sum_{\bk,\bl\in I^{\al}}\sum_{\bi,\bj\in J^d } \dim e(\bk\bg^\bi) R^{\Lambda_0}_{\theta} e(\bl\bg^\bj)\\
&=
\sum_{\bk,\bl\in I^{\al}}
\sum_{\bi,\bj\in J^d } 
2^{\sum_{s=1}^d (i_s-\ell)}2^{\sum_{t=1}^d (j_t-\ell)}
\dim e(\bk\hat\bg^\bi) R^{\Lambda_0}_{\theta} e(\bl\hat\bg^\bj)
\end{align*}
\begin{align*}
&=
\sum_{\bk,\bl\in I^{\al}}
\sum_{\bi,\bj\in J^d } 
2^{\sum_{s=1}^d (i_s-\ell)}2^{\sum_{t=1}^d (j_t-\ell)}
\sum_{\la\in\Par_0}\sum_{\substack{\Stab\in\Std_0(\la,\bk\hat\bg^\bi)\\ \T\in\Std_0(\la,\bl\hat\bg^\bj)}} 2^{m_0^\theta-h(\la)}
\\
&=
\sum_{\bk,\bl\in I^{\al}}
\sum_{\la\in\Par_0}
2^{m_0^\theta-h(\la)}
\sum_{\bi,\bj\in J^d } 
\sum_{\substack{\Stab\in\Std_0(\la,\bk\hat\bg^\bi)\\ \T\in\Std_0(\la,\bl\hat\bg^\bj)}}
2^{\sum_{s=1}^d (i_s-\ell)}2^{\sum_{t=1}^d (j_t-\ell)},
\end{align*}
where we have used Corollary~\ref{CDimAP} for the third  equality. Note by Lemma~\ref{LMichael} that $\Stab_{\leq r},\T_{\leq r}\in \Std_0(\rho)$ whenever $\Stab\in\Std_0(\la,\bk\hat\bg^\bi)$ and $\T\in\Std_0(\la,\bl\hat\bg^\bj)$, so we can rewrite the last line as  
\begin{align*}
&\sum_{\bk,\bl\in I^{\al}}
\sum_{\substack{\U\in\Std_0(\rho,\bk)\\ \V\in\Std_0(\rho,\bl)}}
\sum_{\la\in\Par_0(\rho,d)}
2^{m_0^\theta-h(\la)}
\sum_{\bi,\bj\in J^d } 
\sum_{\substack{\Stab\in\Std_0^{\U,\bi}(\la)\\ \T\in\Std_0^{\V,\bj}(\la)}}
2^{\sum_{s=1}^d (i_s-\ell)}2^{\sum_{t=1}^d (j_t-\ell)}
\\
=\,
&\sum_{\bk,\bl\in I^{\al}}
\sum_{\substack{\U\in\Std_0(\rho,\bk)\\ \V\in\Std_0(\rho,\bl)}}
2^{m^\al_0-h(\rho)}
\sum_{\la\in\Par_0(\rho,d)}
2^{h(\rho)-m^\al_0+m_0^\theta-h(\la)}
\\
&\times\sum_{\bi,\bj\in J^d } 
\sum_{\substack{\Stab\in\Std_0^{\U,\bi}(\la)\\ \T\in\Std_0^{\V,\bj}(\la)}}
2^{\sum_{s=1}^d (i_s-\ell)}2^{\sum_{t=1}^d (j_t-\ell)}.
\end{align*}
Let $(\la^{(0)},\dots,\la^{(\ell)})=\quot(\la)$. Then, by Lemma~\ref{L110921}, we get that, for any fixed $\bk,\bl\in I^{\al}$, $\U\in\Std_0(\rho,\bk)$, $\V\in\Std_0(\rho,\bl)$,
\begin{align*}
\\
&\sum_{\bi,\bj\in J^d } 
\sum_{\substack{\Stab\in\Std_0^{\U,\bi}(\la)\\ \T\in\Std_0^{\V,\bj}(\la)}}
2^{\sum_{s=1}^d (i_s-\ell)}2^{\sum_{t=1}^d (j_t-\ell)}
\\
=\,&
\bigg(\sum_{\bi\in J^d } 
\sum_{\Stab\in\Std_0^{\U,\bi}(\la)}
2^{\sum_{s=1}^d (i_s-\ell)}
\bigg)
\bigg(
\sum_{\bj\in J^d } 
\sum_{\T\in\Std_0^{\V,\bj}(\la)}
2^{\sum_{t=1}^d (j_t-\ell)}\bigg)
\\
=\,&
2^{-2|\la^{(\ell)}|}\binom{d}{|\la^{(0)}|,\dots,|\la^{(\ell)}|}^2\,(K_{\la^{(0)}}')^2\,\prod_{j=1}^\ell K_{\la^{(j)}}^2,
\end{align*}
We conclude that $\dim \ttf R^{\Lambda_0}_{\theta}\ttf$ equals 
\begin{align*}
&
\sum_{\bk,\bl\in I^{\al}}
\sum_{\substack{\U\in\Std_0(\rho,\bk)\\ \V\in\Std_0(\rho,\bl)}}
2^{m^\al_0-h(\rho)}
\\
&\times \sum_{\la\in\Par_0(\rho,d)}
2^{h(\rho)-m^\al_0+m_0^\theta-h(\la)}2^{-2|\la^{(\ell)}|}\binom{d}{|\la^{(0)}|,\dots,|\la^{(\ell)}|}^2\,(K_{\la^{(0)}}')^2\,\prod_{j=1}^\ell K_{\la^{(j)}}^2
\\
=\,&
(\dim R^{\La_0}_\al)
\sum_{\la\in\Par_0(\rho,d)}
2^{2d-h(\la^{(0)})-2|\la^{(\ell)}|}\binom{d}{|\la^{(0)}|,\dots,|\la^{(\ell)}|}^2\,(K_{\la^{(0)}}')^2\,\prod_{j=1}^\ell K_{\la^{(j)}}^2,
\end{align*}
where we have used (\ref{E110921}) and the observation that $h(\la)-h(\rho)=h(\la^{(0)})$ for the last equality. An application of Lemma~\ref{L110921_2} completes the proof.
\end{proof}

\subsection{Basis of $B_d$}
For $1\leq t<d$, we consider the product of transpositions
\begin{equation}\label{EUT}
u_t:=\prod_{a=(t-1)p+1}^{tp}(a,a+p)\in\Si_{dp}.
\index{u@$u_t$}
\end{equation}
The elements $u_1,\dots,u_{d-1}$ generate a subgroup $G_d$ of $\Si_{dp}$ which permutes $d$ blocks of size $p$ and is isomorphic to $\Si_d$. More precisely, if $w=s_{t_1}\cdots s_{t_l}\in\Si_d$ we put $u_w:=u_{t_1}\cdots u_{t_l}\in G_d$,\index{u@$u_w$}
 and then 
$$
\Si_d\iso G_d,\ w\mapsto u_w
$$
is an isomorphism of groups. 

\begin{Lemma} \label{LBlockPermCusp} 
Let $d\in\Z_{\geq 1}$, $\bi^{(1)},\dots,\bi^{(d)},\bj^{(1)},\dots,\bj^{(d)}\in I^\de_\cus$, and  $x\in\D^{(p^d)}$ be such that 
$$x\cdot (\bi^{(1)}\cdots\bi^{(d)})=\bj^{(1)}\dots\bj^{(d)}.$$
Then $x=u_w$ for some $w\in \Si_d$ and $\bj^{(t)}=\bi^{(w^{-1}(t))}$ for all $t=1,\dots,d$.
\end{Lemma}
\begin{proof}
The letter $\ell$ appears in each word $\bi^1,\dots,\bi^d,\bj^{(1)},\dots,\bj^{(d)}$ exactly once and in the first position. In particular, $x^{-1}(1)=kp+1$, for some $0\leq k\leq d-1$ and $x(lp+1)>p$, for any $0\leq l\leq d-1$, with $l\neq k$. If $x^{-1}(m)=lp+n$, for some $2\leq m\leq p$, $l\neq k$ and $1\leq n\leq p$, then, since $x$ is a shuffle, $x(lp+1)\leq m$, a contradiction. Therefore, $x^{-1}(m)=kp+m$, for all $1\leq m\leq p$. Repeating this argument we get that there exists $w\in \Si_d$ such that $x^{-1}(lp+n)=w^{-1}(l)p+n$, for all $1\leq l\leq d-1$ and $1\leq n\leq p$. 
In other words, $x=u_w$ for some $w\in \Si_d$ and $u_w\cdot (\bi^{(1)}\cdots\bi^{(d)})=\bi^{(w^{-1}(1))}\cdots \bi^{(w^{-1}(d))}$.
\end{proof}

\begin{Corollary} \label{C130921_2} 
Let $w\in \Si_{2p}$ and $F$ be a polynomial in $y_1,\dots,y_{2p}$. If $w<u_1$ then in $\bar R_{2\de}$ we have 
$
e_{\de^2}F\psi_w e_{\de^2}\in \bar R_{\de^2}.
$
\end{Corollary}
\begin{proof}
Write $w=xv$ for $x\in\D^{(p^2)}$ and $v\in \Si_p\times \Si_p$ with $\ell(w)=\ell(x)+\ell(v)$. In view of Lemma~\ref{LBKW}(i), we may assume that $\psi_w=\psi_x\psi_v$. Now, 
$$
e_{\de^2}F\psi_we_{\de^2}=\sum_{\bi,\bj,\bk,\bl\in I^\de_\cus}
e(\bi\bj)F\psi_x\psi_ve(\bk\bl).
$$
We may restrict the summation to those $\bk,\bl\in I^\de_\cus$ for which $v\cdot(\bk\bl)\in I^{2\de}_\cus$---otherwise the corresponding summand is $0$. Since $v\in \Si_p\times\Si_p$, by Corollary~\ref{C130921}, for such $\bk,\bl$ we have that  $v\cdot(\bk\bl)$ is of the form $\bk'\bl'$ for $\bk',\bl'\in I^\de_\cus$. Now $x$ satisfies the assumptions of Lemma~\ref{LBlockPermCusp}. So either $x=u_1$ or $x=1$. The assumption that $w<u_1$ now implies that $x=1$. 
\end{proof}

\begin{Corollary} \label{CBruhat}
Let $w\in \Si_{2p}$ and $F$ be a polynomial in $y_1,\dots,y_{2p}$. If $w<u_1$ then in $\bar R_{2\de}$ we have 
$
\ga_{1^2}F\psi_w \ga_{1^2}\in B_{1^2}.
$
\end{Corollary}
\begin{proof}
Comes from Corollary~\ref{C130921_2} multiplying with $\ga_{1^2}$ on both sides. 
\end{proof}

Let $1\leq t<d$. Note that we have a reduced decomposition
$$
u_1=(s_ps_{p+1}\cdots s_{2p-1})(s_{p-1}s_p\cdots s_{2p-2})\cdots(s_{1}s_{2}\cdots s_p)
$$
We always choose this reduced decomposition to define $\psi_{u_1}$, i.e.
\begin{equation}\label{ERedU1}
\psi_{u_1}:=(\psi_p\psi_{p+1}\cdots \psi_{2p-1})(\psi_{p-1}\psi_p\cdots \psi_{2p-2})\cdots(\psi_{1}\psi_{2}\cdots \psi_p).
\end{equation}
In fact, the element $u_1$ is fully commutative, which means that one can go between any of its two reduced decompositions by a series of commutating Coxeter relations $s_rs_t=s_ts_r$ for $|r-t|>1$. So by the relation (\ref{R65}), a different choice of the reduced decomposition could only change $\psi_{u_1}$ by a sign.

Recall that we have identified $\bar R_{\de^2}$ with $\bar R_\de\otimes \bar R_\de$ (via Lemma~\ref{L030216}).  

\begin{Lemma} \label{L130921} 
In $\bar R_{2\de}$, for $\al\otimes \be\in \bar R_\de\otimes \bar R_\de$, we have 
$$
\psi_{u_1}(\al\otimes \be)\equiv(-1)^{|\al||\be|}(\be\otimes \al)\psi_{u_1}\pmod{\bar R_{\de^2}}. 
$$
\end{Lemma}
\begin{proof}
Note that $\psi_{u_1}(\al\otimes \be)\in e_{\de^2}\bar R_{2\de}e_{\de^2}$. Using the commuting relations in $R_{2\de}$ and Lemma~\ref{LBKW}, it is easy to see that 
$$\psi_{u_1}(\al\otimes \be)=(\be\otimes \al)\psi_{u_1}+(*),$$ 
where $(*)$ is a linear combination of the error terms of the form $$e_{\de^2}f_w(y_1,\dots,y_d)\psi_w e_{\de^2}=
f_w(y_1,\dots,y_d)e_{\de^2}\psi_w e_{\de^2},
$$
with $w<u_1$. An application of Corollary~\ref{C130921_2} completes the proof. 
\end{proof}

For $w\in \Si_d$ choose a reduced expression $w=s_{t_1}\dots s_{t_l}$ and define 
\begin{equation*}\label{}
\psi_{u_w}:=\psi_{u_{t_1}}\dots\psi_{u_{t_l}}\in \bar R_{d\de}.
\end{equation*}
As usual, in general, $\psi_w$ depends on the choice of a reduced expression of $w$. 

\begin{Lemma}\label{LBasisFirst}
We have that 
$e_{\de^d}\bar R_{d\de}e_{\de^d}$ is a free left $\bar R_{\de^d}$-module with basis $\{e_{\de^d}\psi_{u_w}\mid w\in\Si_d\}$. 
\end{Lemma}
\begin{proof}
By Lemma~\ref{L030216}(iii), we have that $e_{\de^d}\bar R_{d\de}$ is a free left $\bar R_{\de^d}$-module with basis $\{e_{\de^d}\psi_x\mid x\in {}^{(p^d)}\hspace{-1mm}\D\}$. 
The result now follows from Lemma~\ref{LBlockPermCusp}. 
\end{proof}

Recall the notation (\ref{EWSignAction}). 

\begin{Lemma} \label{L130921_3}
Let $w\in\Si_d$. 
In $e_{\de^d}\bar R_{d\de}e_{\de^d}$, for $\be_1\otimes \dots\otimes \be_d\in \bar R_{\de^d}$, we have 
$$
\psi_{u_w}(\be_1\otimes \dots\otimes\be_d)\equiv{}^w(\be_1\otimes \dots\otimes\be_d)\psi_{u_w}\pmod{\sum_{\substack {v\in\Si_d\\ v< w}} \bar R_{\de^d}\psi_{u_v}}.
$$
\end{Lemma}
\begin{proof}
This follows by repeated application of Lemma~\ref{L130921}.  
\end{proof}

Define 
\begin{equation}\label{ESigma1}
\upsigma:=\ga_{1^2}\psi_{u_1}\ga_{1^2}\in B_{2}.
\index{s@$\upsigma$}
\end{equation}
In terms of Khovanov-Lauda diagrams, we have $\upsigma\ga^{ij}$ is represented by 
$$
\resizebox{126mm}{25mm}{
\begin{braid}\tikzset{baseline=-.3em}
	\draw (0,6) node{\color{blue}$\ell$\color{black}};
	\braidbox{0.6}{2.9}{5.3}{6.6}{};
	\draw (1.8,6) node{$(\ell-1)^2$};
	\draw[dots] (3.5,6)--(5,6);
	\braidbox{5.3}{7.6}{5.3}{6.6}{};
	\draw (6.5,6) node{$(j+1)^2$};
	\draw (8.3,6) node{$j$};
	\draw (9.5,6) node{$j-1$};
	\draw[dots] (10.8,6)--(12.5,6);
	\draw (12.7,6) node{$1$};
	\redbraidbox{13.4}{14.3}{5.3}{6.6}{};
	\draw (13.8,6) node{$\color{red}0\,\,0\color{black}$};
	\draw (15,6) node{$1$};
	\draw[dots] (15.6,6)--(17,6);
	\draw (18.2,6) node{$j-1$};
	\draw (19.5,6) node{$j$};
	\draw (0,-6) node{\color{blue}$\ell$\color{black}};
	\braidbox{0.6}{2.9}{-6.7}{-5.4}{};
	\draw (1.8,-6) node{$(\ell-1)^2$};
	\draw[dots] (3.5,-6)--(5,-6);
	\braidbox{5.3}{7.6}{-6.7}{-5.4}{};
	\draw (6.5,-6) node{$(i+1)^2$};
	\draw (8.3,-6) node{$i$};
	\draw (9.5,-6) node{$i-1$};
	\draw[dots] (10.8,-6)--(12.5,-6);
	\draw (12.7,-6) node{$1$};
	\redbraidbox{13.4}{14.3}{-6.7}{-5.4}{};
	\draw (13.8,-6) node{$\color{red}0\,\,0\color{black}$};
	\draw (15,-6) node{$1$};
	\draw[dots] (15.6,-6)--(17,-6);
	\draw (18.2,-6) node{$i-1$};
	\draw (19.5,-6) node{$i$};

	\draw (20.3,6) node{\color{blue}$\ell$\color{black}};
	\braidbox{20.9}{23.2}{5.3}{6.6}{};
	\draw (22.1,6) node{$(\ell-1)^2$};
	\draw[dots] (23.8,6)--(25.3,6);
	\braidbox{25.6}{27.9}{5.3}{6.6}{};
	\draw (26.8,6) node{$(i+1)^2$};
	\draw (28.6,6) node{$i$};
	\draw (29.8,6) node{$i-1$};
	\draw[dots] (31.1,6)--(32.8,6);
	\draw (33,6) node{$1$};
	\redbraidbox{33.7}{34.6}{5.3}{6.6}{};
	\draw (34.1,6) node{$\color{red}0\,\,0\color{black}$};
	\draw (35.3,6) node{$1$};
	\draw[dots] (35.9,6)--(37.3,6);
	\draw (38.5,6) node{$i-1$};
	\draw (39.8,6) node{$i$};
	\draw (20.3,-6) node{\color{blue}$\ell$\color{black}};
	\braidbox{20.9}{23.2}{-6.7}{-5.4}{};
	\draw (22.1,-6) node{$(\ell-1)^2$};
	\draw[dots] (23.8,-6)--(25.1,-6);
	\braidbox{25.6}{27.9}{-6.7}{-5.4}{};
	\draw (26.8,-6) node{$(j+1)^2$};
	\draw (28.6,-6) node{$j$};
	\draw (29.8,-6) node{$j-1$};
	\draw[dots] (31.1,-6)--(32.8,-6);
	\draw (33,-6) node{$1$};
	\redbraidbox{33.7}{34.6}{-6.7}{-5.4}{};
	\draw (34.1,-6) node{$\color{red}0\,\,0\color{black}$};
	\draw (35.3,-6) node{$1$};
	\draw[dots] (35.9,-6)--(37.3,-6);
	\draw (38.5,-6) node{$j-1$};
	\draw (39.8,-6) node{$j$};
	\draw[blue](0,-5.2)--(20.3,5.2);
	\draw(1.2,-5.2)--(21.3,5.2);
	\draw(2.2,-5.2)--(22.3,5.2);
	\draw(6,-5.2)--(26,5.2);
	\draw(7,-5.2)--(27,5.2);
	\draw(8.3,-5.2)--(28.6,5.2);
	\draw(9.5,-5.2)--(29.8,5.2);
	\draw(12.7,-5.2)--(33.1,5.2);
	\draw[red](13.5,-5.2)--(33.8,5.2);
	\draw[red](14.2,-5.2)--(34.5,5.2);
	\draw(15,-5.2)--(35.3,5.2);
	\draw(18.2,-5.2)--(38.5,5.2);
	\draw(39.8,-5.2)--(19.5,5.2);
	\draw(39.8,5.2)--(19.5,-5.2);
	\draw[blue] (20.3,-5.2)--(0,5.2);
	\draw(21.3,-5.2)--(1.2,5.2);
	\draw(22.3,-5.2)--(2.2,5.2);
	\draw(26,-5.2)--(6,5.2);
	\draw(27,-5.2)--(7,5.2);
	\draw(28.6,-5.2)--(8.3,5.2);
	\draw(29.8,-5.2)--(9.5,5.2);
	\draw(38.5,-5.2)--(18.2,5.2);
	\draw(35.3,-5.2)--(15,5.2);
	\draw[red](34.5,-5.2)--(14.2,5.2);
	\draw[red](33.8,-5.2)--(13.5,5.2);
	\draw(33.1,-5.2)--(12.7,5.2);
	\end{braid},
	}
$$
where the choice of the reduced decomposition that we have made means that we must interpret the ``double odd crossing" as follows:
\begin{equation}\label{EDoubleRed}
\begin{braid}\tikzset{baseline=1.3em}                 
	\draw[red](2,1)node[below]{$0$}--(0,3);
	 \draw[red] (3,1)node[below]{$0$}--(1,3);
	 \draw[red](0,1)node[below]{$0$}--(2,3);
	 \draw[red](1,1)node[below]{$0$}--(3,3);
	  \end{braid}
	  :=
\begin{braid}\tikzset{baseline=1.3em}
                \draw[red](2,1)node[below]{$0$}--(0.7,1.9)--(0,3);
	 \draw[red](3,1)node[below]{$0$}--(2.3,2.2)--(1,3);
	 \draw[red](0,1)node[below]{$0$}--(2,3);
	 \draw[red](1,1)node[below]{$0$}--(3,3);
	  \end{braid},
\end{equation}
while the order of other crossings is unimportant due to the relation (\ref{R65}) and the full commutativity of $u_1$. 
By (\ref{EGrading3}) and the fact that $(\de|\de)=0$, the element $\upsigma\in B_d$ has bidegree $(0,\0)$.

Recalling the embeddings (\ref{EIotaRR+1}), we now define bidegree $(0,\0)$ elements 
\begin{equation}\label{ESigmaR}
\upsigma_r:=\iota_{r,r+1}(\upsigma)\in B_{d}\qquad(1\leq r<d).
\index{s@$\upsigma_r$}
\end{equation}
Note that
$
\upsigma_r=\ga_{1^d}\psi_{u_r}\ga_{1^d}
$
where $\psi_{u_r}$ is defined using the choice of reduced decomposition for $u_r$ similar to (\ref{ERedU1}).

For $w\in \Si_d$ choose a reduced expression $w=s_{t_1}\dots s_{t_l}$ and define the bidegree $(0,\0)$ 
element 
\begin{equation*}\label{}
\upsigma_w:=\upsigma_{t_1}\dots\upsigma_{t_l}\in B_{d}.
\index{s@$\upsigma_w$}
\end{equation*}
In general, $\upsigma_w$ depends on the choice of a reduced expression of $w$.

\begin{Lemma} \label{LIdOneSideSiW} 
Let $w\in\Si_d$. Then $\upsigma_w=\psi_{u_w}\ga_{1^d}$. 
\end{Lemma}
\begin{proof}
This follows by a repeated application of Lemma~\ref{L150921_2}.
\end{proof}

\begin{Lemma} \label{L160921} 
Let $C$ be an algebra and $D$ be a (unital) subalgebra such that $C$ is a free left $D$-module with basis $\{x_1,\dots,x_n\}$. Let $e\in D$ be an idempotent such that for all $i=1,\dots,n$ the following two conditions hold:
\begin{enumerate}
\item[{\rm (1)}] $ex_ie=x_ie$;
\item[{\rm (2)}] $x_ie=ex_i+\sum_{j<i}d_{i,j}x_j$ for some $d_{i,j}\in D$.
\end{enumerate}
Then $eCe$ is a free $eDe$-module with basis $\{ex_1e,\dots,ex_ne\}$. 
\end{Lemma}
\begin{proof}
By assumption, every $c\in C$ can be written uniquely in the form $c=\sum_{i=1}^nd_ix_i$ with all $d_i\in D$. If $c\in eCe$, this implies using (1):
$$
c=\sum_{i=1}^ned_ix_ie=\sum_{i=1}^n(ed_ie)(ex_ie).
$$
Thus every $c\in eCe$ is an $eDe$-linear combination of the elements $ex_ie$. To prove uniqueness, suppose $\sum_{i=1}^nd_iex_ie=0$ with all $d_i\in eDe$ and not all $d_i$ being $0$. Let $j$ be maximal with $d_j\neq 0$.  
Then using (2), we get 
$$
0=\sum_{i=1}^jd_iex_ie=
\sum_{i=1}^jd_ie(x_ie)=\sum_{i=1}^jd_ie(ex_i+\sum_{k<i}d_{i,k}x_k)=d_jx_j+\sum_{i<j}d_i'x_i
$$
for some $d_i'\in D$. Since $\{x_1,\dots,x_n\}$ is a basis of the $D$-module $C$, it follows that $d_j=0$, giving a contradiction. 
\end{proof}

\begin{Proposition} \label{PBdBasis} 
We have that  $B_d$ is a free left $B_{1^d}$-module with basis $\{\upsigma_w\mid w\in \Si_d\}$. In particular, $B_d=\bigoplus_{w\in \Si_d}B_{1^d}\upsigma_w$. 
\end{Proposition}
\begin{proof}
We apply Lemma~\ref{L160921} with $C=e_{\de^d}\bar R_{d\de}e_{\de^d}$, $D=\bar R_{\de^d}$, and $e=\ga_{1^d}$. By Lemma~\ref{LBasisFirst}, $C$ is then a free left $D$-module with basis $\{e_{\de^d}\psi_{u_w}\mid w\in\Si_d\}$. The assumption  (1) of Lemma~\ref{L160921} comes from Lemma~\ref{LIdOneSideSiW}, while the assumption (2) comes from Lemma~\ref{L130921_3}. 
\end{proof} 

Recalling (\ref{EParDecB}), we get:

\begin{Corollary} \label{CBBasis} 
We have as vector spaces  
$$
B_d=\bigoplus_{\substack{c_1,\dots,c_d\in\Z_{\geq 0},\\ w\in {\Si}_d}}z_1^{c_1}\cdots z_d^{c_d} \,\Zig_\ell^{\otimes d}\,\upsigma_w.
$$
Moreover, for $c_1,\dots,c_d\in\Z_{\geq 0}$, $w\in {\Si}_d$ and $\ba\in \Zig_\ell^{\otimes d}$ we have $z_1^{c_1}\cdots z_d^{c_d}\ba\upsigma_w=0$ only if $\ba=0$. 
\end{Corollary}

Since the algebra $\Zig_\ell$ is non-negatively graded, $z_i$'s have degree $4$, and $\upsigma_w$'s have degree $0$, we deduce:

\begin{Corollary} \label{CBDNonNeg} 
The algebra $B_d$ is non-negatively graded. 
\end{Corollary}

For $n\in\Z_{\geq 0}$ we set
\begin{equation}
\label{EBN}
B_d^{(n)}:=\sum_{\substack{c_1+\dots+c_d\leq n,\\ w\in {\Si}_d}}z_1^{c_1}\cdots z_d^{c_d} \,\Zig_\ell^{\otimes d}\,\upsigma_w\subseteq B_d.
\index{b@$B_d^{(n)}$}
\end{equation}

\begin{Lemma} \label{L160921_4}
We have  $\Zig_\ell^{\otimes d}\upsigma_w z_1^{c_1}\cdots z_d^{c_d}\subseteq B_d^{(c_1+\dots+c_d)}$ for any $w\in\Si_d$ and $c_1\dots,c_d\in\Z_{\geq 0}$.
\end{Lemma}
\begin{proof}
Note that $\deg(\upsigma_w z_1^{c_1}\cdots z_d^{c_d})=4c_1+\dots+4c_d$, so by degrees using Proposition~\ref{PBdBasis} and Corollary~\ref{CBDNonNeg} we deduce that $\upsigma_w z_1^{c_1}\cdots z_d^{c_d}\in B_d^{(c_1+\dots+c_d)}$.
Now the lemma follows since by (\ref{ETwistedRelations}), we have $\Zig_\ell^{\otimes d} z_1^{c_1}\cdots z_d^{c_d}=z_1^{c_1}\cdots z_d^{c_d}\Zig_\ell^{\otimes d}$.
\end{proof}

\begin{Lemma} \label{PzFiltSubalg} 
We have that $B_d^{(0)}$ is a subalgebra of $B_d$.
\end{Lemma}
\begin{proof}
It suffices to prove that $\upsigma_w\upsigma_u\in B_d^{(0)}$ and $\upsigma_w\Zig_\ell^{\otimes d}\subseteq B_d^{(0)}$ for all $w,u\in\Si_d$. The first containment is clear using Corollary~\ref{CBDNonNeg} since since $\deg(\upsigma_w)=\deg(\upsigma_u)=0$ while $\deg(z_r)=4$ for all $r$. 

For the second containment, since $\upsigma_w$ is a product of some $\upsigma_r$'s, it suffices to to check that $\upsigma_r\Zig_\ell^{\otimes d}\subseteq  B_d^{(0)}$ for any $r$. Using the embedding $\iota_{r,r+1}$ this reduces to the case $r=1$ and $d=2$, in which case we need to prove  $\upsigma\Zig_\ell^{\otimes 2}\in B_2^{(0)}$. For this it is sufficient to prove that 
$
\upsigma(\zb\otimes 1)\in B_2^{(0)}
$
and
$
\upsigma(1\otimes \zb)\in B_2^{(0)}
$
for any standard generating element $\zb$ of $\Zig_\ell$, i.e. for $\zb$ of the form $\ze^i,\zu,\za^{i,i+1}$ or $\za^{i+1,i}$. Now, by degrees, we only need to consider the generators of degree $4$, i.e. we may assume that $\zb=\za^{i+1,i}$ for some $i$. So it suffices to demonstrate that $\upsigma(\za^{i+1,i}\otimes\ze^j)\in B_2^{(0)}$ and $\upsigma(\ze^j\otimes \za^{i+1,i})\in B_2^{(0)}$ for all admissible $i,j$. The two cases are entirely similar, so we give details for the first one. 

By Corollary~\ref{CBBasis}, $\upsigma(\za^{i+1,i}\otimes\ze^j)$ is a linear combination of terms of the form
$
z_1^{c_1}z_2^{c_2}(\za_1\otimes\za_2)\upsigma_w,
$
and we need to show that the terms corresponding to $c_1+c_2>0$ are zero. If not, taking into account that $\deg(\za^{i+1,i}\otimes\ze^j)=\deg(z_1)=\deg(z_2)=4$, the bad terms can only be of the form $
z_r(\za_1\otimes\za_2)\upsigma_w,
$
for $r\in\{1,2\}$ and $\deg(\za_1)=\deg(\za_2)=0$. Note that
$$
\upsigma(\za^{i+1,i}\otimes\ze^j)=\ga^{j,i+1}\upsigma(\za^{i+1,i}\otimes\ze^j)\ga^{i,j},
$$
so 
$$
z_r(\za_1\otimes\za_2)\upsigma_w=\ga^{j,i+1}z_r(\za_1\otimes\za_2)\upsigma_w\ga^{i,j}.
$$
It follows that $\za_1=\ze^j\za_1$ and $\za_2=\ze^{i+1}\za_2$. 

Suppose $w=1$. Then we also have $\za_1=\za_1\ze^i$ and $\za_2=\za_2\ze^{j}$, i.e. $\za_1=\ze^j\za_1\ze^i$ and $\za_2=\ze^{i+1}\za_2\ze^{j}$, which implies that $\za_1=0$ or $\za_2=0$ since these elements have degree $0$. 

In the case $w=s_1$, we deduce similarly that $\za_1=\ze^j\za_1\ze^j$ and $\za_2=\ze^{i+1}\za_2\ze^i$, whence $\za_2=0$ by degrees.
\end{proof}

\begin{Proposition} \label{PzFilt} 
For all $n,m\in\Z_{\geq 0}$ we have $B_d^{(n)}B_d^{(m)}\subseteq B_d^{(n+m)}$. 
\end{Proposition}
\begin{proof}
This follows from (\ref{EFiltrPar}) and Lemmas~\ref{L160921_4}, \ref{PzFiltSubalg}.  
\end{proof}

\section{Intertwiners}
\label{SIntertwiners}
The elements $\upsigma_w\in B_d$ constricted in the previous section do not satisfy the same relations as the elements $w$ in the symmetric group $\Si_d$, so they do not give us a copy of $\Si_d$ inside $B_d$. 
In this section we construct certain explicit {\em intertwining elements} of the form 
$$\tau_w=\upsigma_w+\sum_{u<w}c_u\upsigma_u\qquad(c_u\in\F)$$ 
which do give us a copy of the symmetric group $\Si_d$ inside  $B_d$. Moreover, we show that this group acts on $\Zig_\ell^{\otimes d}$ by (super)place permutations on tensors. 

The classical analogues of the intertwining elements $\tau _w$ have first appeared in \cite[\S4]{KMR}. But the most effective approach to establishing their key properties is based on the deformation approach going back to \cite{KKK}. The corresponding deformation technique for quiver Hecke superalgebras has been developed by Kang, Kashiwara and Oh in \cite[\S8]{KKO}. We refer the reader to \S\ref{SSKKO} for a review of the key properties of the KKO intertwiners $\phi_u$ which will be used to construct the intertwining elements $\tau_w$.

\subsection{The element {\sansmath$\Phi$}}

Recalling (\ref{EUT}) and (\ref{EPhiW}), we define
\begin{equation}\label{EPhiBig}
\text{\sansmath$\Phi$}:=\ga_{1^2}\phi_{u_1}\ga_{1^2}\in B_2, \index{f@\text{\sansmath$\Phi$}}
\end{equation}
where the reduced decomposition for $u_1$ is chosen as in (\ref{ERedU1}). Thus:
\begin{equation}\label{ERedPhi}
\text{\sansmath$\Phi$}=\ga_{1^2}(\phi_p\phi_{p+1}\cdots \phi_{2p-1})(\phi_{p-1}\phi_p\cdots \phi_{2p-2})\cdots(\phi_{1}\phi_{2}\cdots \phi_p)\ga_{1^2}.
\end{equation}

\begin{Example} \label{ExamplePhiPsi0} 
{\rm 
Using the notation (\ref{EPhiDiagram}), we have for the case $\ell=3$:
$$
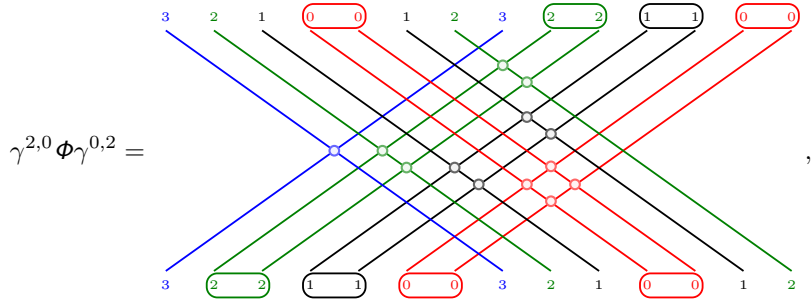
\ga^{2,0}\text{\sansmath$\Phi$}\ga^{0,2}=
\begin{braid}\tikzset{baseline=-.3em}
		\darkgreenbraidbox{1.9}{4.1}{-6.1}{-5.1}{};
		\braidbox{5.9}{8.1}{-6.1}{-5.1}{};
		\redbraidbox{9.9}{12.1}{-6.1}{-5.1}{};
		\redbraidbox{19.9}{22.1}{-6.1}{-5.1}{};
		\redbraidbox{23.9}{26.1}{5.1}{6.1}{};
		\darkgreenbraidbox{15.9}{18.1}{5.1}{6.1}{};
		\braidbox{19.9}{22.1}{5.1}{6.1}{};
\redbraidbox{5.9}{8.1}{5.1}{6.1}{};
	\draw[blue](0,-5)node[below]{$3$}--(14,5)node[above]{$3$};
	\draw[darkgreen](2,-5)node[below]{$2$}--(16,5)node[above]{$2$};
	\draw[darkgreen](4,-5)node[below]{$2$}--(18,5)node[above]{$2$};
	\draw(6,-5)node[below]{$1$}--(20,5)node[above]{$1$};
	\draw(8,-5)node[below]{$1$}--(22,5)node[above]{$1$};
	\draw[red](10,-5)node[below]{$0$}--(24,5)node[above]{$0$};
	\draw[red](12,-5)node[below]{$0$}--(26,5)node[above]{$0$};
	\draw[blue](14,-5)node[below]{$3$}--(0,5)node[above]{$3$};
	\draw[darkgreen](16,-5)node[below]{$2$}--(2,5)node[above]{$2$};
	\draw(18,-5)node[below]{$1$}--(4,5)node[above]{$1$};
\draw[red](20,-5)node[below]{$0$}--(6,5)node[above]{$0$};
\draw[red](22,-5)node[below]{$0$}--(8,5)node[above]{$0$};
\draw(24,-5)node[below]{$1$}--(10,5)node[above]{$1$};
\draw[darkgreen](26,-5)node[below]{$2$}--(12,5)node[above]{$2$};
\filldraw[color=black!60, fill=black!5, thick] (12,-0.7) circle (0.2);
\filldraw[color=black!60, fill=black!5, thick] (13,-1.4) circle (0.2);
\filldraw[color=black!60, fill=black!5, thick] (15,1.4) circle (0.2);
\filldraw[color=black!60, fill=black!5, thick] (16,0.7) circle (0.2);
 \filldraw[color=blue!60, fill=blue!5, thick] (7,0) circle (0.2);
\filldraw[color=darkgreen!60, fill=darkgreen!5, thick] (9,0) circle (0.2);
\filldraw[color=darkgreen!60, fill=darkgreen!5, thick] (10,-0.7) circle (0.2);
\filldraw[color=darkgreen!60, fill=darkgreen!5, thick] (14,3.55) circle (0.2);
\filldraw[color=darkgreen!60, fill=darkgreen!5, thick] (15,2.85) circle (0.2);
\filldraw[color=red!60, fill=red!5, thick] (16,-0.65) circle (0.2);
\filldraw[color=red!60, fill=red!5, thick] (17,-1.4) circle (0.2);
\filldraw[color=red!60, fill=red!5, thick] (15,-1.4) circle (0.2);
\filldraw[color=red!60, fill=red!5, thick] (16,-2.1) circle (0.2);
\end{braid},
$$
where we have used (\ref{EPhiDifCol}) for the different color crossings. 
}
\end{Example}

Here and below, as usual, we interpret the ``double odd crossing" as follows:
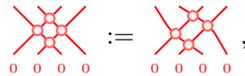
\begin{equation*}\label{EDoubleRedPhi}
\begin{braid}\tikzset{baseline=1.3em}                 
	\draw[red](2,1)node[below]{$0$}--(0,3);
	 \draw[red] (3,1)node[below]{$0$}--(1,3);
	 \draw[red](0,1)node[below]{$0$}--(2,3);
	 \draw[red](1,1)node[below]{$0$}--(3,3);
	 \filldraw[color=red!60, fill=red!5, thick] (1.5,1.5) circle (0.2);
	 \filldraw[color=red!60, fill=red!5, thick] (1.5,2.5) circle (0.2);
	 \filldraw[color=red!60, fill=red!5, thick] (1,2) circle (0.2);
	 \filldraw[color=red!60, fill=red!5, thick] (2,2) circle (0.2);
	  \end{braid}
	  :=
\begin{braid}\tikzset{baseline=1.3em}
                \draw[red](2,1)node[below]{$0$}--(0.7,1.9)--(0,3);
	 \draw[red](3,1)node[below]{$0$}--(2.3,2.2)--(1,3);
	 \draw[red](0,1)node[below]{$0$}--(2,3);
	 \draw[red](1,1)node[below]{$0$}--(3,3);
	 \filldraw[color=red!60, fill=red!5, thick] (1.4,1.4) circle (0.2);
	 \filldraw[color=red!60, fill=red!5, thick] (1.58,2.6) circle (0.2);
	 \filldraw[color=red!60, fill=red!5, thick] (0.86,1.86) circle (0.2);
	 \filldraw[color=red!60, fill=red!5, thick] (2.19,2.19) circle (0.2);
	  \end{braid},
\end{equation*}
while the order of other non-neighboring crossings is unimportant due to Corollary~\ref{C180621}(i). The choice we have made for double odd crossings implies

\begin{Lemma} \label{LQThickRed} 
In $R_{4\al_0}$, we have
$$
\begin{braid}\tikzset{baseline=2.2em}                 
	\draw[red](2,1)node[below]{$0$}--(0,3)--(2,5);
	 \draw[red] (3,1)node[below]{$0$}--(1,3)--(3,5);
	 \draw[red](0,1)node[below]{$0$}--(2,3)--(0,5);
	 \draw[red](1,1)node[below]{$0$}--(3,3)--(1,5);
	 \filldraw[color=red!60, fill=red!5, thick] (1.5,1.5) circle (0.2);
	 \filldraw[color=red!60, fill=red!5, thick] (1.5,2.5) circle (0.2);
	 \filldraw[color=red!60, fill=red!5, thick] (1,2) circle (0.2);
	 \filldraw[color=red!60, fill=red!5, thick] (2,2) circle (0.2);
	 \filldraw[color=red!60, fill=red!5, thick] (1.5,3.5) circle (0.2);
	 \filldraw[color=red!60, fill=red!5, thick] (1.5,4.5) circle (0.2);
	 \filldraw[color=red!60, fill=red!5, thick] (1,4) circle (0.2);
	 \filldraw[color=red!60, fill=red!5, thick] (2,4) circle (0.2);
	  \end{braid}
	  =-(y_1^2+y_3^2)(y_1^2+y_4^2)(y_2^2+y_3^2)(y_2^2+y_4^2).
$$
\end{Lemma}
\begin{proof}
In view of the conventions made, the left hand equals 
\begin{align*}
\begin{braid}\tikzset{baseline=2.2em}
                \draw[red](2,1)node[below]{$0$}--(0.7,1.9)--(0,3)--(2,5);
	 \draw[red](3,1)node[below]{$0$}--(2.3,2.2)--(1,3)--(3,5);
	 \draw[red](0,1)node[below]{$0$}--(2,3)--(0.7,3.9)--(0,5);
	 \draw[red](1,1)node[below]{$0$}--(3,3)--(2.3,4.2)--(1,5);
	 \filldraw[color=red!60, fill=red!5, thick] (1.4,1.4) circle (0.2);
	 \filldraw[color=red!60, fill=red!5, thick] (1.58,2.6) circle (0.2);
	 \filldraw[color=red!60, fill=red!5, thick] (0.86,1.86) circle (0.2);
	 \filldraw[color=red!60, fill=red!5, thick] (2.19,2.19) circle (0.2);
	 \filldraw[color=red!60, fill=red!5, thick] (1.4,3.4) circle (0.2);
	 \filldraw[color=red!60, fill=red!5, thick] (1.58,4.6) circle (0.2);
	 \filldraw[color=red!60, fill=red!5, thick] (0.86,3.86) circle (0.2);
	 \filldraw[color=red!60, fill=red!5, thick] (2.19,4.19) circle (0.2);
	  \end{braid}
&=
\begin{braid}\tikzset{baseline=2.2em}
                \draw[red](2,1)node[below]{$0$}--(0.7,1.9)--(0,3)--(2,5);
	 \draw[red](3,1)node[below]{$0$}--(2.3,2.2)--(1.7,2.7)--(1.7,3.6)--(3,5);
	 \draw[red](0,1)node[below]{$0$}--(1.4,2.5)--(1.4,3.6)--(0.7,3.9)--(0,5);
	 \draw[red](1,1)node[below]{$0$}--(3,3)--(2.3,4.2)--(1,5);
	 \filldraw[color=red!60, fill=red!5, thick] (1.4,1.4) circle (0.2);
	 \filldraw[color=red!60, fill=red!5, thick] (0.86,1.86) circle (0.2);
	 \filldraw[color=red!60, fill=red!5, thick] (2.19,2.19) circle (0.2);
	 \filldraw[color=red!60, fill=red!5, thick] (1.58,4.6) circle (0.2);
	 \filldraw[color=red!60, fill=red!5, thick] (0.86,3.86) circle (0.2);
	 \filldraw[color=red!60, fill=red!5, thick] (2.19,4.19) circle (0.2);
	 \reddot (1.37,3.2);
	  \reddot (1.37,2.8);
	  \end{braid}
+
\begin{braid}\tikzset{baseline=2.2em}
                \draw[red](2,1)node[below]{$0$}--(0.7,1.9)--(0,3)--(2,5);
	 \draw[red](3,1)node[below]{$0$}--(2.3,2.2)--(1.7,2.7)--(1.7,3.6)--(3,5);
	 \draw[red](0,1)node[below]{$0$}--(1.4,2.5)--(1.4,3.6)--(0.7,3.9)--(0,5);
	 \draw[red](1,1)node[below]{$0$}--(3,3)--(2.3,4.2)--(1,5);
	 \filldraw[color=red!60, fill=red!5, thick] (1.4,1.4) circle (0.2);
	 \filldraw[color=red!60, fill=red!5, thick] (0.86,1.86) circle (0.2);
	 \filldraw[color=red!60, fill=red!5, thick] (2.19,2.19) circle (0.2);
	 \filldraw[color=red!60, fill=red!5, thick] (1.58,4.6) circle (0.2);
	 \filldraw[color=red!60, fill=red!5, thick] (0.86,3.86) circle (0.2);
	 \filldraw[color=red!60, fill=red!5, thick] (2.19,4.19) circle (0.2);
	 \reddot (1.7,3.2);
	  \reddot (1.72,2.8);
	  \end{braid}
=
(y_1^2+y_4^2)\begin{braid}\tikzset{baseline=2.2em}
                \draw[red](2,1)node[below]{$0$}--(0.7,1.9)--(0,3)--(2,5);
	 \draw[red](3,1)node[below]{$0$}--(2.3,2.2)--(1.7,2.7)--(1.7,3.6)--(3,5);
	 \draw[red](0,1)node[below]{$0$}--(1.4,2.5)--(1.4,3.6)--(0.7,3.9)--(0,5);
	 \draw[red](1,1)node[below]{$0$}--(3,3)--(2.3,4.2)--(1,5);
	 \filldraw[color=red!60, fill=red!5, thick] (1.4,1.4) circle (0.2);
	 \filldraw[color=red!60, fill=red!5, thick] (0.86,1.86) circle (0.2);
	 \filldraw[color=red!60, fill=red!5, thick] (2.19,2.19) circle (0.2);
	 \filldraw[color=red!60, fill=red!5, thick] (1.58,4.6) circle (0.2);
	 \filldraw[color=red!60, fill=red!5, thick] (0.86,3.86) circle (0.2);
	 \filldraw[color=red!60, fill=red!5, thick] (2.19,4.19) circle (0.2);
	  \end{braid}
	  =
-(y_1^2+y_4^2)\begin{braid}\tikzset{baseline=2.2em}
                \draw[red](2,1)node[below]{$0$}--(0.7,1.9)--(0,3)--(2,5);
	 \draw[red](3,1)node[below]{$0$}--(2.2,2.1)--(1.8,3)--(2.5,3.6)--(3,5);
	 \draw[red](0,1)node[below]{$0$}--(1.4,2.5)--(1.4,3.6)--(0.7,3.9)--(0,5);
	 \draw[red](1,1)node[below]{$0$}--(2.6,2.6)--(1,5);
	 \filldraw[color=red!60, fill=red!5, thick] (1.4,1.4) circle (0.2);
	 \filldraw[color=red!60, fill=red!5, thick] (0.86,1.86) circle (0.2);
	 \filldraw[color=red!60, fill=red!5, thick] (2.19,2.19) circle (0.2);
	 \filldraw[color=red!60, fill=red!5, thick] (1.41,4.41) circle (0.2);
	 \filldraw[color=red!60, fill=red!5, thick] (0.86,3.86) circle (0.2);
	 \filldraw[color=red!60, fill=red!5, thick] (2.15,3.31) circle (0.2);
	  \end{braid},
\end{align*}
where we have used Lemma~\ref{LphiQuadr} for the first equality, (\ref{EDotPastCircle}) for the third equality, and Corollary~\ref{C180621}(i) for the last equality. Using Lemma~\ref{LphiQuadr} and (\ref{EDotPastCircle}) three more times yields the required result.
\end{proof}

Recall the parabolic subalgebra $B\otimes B=B_{1^2}\subseteq B_2$ from (\ref{EBPar}) and the notation (\ref{EAR}). 
 
\begin{Lemma} \label{LPhiGa} 
Let $b,b'\in B$. Then in $B_2$, we have 
$$\text{\sansmath$\Phi$} (b\otimes b')=(-1)^{|b||b'|}(b'\otimes b)\text{\sansmath$\Phi$}.$$ 
In particular, $\text{\sansmath$\Phi$} z_1=z_2\text{\sansmath$\Phi$}$, $\text{\sansmath$\Phi$} z_2=z_1\text{\sansmath$\Phi$}$, and 
\begin{align*}
\text{\sansmath$\Phi$}\ga^{i,j}=\ga^{j,i}\text{\sansmath$\Phi$}=\ga^{j,i}\text{\sansmath$\Phi$}\ga^{i,j}
=
\ga_{1^2}\phi_{u_1}\ga^{i,j}\qquad(i,j\in J).
\end{align*}
\end{Lemma}
\begin{proof}
In view of Theorem~\ref{TCuspIsoRank1}, we may assume that $b$ and $b'$ are of the form $e^i,u,a^{i,j}$ or $z^i$, as defined explicitly in \S\ref{SSSettingd=1}. 
For these elements the first required commutation relation follows from Lemma~\ref{LKKO}. The equalities $\text{\sansmath$\Phi$}\ga^{i,j}=\ga^{j,i}\text{\sansmath$\Phi$}=\ga^{j,i}\text{\sansmath$\Phi$}\ga^{i,j}$ now come as a special case using the fact that $\ga_{1^2}=\sum_{i,j\in J}\ga^{i,j}$. For the last equality, we have
$$
\ga^{j,i}\text{\sansmath$\Phi$}\ga^{i,j}
=\ga^{j,i}\phi_{u_1}\ga^{i,j}
=(\ze^j\otimes \ze^i)\phi_{u_1}\ga^{i,j}
=\ga_{1^2}\phi_{u_1}(\ze^i\otimes \ze^j)\ga^{i,j}
=\ga_{1^2}\phi_{u_1}\ga^{i,j},
$$
where we have used Lemma~\ref{LKKO} for the third equality.
\end{proof}

\subsection{Quadratic relation for $\text{\sansmath$\Phi$}$}\label{sec:quad_rel}

Recall the notation (\ref{ERS}). For $i,k\in J$, denote
\begin{align*}
&R^i_k:=\index{r@$R^i_k$}
\left\{
\begin{array}{ll}
y_{r^i_k} &\hbox{if $k\neq 0$,}\\
y_{r^i_0}^2 &\hbox{if $k=0$,}
\end{array}
\right.
&\tilde R^i_k:=\index{r@$\tilde R^i_k$}
\left\{
\begin{array}{ll}
y_{r^i_k+p} &\hbox{if $k\neq 0$,}\\
y_{r^i_0+p}^2 &\hbox{if $k=0$,}
\end{array}
\right.
\\
&S^i_k:=\index{s@$S^i_k$}
\left\{
\begin{array}{ll}
y_{s^i_k} &\hbox{if $k\neq 0$,}\\
y_{s^i_0}^2 &\hbox{if $k=0$,}
\end{array}
\right.
&\tilde S^i_k:=\index{s@$\tilde S^i_k$}
\left\{
\begin{array}{ll}
y_{s^i_k+p} &\hbox{if $k\neq 0$,}\\
y_{s^i_0+p}^2 &\hbox{if $k=0$,}
\end{array}
\right.
\end{align*}

Let $i,j,k,l\in J$. We define
\begin{equation}\label{EF}
F^{i,j}_{k,l,\pm}:=(R^i_k\pm\tilde R^j_l)(R^i_k\pm\tilde S^j_l)(S^i_k\pm\tilde R^j_l)(S^i_k\pm\tilde S^j_l),
\index{f@$F^{i,j}_{k,l,\pm}$}
\end{equation}
We also define for $k\neq 0$:
\begin{eqnarray}
F^{i,j}_{\ell,k,-}&:=&(y_1- y_{r^j_k+p}^2)(y_1- y_{s^j_k+p}^2),
\\
F^{i,j}_{k,\ell,-}&:=&(y_{r^i_k}^2- y_{p+1})(y_{s^i_k}^2-y_{p+1}).
\end{eqnarray}
(Note that $F^{i,j}_{\ell,k,-}$ does not depend on $i$ and $F^{i,j}_{k,\ell,-}$ does not depend on $j$). 

Recall the filtration (\ref{EBNPar}). 

\begin{Lemma} \label{LF1F2New}
For all admissible $i,j,k,l$, in $B_{1^2}$, we have 
$$
\ga^{i,j} F^{i,j}_{k,l,\pm}\ga^{i,j}=\ga^{i,j} F^{i,j}_{k,l,\pm}\equiv \ga^{i,j}(z_1^2-z_2^2)^2\pmod{B_{1^2}^{(3)}}.
$$ 
\end{Lemma}
\begin{proof}
We give details for the generic case $k,l\neq \ell$, the exceptional case being similar. 
Denote 
$$
\la:=R^i_k+S^i_k,\ \mu:=R^i_kS^i_k, \ \tilde \la:=\tilde R^i_k+\tilde S^i_k,\ \tilde \mu:=\tilde R^i_k\tilde S^i_k.
$$
Note that
$$
(R^i_k\pm\tilde R^j_l)(R^i_k\pm\tilde S^j_l)(S^i_k\pm\tilde R^j_l)(S^i_k\pm\tilde S^j_l)
=(\mu-\tilde\mu)^2\pm(\mu+\tilde\mu)\la\tilde\la+\mu\tilde\la^2+\la^2\tilde\mu.
$$
Since $\la,\tilde\la,\mu,\tilde \mu$ commute with the idempotent $\ga^{i,j}$ thanks to Corollary~\ref{CSumProduct}, it follows that 
$\ga^{i,j} F^{i,j}_{k,l,\pm}\ga^{i,j}=\ga^{i,j} F^{i,j}_{k,l,\pm}$.

Furthermore, by Corollary~\ref{CSumProduct}(i)(iii), $
\ga^{i,j}\la$ and $\ga^{i,j}\tilde\la$ belong to $B_{1^2}^{(0)}$, and   by Corollary~\ref{CSumProduct}(ii),(iv), we have  
$$
\ga^{i,j}\mu\equiv -\ga^{i,j}z_1^2\pmod{B_{1^2}^{(1)}},\quad
\ga^{i,j}\tilde\mu\equiv -\ga^{i,j}z_2^2\pmod{B_{1^2}^{(1)}}.
$$
So, using (\ref{EFiltrPar}), we get 
\begin{align*}
\ga^{i,j}F^{i,j}_{k,l,\pm}=&\,\ga^{i,j}(R^i_k\pm\tilde R^j_l)(R^i_k\pm\tilde S^j_l)(S^i_k\pm\tilde R^j_l)(S^i_k\pm\tilde S^j_l)
\\
=&\,\ga^{i,j}\big((\mu-\tilde\mu)^2\pm(\mu+\tilde\mu)\la\tilde\la+\mu\tilde\la^2+\la^2\tilde\mu\big)
\\
\equiv&\,\ga^{i,j}(z_1^2-z_2^2)^2\pmod{B_{1^2}^{(3)}},
\end{align*}
as required.
\end{proof}

\begin{Proposition} \label{PPhi^2} 
We have $\text{\sansmath$\Phi$}^2\equiv -(z_1^2-z_2^2)^{2p}\pmod{B_{1^2}^{(4p-1)}}$. 
\end{Proposition}
\begin{proof}
We provide details for the generic case $\ell>1$, the case $\ell=1$ being similar. It is easy to see using  Lemmas~\ref{LphiQuadr}, \ref{LQThickRed} and the relation (\ref{EDotPastCircle}), 
that for all $i,j\in J$, we have 
$$
\ga^{i,j}\text{\sansmath$\Phi$}^2=-\ga^{i,j}
F^{i,j}_{0,0,+}\cdot
\prod_{k=1}^\ell (F^{i,j}_{k,k-1,-}\cdot F^{i,j}_{k-1,k,-}).
$$
By Lemma~\ref{LF1F2New}, every factor is congruent to $\ga^{i,j}(z_1^2-z_2^2)^2$ modulo $B_{1^2}^{(3)}$, and the result follows from (\ref{EFiltrPar}). 
\end{proof}

\subsection{The element $\uptau$}\label{sec:uptau}

For $m=1,\dots,p$, we set
$$
P_m:=\phi_m\phi_{m+1}\cdots \phi_{m+p-1} \in R_{2\de}
$$
and, for $0\leq t\leq p$, we also denote 
\begin{equation*}\label{E180921}
P_m^{(t)}:=\phi_m\cdots \phi_{m+t-1}\psi_{m+t}\cdots \psi_{m+p-1},
\end{equation*}
i.e. we have replaced the last $p-t$ factors $\phi_l$ in $P_m$ with $\psi_l$. Set $P^{(t)}:=P_p^{(t)}P_{p-1}^{(t)}\cdots P_1^{(t)}$. In particular, note that $P^{(0)} = \psi_{u_1}$.

We define the intertwining element $\uptau\in B_2$ via: 
\begin{equation}\label{ETau}
\uptau\ga^{i,j}=(\upsigma+(-1)^{i}\de_{i,j})\ga^{i,j}\qquad(i,j\in J).\index{t@$\uptau$}
\end{equation}
Note that $\uptau\in B_2^{(0)}$. 
Recalling Corollary~\ref{CBBasis}, the following theorem shows that $\uptau$ is the `coefficient' of the top degree $z$-term of $\text{\sansmath$\Phi$}$. The proof of the theorem relies on Lemmas~\ref{lem:poly_P},~\ref{LKeySameColors} and \ref{LKeyDifferentColors} which require some lengthy computations and therefore are delegated to the Appendix. 

\begin{Theorem} \label{TLeadingTerm}
We have
$$
\text{\sansmath$\Phi$}\equiv(z_1^2-z_{2}^2)^p\uptau\pmod{B_2^{(2p-1)}}.
$$
\end{Theorem}
\begin{proof}
It is assumed throughout this proof that we are working in ${\bar R}_{2\de}$, since we ultimately want to prove a statement about $B_2$.

We need to prove that for any $i,j\in J$ we have  
$$
\text{\sansmath$\Phi$}\ga^{i,j}\equiv(z_1^2-z_2^2)^p(\upsigma+(-1)^i\de_{i,j})\ga^{i,j}
\pmod{B_2^{(2p-1)}}.
$$
We take the definition of $P_m$, $P_m^{(t)}$ and $P^{(t)}$, from the beginning of this subsection and set $P:=P_pP_{p-1}\cdots P_1$. Then 
by (\ref{ERedPhi}) and Lemma~\ref{LPhiGa}, we have 
$$
\text{\sansmath$\Phi$}\ga^{i,j}=\ga^{j,i}P\ga^{i,j}.
$$

Recall the definition of $F^{j,i}_{k,k,-}$ from $\S$\ref{sec:quad_rel}. Repeated applications of Lemma \ref{lem:poly_P} give
$$
\ga^{j,i}P\ga^{i,j}=\ga^{j,i}P^{(p)}\ga^{i,j}=\ga^{j,i}\bigg(\prod_{k\in J}F^{j,i}_{k,k,-}\bigg)\,P^{(1)}\ga^{i,j}.
$$
(Note that in Lemma \ref{lem:poly_P} $t=r^i_g$ or $s^i_g$, hence the second equality above.)

Furthermore, again by Lemma \ref{lem:poly_P},
\begin{align}
\begin{split}\label{algn:two_summands}
\text{\sansmath$\Phi$}\ga^{i,j} = &\ga^{j,i}\bigg(\prod_{k\in J}F^{j,i}_{k,k,-}\bigg)\,P^{(1)}\ga^{i,j} \\
= & \ga^{j,i}\bigg(\prod_{k\in J}F^{j,i}_{k,k,-}\bigg)\,\Upsilon \ga^{i,j} + \ga^{j,i}\bigg(\prod_{k\in J}F^{j,i}_{k,k,-}\bigg)\,(y_1 - y_{p+1})P^{(0)} \ga^{i,j}.
\end{split}
\end{align}

We first deal with the second summand in (\ref{algn:two_summands}). By (\ref{ESigma1}) and Lemma \ref{LF1F2New},
\begin{align*}
\ga^{j,i}\bigg(\prod_{k\in J}F^{j,i}_{k,k,-}\bigg)(y_1-y_{p+1})\, P^{(0)}\ga^{i,j} = \ga^{j,i}\bigg(\prod_{k\in J}F^{j,i}_{k,k,-}\bigg)\ga^{j,i}(y_1-y_{p+1})\ga^{j,i}\upsigma\ga^{i,j}.
\end{align*}
(Here we are also using the fact that $\ga^{i,j}$ commutes with $y_1$ and $y_{p+1}$.) Now, by Lemma~\ref{LF1F2New} and Corollary~\ref{PzFilt}, we have 
\begin{equation}\label{E180921_2}
\ga^{j,i}\bigg(\prod_{k\in J}F^{j,i}_{k,k,-}\bigg)\ga^{j,i}\equiv\ga^{j,i}(z_1^2-z_2^2)^{p-1}\ga^{j,i}
\pmod{B_2^{(2p-3)}}.
\end{equation}
Moreover, by (\ref{ERelB1}) and (\ref{ERelB3}), we have 
$$
\ga^{j,i}(y_1-y_{p+1})\ga^{j,i}\equiv\ga^{j,i}(z_1^2-z_2^2)\ga^{j,i}
\pmod{B_2^{(1)}}.
$$
So, using Corollary~\ref{PzFilt} again, we get 
\begin{align}\label{algn:sec_summand}
\ga^{j,i}\bigg(\prod_{k\in J}F^{j,i}_{k,k,-}\bigg)(y_1-y_{p+1})\, P^{(0)}\ga^{i,j}\equiv (z_1^2-z_2^2)^{p}\upsigma\ga^{i,j}
\pmod{B_2^{(2p-1)}}.
\end{align}

We now deal with the first summand in (\ref{algn:two_summands}).
By Lemmas~\ref{LKeySameColors} and \ref{LKeyDifferentColors}, we have
\begin{align*}
\ga^{j,i}\Upsilon\ga^{i,j}
\equiv(-1)^i\de_{i,j}(z_1^2-z_2^2)\ga^{i,j}\pmod{B_2^{(1)}}. 
\end{align*}
So, in view of (\ref{E180921_2}), we have 
\begin{align}\label{algn:first_summand}
\ga^{j,i}\bigg(\prod_{k\in J}F^{j,i}_{k,k,-}\bigg)\,\Upsilon \ga^{i,j} \equiv (-1)^i\de_{i,j}(z_1^2-z_2^2)^p\ga^{i,j}\pmod{B_2^{(2p-1)}}.
\end{align}
Therefore, taking (\ref{algn:two_summands}), (\ref{algn:sec_summand}) and (\ref{algn:first_summand}) together gives
$$
\text{\sansmath$\Phi$}\ga^{i,j}\equiv (z_1^2-z_2^2)^p(\upsigma+(-1)^i\de_{i,j})\ga^{i,j}
\pmod{B_2^{(2p-1)}},
$$
as desired.
\end{proof}

\subsection{Quadratic and wreath product relations for $\uptau$}

\begin{Proposition} \label{PTauQuad} 
In $B_2$, we have $\uptau^2=1$. 
\end{Proposition}
\begin{proof}
By Proposition~\ref{PPhi^2}, we have 
$$\text{\sansmath$\Phi$}^2\equiv -(z_1^2-z_2^2)^{2p}\pmod{B_2^{(4p-1)}}.$$ 
On the other hand, by Theorem~\ref{TLeadingTerm}, we have 
$$\text{\sansmath$\Phi$}\equiv(z_1^2-z_2^2)^p\uptau\pmod{B_2^{(2p-1)}},$$ 
so we can write 
$$\text{\sansmath$\Phi$}=(z_1^2-z_2^2)^p\uptau+b$$ for some $b\in B_2^{(2p-1)}$. Therefore, modulo $B_2^{(4p-1)}$, we get:
\begin{align*}
-(z_1^2-z_2^2)^{2p}&\equiv\text{\sansmath$\Phi$}^2
\\&=\text{\sansmath$\Phi$}((z_1^2-z_2^2)^p\uptau+b)
\\&=(z_2^2-z_1^2)^p\text{\sansmath$\Phi$}\uptau+\text{\sansmath$\Phi$} b
\\&\equiv -(z_1^2-z_2^2)^p\text{\sansmath$\Phi$}\uptau
\\&= -(z_1^2-z_2^2)^p\big((z_1^2-z_2^2)^p\uptau+b\big)\uptau
\\&\equiv-(z_1^2-z_2)^{2p}\uptau^2,
\end{align*}
where we have used Lemma~\ref{LPhiGa} for the third equality and Proposition~\ref{PzFilt} for the fourth and sixth equalities. Since $\uptau\in B_2^{(0)}$, we also have $\uptau^2\in B_2^{(0)}$, so from Corollary~\ref{CBBasis}, we now deduce that $\uptau^2=1$. 
\end{proof}

\begin{Corollary} \label{CSiQuad} 
Let $i,j\in J$. In $B_2$, we have $\upsigma^2\ga^{i,j}=\ga^{i,j}$ if $i\neq j$ and $\upsigma^2\ga^{i,i}=(-1)^{i+1}2\upsigma\ga^{i,i}$.
\end{Corollary}
\begin{proof}
By definition, $\uptau\ga^{i,j}=(\upsigma+(-1)^{i}\de_{i,j})\ga^{i,j}$. We have proved that $\uptau^2=1$, so $(\upsigma+(-1)^{i}\de_{i,j})^2\ga^{i,j}=\ga^{i,j}$, which yields the claim. 
\end{proof}

In view of (\ref{EParDecB}), we have a subalgebra $\Zig_\ell\otimes \Zig_\ell\subseteq B_{1^2}\subseteq B_2$.

\begin{Proposition} \label{PTauWreath}
Let $\za,\za'\in \Zig_\ell$. Then in $B_2$, we have 
$$\uptau (\za\otimes \za')=(-1)^{|\za||\za'|}(\za'\otimes \za)\uptau.$$
\end{Proposition}
\begin{proof}
The proof is similar to that of Proposition~\ref{PTauQuad}. 
By Theorem~\ref{TLeadingTerm}, we 
can write 
$\text{\sansmath$\Phi$}=(z_1^2-z_2^2)^p\uptau+b$ for some $b\in B_2^{(2p-1)}$. Now, modulo $B_2^{(2p-1)}$, 
\begin{align*}
(z_1^2-z_2^2)^p(\za'\otimes \za)\uptau
&=
(\za'\otimes \za)(z_1^2-z_2^2)^p\uptau
\\&\equiv 
(\za'\otimes \za)((z_1^2-z_2^2)^p\uptau+b)
\\&=
(\za'\otimes \za)\text{\sansmath$\Phi$}
\\&=(-1)^{|\za||\za'|}
\text{\sansmath$\Phi$}(\za\otimes \za')
\\&=(-1)^{|\za||\za'|}((z_1^2-z_2^2)^p\uptau+b)(\za\otimes \za')
\\&\equiv (-1)^{|\za||\za'|}(z_1^2-z_2^2)^p\uptau(\za\otimes \za'),
\end{align*}
where the first equality follows from the fact that $\zz^2$ is central in $\Zig_\ell[\zz]$, and the fourth equality comes from Lemma~\ref{LPhiGa}. From Corollary~\ref{CBBasis}, we now deduce that $\uptau (\za\otimes \za')=(-1)^{|\za||\za'|}(\za'\otimes \za)\uptau$.
\end{proof}

\subsection{Rough $(\uptau,z)$ commutation relation}
We will later prove the `sharp' $(\uptau,z)$ commutation relation (\ref{ETauZIdTarget}). 
Unfortunately, this will require a lengthy computation. So for now, in the following proposition we show less, which will still suffice to prove our main result on RoCK blocks.

\begin{Proposition}\label{PTauZRough} 
In $B_2$, we have $\uptau z_1-z_2\uptau\in \Zig_\ell\otimes \Zig_\ell$ and $\uptau z_2-z_1\uptau\in \Zig_\ell\otimes \Zig_\ell$. 
\end{Proposition}
\begin{proof}
It suffices to prove the first containment. 
The second containment then follows from the first one by multiplying with $\uptau$ on the left and on the right and using Propositions~\ref{PTauQuad} and \ref{PTauWreath}. 

Since $\Zig_\ell\otimes \Zig_\ell=B_2^{(0)}\cap B_{1^2}$ by Corollary~\ref{CBBasis}, it suffices to prove that $\uptau z_1-z_2\uptau\in B_2^{(0)}$ and $\uptau z_1-z_2\uptau\in B_{1^2}$. 

We first prove that $\uptau z_1-z_2\uptau\in B_2^{(0)}$. 
By Theorem~\ref{TLeadingTerm}, we 
can write 
$$\text{\sansmath$\Phi$}=(z_1^2-z_2^2)^p\uptau+b$$ 
for some $b\in B_2^{(2p-1)}$. By Lemma~\ref{LPhiGa}, we have $\text{\sansmath$\Phi$} z_1=z_2\text{\sansmath$\Phi$}$. So 
$$
(z_1^2-z_2^2)^p\uptau z_1+bz_1=z_2(z_1^2-z_2^2)^p\uptau+z_2b.
$$
Note that $\uptau z_1$ and $z_2\uptau$ belong to $B_2^{(1)}$. So by Corollary~\ref{CBBasis}, we can write  
$\uptau z_1=z_1\al+z_2\be+\ga$ 
for some $\al,\be,\ga\in B_2^{(0)}$.  Hence from the previous equation we get 
$$
z_1(z_1^2-z_2^2)^p\al
+z_2(z_1^2-z_2^2)^p\be+
(z_1^2-z_2^2)^p\ga+
bz_1=z_2(z_1^2-z_2^2)^p\uptau+z_2b.
$$
In view of Corollary~\ref{CBBasis}, we have $\al=0$ and $\be=\uptau$. This proves the containment $\uptau z_1-z_2\uptau\in B_2^{(0)}$. 

To prove that $\uptau z_1-z_2\uptau\in B_{1^2}$, note by  (\ref{EZFormula}) and (\ref{ETau}) that for any $i,j\in J$:
\begin{align*}
\uptau z_1\ga^{i,j}&=(\upsigma+(-1)^i\de_{i,j})y_{p-i-1}y_{p-i}\ga^{i,j}
\\&=\upsigma y_{p-i-1}y_{p-i}\ga^{i,j}+(-1)^i\de_{i,j}y_{p-i-1}y_{p-i}\ga^{i,j}.
\end{align*}
Note that the second summand is in $B_{1^2}$. As for the first summand, use the defining relations in $R_{2\de}$ to pull $y_{p-i-1}$ and $y_{p-i}$ past $\upsigma$ to get 
\begin{align*}
\upsigma y_{p-i-1}y_{p-i}\ga^{i,j}&=
y_{2p-i-1}y_{2p-i}\upsigma \ga^{i,j}+(*)
\\&=z_2\upsigma\ga^{i,j}+(*)\\&=z_2\uptau\ga^{i,j}-(-1)^i\de_{i,j}z_2\ga^{i,j}+(*),
\end{align*}
where the error term $(*)$ belongs to $B_{1^2}$ by  Corollary~\ref{CBruhat}, and $(-1)^i\de_{i,j}z_2\ga^{i,j}$ is clearly in $B_{1^2}$, too.
\end{proof}

\subsection{The elements $\text{\sansmath$\Phi$}_r$ and $\uptau_r$ and braid relations}
Recalling the embeddings (\ref{EIotaRR+1}), we now define for all $r=1,\dots,d-1$:
\begin{equation}\label{EPhiTauR}
\text{\sansmath$\Phi$}_r:=\iota_{r,r+1}(\text{\sansmath$\Phi$})\in B_{d}\qquad\text{and}\qquad
\uptau_r:=\iota_{r,r+1}(\uptau).
\index{f@$\text{\sansmath$\Phi$}_r$}\index{t@$\uptau_r$}
\end{equation}
In view of (\ref{EPhiBig}), we have 
$$\text{\sansmath$\Phi$}_r=\ga_{1^d}\phi_{u_r}\ga_{1^d}$$  
where $\phi_{u_r}$ is defined using the choice of reduced decomposition for $u_r$ similar to (\ref{ERedU1}). 
In view of (\ref{ETau}), we have 
\begin{equation}\label{ETauR}
\uptau\ga^{\bj}=(\upsigma_r+(-1)^{j_r}\de_{j_r,j_{r+1}})\ga^{\bj}\qquad(1\leq r<d,\ \bj\in J^d). 
\end{equation}

It is clear that 
\begin{eqnarray}\label{ETauFarApart}
\uptau_r\uptau_s=\uptau_s\uptau_r\qquad(|r-s|>1),
\\
\uptau_r z_s=z_s\uptau_r\qquad(s\neq r,r+1).
\label{ETauZFarApart}
\end{eqnarray}
Moreover, recalling the notation (\ref{EWSignAction}), 
by Propositions~\ref{PTauQuad}, \ref{PTauWreath} and \ref{PTauZRough}, using $\iota_{r,r+1}$, we deduce for all admissible $r,s$ and all $a_1,\dots,a_d\in\Zig_\ell$:
\begin{eqnarray}
&\uptau_r^2=1,
\label{ETauRQuad}
\\
&\uptau_r(a_1\otimes\dots \otimes a_d)={}^{s_r}(a_1\otimes\dots\otimes a_d)\uptau_r,
\label{ETauRWreath}
\\
\label{ETauRZRough}
&\uptau_rz_r-z_{r+1}\uptau_r\in \Zig_\ell^{\otimes d}\quad \text{and}\quad  \uptau_rz_{r+1}-z_{r}\uptau_r\in \Zig_\ell^{\otimes d}.
\end{eqnarray} 

\begin{Proposition} \label{PTauBraid}
For all\, $1\leq r\leq d-2$, we have $\uptau_r\uptau_{r+1}\uptau_r=\uptau_{r+1}\uptau_r\uptau_{r+1}$.
\end{Proposition}
\begin{proof}
Using Corollary~\ref{C180621}, we deduce 
$$\text{\sansmath$\Phi$}_r\text{\sansmath$\Phi$}_{r+1}\text{\sansmath$\Phi$}_r=\text{\sansmath$\Phi$}_{r+1}\text{\sansmath$\Phi$}_r\text{\sansmath$\Phi$}_{r+1}.$$ 
By Theorem~\ref{TLeadingTerm}, for every $1\leq s<d$, 
we can write 
$$\text{\sansmath$\Phi$}_s=(z_s^2-z_{s+1}^2)^p\uptau_s+b_s$$ 
for some $b_s\in B_d^{(2p-1)}$. 
Now, using Lemma~\ref{LPhiGa}, modulo $B_d^{(6p-1)}$, we get:
\begin{align*}
\text{\sansmath$\Phi$}_{r+1}\text{\sansmath$\Phi$}_r\text{\sansmath$\Phi$}_{r+1}
&\equiv\text{\sansmath$\Phi$}_{r+1}\text{\sansmath$\Phi$}_r(z_{r+1}^2-z_{r+2}^2)^p\uptau_{r+1}
\\&=(z_{r}^2-z_{r+1}^2)^p\text{\sansmath$\Phi$}_{r+1}\text{\sansmath$\Phi$}_r\uptau_{r+1}
\\&\equiv(z_{r}^2-z_{r+1}^2)^p\text{\sansmath$\Phi$}_{r+1}(z_r^2-z_{r+1}^2)^p\uptau_r\uptau_{r+1}
\\&=(z_{r}^2-z_{r+1}^2)^p(z_r^2-z_{r+2}^2)^p\text{\sansmath$\Phi$}_{r+1}\uptau_r\uptau_{r+1}
\\&=(z_{r}^2-z_{r+1}^2)^p(z_r^2-z_{r+2}^2)^p(z_{r+1}^2-z_{r+2}^2)^p\uptau_{r+1}\uptau_r\uptau_{r+1}
\end{align*}
Similarly:
$$
\text{\sansmath$\Phi$}_r\text{\sansmath$\Phi$}_{r+1}\text{\sansmath$\Phi$}_r\equiv (z_{r}^2-z_{r+1}^2)^p(z_r^2-z_{r+2}^2)^p(z_{r+1}^2-z_{r+2}^2)^p\uptau_{r}\uptau_{r+1}\uptau_{r}. 
$$
From Corollary~\ref{CBBasis}, we now deduce that $\uptau_r\uptau_{r+1}\uptau_r=\uptau_{r+1}\uptau_r\uptau_{r+1}$.
\end{proof}

We have proved that $\uptau_r$'s satisfy braid relations, so for every $w\in \Si_d$ with reduced decomposition $w=s_{r_1}\cdots s_{r_m}$, there is a well defined element
$$
\uptau_w:=\uptau_{r_1}\dots\uptau_{r_m}.\index{t@$\uptau_w$}
$$

\begin{Lemma} \label{LTauSi} 
For any $w\in \Si_d$, we have that $\uptau_w=\upsigma_w+\sum_{u<w}c_u\upsigma_u$, where the coefficients $c_u\in \F$. 
\end{Lemma}
\begin{proof}
This follows from the definition of $\uptau$ and Corollary~\ref{CSiQuad}. 
\end{proof}

\subsection[The isomorphism $Y_{\rho,d}\cong \Zig_\ell\swr \Si_d$]{The isomorphism $Y_{\rho,d}\cong \Zig_\ell\swr \Si_d$ and the Morita equivalence $R^{\La_0}_{\cont(\rho)+d\de}\sim \Zig_\ell\swr \Si_d$}

\begin{Proposition} \label{PBWreathIso} 
There is an isomorphism of graded superalgebras 
$$\bar f_d: \Zig_\ell\swr \Si_d\iso B_d^{(0)},\ (a_1\otimes \dots\otimes a_d)\otimes w\mapsto (a_1\otimes \dots\otimes a_d)\uptau_w.\index{f@$\bar f_d$}
$$
\end{Proposition}
\begin{proof}
By Lemma~\ref{LTauSi} and Corollary~\ref{CBBasis}, the map in the lemma is an isomorphism of vector spaces. It is an algebra isomorphism thanks to the relations (\ref{ETauFarApart}), (\ref{ETauRQuad}), 
(\ref{ETauRWreath}) and braid relations of Proposition~\ref{PTauBraid}. 
\end{proof}

Recall the surjection $\Om_{\rho,d}:B_d\to Y_{\rho,d}$ from (\ref{ETrOmd}) and the condition $QF2$ from \S\ref{SSMainD=1}. 

\begin{Theorem}\label{thm:main}
Let $\rho$ be a $d$-Rouquier $\bar p$-core. 
Suppose that $R^{\La_0}_{\cont(\rho)+\de}$ is a $QF2$-algebra. 
The map $\Om_{\rho,d}\circ \bar f_d: \Zig_\ell\swr \Si_d\to Y_{\rho,d}$ is a graded superalgebra isomorphism. 
\end{Theorem}
\begin{proof}
By Proposition~\ref{PDimYd}, $\Zig_\ell\swr \Si_d$ and $Y_{\rho,d}$ have the same dimension, so it suffices to prove that  the map in the statement of the theorem is surjective. By Proposition~\ref{PBWreathIso}, this means $\Om_{\rho,d}(B_d)=\Om_{\rho,d}(B_d^{(0)})$, i.e.$\Om_{\rho,d}(z_r)\in \Om_{\rho,d}(B_d^{(0)})$ for all $r=1,\dots,d$. We prove this by induction on $r$. 
 
Suppose first that $r=1$. If $d=1$, then by Theorem~\ref{TMaind=1}, we have that $\Om_{\rho,1}(z)=\Om_{\rho,1}(f(a))$ for some $a\in\Zig_\ell\subseteq B_1$. 
For an arbitrary $d$, we have $z_1=\iota_1(z)$, so using  
Lemma~\ref{LY_1Y_d}, we deduce that 
$$\Om_{\rho,d}(z_1)=\Om_{\rho,d}(\iota_1(z))=\zeta(\Om_{\rho,1}(a))=\Om_{\rho,d}(\iota_1(a))\in \Om_{\rho,d}(B_d^{(0)}).
$$
 
 Let $r>1$. Then by (\ref{ETauRZRough}), we have $z_r\equiv\uptau_rz_{r-1}\uptau_r\pmod{B_d^{(0)}}$, and so $\Om_{\rho,d}(z_r)\in \Om_{\rho,d}(B_d^{(0)})$ by the inductive assumption.
\end{proof}

\begin{Corollary} \label{CMorWr} 
Let $\rho$ be a $d$-Rouquier $\bar p$-core. 
Suppose that $R^{\La_0}_{\cont(\rho)+\de}$ is a $QF2$-algebra. Then the algebra $Y_{\rho,d}$ is graded Morita superequivalent to the RoCK block $R^{\La_0}_{\cont(\rho)+d\de}$ if and only if $\cha \F>d$.
\end{Corollary}
\begin{proof}
By Lemma~\ref{LNumberOfIrrCyc}, we have $|\Irr(R^{\La_0}_{\cont(\rho)+d\de})|=|\Par^\ell(d)|$.  
On the other hand, by Theorem~\ref{thm:main}, we have 
$Y_{\rho,d}\cong \Zig_\ell\swr \Si_d$. 
Moreover, we also have $|\Irr(\Zig_\ell\swr \Si_d)|=|\Par^\ell(d)|$ if 
$\cha \F>d$, and $|\Irr(\Zig_\ell\swr \Si_d)|<|\Par^\ell(d)|$ otherwise.  So it remains to recall from (\ref{EIdTrd}) that $Y_{\rho,d}$ is an idempotent truncation of $R^{\La_0}_{\cont(\rho)+d\de}$ and use 
Corollary~\ref{cor:idmpt_Mor}. 
\end{proof}

From Theorem~\ref{thm:main} and Corollary~\ref{CMorWr}, we get:

\begin{Theorem} \label{TMorWr} 
Let $\rho$ be a $d$-Rouquier $\bar p$-core. 
Suppose that $R^{\La_0}_{\cont(\rho)+\de}$ is a $QF2$-algebra. Then the RoCK block $R^{\La_0}_{\cont(\rho)+d\de}$ is graded Morita superequivalent to $\Zig_\ell\swr \Si_d$ if and only if $\cha \F>d$. 
\end{Theorem}

\section{Affine zigzag relations}
Let $H_d(\Zig_\ell)$ be the rank $d$ Brauer tree affine Hecke algebra corresponding to $\Zig_\ell$ defined in \S\ref{SSAff}. Recall the notation (\ref{EInsertion}), (\ref{EAR}). 
In this section we extend the isomorphism 
$\bar f_d$ from Proposition~\ref{PBWreathIso} to an isomorphism 
$f_d: H_d(\Zig_\ell)\iso B_{d}$ which maps
$
\zz_r\mapsto z_r$ for all $1\leq r\leq d$. 
In order to prove that such a homomorphism $f_d$ exists, we just need to check the relation 
\begin{equation}\label{ETauZIdTarget}
(\uptau_r z_t - z_{s_r(t)} \uptau_r)\ga^\bi
=
\begin{cases}
\big((\delta_{r,t}- \delta_{r+1,t})(c_r + c_{r+1})
&
\\
\qquad\qquad\qquad+\de_{i_r,0} u_r u_{r+1}\big)\ga^\bi
&
\text{if}\ i_r = i_{r+1},\\
(\delta_{r,t}- \delta_{r+1,t}) a_r^{i_{r+1},i_r} a_{r+1}^{i_r,i_{r+1}} \ga^\bi
&
\text{if}\ |i_r-i_{r+1}|=1,\\
0
& 
\textup{otherwise}.
\end{cases}
\end{equation}
which corresponds to the defining relation (\ref{ESZId}). 
Then the fact that $f_d$ is an isomorphism follows immediately from Corollary~\ref{CBBasis} and Theorem~\ref{TAffBasis}. 

\subsection{Some reductions}
\label{SSRed}
By (\ref{ETauZFarApart}), we have $\uptau_rz_t=z_t\uptau_r$ for $t\neq r,r+1$. So in proving (\ref{ETauZIdTarget}) we may assume that $t=r$ or $r+1$. 

Suppose we have proved the relation (\ref{ETauZIdTarget}) for the case $t=r$, i.e. we have 
\begin{equation*}
(\uptau_r z_r - z_{r+1} \uptau_r)\ga^\bi=
\begin{cases}
(c_r + c_{r+1}+\de_{i_r,0}u_ru_{r+1})\ga^\bi
&
\text{if}\ i_r = i_{r+1},\\
a_r^{i_{r+1},i_r} a_{r+1}^{i_r,i_{r+1}} \ga^\bi
&
\text{if}\ |i_r-i_{r+1}|=1,\\
0
& 
\textup{otherwise}.
\end{cases}
\end{equation*}
Multiplying on both sides with $\uptau_r$ and using the relations (\ref{ETauRQuad}) and (\ref{ETauRWreath}), we get 

\begin{equation*}
(z_r\uptau_r - \uptau_rz_{r+1})\ga^{s_r\cdot\bi}=
\begin{cases}
(c_{r+1} + c_{r}+\de_{i_r,0}u_{r+1}u_{r})\ga^{s_r\cdot \bi}
&
\text{if}\ i_r = i_{r+1},\\
a_{r+1}^{i_{r+1},i_r} a_{r}^{i_r,i_{r+1}} \ga^{s_r\cdot \bi}
&
\text{if}\ |i_r-i_{r+1}|=1,\\
0
& 
\textup{otherwise}.
\end{cases}
\end{equation*}
Since $u_{r+1}u_{r}=-u_{r}u_{r+1}$, these relations 
are equivalent to the relations (\ref{ETauZIdTarget}) with $t=r+1$ (as $\bi$ runs over all $J^d$). 
So we are reduced to the case $t=r$. Moreover, using the embedding $\iota_{r,r+1}$ we may assume that $t=r=1$ and $d=2$, i.e. it suffices to prove 
\begin{equation}\label{ETauZSharp}
(\uptau z_1-z_2\uptau)\ga^{i,j}=
\left\{
\begin{array}{ll}
(c_1+c_2+\de_{i,0}u_1u_2)\ga^{i,j} &\hbox{if $i=j$,}\\
a_1^{j,i}a_2^{i,j}\ga^{i,j} &\hbox{if $|i-j|=1$,}\\
0 &\hbox{otherwise.}
\end{array}
\right.
\end{equation}

Now recall from Proposition~\ref{PTauZRough} that $\uptau z_1-z_2\uptau\in \Zig_\ell\otimes \Zig_\ell$. So for any $i,j\in J$, we have, using (\ref{ETauRWreath}) 
\begin{equation}\label{E080721}
(\uptau z_1-z_2\uptau)\ga^{i,j}=\ga^{j,i}(\uptau z_1-z_2\uptau)\ga^{i,j}\in e^j\Zig_\ell e^i\otimes e^i\Zig_\ell e^j.
\end{equation}
Note that $|i-j|>1$ implies $e^i\Zig_\ell e^j=0$, which yields the third  case of (\ref{ETauZSharp}). Suppose we are in the second case, i.e. $|i-j|=1$. Then $e^i\Zig_\ell e^j$ is spanned by $a^{i,j}$, so it follows from (\ref{E080721}) that 
\begin{equation}\label{E080721_2}
(\uptau z_1-z_2\uptau)\ga^{i,j}=Ca_1^{j,i}a_2^{i,j}
\end{equation} 
for some constant $C$. We show that $C=1$ provided 
\begin{equation}\label{E080721_1}
(\uptau z_1-z_2\uptau)\ga^{i,i}=
(c_1+c_2+\de_{i,0}u_1u_2)\ga^{i,i}
\end{equation}
thus reducing proving (\ref{ETauZSharp}) just to the first case $i=j$. Indeed, multiplying (\ref{E080721_1}) on the left with $a_1^{j,i}$, we get 
$$a_1^{j,i}(\uptau z_1-z_2\uptau)\ga^{i,i}=
a_1^{j,i}(c_1+c_2+\de_{i,0}u_1u_2)\ga^{i,i}.
$$
By (\ref{ETauRWreath}) we have 
$a_1^{j,i}\uptau\ga^{i,i}=\uptau \ga^{i,j} a_2^{j,i}\ga^{i,i}$. On the other hand,  $a_1^{j,i}c_1=a_1^{j,i}u_1=0$, so we deduce
$$
(\uptau z_1-z_2\uptau)a_2^{j,i}\ga^{i,i}=
a_1^{j,i}c_2\ga^{i,i}.
$$
On the other hand, by (\ref{E080721_2}), we have 
$$
(\uptau z_1-z_2\uptau)a_2^{j,i}\ga^{i,i}=(\uptau z_1-z_2\uptau)\ga^{i,j}a_2^{j,i}\ga^{i,i}=Ca_1^{j,i}a_2^{i,j}\ga^{i,j}a_2^{j,i}
=Ca_1^{j,i}c_2\ga^{i,i}.
$$
Comparing the last two equations and using the fact that $a_1^{j,i}c_2\ga^{i,i}\neq 0$ thanks to Corollary~\ref{CBBasis}, we deduce that $C=1$. 

\subsection{The case $i=j$} 
\label{SSI=J}
We have observed in the previous subsection that to prove (\ref{ETauZSharp}) in full generality it suffices to prove just the $i=j$ case of it, which is (\ref{E080721_1}). Note using (\ref{ETau}) that (\ref{E080721_1}) is equivalent to 
\begin{equation}\label{E080721_4}
\ga^{i,i}\upsigma z_1\ga^{i,i}=\ga^{i,i}(z_2\upsigma+(-1)^{i+1}z_1+(-1)^iz_2+c_1+c_2+\de_{i,0}u_1u_2)\ga^{i,i}.
\end{equation}

\begin{Proposition} 
The equality (\ref{E080721_4}) hold for every $i\in J$. 
\end{Proposition}
\begin{proof}
We give details for the generic case $i\neq 0$, the special case $i=0$ being similar. In terms of diagrams, we have
$$
\ga^{i,i}\upsigma z_1\ga^{i,i}=
\resizebox{106mm}{20mm}{
\begin{braid}\tikzset{baseline=-.3em}
	\draw (0,6) node{\color{blue}$\ell$\color{black}};
	\braidbox{0.6}{2.9}{5.3}{6.6}{};
	\draw (1.8,6) node{$(\ell-1)^2$};
	\draw[dots] (3.5,6)--(5,6);
	\braidbox{5.3}{7.6}{5.3}{6.6}{};
	\draw (6.5,6) node{$(j+1)^2$};
	\draw (8.3,6) node{$j$};
	\draw (9.5,6) node{$j-1$};
	\draw[dots] (10.8,6)--(12.5,6);
	\draw (12.7,6) node{$1$};
	\redbraidbox{13.4}{14.3}{5.3}{6.6}{};
	\draw (13.8,6) node{$\color{red}0\,\,0\color{black}$};
	\draw (15,6) node{$1$};
	\draw[dots] (15.6,6)--(17,6);
	\draw (18.2,6) node{$j-1$};
	\draw (19.5,6) node{$j$};
	\draw (0,-6) node{\color{blue}$\ell$\color{black}};
	\braidbox{0.6}{2.9}{-6.7}{-5.4}{};
	\draw (1.8,-6) node{$(\ell-1)^2$};
	\draw[dots] (3.5,-6)--(5,-6);
	\braidbox{5.3}{7.6}{-6.7}{-5.4}{};
	\draw (6.5,-6) node{$(i+1)^2$};
	\draw (8.3,-6) node{$i$};
	\draw (9.5,-6) node{$i-1$};
	\draw[dots] (10.8,-6)--(12.5,-6);
	\draw (12.7,-6) node{$1$};
	\redbraidbox{13.4}{14.3}{-6.7}{-5.4}{};
	\draw (13.8,-6) node{$\color{red}0\,\,0\color{black}$};
	\draw (15,-6) node{$1$};
	\draw[dots] (15.6,-6)--(17,-6);
	\draw (18.2,-6) node{$i-1$};
	\draw (19.5,-6) node{$i$};

	\draw (20.3,6) node{\color{blue}$\ell$\color{black}};
	\braidbox{20.9}{23.2}{5.3}{6.6}{};
	\draw (22.1,6) node{$(\ell-1)^2$};
	\draw[dots] (23.8,6)--(25.3,6);
	\braidbox{25.6}{27.9}{5.3}{6.6}{};
	\draw (26.8,6) node{$(i+1)^2$};
	\draw (28.6,6) node{$i$};
	\draw (29.8,6) node{$i-1$};
	\draw[dots] (31.1,6)--(32.8,6);
	\draw (33,6) node{$1$};
	\redbraidbox{33.7}{34.6}{5.3}{6.6}{};
	\draw (34.1,6) node{$\color{red}0\,\,0\color{black}$};
	\draw (35.3,6) node{$1$};
	\draw[dots] (35.9,6)--(37.3,6);
	\draw (38.5,6) node{$i-1$};
	\draw (39.8,6) node{$i$};
	\draw (20.3,-6) node{\color{blue}$\ell$\color{black}};
	\braidbox{20.9}{23.2}{-6.7}{-5.4}{};
	\draw (22.1,-6) node{$(\ell-1)^2$};
	\draw[dots] (23.8,-6)--(25.1,-6);
	\braidbox{25.6}{27.9}{-6.7}{-5.4}{};
	\draw (26.8,-6) node{$(j+1)^2$};
	\draw (28.6,-6) node{$j$};
	\draw (29.8,-6) node{$j-1$};
	\draw[dots] (31.1,-6)--(32.8,-6);
	\draw (33,-6) node{$1$};
	\redbraidbox{33.7}{34.6}{-6.7}{-5.4}{};
	\draw (34.1,-6) node{$\color{red}0\,\,0\color{black}$};
	\draw (35.3,-6) node{$1$};
	\draw[dots] (35.9,-6)--(37.3,-6);
	\draw (38.5,-6) node{$j-1$};
	\draw (39.8,-6) node{$j$};
	\draw[blue](0,-5.2)--(20.3,5.2);
	\draw(1.2,-5.2)--(21.3,5.2);
	\draw(2.2,-5.2)--(22.3,5.2);
	\draw(6,-5.2)--(26,5.2);
	\draw(7,-5.2)--(27,5.2);
	\draw(8.3,-5.2)--(28.6,5.2);
	\draw(9.5,-5.2)--(29.8,5.2);
	\draw(12.7,-5.2)--(33.1,5.2);
	\draw[red](13.5,-5.2)--(33.8,5.2);
	\draw[red](14.2,-5.2)--(34.5,5.2);
	\draw(15,-5.2)--(35.3,5.2);
	\draw(18.2,-5.2)--(38.5,5.2);
	\draw(39.8,-5.2)--(19.5,5.2);
	\draw(39.8,5.2)--(19.5,-5.2);
	\draw[blue] (20.3,-5.2)--(0,5.2);
	\draw(21.3,-5.2)--(1.2,5.2);
	\draw(22.3,-5.2)--(2.2,5.2);
	\draw(26,-5.2)--(6,5.2);
	\draw(27,-5.2)--(7,5.2);
	\draw(28.6,-5.2)--(8.3,5.2);
	\draw(29.8,-5.2)--(9.5,5.2);
	\draw(38.5,-5.2)--(18.2,5.2);
	\draw(35.3,-5.2)--(15,5.2);
	\draw[red](34.5,-5.2)--(14.2,5.2);
	\draw[red](33.8,-5.2)--(13.5,5.2);
	\draw(33.1,-5.2)--(12.7,5.2);
	\reddot(14.1,-4.85);
	\reddot(14.5,-5.05);
	\end{braid}.
	}
$$
Using Lemma~\ref{L090721}, braid relations (\ref{R7}) and dot-crossing relations (\ref{R5}), we now get 
$$
\ga^{i,i}\upsigma z_1\ga^{i,i}=\ga^{i,i}z_2\upsigma\ga^{i,i}+(*),
$$
where
\begin{align*}
(*)&=
\begin{braid}\tikzset{baseline=-.3em}
	\draw (0,4) node{\color{blue}$\ell$\color{black}};
	\braidbox{0.6}{2.9}{3.3}{4.6}{};
	\draw (1.8,4) node{$(\ell-1)^2$};
	\draw[dots] (3.5,4)--(5,4);
	\braidbox{5.3}{7.6}{3.3}{4.6}{};
	\draw (6.5,4) node{$(i+1)^2$};
	\draw (8.3,4) node{$i$};
	\draw[dots] (9,4)--(10.3,4);
	\draw[dots] (9,-4)--(10.3,-4);
	\draw (11,4) node{$1$};
	\draw (11,-4) node{$1$};
	\redbraidbox{11.7}{12.6}{3.3}{4.6}{};
	\draw (12.1,4) node{$\color{red}0\,\,0\color{black}$};
	\redbraidbox{11.7}{12.6}{-3.3}{-4.6}{};
	\draw (12.1,-4) node{$\color{red}0\,\,0\color{black}$};
	\draw (13.3,4) node{$1$};
	\draw (13.3,-4) node{$1$};
	\draw[dots] (13.9,4)--(15.3,4);
	\draw[dots] (13.9,-4)--(15.3,-4);
	\draw (16,4) node{$i$};
\draw (16,-4) node{$i$};
	\draw (0,-4) node{\color{blue}$\ell$\color{black}};
	\braidbox{0.6}{2.9}{-4.7}{-3.4}{};
	\draw (1.8,-4) node{$(\ell-1)^2$};
	\draw[dots] (3.5,-4)--(5,-4);
	\braidbox{5.3}{7.6}{-4.7}{-3.4}{};
	\draw (6.5,-4) node{$(i+1)^2$};
	\draw (8.3,-4) node{$i$};
	\draw (16.8,4) node{\color{blue}$\ell$\color{black}};
\draw (16.8,-4) node{\color{blue}$\ell$\color{black}};
	\braidbox{17.4}{19.7}{3.3}{4.6}{};
	\draw (18.6,4) node{$(\ell-1)^2$};
	\braidbox{17.4}{19.7}{-3.3}{-4.6}{};
	\draw (18.6,-4) node{$(\ell-1)^2$};
	\draw[dots] (20.3,4)--(21.8,4);
	\draw[dots] (20.3,-4)--(21.8,-4);
	\braidbox{22.2}{24.5}{3.3}{4.6}{};
	\draw (23.4,4) node{$(i+1)^2$};
	\braidbox{22.2}{24.5}{-3.3}{-4.6}{};
	\draw (23.4,-4) node{$(i+1)^2$};
	\draw (25.2,4) node{$i$};
	\draw (25.2,-4) node{$i$};
	\draw[dots] (25.9,4)--(27.6,4);
	\draw[dots] (25.9,-4)--(27.6,-4);
	\draw (27.9,4) node{$1$};
	\draw (27.9,-4) node{$1$};
	\redbraidbox{28.5}{29.4}{3.3}{4.6}{};
	\draw (28.9,4) node{$\color{red}0\,\,0\color{black}$};
	\redbraidbox{28.5}{29.4}{-3.3}{-4.6}{};
	\draw (28.9,-4) node{$\color{red}0\,\,0\color{black}$};
	\draw (30.1,4) node{$1$};
	\draw (30.1,-4) node{$1$};
	\draw[dots] (30.7,4)--(32.1,4);
	\draw[dots] (30.7,-4)--(32.1,-4);
	\draw (32.8,4) node{$i$};
	\draw (32.8,-4) node{$i$};
	\draw[blue](0,-3.2)--(16.8,3.2);
	\draw[blue](0,3.2)--(16.8,-3.2);
	\draw(1.3,-3.2)--(17.9,3.2);
	\draw(1.3,3.2)--(17.9,-3.2);
	\draw(2.7,-3.2)--(19.3,3.2);
	\draw(2.7,3.2)--(19.3,-3.2);
	\draw(5.8,-3.2)--(22.4,3.2);
	\draw(5.8,3.2)--(22.4,-3.2);
	\draw(7.3,-3.2)--(23.9,3.2);
	\draw(7.3,3.2)--(23.9,-3.2);
	\draw(8.3,-3.2)--(25.2,3.2);
	\draw(25.2,-3.2)--(8.3,3.2);
	\draw(11,-3.2)--(27.9,3.2);
	\draw(27.9,-3.2)--(11,3.2);
	\draw[red](11.8,-3.2)--(21.9,-0.1)--(12.5,2.8)--(12.5,3.2);
	\draw[red](12.5,-3.2)--(23,-0.1)--(29.2,3.2);
	\draw[red](28.6,-3.2)--(23,-0.1)--(11.8,3.2);
	\draw[red](29.3,-3.2)--(23,0.4)--(28.6,3.2);
	\draw(13.3,-3.2)--(25.5,0)--(30.1,3.2);
	\draw(13.3,3.2)--(25.5,0)--(30.1,-3.2);
	\draw(32.8,-3.2)--(29,0)--(16,3.2);
	\draw(32.8,3.2)--(29,0)--(16,-3.2);
	\reddot(12.9,-3.1);
	\end{braid}.
\end{align*}
Using the braid relation (\ref{R7}) for 
$\begin{braid}\tikzset{baseline=.4em}
	\draw (0,0) node{$1$};
        \draw[red](1.1,0) node{$0$};
        \draw (2.2,0) node{$1$};
		\draw(0,0.5)--(2.2,1.5);
\draw(2.2,0.5)--(0,1.5);
\draw[red](1.1,0.5)--(2,1)--(1.1,1.5);
	\end{braid}$ 
 and the fact that the word beginning with $\ell(\ell-1)^2\cdots(i+2)^2i\cdots 20$ is not cuspidal by Lemma~\ref{LCuspExpl}, we get 
 \begin{align*}
(*)&=
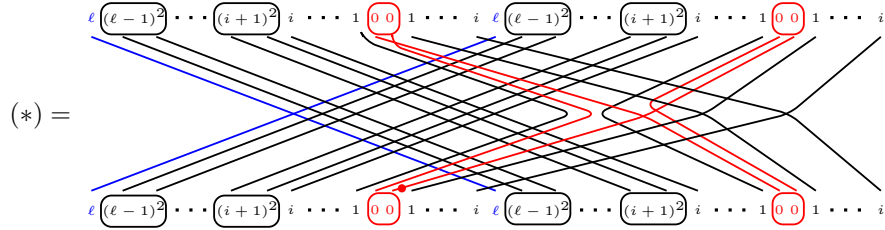
\begin{braid}\tikzset{baseline=-.3em}
	\draw (0,4) node{\color{blue}$\ell$\color{black}};
	\braidbox{0.6}{2.9}{3.3}{4.6}{};
	\draw (1.8,4) node{$(\ell-1)^2$};
	\draw[dots] (3.5,4)--(5,4);
	\braidbox{5.3}{7.6}{3.3}{4.6}{};
	\draw (6.5,4) node{$(i+1)^2$};
	\draw (8.3,4) node{$i$};
	\draw[dots] (9,4)--(10.3,4);
	\draw[dots] (9,-4)--(10.3,-4);
	\draw (11,4) node{$1$};
	\draw (11,-4) node{$1$};
	\redbraidbox{11.7}{12.6}{3.3}{4.6}{};
	\draw (12.1,4) node{$\color{red}0\,\,0\color{black}$};
	\redbraidbox{11.7}{12.6}{-3.3}{-4.6}{};
	\draw (12.1,-4) node{$\color{red}0\,\,0\color{black}$};
	\draw (13.3,4) node{$1$};
	\draw (13.3,-4) node{$1$};
	\draw[dots] (13.9,4)--(15.3,4);
	\draw[dots] (13.9,-4)--(15.3,-4);
	\draw (16,4) node{$i$};
\draw (16,-4) node{$i$};
	\draw (0,-4) node{\color{blue}$\ell$\color{black}};
	\braidbox{0.6}{2.9}{-4.7}{-3.4}{};
	\draw (1.8,-4) node{$(\ell-1)^2$};
	\draw[dots] (3.5,-4)--(5,-4);
	\braidbox{5.3}{7.6}{-4.7}{-3.4}{};
	\draw (6.5,-4) node{$(i+1)^2$};
	\draw (8.3,-4) node{$i$};
	\draw (16.8,4) node{\color{blue}$\ell$\color{black}};
\draw (16.8,-4) node{\color{blue}$\ell$\color{black}};
	\braidbox{17.4}{19.7}{3.3}{4.6}{};
	\draw (18.6,4) node{$(\ell-1)^2$};
	\braidbox{17.4}{19.7}{-3.3}{-4.6}{};
	\draw (18.6,-4) node{$(\ell-1)^2$};
	\draw[dots] (20.3,4)--(21.8,4);
	\draw[dots] (20.3,-4)--(21.8,-4);
	\braidbox{22.2}{24.5}{3.3}{4.6}{};
	\draw (23.4,4) node{$(i+1)^2$};
	\braidbox{22.2}{24.5}{-3.3}{-4.6}{};
	\draw (23.4,-4) node{$(i+1)^2$};
	\draw (25.2,4) node{$i$};
	\draw (25.2,-4) node{$i$};
	\draw[dots] (25.9,4)--(27.6,4);
	\draw[dots] (25.9,-4)--(27.6,-4);
	\draw (27.9,4) node{$1$};
	\draw (27.9,-4) node{$1$};
	\redbraidbox{28.5}{29.4}{3.3}{4.6}{};
	\draw (28.9,4) node{$\color{red}0\,\,0\color{black}$};
	\redbraidbox{28.5}{29.4}{-3.3}{-4.6}{};
	\draw (28.9,-4) node{$\color{red}0\,\,0\color{black}$};
	\draw (30.1,4) node{$1$};
	\draw (30.1,-4) node{$1$};
	\draw[dots] (30.7,4)--(32.1,4);
	\draw[dots] (30.7,-4)--(32.1,-4);
	\draw (32.8,4) node{$i$};
	\draw (32.8,-4) node{$i$};
	\draw[blue](0,-3.2)--(16.8,3.2);
	\draw[blue](0,3.2)--(16.8,-3.2);
	\draw(1.3,-3.2)--(17.9,3.2);
	\draw(1.3,3.2)--(17.9,-3.2);
	\draw(2.7,-3.2)--(19.3,3.2);
	\draw(2.7,3.2)--(19.3,-3.2);
	\draw(5.8,-3.2)--(22.4,3.2);
	\draw(5.8,3.2)--(22.4,-3.2);
	\draw(7.3,-3.2)--(23.9,3.2);
	\draw(7.3,3.2)--(23.9,-3.2);
	\draw(8.3,-3.2)--(25.2,3.2);
	\draw(25.2,-3.2)--(8.3,3.2);
	\draw(11,-3.2)--(20,0)--(11.4,3)--(11.3,3.2);
	\draw(27.9,-3.2)--(21,0)--(27.9,3.2);
	\draw[red](11.8,-3.2)--(21,0)--(12.5,2.85)--(12.5,3.2);
	\draw[red](12.5,-3.2)--(23,-0.1)--(29.2,3.2);
	\draw[red](28.6,-3.2)--(23,-0.1)--(11.8,3.2);
	\draw[red](29.3,-3.2)--(23,0.4)--(28.6,3.2);
	\draw(13.3,-3.2)--(25.5,0)--(30.1,3.2);
	\draw(13.3,3.2)--(25.5,0)--(30.1,-3.2);
	\draw(32.8,-3.2)--(29,0)--(16,3.2);
	\draw(32.8,3.2)--(29,0)--(16,-3.2);
	\reddot(12.9,-3.1);
	\end{braid}.
\end{align*}
Using similarly the braid relations for 
$$\begin{braid}\tikzset{baseline=.4em}
	\draw (0,0) node{$2$};
        \draw(1.1,0) node{$1$};
        \draw (2.2,0) node{$2$};
		\draw(0,0.5)--(2.2,1.5);
\draw(2.2,0.5)--(0,1.5);
\draw(1.1,0.5)--(2,1)--(1.1,1.5);
	\end{braid},\dots,
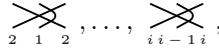
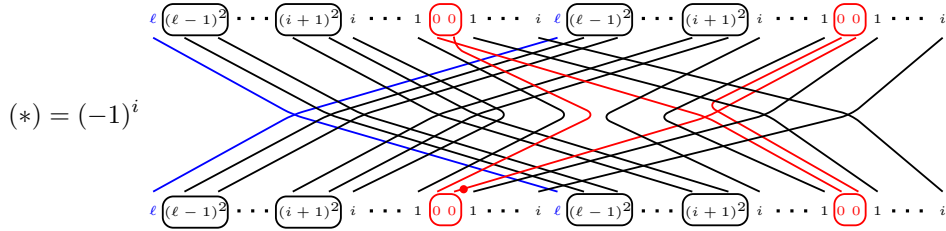
	\begin{braid}\tikzset{baseline=.4em}
	\draw (0,0) node{$i$};
        \draw(1.1,-0.05) node{$i-1$};
        \draw (2.2,0) node{$i$};
		\draw(0,0.5)--(2.2,1.5);
\draw(2.2,0.5)--(0,1.5);
\draw(1.1,0.5)--(2,1)--(1.1,1.5);
	\end{braid},
	$$
we get 
\begin{align*}
(*)&=(-1)^i
\begin{braid}\tikzset{baseline=-.3em}
	\draw (0,4) node{\color{blue}$\ell$\color{black}};
	\braidbox{0.6}{2.9}{3.3}{4.6}{};
	\draw (1.8,4) node{$(\ell-1)^2$};
	\draw[dots] (3.5,4)--(5,4);
	\braidbox{5.3}{7.6}{3.3}{4.6}{};
	\draw (6.5,4) node{$(i+1)^2$};
	\draw (8.3,4) node{$i$};
	\draw[dots] (9,4)--(10.3,4);
	\draw[dots] (9,-4)--(10.3,-4);
	\draw (11,4) node{$1$};
	\draw (11,-4) node{$1$};
	\redbraidbox{11.7}{12.6}{3.3}{4.6}{};
	\draw (12.1,4) node{$\color{red}0\,\,0\color{black}$};
	\redbraidbox{11.7}{12.6}{-3.3}{-4.6}{};
	\draw (12.1,-4) node{$\color{red}0\,\,0\color{black}$};
	\draw (13.3,4) node{$1$};
	\draw (13.3,-4) node{$1$};
	\draw[dots] (13.9,4)--(15.3,4);
	\draw[dots] (13.9,-4)--(15.3,-4);
	\draw (16,4) node{$i$};
\draw (16,-4) node{$i$};
	\draw (0,-4) node{\color{blue}$\ell$\color{black}};
	\braidbox{0.6}{2.9}{-4.7}{-3.4}{};
	\draw (1.8,-4) node{$(\ell-1)^2$};
	\draw[dots] (3.5,-4)--(5,-4);
	\braidbox{5.3}{7.6}{-4.7}{-3.4}{};
	\draw (6.5,-4) node{$(i+1)^2$};
	\draw (8.3,-4) node{$i$};
	\draw (16.8,4) node{\color{blue}$\ell$\color{black}};
\draw (16.8,-4) node{\color{blue}$\ell$\color{black}};
	\braidbox{17.4}{19.7}{3.3}{4.6}{};
	\draw (18.6,4) node{$(\ell-1)^2$};
	\braidbox{17.4}{19.7}{-3.3}{-4.6}{};
	\draw (18.6,-4) node{$(\ell-1)^2$};
	\draw[dots] (20.3,4)--(21.8,4);
	\draw[dots] (20.3,-4)--(21.8,-4);
	\braidbox{22.2}{24.5}{3.3}{4.6}{};
	\draw (23.4,4) node{$(i+1)^2$};
	\braidbox{22.2}{24.5}{-3.3}{-4.6}{};
	\draw (23.4,-4) node{$(i+1)^2$};
	\draw (25.2,4) node{$i$};
	\draw (25.2,-4) node{$i$};
	\draw[dots] (25.9,4)--(27.6,4);
	\draw[dots] (25.9,-4)--(27.6,-4);
	\draw (27.9,4) node{$1$};
	\draw (27.9,-4) node{$1$};
	\redbraidbox{28.5}{29.4}{3.3}{4.6}{};
	\draw (28.9,4) node{$\color{red}0\,\,0\color{black}$};
	\redbraidbox{28.5}{29.4}{-3.3}{-4.6}{};
	\draw (28.9,-4) node{$\color{red}0\,\,0\color{black}$};
	\draw (30.1,4) node{$1$};
	\draw (30.1,-4) node{$1$};
	\draw[dots] (30.7,4)--(32.1,4);
	\draw[dots] (30.7,-4)--(32.1,-4);
	\draw (32.8,4) node{$i$};
	\draw (32.8,-4) node{$i$};
	\draw[blue](0,-3.2)--(5.8,0)--(16.8,3.2);
	\draw[blue](0,3.2)--(5.8,0)--(16.8,-3.2);
	\draw(1.3,-3.2)--(7.3,0)--(17.9,3.2);
	\draw(1.3,3.2)--(7.3,0)--(17.9,-3.2);
	\draw(2.7,-3.2)--(8.7,0)--(19.3,3.2);
	\draw(2.7,3.2)--(8.7,0)--(19.3,-3.2);
	\draw(5.8,-3.2)--(11.8,0)--(22.4,3.2);
	\draw(5.8,3.2)--(11.8,0)--(22.4,-3.2);
	\draw(7.3,-3.2)--(13.4,0)--(23.9,3.2);
	\draw(7.3,3.2)--(13.4,0)--(23.9,-3.2);
	\draw(8.3,-3.2)--(14.8,0)--(8.3,3.2);
	\draw(11,-3.2)--(17.3,0)--(11,3.2);
	\draw[red](11.8,-3.2)--(18.4,0)--(12.5,2.8)--(12.5,3.2);
	\draw[red](12.5,-3.2)--(23,-0.1)--(29.2,3.2);
	\draw[red](28.6,-3.2)--(23,-0.1)--(11.8,3.2);
	\draw[red](29.3,-3.2)--(23,0.4)--(28.6,3.2);
	\draw(13.3,-3.2)--(25.5,0)--(30.1,3.2);
	\draw(13.3,3.2)--(25.5,0)--(30.1,-3.2);
	\draw(32.8,-3.2)--(29,0)--(16,3.2);
	\draw(32.8,3.2)--(29,0)--(16,-3.2);
	\draw(25.2,-3.2)--(18.6,-0.1)--(25.2,3.2);
	\draw(27.9,-3.2)--(21,-0.1)--(27.9,3.2);
	\reddot(12.9,-3.1);
	\end{braid}.
\end{align*}
By Lemma~\ref{C100721}, this equals $(-1)^i\Theta^i_1$, which by Lemma~\ref{L110621Gen}, equals 
$$(-1)^i\ga^{i,i}(-z_1+z_2+(-1)^i(c_1+c_2))\ga^{i,i},$$ proving the proposition.
\end{proof}

\subsection{The isomorphism $B_d\cong H_d(\Zig_\ell)$ and Morita equivalence $\bar R_{d\de}\sim H_d(\Zig_\ell)$}

Recall from Corollary~\ref{CBBasis} that we have identified $\Zig_\ell^{\otimes d}$ with a subalgebra of~$B_d$.

\begin{Theorem} \label{TAffIso}
There exists an isomorphism of graded superalgebras 
\begin{align*}
f_d: H_d(\Zig_\ell)&\iso B_d,\index{f@$f_d$}
\\
\zz_r&\mapsto z_r &(1\leq r\leq d),
\\
a_1\otimes \dots\otimes a_d&\mapsto  a_1\otimes \dots\otimes a_d &(a_1,\dots,a_d\in\Zig_\ell),
\\
w&\mapsto \uptau_w &(w\in\Si_d).
\end{align*}
\end{Theorem}
\begin{proof}
The map $f_d$ is an isomorphism of (graded) superspaces due to Theorem~\ref{TAffBasis} and Corollary~\ref{CBBasis}. It is an algebra map since the images of the generators of $H_d(\Zig_\ell)$ have been checked to satisfy the defining relations of $H_d(\Zig_\ell)$, see (\ref{ETauFarApart})--
(\ref{ETauRWreath}),  Proposition~\ref{PTauBraid}, and (\ref{ETauZIdTarget}) proved in \S\S\ref{SSRed},\ref{SSI=J}.
\end{proof}

\begin{Theorem} \label{T221221} 
Suppose that $\cha \F>d$. Then the algebras $\bar R_{d\de}$ and $H_d(\Zig_\ell)$ are graded Morita superequivalent. 
\end{Theorem}
\begin{proof}
By Theorem~\ref{THeadIrr}(vi), we have $|\Irr(\bar R_{d\de})|=|\Par^\ell(d)|$. On the other hand, $H_d(\Zig_\ell)$ has the algebra $\Zig_\ell\swr \Si_d$ as its quotient, see \cite[Proposition~3.17]{KM}, and under the assumption $\cha \F>d$, we have $|\Irr(\Zig_\ell\swr \Si_d)|=|\Par^\ell(d)|$ irreducible modules. So $|\Irr(H_d(\Zig_\ell)|\geq |\Par^\ell(d)|$. 
So it remains to recall that $B_d=\ga_{1^d}\bar R_{d\de}\ga_{1^d}$ is an idempotent truncation of $\bar R_{d\de}$ and use 
Corollary~\ref{cor:idmpt_Mor}.
\end{proof}


\chapter{RoCK blocks of double covers of symmetric groups}
\label{P4}

\section{Twisted wreath superproducts}\label{sec:equ_super}

In this section, we will work over a field $\F$ containing both $\sqrt{-1}$ a primitive $4^{\nth}$ root of unity and $\sqrt{-2}$, a square root of $-2$. 
All superalgebras will be assumed to be finite dimensional.


We will now develop some theory of equivalences for (twisted) wreath superproducts, building on Morita superequivalence theory discussed in \S\ref{SMoritaSuper}. 

\subsection{Twisted group algebras of symmetric groups}\label{sec:twist_group}
Let $n \in \Z_{>0}$. The {\em twisted group superalgebra}\index{twisted group superalgebra of $\Si_n$} $\cT_n$\index{t@$\cT_n$} of $\Si_n$ is the superalgebra
given by odd generators $\ct_1,\dots,\ct_{n-1}$\index{t@$\ct_r$} subject to the relations  
\begin{equation}\label{ET_nRel}
\ct_r^2=1,\quad \ct_r\ct_s = -\ct_s\ct_r\text{ if }|r-s|>1,\quad (\ct_r\ct_{r+1})^3=1.
\end{equation}
By convention, $\cT_1$ is the totally even superalgebra $\F$.

Choosing for each $w\in \Si_n$ a reduced decomposition $w=s_{r_1}\cdots s_{r_l}$, we define $\ct_w:=\ct_{r_1}\cdots \ct_{r_l}$\index{t@$\ct_w$} (which depends up to a sign on the choice of a reduced decomposition). Then we have a basis $\{\ct_w\mid w\in \Si_n\}$ of $\cT_n$.

Recall from \S\ref{SSBasicRep} the wreath superproduct $A\swr {\Si}_n$, 
for a superalgebra $A$. We also need the {\em twisted wreath superproduct}\index{twisted wreath superproduct}  $A\swr \cT_n$ defined as  the free product $A^{\otimes n}\star \cT_n$ of superalgebras subject to the relations:
\begin{equation}\label{ETwWreathRel}
\ct_r\,(a_1\otimes\dots\otimes a_n) = (-1)^{\sum_{u\neq r,r+1}|a_u|}\left({}^{s_r}(a_1\otimes\dots\otimes a_n)\,\ct_r\right) 
\end{equation}
for all $a_1,\dots,a_n\in A$ and $1\leq r< n$. 

\begin{Proposition}
\label{prop:semi_direct2}
We have $(A\otimes \cC_1)\swr {\Si}_n\cong (A\swr \cT_n)\otimes \cC_n$.
\end{Proposition}
\begin{proof}
Recall the notation (\ref{EInsertion}). It is easy to check directly using the defining relations that the following assignments define mutually inverse isomorphisms of superalgebras $(A\otimes \cC_1)\swr {\Si}_n$ and $(A\swr \cT_n)\otimes \cC_n$:
\begin{align*}
&(a \otimes x)_r \mapsto(-1)^{r|a|}a_r \otimes x_r,\ 
s_r\mapsto\frac{1}{\sqrt{-2}}\ct_r\otimes(\cc_r-\cc_{r+1})
\\
&a_r \otimes x_r \mapsto(-1)^{r|a|}(a \otimes x)_r,\ 
\ct_r \otimes 1 \mapsto-\frac{1}{\sqrt{-2}}s_r(\cc_r-\cc_{r+1})
\end{align*}
for all $a\in A$, $x\in\cC_1$ and $1\leq r< n$.
\end{proof}

\subsection{Extending Morita superequivalence}
Recall the definition of a $G$-graded crossed superproduct from $\S$\ref{sec:more_Morita}. 
Note that $A\swr \Si_n$ is an $\Si_n$-graded crossed superproduct with the $w$-graded component $(A\swr \Si_n)_w=A^{\otimes n}w$, for each $w\in \Si_n$. 
It follows from Proposition~\ref{prop:semi_direct2} that 
 $A\swr \cT_n=\bigoplus_{w\in \Si_n} A^{\otimes n}\ct_w$. In particular, $A\swr \cT_n$ is also an $\Si_n$-graded crossed superproduct with $(A\swr \cT_n)_w = A^{\otimes n}\ct_w$, for each $w\in \Si_n$. 

Recall symmetric superalgebras from $\S$\ref{SSBasicRep} and symmetric graded crossed superproduct from $\S$\ref{sec:more_Morita}. The following lemma is a standard check.

\begin{Lemma}\label{lem:sym_Sn_grdd}
Let $A$ be a symmetric superalgera with symmetrizing form $\mu$. Then $A\swr \Si_n$ and $A\swr \cT_n$ are symmetric $\Si_n$-graded crossed superproduct with symmetrizing forms $\mu_1$ and $\mu_2$, respectively, such that $\mu_1(A^{\otimes n}w) = \mu_2(A^{\otimes n}t_w)=0$ for all $w\neq 1$ and 
\begin{align*}
\mu_1(a_1\otimes\dots\otimes a_n) = \mu_2(a_1\otimes\dots\otimes a_n)= \mu(a_1)\cdots\mu(a_n)\qquad (a_1,\dots,a_n\in A).
\end{align*}
\end{Lemma}

\begin{Proposition}
\label{prop:semi_direct1}
Let $A$ and $B$ be symmetric superalgebras such that $A\sim_{\sM}B$. Then
\begin{enumerate}
\item $A\swr {\Si}_n \sim_{\sM} B\swr {\Si}_n$.
\item $A\swr \cT_n \sim_{\sM} B\swr \cT_n$.
\end{enumerate}
\end{Proposition}
\begin{proof}
Let $A$ and $B$ be Morita superequivalent via an $(A\otimes B^\op)$-supermodule $M$. By Lemma~\ref{lem:MSE} and Remark~\ref{rem:(bi)supmod}, we have that $A^{\otimes n}$ and $B^{\otimes n}$ are Morita superequivalent via the $A^{\otimes n}\otimes (B^\op)^{\otimes n}$-supermodule $M^{\boxtimes n}$. With Lemma \ref{lem:sym_Sn_grdd} in mind, both (i) and (ii) follow from Proposition~\ref{prop:ext_marcus} if we can extend the $A^{\otimes n}\otimes (B^\op)^{\otimes n}$-supermodule structure on $M^{\boxtimes n}$ to $(A\swr \Si_n,B\swr \Si_n)_{\Si_n}$ and $(A\swr \cT_n,B\swr \cT_n)_{\Si_n}$. 

Such extensions are uniquely determined by the formulas 
\begin{align}
(s_r\otimes s_r^{-1}).(m_1\otimes\dots\otimes m_n) &:= {}^{s_r}\,(m_1\otimes\dots\otimes m_n),
\\
(\ct_r\otimes \ct_r^{-1}).(m_1\otimes\dots\otimes m_n)
&:=(-1)^{|m_r|+|m_{r+1}|}\left( {}^{s_r}\,(m_1\otimes \dots \otimes m_n)\right),
\label{E261021}
\end{align}
for all $m_1,\dots,m_n\in M$ and 
 $1\leq r\leq n-1$. 
 
Let us explain the more difficult second case. For $1\leq r<n$, set $T_r:=\ct_r\otimes \ct_r^{-1}=\ct_r\otimes \ct_r$. Note that the subalgebra
 \begin{align*}
(A\swr \cT_n,B\swr \cT_n)_{\Si_n}=\sum_{w\in {\Si}_n}A^{\otimes n}\ct_w \otimes (B^{\otimes n})^{\op}\ct_w^{-1}\subseteq (A\swr \cT_n) \otimes (B\swr \cT_n)^{\op}
\end{align*}
is generated by the subalgebra $A^{\otimes n}\otimes (B^\op)^{\otimes n}$ and the elements $T_1,\dots,T_{n-1}$ subject only to the relations 
\begin{align}
T_r^2=1,\quad T_rT_s = T_sT_r\text{ if }|r-s|>1,\quad (T_rT_{r+1})^3=1,
\label{E261021_1}
\\
T_r(\ua\otimes \ub)
=(-1)^{|a_r|+|a_{r+1}|+|b_r|+|b_{r+1}|}\big({}^{s_r}\ua\otimes
{}^{s_r}\ub)T_r,
\label{E261021_2}
\end{align}
where we have used the notation $\ua:= a_1\otimes \dots\otimes a_n$ and  $\ub:=b_1\otimes\dots\otimes b_n$. 
It is quickly checked that the action (\ref{E261021}) satisfies the relations (\ref{E261021_1}). To check that it also satisfies (\ref{E261021_2}), we denote $\um:=m_1\otimes\dots\otimes m_n$ and show that
\begin{align}\label{algn:twist_rel}
T_r(\ua\otimes \ub)\,\um
=(-1)^{|a_r|+|a_{r+1}|+|b_r|+|b_{r+1}|}\big({}^{s_r}\ua\otimes
{}^{s_r}\ub\big)T_r\,\um,
\end{align}
for all $\ua \in A^{\otimes n}$ and $\ub \in B^{\otimes n}$. It is easily seen that the check boils down to the case $n=2$, where we have to check
\begin{align*}
&T_1((a\otimes a')\otimes (b\otimes b'))(m\otimes m')
\\
=\,&(-1)^{|a|+|a'|+|a||a'|+|b|+|b'|+|b||b'|}((a'\otimes a)\otimes (b'\otimes b))T_1(m\otimes m').
\end{align*}
The left hand side equals
\begin{align*}
\eps_1T_1(amb\otimes a'm'b')
=\eps_1\eps_2(a'm'b'\otimes amb),
\end{align*}
where 
\begin{align*}
\eps_1&:=(-1)^{|b'||m|+|b'||m'|+|b||m|+|a'||m|+|a'||b|},\\
\eps_2&:=(-1)^{|amb|+|a'm'b'|+|amb||a'm'b'|}.
\end{align*}
The right hand side equals 
\begin{align*}
\de_1((a'\otimes a)\otimes (b'\otimes b))T_1(m\otimes m')
&=\de_1\de_2((a'\otimes a)\otimes (b'\otimes b))(m'\otimes m)
\\
&=\de_1\de_2\de_3(a'm'b'\otimes amb),
\end{align*}
where 
\begin{align*}
\de_1&:=(-1)^{|a|+|a'|+|a||a'|+|b|+|b'|+|b||b'|},\\
\de_2&:=(-1)^{|m|+|m'|+|m||m'|},\\
\de_3&:=(-1)^{|b||m'|+|b||m|+|b'||m'|+|a||m'|+|a||b'|}.
\end{align*}
An easy computation with signs now completes the proof. 
\end{proof}

\subsection{Extending complexes}

Let $A$ be a superalgebra with a superunit $u$ (see \S\ref{sec:more_Morita}). If $C_2$ is the cyclic group with generator $g$, we can give $A$ the structure of a $C_2$-graded crossed superproduct  with graded components 
$
A_{g^{\eps}}=
A_{\0}u^{\eps}=A_\eps
$
for $\eps\in\Z/2$. 
For any $n\in\Z_{>0}$, we consider the wreath product 
$$C_2\wr {\Si}_n=\{g_1^{\eps_1}\cdots g_n^{\eps_n}w\mid \eps_1,\dots,\eps_n\in\Z/2,\, w\in \Si_n\}$$
where $g_r$ is the generator of the $r^{\nth}$ $C_2$ factor in the base group $C_2^{\times n}$. Then we can give $A\swr \cT_n$ the structure of a $C_2\wr {\Si}_n$-graded crossed superproduct with 
graded components 
$$
(A\swr \cT_n)_{g_1^{\eps_1}\cdots g_n^{\eps_n}w}=
A_{\0}^{\otimes n}u_1^{\eps_1}\dots u_n^{\eps_n} \ct_w
$$
where each $u_r$ is defined via (\ref{EInsertion}). 

Let $A$ now be a superalgebra with a superunit $u_A$ and 
$B$ be a superalgebra with a superunit $u_B$. 
Recalling (\ref{EDiag}), we have the subsuperalgebras  
\begin{align}
&(A,B)_{C_2}=(A_\0\otimes B_\0^\op)\oplus (A_\1\otimes B_\1^\op)\subseteq A\otimes B^\op,
\label{EC_2}\index{$(A,B)_{C_2}$}
\\
&(A\swr \cT_n,B\swr \cT_n)_{C_2\wr\Si_n}\subseteq (A\swr \cT_n)\otimes (B\swr \cT_n)^\op.\index{$(A\swr \cT_n,B\swr \cT_n)_{C_2\wr\Si_n}$}
\end{align}
Setting $U_i:=(u_A)_i \otimes (u_B^{-1})_i$, for $1\leq i\leq n$, we have that $A_{\0}^{\otimes n}\otimes (B_{\0}^{\otimes n})^{\op}$ together with the $U_i$'s generate a subsuperalgebra of $(A\swr \cT_n,B\swr \cT_n)_{C_2\wr\Si_n}$ isomorphic to $((A,B)_{C_2})^{\otimes n}$. We identify the two algebras. 
Then $(A\swr \cT_n,B\swr \cT_n)_{C_2\wr\Si_n}$ is generated by $((A,B)_{C_2})^{\otimes n}$ and $T_i:=\ct_i \otimes \ct_i^{-1}$, for $1\leq i\leq n-1$ subject to the following relations:
\begin{align}\label{algn:complex_rel}
\begin{split}
&T_i (\ua \otimes \ub) = ({}^{s_i}\ua \otimes {}^{s_i}\ub) T_i, \text{ for all }1\leq i\leq n-1,\\
&T_i^2 =1, \text{ for all }1\leq i\leq n-1,\\
&T_iT_j = T_jT_i, \text{ for all }1\leq i,j\leq n-1\text{ with } |i-j|>1,\\
&T_iT_{i+1}T_i = T_{i+1}T_iT_{i+1}, \text{ for all }1\leq i\leq n-2,\\
&T_i U_j = U_j T_i, \text{ for all }1\leq i \leq n-1, 1\leq j\leq n\text{ and }j\neq i,i+1,\\
&T_i U_i = U_{i+1} T_i, \text{ for all }1\leq i \leq n-1,\\
&T_i U_{i+1} = U_i T_i, \text{ for all }1\leq i \leq n-1,
\end{split}
\end{align}
where $\ua = a_1\otimes \dots \otimes a_n \in A_{\0}^{\otimes n}$ and $\ub = b_1 \otimes \dots \otimes b_n \in (B_{\0}^{\otimes n})^{\op}$.

\begin{Proposition}\label{prop:ext_der_wreath}
Let $n\in\Z_{> 0}$ and $A$ and $B$ be superalgebras with superunit. If $X$ is a complex of totally even $A_{\0}$-$B_{\0}$-bisupermodules that extends to a complex of $(A,B)_{C_2}$-supermodules, then $X^{\otimes n}$ extends to a complex of $(A\swr \cT_n,B\swr \cT_n)_{C_2\wr\Si_n}$-supermodules.
\end{Proposition}

In the above proposition $A_{\0}$, $B_{\0}$ and $(A\swr \cT_n,B\swr \cT_n)_{C_2\wr\Si_n}$ are all totally even superalgebras and so ``supermodule'' and ``bisupermodule'' can really be replaced with ``module'' and ``bimodule''. However, we use the ``super'' notation to better fit with the language of $\S$\ref{sec:more_Morita}.

\begin{proof}
We first note that $X$ is already a complex of $(A,B)_{C_2}$-supermodules and so $X^{\otimes n}$ automatically extends to a complex of $((A,B)_{C_2})^{\otimes n}$-supermodules.

Throughout this proof $x_{i_1}\otimes \dots \otimes x_{i_n}$ will be an element in degree $i_1+\dots+i_n$ of the complex $X^{\otimes n}$, where each $x_{i_j}$ has degree $i_j$. By~\cite[Lemma 4.1(b)]{Mar}, for each $w\in {\Si}_n$ and $i_1,\dots,i_n\in \ZZ$, there exists $\epsilon_w(i_1,\dots,i_n)\in \{\pm1\}$ such that
\begin{align}\label{def_act}
w*(x_{i_1}\otimes \dots \otimes x_{i_n}) = \epsilon_w(i_1,\dots,i_n) (x_{i_{w^{-1}(1)}}\otimes \dots \otimes x_{i_{w^{-1}(n)}})
\end{align}
turns $X^{\otimes n}$ into a complex of $\F {\Si}_n$-modules. We extend $X^{\otimes n}$ to a complex of $(A\swr \cT_n,B\swr \cT_n)_{C_2\wr\Si_n}$-supermodules via
\begin{align*}
T_r.(x_{i_1}\otimes \dots \otimes x_{i_n}) &= s_r*(x_{i_1}\otimes \dots \otimes x_{i_n}),
\end{align*}
for $1\leq r\leq n-1$. First note that, since `$*$' turns $X^{\otimes n}$ into a complex of $\F {\Si}_n$-modules, the action of the $T_r$'s commutes with the differential. To complete the proof one needs to show that all the relations in (\ref{algn:complex_rel}) hold when considering the action of $(A\swr \cT_n,B\swr \cT_n)_{C_2\wr\Si_n}$ on $X^{\otimes n}$. However, they all follow directly from (\ref{def_act}) or from the fact that $*$ turns turns $X^{\otimes n}$ into a complex of $\F {\Si}_n$-modules. Note that these relations are significantly easier to check than those in the proof of Proposition~\ref{prop:semi_direct1}, since $X$ is a complex of totally even modules.
\end{proof}

\section{Blocks of double covers}\label{sec:spin_blocks}

Let $p$ be an odd prime. We can write $p= 2\ell + 1$ for some $\ell\in\Z_{>0}$. Let $\F$ be an algebraically closed field of characteristic $p$. 

\subsection{Double covers of symmetric and alternating groups}
\label{SDoubleCovers}

Let $n \in \Z_{>0}$. We have the double cover
\index{double cover of the symmetric group}
 $\tilde{\Si}_n^+$\index{s@$\tilde{\Si}_n^\pm$} of the symmetric group ${\Si}_n$ given by generators $t_1,\dots,t_{n-1},z$ subject to the relations 
\begin{align*}
z^2=1,\quad t_r z = z t_r,\quad t_r^2=1,
\quad
t_r t_s = z t_s t_r \text{ if } |r-s| > 1,\quad (t_r t_{r+1})^3=1.
\end{align*}
Another double cover $\tilde{\Si}_n^-$ of ${\Si}_n$ is given by generators $t_1,\dots,t_{n-1},z$ subject to the relations 
\begin{align*}
z^2=1,\quad t_r z = z t_r,\quad t_r^2=z,
\quad
t_r t_s = z t_s t_r \text{ if } |r-s| > 1,\quad (t_r t_{r+1})^3=1.
\end{align*}
We have the natural surjections
$$
\pi_n: \tilde{\Si}^\pm_n \to {\Si}_n,\ z\mapsto 1,\ t_r \mapsto s_r.
$$
We denote by $\tilde{\Ai}_n^{\pm}$ the inverse image $\pi_n^{-1}({\Ai}_n)$ of the alternating group ${\Ai}_n\leq \Si_n$.\index{double cover of the alternating group}
\index{a@$\tilde{\Ai}_n^\pm$}

We will often use $\tilde{\Si}_n$\index{s@$\tilde{\Si}_n$} (resp. $\tilde{\Ai}_n$)\index{a@$\tilde{\Ai}_n$} to simultaneously denote both $\tilde{\Si}_n^+$ and $\tilde{\Si}_n^-$ (resp. $\tilde{\Ai}_n^+$ and $\tilde{\Ai}_n^-$). We denote by 
$$|\cdot|:\tilde{\Si}_n \to \ZZ/2$$ the unique group homomorphism with kernel $\tilde{\Ai}_n$.

If $X\subseteq \{1,\dots,n\}$, then we have the subgroup
$$\Si_X:=\{g\in\Si_n\mid g(k)=k\ \text{for all}\ k\in\{1, \dots, n\}\backslash X\} \leq \Si_n.\index{s@$\Si_X$}
$$
We then denote 
$$\tilde{\Si}_X:=\pi_n^{-1}(\Si_X)\leq \tilde \Si_n,\index{s@$\tilde{\Si}_X$}\quad \Ai_X:=\Si_X\cap\Ai_n\index{a@$\Ai_X$} \quad \text{and}\quad \tilde{\Ai}_X:=\pi_n^{-1}(\Ai_X).\index{a@$\tilde{\Ai}_X$}
$$
If $X=\{1,\dots,m\}$, $Y=\{m+1,\dots,n\}$, $x\in\tilde\Si_X$, and $y\in\tilde\Si_Y$, then we have
\begin{equation}\label{EZComm}
xy=z^{|x||y|}yx.
\end{equation}
This follows from the fact that $x=z^\eps t_{r_1}\cdots t_{r_k}$, $y=z^\de t_{s_1}\cdots t_{s_l}$ for $\eps,\de\in\{0,1\}$, $|r_t-s_u|>1$ for all $t,u$, and $k \equiv |x|$, $l\equiv |y|\pmod 2$.

We have the idempotent
\begin{equation}\label{EEZ-}
e_{z}:=(1-z)/2\in\F\tilde{\Si}_n.\index{e@$e_{z}$}
\end{equation}
We consider $\F \tilde{\Si}_n e_{z}$ as a superalgebra by setting all the $t_r e_z$'s to be odd. In particular, $(\F \tilde{\Si}_n e_{z})_{\0} = \F \tilde{\Ai}_n e_{z}$. Adopting the notation of (\ref{EZComm}), it follows that
\begin{equation}\label{EZCommAlg}
\F\tilde{\Si}_{X,Y}e_z \cong \F\tilde{\Si}_X e_z \otimes \F\tilde{\Si}_Y e_z
\end{equation}
as superalgebras, where $\tilde{\Si}_{X,Y}$ is the subgroup of $\tilde{\Si}_n$ generated by $\tilde{\Si}_X$ and $\tilde{\Si}_Y$. By conjugating by an appropriate element of $\tilde{\Si}_n$ (\ref{EZCommAlg}) holds for any $X,Y\subseteq \{1,\dots,n\}$ with $X\cap Y = \varnothing$.

There are superalgebra isomorphisms 
\begin{align*}
\cT_n\cong \F \tilde{\Si}_n^+ e_{z},\ \ct_r\mapsto t_r e_{z}\qquad\text{and}\qquad
\cT_n\cong \F \tilde{\Si}_n^- e_{z},\ \ct_r\mapsto(-1)^r \sqrt{-1}t_r e_{z}.
\end{align*}
We {\em identify}\, $\F \tilde{\Si}_n e_{z}$ 
with  $\cT_n$ via these isomorphisms. In this way $\F \tilde{\Ai}_n e_{z}$ gets identified with $(\cT_n)_\0$. 

We will use the adjective ``spin'' to refer to $\F \tilde{\Si}_n e_{z}= \cT_n$ as opposed to the whole of $\F \tilde{\Si}_n$. For example, the spin characters\index{spin character} of $\tilde{\Si}_n$ refers to the ordinary characters $\chi$ of $\tilde{\Si}_n$ such that $\chi(z) = -\chi(1)$ and 
the spin blocks\index{spin block} of $\F\tilde\Si_n$ (or simply the spin $p$-blocks of $\tilde{\Si}_n$) refers to the blocks of $\F \tilde{\Si}_n e_{z}$.

Following Schur \cite{Schur}, to every partition $\la\in\Par_0(n)$, one associates canonically a spin character $\chi^\la$ of $\tilde\Si_n$ such that $\chi^\la$ is irreducible if $n-h(\la)$ is even and $\chi^\la=\chi^\la_++\chi^\la_-$\index{x@$\chi^\la_{(\pm)}$}\index{h@$\chi^\la_{(\pm)}$}
\index{$\chi^\la_{(\pm)}$} for irreducible characters $\chi^\la_\pm$ if $n-h(\la)$ is odd. Moreover 
$$\{\chi^\la\mid\la\in\Par_0(\la)\ \text{and $n-h(\la)$ is even}\}\cup
\{\chi^\la_\pm\mid\la\in\Par_0(\la)\ \text{and $n-h(\la)$ is odd}\}
$$
is a complete irredundant set of irreducible spin characters of $\tilde\Si_n$.

\subsection{Spin blocks}
\label{SSB}

For a finite group $H$, we say the irreducible character $\chi$ of $H$ lies in a particular block of $\F H$ if its reduction to a Brauer character lies in this block. The spin $p$-blocks of $\tilde \Si_n$ are described explicitly in terms of $\bar p$-cores as follows:

\begin{Theorem} \label{THumphreys} {\rm \cite{Humphreys}}\,\,
Let $p$ be an odd prime. With one exception, 
$\chi^\la_{(\pm)}$ and $\chi^\mu_{(\pm)}$ are in the same $p$-block if and only if $\core(\la)=\core(\mu)$.  The exception is $\la=\mu$ is a $\bar p$-core and $n-h(\la)$ is odd, in which case $\chi^\la_+$ and $\chi^\la_-$ are in different $p$-blocks.
\end{Theorem}

Recalling the notation (\ref{ECores}), we denote 
$$
\Xi(n):=\{(\rho,d)\in \Cores_p\times\Z_{\geq 0}\mid |\rho|+dp=n\}. 
\index{$\Xi(n)$}\index{x@$\Xi(n)$}
$$ 
Let $(\rho,d)\in\Xi(n)$. According to Theorem~\ref{THumphreys}, unless $d=0$ and $n-h(\rho)$ is odd, there is a unique block of $\cT_{n}=\F\tilde\Si_{n}e_z$ corresponding to the characters $\chi^\la$ with $\la \in \Par_0(n)$ and $\core(\la)=\rho$. We denote the corresponding block idempotent by $e_{\rho,d}$.\index{e@$e_{\rho,d}$} In the exceptional case, where $d=0$ and $n-h(\rho)$ is odd, we have two blocks of $\cT_n=\F\tilde\Si_{n}e_z$ corresponding to the characters $\chi^\rho_{\pm}$. We denote the corresponding block  idempotents by $e_{\rho,0}^\pm$\index{e@$e_{\rho,0}^\pm$} and set $e_{\rho,0}:=e_{\rho,0}^++e_{\rho,0}^-$. Thus for any $(\rho,d)\in\Xi(n)$, we have a two-sided superideal
\begin{equation}\label{ERhoDBlocks}
\Blo^{\rho,d}:=\cT_ne_{\rho,d}=\F\tilde\Si_{n}e_{\rho,d}\subseteq \cT_n\subseteq \F\tilde\Si_{n}.\index{b@$\Blo^{\rho,d}$}
\end{equation}
We have a decomposition into two-sided ideals 
$$
\cT_n=\bigoplus_{(\rho,d)\in\Xi(n)}\Blo^{\rho,d}.
$$

Note that in all cases $\Blo^{\rho,d}$ is a {\em superblock}\index{superblock} of $\cT_n$, i.e. an indecomposable two-sided superideal of $\cT_{n}$. Indeed, recalling (\ref{ESi}), this follows from the fact that $\sigma_{\F \tilde{\Si}_n}(e_{\rho,d})=e_{\rho,d}$ (and $\sigma_{\F \tilde{\Si}_n}(e_{\rho,0}^\pm)=e_{\rho,0}^\mp$ in the exceptional case). This in turn is clear from the known fact that $\chi^\la_\pm$ is obtained from $\chi^\la_\mp$ by tensoring with a sign character while $\chi^\la$ is invariant under tensoring with a sign character. 
In particular, we always have that $e_{\rho,d}\in \F\tilde{\Ai}_n=(\F\tilde{\Si}_n)_\0$, and 
$$
\Blo^{\rho,d}_\0=\F\tilde{\Ai}_ne_{\rho,d}\index{b@$\Blo^{\rho,d}_\0$}
$$
is a two-sided ideal of $\F\tilde{\Ai}_n$ (which is a block of $\F\tilde{\Ai}_n$ with one exception, see Theorem~\ref{thm:An_blocks} below). 

Let $(\rho,d)\in\Xi(n)$. Throughout the subsection, we denote
$$U=\{1,\dots,|\rho|\},\quad V=\{|\rho|+1,\dots,n\}.$$
Moreover, for $k=1,\dots,d$, we denote 
$$V_k=\{|\rho|+(k-1)p+1,\dots, |\rho|+kp\}.$$
Thus we have 
\begin{equation}\label{E241221}
\{1,\dots,n\}=U\sqcup V=U\sqcup V_1\sqcup\dots\sqcup V_d.
\end{equation}

The defect groups of spin blocks of $\F \tilde{\Si}_n$ are described as follows:

\begin{Theorem}\label{thm:Sn_blocks}
{\rm \cite[Theorems A, B]{Ca}} \,
Let $(\rho,d)\in\Xi(n)$. If $d=0$ then the block $\F \tilde\Si_n e_{\rho,0}^{(\pm)}$ has trivial defect group. If $d>0$, a Sylow $p$-subgroup $D$ of $\tilde\Si_V$ is a defect group of the block $\Blo^{\rho,d}=\F \tilde{\Si}_n e_{\rho,d}$.
\end{Theorem}

The spin blocks of $\F \tilde{\Ai}_n$ are described as follows:

\begin{Theorem}\label{thm:An_blocks}
{\rm \cite[Proposition 3.16]{Ke}} \,
Let $(\rho,d)\in\Xi(n)$. If $n>1$, then Then $\Blo^{\rho,d}_\0=\F \tilde{\Ai}_n e_{\rho,d}$ is a single block of $\F \tilde{\Ai}_n$, unless $d=0$ and $n-h(\rho)$ is even. In this latter case $\Blo^{\rho,0}_\0$ is a direct sum of two blocks of $\F \tilde{\Ai}_n$. If $n=1$, then $\Blo^{(1),0}_\0=\Blo^{(1),0}=\F e_z$ is a single block of $\F \tilde{\Ai}_1$.

In all cases the defect group of $\F \tilde{\Ai}_n e_{\rho,d}$ is the same as that of $\F \tilde{\Si}_n e_{\rho,d}$.
\end{Theorem}

From now on we often consider only $(\rho,d)\in\Xi(n)$ with  $d>0$. In particular, this ensures that blocks and superblocks coincide. The case $d=0$ corresponds to the trivial defect blocks and is completely elementary.

The Brauer correspondents of spin 
blocks of $\F \tilde{\Si}_n$ and $\F \tilde{\Ai}_n$ are described as follows. Let us always embed $e_{\rho,0}$ into $\F N_{\tilde{\Si}_n}(D)$ via $e_{\rho,0}\in\tilde \F\Si_U\into \F N_{\tilde{\Si}_n}(D)$. (Here $e_{\varnothing,0}$ is interpreted as $e_z$.) Now, denote
\begin{equation}\label{EbloDef}
\blo^{\rho,d}:=\F N_{\tilde{\Si}_n}(D) e_{\rho,0}\quad\text{and}\quad
\blo^{\rho,d}_\0:=\F N_{\tilde{\Ai}_n}(D) e_{\rho,0}.
\index{b@$\blo^{\rho,d}$}\index{b@$\blo^{\rho,d}_\0$}
\end{equation}
Note that $\blo^{\rho,d}_\0$ is indeed the even part of the superalgebra $\blo^{\rho,d}$.  

\begin{Theorem}\label{thm:brauer}
Let $(\rho,d)\in \Xi(n)$ with $d>0$, and $D$ be a Sylow $p$-subgroup of $\tilde{\Si}_V$. Then the Brauer correspondent of\, $\Blo^{\rho,d}$ in $N_{\tilde{\Si}_n}(D)$ is\, $\blo^{\rho,d}$, and the Brauer correspondent of\, $\Blo^{\rho,d}_\0$ in $N_{\tilde{\Ai}_n}(D)$ is\, $\blo^{\rho,d}_\0$. 
\end{Theorem}

\begin{proof}
The first statement is \cite[Theorem A, Corollary 26]{Ca}. For the blocks of $\tilde{\Ai}_n$, note that $e_{\rho,d},e_{\rho,0}\in \F \tilde{\Ai}_n$. Therefore, $e_{\rho,0}$ is the image of $e_{\rho,d}$ under the Brauer homomorphism, whether it is with respect to $\tilde{\Si}_n$ or $\tilde{\Ai}_n$.
\end{proof}

Let $P$ be a Sylow $p$-subgroup of $\tilde{\Si}_p$. Note that $P\leq \tilde{\Ai}_p$. 
A special role will be played by the blocks 
\begin{equation}\label{EBlo}
\Blo^{\varnothing,1}=\F\tilde{\Si}_p e_{\varnothing,1}\quad \text{and} \quad 
\Blo^{\varnothing,1}_\0=\F\tilde{\Ai}_p e_{\varnothing,1}
\end{equation}
of $\F\tilde{\Si}_p$ and $\F\tilde{\Ai}_p$ with Brauer correspondents
\begin{equation}\label{Eblo}
\blo^{\varnothing,1}=\F N_{\tilde{\Si}_p}(P)e_{z}
\quad \text{and} \quad 
\blo^{\varnothing,1}_\0=\F N_{\tilde{\Ai}_p}(P)e_{z},
\end{equation}
respectively. 

The following proposition treats the case of abelian defect group in more detail:

\begin{Proposition}\label{prop:Brauer_corr}
Let $(\rho,d)\in \Xi(n)$ with $d>0$, and $D$ be a Sylow $p$-subgroup of\, $\tilde{\Si}_V$. 
Then $D$ is abelian if and only if $d<p$. In this case we can choose $D=D_1\times \dots \times D_d$, where each $D_k$ is a Sylow $p$-subgroup of $\tilde{\Si}_{V_k}$, and we have an isomorphism of superalgebras 
\begin{align*}
\blo^{\rho,d}\cong \F \tilde{\Si}_U e_{\rho,0}\otimes (\blo^{\varnothing,1}\swr \cT_d)\cong \Blo^{\rho,0}\otimes (\blo^{\varnothing,1}\swr \cT_d).
\end{align*}
\end{Proposition}

\begin{proof}
We denote by $P$ a Sylow $p$-subgroup of $\tilde{\Si}_p$. 
Since $p$ is odd, for every $p$-element $w\in\Si_n$ there exists a unique $p$-element $\tilde w\in\tilde\Si_n$ such that $\pi_n(\tilde g)=g$. In particular, $\pi_n$ induces a bijection between the set of Sylow $p$-subgroups of $\tilde\Si_V \cong \tilde \Si_{dp}$ and those of $\Si_V\cong \Si_{dp}$. Moreover $\pi_n$ induces an isomorphism between corresponding Sylows.
It is well-known that Sylow $p$-subgroups of $\Si_{dp}$ are abelian if and only if $d<p$, in which case Sylow $p$-subgroups are isomorphic to $C_p^{\times d}$. The first claim of the proposition now follows. Let $d<p$ from now on. Then
$$D=\lan D_1,\dots,D_d\ran=D_1\times\dots\times D_d\cong C_p^{\times d}$$
is a Sylow $p$-subgroup of $\tilde\Si_V$.
By an abuse of notation we also write $D=\pi_n(D)$, $D_k=\pi_n(D_k)$, $P=\pi_p(P)$, etc. 

We first describe $N_{\Si_n}(D)$.
Define the group homomorphism $\iota:\Si_d \to \Si_V$ via the natural permutation action of $\Si_d$ on the $V_i$'s. More precisely,
\begin{align*}
\iota(w)\cdot (|\rho|+(k-1)p+t) = |\rho|+(w(k)-1)p+t,
\end{align*}
for all $1\leq k\leq d$, $1\leq t\leq p$ and $w\in \Si_d$. We may choose the Sylow $p$-subgroups $D_k\leq \Si_{V_k}$ so that $\iota(w)(D_k)\iota(w)^{-1} = D_{w(k)}$, for all $1\leq k\leq d$ and $w\in \Si_d$. Now, the non-trivial elements of the $D_k$'s are precisely the elements of $D$ with exactly $n-p$ fixed points when considering the action of $\Si_n$ on $\{1,\dots,n\}$. In particular, $N_{\Si_n}(D)$ must permute the $D_k$'s via conjugation. Therefore,
\begin{align*}
N_{\Si_n}(D) & = C_{\Si_n}(D) \times ((N_{\Si_{V_1}}(D_1)\times \dots \times N_{\Si_{V_d}}(D_d)) \rtimes \iota(\Si_d))\\
& = \Si_U \times ((N_{\Si_{V_1}}(D_1)\times \dots \times N_{\Si_{V_d}}(D_d)) \rtimes \iota(\Si_d)) 
\\
& \cong \Si_U \times (N_{\Si_p}(P)\wr \Si_d).
\end{align*}

It now follows from the unique lifting property of $p$-elements recorded in the first paragraph that $N_{\tilde \Si_n}(D) = \pi_n^{-1}(N_{\Si_n}(D))$ and $N_{\tilde \Si_{V_i}}(D_i) = \pi_n^{-1}(N_{\Si_{V_i}}(D_i))$, for all $1\leq i\leq d$. In particular, $\blo^{\rho,d}=\F N_{\tilde{\Si}_n}(D) e_{\rho,0}$ is an $\Si_d$-graded crossed superproduct via
\begin{align}\label{Norm_crossed_prod}
\F N_{\tilde{\Si}_n}(D) e_{\rho,0} = \bigoplus_{w\in \Si_d} \F H e_{\rho,0} T_w,
\end{align}
where
$$H:=\pi_n^{-1}(\Si_U \times N_{\Si_{V_1}}(D_1)\times \dots \times N_{\Si_{V_d}}(D_d))$$
and each $T_w$ is a lift of $\iota(w)$ to $\tilde{\Si}_V$. Since $e_{\rho,0}e_z=e_{\rho,0}$, it follows from (\ref{EZCommAlg}) that
\begin{align*}
\F H e_{\rho,0} \cong \F \tilde{\Si}_U e_{\rho,0}\otimes (\F N_{\tilde{\Si}_{V_1}}(D_1)e_{z} \otimes \dots\otimes \F N_{\tilde{\Si}_{V_d}}(D_d)e_{z}).
\end{align*}


For each $1 < k\leq d$, we identify $\F N_{\tilde{\Si}_{V_k}}(D_k)e_{z}$ with 
$\F N_{\tilde{\Si}_{V_1}}(D_1)e_{z}\cong \blo^{\varnothing,1}$ via the isomorphism
\begin{align}\label{NSp_ident}
\F N_{\tilde{\Si}_{V_1}}(D_1)e_{z} \to \F N_{\tilde{\Si}_{V_k}}(D_k)e_{z},\quad
a \mapsto (-1)^{k|a|}T_{(1,k)}aT_{(1,k)}^{-1}.
\end{align}
(Note this does not depend on our choice of $T_{(1,k)}$, since $T_{(1,k)}$ is unique up to multiplication by $z$.) 
We can therefore write
\begin{align*}
\F H e_{\rho,0} \cong \F \tilde{\Si}_U e_{\rho,0}\otimes (\blo^{\varnothing,1})^{\otimes d}.
\end{align*}
By (\ref{Norm_crossed_prod}) and (\ref{EZCommAlg}), we even have
\begin{align*}
\F N_{\tilde{\Si}_n}(D) e_{\rho,0} = \F \tilde{\Si}_U e_{\rho,0} \otimes \bigoplus_{w\in \Si_d} (\blo^{\varnothing,1})^{\otimes d}\, T_w.
\end{align*}

We claim we have the isomorphism of superalgebras
\begin{align*}
\blo^{\varnothing,1}\swr \cT_d \cong \bigoplus_{w\in \Si_d} (\blo^{\varnothing,1})^{\otimes d}\, T_w,
\end{align*}
via\, $\ct_w \mapsto \kappa_w T_w e_z$, for all $w\in \Si_d$ and some choices of $\kappa_w \in \F^\times$. We set $T_i := T_{(i,i+1)}$ and $\kappa_i:=\kappa_{(i,i+1)}$ (to be chosen), for all $1\leq i< d$. Let $w\in \Si_d$ and $w=s_{i_1}\dots s_{i_m}$ be a reduced expression for $w$. Certainly $\pi_n(T_w) = \pi_n(T_{i_1}\dots T_{i_m})$ and so $T_w = T_{i_1}\dots T_{i_m}$ or $zT_{i_1}\dots T_{i_m}$. In particular, $T_w e_z = \pm T_{i_1}\dots T_{i_m}e_z$ and we already know that $\ct_w = \pm \ct_{i_1}\dots \ct_{i_m}$. We therefore need only show the $\kappa_i$'s can be chosen such that the relations (\ref{ET_nRel}) and 
(\ref{ETwWreathRel})
hold, i.e. that
\begin{align}\label{norm_wreath_rel}
(\kappa_i T_i e_z)(a_1\otimes \dots \otimes a_d)(\kappa_i^{-1} T_i^{-1} e_z) = (-1)^{\sum_{j\neq i,i+1}|a_j|}({}^{s_i}(a_1\otimes \dots \otimes a_d)),
\end{align}
for all $a_1,\dots,a_d\in \blo^{\varnothing,1}$ and
\begin{subequations}
\begin{align}
(\kappa_i T_i e_z)^2 &= e_z,\label{norm_ti_rel_1}\\
(\kappa_i T_i e_z)(\kappa_j T_j e_z) &= -(\kappa_j T_j e_z)(\kappa_i T_i e_z),\text{ for all }|i-j|>1,\label{norm_ti_rel_2}\\
(\kappa_i T_i e_z)(\kappa_{i+1} T_{i+1} e_z)(\kappa_i T_i e_z) &= (\kappa_{i+1} T_{i+1} e_z)(\kappa_i T_i e_z)(\kappa_{i+1} T_{i+1} e_z).\label{norm_ti_rel_3}
\end{align}
\label{norm_ti_rel}
\end{subequations}

We first prove (\ref{norm_wreath_rel}) (note this relation does not depend on the choice of $\kappa_i$). Let $a\in \blo^{\varnothing,1}$. Using the notation of (\ref{EInsertion}), we consider the $a_1,a_i,a_j\in (\blo^{\varnothing,1})^{\otimes d}$, for $1<i\neq j\leq d$. Then
\begin{align*}
T_{(1,i)}a_1 T_{(1,i)}^{-1} = (-1)^{i|a|}a_i, \quad T_{(1,i)}a_i T_{(1,i)}^{-1} = (-1)^{i|a|}a_1, \quad T_{(1,i)}a_jT_{(1,i)}^{-1} = (-1)^{|a|}a_j,
\end{align*}
where the first two equalities hold by our identification of $\F N_{\tilde{\Si}_{V_1}}(D_1)e_{z}$ with $\F N_{\tilde{\Si}_{V_i}}(D_i)e_{z}$ via (\ref{NSp_ident}) and the third by (\ref{EZCommAlg}). It follows that
\begin{align*}
T_{(1,i)}(a_1\otimes \dots \otimes a_d)T_{(1,i)}^{-1} = (-1)^{i(|a_1|+|a_i|)+\sum_{j\neq 1,i}|a_j|}({}^{(1i)}(a_1\otimes \dots \otimes a_d)).
\end{align*}
Therefore, since $T_i = T_{(1,i)}T_{(1,i+1)}T_{(1,i)}$ or $zT_{(1,i)}T_{(1,i+1)}T_{(1,i)}$, a short calculation shows that
\begin{align*}
T_i(a_1\otimes \dots \otimes a_d)T_i^{-1} = (-1)^{\sum_{j\neq i,i+1}|a_j|}({}^{s_i}(a_1\otimes \dots \otimes a_d)).
\end{align*}

Similarly, (\ref{norm_ti_rel_2}) does not depend on the choice of $\kappa_i$ or $\kappa_j$ and follows from (\ref{EZCommAlg}). As $\pi_n(T_i)$ has order $2$, we have that $T_i$ has order $2$ or $4$ in $\tilde{\Si}_n$. Furthermore, since $\pi_n(T_i)$ is conjugate to $\pi_n(T_j)$ in $\Si_n$, $T_i$ is conjugate to $T_j$ or $zT_j$ in $\tilde \Si_n$, for all $1\leq i,j< d$. In particular, all $T_i$'s have the same order. In addition, since $\pi(T_i)$'s satisfy   braid relations, we have that
\begin{align}\label{braid}
T_iT_{i+1}T_i e_{z} = \pm T_{i+1}T_iT_{i+1} e_{z},
\end{align}
for all $1\leq i< d$. If $T_i$'s all have order $4$, we set all $\kappa_i=\sqrt{-1}$ to ensure that (\ref{norm_ti_rel_1}) holds. It follows from (\ref{braid}) that (\ref{norm_ti_rel_3}) holds up to a sign. We may therefore replace some of the $\kappa_i$'s with $-\kappa_i$ to ensure that (\ref{norm_ti_rel_3}) also holds. Note that this last reassignment will not stop (\ref{norm_ti_rel_1}) being satisfied.
%
%
%
\end{proof}

\subsection{On Brauer correspondents for blocks with abelian defect}
\label{SSBrCorr}

Let $P$ be a Sylow $p$-subgroup of $\tilde{\Si}_p$ and recall the blocks $\Blo^{\varnothing,1}, \blo^{\varnothing,1}$ defined in (\ref{EBlo}) and (\ref{Eblo}). 
In addition, we now need the definition of a $G$-graded crossed superproduct from $\S$\ref{SMoritaSuper}, 
and the notation $(A,B)_{C_2}$ from (\ref{EC_2}).

If $A$ and $B$ are two algebras, we write $$A\sim_{\der} B$$\index{$\sim_{\der}$} to indicate that $A$ and $B$ are derived equivalent, i.e. 
the bounded derived categories of $\D^b(A)$ and $\D^b(B)$ are equivalent as triangulated categories. If $A$ and $B$ are superalgebras, $A\sim_{\der} B$ just means that $|A|\sim_{\der} |B|$, i.e. $A$ and $B$ are derived equivalent as algebras.

\begin{Lemma}\label{lem:ext_rou}
There exists a complex of $(\Blo^{\varnothing,1}_\0\otimes (\blo^{\varnothing,1}_\0)^\op)$-modules inducing a derived equivalence between $\Blo^{\varnothing,1}_\0$ and $\blo^{\varnothing,1}_\0$, that extends to a complex of $(\Blo^{\varnothing,1},\blo^{\varnothing,1})_{C_2}$-modules.
\end{Lemma}

\begin{proof}
By~\cite[Theorem 4.4(b)]{Mu}, $\Blo^{\varnothing,1}_{\0}$ has Brauer tree 
\begin{align*}
\begin{braid}\tikzset{baseline=1.5mm}
\coordinate (1) at (0,0);
\coordinate (2) at (3,0);
\coordinate (3) at (6,0);
\coordinate (4) at (9,0);
\coordinate (5) at (12,0);
\coordinate (6) at (15,0);
\coordinate (7) at (18,0);
\coordinate (8) at (21,0);
\draw[thin] (1) -- (2);
\draw[thin] (2) -- (3);
\draw[thin] (3) -- (4);
\draw(10.5,0) node {$\cdots$};
\draw[thin] (5) -- (6);
\draw[thin] (6) -- (7);
\draw[thin] (7) -- (8);
\node at (1)[circle,fill,inner sep=1.5pt]{};
\node at (2)[circle,fill,inner sep=1.5pt]{};
\node at (3)[circle,fill,inner sep=1.5pt]{};
\node at (6)[circle,fill,inner sep=1.5pt]{};
\node at (7)[circle,fill,inner sep=1.5pt]{};
\node at (8)[circle,fill,inner sep=1.5pt]{};
\draw(14,0) node {$\cdots$};
\draw (1) circle (6mm);
\end{braid}
\end{align*}

\vspace{1.5mm}
\noindent
where we have $\ell+1$ nodes and the exceptional vertex has multiplicity $2$. In particular, since this tree has no non-trivial automorphims, conjugation by $\tilde{\Si}_p$ fixes each isomorphism class of irreducible $\Blo^{\varnothing,1}_{\0}$-module and hence each irreducible $\Blo^{\varnothing,1}_{\0}$-module extends to $\Blo^{\varnothing,1}$. Also, since $N_{\tilde{\Si}_p}(P)/P$ is abelian ($\cong C_{2(p-1)}$ or $C_{p-1}\times C_2$) each irreducible $\blo^{\varnothing,1}$-module is $1$-dimensional and so each irreducible $\blo^{\varnothing,1}_{\0}$-module must extend to $\blo^{\varnothing,1}$.

Throughout the remainder of this proof we make multiple references to~\cite{Rou}. As noted in the introduction, Brou\'e's conjecture was first proved for blocks with cyclic defect groups in \cite[Theorem 4.2]{Rickard_1}. However, the derived equivalence constructed in \cite{Rou} is more useful for our purposes as it is splendid. As well as making the current proof easier, this is also needed when we discuss splendid derived equivalences in more detail in $\S$\ref{sec:Splendid}.

First, note that, by~\cite[10.2.15]{Rou}, we have $\blo^{\varnothing,1} \sim_{\Mor} \F(P\rtimes E)$, where $E\cong C_{p-1}$ is the full automorphism group of $P\cong C_p$. Similarly $\blo^{\varnothing,1}_\0 \sim_{\Mor} \F (P\rtimes E')$, where $E'$ is the unique subgroup of $E$ of index $2$. (In fact, it is not hard to see that we even have isomorphisms
$\blo^{\varnothing,1} \cong \F(P\rtimes E)$ and 
$\blo^{\varnothing,1}_\0 \cong \F (P\rtimes E')$.) 
Using these Morita equivalences, we now follow~\cite[\S 10.3.4]{Rou} to construct our complex.

By~\cite[Proposition 6.1]{Lin2} and \cite[5.4,\,5.5]{Lin4}, we have that the $(\Blo^{\varnothing,1}_\0\otimes(\blo^{\varnothing,1}_\0)^\op)$-module $e_{\varnothing,1} \F\tilde{\Ai}_p e_{z}$ 
has a unique 
non-projective, indecomposable summand $M$ 
which induces a stable equivalence of Morita type between $\Blo^{\varnothing,1}_\0$ and $\blo^{\varnothing,1}_\0$. 
We apply the construction of~\cite[\S 10.3.4]{Rou} to $M$. Indeed it is stated in \cite[\S 4]{Rou} that this construction should be done with respect to such an $M$ to ensure our derived equivalence is splendid.

It is noted at the beginning of \cite[\S 10.3.4]{Rou} that
\begin{align*}
\bP:=\bigoplus_{0\leq i\leq |E'|-1}\bP_{\Omega^{2i}S}\otimes \bP_{\Omega^{2i}T}^*
\end{align*}
is a projective cover of $M$, where $T$ corresponds to $\F$, the trivial $\F(P\rtimes E')$-module through a fixed Morita equivalence between $\blo^{\varnothing,1}_\0$ and $\F(P\rtimes E')$, $S\cong M \otimes_{\blo^{\varnothing,1}_\0}T$, $\Omega$ denotes the Heller translate, and $\bP_V$ denotes a projective cover of $V$. As stated in the proof of~\cite[Proposition 10.3.7]{Rou} and the remarks preceding it, the $\Omega^{2i}S$'s all have a unique irreducible quotient and the $\Omega^{2i}T$'s are all irreducible. In particular, all the $\bP_{\Omega^{2i}S}$'s and $\bP_{\Omega^{2i}T}$'s are indecomposable. By \cite[Theorem 10.2]{Rou}, the complex inducing the derived equivalence between $\Blo^{\varnothing,1}_\0$ and $\blo^{\varnothing,1}_\0$ is now given by
\begin{equation}\label{EComplexX}
X:\,=0\to N\mathrel{\mathop{\to}^{\varphi}} M\to 0,
\end{equation}
where $M$ is in degree zero,
\begin{align*}
N=\bigoplus_{i\in Y}\bP_{\Omega^{2i}S}\otimes \bP_{\Omega^{2i}T}^*,
\end{align*}
for some particular subset $Y\subseteq \{0,1,\dots,|E'|-1\}$ and $\varphi$ is the restriction to $N$ of any surjection $\bP\to M$. 

Before continuing we transport our setup to the group algebra setting. Set
\begin{align*}
\funA = \langle\tilde{\Ai}_p \times N_{\tilde{\Ai}_p}(P),(s,s)\rangle\leq \tilde{\Si}_p \times N_{\tilde{\Si}_p}(P),
\end{align*}
for some fixed $s\in N_{\tilde{\Si}_p}(P)\backslash N_{\tilde{\Ai}_p}(P)$. We have the  isomorphism of superalgebras
\begin{align*}
\F \funA (e_{\varnothing,1}\otimes e_{z})&\iso (\Blo^{\varnothing,1},\blo^{\varnothing,1})_{C_2},\\ 
(a_1, a_2)(e_{\varnothing,1}\otimes e_{z})&\mapsto a_1 e_{\varnothing,1}\otimes a_2^{-1}e_{z}, \\
(s,s)(e_{\varnothing,1}\otimes e_{z})&\mapsto s e_{\varnothing,1}\otimes s^{-1} e_{z},
\end{align*}
for all $a_1\in \tilde{\Ai}_p$ and $a_2\in N_{\tilde{\Ai}_p}(P)$. We identify $\Blo^{\varnothing,1}_\0\otimes (\blo^{\varnothing,1}_\0)^{\op}$ with the subalgebra $\F(\tilde{\Ai}_p \times N_{\tilde{\Ai}_p}(P))(e_{\varnothing,1}\otimes e_{z})$ of $\F \funA (e_{\varnothing,1}\otimes e_{z})$.

In this language our claim is now that we can choose $\varphi$ such that the complex $X$ from (\ref{EComplexX}) extends to a complex of $\F \funA$-modules. Certainly $e_{\varnothing,1}\F\tilde{\Ai}_p e_{z}$ extends to an $\F \funA$-module via $(s,s)\cdot x=sxs^{-1}$, for all $x\in e_{\varnothing,1}\F\tilde{\Ai}_p e_{z}$. In particular, $(s,s)e_{\varnothing,1}\F\tilde{\Ai}_p e_{z}\cong e_{\varnothing,1}\F\tilde{\Ai}_p e_{z}$, as $\F(\tilde{\Ai}_p\times N_{\tilde{\Ai}_p}(P))$-modules. As it is the unique non-projective summand of $e_{\varnothing,1}\F\tilde{\Ai}_p e_{z}$, we have $(s,s)M\cong M$ and $M$ must also extend to an $\F \funA$-module. We denote one such extension $\tilde{M}$.

Recall that all $\bP_{\Omega^{2i}S}$'s and $\bP_{\Omega^{2i}T}$'s are indecomposable and so $\bP_{\Omega^{2i}S}\otimes \bP_{\Omega^{2i}T}^*$ is also indecomposable. Therefore, if $M/\rad{M}\cong S_1\oplus\dots\oplus S_r$, for some $r\in \Z_{> 0}$ and irreducible $\F(\tilde{\Ai}_p \times N_{\tilde{\Ai}_p}(P))$-modules $S_i$, then
\begin{align*}
\{\bP_{S_1},\dots,\bP_{S_r}\}\subseteq\{\bP_{\Omega^{2i}S}\otimes \bP_{\Omega^{2i}T}^*\}_{0\leq i\leq |E|/2-1}
\end{align*}
as multisets. Since $[\funA:\tilde{\Ai}_p\times N_{\tilde{\Ai}_p}(P)]=2$, any $\F \funA$-module is semisimple if and only if its restriction to $\F( \tilde{\Ai}_p\times N_{\tilde{\Ai}_p}(P))$ is. Therefore, $\rad{M} = \rad{\tilde{M}}$ and
\begin{align*}
\Res_{\tilde{\Ai}_p\times N_{\tilde{\Ai}_p}(P)}^{\funA}\left(\tilde{M}/\rad{\tilde{M}}\right)\cong M/\rad{M}.
\end{align*}
Given that each irreducible $\F\tilde{\Ai}_p e_{\varnothing,1}$-module extends to an $\F\tilde{\Si}_p e_{\varnothing,1}$-module and every irreducible $\F N_{\tilde{\Ai}_p}(P)e_{z}$-module extends to an $\F N_{\tilde{\Si}_p}(P)e_z$-module, we can extend each $S_i$ to an $\F \funA$-module. In particular,
\begin{align*}
\tilde{M}/\rad{\tilde{M}} \cong \tilde{S}_1 \oplus \dots \oplus \tilde{S}_r
\end{align*}
where each $\tilde{S}_i$ is an $\F \funA$-module extending $S_i$. It is therefore possible to choose extensions $\bQ_i$ of each individual $\bP_{\Omega^{2i}S}\otimes \bP_{\Omega^{2i}T}^*$ such that
\begin{align*}
\{\bP_{\tilde{S}_1},\dots,\bP_{\tilde{S}_r}\}\subseteq\{\bQ_i\}_{0\leq i\leq |E|/2-1}.
\end{align*}
Consequently $\bigoplus_{0\leq i\leq |E'|-1}\bQ_i$ surjects on to $\tilde{S}_1 \oplus \dots \oplus \tilde{S}_r$ and hence $\tilde{M}$. If we now set $\tilde{N}=\bigoplus_{i\in Y}\bQ_i$ and $\tilde{\varphi}:\tilde{N}\to \tilde{M}$ to be the restriction to $\tilde{N}$ of a surjection $\bigoplus_{0\leq i\leq |E'|-1}\bQ_i \to \tilde{M}$, then
\begin{align*}
\tilde{X}:=0\to \tilde{N}\mathrel{\mathop{\to}^{\tilde{\varphi}}} \tilde{M}\to 0
\end{align*}
is an extension of the complex $X$ from (\ref{EComplexX}) to a complex of $\F \funA$-modules, 
where $\varphi$ is simply $\tilde{\varphi}$ viewed as an $\F(\tilde{\Ai}_p\times N_{\tilde{\Ai}_p}(P))$-module homomorphism.
\end{proof}

Recall from Theorem~\ref{thm:brauer} the Brauer correspondents $\blo^{\rho,d}$ and $\blo^{\rho,d}_\0$ of the blocks $\Blo^{\rho,d}$ and $\Blo^{\rho,d}_\0$, respectively. 

\begin{Proposition}\label{lem:Brauer_der}
Let $0<d<p$ and $(\rho,d)\in\Xi(n)$ with $r:=|\rho|$. 
\begin{enumerate}
\item
\begin{align*}
\blo^{\rho,d} \sim_{\der}
\begin{cases}
\Blo^{\varnothing,1} \swr \cT_d &\text{if }r-h(\rho)\text{ is even,}\\
(\Blo^{\varnothing,1}\swr \cT_d)_{\0} &\text{if }r-h(\rho)\text{ is odd,}
\end{cases}
\end{align*}
and
\begin{align*}
\blo^{\rho,d}_\0 \sim_{\der}
\begin{cases}
(\Blo^{\varnothing,1}\swr \cT_d)_{\0} &\text{if }r-h(\rho)\text{ is even,}\\
\Blo^{\varnothing,1}\swr \cT_d &\text{if }r-h(\rho)\text{ is odd.}
\end{cases}
\end{align*}
\item In particular,
\begin{align*}
\blo^{\rho,d} \sim_{\der} \Blo^{\rho,0}\otimes (\Blo^{\varnothing,1} \swr \cT_d)\quad\text{and}\quad \blo^{\rho,d}_\0 \sim_{\der} \big(\Blo^{\rho,0}\otimes (\Blo^{\varnothing,1} \swr \cT_d)\big)_{\0}.
\end{align*}
\end{enumerate}
\end{Proposition}

\begin{proof}
We first note that, by Proposition~\ref{prop:Brauer_corr} and Theorem \ref{thm:brauer},
\begin{align}\label{BS_BA_isom}
\blo^{\rho,d} \cong \Blo^{\rho,0}\otimes (\blo^{\varnothing,1}\swr \cT_d).
\end{align}
If $r = 0$ or $1$, then $e_{\rho,0}=e_z$ (see the comments preceding Theorem \ref{thm:brauer} and Theorem \ref{thm:An_blocks}). In either case
\begin{align}\label{algn:r_0_1}
\Blo^{\rho,0}\otimes (\blo^{\varnothing,1}\swr \cT_d) \cong (\blo^{\varnothing,1}\swr \cT_d).
\end{align}
Until further notice we assume $r>1$.

By Theorems~\ref{thm:Sn_blocks} and~\ref{thm:An_blocks},
\begin{align*}
|\Blo^{\rho,0}|\cong
\begin{cases}
\cM_{2m}(\F) &\text{if }r-h(\rho)\text{ is even,} \\
\cM_m(\F) \oplus \cM_m(\F) &\text{if }r-h(\rho)\text{ is odd,} \\
\end{cases}
\end{align*}
and
\begin{align*}
|\Blo^{\rho,0}_\0|\cong
\begin{cases}
\cM_m(\F)\oplus \cM_m(\F) &\text{if }r-h(\rho)\text{ is even,} \\
\cM_m(\F) &\text{if }r-h(\rho)\text{ is odd,} \\
\end{cases}
\end{align*}
for some $m\in \Z_{> 0}$.

Let us first assume $r-h(\rho)$ is even. The above shows that any decomposition of $e_{\rho,0}$ into primitive idempotents in $\Blo^{\rho,0}_\0$ remains primitive in $\Blo^{\rho,0}$. Therefore, if $e\in \Blo^{\rho,0}_\0$ is a primitive idempotent, we have $e\Blo^{\rho,0}e\cong \F$ and, by Lemma~\ref{lem:idmpt_Mor}, $\Blo^{\rho,0} \sim_{\sM} \F$. Lemma~\ref{lem:MSE} now gives
\begin{equation}\label{E181221}
\Blo^{\rho,0}\otimes (\blo^{\varnothing,1}\swr \cT_d)\sim_{\sM} \blo^{\varnothing,1}\swr \cT_d.
\end{equation}
Note that, due to (\ref{algn:r_0_1}), the above also holds when $r=0$ or $1$. In this case $r-h(\rho)$ must also be even. Therefore, we cease distinguishing between $r\leq 1$ and $r>1$, since we have shown (\ref{E181221}) to hold whenever $r-h(\rho)$ is even.

When $r-h(\rho)$ is odd, first note that $\sigma_{\F\tilde{\Si}_{r}}$ must swap the two matrix factors in $\Blo^{\rho,0}$, since $\sigma_{\F \tilde{\Si}_n}(e_{\rho,0}^\pm)=e_{\rho,0}^\mp$. Therefore, any primitive idempotent in $\Blo^{\rho,0}_\0$ is the sum of two primitive idempotents, one in each of the matrix factors, in $\Blo^{\rho,0}$. Let $e:=e_1+e_2$ be such a decomposition of a primitive idempotent $e\in \Blo^{\rho,0}_\0$ into orthogonal, primitive idempotents $e_1,e_2\in \Blo^{\rho,0}$. 
Then $e \Blo^{\rho,0} e = \F e_1 \oplus \F e_2$. We can now construct an isomorphism of superalgebras
\begin{align*}
e \Blo^{\rho,0} e  \iso \cC_1,\ e_1\mapsto \frac{1+\cc}{2},\ 
e_2\mapsto \frac{1-\cc}{2}.
\end{align*}
Note this isomorphism must respect the superstructures as $e_1$ and $e_2$ get swapped by $\sigma_{\F\tilde{\Si}_r}$. So, by Lemma~\ref{lem:idmpt_Mor}, $\Blo^{\rho,0} \sim_{\sM} \cC_1$ and, by Lemma~\ref{lem:MSE},
\begin{equation*}\label{E181221_2}
\Blo^{\rho,0}\otimes (\blo^{\varnothing,1} \swr \cT_d)\sim_{\sM} \cC_1\otimes (\blo^{\varnothing,1} \swr \cT_d).
\end{equation*}
Therefore, by (\ref{BS_BA_isom}) and Lemmas~\ref{lem:Mor_A0} and \ref{lem:clif_mat2}, we now have that
\begin{align}\label{mor_even_brauer}
\blo^{\rho,d} \sim_{\Mor}
\begin{cases}
\blo^{\varnothing,1} \swr \cT_d, &\text{if }r-h(\rho)\text{ is even,}\\
(\blo^{\varnothing,1} \swr \cT_d)_{\0}, &\text{if }r-h(\rho)\text{ is odd,}
\end{cases}
\end{align}
and
\begin{align}\label{mor_odd_brauer}
\blo^{\rho,d}_\0 \sim_{\Mor}
\begin{cases}
(\blo^{\varnothing,1} \swr \cT_d)_{\0}, &\text{if }r-h(\rho)\text{ is even,}\\
\blo^{\varnothing,1}\swr \cT_d, &\text{if }r-h(\rho)\text{ is odd.}
\end{cases}
\end{align}

Let $X$ be the complex from Lemma~\ref{lem:ext_rou}. Using that result and Proposition~\ref{prop:ext_der_wreath} we can extend $X^{\otimes d}$ to a complex of $(\Blo^{\varnothing,1} \swr \cT_d,\blo^{\varnothing,1} \swr \cT_d)_{C_2\wr \Si_d}$-modules. It now follows from~\cite[Theorem 3.4(b)]{Mar} that $\Blo^{\varnothing,1}\swr \cT_d$ is derived equivalent to $\blo^{\varnothing,1} \swr \cT_d$. (Note that, since $d<p$, $p\nmid |C_2\wr \Si_d|$ and so we can indeed apply~\cite[Theorem 3.4(b)]{Mar}.)

Next consider the group homomorphism
\begin{align*}
C_2\wr {\Si}_d\to \{\pm 1\},\ g_i\mapsto -1,\ 
s_j\mapsto -1
\end{align*}
for $1\leq i\leq d$ and $1\leq j\leq d-1$, where $g_i$ is the generator of the $i^{\nth}$ $C_2$ and $s_j$ is the appropriate elementary transposition in ${\Si}_d$. Let ${\mathcal K}$ be the kernel of this homomorphism. Note that $(\Blo^{\varnothing,1}\swr \cT_d)_{\0}$ is a ${\mathcal K}$-graded (super) crossed product with identity component $(\Blo^{\varnothing,1}_\0)^{\otimes d}$ and $(\blo^{\varnothing,1}\swr \cT_d)_{\0}$ is a ${\mathcal K}$-graded (super) crossed product with identity component $(\blo^{\varnothing,1}_\0)^{\otimes d}$. As above we can extend $X^{\otimes d}$ to a complex of $((\Blo^{\varnothing,1} \swr \cT_d)_{\0},(\blo^{\varnothing,1} \swr \cT_d)_{\0})_{{\mathcal K}}$-modules. Again, it now follows from~\cite[Theorem 3.4(b)]{Mar} that $(\Blo^{\varnothing,1}\swr \cT_d)_{\0}$ is derived equivalent to $(\blo^{\varnothing,1} \swr \cT_d)_{\0}$. Part (i) now follows.

For part (ii) we simply prove
\begin{align*}
\Blo^{\rho,0} \otimes (\Blo^{\varnothing,1} \swr \cT_d) \sim_{\Mor}
\begin{cases}
\Blo^{\varnothing,1} \swr \cT_d, &\text{if }r-h(\rho)\text{ is even,}\\
(\Blo^{\varnothing,1} \swr \cT_d)_{\0}, &\text{if }r-h(\rho)\text{ is odd,}
\end{cases}
\end{align*}
and
\begin{align*}
\left(\Blo^{\rho,0} \otimes (\Blo^{\varnothing,1} \swr \cT_d)\right)_\0 \sim_{\Mor}
\begin{cases}
(\Blo^{\varnothing,1} \swr \cT_d)_{\0}, &\text{if }r-h(\rho)\text{ is even,}\\
\Blo^{\varnothing,1}\swr \cT_d, &\text{if }r-h(\rho)\text{ is odd,}
\end{cases}
\end{align*}
in exactly the same way (\ref{mor_even_brauer}) and (\ref{mor_odd_brauer}) were proved. The result now follows from part (i).
\end{proof}


\subsection{Splendid derived equivalences}\label{sec:Splendid}

We take a brief aside to discuss splendid derived equivalences. This is not needed anywhere else in this article and can be safely ignored by any reader only interested in the main results of this work. However, in the hope of being useful to future study in the area, we prove that the derived equivalences in Proposition~\ref{lem:Brauer_der}(ii) give rise to splendid derived equivalences between blocks. In particular, this will be helpful for any future work on the strengthening of Brou\'e's conjecture for double covers from a derived equivalence to a splendid derived equivalence. 

Let $G$ be a finite group and $b$ a block idempotent of $\F G$ with corresponding defect group $Q\leq G$. A {\em source idempotent}\index{source idempotent} of $\F Gb$ is a primitive idempotent $i\in (\F Gb)^Q$ such that $\Br_Q(i)\neq 0$. (Since $\Br_Q(b)\neq 0$, source idempotents certainly always exist.) Let $H$ be another finite group and $c$ a block idempotent of $\F H$ with corresponding defect group isomorphic to $Q$ that we identify with $Q$ and $j\in (\F Hc)^Q$ a source idempotent of $\F Hc$. A {\em splendid derived equivalence}\index{splendid derived equivalence} between $\F Gb$ and $\F Hc$ is a derived equivalence induced by a complex $X$ of $(\F Gb,\F Hc)$-bimodules such that in each degree we have a finite direct sum of summands of the $(\F Gb,\F Hc)$-bimodules $\F Gi \otimes_{\F R} j\F H$, where $R$ runs over the subgroups of $Q$, cf. \cite[1.10]{Lin3}.

It will be useful to think of source idempotents in an alternative way outlined in \cite{ALR}. The key point for us is \cite[Remark 3]{ALR}. That is, an idempotent $i\in (\F Gb)^Q$ is a source idempotent of $\F Gb$ if and only if $\F Gi$ is an indecomposable $\F(G\times Q)$-module with vertex $\Delta Q := \{(x,x)|x\in Q\}$.

An {\em interior $Q$-algebra}\index{interior algebra} is an algebra $A$ with a group homomorphism $Q\to A^\times$ called the {\em structural homomorphism}.\index{structural homomorphism} In the set up of the previous two paragraphs, $i\F Gi$ is an interior $Q$-algebra with structural homomorphism $$Q \to (i\F Gi)^\times,\ x\mapsto ixi.$$ An isomorphism of interior $Q$-algebras is an algebra isomorphism that commutes with the respective structural homomorphisms.

\begin{Lemma}\label{lem:source_idmpt}
Let $G$ be a finite group with normal subgroup $H$. Let $c$ be a block idempotent of $\F H$ and $b$ be a block idempotent of $\F G$ such that $bc=c$ and $\F Hc$ and $\F Gb$ have common defect group $Q\leq H$. If $C_G(Q) \leq H$, then any source idempotent $i$ of $\F Hc$ is also a source idempotent of $\F Gb$.
\end{Lemma}

\begin{proof}
Let $i \in (\F Hc)^Q$ be a source idempotent of $\F Hc$, in particular $i\in (\F Gb)^Q$. Then $\F Hi$ is an indecomposable $\F(H \times Q)$-module with vertex $\Delta Q$. Now
\begin{align}\label{ind_res_source}
\F Gi\downarrow^{G\times Q}_{H \times Q} \cong \F Hi \uparrow_{H \times Q}^{G \times Q}\downarrow^{G \times Q}_{H \times Q} \cong \bigoplus_{g\in G/H} g\F Hi,
\end{align}
and each $g\F Hi$ has vertex
\begin{align*}
{}^g\Delta Q:=\{(gxg^{-1},x)|x\in Q\}.
\end{align*}
Let $g_1,g_2\in G$ and suppose ${}^{g_1}\Delta Q$ is conjugate to ${}^{g_2}\Delta Q$ in $H \times Q$. Then
\begin{align*}
\{(h_1g_1xg_1^{-1}h_1^{-1},u_1xu_1^{-1})|x\in Q\} = \{(h_2g_2xg_2^{-1}h_2^{-1},u_2xu_2^{-1})|x\in Q\},
\end{align*}
for some $h_1,h_2 \in H$ and $u_1,u_2 \in Q$. Therefore,
\begin{align*}
\{(u_1^{-1}h_1g_1xg_1^{-1}h_1^{-1}u_1,x)|x\in Q\} = \{(u_2^{-1}h_2g_2xg_2^{-1}h_2^{-1}u_2,x)|x\in Q\}.
\end{align*}
This implies
\begin{align*}
g_2^{-1}h_2^{-1}u_2u_1^{-1}h_1g_1 \in C_G(Q)\leq H.
\end{align*}
Since $H$ is normal in $G$ we must have $g_2^{-1}g_1 \in H$. We have now shown that all the $g\F Hi$'s in (\ref{ind_res_source}) are non-isomorphic. Therefore, by \cite[\S5, Propositon 2]{Ward}, $\F Gi$ is indecomposable as an $\F(G \times Q)$-module. Furthermore, $\F Gi \cong \F Hi \uparrow_{H \times Q}^{G \times Q}$ certainly has vertex $\Delta Q$ and so $i$ is a source idempotent of $\F Gb$.
\end{proof}

We adopt the notation from Propositions~\ref{prop:Brauer_corr} and~\ref{lem:Brauer_der}. In particular, we take the $D\leq \tilde{\Si}_n$ from Proposition \ref{prop:Brauer_corr} as the defect group of $\Blo^{\rho,d}$.

We first need to view $\Blo^{\rho,0}\otimes (\Blo^{\varnothing,1} \swr \cT_d)$ and $\left(\Blo^{\rho,0}\otimes (\Blo^{\varnothing,1} \swr \cT_d)\right)_{\0}$ as blocks of finite groups. Recall the partition (\ref{E241221}). 
Set $$G^{\rho,d} \cong \Si_r \times (\Si_p\wr\Si_d)$$ to be the subgroup of $\Si_n$ permuting the $V_i$'s and $\tilde{G}^{\rho,d}:=\pi_n^{-1}(G^{\rho,d})$. Letting $\tilde{L}^{\rho,d}$ be the normal subgroup of $\tilde{G}^{\rho,d}$ generated by $\tilde{\Si}_U$ and all the $\tilde{\Si}_{V_i}$'s, we have that $\tilde{G}^{\rho,d} = \bigcup_{w\in \Si_d}\tilde{L}^{\rho,d} T_w$ is the coset decomposition of $\tilde{G}^{\rho,d}$ with respect to $\tilde{L}^{\rho,d}$, where the $T_w$'s are those from the proof of Proposition \ref{prop:Brauer_corr}. Now, by (\ref{EZCommAlg}),
$$
\F \tilde{L}^{\rho,d}e_z \cong
\F \tilde{\Si}_U e_z \otimes \F \tilde{\Si}_{V_1} e_z \otimes \dots \otimes \F \tilde{\Si}_{V_d} e_z.
$$
Moreover, as in (\ref{NSp_ident}), we can identify $\F \tilde{\Si}_{V_k}e_{z}$ with
$\F \tilde{\Si}_{V_1}e_{z} \cong \F \tilde{\Si}_p e_z$, for $1<k\leq d$, via the isomorphism
\begin{align}\label{Sp_ident}
\F\tilde{\Si}_{V_1}e_{z} \to \F\tilde{\Si}_{V_k}e_{z},\quad
a \mapsto (-1)^{k|a|}T_{(1,k)}aT_{(1,k)}^{-1}.
\end{align}
We set
\begin{align*}
f_{\rho,d}\index{f@$f_{\rho,d}$}&:= e_{\rho,0}\otimes \underbrace{e_{\varnothing,1}\otimes \dots \otimes e_{\varnothing,1}}_\text{$d$ times}  \in \F \tilde{\Si}_r e_z \otimes (\F \tilde{\Si}_p e_z)^{\otimes d}\\
&\cong \F \tilde{\Si}_U e_z \otimes \F \tilde{\Si}_{V_1} e_z \otimes \dots \otimes \F \tilde{\Si}_{V_d} e_z \cong \F \tilde{L}^{\rho,d}e_z
\end{align*}
and identify $f_{\rho,d}$ with its image in $\F \tilde{L}^{\rho,d}$. Then, since $e_{\rho,0}$ and all the $e_{\varnothing,1}$'s are even,
$$
\F \tilde{L}^{\rho,d} f_{\rho,d} \cong \Blo^{\rho,0} \otimes (\Blo^{\varnothing,1})^{\otimes d}.
$$
We can now proceed exactly as in the proof of Proposition~\ref{prop:Brauer_corr} to extend this to an isomorphism \begin{align}
\label{tildeG_isom}
\F \tilde{G}^{\rho,d}f_{\rho,d} \cong \Blo^{\rho,0}\otimes (\Blo^{\varnothing,1} \swr \cT_d),\ 
T_w f_{\rho,d} \mapsto \kappa_w^{-1}(e_{\rho,0} \otimes t_w)
\end{align}
for all $w\in \Si_d$ and some non-zero scalars $\ka_w$. Now, since by Proposition~\ref{lem:Brauer_der}(ii), we have 
$
\blo^{\rho,d} \sim_{\der} \Blo^{\rho,0}\otimes (\Blo^{\varnothing,1} \swr \cT_d),
$ 
the center of $\F \tilde{G}^{\rho,d}f_{\rho,d}$ is a local $\F$-algebra and hence $f_{\rho,d}$ is a block idempotent of $\F \tilde{G}^{\rho,d}$.


For any subgroup $H\leq \tilde\Si_n$, we denote $H_\0:=H\cap\tilde\Ai_n$. 
For example, we have 
\begin{align*}
\F \tilde{G}^{\rho,d}_{\0}f_{\rho,d} \cong \big(\Blo^{\rho,0}\otimes (\Blo^{\varnothing,1} \swr \cT_d)\big)_{\0} \sim_{\der} \blo^{\rho,d}_\0,
\end{align*}
so $f_{\rho,d}$ is a block idempotent of $\F \tilde{G}^{\rho,d}_{\0}$.

Denote $\tilde{N}^{\rho,d} := N_{\tilde{\Si}_n}(D)$. 
We have the following isomorphisms all of which are obtained by restricting either the isomorphism in Proposition~\ref{prop:Brauer_corr} or that in (\ref{tildeG_isom}) to the appropriate subsuperalgebra:
\begin{align}
\begin{split}\label{subsupalg_isom}
&\F (\tilde{N}^{\rho,d}\cap \tilde{\Si}_V)e_z \cong \blo^{\varnothing,1} \swr \cT_d, \hspace{1.04cm} \F (\tilde{N}^{\rho,d}\cap \tilde{\Ai}_V)e_z \cong (\blo^{\varnothing,1} \swr \cT_d)_{\0},
\\
&\F (\tilde{G}^{\rho,d}\cap \tilde{\Si}_V)e_{\varnothing,1}^{\otimes d} \cong \Blo^{\varnothing,1} \swr \cT_d,\qquad \F (\tilde{G}^{\rho,d}\cap \tilde{\Ai}_V)e_{\varnothing,1}^{\otimes d} \cong (\Blo^{\varnothing,1} \swr \cT_d)_{\0}.
\end{split}
\end{align}

Let $\tilde{H}^{\rho,d}$ be the subgroup of $\tilde{\Si}_n$ generated by $\tilde{\Si}_U$ and all the $\tilde{\Ai}_{V_i}$'s. By (\ref{EZCommAlg}),
$$
\F \tilde{H}^{\rho,d}e_z \cong
\F \tilde{\Si}_U e_z \otimes \F \tilde{\Ai}_{V_1} e_z \otimes \dots \otimes \F \tilde{\Ai}_{V_d} e_z.
$$
Note that this is a super tensor product of superalgebras but all the $\F \tilde{\Ai}_{V_i}e_z$'s are totally even and so all the factors actually commute with one another.

\begin{Lemma}\label{lem:source_idmpts_Sn}
Let $i_\rho$ be a primitive idempotent in $\Blo^{\rho,0}$, $\bar{i}_\rho$ a primitive idempotent in $\Blo^{\rho,0}_{\0}$ and $i\in (\Blo^{\varnothing,1}_{\0})^P$ a source idempotent of $\Blo^{\varnothing,1}_{\0}$. Then:
\begin{enumerate}
\item We have that $e_z$ is a source idempotent of both $\blo^{\varnothing,1}_{\0}$ and $\blo^{\varnothing,1}$.
\item We have that $i$ is a source idempotent of $\Blo^{\varnothing,1}$.
\item We have that
\begin{align*} 
i_\rho \otimes e_z^{\otimes d} \in \Blo^{\rho,0} \otimes (\blo^{\varnothing,1}\swr \cT_d) \cong \F \tilde{N}^{\rho,d}e_{\rho,0}.
\end{align*}
is a source idempotent of $\F \tilde{N}^{\rho,d}e_{\rho,0}$, and
\begin{align*}
\bar{i}_\rho \otimes e_z^{\otimes d} \in \big(\Blo^{\rho,0} \otimes (\blo^{\varnothing,1}\swr \cT_d)\big)_{\0} \cong \F \tilde{N}^{\rho,d}_{\0}e_{\rho,0}.
\end{align*}
is a source idempotent of $\F \tilde{N}^{\rho,d}_{\0}e_{\rho,0}$.
\item We have that
\begin{align*}
i_\rho \otimes i^{\otimes d} \in \Blo^{\rho,0} \otimes (\Blo^{\varnothing,1})^{\otimes d} \cong \F \tilde{L}^{\rho,d}f_{\rho,d}
\end{align*}
is a source idempotent of $\F \tilde{G}^{\rho,d}f_{\rho,d}$, and
\begin{align*}
\bar{i}_\rho \otimes i^{\otimes d} \in \big(\Blo^{\rho,0} \otimes (\Blo^{\varnothing,1})^{\otimes d}\big)_{\0} \cong \F \tilde{L}^{\rho,d}_{\0}f_{\rho,d}
\end{align*}
is a source idempotent of $\F \tilde{G}^{\rho,d}_{\0}f_{\rho,d}$.
\end{enumerate}
\end{Lemma}

\begin{proof}
(i) Since $e_z$ is certainly a source idempotent of $\F(P \times \langle z \rangle)e_z$, Lemma \ref{lem:source_idmpt} now gives the result.

(ii) This is a direct application of Lemma \ref{lem:source_idmpt}.

(iii) Let's first assume $r-h(\rho)$ is even. Certainly $i_\rho$ is a source idempotent for $\Blo^{\rho,0}$. It follows from \cite[Lemma 2.3]{EL} that
\begin{align*}
i_\rho \otimes e_z^{\otimes d} \in \F \left(\tilde{\Si}_U \times N_{\tilde{\Ai}_{V_1}}(P)\times \dots \times N_{\tilde{\Ai}_{V_d}}(P)\right)(e_{\rho,0} \otimes e_z^{\otimes d})
\end{align*}
is a source idempotent. Note we are taking a direct product of groups above and so this is not a subgroup of $\tilde{\Si}_n$. However, as interior $D$-algebras,
\begin{align*}
\F \left(\tilde{\Si}_U \times N_{\tilde{\Ai}_{V_1}}(P)\times \dots \times N_{\tilde{\Ai}_{V_d}}(P)\right)(e_{\rho,0} \otimes e_z^{\otimes d}) \cong \F (\tilde{N}^{\rho,d} \cap \tilde{H}^{\rho,d})e_{\rho,0}.
\end{align*}
It follows from the \cite{ALR} description of source idempotents that $i_\rho \otimes e_z^{\otimes d}$ is a source idempotent of $\F (\tilde{N}^{\rho,d} \cap \tilde{H}^{\rho,d})e_{\rho,0}$. Certainly $C_{\tilde{N}^{\rho,d}}(D) \leq \tilde{N}^{\rho,d} \cap \tilde{H}^{\rho,d}$ and Lemma \ref{lem:source_idmpt} gives the result.

If $r-h(\rho)$ is odd, without loss of generality we assume $i_\rho e_{\rho,0}^+ = i_\rho$. Certainly $i_\rho$ is a source idempotent for $\F \tilde{\Si}_r e_{\rho,0}^+$. Once again, $i_\rho \otimes e_z^{\otimes d}$ is a source idempotent of $\F(\tilde{N}^{\rho,d} \cap \tilde{H}^{\rho,d})e_{\rho,0}^+$. As in the $r-h(\rho)$ even case, the first claim follows from Lemma \ref{lem:source_idmpt}.

The second statement is proved similarly but with $\tilde{\Si}_U$ replaced by $\tilde{\Ai}_U$ and $\tilde{N}^{\rho,d}$ by $\tilde{N}^{\rho,d}_{\0}$.

(iv) The proof is identical to that of (iii) but with $\tilde{N}^{\rho,d}$ replaced by $\tilde{G}^{\rho,d}$, $i_\rho \otimes e_z^{\otimes d}$ by $i_\rho \otimes e_{\varnothing,1}^{\otimes d}$ and $\bar{i}_\rho \otimes e_z^{\otimes d}$ by $\bar{i}_\rho \otimes e_{\varnothing,1}^{\otimes d}$.
\end{proof}

We are now in a position to state and prove the main result of this subsection.

\begin{Proposition}
$\blo^{\rho,d}$ is splendidly derived equivalent to $\F \tilde{G}^{\rho,d}f_{\rho,d}$ and $\blo^{\rho,d}_{\0}$ is splendidly derived equivalent to $\F \tilde{G}^{\rho,d}_{\0}f_{\rho,d}$.
\end{Proposition}

\begin{proof}
The idea of the proof is that, through the isomorphisms 
$$\F \tilde{G}^{\rho,d}f_{\rho,d} \cong \Blo^{\rho,0}\otimes (\Blo^{\varnothing,1} \swr \cT_d)\quad \text{and} \quad\F \tilde{G}^{\rho,d}_{\0}f_{\rho,d} \cong \big(\Blo^{\rho,0}\otimes (\Blo^{\varnothing,1} \swr \cT_d)\big)_{\0},$$ we show that the derived equivalence between $\blo^{\rho,d}$ and $\Blo^{\rho,0}\otimes (\Blo^{\varnothing,1} \swr \cT_d)$ constructed in Proposition~\ref{lem:Brauer_der}(ii) is splendid.

We first note that the complex $X$ constructed in Lemma \ref{lem:ext_rou} induces a splendid derived equivalence between $\blo^{\varnothing,1}_{\0}$ and $\Blo^{\varnothing,1}_{\0}$. Indeed, this is stated in \cite[Theorem 1.1]{Rou}, our main reference for the proof of Lemma \ref{lem:ext_rou}. 

Now, $X^{\otimes d}$ induces a derived equivalence between
\begin{align*}
\F (\tilde{N}^{\rho,d}\cap \tilde{H}^{\rho,d} \cap \tilde{\Si}_V)e_z \cong \F N_{\tilde{\Ai}_{V_1}}(D_1)e_z\otimes \dots \otimes \F N_{\tilde{\Ai}_{V_d}}(D_d)e_z \cong (\blo^{\varnothing,1}_{\0})^{\otimes d}
\end{align*}
and
\begin{align*}
\F (\tilde{H}^{\rho,d} \cap \tilde{\Si}_V)e_{\varnothing,1}^{\otimes d} \cong \F \tilde{\Ai}_{V_1} e_{\varnothing,1}\otimes \dots \otimes \F \tilde{\Ai}_{V_d} e_{\varnothing,1} \cong (\Blo^{\varnothing,1}_{\0})^{\otimes d}.
\end{align*}
We adopt the labeling of idempotents from Lemma \ref{lem:source_idmpts_Sn}. Since every bimodule in each degree of $X$ is of the form
\begin{align*}
\F N_{\tilde{\Ai}_p}(P) e_z \otimes_{\F R} i\F\tilde{\Ai}_p,
\end{align*}
as $R$ ranges over the subgroups of $P$, every bimodule in each degree of $X^{\otimes d}$ is a direct summand of
\begin{align}\label{X^d:idmpts}
\F (\tilde{N}^{\rho,d}\cap \tilde{H}^{\rho,d} \cap \tilde{\Si}_V)e_z \otimes_{\F R} i^{\otimes d}\F (\tilde{H}^{\rho,d} \cap \tilde{\Si}_V),
\end{align}
as $R$ ranges over the subgroups of $D$. (Note that the $e_z$ above is really the $e_z^{\otimes d} \in \blo^{\varnothing,1}\swr \cT_d$ from Lemma \ref{lem:source_idmpts_Sn} interpreted as an element of $\F (\tilde{N}^{\rho,d}\cap \tilde{\Si}_V)$.) Now, using the identifications in (\ref{subsupalg_isom}), the proof of Proposition~\ref{lem:Brauer_der}(i) details a derived equivalence between $\F (\tilde{N}^{\rho,d}\cap \tilde{\Si}_V)e_z$ and $\F (\tilde{G}^{\rho,d}\cap \tilde{\Si}_V)e_{\varnothing,1}^{\otimes d}$. We describe how to construct the appropriate complex $Y$ in this group algebra setting. This involves looking in more detail at~\cite[Theorem 3.4(b)]{Mar}.

We first extend $X^{\otimes d}$ to a complex $\tilde{X}^{\otimes d}$ of modules for the group algebra of
\begin{align*}
\left\langle (\tilde{N}^{\rho,d}\cap \tilde{H}^{\rho,d} \cap \tilde{\Si}_V) \times (\tilde{H}^{\rho,d} \cap \tilde{\Si}_V), (u_i,u_i), (T_w,T_w)\right\rangle& \\
\leq (\tilde{N}^{\rho,d}\cap \tilde{\Si}_V) \times (\tilde{G}^{\rho,d}\cap \tilde{\Si}_V)&,
\end{align*}
as $i$ ranges over $\{1,\dots,d\}$ with each $u_i \in N_{\tilde{\Si}_{V_i}}(D_i) \backslash N_{\tilde{\Ai}_{V_i}}(D_i)$ and $w$ ranges over $\Si_d$ with the $T_w$'s those from the proof of Proposition~\ref{prop:Brauer_corr}. We then induce $\tilde{X}^{\otimes d}$ up to
\begin{align*}
(\tilde{N}^{\rho,d}\cap \tilde{\Si}_V) \times (\tilde{G}^{\rho,d}\cap \tilde{\Si}_V)
\end{align*}
to obtain $Y$. Since
\begin{align*}
p\nmid 2^{2d}(d!)^2 = [(\tilde{N}^{\rho,d}\cap \tilde{\Si}_V) \times (\tilde{G}^{\rho,d}\cap \tilde{\Si}_V):(\tilde{N}^{\rho,d}\cap \tilde{H}^{\rho,d} \cap \tilde{\Si}_V)\times (\tilde{H}^{\rho,d} \cap \tilde{\Si}_V)],
\end{align*}
it follows from (\ref{X^d:idmpts}) that any indecomposable summand appearing in each degree of $Y$ is a direct summand of
\begin{align*}
\F (\tilde{N}^{\rho,d} \cap \tilde{\Si}_V)e_z \otimes_{\F R} i^{\otimes d} \F (\tilde{G}^{\rho,d} \cap \tilde{\Si}_V),
\end{align*}
as $R$ ranges over the subgroups of $D$. Note that, by taking $\rho = \varnothing$ in Lemma \ref{lem:source_idmpts_Sn}(iii) and (iv), $e_z$ is a source idempotent of $\F (\tilde{N}^{\rho,d} \cap \tilde{\Si}_V)e_z$ and $i^{\otimes d}$ of $\F (\tilde{G}^{\rho,d} \cap \tilde{\Si}_V) e_{\varnothing,1}^{\otimes d}$. Therefore, $Y$ induces a splendid derived equivalence between $\F (\tilde{N}^{\rho,d}\cap \tilde{\Si}_V)e_z$ and $\F (\tilde{G}^{\rho,d}\cap \tilde{\Si}_V)e_{\varnothing,1}^{\otimes d}$.

We first assume $r-h(\rho)$ is even and show that $\F (\tilde{N}^{\rho,d}\cap \tilde{\Si}_V)e_z$ and $\F \tilde{N}^{\rho,d} e_{\rho,0}$ are splendidly derived equivalent. They are even splendidly Morita equivalent. That is, our complex can be chosen to be concentrated in degree zero. This, in turn, is equivalent to their source algebras being isomorphic as interior $D$-algebras, cf. \cite{Puig}. Since $r-h(\rho)$ is even, as in the proof of Proposition~\ref{lem:Brauer_der}(i), we can choose $i_\rho \in \F\tilde{\Ai}_U$. Then
\begin{align}
\label{D-alg}
\F (\tilde{N}^{\rho,d}\cap \tilde{\Si}_V)e_z \to i_\rho\F \tilde{N}^{\rho,d}i_\rho,\ 
x\mapsto i_\rho x i_\rho,
\end{align}
is such an isomorphism of source algebras. The only non-trivial things to check are injectivity and surjectivity. However, both follow immediately from the natural isomorphism
\begin{align*}
\F \tilde{N}^{\rho,d}e_z \cong \F \tilde{\Si}_U e_z \otimes \F (\tilde{N}^{\rho,d}\cap \tilde{\Si}_V) e_z.
\end{align*}
Similarly, the isomorphism
\begin{align*}
i^{\otimes d}\F (\tilde{G}^{\rho,d}\cap \tilde{\Si}_V)i^{\otimes d} \to (i_\rho \otimes i^{\otimes d})\F \tilde{G}^{\rho,d}(i_\rho \otimes i^{\otimes d}),\ 
x\mapsto i_\rho x i_\rho
\end{align*}
shows that $\F (\tilde{G}^{\rho,d}\cap \tilde{\Si}_V)e_{\varnothing,1}^{\otimes d}$ and $\F \tilde{G}^{\rho,d} f_{\rho,d}$ are splendidly Morita equivalent. We have now proved that $\F \tilde{G}^{\rho,d}f_{\rho,d}$ and $\F \tilde{N}^{\rho,d}e_{\rho,0}$ are splendidly derived equivalent when $r-h(\rho)$ is even.

To prove the remaining cases we must first prove that the derived equivalence between $(\blo^{\varnothing,1} \swr \cT_d)_{\0}$ and $(\Blo^{\varnothing,1} \swr \cT_d)_{\0}$, constructed in the proof of Proposition~\ref{lem:Brauer_der}(i), is splendid when viewed as a derived equivalence between $\F (\tilde{N}^{\rho,d}\cap \tilde{\Ai}_V)e_z$ and $\F (\tilde{G}^{\rho,d}\cap \tilde{\Ai}_V)e_{\varnothing,1}^{\otimes d}$. However, this follows in exactly the same manner as for the derived equivalence between $\F (\tilde{N}^{\rho,d}\cap \tilde{\Si}_V)e_z$ and $\F (\tilde{G}^{\rho,d}\cap \tilde{\Si}_V)e_{\varnothing,1}^{\otimes d}$.

Next we have the interior $D$-algebra isomorphism:
\begin{align}\label{D-alg1}
\F (\tilde{N}^{\rho,d}\cap \tilde{\Ai}_V)e_z \to \bar{i}_\rho\F \tilde{N}^{\rho,d}_{\0}\bar{i}_\rho,\ 
x\mapsto \bar{i}_\rho x \bar{i}_\rho,
\end{align}
when $r-h(\rho)$ is even. Since we are already choosing $i_\rho \in \F\tilde{\Ai}_U$, we can assume $\bar{i}_\rho = i_\rho$. Therefore, (\ref{D-alg1}) just follows from taking the even part of both sides of (\ref{D-alg}).

When $r-h(\rho)$ is odd, as noted in the proof of Proposition~\ref{lem:Brauer_der}(i), $\cC_1 \cong \bar{i}_\rho \F\tilde{\Si}_U \bar{i}_\rho$ as superalgebras. Let $\cc_\rho\in \bar{i}_\rho \F\tilde{\Si}_U \bar{i}_\rho$ be the image of $\cc$ under this isomorphism. We now define the $D$-algebra isomorphisms
\begin{align}\label{D-alg2}
\F (\tilde{N}^{\rho,d}\cap \tilde{\Ai}_V)e_z \to i_\rho\F \tilde{N}^{\rho,d}i_\rho,\ 
x\mapsto i_\rho x i_\rho
\end{align}
and
\begin{align}\label{D-alg3}
\begin{split}
\F (\tilde{N}^{\rho,d}\cap \tilde{\Si}_V)e_z&\to \bar{i}_\rho\F \tilde{N}^{\rho,d}_{\0}\bar{i}_\rho,\\ 
x&\mapsto 
\left\{
\begin{array}{ll}
\bar{i}_\rho x \bar{i}_\rho &\hbox{if $x\in \big(\F (\tilde{N}^{\rho,d}\cap \tilde{\Si}_V)e_z\big)_{\0}$,}\\
\sqrt{-1}\cc_\rho \bar{i}_\rho x \bar{i}_\rho &\hbox{if $x\in \big(\F (\tilde{N}^{\rho,d}\cap \tilde{\Si}_V)e_z\big)_{\1}$.}
\end{array}
\right.
\end{split}
\end{align}

Indeed, without loss of generality we assume $i_\rho e_{\rho,0}^+ = i_\rho$. Therefore,
\begin{align*}
i_\rho v i_\rho = i_\rho e_{\rho,0}^+ v e_{\rho,0}^+ i_\rho = i_\rho e_{\rho,0}^+ e_{\rho,0}^- v i_\rho = 0,
\end{align*}
for all $v\in \tilde{\Si}_V \backslash \tilde{\Ai}_V$. Since $i_\rho\F\tilde{\Si}_U i_\rho = \F i_\rho$, we have 
\begin{align*}
i_{\rho}\F (\tilde{N}^{\rho,d}\cap \tilde{\Ai}_V)i_\rho = i_{\rho}\F (\tilde{N}^{\rho,d}\cap \tilde{\Si}_V)i_\rho = i_{\rho}\F \tilde{N}^{\rho,d}i_\rho,
\end{align*}
and (\ref{D-alg2}) is surjective. Injectivity of (\ref{D-alg2}) follows since $i_\rho x i_\rho$ is precisely the image of $i_\rho \otimes x$ under the natural algebra isomorphism
\begin{align*}
\F\tilde{\Si}_U e_z \otimes \F (\tilde{N}^{\rho,d} \cap \tilde{\Ai}_V)e_z \to \F \tilde{N}^{\rho,d}e_z.
\end{align*}

For (\ref{D-alg3}), one can readily check that the homomorphism  is well-defined. To show its injectivity and surjectivity we identify $\bar{i}_\rho\F \tilde{N}^{\rho,d}\bar{i}_\rho$ with
\begin{align*}
\bar{i}_\rho\F \tilde{\Si}_U \bar{i}_\rho \otimes \F(\tilde{N}^{\rho,d} \cap \tilde{\Si}_V)e_z
\end{align*}
using (\ref{EZCommAlg}). Next we decompose $\bar{i}_\rho\F \tilde{N}^{\rho,d}_{\0}\bar{i}_\rho$ as
\begin{align*}
&\Big[\big(\bar{i}_\rho\F \tilde{\Si}_U \bar{i}_\rho\big)_{\0} \otimes \big(\F(\tilde{N}^{\rho,d} \cap \tilde{\Si}_V)e_z\big)_{\0}\Big] \oplus \Big[\big(\bar{i}_\rho\F \tilde{\Si}_U \bar{i}_\rho\big)_{\1} \otimes \big(\F(\tilde{N}^{\rho,d} \cap \tilde{\Si}_V)e_z\big)_{\1}\Big]\\
=& \Big[\F \bar{i}_\rho \otimes \big(\F(\tilde{N}^{\rho,d} \cap \tilde{\Si}_V)e_z\big)_{\0}\Big] \oplus \Big[\F \cc_\rho \otimes \big(\F(\tilde{N}^{\rho,d} \cap \tilde{\Si}_V)e_z\big)_{\1}\Big].
\end{align*}
One can now readily check that (\ref{D-alg3}) induces bijections
\begin{align*}
\big(\F (\tilde{N}^{\rho,d}\cap \tilde{\Si}_V)e_z\big)_{\0} \to \F \bar{i}_\rho \otimes \big(\F(\tilde{N}^{\rho,d} \cap \tilde{\Si}_V)e_z\big)_{\0}
\end{align*}
and
\begin{align*}
\big(\F (\tilde{N}^{\rho,d}\cap \tilde{\Si}_V)e_z\big)_{\1} \to \F \cc_\rho \otimes \big(\F(\tilde{N}^{\rho,d} \cap \tilde{\Si}_V)e_z\big)_{\1}.
\end{align*}
Therefore, (\ref{D-alg3}) is both injective and surjective.

The isomorphisms
\begin{align*}
i^{\otimes d} \F (\tilde{G}^{\rho,d}\cap \tilde{\Ai}_V) i^{\otimes d} \to (\bar{i}_\rho \otimes i^{\otimes d})\F \tilde{G}^{\rho,d}_{\0}(\bar{i}_\rho \otimes i^{\otimes d}),\ 
x\mapsto \bar{i}_\rho x \bar{i}_\rho,
\end{align*}
when $r-h(\rho)$ is even, and 
\begin{align*}
i^{\otimes d}\F (\tilde{G}^{\rho,d}\cap \tilde{\Ai}_V)i^{\otimes d} \to (i_\rho \otimes i^{\otimes d})\F \tilde{G}^{\rho,d}(i_\rho \otimes i^{\otimes d}),\ 
x\mapsto i_\rho x i_\rho
\end{align*}
and
\begin{align*}
i^{\otimes d}\F (\tilde{G}^{\rho,d}\cap \tilde{\Si}_V)i^{\otimes d}&\to (\bar{i}_\rho \otimes i^{\otimes d})\F \tilde{G}^{\rho,d}_{\0}(\bar{i}_\rho \otimes i^{\otimes d})\\
x&\mapsto 
\left\{
\begin{array}{ll}
\bar{i}_\rho x \bar{i}_\rho &\hbox{if $x\in \big(i^{\otimes d}\F (\tilde{G}^{\rho,d}\cap \tilde{\Si}_V)i^{\otimes d}\big)_{\0}$,}\\
\sqrt{-1} \cc_\rho \bar{i}_\rho x \bar{i}_\rho &\hbox{if $x\in \big(i^{\otimes d}\F (\tilde{G}^{\rho,d}\cap \tilde{\Si}_V)i^{\otimes d}\big)_{\1}$}
\end{array}
\right.
\end{align*}
when $r-h(\rho)$ is odd, are proved in a completely analogous way to (\ref{D-alg1}), (\ref{D-alg2}) and (\ref{D-alg3}).

With Lemma \ref{lem:source_idmpts_Sn} in mind, we have now shown that, when $r-h(\rho)$ is even, the following chain of Morita/derived equivalences are actually splendid
\begin{align*}
\blo^{\rho,d}_{\0} = \F \tilde{N}^{\rho,d}_{\0}e_{\rho,0} \sim_{\Mor} \F (\tilde{N}^{\rho,d} \cap \tilde{\Ai}_V)e_z \sim_{\der} \F (\tilde{G}^{\rho,d} \cap \tilde{\Ai}_V)e_{\varnothing,1}^{\otimes d} \sim_{\Mor} \F \tilde{G}^{\rho,d}_{\0} f_{\rho,d}.
\end{align*}
Similarly, when $r-h(\rho)$ is odd
\begin{align*}
\blo^{\rho,d} = \F \tilde{N}^{\rho,d}e_{\rho,0} \sim_{\Mor} \F (\tilde{N}^{\rho,d} \cap \tilde{\Ai}_V)e_z \sim_{\der} \F (\tilde{G}^{\rho,d} \cap \tilde{\Ai}_V)e_{\varnothing,1}^{\otimes d} \sim_{\Mor} \F \tilde{G}^{\rho,d} f_{\rho,d}
\end{align*}
and
\begin{align*}
\blo^{\rho,d}_{\0} = \F \tilde{N}^{\rho,d}_{\0}e_{\rho,0} \sim_{\Mor} \F (\tilde{N}^{\rho,d} \cap \tilde{\Si}_V)e_z \sim_{\der} \F (\tilde{G}^{\rho,d} \cap \tilde{\Si}_V)e_{\varnothing,1}^{\otimes d} \sim_{\Mor} \F \tilde{G}^{\rho,d}_{\0} f_{\rho,d}
\end{align*}
are both chains of splendid Morita/derived equivalences. This completes the proof.
\end{proof}

\section{Kang-Kashiwara-Tsuchioka isomorphism}
\label{SKKT}

We continue with our assumptions on $\F$ from Section \ref{sec:spin_blocks}. 
The Kang-Kashiwara-Tsuchioka isomorphism allows us to 
relate the spin blocks of $\cT_n$ to the quiver Hecke superalgebras $R^{\La_0}_\theta$ studied in Chapters~\ref{Part2} and \ref{ChRockqHs}. The relationship goes through some auxiliary algebras which are of independent interest and which we introduce in \S\ref{SSS}.

\subsection{Alternative description of superblocks of $\cT_n$}\label{sec:alt_descr}

We now give an alternative description of the superblocks of $\cT_n$ from~\cite{BKdurham} (see also \cite{BK1,Kbook}). In this subsection all our superalgebras are $\F$-algebras.
For $1\leq s<r\leq n$, we define
$$
[s,r]:=(-1)^{r-s-1}\ct_{r-1}\cdots \ct_{s+1}\ct_s\ct_{s+1}\cdots \ct_{r-1}\in\cT_n.\index{$[s,r]$}
$$
Now, the {\em spin Jucys-Murphy elements}\index{spin Jucys-Murphy elements} are defined as follows 
$$
\cm_r:=\sum_{s=1}^{r-1}[s,r] \in\cT_n\qquad(1\leq r\leq n).
\index{m@$\cm_r$}
$$
Alternatively, the elements $\cm_r$ can be defined inductively via:
\begin{equation}\label{ERecY}
\cm_1=0,\quad \cm_{r+1}=-\ct_r\cm_r\ct_r+\ct_r.
\end{equation}
Note that the elements $\cm_r$ are odd, so the elements $\cm_r^2$ are even. 

\begin{Lemma} \label{SMurphyCenter} {\rm \cite[Theorem 3.2]{BKdurham}} The elements $\cm_1^2,\dots,\cm_n^2$ commute, and 
$Z(\cT_n)_\0$ consists of all symmetric polynomials in $\cm_1^2,\dots, \cm_n^2$.
\end{Lemma}

Now recall that $\ell=(p-1)/2$ and we have the notation $I = \{0,1,\dots,\ell\}$, see \S\ref{ChBasicNotGen}. We consider the elements $i\in I$ as elements of $\F$, i.e. we identify $i=i\cdot 1_\F$. Let $V$ be a finite dimensional $\cT_n$-supermodule. For any $\bi=i_1\cdots i_n\in I^n$, we consider the simultaneous generalized eigenspace 
\begin{equation}\label{ETWtSp}
\hspace{3mm}
V_\bi:=\{v\in V\mid (\cm_r^2-i_r(i_r+1)/2)^Nv=0\ \text{for $N\gg0$ and $r=1,\dots,n$}\}.\index{$V_\bi$}
\end{equation}
By \cite[Lemma 3.3]{BKdurham}, we then have the `weight space decomposition' of $V$:
$$
V=\bigoplus_{\bi\in I^n}V_\bi.
$$
Considering the `weight space decomposition' of the regular $\cT_n$-module, we deduce that there is a system $\{e(\bi) \mid \bi \in I^n\}$ of mutually orthogonal idempotents in $\cT_n$ summing to the identity, uniquely determined by the property that $e(\bi)V = V_\bi$ for each $V\in\mod{\cT_n}$ (some of the $e(\bi)$ might be zero). 
In fact, each $e(\bi)$ lies in the commutative subalgebra generated by $\cm_1^2,\dots,\cm_n^2$.

Recalling the notation of \S\ref{ChBasicNotGen}, let $\theta\in Q_+$ with $\height(\theta)=n$. Set 
$$
e_\theta:=\sum_{\bi\in I^\theta}e(\bi)\in\cT_n.
$$

\begin{Lemma} \label{LSuperblocksDurham} {\rm \cite[Lemma 3.4]{BKdurham}}
Let $V\in\mod{\cT_n}$. 
The decomposition 
$$
V=\bigoplus_{\theta\in Q_+\, \text{\rm with}\,\, \height(\theta)=n}e_\theta V
$$
is precisely the decomposition of $V$ into superblocks.
\end{Lemma}

In particular, $e_\theta\in \cT_n$ is a central idempotent, and, setting $\cT_\theta:=e_\theta \cT_n$,\index{t@$\cT_\theta$} we have that 
\begin{equation}\label{EThetaBlocks}
\cT_n=\bigoplus_{\theta\in Q_+\, \text{\rm with}\,\, \height(\theta)=n}\cT_\theta
\end{equation}
is the superblock decomposition of $\cT_n$ (in general $\cT_\theta=0$ for some $\theta$'s).

In~\cite{BK1,BK2}, the irreducible $\cT_n$-supermodules were classified. Recall the set $\Par_p^{\res}(n)$ of $p$-restricted $p$-strict partitions of $n$ from \S\ref{SSGen}. In fact, \cite[Theorem 9.10]{BK1} and \cite[Theorem 10.3]{BK2} yield two different  canonical ways to label the irreducible $\cT_n$-supermodules by the elements of $\Par_p^{\res}(n)$. In \cite[Theorem B]{KShch}, it is shown that the two parametrizations of irreducible $\cT_n$-supermodules agree with each other. Thus using any of the parametrizations, we have a complete and irredundant set of irreducible $\cT_n$-supermodules
$$
\{D(\la)\mid \la\in \Par_p^{\res}(n)\}. \index{d@$D(\la)$}
$$

For $\la\in\Par_p(n)$, we denote by $\la^R\in\Par_p^{\res}(n)$\index{$\la^R$} the {\em regularization of $\la$},\index{regularization} as defined in \cite[\S2]{BKReg}.

\begin{Lemma} \label{LReg} 
For any $\la\in\Par_0(n)$, reduction modulo $p$ of the complex irreducible $\C \tilde{\Si}_n$-supermodule with character $\chi^\la$ contains $D(\la^R)$ as a composition factor. 
\end{Lemma}
\begin{proof}
This follows easily from \cite{BKReg} and \cite{BK2}. To give more details, let $E_\la$ be the symmetric polynomial and  $d_{\la,\mu}$ the decomposition numbers defined in \cite[\S10]{BK2}. In view of 
\cite[4.3]{BrundanQ}, this is the same $E_\la$ as used in \cite{BKReg}. By \cite[Theorem 4.4(ii)]{BKReg}, we have $d_{\la,\la^R}\neq 0$, which implies the lemma in view of \cite[Theorem 10.8]{BK2}.
\end{proof}


Since $\la^R$ is obtained from $\la$ by moving nodes along so-called ladders, and since all the nodes in a ladder have the same residue, we have $\cont(\la^R)=\cont(\la)$. It follows from Lemma~\ref{LReg} that the superblocks defined in (\ref{ERhoDBlocks}) and (\ref{EThetaBlocks}) match as follows:
$$
\Blo^{\rho,d}=\cT_{\cont(\rho)+d\de}.
$$

\subsection{Sergeev and related algebras}\label{sec:Sergeev}
\label{SSS}
The {\em affine Sergeev superalgebra}\index{affine Sergeev superalgebra} $\cX_n$\index{x@$\cX_n$}, defined by Nazarov \cite{Naz} (see also \cite[\S14.1]{Kbook}), has even generators $x_1,\dots,x_n$, $s_1,\dots,s_{n-1}$ and odd
generators $\cc_1,\dots,\cc_n$ subject only to the following relations for all admissible $i,j$:
\begin{eqnarray}
&&x_ix_j = x_jx_i,
\\
&&\cc_i^2=1,\ \cc_i\cc_j = -\cc_j\cc_i,
\\
&&s_i^2=1,\ s_is_{i+1}s_i=s_{i+1}s_is_{i+1}, 
\ \text{and}\ 
s_is_j=s_js_i\ \text{for}\ |i-j|>1,
\\
&&s_i\cc_i=\cc_{i+1}s_i,\ s_i\cc_{i+1}=\cc_is_i,
\ \text{and}\ 
s_i\cc_j=\cc_js_i\ \text{for}\ j\neq i,i+1, 
\\
&&\cc_ix_i = -x_i\cc_i,\ \text{and}\ \cc_ix_j = x_j\cc_i\ \text{for}\ i\neq j,
\\
&&s_ix_i = x_{i+1}s_i-1-\cc_i\cc_{i+1},\ \text{and}\ s_ix_j = x_js_i\ \text{for}\ j \neq i,i+1.
\end{eqnarray}

\begin{Lemma} \label{LAffSergCenter} {\rm \cite[Theorem 14.3.1]{Kbook}} 
The center of $\cX_n$ consists of all symmetric polynomials in $x_1^2,\dots,x_n^2$.
\end{Lemma}

We now suppose that $\cha \F=2\ell+1$. 
Let $\La=\sum_{i\in I}a_i\La_i\in P_+$ (dominant weight of Lie type $A_{2\ell}^{(2)}$). The {\em cyclotomic Sergeev superalgebra}\index{cyclotomic Sergeev superalgebra} $\cX_n^{\La}$,\index{x@$\cX_n^{\La}$} see \cite[\S15.3]{Kbook}, is defined to be $\cX_n$ modulo the relation
\begin{equation}
x_1^{a_0}\prod_{i\in I\backslash\{0\}} (x_1^2-i(i+1))^{a_i}=0.
\end{equation}


Fix $\La$ and $n$. Given $V\in\mod{\cX_n^\La}$ and $\bi = i_1\cdots i_n\in I^n$, we consider the simultaneous generalized eigenspace 
\begin{equation}\label{ETWtXLa}
\hspace{3mm}
V_\bi:=\{v\in V\mid (x_r^2-i_r(i_r+1))^Nv=0\ \text{for $N\gg0$ and $r=1,\dots,n$}\}.\index{$V_\bi$}
\end{equation}
By \cite[Definition 15.1.1, Lemma 15.3.1]{Kbook}, we then have a `weight space decomposition' 
$
V=\bigoplus_{\bi\in I^n}V_\bi.
$
Considering the `weight space decomposition' of the regular $\cX_n^\La$-module, we deduce that there is a system $\{e(\bi) \mid \bi \in I^n\}$ of mutually orthogonal idempotents in $\cX_n^\La$ summing to the identity uniquely determined by the property that $e(\bi)V = V_\bi$ for each $V\in\mod{\cX_n^\La}$ (some of the $e(\bi)$ might be zero). 
In fact, each $e(\bi)$ lies in the commutative subalgebra generated by $x_1^2,\dots,x_n^2$.

For $\theta\in Q_+$ be of height $n$, we set
$$
e_\theta:=\sum_{\bi\in I^\theta}e(\bi)\in \cX_n^\La.
$$
Let $\Ga_n$ denote the algebra of symmetric polynomials in $n$ variables. Then by Lemma~\ref{LAffSergCenter}, we have that $$A:=\{f(x_1^2,\dots,x_n^2)\mid f\in\Ga_n\}\subseteq \cX_n^\La$$ is a  central subalgebra. Note that for any $f\in\Ga_n$, $V\in\mod{\cX_n^\La}$ and a simultaneous eigenvector $v\in V$ for $x_1^2,\dots,x_n^2$ corresponding to the eigenvalues $i_1(i_1+1),\dots,i_n(i_n+1)$, we have 
$$f(x_1^2,\dots,x_n^2)v=f(i_1(i_1+1),\dots,i_n(i_n+1))v.$$ Moreover, if $\bi,\bj\in I^\theta$, then we have 
$$f(i_1(i_1+1),\dots,i_n(i_n+1))=f(j_1(j_1+1),\dots,j_n(j_n+1))$$ for all $f\in\Ga_n$, so there is a well-defined character 
$$\chi_\theta:A\to\F,\ f(x_1^2,\dots,x_n^2)\mapsto f(i_1(i_1+1),\dots,i_n(i_n+1)),$$ 
where $\bi\in I^\theta$. Note that 
$$e_\theta V=\{v\in V\mid (a-\chi_\theta(a))^Nv=0\ \text{for $N\gg0$ and all $a\in A$}\}.
$$ 
Applying this to the case where $V$ is the left regular $\cX_n^\La$-module, we deduce that $e_\theta\in A$, in particular, $e_\theta$ is a central idempotent in $\cX_n^\La$ (possibly zero). Set 
$$
\cX^\La_\theta:=e_\theta \cX_n^\La. \index{x@$\cX^\La_\theta$}
$$

We define the {\em Sergeev superalgebra}\index{Sergeev superalgebra} $\cY_n$,\index{t@$\cY_n$} see \cite[\S13.2]{Kbook}, to be 
$$
\cY_n:= \cT_n\otimes \cC_n.
$$
Taking $A=\F$ in Proposition~\ref{prop:semi_direct2}, yields the isomorphism of superalgebras  
\begin{equation}\label{Se_isom}
{\tt Se}:\cC_1\swr {\Si}_n\iso \cT_n\otimes \cC_n = \cY_n,
\end{equation}
first established by Sergeev \cite[\S2]{Ser} and Yamaguchi \cite{Yam}. 

\begin{Lemma}\label{LLevelOne}
There is an isomorphism of superalgebras 
$$\phi:\cX_n^{\La_0} \iso \cC_1 \swr {\Si}_n$$ 
such that 
 $\phi(\cc_i) = \cc_i$, $\phi(s_i) = s_i$, 
$\phi(x_1)=0$ and $\phi(x_{i+1}) = s_i \phi(x_i) s_i + s_i + \cc_i\cc_{i+1}s_i$ for all $i=0,\dots,n-1$. Moreover, for all $r=1,\dots,n$, $\bi\in I^n$, and $\theta\in Q_+$, we have 
$${\tt Se}(\phi(x_i^2))=2\cm_i^2\otimes 1,\quad 
{\tt Se}(\phi(e(\bi)))=e(\bi)\otimes 1\quad 
\text{and}\quad  {\tt Se}(\phi(e_\theta))=e_\theta\otimes 1.
$$
\end{Lemma}
\begin{proof}
The isomorphism $\phi$ is well-known, see e.g.  \cite[Remark 15.4.7]{Kbook}. The equality ${\tt Se}(\phi(x_i^2))=2\cm_i^2\otimes 1$ is easy to check and has been 
noted in the proof of~\cite[Lemma 3.3]{BKdurham}. The equality ${\tt Se}(\phi(e(\bi)))=e(\bi)\otimes 1
$ follows from the equalities ${\tt Se}(\phi(x_i^2))=2\cm_i^2\otimes 1$ comparing (\ref{ETWtSp}) and (\ref{ETWtXLa}). 
The final equality now follows from the definitions of $e_\theta$'s. 
\end{proof}

The $q$-analogue of the affine Sergeev superalgebra $\cX_n$ is the {\em affine Hecke-Clifford superalgebra}\index{affine Hecke-Clifford superalgebra} $\cX_n(q)$\index{x@$\cX_n(q)$} defined by Jones and Nazarov \cite{JN}. Now let $\F$ be an algebraically closed field of characteristic different from $2$, and the parameter $q\in\F^\times$ be a primitive $(2\ell+1)$st root of unity. Set $\xi:=q-q^{-1}$. We define $\cX_n(q)$ to be the superalgebra with even generators $X_1^{\pm},\dots,X_n^{\pm}$, $T_1,\dots,T_{n-1}$ and odd generators $C_1,\dots,C_n$ subject only to the following relations for all admissible $i,j$:
\begin{eqnarray}
&&
\hspace{5mm}X_iX_i^{-1} = X_i^{-1}X_i = 1,\ X_iX_j = X_jX_i,
\label{AHC1}
\\
&&
\hspace{5mm}C_i^2=1 \ \text{and}\ 
C_iC_j = -C_jC_i\ \text{for}\ i\neq j,
\label{AHC2}
\\
&&
\hspace{5mm}T_i^2 = \xi T_i + 1,\ T_iT_{i+1}T_i = T_{i+1}T_iT_{i+1},\ \text{and}\  T_iT_j = T_j T_i\ \text{for}\ |i - j| > 1,  
\label{AHC3}
\\
&&
\hspace{5mm}T_iC_i = C_{i+1}T_i\ \text{and}\  
T_iC_j = C_j T_i\ \text{for}\ j \neq i,i+1,
\label{AHC4}
\\
&&
\hspace{5mm}
C_iX_i^{\pm} = X_i^{\mp}C_i\ \text{and}\ 
C_iX_j = X_jC_i\ \text{for}\ i \neq j,
\label{AHC5}
\\
&&
\hspace{5mm}
(T_i + \xi C_iC_{i+1})X_iT_i = X_{i+1}\ \text{and}\ 
T_iX_j = X_j T_i\ \text{for}\ j\neq i,i+1.
\label{AHC6}
\end{eqnarray}

\begin{Lemma} \label{LAffSergCenterQ} {\rm \cite[Theorem 2.3]{BK1}} 
The center of $\cX_n(q)$ consists of all symmetric polynomials in $X_1+X_1^{-1},\dots,X_n+X_n^{-1}$.
\end{Lemma}

For any $i\in I$, we define
$$
q(i):=2\frac{q^{2i+1}+q^{-2i-1}}{q+q^{-1}}\in\F.
$$
For $\La=\sum_{i\in I}a_i\La_i\in P_+$, the {\em cyclotomic Hecke-Clifford superalgebra} $\cX_n^{\La}(q)$, see \cite[\S\S3-a,4-b]{BK1}, is defined to be $\cX_n(q)$ modulo the relation
\begin{equation}
(X_1-1)^{a_0}\prod_{i\in I\backslash\{0\}} \big(X_1+X_1^{-1}-q(i)\big)^{a_i}=0.
\end{equation}
The special case $\cX^{\La_0}_n(q)$ yields the Olshanski's {\em Hecke-Clifford superalgebra}\index{Hecke-Clifford superalgebra} which is a $q$-analogue $\cY_n(q)$\index{t@$\cY_n(q)$} of $\cY_n$ first considered in \cite{Ol}. As pointed out in the last paragraph of \S3-d in \cite{BK1}, alternatively, one can define $\cY_n(q)$ as the superalgebra generated by elements $C_i,T_j$ subject only to the relations (\ref{AHC2}), (\ref{AHC3}) and (\ref{AHC4}).

Fix $\La$ and $n$. Given $V\in\mod{\cX_n^\La(q)}$ and $i_1,\dots,i_n\in I$, we consider the simultaneous generalized eigenspace 
\begin{equation*}
\hspace{3mm}
V_\bi:=\{v\in V\mid (X_r+X_{r}^{-1}-q(i_r))^Nv=0\ \text{for $N\gg0$ and $r=1,\dots,n$}\}.\index{$V_\bi$}
\end{equation*}
By \cite[Corollary 4.5, Lemma 4.9]{BK1}, we then have a `weight space decomposition' 
$
V=\bigoplus_{\bi\in I^n}V_\bi.
$
Considering the `weight space decomposition' of the regular $\cX_n^\La(q)$-module, we deduce that there is a system $\{e(\bi) \mid \bi \in I^n\}$ of mutually orthogonal idempotents in $\cX_n^\La(q)$ summing to the identity uniquely determined by the property that $e(\bi)V = V_\bi$ for each $V\in\mod{\cX_n^\La(q)}$ (some of the $e(\bi)$ might be zero). 
In fact, each $e(\bi)$ lies in the commutative subalgebra generated by $X_1+X_1^{-1},\dots,X_n+X_n^{-1}$.  

For $\theta\in Q_+$ of height $n$, we set 
$$
e_\theta:=\sum_{\bi\in I^\theta}e(\bi)\in \cX_n^\La(q).
$$
Similarly to the $\cX_n^\La$ case, using the previous paragraph and Lemma~\ref{LAffSergCenterQ}, we have that $e_\theta$ is a central idempotent in $\cX_n^\La(q)$ (possibly zero). Set 
$$
\cX^\La_\theta(q):=e_\theta \cX_n^\La(q). \index{x@$\cX^\La_\theta(q)$}
$$

\subsection{Kang-Kashiwara-Tsuchioka isomorphisms}
\label{SSKKT}
Let 
$$\hat I:=\{0,1,1',2,2',\dots,\ell,\ell'\}\index{i@$\hat I$}$$ 
We have $I\subset \hat I$, and we also have the surjection 
$$
\pr:\hat I\to I,\ 0\mapsto 0,\ i\mapsto i,\ i'\mapsto i
\qquad(\text{for}\ i=1,\dots,\ell),
$$
and the involution 
$$
\ttc:\hat I\to\hat I,\ 0\mapsto 0,\ i\mapsto i',\ i'\mapsto i\qquad(\text{for}\ i=1,\dots,\ell).
$$
These induce the maps
\begin{align*}
&\pr:\hat I^n\to I^n,\ i_1\cdots i_n\mapsto \pr(i_1)\cdots \pr(i_n),
\\
&\ttc_r:\hat I^n\to \hat I^n,\ i_1\cdots i_n\mapsto i_1\cdots i_{r-1}(\ttc i_r)i_{r+1}\cdots i_n\qquad(1\leq r\leq n).
\end{align*}
For $\theta\in Q_+$ with $\height(\theta)=n$, we define 
$
\hat I^\theta=\{\bi\in \hat I^n\mid \pr(\bi)\in I^\theta\}.
$

Recalling the family of polynomials $\{Q_{i,j}(u,v)\mid i,j\in I\}$ from \S\ref{SQHDef}, we extend this to the larger family of polynomials $\{Q_{i,j}(u,v)\mid i,j\in \hat I\}$ using the rule
$$
Q_{\ttc i,j}(-u,v)=Q_{i,\ttc j}(u,-v)=Q_{i,j}(u,v) \qquad (i,j\in\hat I), 
$$
cf. \cite[(3.5)]{KKT}.

We now define the {\em quiver Hecke-Clifford superalgebra} \index{quiver Hecke-Clifford superalgebra}$RC_n$\index{r@$RC_n$} as the (unital) $\F$-algebra generated by the even generators $\{\tty_1,\dots,\tty_n\}\cup\{\si_1,\dots,\si_{n-1}\}\cup\{e(\bi)\mid \bi\in\hat I^n\}$ and the odd generators $\{\cc_1,\dots\cc_n\}$ subject only to the relations (i)--(x) of \cite[Definition 3.5]{KKT}. (Note that in \cite{KKT} the $Q_{i,j}$'s are denoted by $\tilde{Q}_{i,j}$.) Moreover, for $\La=\sum_{i\in I}a_i\La_i\in P_+$, the {\em cyclotomic quiver Hecke-Clifford superalgebra}\index{cyclotomic quiver Hecke-Clifford superalgebra} $RC_n^{\La}$\index{r@$RC_n^{\La}$} is defined as the quotient of $RC_n$ by the additional relations $\tty_1^{a_{i_1}}e(\bi)=0$ for all $\bi\in\hat I^n$. 

For any $\theta\in Q_+$ with $\height(\theta)=n$, we define a central idempotent 
$$
e_\theta:=\sum_{\bi\in \hat I^\theta}e(\bi)\in RC_n,
$$
and set
$$
RC_\theta=e_\theta RC_n\quad \text{and}\quad RC_\theta^{\La}=e_\theta RC_n^{\La}.
\index{r@$RC_\theta$}\index{r@$RC_\theta^{\La}$}
$$

Define also the idempotent 
$$e_I:=\sum_{\bi\in I^n}e(\bi)\in RC_n.$$ 
It follows from the comments just after~\cite[Definition 3.10]{KKT}, that $RC_\theta e_I RC_\theta = RC_\theta$ and $RC_\theta^{\La} e_I RC_\theta^{\La}  = RC_\theta^{\La}$. So, by Lemma~\ref{lem:idmpt_Mor}, we have:

\begin{Lemma}\label{lem:RC_trunc}
We have $RC_\theta\sim_{\sM} e_I RC_\theta e_I$ and $RC_\theta^{\La}\sim_{\sM}e_I RC_\theta^{\La} e_I$. 
\end{Lemma}

The following now readily follows from \cite[Theorem 3.13]{KKT}:

\begin{Proposition}\label{prop:RCtoR}
Let $\theta=\sum_{i\in I}m_i\al_i\in Q_+$.  
We have isomorphisms of superalgebras $e_I RC_\theta e_I\cong R_\theta\otimes \cC_{m_0}$ and $e_I RC_\theta^{\La} e_I\cong R_\theta^{\La}\otimes \cC_{m_0}$. 
\end{Proposition}
\begin{proof}
The isomorphism $e_I RC_\theta e_I\cong R_\theta\otimes \cC_{m_0}$ is established in \cite[Theorem 3.13]{KKT}. Note that in the proof of~\cite[Theorem 3.13]{KKT}, each $\tty_1 e(\bi)\in R_\theta$ gets sent to either $\tty_1 e(\bi)$ or $\cc_1 \tty_1 e(\bi)$ in $e_IRC_\theta e_I$. Since $\cc_1$ is invertible, passing this isomorphism to the cyclotomic setting is valid. 
\end{proof}

Recall the cyclotomic Sergeev superalgebras $\cX_n^\La$ and the cyclotomic Hecke Clifford superalgebras $\cX_n^\La(q)$ from \S\ref{SSS}. We now cite the key results of \cite{KKT}.


\begin{Theorem} \label{TKKTQ2} {\rm \cite[Theorem 5.4]{KKT}}
Let $\F$ be an algebraically closed field of odd characteristic $2\ell+1$.  
Then there is an isomorphisms of superalgebras $RC_n^{\La} \iso \cX^{\La}_n$ which maps 
$$\sum_{\hat{\bi}\in \hat I^n,\, \pr(\hat{\bi})=\bi}e(\hat{\bi})\mapsto e(\bi)\qquad(\bi\in I^n).
$$
In particular, $RC_\theta^{\La}$ and $\cX^{\La}_\theta$ correspond under this isomorphism, for any $\theta\in Q_+$ with $\height(\theta)=n$.
\end{Theorem}

\begin{proof}
It is stated in \cite[Theorem 5.4]{KKT} and \cite[p. 51]{KKT} 
that there is an isomorphism $\om:RC_n^{\Lambda} \iso  \cX^{\Lambda}_n$. To prove the correspondence of $e(\bi)$'s we must examine the explicit form of the isomorphism $\om$ constructed in the proof of \cite[Theorem 5.4]{KKT}.

We first note that in \cite[\S5.2]{KKT} the set $I$, which labels of the nodes of the Dynkin diagram of type $A_{2\ell}^{(2)}$, is given by the sequence $g(\frac{1}{2}),g(\frac{3}{2}),\dots,g(\frac{2\ell+1}{2})$ rather than $0,1,\dots,\ell$, where $g(k):=k^2-\frac{1}{4}$. In other words, $I$ defined in \cite{KKT} is obtained from the $I$ defined in this article via the transformation $i\mapsto i(i+1)$. Then $J$, denoted $\hat{I}$ in this article, is given by
\begin{align*}
J = \{0,\pm\sqrt{2},\dots,\pm\sqrt{i(i+1)},\dots,\pm\sqrt{\ell(\ell+1)}\}.
\end{align*}
Now, in \cite[Definition 5.3]{KKT} it states that $\om(e(\hat{\bi}))$ is a simultaneous generalized eigenvector for the $x_i$'s with corresponding eigenvalues
\begin{align*}
x(\hat{\bi}) := (x(i_1),\dots,x(i_n)),
\end{align*}
for all $\hat{\bi} \in \hat{I}^n$, where $x:J\to \mathbb{A}_x^1$ is just the inclusion function. (This notation of eigenspaces is described in more detail in \cite[\S4.3]{KKT}.) Since the simultaneous generalized eigenspaces for the $x_i$'s are certainly invariant under right multiplication by $\cX^{\La}_n$, it follows from the definition of $e(\bi)\in\cX^{\La}_n$ that
\begin{align*}
\om\Bigg(\sum_{\hat{\bi}\in \hat{I}^n,\, \pr(\hat{\bi})=\bi}e(\hat{\bi})\Bigg)\cX^{\La}_n \subseteq e(\bi)\cX_n^{\La}
\end{align*}
for all $\bi \in I^n$. Since the sum of the $e(\hat{\bi})$'s is $1\in RC_n^{\La}$, the inclusion above is actually an equality and we have the desired correspondence of idempotents.

The correspondence of blocks follows immediately.
\end{proof}

The following theorem is similar to Theorem~\ref{TKKTQ2} but uses \cite[Corollary 4.8]{KKT} instead of \cite[Theorem 5.4]{KKT} and \cite[p. 51]{KKT}.

\begin{Theorem} \label{TKKTQ1} {\rm \cite[Corollary 4.8]{KKT}} Let $\F$ be an algebraically closed field of characteristic different from $2$, and the parameter $q\in\F^\times$ be a primitive $(2\ell+1)$st root of unity. 
Then there is an isomorphisms of superalgebras $RC_n^{\La} \iso \cX^{\La}_n(q)$ which maps 
$$\sum_{\hat{\bi}\in \hat I^n, \pr(\hat{\bi})=\bi}e(\hat{\bi})\mapsto e(\bi)\qquad(\bi\in I^n).
$$
In particular, $RC_\theta^{\La}$ and $\cX^{\La}_\theta(q)$ correspond under this isomorphism, for any $\theta\in Q_+$ with $\height(\theta)=n$.
\end{Theorem}

\begin{Corollary}\label{Cor:QF2}
In the settings of Theorems \ref{TKKTQ2} and \ref{TKKTQ1}, $R_\theta^{\La}$ and $R_n^\La$ are $QF2$-algebras.
\end{Corollary}

\begin{proof}
It follows immediately from the definition of a $QF2$-algebra that $A$ is $QF2$ if $A/\rad(A)\cong \soc(A)$. We begin by showing that if $A$ is a Frobenius superalgebra, then $A\otimes \cC_m$ is a $QF2$-algebra, for all $m\in \Z_{> 0}$.

Let $A$ be a Frobenius algebra. First note that $\rad(A\otimes \cC_m) = \rad(A)\otimes \cC_m$, since $\rad(A)\otimes \cC_m$ is a nilpotent ideal of $A\otimes \cC_m$ with a semisimple quotient. Therefore,
\begin{align*}
(A\otimes \cC_m)/\rad(A\otimes \cC_m) \cong (A/\rad(A)) \otimes \cC_m \cong \soc(A) \otimes \cC_m \cong \soc(A\otimes \cC_m),
\end{align*}
where the second isomorphism follows from the fact that $A$ is Frobenius. Hence $A\otimes \cC_m$ is a $QF2$-algebra.


Now,
\begin{align*}
RC_{\theta}^{\La}\otimes \cC_{m_0} \mathrel{\mathop{\sim_{\sM}}^{\ref{lem:RC_trunc}}} e_I RC_{\theta}^{\La}e_I \otimes \cC_{m_0} \mathrel{\mathop{\cong}^{\ref{prop:RCtoR}}}R_\theta^{\La} \otimes \cC_{2m_0}\mathrel{\mathop{\sim_{\sM}}^{\ref{lem:clif_mat1}}}R_\theta^{\La},
\end{align*}
for all $\theta = \sum_{i\in I}m_i\al_i \in Q_+$.

The algebra $\cX^{\La}_n$, in the setting of Theorem \ref{TKKTQ2}, is Frobenius due to \cite[Corollary 15.6.2]{Kbook} and $\cX^{\La}_n(q)$, in the setting of Theorem \ref{TKKTQ1}, is Frobenius due to \cite[Corollary 3.14]{BK1}. It follows
from the isomorphisms of Theorems~\ref{TKKTQ2} and \ref{TKKTQ1} that in both cases, the algebra $RC_n^{\La}$ and hence each $RC_\theta^{\La}$ is also Frobenius, for all $\theta = \sum_{i\in I}m_i\al_i \in Q_+$, with $\height(\theta)=n$. Therefore, since being $QF2$ is a Morita invariant, each $R_\theta^{\La}$, and hence $R_n^{\La}$, is $QF2$.
\end{proof}

\begin{Theorem}\label{thm:label_match}
Let $\F$ be an algebraically closed field of odd characteristic $2\ell+1$ and $\theta = \sum_{i\in I}m_i\al_i \in Q_+$, with $\height(\theta)=n$. Then
\begin{align*}
\cT_\theta\sim_{\sM} R_\theta^{\Lambda_0}\otimes \cC_{n-m_0}.
\end{align*}
In particular, if $(\rho,d)\in\Xi(n)$ with $\cont(\rho) = \sum_{i\in I}r_i\al_i$, then
\begin{align*}
\Blo^{\rho,d}\sim_{\sM} R_{\cont(\rho)+d\de}^{\Lambda_0}\otimes \cC_{n-r_0}.
\end{align*}
\end{Theorem}

\begin{proof}
First note that
\begin{align*}
RC_{\theta}^{\La_0} \mathrel{\mathop{\sim_{\sM}}^{\ref{lem:RC_trunc}}} e_IRC_{\theta}^{\La_0}e_I \mathrel{\mathop{\cong}^{\ref{prop:RCtoR}}} R_{\theta}^{\La_0} \otimes \cC_{m_0}.
\end{align*}

Now, by Theorem \ref{TKKTQ2}, $RC_\theta^{\La_0}\cong \cX^{\Lambda_0}_\theta$ which, due to Lemma \ref{LLevelOne} and (\ref{Se_isom}), is isomorphic to $$e_\theta\cY_n = e_\theta\cT_n \otimes \cC_n = \cT_\theta \otimes \cC_n.$$
Therefore,
\begin{align*}
\cT_\theta \mathrel{\mathop{\sim_{\sM}}^{\ref{lem:clif_mat1}}} \cT_\theta \otimes \cC_{2n} \cong RC^{\La_0}_{\theta} \otimes \cC_n \sim_{\sM} R^{\La_0}_{\theta} \otimes \cC_{n+m_0} \mathrel{\mathop{\sim_{\sM}}^{\ref{lem:clif_mat1}}} R^{\La_0}_{\theta} \otimes \cC_{n-m_0}.
\end{align*}
For the second statement
\begin{align*}
\Blo^{\rho,d} = \cT_{\cont(\rho)+d\de} \sim_{\sM} R^{\La_0}_{\cont(\rho)+d\de} \otimes \cC_{n-m_0}\mathrel{\mathop{\sim_{\sM}}^{\ref{lem:clif_mat1}}} R^{\La_0}_{\cont(\rho)+d\de} \otimes \cC_{n-r_0},
\end{align*}
since $m_0=r_0+2d$. 
\end{proof}

\section{Brou\'e's conjecture for RoCK blocks of double covers}

We continue to work over the field $\F$ as in Section \ref{SKKT}. In particular, $$\cha\F=p=2\ell+1.$$ 

\subsection{The superblock $\Blo^{\varnothing,1}$}\label{sec:supblo_B}
In this section we investigate the superblock $\Blo^{\varnothing,1}=\F \tilde{\Si}_p e_{\varnothing,1}$, cf. (\ref{EBlo}). 
Recall the notation $I,\, J,\, K$ from \S\ref{ChBasicNot}. 
We define the $\F$-superalgebra $\Zag_\ell$\index{b@$\Zag_\ell$} to be the path algebra of the quiver
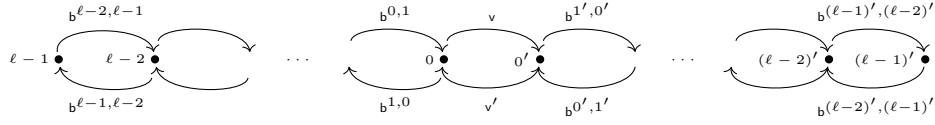
\begin{align*}
\begin{braid}\tikzset{baseline=3mm}
\coordinate (1) at (0,0);
\coordinate (2) at (4,0);
\coordinate (3) at (8,0);
\coordinate (4) at (12,0);
\coordinate (5) at (16,0);
\coordinate (6) at (20,0);
\coordinate (7) at (24,0);
\coordinate (8) at (28,0);
\coordinate (9) at (32,0);
\coordinate (10) at (36,0);
\draw [thin, black,->,shorten <= 0.1cm, shorten >= 0.1cm]   (1) to[distance=1.5cm,out=100, in=100] (2);
\draw [thin,black,->,shorten <= 0.25cm, shorten >= 0.1cm]   (2) to[distance=1.5cm,out=-100, in=-80] (1);
\draw [thin,black,->,shorten <= 0.25cm, shorten >= 0.1cm]   (2) to[distance=1.5cm,out=80, in=100] (3);
\draw [thin,black,->,shorten <= 0.25cm, shorten >= 0.1cm]   (3) to[distance=1.5cm,out=-100, in=-80] (2);
\draw [thin,black,->,shorten <= 0.25cm, shorten >= 0.1cm]   (4) to[distance=1.5cm,out=80, in=100] (5);
\draw [thin,black,->,shorten <= 0.25cm, shorten >= 0.1cm]   (5) to[distance=1.5cm,out=-100, in=-80] (4);
\draw [thin,black,->,shorten <= 0.25cm, shorten >= 0.1cm]   (5) to[distance=1.5cm,out=80, in=100] (6);
\draw [thin,black,->,shorten <= 0.1cm, shorten >= 0.1cm]   (6) to[distance=1.5cm,out=-100, in=-100] (5);
\draw [thin,black,->,shorten <= 0.25cm, shorten >= 0.1cm]   (6) to[distance=1.5cm,out=80, in=100] (7);
\draw [thin,black,->,shorten <= 0.1cm, shorten >= 0.1cm]   (7) to[distance=1.5cm,out=-100, in=-100] (6);
\draw [thin,black,->,shorten <= 0.25cm, shorten >= 0.1cm]   (8) to[distance=1.5cm,out=80, in=100] (9);
\draw [thin,black,->,shorten <= 0.1cm, shorten >= 0.1cm]   (9) to[distance=1.5cm,out=-100, in=-100] (8);
\draw [thin,black,->,shorten <= 0.25cm, shorten >= 0.1cm]   (9) to[distance=1.5cm,out=80, in=100] (10);
\draw [thin,black,->,shorten <= 0.1cm, shorten >= 0.1cm]   (10) to[distance=1.5cm,out=-100, in=-100] (9);
\blackdot(0,0);
\blackdot(4,0);
\blackdot(16,0);
\blackdot(20,0);
\blackdot(32,0);
\blackdot(36,0);
\draw(0,0) node[left]{$\ell-1$};
\draw(4,0) node[left]{$\ell-2$};
\draw(10,0) node {$\cdots$};
\draw(16,0) node[left]{$0$};
\draw(20,0) node[left]{$0'$};
\draw(26,0) node {$\cdots$};
\draw(32,0) node[left]{$(\ell-2)'$};
\draw(36,0) node[left]{$(\ell-1)'$};
\draw(2,1.2) node[above]{$\zb^{\ell-2,\ell-1}$};
\draw(6,1.2); 
\draw(2,-1.2) node[below]{$\zb^{\ell-1,\ell-2}$};
\draw(6,-1.2); 
\draw(14,1.2) node[above]{$\zb^{0,1}$};
\draw(18,1.2) node[above]{$\zv$};
\draw(22,1.2) node[above]{$\zb^{1',0'}$};
\draw(14,-1.2) node[below]{$\zb^{1,0}$};
\draw(18,-1.2) node[below]{$\zv'$};
\draw(22,-1.2) node[below]{$\zb^{0',1'}$};
\draw(30,1.2); 
\draw(34,1.2) node[above]{$\zb^{(\ell-1)',(\ell-2)'}$};
\draw(30,-1.2); 
\draw(34,-1.2) node[below]{$\zb^{(\ell-2)',(\ell-1)'}$};
\end{braid}
\end{align*}
generated by length $1$ paths 
$$\{\zv,\zv'\}\cup\{\zb^{k,k+1},\zb^{k',(k+1)'},\zb^{k+1,k},\zb^{(k+1)',k'}\mid k\in K\}$$  and  length $0$ paths 
$$\{\zf^j,\zf^{j'}\mid j\in J\},$$ modulo the following relations:
\begin{enumerate}
\item All paths of length three or greater are zero.
\item All paths of length two that are not cycles are zero.
\item The length-two cycles based at any fixed vertex are equal.
\end{enumerate}
The superstructure on $\Zag_{\ell}$ is determined by the fact that
\begin{align*}
\zf^j+\zf^{j'},\quad \zb^{k,k+1}+\zb^{k',(k+1)'},\quad \zb^{k+1,k}+\zb^{(k+1)',k'}, \quad \zv+\zv',
\end{align*}
for all $j\in J$ and $k\in K$ generate $(\Zag_{\ell})_{\0}$ and $\Zag_{\ell}$ has superunit
$
\sum_{j\in J}(\zf^j-\zf^{j'}).
$
Equivalently, the involution $\sigma_{\Zag_{\ell}}$ of (\ref{ESi}) is given by
\begin{align*}
\zf^j\leftrightarrow\zf^{j'},\quad \zb^{k,k+1}\leftrightarrow\zb^{k',(k+1)'},\quad \zb^{k+1,k}\leftrightarrow\zb^{(k+1)',k'}, \quad \zv\leftrightarrow\zv',
\end{align*}
for all $j\in J$ and $k\in K$.

\begin{Lemma}\label{lem:zigzag}
$\Zig_{\ell} \otimes \cC_1 \cong \Zag_{\ell}$.
\end{Lemma}

\begin{proof}
It is easily checked that the following assignments define mutually inverse isomorphisms of superalgebras:
\begin{align*}
&\Zig_{\ell}\otimes \cC_1 \to\Zag_{\ell} && \Zag_{\ell}\to\Zig_{\ell}\otimes \cC_1\\
&\ze^j\otimes 1 \mapsto \zf^j+\zf^{j'} && \zf^j\mapsto \ze^j \otimes (1+\cc)/2, \quad \zf^{j'}\mapsto \ze^j \otimes (1-\cc)/2\\
&\zu\otimes 1 \mapsto \zv - \zv' && \zv\mapsto \zu \otimes (1+\cc)/2, \quad \zv'\mapsto -\zu \otimes (1-\cc)/2\\
&\za^{k,k+1}\otimes 1 \mapsto \zb^{k,k+1} + \zb^{k',(k+1)'} && \zb^{k,k+1}\mapsto \za^{k,k+1} \otimes (1+\cc)/2 \\
&\za^{k+1,k}\otimes 1 \mapsto - \zb^{k+1,k} - \zb^{(k+1)',k'} && \zb^{k',(k+1)'}\mapsto \za^{k,k+1} \otimes (1-\cc)/2\\
&1\otimes \cc \mapsto {\textstyle\sum_{i\in J}(\zf^i - \zf^{i'})} && \zb^{k+1,k}\mapsto -\za^{k+1,k} \otimes (1+\cc)/2 \\
& && \zb^{(k+1)',k'}\mapsto -\za^{k+1,k} \otimes (1-\cc)/2,
\end{align*}
for all $j\in J$ and $k\in K$.
\end{proof}

\begin{Lemma}\label{lem:sym_S}
$\Blo^{\varnothing,1}\sim_{\sM}\Zag_{\ell}$.
\end{Lemma}

\begin{proof}
It is enough to prove the statement for $\F$ an algebraic closure of $\F_p$ and we make this assumption throughout the proof.

By~\cite[Theorem 4.4(a)]{Mu}, $\Blo^{\varnothing,1}$ has Brauer tree
\vspace{2mm}
\begin{align}\label{Brauer_tree_B_l}
\begin{braid}\tikzset{baseline=-0.5mm}
\coordinate (1) at (0,0);
\coordinate (2) at (3,0);
\coordinate (3) at (6,0);
\coordinate (4) at (9,0);
\coordinate (5) at (12,0);
\coordinate (6) at (15,0);
\coordinate (7) at (18,0);
\coordinate (8) at (21,0);
\draw[thin] (1) -- (2);
\draw[thin] (2) -- (3);
\draw[thin] (3) -- (4);
\draw(10.5,0) node {$\cdots$};
\draw[thin] (5) -- (6);
\draw[thin] (6) -- (7);
\draw[thin] (7) -- (8);
\node at (1)[circle,fill,inner sep=1.5pt]{};
\node at (2)[circle,fill,inner sep=1.5pt]{};
\node at (3)[circle,fill,inner sep=1.5pt]{};
\node at (6)[circle,fill,inner sep=1.5pt]{};
\node at (7)[circle,fill,inner sep=1.5pt]{};
\node at (8)[circle,fill,inner sep=1.5pt]{};
\draw(14,0) node {$\cdots$};
\end{braid}
\end{align}

\vspace{2.5mm}
\noindent
where we have $p$ nodes. As an algebra, $\Zag_{\ell}$ is the Brauer tree algebra with precisely this tree. In particular, $\Blo^{\varnothing,1}$ is Morita equivalent to $\Zag_{\ell}$. The content of this proof is in showing that they are Morita superequivalent. Recall (\ref{ESi}) and note that the superstructures of 
$\Blo^{\varnothing,1}$ and 
$\Zag_{\ell}$ are determined by $\sigma:=\sigma_{\Blo^{\varnothing,1}}$ and $\sigma':=
\sigma_{\Zag_{\ell}}$, respectively. Therefore, by Lemma~\ref{lem:idmpt_Mor}, it suffices to find an  idempotent $e\in \Blo^{\varnothing,1}_{\0}=\F\tilde{\Ai}_p e_{\varnothing,1}$ such that $eL\neq 0$ for any irreducible $\Blo^{\varnothing,1}$-supermodule $L$ and there is an isomorphism of algebras $\chi:\Zag_\ell\to e\Blo^{\varnothing,1} e$ such that $\si\circ \chi=\chi\circ \si'$.

By~\cite[Theorem 4.4(b)]{Mu}, $\Blo^{\varnothing,1}_{\0}$ has Brauer tree
\begin{align*}
\begin{braid}\tikzset{baseline=1.5mm}
\coordinate (1) at (0,0);
\coordinate (2) at (3,0);
\coordinate (3) at (6,0);
\coordinate (4) at (9,0);
\coordinate (5) at (12,0);
\coordinate (6) at (15,0);
\coordinate (7) at (18,0);
\coordinate (8) at (21,0);
\draw[thin] (1) -- (2);
\draw[thin] (2) -- (3);
\draw[thin] (3) -- (4);
\draw(10.5,0) node {$\cdots$};
\draw[thin] (5) -- (6);
\draw[thin] (6) -- (7);
\draw[thin] (7) -- (8);
\node at (1)[circle,fill,inner sep=1.5pt]{};
\node at (2)[circle,fill,inner sep=1.5pt]{};
\node at (3)[circle,fill,inner sep=1.5pt]{};
\node at (6)[circle,fill,inner sep=1.5pt]{};
\node at (7)[circle,fill,inner sep=1.5pt]{};
\node at (8)[circle,fill,inner sep=1.5pt]{};
\draw(14,0) node {$\cdots$};
\draw (1) circle (6mm);
\end{braid}
\end{align*}

\vspace{2mm}
\noindent
where we have $\ell+1$ nodes and the exceptional vertex has multiplicity $2$. 


Let $\{e_0,\dots,e_{\ell-1}\}$ a set of orthogonal, primitive idempotents in $\Blo^{\varnothing,1}_{\0}$ corresponding to the $\ell$ pairwise non-isomorphic simple $\Blo^{\varnothing,1}_{\0}$-modules. In particular, setting $e:=e_0+\dots+e_{\ell-1}$, we have that $e \Blo^{\varnothing,1}_{\0} e$ is basic and Morita equivalent to $\Blo^{\varnothing,1}_{\0}$.
Pick $w\in \tilde{\Si}_p \backslash \tilde{\Ai}_p$. Then ${}^we:=wew^{-1}={}^we_0+\dots+{}^we_{\ell-1}$ is, again, a sum of orthogonal, primitive idempotents in $\Blo^{\varnothing,1}_{\0}$. %
%
%
Therefore, there exists $\alpha\in (\Blo^{\varnothing,1}_{\0})^\times$ such that
\begin{align*}
\{{}^{\alpha w}e_0,\dots,{}^{\alpha w}e_{\ell-1}\}=\{e_0,\dots,e_{\ell-1}\}.
\end{align*}
However, by inspection of the Brauer tree, every algebra automorphism of $e\Blo^{\varnothing,1}_{\0} e$ must fix the isomorphism class of every irreducible module. Therefore, we even have ${}^{\alpha w}e_j=e_j$, for all $j\in J$. 
%
%
So, by possibly raising $\alpha w e_{\varnothing,1}$ to an odd power, we have a superunit $u\in \Blo^{\varnothing,1}$ of order a power of $2$ that centralizes each $e_i$. (Since $\F$ is an algebraic closure of $\F_p$, $\al w$ lives in a finite subring of $\F\tilde{\Ai}_p$ and hence has finite order.) Say $u$ has order $2^n$. For each $0\leq m\leq 2^n-1$, set
\begin{align*}
f_m=\frac{1}{2^n}\sum_{t=0}^{2^n-1}\zeta^{mt}u^t
\in \Blo^{\varnothing,1},
\end{align*}
where $\zeta\in \F$ is some fixed, primitive $(2^n)^{\nth}$ root of unity. We have
\begin{align*}
\sigma(f_m) &= \frac{1}{2^n}\sum_{t=0}^{2^n-1}\zeta^{mt}\sigma(u^t) 
\\&= \frac{1}{2^n}\sum_{t=0}^{2^n-1}(-1)^t\zeta^{mt}u^t 
\\&= \frac{1}{2^n}\sum_{t=0}^{2^n-1}\zeta^{mt+2^{n-1}t}u^t\\
&= \frac{1}{2^n}\sum_{t=0}^{2^n-1}\zeta^{(m+2^{n-1})t}u^t 
\\&=  f_{m+2^{n-1}},
\end{align*}
where the subscript of $f$ is considered modulo $2^n$. Note that $\sum_{m=0}^{2^n-1} f_m = e_{\varnothing,1}$. Therefore, for every $j\in J$, there must exist $0\leq m_j\leq 2^n-1$ such that $e_jf_{m_j}\neq0$ and so
$$0\neq \sigma(e_jf_{m_j})=e_j\sigma(f_{m_j})=e_jf_{m_j+2^{n-1}}.$$ 
Now, since $e_j$ is primitive and $[\tilde{\Si}_p:\tilde{\Ai}_p]=2$, $\F\tilde{\Si}_p e_j$ is the direct sum of at most two projective, indecomposable $\F\tilde{\Si}_p$-modules. Therefore, $e_j$ must decompose into the sum of at most two orthogonal, primitive idempotents in $\F\tilde{\Si}_p$. However, $e_j = \sum_{m=0}^{2^n-1}e_j f_m$ and so, in fact, $e_j=e_jf_{m_j}+e_jf_{m_j+2^{n-1}}$ is a decomposition of $e_j$ into precisely two orthogonal, primitive idempotents. We now have
\begin{align*}
e=e_0f_{m_0}+\sigma(e_0f_{m_0})+\dots+e_{\ell-1}f_{m_{\ell-1}}+\sigma(e_{\ell-1}f_{m_{\ell-1}})
\end{align*}
is a sum of orthogonal, primitive idempotents in $\Blo^{\varnothing,1}$. Moreover, since
\begin{align*}
\Blo^{\varnothing,1} e\Blo^{\varnothing,1}  = \F\tilde{\Si}_p\Blo^{\varnothing,1}_{\0}e\Blo^{\varnothing,1}_{\0}\F\tilde{\Si}_p= \F\tilde{\Si}_p\Blo^{\varnothing,1}_{\0}\F\tilde{\Si}_p
= \Blo^{\varnothing,1} 
\end{align*}
and $\Blo^{\varnothing,1} \sim_{\Mor} \Zag_{\ell}$ has exactly $p-1=2\ell$ pairwise non-isomorphic irreducible modules, $e\Blo^{\varnothing,1} e$ is basic and therefore isomorphic, as an algebra, to $\Zag_{\ell}$. In particular, there exists an algebra isomorphism $\varphi:\Zag_{\ell}\to e\Blo^{\varnothing,1} e$ such that
\begin{align*}
\varphi(\{\zf^0,\zf^{0'},\dots,\zf^{\ell},\zf^{\ell'}\}) = \{e_0f_{m_0},\sigma(e_0f_{m_0}),\dots,e_{\ell}f_{m_{\ell}},\sigma(e_{\ell}f_{m_{\ell}})\}.
\end{align*}
This implies the algebra involution $\sigma^{\varphi}=\varphi^{-1}\circ \sigma \circ \varphi$ of $\Zag_{\ell}$ restricts to a non-trivial permutation of the $\zf^j$'s and $\zf^{j'}$'s. Certainly the only non-trivial automorphism of the Brauer tree (\ref{Brauer_tree_B_l}) of $\Zag_{\ell}$ is given by reflecting along its vertical line of symmetry. Therefore, $\sigma^{\varphi}$ must swap $\zf^j$ and $\zf^{j'}$, for all $j\in J$. In other words,
\begin{align}\label{sigma_phi}
\sigma^{\varphi}(\zf^j) = \zf^{j'} = \sigma'(\zf^j),\quad \sigma^{\varphi}(\zf^{j'}) = \zf^j = \sigma'(\zf^{j'}),
\end{align}
for all $j\in J$. Since $\dim_k(\zf \Zag_{\ell}\zf')\leq 1$, for all distinct $\zf,\zf'\in \{\zf^0,\zf^{0'},\dots,\zf^{\ell},\zf^{\ell'}\}$, (\ref{sigma_phi}) implies there exist $\lambda_{k,k+1}$, $\lambda_{k+1,k}$, $\mu\in \F^{\times}$, such that
\begin{align*}
\sigma^{\varphi}(\zv) = \mu\zv', && \sigma^{\varphi}(\zb^{k,k+1}) = \lambda_{k,k+1} \zb^{k',(k+1)'}, && \sigma^{\varphi}(\zb^{k+1,k}) = \lambda_{k+1,k} \zb^{(k+1)',k'},\\
\sigma^{\varphi}(\zv') = \mu^{-1}\zv, && \sigma^{\varphi}(\zb^{k',(k+1)'}) = \lambda_{k,k+1}^{-1} \zb^{k,k+1}, && \sigma^{\varphi}(\zb^{(k+1)',k'}) = \lambda_{k+1,k}^{-1} \zb^{k+1,k},
\end{align*}
for all $k\in K$. Note that $\sigma^{\varphi}(\zv' \zv)=\zv\zv'$. Therefore $\sigma^{\varphi}(\zb^{0,1}\zb^{1,0})=\zb^{0',1'}\zb^{1',0'}$ and we can continue by induction to show that $\lambda_{k,k+1}\lambda_{k+1,k}=1$, for all $k\in K$. Now define the algebra automorphism $\vartheta:\Zag_{\ell}\to \Zag_{\ell}$ via
\begin{align*}
&\vartheta(\zf^j) = \zf^j, \qquad \vartheta(\zf^{j'}) = \zf^{j'}, \qquad \vartheta(\zv) = \zv, \qquad \vartheta(\zv') = \mu\zv', \\
&\vartheta(\zb^{k,k+1}) = \mu \zb^{k,k+1}, \qquad \vartheta(\zb^{k+1,k}) = \zb^{k+1,k}, \\ &\vartheta(\zb^{k',(k+1)'}) = \mu \lambda_{k,k+1} \zb^{k',(k+1)'}, \qquad \vartheta(\zb^{(k+1)',k'}) = \lambda_{k+1,k} \zb^{(k+1)',k'},
\end{align*}
for all $j\in J$ and $k\in K$. One can now check that $\sigma^{\varphi \circ \vartheta} = \vartheta^{-1} \circ \sigma^{\varphi} \circ \vartheta = \sigma'$. 
So we can take $\chi:=\varphi \circ \vartheta$. 
\end{proof}

\begin{Corollary} \label{C021121} 
$\Blo^{\varnothing,1} \otimes \cC_1$ is Morita superequivalent to $\Zig_{\ell}$.
\end{Corollary}
\begin{proof}
By Lemmas~\ref{lem:sym_S},\,\ref{lem:zigzag}, we have 
$\Blo^{\varnothing,1}  \otimes \cC_1 \mathrel{\mathop{\sim_{\sM}}} \Zag_{\ell} \otimes \cC_1\cong \Zig_{\ell} \otimes \cC_2$, which is Morita superequivalent to $\Zig_{\ell}$ thanks to Lemma~\ref{lem:clif_mat1}. 
\end{proof}

\begin{Corollary} \label{C021121_2} 
$(\Blo^{\varnothing,1} \otimes \cC_1)\swr {\Si}_d$ is Morita superequivalent to $\Zig_{\ell}\swr {\Si}_d$.
\end{Corollary}
\begin{proof}
As noted in \S\ref{SSZig}, $\Zig_{\ell}$ is a symmetric superalgebra. In addition, since $\cC_1$ is symmetric via $1\mapsto 1$, $\cc\mapsto 0$, $\Blo^{\varnothing,1} \otimes \cC_1$ is also a symmetric superalgebra by (\ref{tensor_symm_algs}). Corollary~\ref{C021121} gives that $\Blo^{\varnothing,1} \otimes \cC_1\sim_{\sM}\Zig_{\ell}$. We can now apply Proposition~\ref{prop:semi_direct1}(i).
\end{proof}

\subsection{Local description of RoCK blocks}
Let $d$ be a positive integer less than $p$ and $\rho$ a $d$-Rouquier ${\bar p}$-core. Set $r:=|\rho|$, $n:=|\rho|+dp$. By Corollary~\ref{Cor:QF2}, $R^{\La_0}_{\cont(\rho)+\de}$ is a $QF2$-algebra and we can apply Theorem~\ref{TMorWr} to get 
\begin{align}\label{algn:Rn_sM_ASd}
R^{\La_0}_{\cont(\rho)+d\de}\sim_{\sM} \Zig_{\ell}\swr {\Si}_d.
\end{align}

\begin{Lemma}\label{lem:Yd_N}
$R^{\La_0}_{\cont(\rho)+d\de}$ is Morita superequivalent to $(\Blo^{\varnothing,1} \swr \cT_d)\otimes \cC_d$.
\end{Lemma}

\begin{proof}
By (\ref{algn:Rn_sM_ASd}) and Corollary~\ref{C021121_2}, we have
$$R^{\La_0}_{\cont(\rho)+d\de} \sim_{\sM}\Zig_{\ell}\swr {\Si}_d\sim_{\sM} (\Blo^{\varnothing,1} \otimes \cC_1)\swr {\Si}_d.$$
Therefore, by Proposition~\ref{prop:semi_direct2}, $R^{\La_0}_{\cont(\rho)+d\de} \sim_{\sM} (\Blo^{\varnothing,1} \swr \cT_d)\otimes \cC_d$.

\end{proof}

Let $\rho'$ be an arbitrary ${\bar p}$-core, and set $r':=|\rho'|$ and $n':=|\rho'|+dp$. Write $\operatorname{cont}(\rho')=\sum_{i\in I}r_i\alpha_i$. We now have, by Theorem~\ref{thm:label_match}: 
 $${\mathfrak B}^{\rho',d}\sim_{\operatorname{sMor}}R^{\Lambda_0}_{\operatorname{cont}(\rho')+d\delta}\otimes {\mathcal C}_{n'-r_0}.
 $$
Note that
\begin{align*}
n'-r_0+d = r' - r_0 + d(p+1) \equiv r'-r_0 \pmod 2\equiv r'-h(\rho') \pmod 2,
\end{align*}
where the final congruence follows since $\rho'$ has no parts divisible by $p$ and so each row of $\rho'$ has an odd number of boxes with residue $0$. Taking into account Lemma~\ref{lem:clif_mat1} we get
\begin{align}\label{algn:Blo_R_sMor}
{\mathfrak B}^{\rho',d}&\sim_{\operatorname{sMor}}
\begin{cases}
R^{\Lambda_0}_{\operatorname{cont}(\rho')+d\delta} & \text{if $r-h(\rho')+d$ is even},\\
R^{\Lambda_0}_{\operatorname{cont}(\rho')+d\delta} \otimes \cC_1 & \text{if $r-h(\rho')+d$ is odd}.
\end{cases}
\end{align}

\begin{Proposition} \label{P021221_3}
We have the following Morita superequivalences:
\begin{enumerate}
\item 
\begin{align*}
\Blo^{\rho,d} \sim_{\sM}
\begin{cases}
\Zig_\ell \swr \Si_d &\text{if }r-h(\rho)\text{ is even,}\\
(\Zig_\ell \swr \Si_d)\otimes \cC_1&\text{if }r-h(\rho)\text{ is odd.}
\end{cases}
\end{align*}
\item
\begin{align*}
\Blo^{\rho,d} \sim_{\sM}
\begin{cases}
\Blo^{\varnothing,1} \swr \cT_d&\text{if }r-h(\rho)\text{ is even,}\\
(\Blo^{\varnothing,1} \swr \cT_d)\otimes \cC_1&\text{if }r-h(\rho)\text{ is odd.}
\end{cases}
\end{align*}
\end{enumerate}
\end{Proposition}

\begin{proof}
The first part is an immediate consequence of (\ref{algn:Rn_sM_ASd}) and Lemma \ref{lem:clif_mat1}.

For the second part, by Lemma \ref{lem:Yd_N} and (\ref{algn:Blo_R_sMor}),
\begin{align*}
{\mathfrak B}^{\rho,d}&\sim_{\operatorname{sMor}}
\begin{cases}
(\Blo^{\varnothing,1} \swr \cT_d) \otimes \cC_d & \text{if $r-h(\rho)+d$ is even},\\
(\Blo^{\varnothing,1} \swr \cT_d) \otimes \cC_{d+1} & \text{if $r-h(\rho)+d$ is odd}.
\end{cases}
\end{align*}
The claim now follows from Lemma \ref{lem:clif_mat1}.
\end{proof}

We now `desuperize' Proposition~\ref{P021221_3}(ii):

\begin{Corollary} \label{C021221_4}
We have 
\begin{align*}
\Blo^{\rho,d} \sim_{\Mor}
\begin{cases}
\Blo^{\varnothing,1} \swr \cT_d&\text{if }r-h(\rho)\text{ is even,}\\
(\Blo^{\varnothing,1} \swr \cT_d)_{\0}&\text{if }r-h(\rho)\text{ is odd,}
\end{cases}
\end{align*}
and
\begin{align*}
\Blo^{\rho,d}_\0 \sim_{\Mor}
\begin{cases}
(\Blo^{\varnothing,1} \swr \cT_d)_{\0}&\text{if }r-h(\rho)\text{ is even,}\\
\Blo^{\varnothing,1} \swr \cT_d&\text{if }r-h(\rho)\text{ is odd.}
\end{cases}
\end{align*}
\end{Corollary}
\begin{proof}
Follows from Proposition~\ref{P021221_3}(ii) using Lemmas~\ref{lem:Mor_A0} and~\ref{lem:clif_mat2}.
\end{proof}

\begin{Theorem}
Brou\'e's abelian defect group conjecture holds for the RoCK blocks 
$\Blo^{\rho,d}$ of 
$\F \tilde{\Si}_n$ and $\Blo^{\rho,d}_\0$ of $\F \tilde{\Ai}_n$.
\end{Theorem}

\begin{proof}
Follows from Corollary~\ref{C021221_4} and Proposition~\ref{lem:Brauer_der}(i).
\end{proof}

We now explain how Conjecture 2 from the Introduction implies Kessar-Schaps' conjecture from the introduction and Brou\'e's abelian defect group conjecture for double covers of symmetric and alternating groups.

\begin{Theorem}\label{TBroue} 
Suppose that Conjecture 2 from the Introduction holds. 
Let $\F$ be an algebraically closed field of positive characteristic. 
Then Brou\'e's abelian defect group conjecture for blocks of double covers of symmetric and alternating groups holds.
\end{Theorem}
\begin{proof}
We first deal with the $p=2$ case. Let $n \in \Z_{>0}$. Since $\langle z \rangle$ is a central $2$-subgroup of $\tilde{\Si}_n$, the blocks of $\F \tilde{\Si}_n$ are in one-to-one correspondence with those of $\F \Si_n$. Moreover, if $D$ is the defect group of a block of $\F \Si_n$, then $\pi_n^{-1}(D)$ is defect group of the corresponding block of $\F \tilde{\Si}_n$.

Now, by \cite[Theorem 6.2.45]{JK}, a defect group of a block of $\F \Si_n$ of `weight' $d$ is given by a Sylow $2$-subgroup of $\Si_V \leq \Si_n$ where $V := \{n-2d+1, \dots, n\}$. Therefore, a defect group of the corresponding block of $\F \tilde{\Si}_n$ is a Sylow $2$-subgroup of $\tilde{\Si}_V \leq \tilde{\Si}_n$. One can now quickly check that any abelian defect group of a block of $\F \tilde{\Si}_n$ must be isomorphic to either $C_2$, $C_4$ or $C_2 \times C_2$.

If the defect group is cyclic, then Brou\'e's conjecture is known to hold, as noted in the introduction, by \cite[Theorem 4.2]{Rickard_1}. If the defect group is isomorphic to $C_2 \times C_2$, then Brou\'e's conjecture is again known to hold, this time due to \cite{Lin}.

From now on we assume $p>2$. As explained in \S\ref{SSSpinBlocks}, it suffices to prove Brou\'e's abelian defect group conjecture for spin blocks of $\F\tilde\Si_n$ and $\F\tilde\Ai_n$, since for blocks of $\F\Si_n$ and $\F\Ai_n$ it has been proved in \cite{CR,MarAlt}. Moreover, the case of defect zero blocks is trivial, so we may assume that $0<d<p$, $\rho'$ is any $\bar p$-core and prove that  
$
{\mathfrak B}^{\rho',d}
\sim_{\operatorname{der}} {\mathfrak b}^{\rho',d}$ and 
$
{\mathfrak B}^{\rho',d}_\0
\sim_{\operatorname{der}} {\mathfrak b}^{\rho',d}_\0$. 

By Proposition~\ref{lem:Brauer_der}(i), we have 
$$
{\mathfrak b}^{\rho',d}
\sim_{\operatorname{der}}
\begin{cases}
{\mathfrak B}^{\varnothing,1}\wr_{\tt s} {\mathcal T}_d&\text{if }|\rho'|-h(\rho')\text{ is even,}\\
({\mathfrak B}^{\varnothing,1} \wr_{\tt s} {\mathcal T}_d)_{\0}&\text{if }|\rho'|-h(\rho')\text{ is odd,}
\end{cases}
$$
and
$$
{\mathfrak b}^{\rho',d}_\0
\sim_{\operatorname{der}}
\begin{cases}
({\mathfrak B}^{\varnothing,1} \wr_{\tt s} {\mathcal T}_d)_{\0}&\text{if }|\rho'|-h(\rho')\text{ is even,}\\
{\mathfrak B}^{\varnothing,1} \wr_{\tt s} {\mathcal T}_d&\text{if }|\rho'|-h(\rho')\text{ is odd.}
\end{cases}
$$
On the other hand, let $\rho$ be a $d$-Rouquier core with $|\rho|-h(\rho)\equiv|\rho'|-h(\rho')\pmod{2}$. (If $|\rho|-h(\rho)|$ has the wrong parity we can always change its parity by adding $1$ or $2$ beads to the $\ell^{\nth}$ runner.) By Corollary~\ref{C021221_4}, we have
$$
{\mathfrak B}^{\rho,d}
\sim_{\operatorname{Mor}}
\begin{cases}
{\mathfrak B}^{\varnothing,1}\wr_{\tt s} {\mathcal T}_d&\text{if }|\rho|-h(\rho)\text{ is even,}\\
({\mathfrak B}^{\varnothing,1} \wr_{\tt s} {\mathcal T}_d)_{\0}&\text{if }|\rho|-h(\rho)\text{ is odd,}
\end{cases}
$$
and
$$
{\mathfrak B}^{\rho,d}_\0
\sim_{\operatorname{Mor}}
\begin{cases}
({\mathfrak B}^{\varnothing,1} \wr_{\tt s} {\mathcal T}_d)_{\0}&\text{if }|\rho|-h(\rho)\text{ is even,}\\
{\mathfrak B}^{\varnothing,1} \wr_{\tt s} {\mathcal T}_d&\text{if }|\rho|-h(\rho)\text{ is odd.}
\end{cases}
$$

So we have
\begin{equation}\label{eqn:block_brauer_even}
{\mathfrak B}^{\rho,d}\sim_{\operatorname{der}}{\mathfrak b}^{\rho,d}\sim_{\operatorname{Mor}}{\mathfrak b}^{\rho',d}
\end{equation}
and
\begin{equation}\label{eqn:block_brauer_odd}
{\mathfrak B}^{\rho,d}_{\0}\sim_{\operatorname{der}}{\mathfrak b}^{\rho,d}_{\0}\sim_{\operatorname{Mor}}{\mathfrak b}^{\rho',d}_{\0},
\end{equation}
where the Morita equivalences between Brauer correspondents follow from (\ref{mor_even_brauer}) and (\ref{mor_odd_brauer}) in the proof of Proposition \ref{lem:Brauer_der}.

We now assume that Conjecture 2 from the Introduction holds. Let's first assume $|\rho'|-h(\rho')$ is even. Using (\ref{algn:Blo_R_sMor}) we get
\begin{align*}
{\mathfrak B}^{\rho',d} \sim_{\operatorname{sMor}}R^{\Lambda_0}_{\operatorname{cont}(\rho')+d\delta} \sim_{\operatorname{der}}R^{\Lambda_0}_{\operatorname{cont}(\rho)+d\delta} \sim_{\operatorname{sMor}}{\mathfrak B}^{\rho,d} \sim_{\operatorname{der}}{\mathfrak b}^{\rho',d},
\end{align*}
as required. Similarly,
\begin{align*}
{\mathfrak B}^{\rho',d}_{\0} & \sim_{\operatorname{Mor}} {\mathfrak B}^{\rho',d} \otimes \cC_1 \sim_{\operatorname{sMor}}R^{\Lambda_0}_{\operatorname{cont}(\rho')+d\delta} \otimes \cC_1 \sim_{\operatorname{der}}R^{\Lambda_0}_{\operatorname{cont}(\rho)+d\delta} \otimes \cC_1\\
&\sim_{\operatorname{sMor}}{\mathfrak B}^{\rho,d} \otimes \cC_1 \sim_{\operatorname{Mor}} {\mathfrak B}^{\rho,d}_{\0} \sim_{\operatorname{der}}{\mathfrak b}^{\rho',d}_{\0},
\end{align*}
where the first and fifth equivalences follows from Lemma \ref{lem:clif_mat2}.

Let's now assume $|\rho'|-h(\rho')$ is odd. We now have
\begin{align*}
{\mathfrak B}^{\rho',d} \sim_{\operatorname{sMor}}R^{\Lambda_0}_{\operatorname{cont}(\rho')+d\delta} \otimes \cC_1 \sim_{\operatorname{der}}R^{\Lambda_0}_{\operatorname{cont}(\rho)+d\delta} \otimes \cC_1 \sim_{\operatorname{sMor}}{\mathfrak B}^{\rho,d} \sim_{\operatorname{der}}{\mathfrak b}^{\rho',d},
\end{align*}
as required. Similarly,
\begin{align*}
{\mathfrak B}^{\rho',d}_{\0} & \sim_{\operatorname{Mor}} {\mathfrak B}^{\rho',d} \otimes \cC_1 \sim_{\operatorname{sMor}}R^{\Lambda_0}_{\operatorname{cont}(\rho')+d\delta} \sim_{\operatorname{der}}R^{\Lambda_0}_{\operatorname{cont}(\rho)+d\delta} \\ & \sim_{\operatorname{sMor}}{\mathfrak B}^{\rho,d} \otimes \cC_1 \sim_{\operatorname{Mor}} {\mathfrak B}^{\rho,d}_{\0} \sim_{\operatorname{der}}{\mathfrak b}^{\rho',d}_{\0},
\end{align*}
where, again, the first and fifth equivalences follows from Lemma \ref{lem:clif_mat2}.
\end{proof}

For the following proof, $\rho$ and $\rho'$ will just denote arbitrary ${\bar p}$-cores, as in the statement of Kessar-Schaps' conjecture.

\begin{Theorem}\label{TKS} 
Suppose that Conjecture 2 from the Introduction holds. Then Kessar-Schaps' conjecture from the Introduction holds. 
\end{Theorem}
\begin{proof}
We assume throughout the proof that Conjecture 2 holds. Let $d \in \Z_{>0}$. If $|\rho| - h(\rho)$ is even and $|\rho'| - h(\rho')$ is odd, it follows from (\ref{algn:Blo_R_sMor}) that
\begin{align*}
\Blo^{\rho,d} \sim_{\sM} R^{\Lambda_0}_{\cont(\rho)+d\delta} \sim_{\der} R^{\Lambda_0}_{\cont(\rho')+d\delta} \sim_{\sM} \Blo^{\rho',d} \otimes \cC_1 \sim_{\Mor} \Blo^{\rho',d}_{\0},
\end{align*}
where the final equivalence follows from Lemma \ref{lem:clif_mat2}. Similarly, if $|\rho| - h(\rho)$ is odd and $|\rho'| - h(\rho')$ is even,
\begin{align*}
\Blo^{\rho,d} \sim_{\sM} R^{\Lambda_0}_{\cont(\rho)+d\delta} \otimes \cC_1 \sim_{\der} R^{\Lambda_0}_{\cont(\rho')+d\delta} \otimes \cC_1 \sim_{\sM} \Blo^{\rho',d} \otimes \cC_1 \sim_{\Mor} \Blo^{\rho',d}_{\0}.
\end{align*}
The other `if' statements of the Kessar-Schaps conjecture are proved very similarly.

Now let $d' \in \Z_{>0}$ and suppose $\Blo^{\rho,d}$ is derived equivalent to $\Blo^{\rho',d'}_{\0}$. Recall the elementary divisors of the Cartan matrix of a block are invariant under derived equivalence. (See \cite[Proposition 5.1]{Xi} and its proof.) Since the largest elementary divisor of the Cartan matrix coincides with the order of the defect group, $\Blo^{\rho,d}$ and $\Blo^{\rho',d'}$ must have defect groups of the same order. It now follows from Theorems \ref{thm:Sn_blocks} and \ref{thm:An_blocks} that $d=d'$.

Since $\Blo^{\rho,d}$ and $\Blo^{\rho',d'}_{\0}$ are derived equivalent they must have the same number of isomorphism classes of irreducible modules. It now follows from \cite[Proposition 13.17]{Ols} that $|\rho| - h(\rho)$ and $|\rho'| - h(\rho')$ must have opposite parity, as required. Again, the other `only if' statements of the Kessar-Schaps conjecture are proved very similarly. 
\end{proof}

\begin{Remark}
We note that we have not actually used Conjecture 2 to prove the `only if' direction of Kessar-Schaps' conjecture. Consequently this direction holds with no extra assumptions.
\end{Remark}

\subsection{Morita equivalences for more general algebras}
In this subsection, we do not assume that $d<p$ and consider a $d$-Rouquier core $\rho$. We write $\cont(\rho)+d\de=\sum_{i\in I}m_i\al_i$. We also assume that we are in the setting of Theorems \ref{TKKTQ2} or \ref{TKKTQ1}; in particular, 
either $q=1$ and $\cha\F=p=2\ell-1$ or $q$ is a primitive $(2\ell+1)$st toot of unity in a field $\F$ of characteristic different from $2$. 

\begin{Theorem} \label{T231221}
There exists an idempotent $f\in\cX^{\La_0}_{\cont(\rho)+d\de}(q)$ such that 
$$
f\cX^{\La_0}_{\cont(\rho)+d\de}(q)f\cong (\Zig_\ell\swr\Si_d)\otimes \cC_{m_0}.
$$
Moreover, $\cX^{\La_0}_{\cont(\rho)+d\de}(q)\sim_{\sM}f\cX^{\La_0}_{\cont(\rho)+d\de}(q)f$ if and only if $d<p$. 
\end{Theorem}
\begin{proof}
By Lemma~\ref{LIdSym}, there exists an idempotent $\eps_d\in R^{\La_0}_{\cont(\rho)+d\de}$ such that 
$$\eps_d R^{\La_0}_{\cont(\rho)+d\de}\eps_d\cong Z_{\rho,d}.$$ Take $$e:=\ga_{1^d}\eps_d\in R^{\La_0}_{\cont(\rho)+d\de}.$$  It is clear from Lemma~\ref{LIdSymNew} that $\eps_d$ and $\ga_{1^d}$ commute, so $e$ is an idempotent and 
$$eR^{\La_0}_{\cont(\rho)+d\de}e\cong Y_{\rho,d},$$ see (\ref{ECTrBold}). Moreover, by Corollary~\ref{CMorWr}, we have $R^{\La_0}_{\cont(\rho)+d\de}\sim_{\sM} Y_{\rho,d}$ if and only if $d<p$. By Corollary~\ref{Cor:QF2} and Theorem~\ref{thm:main}, we now deduce that $$eR^{\La_0}_{\cont(\rho)+d\de}e\cong \Zig_\ell\swr\Si_d$$ and 
$$(e\otimes 1) (R^{\La_0}_{\cont(\rho)+d\de}\otimes\cC_{m_0})(e\otimes 1)\cong (\Zig_\ell\swr\Si_d)\otimes \cC_{m_0}.$$
Let $f\in \cX^{\La_0}_{\cont(\rho)+d\de}$ be the image of $e\otimes 1$ under the maps
$$
R^{\La_0}_{\cont(\rho)+d\de}\otimes\cC_{m_0}\iso e_IRC^{\La_0}_{\cont(\rho)+d\de}e_I\into RC^{\La_0}_{\cont(\rho)+d\de}\iso \cX^{\La_0}_{\cont(\rho)+d\de}(q),
$$
where the first isomorphism comes from Proposition~\ref{prop:RCtoR} and the second isomorphism comes from Theorems~\ref{TKKTQ2}, \ref{TKKTQ1}. 
Then 
$$
f\cX^{\La_0}_{\cont(\rho)+d\de}(q)f\cong (\Zig_\ell\swr\Si_d)\otimes \cC_{m_0}.
$$ 
Moreover, taking into account Lemma~\ref{lem:RC_trunc}, we deduce that $\cX^{\La_0}_{\cont(\rho)+d\de}(q)\sim_{\sM}f\cX^{\La_0}_{\cont(\rho)+d\de}(q)f$ if and only if $d<p$. 
\end{proof}


\chapter{Appendix. Some calculations in \texorpdfstring{$B_2$}{}}

\section{Some small rank computations}
The lemmas in this subsection follow by applying braid relations (\ref{R7}), dot-crossing relations (\ref{R5}), and quadratic relations (\ref{R6}) (which in particular implies  
$
	  \begin{braid}\tikzset{baseline=3em}
	\braidbox{0}{1.3}{3}{4}{};
	\draw (0.7,3.5) node{$j\ j$};
	\draw(0.3,4)--(1.1,4.7);
	\draw(1.1,4)--(0.3,4.7);
	  \end{braid}=0$).

\vspace{0.5mm}\begin{Lemma} \label{LA3/4}
Let $1\leq i\leq \ell-2$. Then in $R_{2\al_{i}+2\al_{i+1}}$ we have 
$$
\resizebox{126mm}{8mm}{
\begin{braid}\tikzset{baseline=1.6em}
	\draw (0,0) node{\normalsize$i$};
	\braidbox{1}{4.2}{-0.8}{0.7}{$(i+1)^2$};
	\draw (5.2,0) node{\normalsize$i$};
	\draw (0,4) node{\normalsize$i$};
	\braidbox{1}{4.2}{3.2}{4.7}{$(i+1)^2$};
	\draw (5.2,4) node{\normalsize$i$};
	\draw(0,3)--(5.2,1);
	 \draw(5.2,3)--(0,1);
	 \draw(1.5,3)--(0,2)--(1.5,1);
	 \draw(3.3,3)--(1.8,2)--(3.3,1);
	  \end{braid}
	  =
\begin{braid}\tikzset{baseline=1.6em}
	\draw (0,0) node{\normalsize$i$};
	\braidbox{1}{4.2}{-0.8}{0.7}{$(i+1)^2$};
	\draw (5.2,0) node{\normalsize$i$};
	\draw (0,4) node{\normalsize$i$};
	\braidbox{1}{4.2}{3.2}{4.7}{$(i+1)^2$};
	\draw (5.2,4) node{\normalsize$i$};
	\draw(0,3)--(5.2,1);
	 \draw(5.2,3)--(0,1);
	 \draw(1.5,3)--(3.5,2)--(1.5,1);
	 \draw(3.3,3)--(4.8,2)--(3.3,1);
	  \end{braid}
	  --
	  \begin{braid}\tikzset{baseline=1.6em}
	\draw (0,0) node{\normalsize$i$};
	\braidbox{1}{4.2}{-0.8}{0.7}{$(i+1)^2$};
	\draw (5.2,0) node{\normalsize$i$};
	\draw (0,4) node{\normalsize$i$};
	\braidbox{1}{4.2}{3.2}{4.7}{$(i+1)^2$};
	\draw (5.2,4) node{\normalsize$i$};
	\draw(0,3)--(0,1);
	 \draw(5.2,3)--(5.2,1);
	 \draw(1.5,3)--(1.5,1);
	 \draw(3.3,3)--(3.3,1);
	 \blackdot (0,2);
	  \end{braid}
	  +
	  \begin{braid}\tikzset{baseline=1.6em}
	\draw (0,0) node{\normalsize$i$};
	\braidbox{1}{4.2}{-0.8}{0.7}{$(i+1)^2$};
	\draw (5.2,0) node{\normalsize$i$};
	\draw (0,4) node{\normalsize$i$};
	\braidbox{1}{4.2}{3.2}{4.7}{$(i+1)^2$};
	\draw (5.2,4) node{\normalsize$i$};
	\draw(0,3)--(0,1);
	 \draw(5.2,3)--(5.2,1);
	 \draw(1.5,3)--(1.5,1);
	 \draw(3.3,3)--(3.3,1);
	 \blackdot (1.5,2);
	  \end{braid}
	  +
	  \begin{braid}\tikzset{baseline=1.6em}
	\draw (0,0) node{\normalsize$i$};
	\braidbox{1}{4.2}{-0.8}{0.7}{$(i+1)^2$};
	\draw (5.2,0) node{\normalsize$i$};
	\draw (0,4) node{\normalsize$i$};
	\braidbox{1}{4.2}{3.2}{4.7}{$(i+1)^2$};
	\draw (5.2,4) node{\normalsize$i$};
	\draw(0,3)--(0,1);
	 \draw(5.2,3)--(5.2,1);
	 \draw(1.5,3)--(1.5,1);
	 \draw(3.3,3)--(3.3,1);
	 \blackdot (3.3,2);
	  \end{braid}
	  --
	  \begin{braid}\tikzset{baseline=1.6em}
	\draw (0,0) node{\normalsize$i$};
	\braidbox{1}{4.2}{-0.8}{0.7}{$(i+1)^2$};
	\draw (5.2,0) node{\normalsize$i$};
	\draw (0,4) node{\normalsize$i$};
	\braidbox{1}{4.2}{3.2}{4.7}{$(i+1)^2$};
	\draw (5.2,4) node{\normalsize$i$};
	\draw(0,3)--(0,1);
	 \draw(5.2,3)--(5.2,1);
	 \draw(1.5,3)--(1.5,1);
	 \draw(3.3,3)--(3.3,1);
	 \blackdot (5.2,2);
	  \end{braid}.
	  }
$$
\end{Lemma}

\vspace{0.5mm}\begin{Lemma} \label{LA4/5}
Let $\ell>1$. Then in $R_{2\al_{0}+2\al_{1}}$ we have 
$$
\resizebox{126mm}{8mm}{
\begin{braid}\tikzset{baseline=1.6em}
	\draw (0,0) node{\normalsize$1$};
	\redbraidbox{1}{4.2}{-0.8}{0.7}{$\color{red}0\,\,0\color{black}$};
	\draw (5.2,0) node{\normalsize$1$};
	\draw (0,4) node{\normalsize$1$};
	\redbraidbox{1}{4.2}{3.2}{4.7}{$\color{red}0\,\,0\color{black}$};
	\draw (5.2,4) node{\normalsize$1$};
	\draw(0,3)--(5.2,1);
	 \draw(5.2,3)--(0,1);
	 \draw[red](1.5,3)--(0,2)--(1.5,1);
	 \draw[red](3.3,3)--(1.8,2)--(3.3,1);
	  \end{braid}
	  =
\begin{braid}\tikzset{baseline=1.6em}
	\draw (0,0) node{\normalsize$1$};
	\redbraidbox{1}{4.2}{-0.8}{0.7}{$\color{red}0\,\,0\color{black}$};
	\draw (5.2,0) node{\normalsize$1$};
	\draw (0,4) node{\normalsize$1$};
	\redbraidbox{1}{4.2}{3.2}{4.7}{$\color{red}0\,\,0\color{black}$};
	\draw (5.2,4) node{\normalsize$1$};
	\draw(0,3)--(5.2,1);
	 \draw(5.2,3)--(0,1);
	 \draw[red](1.5,3)--(3.5,2)--(1.5,1);
	 \draw[red](3.3,3)--(4.8,2)--(3.3,1);
	  \end{braid}
	  --
	  \begin{braid}\tikzset{baseline=1.6em}
	\draw (0,0) node{\normalsize$1$};
	\redbraidbox{1}{4.2}{-0.8}{0.7}{$\color{red}0\,\,0\color{black}$};
	\draw (5.2,0) node{\normalsize$1$};
	\draw (0,4) node{\normalsize$1$};
	\redbraidbox{1}{4.2}{3.2}{4.7}{$\color{red}0\,\,0\color{black}$};
	\draw (5.2,4) node{\normalsize$1$};
	\draw(0,3)--(0,1);
	 \draw(5.2,3)--(5.2,1);
	 \draw[red](1.5,3)--(1.5,1);
	 \draw[red](3.3,3)--(3.3,1);
	 \blackdot (0,2);
	  \end{braid}
	  +
	  \begin{braid}\tikzset{baseline=1.6em}
	\draw (0,0) node{\normalsize$1$};
	\redbraidbox{1}{4.2}{-0.8}{0.7}{$\color{red}0\,\,0\color{black}$};
	\draw (5.2,0) node{\normalsize$1$};
	\draw (0,4) node{\normalsize$1$};
	\redbraidbox{1}{4.2}{3.2}{4.7}{$\color{red}0\,\,0\color{black}$};
	\draw (5.2,4) node{\normalsize$1$};
	\draw(0,3)--(0,1);
	 \draw(5.2,3)--(5.2,1);
	 \draw[red](1.5,3)--(1.5,1);
	 \draw[red](3.3,3)--(3.3,1);
	 \reddot (1.5,2.2);
	 \reddot (1.5,1.8);
	  \end{braid}
	  +
	  \begin{braid}\tikzset{baseline=1.6em}
	\draw (0,0) node{\normalsize$1$};
	\redbraidbox{1}{4.2}{-0.8}{0.7}{$\color{red}0\,\,0\color{black}$};
	\draw (5.2,0) node{\normalsize$1$};
	\draw (0,4) node{\normalsize$1$};
	\redbraidbox{1}{4.2}{3.2}{4.7}{$\color{red}0\,\,0\color{black}$};
	\draw (5.2,4) node{\normalsize$1$};
	\draw(0,3)--(0,1);
	 \draw(5.2,3)--(5.2,1);
	 \draw[red](1.5,3)--(1.5,1);
	 \draw[red](3.3,3)--(3.3,1);
	 \reddot (3.3,2.2);
	 \reddot (3.3,1.8);
	  \end{braid}
	  --
	  \begin{braid}\tikzset{baseline=1.6em}
	\draw (0,0) node{\normalsize$1$};
	\redbraidbox{1}{4.2}{-0.8}{0.7}{$\color{red}0\,\,0\color{black}$};
	\draw (5.2,0) node{\normalsize$1$};
	\draw (0,4) node{\normalsize$1$};
	\redbraidbox{1}{4.2}{3.2}{4.7}{$\color{red}0\,\,0\color{black}$};
	\draw (5.2,4) node{\normalsize$1$};
	\draw(0,3)--(0,1);
	 \draw(5.2,3)--(5.2,1);
	 \draw[red](1.5,3)--(1.5,1);
	 \draw[red](3.3,3)--(3.3,1);
	 \blackdot (5.2,2);
	  \end{braid}.
	  }
$$
\end{Lemma}

\vspace{0.5mm}\begin{Lemma} \label{LA} 
Let $\ell>1$. Then in $R_{4\al_{\ell-1}+\al_\ell}$ we have 
$$
\resizebox{126mm}{8mm}{
\begin{braid}\tikzset{baseline=1.6em}
	\braidbox{0}{3.2}{-0.8}{0.7}{$(\ell-1)^2$};
	\draw (4.1,0) node{\color{blue}\normalsize$\ell$\color{black}};
	\braidbox{5}{8.1}{-0.8}{0.7}{$(\ell-1)^2$};
	\braidbox{0}{3.2}{3.2}{4.7}{$(\ell-1)^2$};
	\draw (4.1,4) node{\color{blue}\normalsize$\ell$\color{black}};
	\braidbox{5}{8.1}{3.2}{4.7}{$(\ell-1)^2$};
         \draw(0.5,3)--(5.5,1);
         \draw(2.5,3)--(7.5,1);
	 \draw(5.2,3)--(0.5,1);
	 \draw(2.5,1)--(7.5,3);
	 \draw[blue](4.1,3)--(1,2)--(4.1,1);
	  \end{braid}
	  =
\begin{braid}\tikzset{baseline=1.6em}
	\braidbox{0}{3.2}{-0.8}{0.7}{$(\ell-1)^2$};
	\draw (4.1,0) node{\color{blue}\normalsize$\ell$\color{black}};
	\braidbox{5}{8.1}{-0.8}{0.7}{$(\ell-1)^2$};
	\braidbox{0}{3.2}{3.2}{4.7}{$(\ell-1)^2$};
	\draw (4.1,4) node{\color{blue}\normalsize$\ell$\color{black}};
	\braidbox{5}{8.1}{3.2}{4.7}{$(\ell-1)^2$};
         \draw(0.5,3)--(5.5,1);
         \draw(2.5,3)--(7.5,1);
	 \draw(5.2,3)--(0.5,1);
	 \draw(2.5,1)--(7.5,3);
	 \draw[blue](4.1,3)--(2.5,2.5)--(5.4,1.5)--(4.1,1);
	  \end{braid}
	  +\begin{braid}\tikzset{baseline=1.6em}
	\braidbox{0}{3.2}{-0.8}{0.7}{$(\ell-1)^2$};
	\draw (4.1,0) node{\color{blue}\normalsize$\ell$\color{black}};
	\braidbox{5}{8.1}{-0.8}{0.7}{$(\ell-1)^2$};
	\braidbox{0}{3.2}{3.2}{4.7}{$(\ell-1)^2$};
	\draw (4.1,4) node{\color{blue}\normalsize$\ell$\color{black}};
	\braidbox{5}{8.1}{3.2}{4.7}{$(\ell-1)^2$};
         \draw(0.5,3)--(0.5,1);
         \draw(2.5,3)--(5.2,1);
	 \draw(5.2,3)--(7.5,1);
	 \draw(2.5,1)--(7.5,3);
	 \draw[blue](4.1,3)--(5.6,1.5)--(4.1,1);
	  \end{braid}
	  +\begin{braid}\tikzset{baseline=1.6em}
	\braidbox{0}{3.2}{-0.8}{0.7}{$(\ell-1)^2$};
	\draw (4.1,0) node{\color{blue}\normalsize$\ell$\color{black}};
	\braidbox{5}{8.1}{-0.8}{0.7}{$(\ell-1)^2$};
	\braidbox{0}{3.2}{3.2}{4.7}{$(\ell-1)^2$};
	\draw (4.1,4) node{\color{blue}\normalsize$\ell$\color{black}};
	\braidbox{5}{8.1}{3.2}{4.7}{$(\ell-1)^2$};
         \draw(0.5,3)--(0.5,1);
         \draw(2.5,3)--(2.5,1);
	 \draw(5.2,3)--(5.2,1);
	 \draw(7.5,1)--(7.5,3);
	 \draw[blue](4.1,3)--(4.1,1);
	 \end{braid}.
	 }
	  $$
\end{Lemma}

\vspace{0.5mm}\begin{Lemma} \label{L090721_3} 
Let $1\leq i<\ell$. Then 
In $R_{\al_{i}+4\al_{i+1}}$ we have 
\begin{align*}
\resizebox{95mm}{5.5mm}{
\begin{braid}\tikzset{baseline=1em}
	\braidbox{-.5}{2.5}{-0.8}{0.7}{$(i+1)^2$};
	\draw (3.3,0) node{\normalsize$i$};
	\braidbox{4.1}{7.1}{-0.8}{0.7}{$(i+1)^2$};
	\draw(3.3,3)--(5.5,2)--(3.3,1);
	 \draw(0.2,3)--(4.6,1);
	 \draw(1.8,3)--(6.2,1);
	  \draw(0.2,1)--(4.6,3);
	 \draw(1.8,1)--(6.2,3);
	  \end{braid}
	  =
\begin{braid}\tikzset{baseline=1em}
	\braidbox{-.5}{2.5}{-0.8}{0.7}{$(i+1)^2$};
	\draw (3.3,0) node{\normalsize$i$};
	\braidbox{4.1}{7.1}{-0.8}{0.7}{$(i+1)^2$};
	\draw(3.3,3)--(1.5,2)--(3.3,1);
	 \draw(0.2,3)--(4.6,1);
	 \draw(1.8,3)--(6.2,1);
	  \draw(0.2,1)--(4.6,3);
	 \draw(1.8,1)--(6.2,3);
	  \end{braid}
	  --
	 \begin{braid}\tikzset{baseline=1em}
	\braidbox{-.5}{2.5}{-0.8}{0.7}{$(i+1)^2$};
	\draw (3.3,0) node{\normalsize$i$};
	\braidbox{4.1}{7.1}{-0.8}{0.7}{$(i+1)^2$};
	\draw(3.3,3)--(2.1,1.6)--(3.3,1);
	 \draw(0.2,3)--(4.6,1);
	 \draw(4.6,3)--(6.2,1);
	  \draw(0.2,1)--(1.8,3);
	 \draw(1.8,1)--(6.2,3);
	  \end{braid}.
	  }
\end{align*}
\end{Lemma}

\vspace{0.5mm}\begin{Lemma} \label{LA1/2}
Let $1\leq i\leq \ell-2$. Then in $R_{4\al_{i}+2\al_{i+1}}$ we have 
\begin{align*}
\resizebox{126mm}{12mm}{
\begin{braid}\tikzset{baseline=3em}
	\braidbox{0}{2.5}{-0.8}{0.7}{$i\ \ i$};
	\braidbox{3.3}{7}{-0.8}{0.7}{$(i+1)^2$};
	\braidbox{7.8}{10.3}{-0.8}{0.7}{$i\ \ i$};
	\braidbox{0}{2.5}{7.2}{8.7}{$i\ \ i$};
	\braidbox{3.3}{7}{7.2}{8.7}{$(i+1)^2$};
	\braidbox{7.8}{10.3}{7.2}{8.7}{$i\ \ i$};
	\draw(0.4,7)--(8.3,1);
	 \draw(2,7)--(9.8,1);
	 \draw(0.4,1)--(8.3,7);
	 \draw(2,1)--(9.8,7);
	 \draw(4.2,7)--(1,4)--(4.2,1);
	 \draw(6,7)--(2.8,4)--(6,1);
	  \end{braid}
	  =
	  \begin{braid}\tikzset{baseline=3em}
	\braidbox{0}{2.5}{-0.8}{0.7}{$i\ \ i$};
	\braidbox{3.3}{7}{-0.8}{0.7}{$(i+1)^2$};
	\braidbox{7.8}{10.3}{-0.8}{0.7}{$i\ \ i$};
	\braidbox{0}{2.5}{7.2}{8.7}{$i\ \ i$};
	\braidbox{3.3}{7}{7.2}{8.7}{$(i+1)^2$};
	\braidbox{7.8}{10.3}{7.2}{8.7}{$i\ \ i$};
	\draw(0.4,7)--(8.3,1);
	 \draw(2,7)--(9.8,1);
	 \draw(0.4,1)--(8.3,7);
	 \draw(2,1)--(9.8,7);
	 \draw(4.2,7)--(3,4)--(4.2,1);
	 \draw(6,7)--(7,4)--(6,1);
	  \end{braid}
	  --
	  \begin{braid}\tikzset{baseline=3em}
	\braidbox{0}{2.5}{-0.8}{0.7}{$i\ \ i$};
	\braidbox{3.3}{7}{-0.8}{0.7}{$(i+1)^2$};
	\braidbox{7.8}{10.3}{-0.8}{0.7}{$i\ \ i$};
	\braidbox{0}{2.5}{7.2}{8.7}{$i\ \ i$};
	\braidbox{3.3}{7}{7.2}{8.7}{$(i+1)^2$};
	\braidbox{7.8}{10.3}{7.2}{8.7}{$i\ \ i$};
	\draw(0.4,7)--(8.3,1);
	 \draw(2,7)--(0.4,1);
	 \draw(9.8,1)--(8.3,7);
	 \draw(2,1)--(9.8,7);
	 \draw(4.2,7)--(6.2,3.3)--(4.2,1);
	 \draw(6,7)--(8,3.3)--(6,1);
	 \blackdot (0.5,1.4);
	  \end{braid}
 +
	  \begin{braid}\tikzset{baseline=3em}
	\braidbox{0}{2.5}{-0.8}{0.7}{$i\ \ i$};
	\braidbox{3.3}{7}{-0.8}{0.7}{$(i+1)^2$};
	\braidbox{7.8}{10.3}{-0.8}{0.7}{$i\ \ i$};
	\braidbox{0}{2.5}{7.2}{8.7}{$i\ \ i$};
	\braidbox{3.3}{7}{7.2}{8.7}{$(i+1)^2$};
	\braidbox{7.8}{10.3}{7.2}{8.7}{$i\ \ i$};
	\draw(0.4,7)--(8.3,1);
	 \draw(2,7)--(0.4,1);
	 \draw(9.8,1)--(8.3,7);
	 \draw(2,1)--(9.8,7);
	 \draw(4.2,7)--(6.2,3.3)--(4.2,1);
	 \draw(6,7)--(8,3.3)--(6,1);
	 \blackdot (4.55,1.4);
	  \end{braid}
	  +
	  \begin{braid}\tikzset{baseline=3em}
	\braidbox{0}{2.5}{-0.8}{0.7}{$i\ \ i$};
	\braidbox{3.3}{7}{-0.8}{0.7}{$(i+1)^2$};
	\braidbox{7.8}{10.3}{-0.8}{0.7}{$i\ \ i$};
	\braidbox{0}{2.5}{7.2}{8.7}{$i\ \ i$};
	\braidbox{3.3}{7}{7.2}{8.7}{$(i+1)^2$};
	\braidbox{7.8}{10.3}{7.2}{8.7}{$i\ \ i$};
	\draw(0.4,7)--(0.4,1);
	 \draw(2,7)--(2,1);
	 \draw(8.3,1)--(8.3,7);
	 \draw(9.8,1)--(9.8,7);
	 \draw(4.2,7)--(4.2,1);
	 \draw(6,7)--(6,1);
	  \end{braid}.
	  }
\end{align*}
\end{Lemma}

\vspace{1mm}
In the following lemma and below, as usual, to have certainty with signs, we interpret the ``double odd crossing" as in (\ref{EDoubleRed}).

\begin{Lemma} \label{LB} 
Let $\ell>1$. Then in $R_{4\al_{0}+\al_{1}}$ we have 
$$
\resizebox{126mm}{8mm}{
\begin{braid}\tikzset{baseline=1.6em}
	\redbraidbox{0}{2.1}{-0.8}{0.7}{$\color{red}0\,\,0\color{black}$};
	\draw (3.1,0) node{\normalsize$1$};
	\redbraidbox{4.1}{6.2}{-0.8}{0.7}{$\color{red}0\,\,0\color{black}$};
	\redbraidbox{0}{2.1}{3.2}{4.7}{$\color{red}0\,\,0\color{black}$};
	\draw (3.1,4) node{\normalsize$1$};
	\redbraidbox{4.1}{6.2}{3.2}{4.7}{$\color{red}0\,\,0\color{black}$};
	\draw(3.1,3)--(1,2)--(3.1,1);
	 \draw[red](0.4,3)--(4.4,1);
	 \draw[red](1.8,3)--(5.8,1);
	  \draw[red](0.4,1)--(4.4,3);
	 \draw[red](1.8,1)--(5.8,3);
	  \end{braid}
	  =
\begin{braid}\tikzset{baseline=1.6em}
	\redbraidbox{0}{2.1}{-0.8}{0.7}{$\color{red}0\,\,0\color{black}$};
	\draw (3.1,0) node{\normalsize$1$};
	\redbraidbox{4.1}{6.2}{-0.8}{0.7}{$\color{red}0\,\,0\color{black}$};
	\redbraidbox{0}{2.1}{3.2}{4.7}{$\color{red}0\,\,0\color{black}$};
	\draw (3.1,4) node{\normalsize$1$};
	\redbraidbox{4.1}{6.2}{3.2}{4.7}{$\color{red}0\,\,0\color{black}$};
	\draw(3.1,3)--(5.4,2)--(3.1,1);
	 \draw[red](0.4,3)--(4.4,1);
	 \draw[red](1.8,3)--(5.8,1);
	  \draw[red](0.4,1)--(4.4,3);
	 \draw[red](1.8,1)--(5.8,3);
	  \end{braid}
	  +
	  \begin{braid}\tikzset{baseline=1.6em}
	\redbraidbox{0}{2.1}{-0.8}{0.7}{$\color{red}0\,\,0\color{black}$};
	\draw (3.1,0) node{\normalsize$1$};
	\redbraidbox{4.1}{6.2}{-0.8}{0.7}{$\color{red}0\,\,0\color{black}$};
	\redbraidbox{0}{2.1}{3.2}{4.7}{$\color{red}0\,\,0\color{black}$};
	\draw (3.1,4) node{\normalsize$1$};
	\redbraidbox{4.1}{6.2}{3.2}{4.7}{$\color{red}0\,\,0\color{black}$};
	\draw(3.1,3)--(4,2)--(3.1,1);
	 \draw[red](0.4,3)--(4.4,1);
	 \draw[red](1.8,3)--(0.4,1);
	  \draw[red](1.8,1)--(4.4,3);
	 \draw[red](5.8,1)--(5.8,3);
	  \end{braid}
	  -- \begin{braid}\tikzset{baseline=1.6em}
	\redbraidbox{0}{2.1}{-0.8}{0.7}{$\color{red}0\,\,0\color{black}$};
	\draw (3.1,0) node{\normalsize$1$};
	\redbraidbox{4.1}{6.2}{-0.8}{0.7}{$\color{red}0\,\,0\color{black}$};
	\redbraidbox{0}{2.1}{3.2}{4.7}{$\color{red}0\,\,0\color{black}$};
	\draw (3.1,4) node{\normalsize$1$};
	\redbraidbox{4.1}{6.2}{3.2}{4.7}{$\color{red}0\,\,0\color{black}$};
	\draw(3.1,3)--(4.2,1.6)--(3.1,1);
	 \draw[red](0.4,3)--(1,2)--(4.4,1);
	 \draw[red](1.8,3)--(0.4,1);
	  \draw[red](1.8,1)--(5.8,3);
	 \draw[red](5.8,1)--(4.4,3);
	 \reddot (0.5,1.2);
	  \end{braid}
	  --
	  \begin{braid}\tikzset{baseline=1.6em}
	\redbraidbox{0}{2.1}{-0.8}{0.7}{$\color{red}0\,\,0\color{black}$};
	\draw (3.1,0) node{\normalsize$1$};
	\redbraidbox{4.1}{6.2}{-0.8}{0.7}{$\color{red}0\,\,0\color{black}$};
	\redbraidbox{0}{2.1}{3.2}{4.7}{$\color{red}0\,\,0\color{black}$};
	\draw (3.1,4) node{\normalsize$1$};
	\redbraidbox{4.1}{6.2}{3.2}{4.7}{$\color{red}0\,\,0\color{black}$};
	\draw(3.1,3)--(3.1,1);
	 \draw[red](0.4,3)--(0.4,1);
	 \draw[red](1.8,3)--(1.8,1);
	  \draw[red](4.4,1)--(4.4,3);
	 \draw[red](5.8,1)--(5.8,3);
	  \end{braid}.
	  }
$$
\end{Lemma}

\vspace{1mm}
\begin{Lemma} \label{LBSpecial} 
Let $\ell=1$. Then in $R_{4\al_{0}+\al_{1}}$ we have 
$$
\resizebox{120mm}{8mm}{
\begin{braid}\tikzset{baseline=1.6em}
	\redbraidbox{0}{2.1}{-0.8}{0.7}{$\color{red}0\,\,0\color{black}$};
	\draw (3.1,0) node{\color{blue}\normalsize$1$\color{black}};
	\redbraidbox{4.1}{6.2}{-0.8}{0.7}{$\color{red}0\,\,0\color{black}$};
	\redbraidbox{0}{2.1}{3.2}{4.7}{$\color{red}0\,\,0\color{black}$};
	\draw (3.1,4) node{\color{blue}\normalsize$1$\color{black}};
	\redbraidbox{4.1}{6.2}{3.2}{4.7}{$\color{red}0\,\,0\color{black}$};
	\draw[blue](3.1,3)--(1,2)--(3.1,1);
	 \draw[red](0.4,3)--(4.4,1);
	 \draw[red](1.8,3)--(5.8,1);
	  \draw[red](0.4,1)--(4.4,3);
	 \draw[red](1.8,1)--(5.8,3);
	  \end{braid}
	  =
\begin{braid}\tikzset{baseline=1.6em}
	\redbraidbox{0}{2.1}{-0.8}{0.7}{$\color{red}0\,\,0\color{black}$};
	\draw (3.1,0) node{\color{blue}\normalsize$1$\color{black}};
	\redbraidbox{4.1}{6.2}{-0.8}{0.7}{$\color{red}0\,\,0\color{black}$};
	\redbraidbox{0}{2.1}{3.2}{4.7}{$\color{red}0\,\,0\color{black}$};
	\draw (3.1,4) node{\color{blue}\normalsize$1$\color{black}};
	\redbraidbox{4.1}{6.2}{3.2}{4.7}{$\color{red}0\,\,0\color{black}$};
	\draw[blue](3.1,3)--(2,2.5)--(4.15,1.5)--(3.1,1);
	 \draw[red](0.4,3)--(4.4,1);
	 \draw[red](1.8,3)--(5.8,1);
	  \draw[red](0.4,1)--(4.4,3);
	 \draw[red](1.8,1)--(5.8,3);
	  \end{braid}
	  +
	  \begin{braid}\tikzset{baseline=1.6em}
	\redbraidbox{0}{2.1}{-0.8}{0.7}{$\color{red}0\,\,0\color{black}$};
	\draw (3.1,0) node{\color{blue}\normalsize$1$\color{black}};
	\redbraidbox{4.1}{6.2}{-0.8}{0.7}{$\color{red}0\,\,0\color{black}$};
	\redbraidbox{0}{2.1}{3.2}{4.7}{$\color{red}0\,\,0\color{black}$};
	\draw (3.1,4) node{\color{blue}\normalsize$1$\color{black}};
	\redbraidbox{4.1}{6.2}{3.2}{4.7}{$\color{red}0\,\,0\color{black}$};
	\draw[blue](3.1,3)--(4.4,2)--(3.1,1);
	 \draw[red](0.4,3)--(0.4,1);
	 \draw[red](1.8,3)--(4.4,1);
	  \draw[red](1.8,1)--(5.8,3);
	 \draw[red](5.8,1)--(4.4,3);
	  \reddot (0.4,2.9);
	   \reddot (0.4,2.6);
	  \end{braid}
	  +  \begin{braid}\tikzset{baseline=1.6em}
	\redbraidbox{0}{2.1}{-0.8}{0.7}{$\color{red}0\,\,0\color{black}$};
	\draw (3.1,0) node{\color{blue}\normalsize$1$\color{black}};
	\redbraidbox{4.1}{6.2}{-0.8}{0.7}{$\color{red}0\,\,0\color{black}$};
	\redbraidbox{0}{2.1}{3.2}{4.7}{$\color{red}0\,\,0\color{black}$};
	\draw (3.1,4) node{\color{blue}\normalsize$1$\color{black}};
	\redbraidbox{4.1}{6.2}{3.2}{4.7}{$\color{red}0\,\,0\color{black}$};
	\draw[blue](3.1,3)--(5,2)--(3.1,1);
	 \draw[red](0.4,3)--(0.4,1);
	 \draw[red](1.8,3)--(4.4,2)--(5.8,1);
	  \draw[red](1.8,1)--(4.4,2)--(5.8,3);
	 \draw[red](4.4,1)--(5.8,2)--(4.4,3);
	  \reddot (0.4,2.9);
	   \reddot (4.5,1.15);
	    \reddot (4.85,1.35);
	  \end{braid}
	 +  \begin{braid}\tikzset{baseline=1.6em}
	\redbraidbox{0}{2.1}{-0.8}{0.7}{$\color{red}0\,\,0\color{black}$};
	\draw (3.1,0) node{\color{blue}\normalsize$1$\color{black}};
	\redbraidbox{4.1}{6.2}{-0.8}{0.7}{$\color{red}0\,\,0\color{black}$};
	\redbraidbox{0}{2.1}{3.2}{4.7}{$\color{red}0\,\,0\color{black}$};
	\draw (3.1,4) node{\color{blue}\normalsize$1$\color{black}};
	\redbraidbox{4.1}{6.2}{3.2}{4.7}{$\color{red}0\,\,0\color{black}$};
	\draw[blue](3.1,3)--(5,2)--(3.1,1);
	 \draw[red](0.4,3)--(0.4,1);
	 \draw[red](1.8,3)--(4.4,2)--(5.8,1);
	  \draw[red](1.8,1)--(4.4,2)--(5.8,3);
	 \draw[red](4,1)--(5.8,2)--(4.4,3);
	  \reddot (4.71,1.45);
	   \reddot (4.1,1.1);
	    \reddot (4.4,1.28);
	  \end{braid}+}F,
$$
where
$$
F:=z_1^2-z_2^2+z_1u_1u_2-u_1^2u_2^2-u_1^2z_2-z_1u_2^2+u_1u_2z_2+u_1u_2^3-2u_2^2z_2-u_2^4
$$
and we have put
$$
z_1:=
\resizebox{21mm}{6mm}{
\begin{braid}\tikzset{baseline=1.6em}
	\redbraidbox{0}{2.1}{-0.8}{0.7}{$\color{red}0\,\,0\color{black}$};
	\draw (3.1,0) node{\color{blue}\normalsize$1$\color{black}};
	\redbraidbox{4.1}{6.2}{-0.8}{0.7}{$\color{red}0\,\,0\color{black}$};
	\redbraidbox{0}{2.1}{3.2}{4.7}{$\color{red}0\,\,0\color{black}$};
	\draw (3.1,4) node{\color{blue}\normalsize$1$\color{black}};
	\redbraidbox{4.1}{6.2}{3.2}{4.7}{$\color{red}0\,\,0\color{black}$};
	\draw[blue](3.1,3)--(3.1,1);
	 \draw[red](0.4,3)--(0.4,1);
	 \draw[red](1.8,3)--(1.8,1);
	  \draw[red](4.4,1)--(4.4,3);
	 \draw[red](5.8,1)--(5.8,3);
	  \reddot (0.4,2.2);
	   \reddot (1.8,1.8);
	  \end{braid},
	  }
	  \  
	  z_2:=
	  \resizebox{21mm}{6mm}{
	  \begin{braid}\tikzset{baseline=1.6em}
	\redbraidbox{0}{2.1}{-0.8}{0.7}{$\color{red}0\,\,0\color{black}$};
	\draw (3.1,0) node{\color{blue}\normalsize$1$\color{black}};
	\redbraidbox{4.1}{6.2}{-0.8}{0.7}{$\color{red}0\,\,0\color{black}$};
	\redbraidbox{0}{2.1}{3.2}{4.7}{$\color{red}0\,\,0\color{black}$};
	\draw (3.1,4) node{\color{blue}\normalsize$1$\color{black}};
	\redbraidbox{4.1}{6.2}{3.2}{4.7}{$\color{red}0\,\,0\color{black}$};
	\draw[blue](3.1,3)--(3.1,1);
	 \draw[red](0.4,3)--(0.4,1);
	 \draw[red](1.8,3)--(1.8,1);
	  \draw[red](4.4,1)--(4.4,3);
	 \draw[red](5.8,1)--(5.8,3);
	  \reddot (4.4,2.2);
	   \reddot (5.8,1.8);
	  \end{braid},
	  }
	  \  
	  u_1:=
	  \resizebox{21mm}{6mm}{
	  \begin{braid}\tikzset{baseline=1.6em}
	\redbraidbox{0}{2.1}{-0.8}{0.7}{$\color{red}0\,\,0\color{black}$};
	\draw (3.1,0) node{\color{blue}\normalsize$1$\color{black}};
	\redbraidbox{4.1}{6.2}{-0.8}{0.7}{$\color{red}0\,\,0\color{black}$};
	\redbraidbox{0}{2.1}{3.2}{4.7}{$\color{red}0\,\,0\color{black}$};
	\draw (3.1,4) node{\color{blue}\normalsize$1$\color{black}};
	\redbraidbox{4.1}{6.2}{3.2}{4.7}{$\color{red}0\,\,0\color{black}$};
	\draw[blue](3.1,3)--(3.1,1);
	 \draw[red](0.4,3)--(0.4,1);
	 \draw[red](1.8,3)--(1.8,1);
	  \draw[red](4.4,1)--(4.4,3);
	 \draw[red](5.8,1)--(5.8,3);
	  \reddot (1.8,2);
	  \end{braid}
	  }
,
	  \  
	  u_2:=
	  \resizebox{21mm}{6mm}{
	  \begin{braid}\tikzset{baseline=1.6em}
	\redbraidbox{0}{2.1}{-0.8}{0.7}{$\color{red}0\,\,0\color{black}$};
	\draw (3.1,0) node{\color{blue}\normalsize$1$\color{black}};
	\redbraidbox{4.1}{6.2}{-0.8}{0.7}{$\color{red}0\,\,0\color{black}$};
	\redbraidbox{0}{2.1}{3.2}{4.7}{$\color{red}0\,\,0\color{black}$};
	\draw (3.1,4) node{\color{blue}\normalsize$1$\color{black}};
	\redbraidbox{4.1}{6.2}{3.2}{4.7}{$\color{red}0\,\,0\color{black}$};
	\draw[blue](3.1,3)--(3.1,1);
	 \draw[red](0.4,3)--(0.4,1);
	 \draw[red](1.8,3)--(1.8,1);
	  \draw[red](4.4,1)--(4.4,3);
	 \draw[red](5.8,1)--(5.8,3);
	  \reddot (5.8,2);
	  \end{braid}
	  }.
$$
\end{Lemma}

\vspace{1mm}
\begin{Lemma} \label{L090721} 
In $R_{4\al_{0}}$ we have 
$$
\begin{braid}\tikzset{baseline=1.3em}
                \draw[red](2,1)node[below]{$0$}--(0.7,1.9)--(0,3);
	 \draw[red](3,1)node[below]{$0$}--(2.3,2.2)--(1,3);
	 \draw[red](0,1)node[below]{$0$}--(2,3);
	 \draw[red](1,1)node[below]{$0$}--(3,3);
	 \reddot (.25,1.25);
	 \reddot (1.1,1.1);
	  \end{braid}
	  =
	  \begin{braid}\tikzset{baseline=1.3em}
                \draw[red](2,1)node[below]{$0$}--(0.7,1.9)--(0,3);
	 \draw[red](3,1)node[below]{$0$}--(2.3,2.2)--(1,3);
	 \draw[red](0,1)node[below]{$0$}--(2,3);
	 \draw[red](1,1)node[below]{$0$}--(3,3);
	 \reddot (1.9,2.9);
	 \reddot (2.75,2.75);
	  \end{braid}
	  +\begin{braid}\tikzset{baseline=1.3em}
                \draw[red](2,1)node[below]{$0$}--(0,3);
	 \draw[red](3,1)node[below]{$0$}--(2.3,1.4)--(2,3);
	 \draw[red](0,1)node[below]{$0$}--(0.5,2.7)--(1,3);
	 \draw[red](1,1)node[below]{$0$}--(3,3);
	 \reddot (1.2,1.2);
	  \end{braid}.
$$
\end{Lemma}

\section{A commutation lemma}
In this subsection we set 
$$
\theta:=\al_i+4(\al_{i+1}+\dots+\al_{\ell-1})+2\al_\ell
$$
for $0\leq i<\ell$. 
We consider the divided power idempotent $e\in R_\theta$ defined as follows: 
\begin{align*}
e&:=e(\ell\,(\ell-1)^{(2)}\cdots(i+2)^{(2)}(i+1)^{(2)}\,i\,\ell\,(\ell-1)^{(2)}\cdots(i+2)^{(2)}(i+1)^{(2)})
\\&=
\begin{braid}\tikzset{baseline=-.3em}
	\draw (-0.2,0) node{\color{blue}$\ell$\color{black}};
	\braidbox{0.6}{2.95}{-0.7}{.6}{};
	\draw (1.8,0) node{$(\ell-1)^2$};
	\draw[dots] (3.5,0)--(5,0);
	\braidbox{5.3}{7.6}{-.7}{.6}{};
	\draw (6.5,0) node{$(i+2)^2$};
	\braidbox{8.3}{10.6}{-.7}{.6}{};
	\draw (9.5,0) node{$(i+1)^2$};
	\draw (11.5,0) node{$i$};
\draw (12.3,0) node{\color{blue}$\ell$\color{black}};
\braidbox{13.1}{15.44}{-0.7}{.6}{};
	\draw (14.3,0) node{$(\ell-1)^2$};
	\draw[dots] (16,0)--(17.4,0);
	\braidbox{17.8}{20.1}{-.7}{.6}{};
	\draw (19,0) node{$(i+2)^2$};
	\braidbox{20.8}{23.1}{-.7}{.6}{};
	\draw (22,0) node{$(i+1)^2$};	
	\end{braid}.
\end{align*}

We consider the left ideal $\I$ in $R_{\theta}$ generated by the elements 
\begin{align*}
&\psi_1 e,\, \psi_2\psi_3 e,\,\dots,\,\psi_{2\ell-2i-4}\psi_{2\ell-2i-3}e,\, \psi_{2\ell-2i-2}\psi_{2\ell-2i-1}e,
\\ 
&\psi_{2\ell-2i}e,\, \psi_{2\ell-2i+2}\psi_{2\ell-2i+3}e,\,\dots,\, \psi_{4\ell-4i-3}\psi_{4\ell-4i-2}e.
\end{align*}
In terms of the diagrams, we have:
\begin{align*}
\psi_1 e&=\begin{braid}\tikzset{baseline=0.3em}
	\draw (-0.2,0) node{\color{blue}$\ell$\color{black}};
	\braidbox{0.6}{2.95}{-0.7}{.6}{};
	\draw (1.8,0) node{$(\ell-1)^2$};
	\draw[dots] (3.5,0)--(5,0);
	\braidbox{5.3}{7.6}{-.7}{.6}{};
	\draw (6.5,0) node{$(i+2)^2$};
	\braidbox{8.3}{10.6}{-.7}{.6}{};
	\draw (9.5,0) node{$(i+1)^2$};
	\draw (11.5,0) node{$i$};
\draw (12.3,0) node{\color{blue}$\ell$\color{black}};
\braidbox{13.1}{15.44}{-0.7}{.6}{};
	\draw (14.3,0) node{$(\ell-1)^2$};
	\draw[dots] (16,0)--(17.4,0);
	\braidbox{17.8}{20.1}{-.7}{.6}{};
	\draw (19,0) node{$(i+2)^2$};
	\braidbox{20.8}{23.1}{-.7}{.6}{};
	\draw (22,0) node{$(i+1)^2$};
	\draw[blue](-0.2,0.8)--(1,1.8);
	\draw(1,0.8)--(-0.2,1.8);	
	\end{braid},
	\\
	\psi_2\psi_3 e&=
	\begin{braid}\tikzset{baseline=.3em}
	\draw (-3.2,0) node{\color{blue}$\ell$\color{black}};
	\braidbox{-2.4}{-.05}{-0.7}{.6}{};
	\draw (-1.2,0) node{$(\ell-1)^2$};
	\braidbox{0.6}{2.95}{-0.7}{.6}{};
	\draw (1.8,0) node{$(\ell-2)^2$};
	\draw[dots] (3.5,0)--(5,0);
	\braidbox{5.3}{7.6}{-.7}{.6}{};
	\draw (6.5,0) node{$(i+2)^2$};
	\braidbox{8.3}{10.6}{-.7}{.6}{};
	\draw (9.5,0) node{$(i+1)^2$};
	\draw (11.5,0) node{$i$};
\draw (12.3,0) node{\color{blue}$\ell$\color{black}};
\braidbox{13.1}{15.44}{-0.7}{.6}{};
	\draw (14.3,0) node{$(\ell-1)^2$};
	\draw[dots] (16,0)--(17.4,0);
	\braidbox{17.8}{20.1}{-.7}{.6}{};
	\draw (19,0) node{$(i+2)^2$};
	\braidbox{20.8}{23.1}{-.7}{.6}{};
	\draw (22,0) node{$(i+1)^2$};
	\draw(-0.5,0.8)--(1,1.8);
	\draw(1,0.8)--(-2,1.8);	
	\draw(-2,0.8)--(-0.5,1.8);
	\end{braid},
	\\&\ \, \vdots
	\\
	\psi_{2\ell-2i-4}\psi_{2\ell-2i-3}e&=
	\begin{braid}\tikzset{baseline=.3em}
	\draw (-3.2,0) node{\color{blue}$\ell$\color{black}};
	\braidbox{-2.4}{-.05}{-0.7}{.6}{};
	\draw (-1.2,0) node{$(\ell-1)^2$};
	\braidbox{0.6}{2.95}{-0.7}{.6}{};
	\draw (1.8,0) node{$(\ell-2)^2$};
	\draw[dots] (3.5,0)--(5,0);
	\braidbox{5.3}{7.6}{-.7}{.6}{};
	\draw (6.5,0) node{$(i+2)^2$};
	\braidbox{8.3}{10.6}{-.7}{.6}{};
	\draw (9.5,0) node{$(i+1)^2$};
	\draw (11.5,0) node{$i$};
\draw (12.3,0) node{\color{blue}$\ell$\color{black}};
\braidbox{13.1}{15.44}{-0.7}{.6}{};
	\draw (14.3,0) node{$(\ell-1)^2$};
	\draw[dots] (16,0)--(17.4,0);
	\braidbox{17.8}{20.1}{-.7}{.6}{};
	\draw (19,0) node{$(i+2)^2$};
	\braidbox{20.8}{23.1}{-.7}{.6}{};
	\draw (22,0) node{$(i+1)^2$};
	\draw(7.3,0.8)--(8.8,1.8);
	\draw(8.8,0.8)--(5.8,1.8);	
	\draw(5.8,0.8)--(7.3,1.8);
	\end{braid},
	\\
	\psi_{2\ell-2i-2}\psi_{2\ell-2i-1}e&=
	\begin{braid}\tikzset{baseline=.3em}
	\draw (-3.2,0) node{\color{blue}$\ell$\color{black}};
	\braidbox{-2.4}{-.05}{-0.7}{.6}{};
	\draw (-1.2,0) node{$(\ell-1)^2$};
	\braidbox{0.6}{2.95}{-0.7}{.6}{};
	\draw (1.8,0) node{$(\ell-2)^2$};
	\draw[dots] (3.5,0)--(5,0);
	\braidbox{5.3}{7.6}{-.7}{.6}{};
	\draw (6.5,0) node{$(i+2)^2$};
	\braidbox{8.3}{10.6}{-.7}{.6}{};
	\draw (9.5,0) node{$(i+1)^2$};
	\draw (11.5,0) node{$i$};
\draw (12.3,0) node{\color{blue}$\ell$\color{black}};
\braidbox{13.1}{15.44}{-0.7}{.6}{};
	\draw (14.3,0) node{$(\ell-1)^2$};
	\draw[dots] (16,0)--(17.4,0);
	\braidbox{17.8}{20.1}{-.7}{.6}{};
	\draw (19,0) node{$(i+2)^2$};
	\braidbox{20.8}{23.1}{-.7}{.6}{};
	\draw (22,0) node{$(i+1)^2$};
	\draw(10.2,0.8)--(11.5,1.8);
	\draw(11.5,0.8)--(8.7,1.8);	
	\draw(8.7,0.8)--(10.2,1.8);
	\end{braid},
	\end{align*}
	\begin{align*}
	\psi_{2\ell-2i}e&=
	\begin{braid}\tikzset{baseline=.3em}
	\draw (-0.2,0) node{\color{blue}$\ell$\color{black}};
	\braidbox{0.6}{2.95}{-0.7}{.6}{};
	\draw (1.8,0) node{$(\ell-1)^2$};
	\draw[dots] (3.5,0)--(5,0);
	\braidbox{5.3}{7.6}{-.7}{.6}{};
	\draw (6.5,0) node{$(i+2)^2$};
	\braidbox{8.3}{10.6}{-.7}{.6}{};
	\draw (9.5,0) node{$(i+1)^2$};
	\draw (11.5,0) node{$i$};
\draw (12.3,0) node{\color{blue}$\ell$\color{black}};
\braidbox{13.1}{15.44}{-0.7}{.6}{};
	\draw (14.3,0) node{$(\ell-1)^2$};
	\braidbox{16.1}{18.44}{-0.7}{.6}{};
	\draw (17.3,0) node{$(\ell-2)^2$};
	\draw[dots] (19,0)--(20.4,0);
	\braidbox{20.8}{23.1}{-.7}{.6}{};
	\draw (22,0) node{$(i+2)^2$};
	\braidbox{23.8}{26.1}{-.7}{.6}{};
	\draw (25,0) node{$(i+1)^2$};
	\draw[blue](12.3,0.8)--(13.5,1.8);
	\draw(13.5,0.8)--(12.3,1.8);	
	\end{braid},
	\\
	\psi_{2\ell-2i+2}\psi_{2\ell-2i+3}e&=
	\begin{braid}\tikzset{baseline=.3em}
	\draw (-0.2,0) node{\color{blue}$\ell$\color{black}};
	\braidbox{0.6}{2.95}{-0.7}{.6}{};
	\draw (1.8,0) node{$(\ell-1)^2$};
	\draw[dots] (3.5,0)--(5,0);
	\braidbox{5.3}{7.6}{-.7}{.6}{};
	\draw (6.5,0) node{$(i+2)^2$};
	\braidbox{8.3}{10.6}{-.7}{.6}{};
	\draw (9.5,0) node{$(i+1)^2$};
	\draw (11.5,0) node{$i$};
\draw (12.3,0) node{\color{blue}$\ell$\color{black}};
\braidbox{13.1}{15.44}{-0.7}{.6}{};
	\draw (14.3,0) node{$(\ell-1)^2$};
	\braidbox{16.1}{18.44}{-0.7}{.6}{};
	\draw (17.3,0) node{$(\ell-2)^2$};
	\draw[dots] (19,0)--(20.4,0);
	\braidbox{20.8}{23.1}{-.7}{.6}{};
	\draw (22,0) node{$(i+2)^2$};
	\braidbox{23.8}{26.1}{-.7}{.6}{};
	\draw (25,0) node{$(i+1)^2$};
	\draw(15.2,0.8)--(16.5,1.8);
	\draw(16.5,0.8)--(13.7,1.8);	
	\draw(13.7,0.8)--(15.2,1.8);
	\end{braid},
	\\
	&\ \,\vdots
	\\
	\psi_{4\ell-4i-3}\psi_{4\ell-4i-2}e&=
	\begin{braid}\tikzset{baseline=.3em}
	\draw (-0.2,0) node{\color{blue}$\ell$\color{black}};
	\braidbox{0.6}{2.95}{-0.7}{.6}{};
	\draw (1.8,0) node{$(\ell-1)^2$};
	\draw[dots] (3.5,0)--(5,0);
	\braidbox{5.3}{7.6}{-.7}{.6}{};
	\draw (6.5,0) node{$(i+2)^2$};
	\braidbox{8.3}{10.6}{-.7}{.6}{};
	\draw (9.5,0) node{$(i+1)^2$};
	\draw (11.5,0) node{$i$};
\draw (12.3,0) node{\color{blue}$\ell$\color{black}};
\braidbox{13.1}{15.44}{-0.7}{.6}{};
	\draw (14.3,0) node{$(\ell-1)^2$};
	\braidbox{16.1}{18.44}{-0.7}{.6}{};
	\draw (17.3,0) node{$(\ell-2)^2$};
	\draw[dots] (19,0)--(20.4,0);
	\braidbox{20.8}{23.1}{-.7}{.6}{};
	\draw (22,0) node{$(i+2)^2$};
	\braidbox{23.8}{26.1}{-.7}{.6}{};
	\draw (25,0) node{$(i+1)^2$};
	\draw(22.8,0.8)--(24.1,1.8);
	\draw(24.1,0.8)--(21.3,1.8);	
	\draw(21.3,0.8)--(22.8,1.8);
	\end{braid}.
\end{align*}

\vspace{1mm}
\begin{Lemma} \label{L100721_0} 
In $R_\theta$ we have 
$$
\begin{braid}\tikzset{baseline=-.3em}
		\draw (-0.2,-4) node{\color{blue}$\ell$\color{black}};
	\braidbox{0.6}{2.95}{-3.3}{-4.6}{};
	\draw (1.8,-4) node{$(\ell-1)^2$};
	\draw[dots] (3.5,-4)--(5,-4);
	\braidbox{5.3}{7.6}{-3.3}{-4.6}{};
	\draw (6.5,-4) node{$(i+2)^2$};
	\braidbox{8.3}{10.6}{-3.3}{-4.6}{};
	\draw (9.5,-4) node{$(i+1)^2$};
	\draw (11.5,-4) node{$i$};
\draw (12.3,-4) node{\color{blue}$\ell$\color{black}};
\braidbox{13.1}{15.44}{-3.3}{-4.6}{};
	\draw (14.3,-4) node{$(\ell-1)^2$};
	\draw[dots] (16,-4)--(17.4,-4);
	\braidbox{17.8}{20.1}{-3.3}{-4.6}{};
	\draw (19,-4) node{$(i+2)^2$};
	\braidbox{20.8}{23.1}{-3.3}{-4.6}{};
	\draw (22,-4) node{$(i+1)^2$};	
	\draw[blue](-0.2,-3.2)--(12.3,3.2);
	\draw[blue](-0.2,3.2)--(12.3,-3.2);
	\draw(1,-3.2)--(13.5,3.2);
	\draw(1,3.2)--(13.5,-3.2);
	\draw(2.5,-3.2)--(15,3.2);
	\draw(2.5,3.2)--(15,-3.2);
	\draw(5.7,-3.2)--(18.2,3.2);
	\draw(5.7,3.2)--(18.2,-3.2);
	\draw(7.2,-3.2)--(19.7,3.2);
	\draw(7.2,3.2)--(19.7,-3.2);
	\draw(8.7,-3.2)--(21.2,3.2);
	\draw(8.7,3.2)--(21.2,-3.2);
	\draw(10.2,-3.2)--(22.7,3.2);
	\draw(10.2,3.2)--(22.7,-3.2);
	\draw(11.5,-3.2)--(18,0)--(11.5,3.2);
	\end{braid}
	\ \equiv-e\pmod{\I}.
$$
\end{Lemma}
\begin{proof}
We will write `$\equiv$' for `$\equiv\pmod{\I}$'. 
We induct on $\ell-i$. If $\ell-i=1$, we have using braid relations (\ref{R7}): 
$$
\begin{braid}\tikzset{baseline=-.3em}
	\draw (0,-2) node{\color{blue}$\ell$\color{black}};
	\draw (1.5,-2) node{$\ell-1$};
	\draw (3,-2) node{\color{blue}$\ell$\color{black}};
	\draw[blue](0,-1.2)--(3,1.2);
	\draw[blue](0,1.2)--(3,-1.2);
	\draw(1.5,-1.2)--(2.5,0)--(1.5,1.2);
	\end{braid}
	=\begin{braid}\tikzset{baseline=-.3em}
	\draw (0,-2) node{\color{blue}$\ell$\color{black}};
	\draw (1.5,-2) node{$\ell-1$};
	\draw (3,-2) node{\color{blue}$\ell$\color{black}};
	\draw[blue](0,-1.2)--(3,1.2);
	\draw[blue](0,1.2)--(3,-1.2);
	\draw(1.5,-1.2)--(0.5,0)--(1.5,1.2);
	\end{braid}-
	\begin{braid}\tikzset{baseline=-.3em}
		\draw (0,-2) node{\color{blue}$\ell$\color{black}};
	\draw (1.5,-2) node{$\ell-1$};
	\draw (3,-2) node{\color{blue}$\ell$\color{black}};
	\draw[blue](0,-1.2)--(0,1.2);
	\draw[blue](3,1.2)--(3,-1.2);
	\draw(1.5,-1.2)--(1.5,1.2);
	\end{braid}
	\equiv
	-
	\begin{braid}\tikzset{baseline=-.3em}
		\draw (0,-2) node{\color{blue}$\ell$\color{black}};
	\draw (1.5,-2) node{$\ell-1$};
	\draw (3,-2) node{\color{blue}$\ell$\color{black}};
	\draw[blue](0,-1.2)--(0,1.2);
	\draw[blue](3,1.2)--(3,-1.2);
	\draw(1.5,-1.2)--(1.5,1.2);
	\end{braid}
=-e.
$$
For the inductive step, let $\ell-i>1$. We have using Lemma~\ref{L090721_3} and braid relations:
\begin{align*}
&\,\begin{braid}\tikzset{baseline=-.3em}
	\draw (-0.2,-4) node{\color{blue}$\ell$\color{black}};
	\braidbox{0.6}{2.95}{-3.3}{-4.6}{};
	\draw (1.8,-4) node{$(\ell-1)^2$};
	\draw[dots] (3.5,-4)--(5,-4);
	\braidbox{5.3}{7.6}{-3.3}{-4.6}{};
	\draw (6.5,-4) node{$(i+2)^2$};
	\braidbox{8.3}{10.6}{-3.3}{-4.6}{};
	\draw (9.5,-4) node{$(i+1)^2$};
	\draw (11.5,-4) node{$i$};
\draw (12.3,-4) node{\color{blue}$\ell$\color{black}};
\braidbox{13.1}{15.44}{-3.3}{-4.6}{};
	\draw (14.3,-4) node{$(\ell-1)^2$};
	\draw[dots] (16,-4)--(17.4,-4);
	\braidbox{17.8}{20.1}{-3.3}{-4.6}{};
	\draw (19,-4) node{$(i+2)^2$};
	\braidbox{20.8}{23.1}{-3.3}{-4.6}{};
	\draw (22,-4) node{$(i+1)^2$};	
	\draw[blue](-0.2,-3.2)--(12.3,3.2);
	\draw[blue](-0.2,3.2)--(12.3,-3.2);
	\draw(1,-3.2)--(13.5,3.2);
	\draw(1,3.2)--(13.5,-3.2);
	\draw(2.5,-3.2)--(15,3.2);
	\draw(2.5,3.2)--(15,-3.2);
	\draw(5.7,-3.2)--(18.2,3.2);
	\draw(5.7,3.2)--(18.2,-3.2);
	\draw(7.2,-3.2)--(19.7,3.2);
	\draw(7.2,3.2)--(19.7,-3.2);
	\draw(8.7,-3.2)--(21.2,3.2);
	\draw(8.7,3.2)--(21.2,-3.2);
	\draw(10.2,-3.2)--(22.7,3.2);
	\draw(10.2,3.2)--(22.7,-3.2);
	\draw(11.5,-3.2)--(18,0)--(11.5,3.2);
	\end{braid}
	\\
	\equiv&\,-
	\begin{braid}\tikzset{baseline=-.3em}
	\draw (-0.2,-4) node{\color{blue}$\ell$\color{black}};
	\braidbox{0.6}{2.95}{-3.3}{-4.6}{};
	\draw (1.8,-4) node{$(\ell-1)^2$};
	\draw[dots] (3.5,-4)--(5,-4);
	\braidbox{5.3}{7.6}{-3.3}{-4.6}{};
	\draw (6.5,-4) node{$(i+2)^2$};
	\braidbox{8.3}{10.6}{-3.3}{-4.6}{};
	\draw (9.5,-4) node{$(i+1)^2$};
	\draw (11.5,-4) node{$i$};
\draw (12.3,-4) node{\color{blue}$\ell$\color{black}};
\braidbox{13.1}{15.44}{-3.3}{-4.6}{};
	\draw (14.3,-4) node{$(\ell-1)^2$};
	\draw[dots] (16,-4)--(17.4,-4);
	\braidbox{17.8}{20.1}{-3.3}{-4.6}{};
	\draw (19,-4) node{$(i+2)^2$};
	\braidbox{20.8}{23.1}{-3.3}{-4.6}{};
	\draw (22,-4) node{$(i+1)^2$};	
	\draw[blue](-0.2,-3.2)--(12.3,3.2);
	\draw[blue](-0.2,3.2)--(12.3,-3.2);
	\draw(1,-3.2)--(13.5,3.2);
	\draw(1,3.2)--(13.5,-3.2);
	\draw(2.5,-3.2)--(15,3.2);
	\draw(2.5,3.2)--(15,-3.2);
	\draw(5.7,-3.2)--(18.2,3.2);
	\draw(5.7,3.2)--(18.2,-3.2);
	\draw(7.2,-3.2)--(19.7,3.2);
	\draw(7.2,3.2)--(19.7,-3.2);
	\draw(8.7,-3.2)--(15.6,0.4)--(10.2,3.2);
	\draw(8.7,3.2)--(21.2,-3.2);
	\draw(10.2,-3.2)--(22.7,3.2);
	\draw(21.3,3.2)--(16.9,0.8)--(22.7,-3.2);
	\draw(11.5,-3.2)--(16.5,-1.1)--(14.8,-0.4)--(16.6,0.5)--(11.5,3.2);
	\end{braid}
	\\
	=&\,
	-\begin{braid}\tikzset{baseline=-.3em}
	\draw (-0.2,-4) node{\color{blue}$\ell$\color{black}};
	\braidbox{0.6}{2.95}{-3.3}{-4.6}{};
	\draw (1.8,-4) node{$(\ell-1)^2$};
	\draw[dots] (3.5,-4)--(5,-4);
	\braidbox{5.3}{7.6}{-3.3}{-4.6}{};
	\draw (6.5,-4) node{$(i+2)^2$};
	\braidbox{8.3}{10.6}{-3.3}{-4.6}{};
	\draw (9.5,-4) node{$(i+1)^2$};
	\draw (11.5,-4) node{$i$};
\draw (12.3,-4) node{\color{blue}$\ell$\color{black}};
\braidbox{13.1}{15.44}{-3.3}{-4.6}{};
	\draw (14.3,-4) node{$(\ell-1)^2$};
	\draw[dots] (16,-4)--(17.4,-4);
	\braidbox{17.8}{20.1}{-3.3}{-4.6}{};
	\draw (19,-4) node{$(i+2)^2$};
	\braidbox{20.8}{23.1}{-3.3}{-4.6}{};
	\draw (22,-4) node{$(i+1)^2$};	
	\draw[blue](-0.2,-3.2)--(12.3,3.2);
	\draw[blue](-0.2,3.2)--(12.3,-3.2);
	\draw(1,-3.2)--(13.5,3.2);
	\draw(1,3.2)--(13.5,-3.2);
	\draw(2.5,-3.2)--(15,3.2);
	\draw(2.5,3.2)--(15,-3.2);
	\draw(5.7,-3.2)--(18.2,3.2);
	\draw(5.7,3.2)--(18.2,-3.2);
	\draw(7.2,-3.2)--(19.7,3.2);
	\draw(7.2,3.2)--(19.7,-3.2);
	\draw(8.7,-3.2)--(14.7,0)--(9.3,2.7)--(10.2,3.2);
	\draw(9,3.2)--(10.5,2.8)--(15,0.6)--(21.2,-3.2);
	\draw(10.2,-3.2)--(22.7,3.2);
	\draw(21.3,3.2)--(22.7,-3.2);
	\draw(11.5,-3.2)--(16.4,-0.9)--(15,-0.2)--(16.4,0.6)--(11.5,3.2);
	\end{braid}.
\end{align*}
By the inductive assumption, the last expression equals
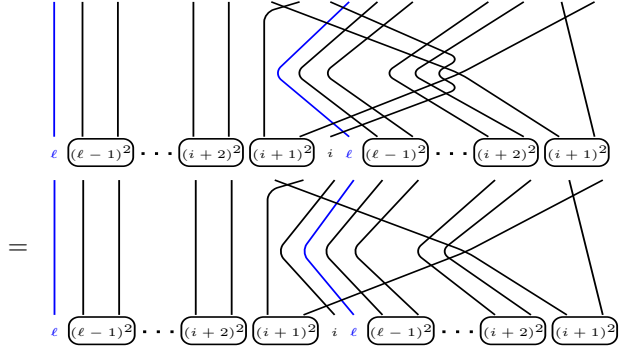
\begin{align*}
	&
	\resizebox{80mm}{10mm}{
	\begin{braid}\tikzset{baseline=-.3em}
	\draw (-0.2,-4) node{\color{blue}$\ell$\color{black}};
	\braidbox{0.6}{2.95}{-3.3}{-4.6}{};
	\draw (1.8,-4) node{$(\ell-1)^2$};
	\draw[dots] (3.5,-4)--(5,-4);
	\braidbox{5.3}{7.6}{-3.3}{-4.6}{};
	\draw (6.5,-4) node{$(i+2)^2$};
	\braidbox{8.3}{10.6}{-3.3}{-4.6}{};
	\draw (9.5,-4) node{$(i+1)^2$};
	\draw (11.5,-4) node{$i$};
\draw (12.3,-4) node{\color{blue}$\ell$\color{black}};
\braidbox{13.1}{15.44}{-3.3}{-4.6}{};
	\draw (14.3,-4) node{$(\ell-1)^2$};
	\draw[dots] (16,-4)--(17.4,-4);
	\braidbox{17.8}{20.1}{-3.3}{-4.6}{};
	\draw (19,-4) node{$(i+2)^2$};
	\braidbox{20.8}{23.1}{-3.3}{-4.6}{};
	\draw (22,-4) node{$(i+1)^2$};	
	\draw[blue](-0.2,-3.2)--(-0.2,3.2);
	\draw[blue](12.3,3.2)--(9.1,-0.2)--(12.3,-3.2);
	\draw(1,-3.2)--(1,3.2);
	\draw(13.5,3.2)--(10,-0.2)--(13.5,-3.2);
	\draw(2.5,-3.2)--(2.5,3.2);
	\draw(15,3.2)--(11.2,-0.2)--(15,-3.2);
	\draw(5.7,-3.2)--(5.7,3.2);
	\draw(18.2,3.2)--(13.8,-0.2)--(18.2,-3.2);
	\draw(7.2,-3.2)--(7.2,3.2);
	\draw(19.7,3.2)--(14.9,-0.2)--(19.7,-3.2);
	\draw(8.7,-3.2)--(8.7,2.7)--(10.2,3.2);
	\draw(9,3.2)--(17,-0.2)--(21.2,-3.2);
	\draw(10.2,-3.2)--(17,-0.2)--(22.7,3.2);
	\draw(21.3,3.2)--(22.7,-3.2);
	\draw(11.5,-3.2)--(17,-0.9)--(15.9,-0.2)--(17,0.5)--(11.5,3.2);
	\end{braid}
	}
	\\
	=\,
	&
	\resizebox{80mm}{10mm}{
	\begin{braid}\tikzset{baseline=-.3em}
	\draw (-0.2,-4) node{\color{blue}$\ell$\color{black}};
	\braidbox{0.6}{2.95}{-3.3}{-4.6}{};
	\draw (1.8,-4) node{$(\ell-1)^2$};
	\draw[dots] (3.5,-4)--(5,-4);
	\braidbox{5.3}{7.6}{-3.3}{-4.6}{};
	\draw (6.5,-4) node{$(i+2)^2$};
	\braidbox{8.3}{10.6}{-3.3}{-4.6}{};
	\draw (9.5,-4) node{$(i+1)^2$};
	\draw (11.5,-4) node{$i$};
\draw (12.3,-4) node{\color{blue}$\ell$\color{black}};
\braidbox{13.1}{15.44}{-3.3}{-4.6}{};
	\draw (14.3,-4) node{$(\ell-1)^2$};
	\draw[dots] (16,-4)--(17.4,-4);
	\braidbox{17.8}{20.1}{-3.3}{-4.6}{};
	\draw (19,-4) node{$(i+2)^2$};
	\braidbox{20.8}{23.1}{-3.3}{-4.6}{};
	\draw (22,-4) node{$(i+1)^2$};	
	\draw[blue](-0.2,-3.2)--(-0.2,3.2);
	\draw[blue](12.3,3.2)--(10.1,-0.2)--(12.3,-3.2);
	\draw(1,-3.2)--(1,3.2);
	\draw(13.5,3.2)--(11,-0.2)--(13.5,-3.2);
	\draw(2.5,-3.2)--(2.5,3.2);
	\draw(15,3.2)--(12.2,-0.2)--(15,-3.2);
	\draw(5.7,-3.2)--(5.7,3.2);
	\draw(18.2,3.2)--(14.8,-0.2)--(18.2,-3.2);
	\draw(7.2,-3.2)--(7.2,3.2);
	\draw(19.7,3.2)--(15.9,-0.2)--(19.7,-3.2);
	\draw(8.7,-3.2)--(8.7,2.7)--(10.2,3.2);
	\draw(9,3.2)--(17,-0.2)--(21.2,-3.2);
	\draw(10.2,-3.2)--(17,-0.2)--(22.7,3.2);
	\draw(21.3,3.2)--(22.7,-3.2);
	\draw(11.5,-3.2)--(9.1,-0.2)--(11.5,3.2);
	\end{braid}}.
\end{align*}
Applying Lemma~\ref{LA3/4} (with `$i=i+1$') and ignoring the summands which become zero due to quadratic relations, we get modulo $\I$:
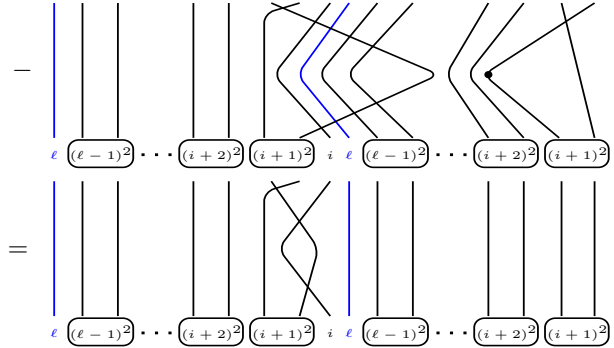
\begin{align*}
-&
\resizebox{80mm}{10mm}{
\begin{braid}\tikzset{baseline=-.3em}
	\draw (-0.2,-4) node{\color{blue}$\ell$\color{black}};
	\braidbox{0.6}{2.95}{-3.3}{-4.6}{};
	\draw (1.8,-4) node{$(\ell-1)^2$};
	\draw[dots] (3.5,-4)--(5,-4);
	\braidbox{5.3}{7.6}{-3.3}{-4.6}{};
	\draw (6.5,-4) node{$(i+2)^2$};
	\braidbox{8.3}{10.6}{-3.3}{-4.6}{};
	\draw (9.5,-4) node{$(i+1)^2$};
	\draw (11.5,-4) node{$i$};
\draw (12.3,-4) node{\color{blue}$\ell$\color{black}};
\braidbox{13.1}{15.44}{-3.3}{-4.6}{};
	\draw (14.3,-4) node{$(\ell-1)^2$};
	\draw[dots] (16,-4)--(17.4,-4);
	\braidbox{17.8}{20.1}{-3.3}{-4.6}{};
	\draw (19,-4) node{$(i+2)^2$};
	\braidbox{20.8}{23.1}{-3.3}{-4.6}{};
	\draw (22,-4) node{$(i+1)^2$};	
	\draw[blue](-0.2,-3.2)--(-0.2,3.2);
	\draw[blue](12.3,3.2)--(10.1,-0.2)--(12.3,-3.2);
	\draw(1,-3.2)--(1,3.2);
	\draw(13.5,3.2)--(11,-0.2)--(13.5,-3.2);
	\draw(2.5,-3.2)--(2.5,3.2);
	\draw(15,3.2)--(12.2,-0.2)--(15,-3.2);
	\draw(5.7,-3.2)--(5.7,3.2);
	\draw(18.2,3.2)--(16.4,-0.2)--(18.2,-3.2);
	\draw(7.2,-3.2)--(7.2,3.2);
	\draw(19.7,3.2)--(17.3,-0.2)--(19.7,-3.2);
	\draw(8.7,-3.2)--(8.7,2.7)--(10.2,3.2);
	\draw(9,3.2)--(16.1,-0.2)--(10.2,-3.2);
	\draw(21.2,-3.2)--(18,-0.2)--(22.7,3.2);
	\draw(21.3,3.2)--(22.7,-3.2);
	\draw(11.5,-3.2)--(9.1,-0.2)--(11.5,3.2);
	\blackdot(18.2,-.2);
	\end{braid}
	}
	\\=\,
	&
	\resizebox{80mm}{10mm}{
	\begin{braid}\tikzset{baseline=-.3em}
	\draw (-0.2,-4) node{\color{blue}$\ell$\color{black}};
	\braidbox{0.6}{2.95}{-3.3}{-4.6}{};
	\draw (1.8,-4) node{$(\ell-1)^2$};
	\draw[dots] (3.5,-4)--(5,-4);
	\braidbox{5.3}{7.6}{-3.3}{-4.6}{};
	\draw (6.5,-4) node{$(i+2)^2$};
	\braidbox{8.3}{10.6}{-3.3}{-4.6}{};
	\draw (9.5,-4) node{$(i+1)^2$};
	\draw (11.5,-4) node{$i$};
\draw (12.3,-4) node{\color{blue}$\ell$\color{black}};
\braidbox{13.1}{15.44}{-3.3}{-4.6}{};
	\draw (14.3,-4) node{$(\ell-1)^2$};
	\draw[dots] (16,-4)--(17.4,-4);
	\braidbox{17.8}{20.1}{-3.3}{-4.6}{};
	\draw (19,-4) node{$(i+2)^2$};
	\braidbox{20.8}{23.1}{-3.3}{-4.6}{};
	\draw (22,-4) node{$(i+1)^2$};	
	\draw[blue](-0.2,-3.2)--(-0.2,3.2);
	\draw[blue](12.3,3.2)--(12.3,-3.2);
	\draw(1,-3.2)--(1,3.2);
	\draw(13.5,3.2)--(13.5,-3.2);
	\draw(2.5,-3.2)--(2.5,3.2);
	\draw(15,3.2)--(15,-3.2);
	\draw(5.7,-3.2)--(5.7,3.2);
	\draw(18.2,3.2)--(18.2,-3.2);
	\draw(7.2,-3.2)--(7.2,3.2);
	\draw(19.7,3.2)--(19.7,-3.2);
	\draw(8.7,-3.2)--(8.7,2.7)--(10.2,3.2);
	\draw(9,3.2)--(11,0)--(10.2,-3.2);
	\draw(22.7,-3.2)--(22.7,3.2);
	\draw(21.3,3.2)--(21.3,-3.2);
	\draw(11.5,-3.2)--(9.3,0)--(11.5,3.2);
	\end{braid}
	}.
\end{align*}
Using the quadratic and dot-crossing relations completes the proof. 
\end{proof}

\section{The elements \texorpdfstring{$\Theta^i_j$}{}}
Throughout the subsection, we suppose that $i\in J\setminus\{0\}$. For $1\leq j\leq i$, we define the element $\Theta^i_j\in B_2$ to be 
$$
\Theta^i_j:=\index{$\Theta^i_j$}
\resizebox{116mm}{15mm}{
\begin{braid}\tikzset{baseline=-.3em}
	\draw (0,4) node{\color{blue}$\ell$\color{black}};
	\draw (0,-4) node{\color{blue}$\ell$\color{black}};
	\draw[blue](0,3.2)--(0,-3.2);
	\braidbox{0.6}{2.9}{3.3}{4.6}{};
	\draw (1.8,4) node{$(\ell-1)^2$};
	\braidbox{0.6}{2.9}{-4.7}{-3.4}{};
	\draw (1.8,-4) node{$(\ell-1)^2$};
	\draw(1,-3.2)--(1,3.2);
	\draw(2.6,-3.2)--(2.6,3.2);
%
	\draw[dots] (3.5,4)--(5,4);
	\braidbox{5.3}{7.6}{3.3}{4.6}{};
	\draw (6.5,4) node{$(i+1)^2$};
	\braidbox{5.3}{7.6}{-4.7}{-3.4}{};
	\draw (6.5,-4) node{$(i+1)^2$};
	\draw(5.8,-3.2)--(5.8,3.2);
	\draw(7.3,-3.2)--(7.3,3.2);
%
	\draw (8.3,4) node{$i$};
	\draw (8.3,-4) node{$i$};
	\draw(8.3,-3.2)--(8.3,3.2);
	\draw[dots] (9,4)--(10.3,4);
	\draw[dots] (9,-4)--(10.3,-4);
	\draw (11,4) node{$1$};
	\draw (11,-4) node{$1$};
	\redbraidbox{11.7}{12.6}{3.3}{4.6}{};
	\draw (12.1,4) node{$\color{red}0\,\,0\color{black}$};
	\redbraidbox{11.7}{12.6}{-3.3}{-4.6}{};
	\draw (12.1,-4) node{$\color{red}0\,\,0\color{black}$};
	\draw (13.3,4) node{$1$};
	\draw (13.3,-4) node{$1$};
	\draw(13.3,-3.2)--(13.3,3.2);
	\draw[dots] (13.9,4)--(15.3,4);
	\draw[dots] (13.9,-4)--(15.3,-4);
	\draw (16.3,4) node{$j-1$};
\draw (16.3,-4) node{$j-1$};
\draw(16.3,-3.2)--(16.3,3.2);
\draw[dots] (3.5,-4)--(5,-4);
	\draw(11,-3.2)--(11,3.2);
	\draw[red](11.8,-3.2)--(11.8,3.2);
	\draw[red](12.5,-3.2)--(12.5,3.2);
	\draw (17.6,4) node{$j$};
	\draw (17.6,-4) node{$j$};
	\draw (38.4,4) node{$j$};
	\draw (38.4,-4) node{$j$};
	\draw(17.6,-3.2)--(34.5,0)--(38.4,3.2);
	\draw(17.6,3.2)--(34.5,0)--(38.4,-3.2);

	\draw[dots] (18.1,4)--(19.5,4);
	\draw[dots] (18.1,-4)--(19.5,-4);
	\draw (19.9,4) node{$i$};
	\draw (19.9,-4) node{$i$};
	\draw (40.7,4) node{$i$};
	\draw (40.7,-4) node{$i$};
	\draw(19.9,-3.2)--(37,0)--(40.7,3.2);
	\draw(19.9,3.2)--(37,0)--(40.7,-3.2);
	\draw (20.8,4) node{\color{blue}$\ell$\color{black}};
	\draw (20.8,-4) node{\color{blue}$\ell$\color{black}};
	\draw[blue](20.8,-3.2)--(17,0)--(20.8,3.2);
	\braidbox{21.4}{23.7}{3.3}{4.6}{};
	\draw (22.6,4) node{$(\ell-1)^2$};
	\braidbox{21.4}{23.7}{-3.3}{-4.6}{};
	\draw (22.6,-4) node{$(\ell-1)^2$};
	\draw(22.2,-3.2)--(18.4,0)--(22.2,3.2);
	\draw(23.6,-3.2)--(19.8,0)--(23.6,3.2);
\draw[dots] (24.3,4)--(25.8,4);
\draw[dots] (24.3,-4)--(25.8,-4);
\braidbox{26.1}{28.4}{3.3}{4.6}{};
	\draw (27.3,4) node{$(i+1)^2$};
	\braidbox{26.1}{28.4}{-3.3}{-4.6}{};
	\draw (27.3,-4) node{$(i+1)^2$};
	\draw(26.4,-3.2)--(22.6,0)--(26.4,3.2);
	\draw(27.8,-3.2)--(24,0)--(27.8,3.2);
	\draw (29.1,4) node{$i$};
	\draw(29.1,-3.2)--(25.3,0)--(29.1,3.2);
	\draw[dots] (29.8,4)--(31.1,4);
	\draw (29.1,-4) node{$i$};
	\draw[dots] (29.8,-4)--(31.1,-4);
	\draw (31.8,4) node{$1$};
	\draw (31.8,-4) node{$1$};
	\draw(31.8,-3.2)--(28,0)--(31.8,3.2);
	\redbraidbox{32.5}{33.4}{3.3}{4.6}{};
	\draw (32.9,4) node{$\color{red}0\,\,0\color{black}$};
	\redbraidbox{32.5}{33.4}{-3.3}{-4.6}{};
	\draw (32.9,-4) node{$\color{red}0\,\,0\color{black}$};
	\draw[red](32.6,-3.2)--(28.8,0)--(32.6,3.2);
	\draw[red](33.3,-3.2)--(29.5,0)--(33.3,3.2);
	\draw (34.1,4) node{$1$};
	\draw (34.1,-4) node{$1$};
	\draw(34.1,-3.2)--(30.3,0)--(34.1,3.2);
	\draw[dots] (34.7,4)--(36.1,4);
	\draw[dots] (34.7,-4)--(36.1,-4);
	\draw (37.1,4) node{$j-1$};
	\draw (37.1,-4) node{$j-1$};
	\draw(37.1,-3.2)--(33.3,0)--(37.1,3.2);
		\draw[dots] (38.9,4)--(40.3,4);
	\draw[dots] (38.9,-4)--(40.3,-4);
\end{braid}.
}
$$

\begin{Lemma} \label{C100721}
Let $1\leq i<\ell$. Then in $B_2$ we have
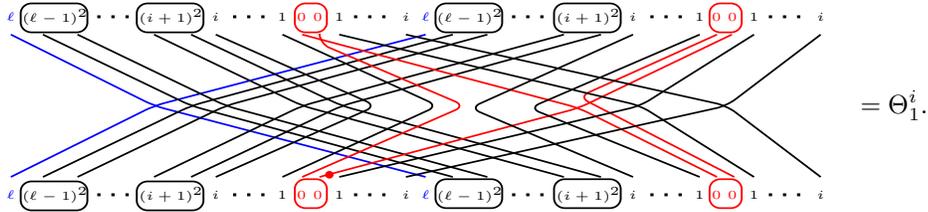
\begin{align*}
\resizebox{115mm}{15mm}{
\begin{braid}\tikzset{baseline=-.3em}
	\draw (0,4) node{\color{blue}$\ell$\color{black}};
	\braidbox{0.6}{2.9}{3.3}{4.6}{};
	\draw (1.8,4) node{$(\ell-1)^2$};
	\draw[dots] (3.5,4)--(5,4);
	\braidbox{5.3}{7.6}{3.3}{4.6}{};
	\draw (6.5,4) node{$(i+1)^2$};
	\draw (8.3,4) node{$i$};
	\draw[dots] (9,4)--(10.3,4);
	\draw[dots] (9,-4)--(10.3,-4);
	\draw (11,4) node{$1$};
	\draw (11,-4) node{$1$};
	\redbraidbox{11.7}{12.6}{3.3}{4.6}{};
	\draw (12.1,4) node{$\color{red}0\,\,0\color{black}$};
	\redbraidbox{11.7}{12.6}{-3.3}{-4.6}{};
	\draw (12.1,-4) node{$\color{red}0\,\,0\color{black}$};
	\draw (13.3,4) node{$1$};
	\draw (13.3,-4) node{$1$};
	\draw[dots] (13.9,4)--(15.3,4);
	\draw[dots] (13.9,-4)--(15.3,-4);
	\draw (16,4) node{$i$};
\draw (16,-4) node{$i$};
	\draw (0,-4) node{\color{blue}$\ell$\color{black}};
	\braidbox{0.6}{2.9}{-4.7}{-3.4}{};
	\draw (1.8,-4) node{$(\ell-1)^2$};
	\draw[dots] (3.5,-4)--(5,-4);
	\braidbox{5.3}{7.6}{-4.7}{-3.4}{};
	\draw (6.5,-4) node{$(i+1)^2$};
	\draw (8.3,-4) node{$i$};
\draw (16.8,4) node{\color{blue}$\ell$\color{black}};
\draw (16.8,-4) node{\color{blue}$\ell$\color{black}};
	\braidbox{17.4}{19.7}{3.3}{4.6}{};
	\draw (18.6,4) node{$(\ell-1)^2$};
	\braidbox{17.4}{19.7}{-3.3}{-4.6}{};
	\draw (18.6,-4) node{$(\ell-1)^2$};
	\draw[dots] (20.3,4)--(21.8,4);
	\draw[dots] (20.3,-4)--(21.8,-4);
	\braidbox{22.2}{24.5}{3.3}{4.6}{};
	\draw (23.4,4) node{$(i+1)^2$};
	\braidbox{22.2}{24.5}{-3.3}{-4.6}{};
	\draw (23.4,-4) node{$(i+1)^2$};
	\draw (25.2,4) node{$i$};
	\draw (25.2,-4) node{$i$};
	\draw[dots] (25.9,4)--(27.6,4);
	\draw[dots] (25.9,-4)--(27.6,-4);
	\draw (27.9,4) node{$1$};
	\draw (27.9,-4) node{$1$};
	\redbraidbox{28.5}{29.4}{3.3}{4.6}{};
	\draw (28.9,4) node{$\color{red}0\,\,0\color{black}$};
	\redbraidbox{28.5}{29.4}{-3.3}{-4.6}{};
	\draw (28.9,-4) node{$\color{red}0\,\,0\color{black}$};
	\draw (30.1,4) node{$1$};
	\draw (30.1,-4) node{$1$};
	\draw[dots] (30.7,4)--(32.1,4);
	\draw[dots] (30.7,-4)--(32.1,-4);
	\draw (32.8,4) node{$i$};
	\draw (32.8,-4) node{$i$};
	\draw[blue](0,-3.2)--(5.8,0)--(16.8,3.2);
	\draw[blue](0,3.2)--(5.8,0)--(16.8,-3.2);
	\draw(1.3,-3.2)--(7.3,0)--(17.9,3.2);
	\draw(1.3,3.2)--(7.3,0)--(17.9,-3.2);
	\draw(2.7,-3.2)--(8.7,0)--(19.3,3.2);
	\draw(2.7,3.2)--(8.7,0)--(19.3,-3.2);
	\draw(5.8,-3.2)--(11.8,0)--(22.4,3.2);
	\draw(5.8,3.2)--(11.8,0)--(22.4,-3.2);
	\draw(7.3,-3.2)--(13.4,0)--(23.9,3.2);
	\draw(7.3,3.2)--(13.4,0)--(23.9,-3.2);
	\draw(8.3,-3.2)--(14.8,0)--(8.3,3.2);
	\draw(11,-3.2)--(17.3,0)--(11,3.2);
	\draw[red](11.8,-3.2)--(18.4,0)--(12.5,2.8)--(12.5,3.2);
	\draw[red](12.5,-3.2)--(23,-0.1)--(29.2,3.2);
	\draw[red](28.6,-3.2)--(23,-0.1)--(11.8,3.2);
	\draw[red](29.3,-3.2)--(23,0.4)--(28.6,3.2);
	\draw(13.3,-3.2)--(25.5,0)--(30.1,3.2);
	\draw(13.3,3.2)--(25.5,0)--(30.1,-3.2);
	\draw(32.8,-3.2)--(29,0)--(16,3.2);
	\draw(32.8,3.2)--(29,0)--(16,-3.2);
	\draw(25.2,-3.2)--(18.6,-0.1)--(25.2,3.2);
	\draw(27.9,-3.2)--(21,-0.1)--(27.9,3.2);
	\reddot(12.9,-3.1);
	\end{braid}
	}
	=\Theta^i_1.
\end{align*}
\end{Lemma}
\begin{proof}
We apply Lemma~\ref{L100721_0} to the element in the left hand side. Note that the elements of $\I$ in Lemma~\ref{L100721_0} create non-cuspidal words, so the corresponding terms are zero, and we get 
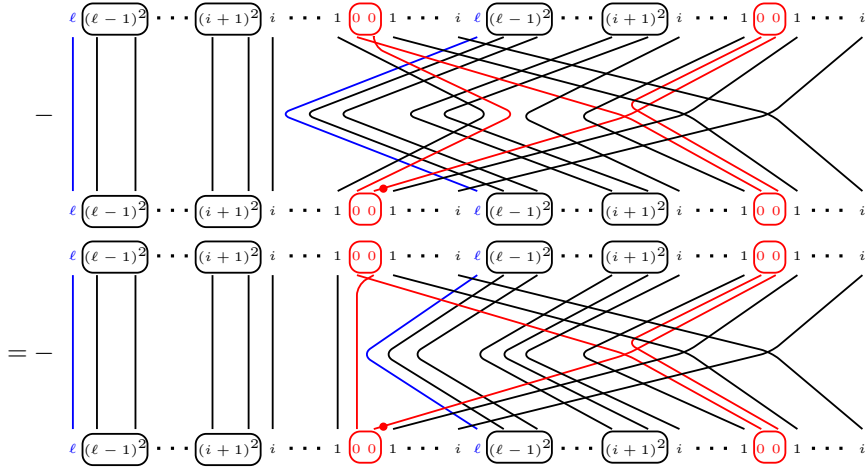
\begin{align*}
&-\begin{braid}\tikzset{baseline=-.3em}
	\draw (0,4) node{\color{blue}$\ell$\color{black}};
	\braidbox{0.6}{2.9}{3.3}{4.6}{};
	\draw (1.8,4) node{$(\ell-1)^2$};
	\draw[dots] (3.5,4)--(5,4);
	\braidbox{5.3}{7.6}{3.3}{4.6}{};
	\draw (6.5,4) node{$(i+1)^2$};
	\draw (8.3,4) node{$i$};
	\draw[dots] (9,4)--(10.3,4);
	\draw[dots] (9,-4)--(10.3,-4);
	\draw (11,4) node{$1$};
	\draw (11,-4) node{$1$};
	\redbraidbox{11.7}{12.6}{3.3}{4.6}{};
	\draw (12.1,4) node{$\color{red}0\,\,0\color{black}$};
	\redbraidbox{11.7}{12.6}{-3.3}{-4.6}{};
	\draw (12.1,-4) node{$\color{red}0\,\,0\color{black}$};
	\draw (13.3,4) node{$1$};
	\draw (13.3,-4) node{$1$};
	\draw[dots] (13.9,4)--(15.3,4);
	\draw[dots] (13.9,-4)--(15.3,-4);
	\draw (16,4) node{$i$};
\draw (16,-4) node{$i$};
	\draw (0,-4) node{\color{blue}$\ell$\color{black}};
	\braidbox{0.6}{2.9}{-4.7}{-3.4}{};
	\draw (1.8,-4) node{$(\ell-1)^2$};
	\draw[dots] (3.5,-4)--(5,-4);
	\braidbox{5.3}{7.6}{-4.7}{-3.4}{};
	\draw (6.5,-4) node{$(i+1)^2$};
	\draw (8.3,-4) node{$i$};
	\draw (16.8,4) node{\color{blue}$\ell$\color{black}};
\draw (16.8,-4) node{\color{blue}$\ell$\color{black}};
	\braidbox{17.4}{19.7}{3.3}{4.6}{};
	\draw (18.6,4) node{$(\ell-1)^2$};
	\braidbox{17.4}{19.7}{-3.3}{-4.6}{};
	\draw (18.6,-4) node{$(\ell-1)^2$};
	\draw[dots] (20.3,4)--(21.8,4);
	\draw[dots] (20.3,-4)--(21.8,-4);
	\braidbox{22.2}{24.5}{3.3}{4.6}{};
	\draw (23.4,4) node{$(i+1)^2$};
	\braidbox{22.2}{24.5}{-3.3}{-4.6}{};
	\draw (23.4,-4) node{$(i+1)^2$};
	\draw (25.2,4) node{$i$};
	\draw (25.2,-4) node{$i$};
	\draw[dots] (25.9,4)--(27.6,4);
	\draw[dots] (25.9,-4)--(27.6,-4);
	\draw (27.9,4) node{$1$};
	\draw (27.9,-4) node{$1$};
	\redbraidbox{28.5}{29.4}{3.3}{4.6}{};
	\draw (28.9,4) node{$\color{red}0\,\,0\color{black}$};
	\redbraidbox{28.5}{29.4}{-3.3}{-4.6}{};
	\draw (28.9,-4) node{$\color{red}0\,\,0\color{black}$};
	\draw (30.1,4) node{$1$};
	\draw (30.1,-4) node{$1$};
	\draw[dots] (30.7,4)--(32.1,4);
	\draw[dots] (30.7,-4)--(32.1,-4);
	\draw (32.8,4) node{$i$};
	\draw (32.8,-4) node{$i$};
	\draw[blue](16.8,-3.2)--(8.6,0)--(16.8,3.2);
	\draw[blue](0,3.2)--(0,-3.2);
	\draw(1,-3.2)--(1,3.2);
	\draw(17.9,3.2)--(9.6,0)--(17.9,-3.2);
	\draw(2.6,-3.2)--(2.6,3.2);
	\draw(19.3,3.2)--(11,0)--(19.3,-3.2);
	\draw(5.8,-3.2)--(5.8,3.2);
	\draw(22.4,3.2)--(13.8,0)--(22.4,-3.2);
	\draw(7.3,-3.2)--(7.3,3.2);
	\draw(23.9,3.2)--(15.2,0)--(23.9,-3.2);
	\draw(8.3,-3.2)--(8.3,3.2);
	\draw(11,-3.2)--(17.3,0)--(11,3.2);
	\draw[red](11.8,-3.2)--(18.4,0)--(12.5,2.8)--(12.5,3.2);
	\draw[red](12.5,-3.2)--(23,-0.1)--(29.2,3.2);
	\draw[red](28.6,-3.2)--(23,-0.1)--(11.8,3.2);
	\draw[red](29.3,-3.2)--(23,0.4)--(28.6,3.2);
	\draw(13.3,-3.2)--(25.5,0)--(30.1,3.2);
	\draw(13.3,3.2)--(25.5,0)--(30.1,-3.2);
	\draw(32.8,-3.2)--(29,0)--(16,3.2);
	\draw(32.8,3.2)--(29,0)--(16,-3.2);
	\draw(25.2,-3.2)--(18.6,-0.1)--(25.2,3.2);
	\draw(27.9,-3.2)--(21,-0.1)--(27.9,3.2);
	\reddot(12.9,-3.1);
	\end{braid}
	\\
	=&
	-\begin{braid}\tikzset{baseline=-.3em}
	\draw (0,4) node{\color{blue}$\ell$\color{black}};
	\braidbox{0.6}{2.9}{3.3}{4.6}{};
	\draw (1.8,4) node{$(\ell-1)^2$};
	\draw[dots] (3.5,4)--(5,4);
	\braidbox{5.3}{7.6}{3.3}{4.6}{};
	\draw (6.5,4) node{$(i+1)^2$};
	\draw (8.3,4) node{$i$};
	\draw[dots] (9,4)--(10.3,4);
	\draw[dots] (9,-4)--(10.3,-4);
	\draw (11,4) node{$1$};
	\draw (11,-4) node{$1$};
	\redbraidbox{11.7}{12.6}{3.3}{4.6}{};
	\draw (12.1,4) node{$\color{red}0\,\,0\color{black}$};
	\redbraidbox{11.7}{12.6}{-3.3}{-4.6}{};
	\draw (12.1,-4) node{$\color{red}0\,\,0\color{black}$};
	\draw (13.3,4) node{$1$};
	\draw (13.3,-4) node{$1$};
	\draw[dots] (13.9,4)--(15.3,4);
	\draw[dots] (13.9,-4)--(15.3,-4);
	\draw (16,4) node{$i$};
\draw (16,-4) node{$i$};
	\draw (0,-4) node{\color{blue}$\ell$\color{black}};
	\braidbox{0.6}{2.9}{-4.7}{-3.4}{};
	\draw (1.8,-4) node{$(\ell-1)^2$};
	\draw[dots] (3.5,-4)--(5,-4);
	\braidbox{5.3}{7.6}{-4.7}{-3.4}{};
	\draw (6.5,-4) node{$(i+1)^2$};
	\draw (8.3,-4) node{$i$};
	\draw (16.8,4) node{\color{blue}$\ell$\color{black}};
\draw (16.8,-4) node{\color{blue}$\ell$\color{black}};
	\braidbox{17.4}{19.7}{3.3}{4.6}{};
	\draw (18.6,4) node{$(\ell-1)^2$};
	\braidbox{17.4}{19.7}{-3.3}{-4.6}{};
	\draw (18.6,-4) node{$(\ell-1)^2$};
	\draw[dots] (20.3,4)--(21.8,4);
	\draw[dots] (20.3,-4)--(21.8,-4);
	\braidbox{22.2}{24.5}{3.3}{4.6}{};
	\draw (23.4,4) node{$(i+1)^2$};
	\braidbox{22.2}{24.5}{-3.3}{-4.6}{};
	\draw (23.4,-4) node{$(i+1)^2$};
	\draw (25.2,4) node{$i$};
	\draw (25.2,-4) node{$i$};
	\draw[dots] (25.9,4)--(27.6,4);
	\draw[dots] (25.9,-4)--(27.6,-4);
	\draw (27.9,4) node{$1$};
	\draw (27.9,-4) node{$1$};
	\redbraidbox{28.5}{29.4}{3.3}{4.6}{};
	\draw (28.9,4) node{$\color{red}0\,\,0\color{black}$};
	\redbraidbox{28.5}{29.4}{-3.3}{-4.6}{};
	\draw (28.9,-4) node{$\color{red}0\,\,0\color{black}$};
	\draw (30.1,4) node{$1$};
	\draw (30.1,-4) node{$1$};
	\draw[dots] (30.7,4)--(32.1,4);
	\draw[dots] (30.7,-4)--(32.1,-4);
	\draw (32.8,4) node{$i$};
	\draw (32.8,-4) node{$i$};
	\draw[blue](16.8,-3.2)--(12,-.10)--(16.8,3.2);
	\draw[blue](0,3.2)--(0,-3.2);
	\draw(1,-3.2)--(1,3.2);
	\draw(17.9,3.2)--(12.9,-0.1)--(17.9,-3.2);
	\draw(2.6,-3.2)--(2.6,3.2);
	\draw(19.3,3.2)--(14.1,-0.1)--(19.3,-3.2);
	\draw(5.8,-3.2)--(5.8,3.2);
	\draw(22.4,3.2)--(16.7,-.1)--(22.4,-3.2);
	\draw(7.3,-3.2)--(7.3,3.2);
	\draw(23.9,3.2)--(17.7,-0.1)--(23.9,-3.2);
	\draw(8.3,-3.2)--(8.3,3.2);
	\draw(11,-3.2)--(11,3.2);
	\draw[red](11.8,-3.2)--(11.8,2.8)--(12.5,3.2);
	\draw[red](12.5,-3.2)--(23,-0.1)--(29.2,3.2);
	\draw[red](28.6,-3.2)--(23,-0.1)--(11.8,3.2);
	\draw[red](29.3,-3.2)--(23,0.4)--(28.6,3.2);
	\draw(13.3,-3.2)--(25.5,0)--(30.1,3.2);
	\draw(13.3,3.2)--(25.5,0)--(30.1,-3.2);
	\draw(32.8,-3.2)--(29,0)--(16,3.2);
	\draw(32.8,3.2)--(29,0)--(16,-3.2);
	\draw(25.2,-3.2)--(18.6,-0.1)--(25.2,3.2);
	\draw(27.9,-3.2)--(21,-0.1)--(27.9,3.2);
	\reddot(12.9,-3.1);
	\end{braid}
\end{align*}
Applying the braid relation 
$$
\begin{braid}\tikzset{baseline=2mm}
  \draw[red] (0,0)node[below]{$0$}--(2,2);
  \draw[red] (2,0)node[below]{$0$}--(0,2);
  \draw (1,0)node[below]{$1$}--(0,1)--(1,2);
\end{braid}
=\begin{braid}\tikzset{baseline=2mm}
  \draw[red] (0,0)node[below]{$0$}--(2,2);
  \draw[red] (2,0)node[below]{$0$}--(0,2);
  \draw (1,0)node[below]{$1$}--(2,1)--(1,2);
\end{braid}
+
\begin{braid}\tikzset{baseline=2mm}
  \draw[red] (0,0)node[below]{$0$}--(0,2);
  \draw[red] (2,0)node[below]{$0$}--(2,2);
  \draw (1,0)node[below]{$1$}--(1,2);
  \reddot(0,1);
\end{braid}
-\begin{braid}\tikzset{baseline=2mm}
  \draw[red] (0,0)node[below]{$0$}--(0,2);
  \draw[red] (2,0)node[below]{$0$}--(2,2);
  \draw (1,0)node[below]{$1$}--(1,2);
  \reddot(2,1);
\end{braid}
$$
to the last diagram we get three summands. The first summand is zero since the word starting with $\ell(\ell-1)^2\cdots(i+1)^2i\cdots20$ is not cuspidal by Lemma~\ref{LCuspExpl}. The second summand is also zero since 
$\begin{braid}\tikzset{baseline=.5em}
\redbraidbox{-.3}{.6}{-.7}{.6}{};
	\draw (0.1,0) node{$\color{red}0\,\,0\color{black}$};
		\draw[red](-0.2,0.6)--(.5,1.3);
\draw[red](0.6,0.6)--(-.2,1.3);
	\end{braid}=0$. 
The third summand is 
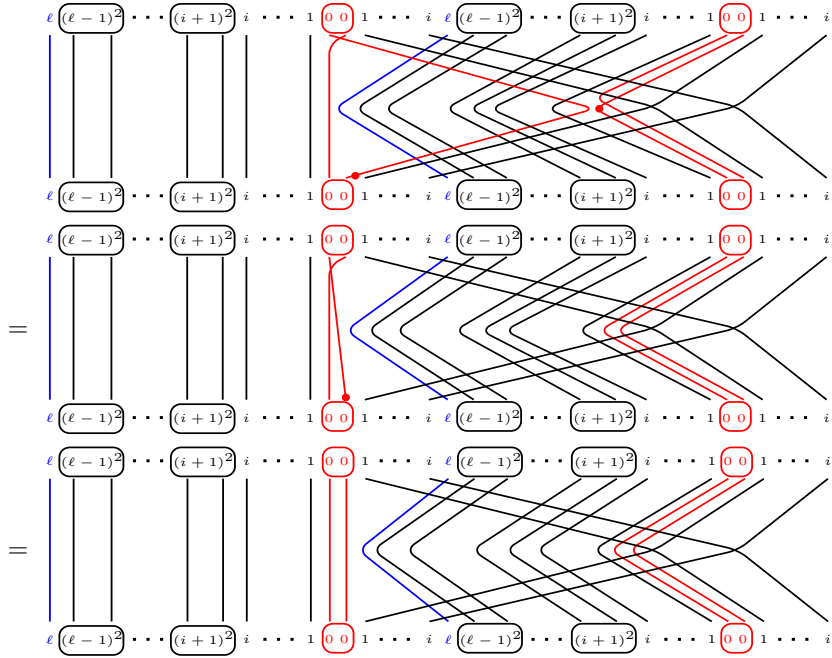
\begin{align*}
&
\resizebox{110mm}{15mm}{
\begin{braid}\tikzset{baseline=-.3em}
	\draw (0,4) node{\color{blue}$\ell$\color{black}};
	\braidbox{0.6}{2.9}{3.3}{4.6}{};
	\draw (1.8,4) node{$(\ell-1)^2$};
	\draw[dots] (3.5,4)--(5,4);
	\braidbox{5.3}{7.6}{3.3}{4.6}{};
	\draw (6.5,4) node{$(i+1)^2$};
	\draw (8.3,4) node{$i$};
	\draw[dots] (9,4)--(10.3,4);
	\draw[dots] (9,-4)--(10.3,-4);
	\draw (11,4) node{$1$};
	\draw (11,-4) node{$1$};
	\redbraidbox{11.7}{12.6}{3.3}{4.6}{};
	\draw (12.1,4) node{$\color{red}0\,\,0\color{black}$};
	\redbraidbox{11.7}{12.6}{-3.3}{-4.6}{};
	\draw (12.1,-4) node{$\color{red}0\,\,0\color{black}$};
	\draw (13.3,4) node{$1$};
	\draw (13.3,-4) node{$1$};
	\draw[dots] (13.9,4)--(15.3,4);
	\draw[dots] (13.9,-4)--(15.3,-4);
	\draw (16,4) node{$i$};
\draw (16,-4) node{$i$};
	\draw (0,-4) node{\color{blue}$\ell$\color{black}};
	\braidbox{0.6}{2.9}{-4.7}{-3.4}{};
	\draw (1.8,-4) node{$(\ell-1)^2$};
	\draw[dots] (3.5,-4)--(5,-4);
	\braidbox{5.3}{7.6}{-4.7}{-3.4}{};
	\draw (6.5,-4) node{$(i+1)^2$};
	\draw (8.3,-4) node{$i$};
	\draw (16.8,4) node{\color{blue}$\ell$\color{black}};
\draw (16.8,-4) node{\color{blue}$\ell$\color{black}};
	\braidbox{17.4}{19.7}{3.3}{4.6}{};
	\draw (18.6,4) node{$(\ell-1)^2$};
	\braidbox{17.4}{19.7}{-3.3}{-4.6}{};
	\draw (18.6,-4) node{$(\ell-1)^2$};
	\draw[dots] (20.3,4)--(21.8,4);
	\draw[dots] (20.3,-4)--(21.8,-4);
	\braidbox{22.2}{24.5}{3.3}{4.6}{};
	\draw (23.4,4) node{$(i+1)^2$};
	\braidbox{22.2}{24.5}{-3.3}{-4.6}{};
	\draw (23.4,-4) node{$(i+1)^2$};
	\draw (25.2,4) node{$i$};
	\draw (25.2,-4) node{$i$};
	\draw[dots] (25.9,4)--(27.6,4);
	\draw[dots] (25.9,-4)--(27.6,-4);
	\draw (27.9,4) node{$1$};
	\draw (27.9,-4) node{$1$};
	\redbraidbox{28.5}{29.4}{3.3}{4.6}{};
	\draw (28.9,4) node{$\color{red}0\,\,0\color{black}$};
	\redbraidbox{28.5}{29.4}{-3.3}{-4.6}{};
	\draw (28.9,-4) node{$\color{red}0\,\,0\color{black}$};
	\draw (30.1,4) node{$1$};
	\draw (30.1,-4) node{$1$};
	\draw[dots] (30.7,4)--(32.1,4);
	\draw[dots] (30.7,-4)--(32.1,-4);
	\draw (32.8,4) node{$i$};
	\draw (32.8,-4) node{$i$};
	\draw[blue](16.8,-3.2)--(12,-.10)--(16.8,3.2);
	\draw[blue](0,3.2)--(0,-3.2);
	\draw(1,-3.2)--(1,3.2);
	\draw(17.9,3.2)--(12.9,-0.1)--(17.9,-3.2);
	\draw(2.6,-3.2)--(2.6,3.2);
	\draw(19.3,3.2)--(14.1,-0.1)--(19.3,-3.2);
	\draw(5.8,-3.2)--(5.8,3.2);
	\draw(22.4,3.2)--(16.7,-.1)--(22.4,-3.2);
	\draw(7.3,-3.2)--(7.3,3.2);
	\draw(23.9,3.2)--(17.7,-0.1)--(23.9,-3.2);
	\draw(8.3,-3.2)--(8.3,3.2);
	\draw(11,-3.2)--(11,3.2);
	\draw[red](11.8,-3.2)--(11.8,2.8)--(12.5,3.2);
	\draw[red](12.5,-3.2)--(23,-0.1)--(11.8,3.2);
	\draw[red](28.6,-3.2)--(23,-0.1)--(29.3,3.2);
	\reddot(23.2,-.1);
	\draw[red](29.3,-3.2)--(23,0.4)--(28.6,3.2);
	\draw(13.3,-3.2)--(25.5,0)--(30.1,3.2);
	\draw(13.3,3.2)--(25.5,0)--(30.1,-3.2);
	\draw(32.8,-3.2)--(29,0)--(16,3.2);
	\draw(32.8,3.2)--(29,0)--(16,-3.2);
	\draw(25.2,-3.2)--(18.6,-0.1)--(25.2,3.2);
	\draw(27.9,-3.2)--(21,-0.1)--(27.9,3.2);
	\reddot(12.9,-3.1);
	\end{braid}
	}
	\\
	=&
	\resizebox{110mm}{15mm}{
	\begin{braid}\tikzset{baseline=-.3em}
	\draw (0,4) node{\color{blue}$\ell$\color{black}};
	\braidbox{0.6}{2.9}{3.3}{4.6}{};
	\draw (1.8,4) node{$(\ell-1)^2$};
	\draw[dots] (3.5,4)--(5,4);
	\braidbox{5.3}{7.6}{3.3}{4.6}{};
	\draw (6.5,4) node{$(i+1)^2$};
	\draw (8.3,4) node{$i$};
	\draw[dots] (9,4)--(10.3,4);
	\draw[dots] (9,-4)--(10.3,-4);
	\draw (11,4) node{$1$};
	\draw (11,-4) node{$1$};
	\redbraidbox{11.7}{12.6}{3.3}{4.6}{};
	\draw (12.1,4) node{$\color{red}0\,\,0\color{black}$};
	\redbraidbox{11.7}{12.6}{-3.3}{-4.6}{};
	\draw (12.1,-4) node{$\color{red}0\,\,0\color{black}$};
	\draw (13.3,4) node{$1$};
	\draw (13.3,-4) node{$1$};
	\draw[dots] (13.9,4)--(15.3,4);
	\draw[dots] (13.9,-4)--(15.3,-4);
	\draw (16,4) node{$i$};
\draw (16,-4) node{$i$};
	\draw (0,-4) node{\color{blue}$\ell$\color{black}};
	\braidbox{0.6}{2.9}{-4.7}{-3.4}{};
	\draw (1.8,-4) node{$(\ell-1)^2$};
	\draw[dots] (3.5,-4)--(5,-4);
	\braidbox{5.3}{7.6}{-4.7}{-3.4}{};
	\draw (6.5,-4) node{$(i+1)^2$};
	\draw (8.3,-4) node{$i$};
	\draw (16.8,4) node{\color{blue}$\ell$\color{black}};
\draw (16.8,-4) node{\color{blue}$\ell$\color{black}};
	\braidbox{17.4}{19.7}{3.3}{4.6}{};
	\draw (18.6,4) node{$(\ell-1)^2$};
	\braidbox{17.4}{19.7}{-3.3}{-4.6}{};
	\draw (18.6,-4) node{$(\ell-1)^2$};
	\draw[dots] (20.3,4)--(21.8,4);
	\draw[dots] (20.3,-4)--(21.8,-4);
	\braidbox{22.2}{24.5}{3.3}{4.6}{};
	\draw (23.4,4) node{$(i+1)^2$};
	\braidbox{22.2}{24.5}{-3.3}{-4.6}{};
	\draw (23.4,-4) node{$(i+1)^2$};
	\draw (25.2,4) node{$i$};
	\draw (25.2,-4) node{$i$};
	\draw[dots] (25.9,4)--(27.6,4);
	\draw[dots] (25.9,-4)--(27.6,-4);
	\draw (27.9,4) node{$1$};
	\draw (27.9,-4) node{$1$};
	\redbraidbox{28.5}{29.4}{3.3}{4.6}{};
	\draw (28.9,4) node{$\color{red}0\,\,0\color{black}$};
	\redbraidbox{28.5}{29.4}{-3.3}{-4.6}{};
	\draw (28.9,-4) node{$\color{red}0\,\,0\color{black}$};
	\draw (30.1,4) node{$1$};
	\draw (30.1,-4) node{$1$};
	\draw[dots] (30.7,4)--(32.1,4);
	\draw[dots] (30.7,-4)--(32.1,-4);
	\draw (32.8,4) node{$i$};
	\draw (32.8,-4) node{$i$};
	\draw[blue](16.8,-3.2)--(12.5,-.1)--(16.8,3.2);
	\draw[blue](0,3.2)--(0,-3.2);
	\draw(1,-3.2)--(1,3.2);
	\draw(17.9,3.2)--(13.4,-0.1)--(17.9,-3.2);
	\draw(2.6,-3.2)--(2.6,3.2);
	\draw(19.3,3.2)--(14.6,-0.1)--(19.3,-3.2);
	\draw(5.8,-3.2)--(5.8,3.2);
	\draw(22.4,3.2)--(17.1,-.1)--(22.4,-3.2);
	\draw(7.3,-3.2)--(7.3,3.2);
	\draw(23.9,3.2)--(18.2,-0.1)--(23.9,-3.2);
	\draw(8.3,-3.2)--(8.3,3.2);
	\draw(11,-3.2)--(11,3.2);
	\draw[red](11.8,-3.2)--(11.8,2.8)--(12.5,3.2);
	\draw[red](12.5,-3.2)--(11.8,3.2);
	\draw[red](28.6,-3.2)--(23.2,-0.1)--(28.6,3.2);
	\draw[red](29.3,-3.2)--(23.9,-0.1)--(29.3,3.2);
	\draw(13.3,-3.2)--(25.5,0)--(30.1,3.2);
	\draw(13.3,3.2)--(25.5,0)--(30.1,-3.2);
	\draw(32.8,-3.2)--(29,0)--(16,3.2);
	\draw(32.8,3.2)--(29,0)--(16,-3.2);
	\draw(25.2,-3.2)--(19.2,-0.1)--(25.2,3.2);
	\draw(27.9,-3.2)--(22.3,-0.1)--(27.9,3.2);
	\reddot(12.5,-3.1);
	\end{braid}
	}
	\\
	=&
	\resizebox{110mm}{15mm}{
	\begin{braid}\tikzset{baseline=-.3em}
	\draw (0,4) node{\color{blue}$\ell$\color{black}};
	\draw (0,-4) node{\color{blue}$\ell$\color{black}};
	\draw (16.8,4) node{\color{blue}$\ell$\color{black}};
\draw (16.8,-4) node{\color{blue}$\ell$\color{black}};
\draw[blue](16.8,-3.2)--(13,0)--(16.8,3.2);
	\draw[blue](0,3.2)--(0,-3.2);
	\braidbox{0.6}{2.9}{3.3}{4.6}{};
	\draw (1.8,4) node{$(\ell-1)^2$};
	\braidbox{0.6}{2.9}{-4.7}{-3.4}{};
	\draw (1.8,-4) node{$(\ell-1)^2$};
	\braidbox{17.4}{19.7}{3.3}{4.6}{};
	\draw (18.6,4) node{$(\ell-1)^2$};
	\braidbox{17.4}{19.7}{-3.3}{-4.6}{};
	\draw (18.6,-4) node{$(\ell-1)^2$};
	\draw(1,-3.2)--(1,3.2);
	\draw(17.9,3.2)--(13.6,0)--(17.9,-3.2);
	\draw(2.6,-3.2)--(2.6,3.2);
	\draw(19.3,3.2)--(15,0)--(19.3,-3.2);
	\draw[dots] (3.5,4)--(5,4);
	\braidbox{5.3}{7.6}{3.3}{4.6}{};
	\draw (6.5,4) node{$(i+1)^2$};
	\braidbox{5.3}{7.6}{-4.7}{-3.4}{};
	\draw (6.5,-4) node{$(i+1)^2$};
	\braidbox{22.2}{24.5}{3.3}{4.6}{};
	\draw (23.4,4) node{$(i+1)^2$};
	\braidbox{22.2}{24.5}{-3.3}{-4.6}{};
	\draw (23.4,-4) node{$(i+1)^2$};
	\draw(5.8,-3.2)--(5.8,3.2);
	\draw(22.4,3.2)--(17.8,0)--(22.4,-3.2);
	\draw(7.3,-3.2)--(7.3,3.2);
	\draw(23.9,3.2)--(19.2,0)--(23.9,-3.2);
	\draw (8.3,4) node{$i$};
	\draw[dots] (9,4)--(10.3,4);
	\draw[dots] (9,-4)--(10.3,-4);
	\draw (11,4) node{$1$};
	\draw (11,-4) node{$1$};
	\redbraidbox{11.7}{12.6}{3.3}{4.6}{};
	\draw (12.1,4) node{$\color{red}0\,\,0\color{black}$};
	\redbraidbox{11.7}{12.6}{-3.3}{-4.6}{};
	\draw (12.1,-4) node{$\color{red}0\,\,0\color{black}$};
	\draw (13.3,4) node{$1$};
	\draw (13.3,-4) node{$1$};
	\draw[dots] (13.9,4)--(15.3,4);
	\draw[dots] (13.9,-4)--(15.3,-4);
	\draw (16,4) node{$i$};
\draw (16,-4) node{$i$};
			\draw[dots] (3.5,-4)--(5,-4);
	\draw (8.3,-4) node{$i$};
	\draw[dots] (20.3,4)--(21.8,4);
	\draw[dots] (20.3,-4)--(21.8,-4);
	\draw (25.2,4) node{$i$};
	\draw (25.2,-4) node{$i$};
	\draw[dots] (25.9,4)--(27.6,4);
	\draw[dots] (25.9,-4)--(27.6,-4);
	\draw (27.9,4) node{$1$};
	\draw (27.9,-4) node{$1$};
	\redbraidbox{28.5}{29.4}{3.3}{4.6}{};
	\draw (28.9,4) node{$\color{red}0\,\,0\color{black}$};
	\redbraidbox{28.5}{29.4}{-3.3}{-4.6}{};
	\draw (28.9,-4) node{$\color{red}0\,\,0\color{black}$};
	\draw (30.1,4) node{$1$};
	\draw (30.1,-4) node{$1$};
	\draw[dots] (30.7,4)--(32.1,4);
	\draw[dots] (30.7,-4)--(32.1,-4);
	\draw (32.8,4) node{$i$};
	\draw (32.8,-4) node{$i$};
	\draw(8.3,-3.2)--(8.3,3.2);
	\draw(11,-3.2)--(11,3.2);
	\draw[red](11.8,-3.2)--(11.8,3.2);
	\draw[red](12.5,-3.2)--(12.5,3.2);
	\draw[red](28.6,-3.2)--(23.6,0)--(28.6,3.2);
	\draw[red](29.3,-3.2)--(24.4,0)--(29.3,3.2);
	\draw(13.3,-3.2)--(25.5,0)--(30.1,3.2);
	\draw(13.3,3.2)--(25.5,0)--(30.1,-3.2);
	\draw(32.8,-3.2)--(29,0)--(16,3.2);
	\draw(32.8,3.2)--(29,0)--(16,-3.2);
	\draw(25.2,-3.2)--(20.5,0)--(25.2,3.2);
	\draw(27.9,-3.2)--(22.9,0)--(27.9,3.2);
	\end{braid},}
\end{align*}
as required. 
\end{proof}

\begin{Lemma} \label{L110621} 
Let $1\leq i\leq \ell-1$. Then 
$\Theta^i_i=\ga^{i,i}(-z_1+z_2+(-1)^i(c_1+c_2))\ga^{i,i}$. 
\end{Lemma}
\begin{proof}
We provide details for the generic case $1<i<\ell-1$, the special cases $i=1$ and $i=\ell-1$, being similar, are left to the reader. To simplify the notation we drop the first part of the diagram with vertical strings, since it is not going to change throughout the computation. Using braid relations (\ref{R7}) and quadratic relations (\ref{R6}), we get that $\Theta^i_i$ equals 
\begin{align*}
&
\resizebox{59mm}{15mm}{
\begin{braid}\tikzset{baseline=-.3em}
	\draw (19.5,4) node{$i$};
		\draw (19.5,-4) node{$i$};
	\draw (20.3,4) node{\color{blue}$\ell$\color{black}};
	\braidbox{20.9}{23.2}{3.3}{4.6}{};
	\draw (22.1,4) node{$(\ell-1)^2$};
	\draw[dots] (23.8,4)--(25.3,4);
	\braidbox{25.6}{27.9}{3.3}{4.6}{};
	\draw (26.8,4) node{$(i+1)^2$};
	\draw (28.6,4) node{$i$};
	\draw (29.8,4) node{$i-1$};
	\draw[dots] (31.1,4)--(32.8,4);
	\draw (33,-4) node{$1$};
	\redbraidbox{33.7}{34.6}{3.3}{4.6}{};
	\draw (34.1,4) node{$\color{red}0\,\,0\color{black}$};
	\draw (35.3,4) node{$1$};
	\draw[dots] (35.9,4)--(37.3,4);
	\draw (38.5,4) node{$i-1$};
	\draw (39.8,4) node{$i$};
	\draw (20.3,-4) node{\color{blue}$\ell$\color{black}};
	\braidbox{20.9}{23.2}{-4.7}{-3.4}{};
	\draw (22.1,-4) node{$(\ell-1)^2$};
	\draw[dots] (23.8,-4)--(25.1,-4);
	\braidbox{25.6}{27.9}{-4.7}{-3.4}{};
	\draw (26.8,-4) node{$(i+1)^2$};
	\draw (28.6,-4) node{$i$};
	\draw (29.8,-4) node{$i-1$};
	\draw[dots] (31.1,-4)--(32.8,-4);
	\draw (33,-4) node{$1$};
	\redbraidbox{33.7}{34.6}{-4.7}{-3.4}{};
	\draw (34.1,-4) node{$\color{red}0\,\,0\color{black}$};
	\draw (35.3,-4) node{$1$};
	\draw[dots] (35.9,-4)--(37.3,-4);
	\draw (38.5,-4) node{$i-1$};
	\draw (39.8,-4) node{$i$};
	\draw(39.8,-3.2)--(36,0)--(19.5,3.2);
	\draw(39.8,3.2)--(36,0)--(19.5,-3.2);
	\draw[blue](20.3,-3.2)--(18.5,0)--(20.3,3.2);
	\draw(21.3,-3.2)--(19.5,0)--(21.3,3.2);
	\draw(22.3,-3.2)--(20.5,0)--(22.3,3.2);
	\draw(26,-3.2)--(24.1,0)--(26,3.2);
	\draw(27,-3.2)--(25.1,0)--(27,3.2);
	\draw(28.6,-3.2)--(26.7,0)--(28.6,3.2);
	\draw(29.8,-3.2)--(27.9,0)--(29.8,3.2);
	\draw(38.5,-3.2)--(35,0)--(38.5,3.2);
	\draw(35.3,-3.2)--(32.5,0)--(35.3,3.2);
	\draw[red](34.5,-3.2)--(31.8,0)--(34.5,3.2);
	\draw[red](33.8,-3.2)--(31.1,0)--(33.8,3.2);
	\draw(33.1,-3.2)--(30.4,0)--(33.1,3.2);
	\end{braid}
	}
	\\
	=&
	\resizebox{59mm}{15mm}{
	\begin{braid}\tikzset{baseline=-.3em}
	\draw (19.5,4) node{$i$};
		\draw (19.5,-4) node{$i$};
	\draw (20.3,4) node{\color{blue}$\ell$\color{black}};
	\braidbox{20.9}{23.2}{3.3}{4.6}{};
	\draw (22.1,4) node{$(\ell-1)^2$};
	\draw[dots] (23.8,4)--(25.3,4);
	\braidbox{25.6}{27.9}{3.3}{4.6}{};
	\draw (26.8,4) node{$(i+1)^2$};
	\draw (28.6,4) node{$i$};
	\draw (29.8,4) node{$i-1$};
	\draw[dots] (31.1,4)--(32.8,4);
	\draw (33,4) node{$1$};
	\redbraidbox{33.7}{34.6}{3.3}{4.6}{};
	\draw (34.1,4) node{$\color{red}0\,\,0\color{black}$};
	\draw (35.3,4) node{$1$};
	\draw[dots] (35.9,4)--(37.3,4);
	\draw (38.5,4) node{$i-1$};
	\draw (39.8,4) node{$i$};
	\draw (20.3,-4) node{\color{blue}$\ell$\color{black}};
	\braidbox{20.9}{23.2}{-4.7}{-3.4}{};
	\draw (22.1,-4) node{$(\ell-1)^2$};
	\draw[dots] (23.8,-4)--(25.1,-4);
	\braidbox{25.6}{27.9}{-4.7}{-3.4}{};
	\draw (26.8,-4) node{$(i+1)^2$};
	\draw (28.6,-4) node{$i$};
	\draw (29.8,-4) node{$i-1$};
	\draw[dots] (31.1,-4)--(32.8,-4);
	\draw (33,-4) node{$1$};
	\redbraidbox{33.7}{34.6}{-4.7}{-3.4}{};
	\draw (34.1,-4) node{$\color{red}0\,\,0\color{black}$};
	\draw (35.3,-4) node{$1$};
	\draw[dots] (35.9,-4)--(37.3,-4);
	\draw (38.5,-4) node{$i-1$};
	\draw (39.8,-4) node{$i$};
	\draw(39.8,-3.2)--(36,0)--(19.5,3.2);
	\draw(39.8,3.2)--(36,0)--(19.5,-3.2);
	\draw[blue](20.3,-3.2)--(18.5,0)--(20.3,3.2);
	\draw(21.3,-3.2)--(19.5,0)--(21.3,3.2);
	\draw(22.3,-3.2)--(20.5,0)--(22.3,3.2);
	\draw(26,-3.2)--(24.1,0)--(26,3.2);
	\draw(27,-3.2)--(25.1,0)--(27,3.2);
	\draw(28.6,-3.2)--(26.7,0)--(28.6,3.2);
	\draw(29.8,-3.2)--(27.9,0)--(29.8,3.2);
	\draw(38.5,-3.2)--(40,0)--(38.5,3.2);
	\draw(35.3,-3.2)--(32.5,0)--(35.3,3.2);
	\draw[red](34.5,-3.2)--(31.8,0)--(34.5,3.2);
	\draw[red](33.8,-3.2)--(31.1,0)--(33.8,3.2);
	\draw(33.1,-3.2)--(30.4,0)--(33.1,3.2);
	\end{braid}
	}
	+
	\resizebox{59mm}{15mm}{
	\begin{braid}\tikzset{baseline=-.3em}
	\draw (19.5,4) node{$i$};
		\draw (19.5,-4) node{$i$};
	\draw (20.3,4) node{\color{blue}$\ell$\color{black}};
	\braidbox{20.9}{23.2}{3.3}{4.6}{};
	\draw (22.1,4) node{$(\ell-1)^2$};
	\draw[dots] (23.8,4)--(25.3,4);
	\braidbox{25.6}{27.9}{3.3}{4.6}{};
	\draw (26.8,4) node{$(i+1)^2$};
	\draw (28.6,4) node{$i$};
	\draw (29.8,4) node{$i-1$};
	\draw[dots] (31.1,4)--(32.8,4);
	\draw (33,4) node{$1$};
	\redbraidbox{33.7}{34.6}{3.3}{4.6}{};
	\draw (34.1,4) node{$\color{red}0\,\,0\color{black}$};
	\draw (35.3,4) node{$1$};
	\draw[dots] (35.9,4)--(37.3,4);
	\draw (38.5,4) node{$i-1$};
	\draw (39.8,4) node{$i$};
	\draw (20.3,-4) node{\color{blue}$\ell$\color{black}};
	\braidbox{20.9}{23.2}{-4.7}{-3.4}{};
	\draw (22.1,-4) node{$(\ell-1)^2$};
	\draw[dots] (23.8,-4)--(25.1,-4);
	\braidbox{25.6}{27.9}{-4.7}{-3.4}{};
	\draw (26.8,-4) node{$(i+1)^2$};
	\draw (28.6,-4) node{$i$};
	\draw (29.8,-4) node{$i-1$};
	\draw[dots] (31.1,-4)--(32.8,-4);
	\draw (33,-4) node{$1$};
	\redbraidbox{33.7}{34.6}{-4.7}{-3.4}{};
	\draw (34.1,-4) node{$\color{red}0\,\,0\color{black}$};
	\draw (35.3,-4) node{$1$};
	\draw[dots] (35.9,-4)--(37.3,-4);
	\draw (38.5,-4) node{$i-1$};
	\draw (39.8,-4) node{$i$};
	\draw(39.8,-3.2)--(36.5,0)--(39.8,3.2);
	\draw(19.5,3.2)--(35.5,0)--(19.5,-3.2);
	\draw[blue](20.3,-3.2)--(18.5,0)--(20.3,3.2);
	\draw(21.3,-3.2)--(19.5,0)--(21.3,3.2);
	\draw(22.3,-3.2)--(20.5,0)--(22.3,3.2);
	\draw(26,-3.2)--(24.1,0)--(26,3.2);
	\draw(27,-3.2)--(25.1,0)--(27,3.2);
	\draw(28.6,-3.2)--(26.7,0)--(28.6,3.2);
	\draw(29.8,-3.2)--(27.9,0)--(29.8,3.2);
	\draw(38.5,-3.2)--(35.5,0)--(38.5,3.2);
	\draw(35.3,-3.2)--(32.5,0)--(35.3,3.2);
	\draw[red](34.5,-3.2)--(31.8,0)--(34.5,3.2);
	\draw[red](33.8,-3.2)--(31.1,0)--(33.8,3.2);
	\draw(33.1,-3.2)--(30.4,0)--(33.1,3.2);
	\end{braid}
	}
	\\
	=&
	\resizebox{59mm}{15mm}{
	\begin{braid}\tikzset{baseline=-.3em}
	\draw (19.5,4) node{$i$};
		\draw (19.5,-4) node{$i$};
	\draw (20.3,4) node{\color{blue}$\ell$\color{black}};
	\braidbox{20.9}{23.2}{3.3}{4.6}{};
	\draw (22.1,4) node{$(\ell-1)^2$};
	\draw[dots] (23.8,4)--(25.3,4);
	\braidbox{25.6}{27.9}{3.3}{4.6}{};
	\draw (26.8,4) node{$(i+1)^2$};
	\draw (28.6,4) node{$i$};
	\draw (29.8,4) node{$i-1$};
	\draw[dots] (31.1,4)--(32.8,4);
	\draw (33,4) node{$1$};
	\redbraidbox{33.7}{34.6}{3.3}{4.6}{};
	\draw (34.1,4) node{$\color{red}0\,\,0\color{black}$};
	\draw (35.3,4) node{$1$};
	\draw[dots] (35.9,4)--(37.3,4);
	\draw (38.5,4) node{$i-1$};
	\draw (39.8,4) node{$i$};
	\draw (20.3,-4) node{\color{blue}$\ell$\color{black}};
	\braidbox{20.9}{23.2}{-4.7}{-3.4}{};
	\draw (22.1,-4) node{$(\ell-1)^2$};
	\draw[dots] (23.8,-4)--(25.1,-4);
	\braidbox{25.6}{27.9}{-4.7}{-3.4}{};
	\draw (26.8,-4) node{$(i+1)^2$};
	\draw (28.6,-4) node{$i$};
	\draw (29.8,-4) node{$i-1$};
	\draw[dots] (31.1,-4)--(32.8,-4);
	\draw (33,-4) node{$1$};
	\redbraidbox{33.7}{34.6}{-4.7}{-3.4}{};
	\draw (34.1,-4) node{$\color{red}0\,\,0\color{black}$};
	\draw (35.3,-4) node{$1$};
	\draw[dots] (35.9,-4)--(37.3,-4);
	\draw (38.5,-4) node{$i-1$};
	\draw (39.8,-4) node{$i$};
	\draw(39.8,-3.2)--(19.5,3.2);
	\draw(39.8,3.2)--(19.5,-3.2);
	\draw[blue](20.3,-3.2)--(18.5,0)--(20.3,3.2);
	\draw(21.3,-3.2)--(19.5,0)--(21.3,3.2);
	\draw(22.3,-3.2)--(20.5,0)--(22.3,3.2);
	\draw(26,-3.2)--(24.1,0)--(26,3.2);
	\draw(27,-3.2)--(25.1,0)--(27,3.2);
	\draw(28.6,-3.2)--(26.7,0)--(28.6,3.2);
	\draw(29.8,-3.2)--(28.1,0)--(29.8,3.2);
	\draw(38.5,-3.2)--(40,0)--(38.5,3.2);
	\draw(35.3,-3.2)--(37.5,0)--(35.3,3.2);
	\draw[red](34.5,-3.2)--(36.8,0)--(34.5,3.2);
	\draw[red](33.8,-3.2)--(36.1,0)--(33.8,3.2);
	\draw(33.1,-3.2)--(35.5,0)--(33.1,3.2);
	\end{braid}
	}
	+
	\resizebox{59mm}{15mm}{
	\begin{braid}\tikzset{baseline=-.3em}
	\draw (19.5,4) node{$i$};
		\draw (19.5,-4) node{$i$};
	\draw (20.3,4) node{\color{blue}$\ell$\color{black}};
	\braidbox{20.9}{23.2}{3.3}{4.6}{};
	\draw (22.1,4) node{$(\ell-1)^2$};
	\draw[dots] (23.8,4)--(25.3,4);
	\braidbox{25.6}{27.9}{3.3}{4.6}{};
	\draw (26.8,4) node{$(i+1)^2$};
	\draw (28.6,4) node{$i$};
	\draw (29.8,4) node{$i-1$};
	\draw[dots] (31.1,4)--(32.8,4);
	\draw (33,4) node{$1$};
	\redbraidbox{33.7}{34.6}{3.3}{4.6}{};
	\draw (34.1,4) node{$\color{red}0\,\,0\color{black}$};
	\draw (35.3,4) node{$1$};
	\draw[dots] (35.9,4)--(37.3,4);
	\draw (38.5,4) node{$i-1$};
	\draw (39.8,4) node{$i$};
	\draw (20.3,-4) node{\color{blue}$\ell$\color{black}};
	\braidbox{20.9}{23.2}{-4.7}{-3.4}{};
	\draw (22.1,-4) node{$(\ell-1)^2$};
	\draw[dots] (23.8,-4)--(25.1,-4);
	\braidbox{25.6}{27.9}{-4.7}{-3.4}{};
	\draw (26.8,-4) node{$(i+1)^2$};
	\draw (28.6,-4) node{$i$};
	\draw (29.8,-4) node{$i-1$};
	\draw[dots] (31.1,-4)--(32.8,-4);
	\draw (33,-4) node{$1$};
	\redbraidbox{33.7}{34.6}{-4.7}{-3.4}{};
	\draw (34.1,-4) node{$\color{red}0\,\,0\color{black}$};
	\draw (35.3,-4) node{$1$};
	\draw[dots] (35.9,-4)--(37.3,-4);
	\draw (38.5,-4) node{$i-1$};
	\draw (39.8,-4) node{$i$};
	\draw(39.8,-3.2)--(39.8,3.2);
	\draw(19.5,3.2)--(28.6,0)--(19.5,-3.2);
	\draw[blue](20.3,-3.2)--(18.5,0)--(20.3,3.2);
	\draw(21.3,-3.2)--(19.5,0)--(21.3,3.2);
	\draw(22.3,-3.2)--(20.5,0)--(22.3,3.2);
	\draw(26,-3.2)--(24.1,0)--(26,3.2);
	\draw(27,-3.2)--(25.1,0)--(27,3.2);
	\draw(28.6,-3.2)--(26.7,0)--(28.6,3.2);
	\draw(29.8,-3.2)--(27.7,0)--(29.8,3.2);
	\draw(38.5,-3.2)--(38.5,3.2);
	\draw(35.3,-3.2)--(35.3,3.2);
	\draw[red](34.5,-3.2)--(34.5,3.2);
	\draw[red](33.8,-3.2)--(33.8,3.2);
	\draw(33.1,-3.2)--(33.1,3.2);
	\end{braid}.
	}
\end{align*}
Denote the first summand by $I$ and the second summand by $II$. 
The error term arising from the application of the braid relation 
$\begin{braid}\tikzset{baseline=.4em}
	\draw (0,0) node{$i$};
        \draw (1.1,0) node{$i-1$};
        \draw (2.2,0) node{$i$};
		\draw(0,0.5)--(2.2,1.5);
\draw(2.2,0.5)--(0,1.5);
\draw(1.1,0.5)--(0.2,1)--(1.1,1.5);
	\end{braid}$ 
in $I$ equals $0$ due to the quadratic relation $\begin{braid}\tikzset{baseline=.5em}
	\draw (0,0) node{$i$};
        \draw (1,0) node{$i$};
		\draw(0,0.5)--(1,1)--(0,1.5);
\draw(1,0.5)--(0,1)--(1,1.5);
	\end{braid}=0$. 
	So 
\begin{align*}
I=
\resizebox{59mm}{14mm}{
\begin{braid}\tikzset{baseline=-.3em}
	\draw (19.5,4) node{$i$};
		\draw (19.5,-4) node{$i$};
	\draw (20.3,4) node{\color{blue}$\ell$\color{black}};
	\braidbox{20.9}{23.2}{3.3}{4.6}{};
	\draw (22.1,4) node{$(\ell-1)^2$};
	\draw[dots] (23.8,4)--(25.3,4);
	\braidbox{25.6}{27.9}{3.3}{4.6}{};
	\draw (26.8,4) node{$(i+1)^2$};
	\draw (28.6,4) node{$i$};
	\draw (29.8,4) node{$i-1$};
	\draw[dots] (31.1,4)--(32.8,4);
	\draw (33,4) node{$1$};
	\redbraidbox{33.7}{34.6}{3.3}{4.6}{};
	\draw (34.1,4) node{$\color{red}0\,\,0\color{black}$};
	\draw (35.3,4) node{$1$};
	\draw[dots] (35.9,4)--(37.3,4);
	\draw (38.5,4) node{$i-1$};
	\draw (39.8,4) node{$i$};
	\draw (20.3,-4) node{\color{blue}$\ell$\color{black}};
	\braidbox{20.9}{23.2}{-4.7}{-3.4}{};
	\draw (22.1,-4) node{$(\ell-1)^2$};
	\draw[dots] (23.8,-4)--(25.1,-4);
	\braidbox{25.6}{27.9}{-4.7}{-3.4}{};
	\draw (26.8,-4) node{$(i+1)^2$};
	\draw (28.6,-4) node{$i$};
	\draw (29.8,-4) node{$i-1$};
	\draw[dots] (31.1,-4)--(32.8,-4);
	\draw (33,-4) node{$1$};
	\redbraidbox{33.7}{34.6}{-4.7}{-3.4}{};
	\draw (34.1,-4) node{$\color{red}0\,\,0\color{black}$};
	\draw (35.3,-4) node{$1$};
	\draw[dots] (35.9,-4)--(37.3,-4);
	\draw (38.5,-4) node{$i-1$};
	\draw (39.8,-4) node{$i$};
	\draw(39.8,-3.2)--(26,0)--(19.5,3.2);
	\draw(39.8,3.2)--(26,0)--(19.5,-3.2);
	\draw[blue](20.3,-3.2)--(18.5,0)--(20.3,3.2);
	\draw(21.3,-3.2)--(19.5,0)--(21.3,3.2);
	\draw(22.3,-3.2)--(20.5,0)--(22.3,3.2);
	\draw(26,-3.2)--(24.1,0)--(26,3.2);
	\draw(27,-3.2)--(25.1,0)--(27,3.2);
	\draw(28.6,-3.2)--(30.7,0)--(28.6,3.2);
	\draw(29.8,-3.2)--(32,0)--(29.8,3.2);
	\draw(38.5,-3.2)--(40,0)--(38.5,3.2);
	\draw(35.3,-3.2)--(37.5,0)--(35.3,3.2);
	\draw[red](34.5,-3.2)--(36.8,0)--(34.5,3.2);
	\draw[red](33.8,-3.2)--(36.1,0)--(33.8,3.2);
	\draw(33.1,-3.2)--(35.5,0)--(33.1,3.2);
	\end{braid}.}
\end{align*}
Applying Lemma~\ref{LA3/4}, using again $\begin{braid}\tikzset{baseline=.5em}
	\draw (0,0) node{$i$};
        \draw (1,0) node{$i$};
		\draw(0,0.5)--(1,1)--(0,1.5);
\draw(1,0.5)--(0,1)--(1,1.5);
	\end{braid}=0$
, and the fact that the word beginning with $\ell(\ell-1)^2\cdots(i+2)^2i$ is not cuspidal by Lemma~\ref{LCuspExpl}, we get 
\begin{align*}
I=-
\resizebox{59mm}{14mm}{\begin{braid}\tikzset{baseline=-.3em}
	\draw (19.5,4) node{$i$};
		\draw (19.5,-4) node{$i$};
	\draw (20.3,4) node{\color{blue}$\ell$\color{black}};
	\braidbox{20.9}{23.2}{3.3}{4.6}{};
	\draw (22.1,4) node{$(\ell-1)^2$};
	\draw[dots] (23.8,4)--(25.3,4);
	\braidbox{25.6}{27.9}{3.3}{4.6}{};
	\draw (26.8,4) node{$(i+1)^2$};
	\draw (28.6,4) node{$i$};
	\draw (29.8,4) node{$i-1$};
	\draw[dots] (31.1,4)--(32.8,4);
	\draw (33,4) node{$1$};
	\redbraidbox{33.7}{34.6}{3.3}{4.6}{};
	\draw (34.1,4) node{$\color{red}0\,\,0\color{black}$};
	\draw (35.3,4) node{$1$};
	\draw[dots] (35.9,4)--(37.3,4);
	\draw (38.5,4) node{$i-1$};
	\draw (39.8,4) node{$i$};
	\draw (20.3,-4) node{\color{blue}$\ell$\color{black}};
	\braidbox{20.9}{23.2}{-4.7}{-3.4}{};
	\draw (22.1,-4) node{$(\ell-1)^2$};
	\draw[dots] (23.8,-4)--(25.1,-4);
	\braidbox{25.6}{27.9}{-4.7}{-3.4}{};
	\draw (26.8,-4) node{$(i+1)^2$};
	\draw (28.6,-4) node{$i$};
	\draw (29.8,-4) node{$i-1$};
	\draw[dots] (31.1,-4)--(32.8,-4);
	\draw (33,-4) node{$1$};
	\redbraidbox{33.7}{34.6}{-4.7}{-3.4}{};
	\draw (34.1,-4) node{$\color{red}0\,\,0\color{black}$};
	\draw (35.3,-4) node{$1$};
	\draw[dots] (35.9,-4)--(37.3,-4);
	\draw (38.5,-4) node{$i-1$};
	\draw (39.8,-4) node{$i$};
	\draw(19.5,-3.2)--(25.5,0)--(19.5,3.2);
	\draw(39.8,3.2)--(27.5,0)--(39.8,-3.2);
	\draw[blue](20.3,-3.2)--(18.5,0)--(20.3,3.2);
	\draw(21.3,-3.2)--(19.5,0)--(21.3,3.2);
	\draw(22.3,-3.2)--(20.5,0)--(22.3,3.2);
	\draw(26,-3.2)--(26,3.2);
	\draw(27,-3.2)--(27,3.2);
	\draw(28.6,-3.2)--(30.7,0)--(28.6,3.2);
	\draw(29.8,-3.2)--(32,0)--(29.8,3.2);
	\draw(38.5,-3.2)--(40,0)--(38.5,3.2);
	\draw(35.3,-3.2)--(37.5,0)--(35.3,3.2);
	\draw[red](34.5,-3.2)--(36.8,0)--(34.5,3.2);
	\draw[red](33.8,-3.2)--(36.1,0)--(33.8,3.2);
	\draw(33.1,-3.2)--(35.5,0)--(33.1,3.2);
	\blackdot (27.7,0);
	\end{braid}.
	}
\end{align*}
Using the relation (\ref{R5}) and quadratic relations, we get
$$
I=
\resizebox{70mm}{19mm}{
\begin{braid}\tikzset{baseline=-.3em}
	\draw (19.5,4) node{$i$};
		\draw (19.5,-4) node{$i$};
	\draw (20.3,4) node{\color{blue}$\ell$\color{black}};
	\braidbox{20.9}{23.2}{3.3}{4.6}{};
	\draw (22.1,4) node{$(\ell-1)^2$};
	\draw[dots] (23.8,4)--(25.3,4);
	\braidbox{25.6}{27.9}{3.3}{4.6}{};
	\draw (26.8,4) node{$(i+1)^2$};
	\draw (28.6,4) node{$i$};
	\draw (29.8,4) node{$i-1$};
	\draw[dots] (31.1,4)--(32.8,4);
	\draw (33,4) node{$1$};
	\redbraidbox{33.7}{34.6}{3.3}{4.6}{};
	\draw (34.1,4) node{$\color{red}0\,\,0\color{black}$};
	\draw (35.3,4) node{$1$};
	\draw[dots] (35.9,4)--(37.3,4);
	\draw (38.5,4) node{$i-1$};
	\draw (39.8,4) node{$i$};
	\draw (20.3,-4) node{\color{blue}$\ell$\color{black}};
	\braidbox{20.9}{23.2}{-4.7}{-3.4}{};
	\draw (22.1,-4) node{$(\ell-1)^2$};
	\draw[dots] (23.8,-4)--(25.1,-4);
	\braidbox{25.6}{27.9}{-4.7}{-3.4}{};
	\draw (26.8,-4) node{$(i+1)^2$};
	\draw (28.6,-4) node{$i$};
	\draw (29.8,-4) node{$i-1$};
	\draw[dots] (31.1,-4)--(32.8,-4);
	\draw (33,-4) node{$1$};
	\redbraidbox{33.7}{34.6}{-4.7}{-3.4}{};
	\draw (34.1,-4) node{$\color{red}0\,\,0\color{black}$};
	\draw (35.3,-4) node{$1$};
	\draw[dots] (35.9,-4)--(37.3,-4);
	\draw (38.5,-4) node{$i-1$};
	\draw (39.8,-4) node{$i$};
	\draw(19.5,-3.2)--(19.5,3.2);
	\draw(28.6,3.2)--(30,0)--(39.8,-3.2);
	\draw[blue](20.3,-3.2)--(20.3,3.2);
	\draw(21.3,-3.2)--(21.3,3.2);
	\draw(22.3,-3.2)--(22.3,3.2);
	\draw(26,-3.2)--(26,3.2);
	\draw(27,-3.2)--(27,3.2);
	\draw(28.6,-3.2)--(30,0)--(39.8,3.2);
	\draw(29.8,-3.2)--(32,0)--(29.8,3.2);
	\draw(38.5,-3.2)--(40,0)--(38.5,3.2);
	\draw(35.3,-3.2)--(37.5,0)--(35.3,3.2);
	\draw[red](34.5,-3.2)--(36.8,0)--(34.5,3.2);
	\draw[red](33.8,-3.2)--(36.1,0)--(33.8,3.2);
	\draw(33.1,-3.2)--(35.5,0)--(33.1,3.2);
	\end{braid}.
	}
$$
Using the braid relation for 
$\begin{braid}\tikzset{baseline=.4em}
	\draw (0,0) node{$i$};
        \draw (1.1,0) node{$i-1$};
        \draw (2.2,0) node{$i$};
		\draw(0,0.5)--(2.2,1.5);
\draw(2.2,0.5)--(0,1.5);
\draw(1.1,0.5)--(2,1)--(1.1,1.5);
	\end{braid}$ 
 and the fact that the word beginning with $\ell(\ell-1)^2\cdots(i+2)^2(i-1)$ is not cuspidal by Lemma~\ref{LCuspExpl}, we get 
$$
I=-
\resizebox{70mm}{19mm}{
\begin{braid}\tikzset{baseline=-.3em}
	\draw (19.5,4) node{$i$};
		\draw (19.5,-4) node{$i$};
	\draw (20.3,4) node{\color{blue}$\ell$\color{black}};
	\braidbox{20.9}{23.2}{3.3}{4.6}{};
	\draw (22.1,4) node{$(\ell-1)^2$};
	\draw[dots] (23.8,4)--(25.3,4);
	\braidbox{25.6}{27.9}{3.3}{4.6}{};
	\draw (26.8,4) node{$(i+1)^2$};
	\draw (28.6,4) node{$i$};
	\draw (29.8,4) node{$i-1$};
	\draw[dots] (31.1,4)--(32.8,4);
	\draw (33,4) node{$1$};
	\redbraidbox{33.7}{34.6}{3.3}{4.6}{};
	\draw (34.1,4) node{$\color{red}0\,\,0\color{black}$};
	\draw (35.3,4) node{$1$};
	\draw[dots] (35.9,4)--(37.3,4);
	\draw (38.5,4) node{$i-1$};
	\draw (39.8,4) node{$i$};
	\draw (20.3,-4) node{\color{blue}$\ell$\color{black}};
	\braidbox{20.9}{23.2}{-4.7}{-3.4}{};
	\draw (22.1,-4) node{$(\ell-1)^2$};
	\draw[dots] (23.8,-4)--(25.1,-4);
	\braidbox{25.6}{27.9}{-4.7}{-3.4}{};
	\draw (26.8,-4) node{$(i+1)^2$};
	\draw (28.6,-4) node{$i$};
	\draw (29.8,-4) node{$i-1$};
	\draw[dots] (31.1,-4)--(32.8,-4);
	\draw (33,-4) node{$1$};
	\redbraidbox{33.7}{34.6}{-4.7}{-3.4}{};
	\draw (34.1,-4) node{$\color{red}0\,\,0\color{black}$};
	\draw (35.3,-4) node{$1$};
	\draw[dots] (35.9,-4)--(37.3,-4);
	\draw (38.5,-4) node{$i-1$};
	\draw (39.8,-4) node{$i$};
	\draw(19.5,-3.2)--(19.5,3.2);
	\draw(39.8,3.2)--(30.5,0)--(39.8,-3.2);
	\draw[blue](20.3,-3.2)--(20.3,3.2);
	\draw(21.3,-3.2)--(21.3,3.2);
	\draw(22.3,-3.2)--(22.3,3.2);
	\draw(26,-3.2)--(26,3.2);
	\draw(27,-3.2)--(27,3.2);
	\draw(28.6,-3.2)--(28.6,3.2);
	\draw(29.8,-3.2)--(29.8,3.2);
	\draw(38.5,-3.2)--(40,0)--(38.5,3.2);
	\draw(35.3,-3.2)--(37.5,0)--(35.3,3.2);
	\draw[red](34.5,-3.2)--(36.8,0)--(34.5,3.2);
	\draw[red](33.8,-3.2)--(36.1,0)--(33.8,3.2);
	\draw(33.1,-3.2)--(35.5,0)--(33.1,3.2);
	\end{braid}.
	}
$$
Using quadratic relations we now get
$$
I=-
\resizebox{70mm}{19mm}{
\begin{braid}\tikzset{baseline=-.3em}
	\draw (19.5,4) node{$i$};
		\draw (19.5,-4) node{$i$};
	\draw (20.3,4) node{\color{blue}$\ell$\color{black}};
	\braidbox{20.9}{23.2}{3.3}{4.6}{};
	\draw (22.1,4) node{$(\ell-1)^2$};
	\draw[dots] (23.8,4)--(25.3,4);
	\braidbox{25.6}{27.9}{3.3}{4.6}{};
	\draw (26.8,4) node{$(i+1)^2$};
	\draw (28.6,4) node{$i$};
	\draw (29.8,4) node{$i-1$};
	\draw[dots] (31.1,4)--(32.8,4);
	\draw (33,4) node{$1$};
	\redbraidbox{33.7}{34.6}{3.3}{4.6}{};
	\draw (34.1,4) node{$\color{red}0\,\,0\color{black}$};
	\draw (35.3,4) node{$1$};
	\draw[dots] (35.9,4)--(37.3,4);
	\draw (38.5,4) node{$i-1$};
	\draw (39.8,4) node{$i$};
	\draw (20.3,-4) node{\color{blue}$\ell$\color{black}};
	\braidbox{20.9}{23.2}{-4.7}{-3.4}{};
	\draw (22.1,-4) node{$(\ell-1)^2$};
	\draw[dots] (23.8,-4)--(25.1,-4);
	\braidbox{25.6}{27.9}{-4.7}{-3.4}{};
	\draw (26.8,-4) node{$(i+1)^2$};
	\draw (28.6,-4) node{$i$};
	\draw (29.8,-4) node{$i-1$};
	\draw[dots] (31.1,-4)--(32.8,-4);
	\draw (33,-4) node{$1$};
	\redbraidbox{33.7}{34.6}{-4.7}{-3.4}{};
	\draw (34.1,-4) node{$\color{red}0\,\,0\color{black}$};
	\draw (35.3,-4) node{$1$};
	\draw[dots] (35.9,-4)--(37.3,-4);
	\draw (38.5,-4) node{$i-1$};
	\draw (39.8,-4) node{$i$};
	\draw(19.5,-3.2)--(19.5,3.2);
	\draw(39.8,3.2)--(38.5,0)--(39.8,-3.2);
	\draw[blue](20.3,-3.2)--(20.3,3.2);
	\draw(21.3,-3.2)--(21.3,3.2);
	\draw(22.3,-3.2)--(22.3,3.2);
	\draw(26,-3.2)--(26,3.2);
	\draw(27,-3.2)--(27,3.2);
	\draw(28.6,-3.2)--(28.6,3.2);
	\draw(29.8,-3.2)--(29.8,3.2);
	\draw(38.5,-3.2)--(39.5,0)--(38.5,3.2);
	\draw(35.3,-3.2)--(35.3,3.2);
	\draw[red](34.5,-3.2)--(34.5,3.2);
	\draw[red](33.8,-3.2)--(33.8,3.2);
	\draw(33.1,-3.2)--(33.1,3.2);
	\end{braid}
	}
	=\ga^{i,i}(y_{2p}-y_{2p-1})\ga^{i,i}.
$$
By (\ref{ECFormula}), we have 
$\ga^{i,i}y_{2p}\ga^{i,i}=\ga^{i,i}((-1)^ic_2-z_2)\ga^{i,i},$ 
while by (\ref{ERelB3}), we have $\ga^{i,i}y_{2p-1}\ga^{i,i}=-\ga^{i,i}z_2\ga^{i,i}$. So 
$$I=(-1)^i\ga^{i,i}c_2\ga^{i,i}.$$ 

On the other hand, using (\ref{R5}) and (\ref{R6}), we get
$$
II=-
\resizebox{70mm}{19mm}{
\begin{braid}\tikzset{baseline=-.3em}
	\draw (19.5,4) node{$i$};
		\draw (19.5,-4) node{$i$};
	\draw (20.3,4) node{\color{blue}$\ell$\color{black}};
	\braidbox{20.9}{23.2}{3.3}{4.6}{};
	\draw (22.1,4) node{$(\ell-1)^2$};
	\draw[dots] (23.8,4)--(25.3,4);
	\braidbox{25.6}{27.9}{3.3}{4.6}{};
	\draw (26.8,4) node{$(i+1)^2$};
	\draw (28.6,4) node{$i$};
	\draw (29.8,4) node{$i-1$};
	\draw[dots] (31.1,4)--(32.8,4);
	\draw (33,4) node{$1$};
	\redbraidbox{33.7}{34.6}{3.3}{4.6}{};
	\draw (34.1,4) node{$\color{red}0\,\,0\color{black}$};
	\draw (35.3,4) node{$1$};
	\draw[dots] (35.9,4)--(37.3,4);
	\draw (38.5,4) node{$i-1$};
	\draw (39.8,4) node{$i$};
	\draw (20.3,-4) node{\color{blue}$\ell$\color{black}};
	\braidbox{20.9}{23.2}{-4.7}{-3.4}{};
	\draw (22.1,-4) node{$(\ell-1)^2$};
	\draw[dots] (23.8,-4)--(25.1,-4);
	\braidbox{25.6}{27.9}{-4.7}{-3.4}{};
	\draw (26.8,-4) node{$(i+1)^2$};
	\draw (28.6,-4) node{$i$};
	\draw (29.8,-4) node{$i-1$};
	\draw[dots] (31.1,-4)--(32.8,-4);
	\draw (33,-4) node{$1$};
	\redbraidbox{33.7}{34.6}{-4.7}{-3.4}{};
	\draw (34.1,-4) node{$\color{red}0\,\,0\color{black}$};
	\draw (35.3,-4) node{$1$};
	\draw[dots] (35.9,-4)--(37.3,-4);
	\draw (38.5,-4) node{$i-1$};
	\draw (39.8,-4) node{$i$};
	\draw(39.8,-3.2)--(39.8,3.2);
	\draw(19.5,3.2)--(26,0)--(28.6,-3.2);
	\draw[blue](20.3,-3.2)--(18.5,0)--(20.3,3.2);
	\draw(21.3,-3.2)--(19.5,0)--(21.3,3.2);
	\draw(22.3,-3.2)--(20.5,0)--(22.3,3.2);
	\draw(26,-3.2)--(24.1,0)--(26,3.2);
	\draw(27,-3.2)--(25.1,0)--(27,3.2);
	\draw(19.5,-3.2)--(26,0)--(28.6,3.2);
	\draw(29.8,-3.2)--(29.8,3.2);
	\draw(38.5,-3.2)--(38.5,3.2);
	\draw(35.3,-3.2)--(35.3,3.2);
	\draw[red](34.5,-3.2)--(34.5,3.2);
	\draw[red](33.8,-3.2)--(33.8,3.2);
	\draw(33.1,-3.2)--(33.1,3.2);
	\end{braid}.
	}
$$
Applying Lemma~\ref{LA3/4} and (\ref{R5}), (\ref{R6}), we deduce that $II$ equals 
\begin{align*}
&
\resizebox{59mm}{15mm}{
\begin{braid}\tikzset{baseline=-.3em}
	\draw (19.5,4) node{$i$};
		\draw (19.5,-4) node{$i$};
	\draw (20.3,4) node{\color{blue}$\ell$\color{black}};
	\braidbox{20.9}{23.2}{3.3}{4.6}{};
	\draw (22.1,4) node{$(\ell-1)^2$};
	\draw[dots] (23.8,4)--(25.3,4);
	\braidbox{25.6}{27.9}{3.3}{4.6}{};
	\draw (26.8,4) node{$(i+1)^2$};
	\draw (28.6,4) node{$i$};
	\draw (29.8,4) node{$i-1$};
	\draw[dots] (31.1,4)--(32.8,4);
	\draw (33,4) node{$1$};
	\redbraidbox{33.7}{34.6}{3.3}{4.6}{};
	\draw (34.1,4) node{$\color{red}0\,\,0\color{black}$};
	\draw (35.3,4) node{$1$};
	\draw[dots] (35.9,4)--(37.3,4);
	\draw (38.5,4) node{$i-1$};
	\draw (39.8,4) node{$i$};
	\draw (20.3,-4) node{\color{blue}$\ell$\color{black}};
	\braidbox{20.9}{23.2}{-4.7}{-3.4}{};
	\draw (22.1,-4) node{$(\ell-1)^2$};
	\draw[dots] (23.8,-4)--(25.1,-4);
	\braidbox{25.6}{27.9}{-4.7}{-3.4}{};
	\draw (26.8,-4) node{$(i+1)^2$};
	\draw (28.6,-4) node{$i$};
	\draw (29.8,-4) node{$i-1$};
	\draw[dots] (31.1,-4)--(32.8,-4);
	\draw (33,-4) node{$1$};
	\redbraidbox{33.7}{34.6}{-4.7}{-3.4}{};
	\draw (34.1,-4) node{$\color{red}0\,\,0\color{black}$};
	\draw (35.3,-4) node{$1$};
	\draw[dots] (35.9,-4)--(37.3,-4);
	\draw (38.5,-4) node{$i-1$};
	\draw (39.8,-4) node{$i$};
	\draw(39.8,-3.2)--(39.8,3.2);
	\draw(28.6,3.2)--(28.6,-3.2);
	\draw[blue](20.3,-3.2)--(20.3,3.2);
	\draw(21.3,-3.2)--(21.3,3.2);
	\draw(22.3,-3.2)--(22.3,3.2);
	\draw(26,-3.2)--(26,3.2);
	\draw(27,-3.2)--(27,3.2);
	\draw(19.5,-3.2)--(19.5,3.2);
	\draw(29.8,-3.2)--(29.8,3.2);
	\draw(38.5,-3.2)--(38.5,3.2);
	\draw(35.3,-3.2)--(35.3,3.2);
	\draw[red](34.5,-3.2)--(34.5,3.2);
	\draw[red](33.8,-3.2)--(33.8,3.2);
	\draw(33.1,-3.2)--(33.1,3.2);
	\blackdot (19.5,0);
	\end{braid}
	}
	-
	\resizebox{59mm}{15mm}{
	\begin{braid}\tikzset{baseline=-.3em}
	\draw (19.5,4) node{$i$};
		\draw (19.5,-4) node{$i$};
	\draw (20.3,4) node{\color{blue}$\ell$\color{black}};
	\braidbox{20.9}{23.2}{3.3}{4.6}{};
	\draw (22.1,4) node{$(\ell-1)^2$};
	\draw[dots] (23.8,4)--(25.3,4);
	\braidbox{25.6}{27.9}{3.3}{4.6}{};
	\draw (26.8,4) node{$(i+1)^2$};
	\draw (28.6,4) node{$i$};
	\draw (29.8,4) node{$i-1$};
	\draw[dots] (31.1,4)--(32.8,4);
	\draw (33,4) node{$1$};
	\redbraidbox{33.7}{34.6}{3.3}{4.6}{};
	\draw (34.1,4) node{$\color{red}0\,\,0\color{black}$};
	\draw (35.3,4) node{$1$};
	\draw[dots] (35.9,4)--(37.3,4);
	\draw (38.5,4) node{$i-1$};
	\draw (39.8,4) node{$i$};
	\draw (20.3,-4) node{\color{blue}$\ell$\color{black}};
	\braidbox{20.9}{23.2}{-4.7}{-3.4}{};
	\draw (22.1,-4) node{$(\ell-1)^2$};
	\draw[dots] (23.8,-4)--(25.1,-4);
	\braidbox{25.6}{27.9}{-4.7}{-3.4}{};
	\draw (26.8,-4) node{$(i+1)^2$};
	\draw (28.6,-4) node{$i$};
	\draw (29.8,-4) node{$i-1$};
	\draw[dots] (31.1,-4)--(32.8,-4);
	\draw (33,-4) node{$1$};
	\redbraidbox{33.7}{34.6}{-4.7}{-3.4}{};
	\draw (34.1,-4) node{$\color{red}0\,\,0\color{black}$};
	\draw (35.3,-4) node{$1$};
	\draw[dots] (35.9,-4)--(37.3,-4);
	\draw (38.5,-4) node{$i-1$};
	\draw (39.8,-4) node{$i$};
	\draw(39.8,-3.2)--(39.8,3.2);
	\draw(28.6,3.2)--(28.6,-3.2);
	\draw[blue](20.3,-3.2)--(20.3,3.2);
	\draw(21.3,-3.2)--(21.3,3.2);
	\draw(22.3,-3.2)--(22.3,3.2);
	\draw(26,-3.2)--(26,3.2);
	\draw(27,-3.2)--(27,3.2);
	\draw(19.5,-3.2)--(19.5,3.2);
	\draw(29.8,-3.2)--(29.8,3.2);
	\draw(38.5,-3.2)--(38.5,3.2);
	\draw(35.3,-3.2)--(35.3,3.2);
	\draw[red](34.5,-3.2)--(34.5,3.2);
	\draw[red](33.8,-3.2)--(33.8,3.2);
	\draw(33.1,-3.2)--(33.1,3.2);
	\blackdot (26,0);
	\end{braid}
	}
	\\
	&-
	\resizebox{59mm}{15mm}{
	\begin{braid}\tikzset{baseline=-.3em}
	\draw (19.5,4) node{$i$};
		\draw (19.5,-4) node{$i$};
	\draw (20.3,4) node{\color{blue}$\ell$\color{black}};
	\braidbox{20.9}{23.2}{3.3}{4.6}{};
	\draw (22.1,4) node{$(\ell-1)^2$};
	\draw[dots] (23.8,4)--(25.3,4);
	\braidbox{25.6}{27.9}{3.3}{4.6}{};
	\draw (26.8,4) node{$(i+1)^2$};
	\draw (28.6,4) node{$i$};
	\draw (29.8,4) node{$i-1$};
	\draw[dots] (31.1,4)--(32.8,4);
	\draw (33,4) node{$1$};
	\redbraidbox{33.7}{34.6}{3.3}{4.6}{};
	\draw (34.1,4) node{$\color{red}0\,\,0\color{black}$};
	\draw (35.3,4) node{$1$};
	\draw[dots] (35.9,4)--(37.3,4);
	\draw (38.5,4) node{$i-1$};
	\draw (39.8,4) node{$i$};
	\draw (20.3,-4) node{\color{blue}$\ell$\color{black}};
	\braidbox{20.9}{23.2}{-4.7}{-3.4}{};
	\draw (22.1,-4) node{$(\ell-1)^2$};
	\draw[dots] (23.8,-4)--(25.1,-4);
	\braidbox{25.6}{27.9}{-4.7}{-3.4}{};
	\draw (26.8,-4) node{$(i+1)^2$};
	\draw (28.6,-4) node{$i$};
	\draw (29.8,-4) node{$i-1$};
	\draw[dots] (31.1,-4)--(32.8,-4);
	\draw (33,-4) node{$1$};
	\redbraidbox{33.7}{34.6}{-4.7}{-3.4}{};
	\draw (34.1,-4) node{$\color{red}0\,\,0\color{black}$};
	\draw (35.3,-4) node{$1$};
	\draw[dots] (35.9,-4)--(37.3,-4);
	\draw (38.5,-4) node{$i-1$};
	\draw (39.8,-4) node{$i$};
	\draw(39.8,-3.2)--(39.8,3.2);
	\draw(28.6,3.2)--(28.6,-3.2);
	\draw[blue](20.3,-3.2)--(20.3,3.2);
	\draw(21.3,-3.2)--(21.3,3.2);
	\draw(22.3,-3.2)--(22.3,3.2);
	\draw(26,-3.2)--(26,3.2);
	\draw(27,-3.2)--(27,3.2);
	\draw(19.5,-3.2)--(19.5,3.2);
	\draw(29.8,-3.2)--(29.8,3.2);
	\draw(38.5,-3.2)--(38.5,3.2);
	\draw(35.3,-3.2)--(35.3,3.2);
	\draw[red](34.5,-3.2)--(34.5,3.2);
	\draw[red](33.8,-3.2)--(33.8,3.2);
	\draw(33.1,-3.2)--(33.1,3.2);
	\blackdot (27,0);
	\end{braid}
	}
	+
	\resizebox{59mm}{15mm}{
	\begin{braid}\tikzset{baseline=-.3em}
	\draw (19.5,4) node{$i$};
		\draw (19.5,-4) node{$i$};
	\draw (20.3,4) node{\color{blue}$\ell$\color{black}};
	\braidbox{20.9}{23.2}{3.3}{4.6}{};
	\draw (22.1,4) node{$(\ell-1)^2$};
	\draw[dots] (23.8,4)--(25.3,4);
	\braidbox{25.6}{27.9}{3.3}{4.6}{};
	\draw (26.8,4) node{$(i+1)^2$};
	\draw (28.6,4) node{$i$};
	\draw (29.8,4) node{$i-1$};
	\draw[dots] (31.1,4)--(32.8,4);
	\draw (33,4) node{$1$};
	\redbraidbox{33.7}{34.6}{3.3}{4.6}{};
	\draw (34.1,4) node{$\color{red}0\,\,0\color{black}$};
	\draw (35.3,4) node{$1$};
	\draw[dots] (35.9,4)--(37.3,4);
	\draw (38.5,4) node{$i-1$};
	\draw (39.8,4) node{$i$};
	\draw (20.3,-4) node{\color{blue}$\ell$\color{black}};
	\braidbox{20.9}{23.2}{-4.7}{-3.4}{};
	\draw (22.1,-4) node{$(\ell-1)^2$};
	\draw[dots] (23.8,-4)--(25.1,-4);
	\braidbox{25.6}{27.9}{-4.7}{-3.4}{};
	\draw (26.8,-4) node{$(i+1)^2$};
	\draw (28.6,-4) node{$i$};
	\draw (29.8,-4) node{$i-1$};
	\draw[dots] (31.1,-4)--(32.8,-4);
	\draw (33,-4) node{$1$};
	\redbraidbox{33.7}{34.6}{-4.7}{-3.4}{};
	\draw (34.1,-4) node{$\color{red}0\,\,0\color{black}$};
	\draw (35.3,-4) node{$1$};
	\draw[dots] (35.9,-4)--(37.3,-4);
	\draw (38.5,-4) node{$i-1$};
	\draw (39.8,-4) node{$i$};
	\draw(39.8,-3.2)--(39.8,3.2);
	\draw(28.6,3.2)--(28.6,-3.2);
	\draw[blue](20.3,-3.2)--(20.3,3.2);
	\draw(21.3,-3.2)--(21.3,3.2);
	\draw(22.3,-3.2)--(22.3,3.2);
	\draw(26,-3.2)--(26,3.2);
	\draw(27,-3.2)--(27,3.2);
	\draw(19.5,-3.2)--(19.5,3.2);
	\draw(29.8,-3.2)--(29.8,3.2);
	\draw(38.5,-3.2)--(38.5,3.2);
	\draw(35.3,-3.2)--(35.3,3.2);
	\draw[red](34.5,-3.2)--(34.5,3.2);
	\draw[red](33.8,-3.2)--(33.8,3.2);
	\draw(33.1,-3.2)--(33.1,3.2);
	\blackdot (28.6,0);
	\end{braid}.
	}
\end{align*}
By (\ref{ECFormula}), (\ref{ERelB2}) and  (\ref{ERelB3}), we now deduce 
$$II=\ga^{i,i}((-1)^ic_1-z_1+z_2)\ga^{i,i}.$$
So $\Theta^i_i=I+II$ is as required. 
\end{proof}

\begin{Lemma} \label{L110621Gen} 
Let $1\leq j\leq  i\leq \ell-1$. Then 
$$\Theta^i_j=\ga^{i,i}(-z_1+z_2+(-1)^i(c_1+c_2))\ga^{i,i}.$$ 
\end{Lemma}
\begin{proof}
In view of Lemma~\ref{L110621}, it suffices to prove that for $j<i$ we have $\Theta^i_j=\Theta^i_{j+1}$. We provide details for the generic case $j>1$. Dropping the left trivial part of the diagrams, we have
\begin{align*}
\Theta^i_j&=
\resizebox{117mm}{22mm}{
\begin{braid}\tikzset{baseline=-.3em}
        \draw[dots] (17.5,7)--(19.2,7);
	\draw[dots] (17.5,-7)--(19.2,-7);
	\draw[dots] (24.2,7)--(25.7,7);
	\draw[dots] (24.2,-7)--(25.7,-7);
	\draw[dots] (29.9,7)--(31.6,7);
	\draw[dots] (29.9,-7)--(31.6,-7);
	\draw[dots] (36.8,7)--(38.5,7);
	\draw[dots] (36.8,-7)--(38.5,-7);
	\draw[dots] (42,7)--(43.7,7);
	\draw[dots] (42,-7)--(43.7,-7);
	\draw[dots] (48.9,7)--(50.6,7);
	\draw[dots] (48.9,-7)--(50.6,-7);
	\draw (14.7,7) node{$j$};
	\draw (14.7,-7) node{$j$};
	\draw (16.2,7) node{$j+1$};
	\draw (16.2,-7) node{$j+1$};
	\draw (19.5,7) node{$i$};
	\draw (19.5,-7) node{$i$};
	\draw (20.4,7) node{\color{blue}$\ell$\color{black}};
	\draw (20.4,-7) node{\color{blue}$\ell$\color{black}};
	\braidbox{21.2}{23.6}{6.3}{7.6}{};
	\draw (22.4,7) node{$(\ell-1)^2$};
	\braidbox{21.2}{23.6}{-6.3}{-7.6}{};
	\draw (22.4,-7) node{$(\ell-1)^2$};
	\braidbox{26}{28.3}{6.3}{7.6}{};
	\draw (27.2,7) node{$(i+1)^2$};
	\braidbox{26}{28.3}{-6.3}{-7.6}{};
	\draw (27.2,-7) node{$(i+1)^2$};
	\draw (29.2,7) node{$i$};
	\draw (29.2,-7) node{$i$};
	\draw (32.5,7) node{$j+1$};
	\draw (32.5,-7) node{$j+1$};
	\draw (34,7) node{$j$};
	\draw (34,-7) node{$j$};
	\draw (35.5,7) node{$j-1$};
	\draw (35.5,-7) node{$j-1$};
	\draw (38.7,7) node{$1$};
	\draw (38.7,-7) node{$1$};
	\redbraidbox{39.5}{40.5}{6.3}{7.6}{};
	\draw (40,7) node{$\color{red}0\,\,0\color{black}$};
	\redbraidbox{39.5}{40.5}{-6.3}{-7.6}{};
	\draw (40,-7) node{$\color{red}0\,\,0\color{black}$};
	\draw (41.3,7) node{$1$};
	\draw (41.3,-7) node{$1$};
	\draw (44.7,7) node{$j-1$};
	\draw (44.7,-7) node{$j-1$};
	\draw (46.2,7) node{$j$};
	\draw (46.2,-7) node{$j$};
	\draw (47.7,7) node{$j+1$};
	\draw (47.7,-7) node{$j+1$};
	\draw (50.9,7) node{$i$};
	\draw (50.9,-7) node{$i$};
	\draw (14.7,6.3)--(40,0)--(46.2,-6.3);
	\draw (14.7,-6.3)--(40,0)--(46.2,6.3);
	\draw (16.2,6.3)--(41.4,0)--(47.7,-6.3);
	\draw (16.2,-6.3)--(41.4,0)--(47.7,6.3);
	\draw(19.5,6.3)--(44.7,0)--(50.9,-6.3);
	\draw(19.5,-6.3)--(44.7,0)--(50.9,6.3);
	\draw (44.7,6.3)--(38.7,0)--(44.7,-6.3);
	\draw (41.3,6.3)--(35.6,0)--(41.3,-6.3);
	\draw[red](40.4,6.3)--(34.7,0)--(40.4,-6.3);
	\draw[red](39.6,6.3)--(33.9,0)--(39.6,-6.3);
	\draw (38.7,6.3)--(32.9,0)--(38.7,-6.3);
	\draw (35.5,6.3)--(30.1,0)--(35.5,-6.3);
	\draw (34,6.3)--(28.8,0)--(34,-6.3);
	\draw (32.5,6.3)--(27.4,0)--(32.5,-6.3);
	\draw (29.2,6.3)--(24.3,0)--(29.2,-6.3);
	\draw (28,6.3)--(23,0)--(28,-6.3);
	\draw (26.5,6.3)--(21.5,0)--(26.5,-6.3);
	\draw (23.2,6.3)--(18.2,0)--(23.2,-6.3);
	\draw (21.7,6.3)--(16.7,0)--(21.7,-6.3);
	\draw[blue] (20.3,6.3)--(15.4,0)--(20.3,-6.3);
	\end{braid}.
	}
\end{align*}
We apply the braid relation (\ref{R7}) for 
$\begin{braid}\tikzset{baseline=.4em}
	\draw (0,0) node{$j$};
        \draw (1.1,0) node{$j-1$};
        \draw (2.2,0) node{$j$};
		\draw(0,0.5)--(2.2,1.5);
\draw(2.2,0.5)--(0,1.5);
\draw(1.1,0.5)--(0,1)--(1.1,1.5);
	\end{braid}$. 
The term corresponding to 
$\begin{braid}\tikzset{baseline=.4em}
	\draw (0,0) node{$j$};
        \draw (1.1,0) node{$j-1$};
        \draw (2.2,0) node{$j$};
		\draw(0,0.5)--(2.2,1.5);
\draw(2.2,0.5)--(0,1.5);
\draw(1.1,0.5)--(2.2,1)--(1.1,1.5);
	\end{braid}$
is $0$ since the word beginning with $\ell(\ell-1)^2\cdots(i+2)^2i\cdots 10^21\cdots (j-2)j$ is not cuspidal by Lemma~\ref{LCuspExpl}. So 
\begin{align*}
\Theta^i_j&=
\resizebox{117mm}{22mm}{
\begin{braid}\tikzset{baseline=-.3em}
        \draw[dots] (17.5,7)--(19.2,7);
	\draw[dots] (17.5,-7)--(19.2,-7);
	\draw[dots] (24.2,7)--(25.7,7);
	\draw[dots] (24.2,-7)--(25.7,-7);
	\draw[dots] (29.9,7)--(31.6,7);
	\draw[dots] (29.9,-7)--(31.6,-7);
	\draw[dots] (36.8,7)--(38.5,7);
	\draw[dots] (36.8,-7)--(38.5,-7);
	\draw[dots] (42,7)--(43.7,7);
	\draw[dots] (42,-7)--(43.7,-7);
	\draw[dots] (48.9,7)--(50.6,7);
	\draw[dots] (48.9,-7)--(50.6,-7);
	\draw (14.7,7) node{$j$};
	\draw (14.7,-7) node{$j$};
	\draw (16.2,7) node{$j+1$};
	\draw (16.2,-7) node{$j+1$};
	\draw (19.5,7) node{$i$};
	\draw (19.5,-7) node{$i$};
	\draw (20.4,7) node{\color{blue}$\ell$\color{black}};
	\draw (20.4,-7) node{\color{blue}$\ell$\color{black}};
	\braidbox{21.2}{23.6}{6.3}{7.6}{};
	\draw (22.4,7) node{$(\ell-1)^2$};
	\braidbox{21.2}{23.6}{-6.3}{-7.6}{};
	\draw (22.4,-7) node{$(\ell-1)^2$};
	\braidbox{26}{28.3}{6.3}{7.6}{};
	\draw (27.2,7) node{$(i+1)^2$};
	\braidbox{26}{28.3}{-6.3}{-7.6}{};
	\draw (27.2,-7) node{$(i+1)^2$};
	\draw (29.2,7) node{$i$};
	\draw (29.2,-7) node{$i$};
	\draw (32.5,7) node{$j+1$};
	\draw (32.5,-7) node{$j+1$};
	\draw (34,7) node{$j$};
	\draw (34,-7) node{$j$};
	\draw (35.5,7) node{$j-1$};
	\draw (35.5,-7) node{$j-1$};
	\draw (38.7,7) node{$1$};
	\draw (38.7,-7) node{$1$};
	\redbraidbox{39.5}{40.5}{6.3}{7.6}{};
	\draw (40,7) node{$\color{red}0\,\,0\color{black}$};
	\redbraidbox{39.5}{40.5}{-6.3}{-7.6}{};
	\draw (40,-7) node{$\color{red}0\,\,0\color{black}$};
	\draw (41.3,7) node{$1$};
	\draw (41.3,-7) node{$1$};
	\draw (44.7,7) node{$j-1$};
	\draw (44.7,-7) node{$j-1$};
	\draw (46.2,7) node{$j$};
	\draw (46.2,-7) node{$j$};
	\draw (47.7,7) node{$j+1$};
	\draw (47.7,-7) node{$j+1$};
	\draw (50.9,7) node{$i$};
	\draw (50.9,-7) node{$i$};
	\draw (14.7,6.3)--(39.5,0)--(14.7,-6.3);
	\draw (46.2,-6.3)--(40.4,0)--(46.2,6.3);
	\draw (16.2,6.3)--(41.4,0)--(47.7,-6.3);
	\draw (16.2,-6.3)--(41.4,0)--(47.7,6.3);
	\draw(19.5,6.3)--(44.7,0)--(50.9,-6.3);
	\draw(19.5,-6.3)--(44.7,0)--(50.9,6.3);
	\draw (44.7,6.3)--(39.5,0)--(44.7,-6.3);
	\draw (41.3,6.3)--(35.6,0)--(41.3,-6.3);
	\draw[red](40.4,6.3)--(34.7,0)--(40.4,-6.3);
	\draw[red](39.6,6.3)--(33.9,0)--(39.6,-6.3);
	\draw (38.7,6.3)--(32.9,0)--(38.7,-6.3);
	\draw (35.5,6.3)--(30.1,0)--(35.5,-6.3);
	\draw (34,6.3)--(28.8,0)--(34,-6.3);
	\draw (32.5,6.3)--(27.4,0)--(32.5,-6.3);
	\draw (29.2,6.3)--(24.3,0)--(29.2,-6.3);
	\draw (28,6.3)--(23,0)--(28,-6.3);
	\draw (26.5,6.3)--(21.5,0)--(26.5,-6.3);
	\draw (23.2,6.3)--(18.2,0)--(23.2,-6.3);
	\draw (21.7,6.3)--(16.7,0)--(21.7,-6.3);
	\draw[blue] (20.3,6.3)--(15.4,0)--(20.3,-6.3);
	\end{braid}.
	}
\end{align*}
Using quadratic and dot-crossing relations we now get
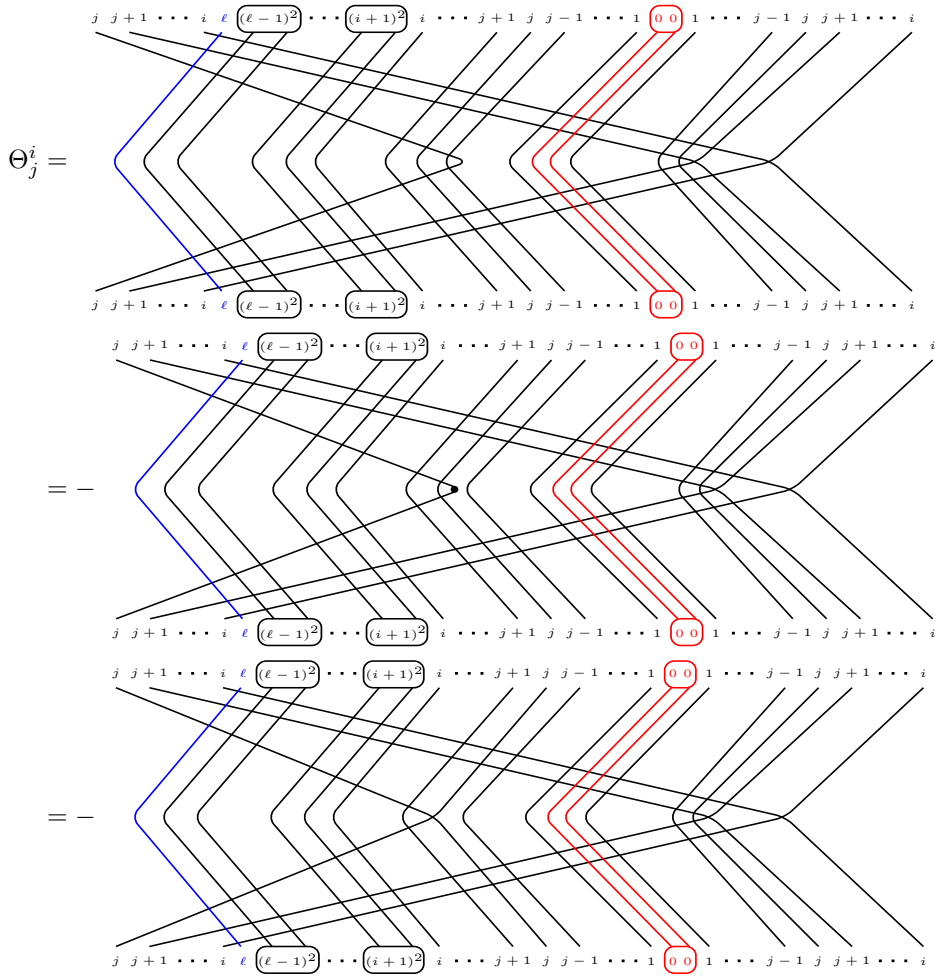
\begin{align*}
\Theta^i_j&=
\resizebox{115mm}{22mm}{
\begin{braid}\tikzset{baseline=-.3em}
        \draw[dots] (17.5,7)--(19.2,7);
	\draw[dots] (17.5,-7)--(19.2,-7);
	\draw[dots] (24.2,7)--(25.7,7);
	\draw[dots] (24.2,-7)--(25.7,-7);
	\draw[dots] (29.9,7)--(31.6,7);
	\draw[dots] (29.9,-7)--(31.6,-7);
	\draw[dots] (36.8,7)--(38.5,7);
	\draw[dots] (36.8,-7)--(38.5,-7);
	\draw[dots] (42,7)--(43.7,7);
	\draw[dots] (42,-7)--(43.7,-7);
	\draw[dots] (48.9,7)--(50.6,7);
	\draw[dots] (48.9,-7)--(50.6,-7);
	\draw (14.7,7) node{$j$};
	\draw (14.7,-7) node{$j$};
	\draw (16.2,7) node{$j+1$};
	\draw (16.2,-7) node{$j+1$};
	\draw (19.5,7) node{$i$};
	\draw (19.5,-7) node{$i$};
	\draw (20.4,7) node{\color{blue}$\ell$\color{black}};
	\draw (20.4,-7) node{\color{blue}$\ell$\color{black}};
	\braidbox{21.2}{23.6}{6.3}{7.6}{};
	\draw (22.4,7) node{$(\ell-1)^2$};
	\braidbox{21.2}{23.6}{-6.3}{-7.6}{};
	\draw (22.4,-7) node{$(\ell-1)^2$};
	\braidbox{26}{28.3}{6.3}{7.6}{};
	\draw (27.2,7) node{$(i+1)^2$};
	\braidbox{26}{28.3}{-6.3}{-7.6}{};
	\draw (27.2,-7) node{$(i+1)^2$};
	\draw (29.2,7) node{$i$};
	\draw (29.2,-7) node{$i$};
	\draw (32.5,7) node{$j+1$};
	\draw (32.5,-7) node{$j+1$};
	\draw (34,7) node{$j$};
	\draw (34,-7) node{$j$};
	\draw (35.5,7) node{$j-1$};
	\draw (35.5,-7) node{$j-1$};
	\draw (38.7,7) node{$1$};
	\draw (38.7,-7) node{$1$};
	\redbraidbox{39.5}{40.5}{6.3}{7.6}{};
	\draw (40,7) node{$\color{red}0\,\,0\color{black}$};
	\redbraidbox{39.5}{40.5}{-6.3}{-7.6}{};
	\draw (40,-7) node{$\color{red}0\,\,0\color{black}$};
	\draw (41.3,7) node{$1$};
	\draw (41.3,-7) node{$1$};
	\draw (44.7,7) node{$j-1$};
	\draw (44.7,-7) node{$j-1$};
	\draw (46.2,7) node{$j$};
	\draw (46.2,-7) node{$j$};
	\draw (47.7,7) node{$j+1$};
	\draw (47.7,-7) node{$j+1$};
	\draw (50.9,7) node{$i$};
	\draw (50.9,-7) node{$i$};
	\draw (14.7,6.3)--(31.2,0)--(14.7,-6.3);
	\draw (46.2,-6.3)--(40.4,0)--(46.2,6.3);
	\draw (16.2,6.3)--(41.4,0)--(47.7,-6.3);
	\draw (16.2,-6.3)--(41.4,0)--(47.7,6.3);
	\draw(19.5,6.3)--(44.7,0)--(50.9,-6.3);
	\draw(19.5,-6.3)--(44.7,0)--(50.9,6.3);
	\draw (44.7,6.3)--(39.5,0)--(44.7,-6.3);
	\draw (41.3,6.3)--(35.6,0)--(41.3,-6.3);
	\draw[red](40.4,6.3)--(34.7,0)--(40.4,-6.3);
	\draw[red](39.6,6.3)--(33.9,0)--(39.6,-6.3);
	\draw (38.7,6.3)--(32.9,0)--(38.7,-6.3);
	\draw (35.5,6.3)--(30.1,0)--(35.5,-6.3);
	\draw (34,6.3)--(28.8,0)--(34,-6.3);
	\draw (32.5,6.3)--(27.4,0)--(32.5,-6.3);
	\draw (29.2,6.3)--(24.3,0)--(29.2,-6.3);
	\draw (28,6.3)--(23,0)--(28,-6.3);
	\draw (26.5,6.3)--(21.5,0)--(26.5,-6.3);
	\draw (23.2,6.3)--(18.2,0)--(23.2,-6.3);
	\draw (21.7,6.3)--(16.7,0)--(21.7,-6.3);
	\draw[blue] (20.3,6.3)--(15.4,0)--(20.3,-6.3);
	\end{braid}
	}
	\\
	&=-
\resizebox{115mm}{22mm}{
\begin{braid}\tikzset{baseline=-.3em}
        \draw[dots] (17.5,7)--(19.2,7);
	\draw[dots] (17.5,-7)--(19.2,-7);
	\draw[dots] (24.2,7)--(25.7,7);
	\draw[dots] (24.2,-7)--(25.7,-7);
	\draw[dots] (29.9,7)--(31.6,7);
	\draw[dots] (29.9,-7)--(31.6,-7);
	\draw[dots] (36.8,7)--(38.5,7);
	\draw[dots] (36.8,-7)--(38.5,-7);
	\draw[dots] (42,7)--(43.7,7);
	\draw[dots] (42,-7)--(43.7,-7);
	\draw[dots] (48.9,7)--(50.6,7);
	\draw[dots] (48.9,-7)--(50.6,-7);
	\draw (14.7,7) node{$j$};
	\draw (14.7,-7) node{$j$};
	\draw (16.2,7) node{$j+1$};
	\draw (16.2,-7) node{$j+1$};
	\draw (19.5,7) node{$i$};
	\draw (19.5,-7) node{$i$};
	\draw (20.4,7) node{\color{blue}$\ell$\color{black}};
	\draw (20.4,-7) node{\color{blue}$\ell$\color{black}};
	\braidbox{21.2}{23.6}{6.3}{7.6}{};
	\draw (22.4,7) node{$(\ell-1)^2$};
	\braidbox{21.2}{23.6}{-6.3}{-7.6}{};
	\draw (22.4,-7) node{$(\ell-1)^2$};
	\braidbox{26}{28.3}{6.3}{7.6}{};
	\draw (27.2,7) node{$(i+1)^2$};
	\braidbox{26}{28.3}{-6.3}{-7.6}{};
	\draw (27.2,-7) node{$(i+1)^2$};
	\draw (29.2,7) node{$i$};
	\draw (29.2,-7) node{$i$};
	\draw (32.5,7) node{$j+1$};
	\draw (32.5,-7) node{$j+1$};
	\draw (34,7) node{$j$};
	\draw (34,-7) node{$j$};
	\draw (35.5,7) node{$j-1$};
	\draw (35.5,-7) node{$j-1$};
	\draw (38.7,7) node{$1$};
	\draw (38.7,-7) node{$1$};
	\redbraidbox{39.5}{40.5}{6.3}{7.6}{};
	\draw (40,7) node{$\color{red}0\,\,0\color{black}$};
	\redbraidbox{39.5}{40.5}{-6.3}{-7.6}{};
	\draw (40,-7) node{$\color{red}0\,\,0\color{black}$};
	\draw (41.3,7) node{$1$};
	\draw (41.3,-7) node{$1$};
	\draw (44.7,7) node{$j-1$};
	\draw (44.7,-7) node{$j-1$};
	\draw (46.2,7) node{$j$};
	\draw (46.2,-7) node{$j$};
	\draw (47.7,7) node{$j+1$};
	\draw (47.7,-7) node{$j+1$};
	\draw (50.9,7) node{$i$};
	\draw (50.9,-7) node{$i$};
	\draw (14.7,6.3)--(30,0)--(14.7,-6.3);
	\blackdot (29.7,0);
	\draw (46.2,-6.3)--(40.4,0)--(46.2,6.3);
	\draw (16.2,6.3)--(41.4,0)--(47.7,-6.3);
	\draw (16.2,-6.3)--(41.4,0)--(47.7,6.3);
	\draw(19.5,6.3)--(44.7,0)--(50.9,-6.3);
	\draw(19.5,-6.3)--(44.7,0)--(50.9,6.3);
	\draw (44.7,6.3)--(39.5,0)--(44.7,-6.3);
	\draw (41.3,6.3)--(35.6,0)--(41.3,-6.3);
	\draw[red](40.4,6.3)--(34.7,0)--(40.4,-6.3);
	\draw[red](39.6,6.3)--(33.9,0)--(39.6,-6.3);
	\draw (38.7,6.3)--(32.9,0)--(38.7,-6.3);
	\draw (35.5,6.3)--(30.1,0)--(35.5,-6.3);
	\draw (34,6.3)--(28.8,0)--(34,-6.3);
	\draw (32.5,6.3)--(27.4,0)--(32.5,-6.3);
	\draw (29.2,6.3)--(24.3,0)--(29.2,-6.3);
	\draw (28,6.3)--(23,0)--(28,-6.3);
	\draw (26.5,6.3)--(21.5,0)--(26.5,-6.3);
	\draw (23.2,6.3)--(18.2,0)--(23.2,-6.3);
	\draw (21.7,6.3)--(16.7,0)--(21.7,-6.3);
	\draw[blue] (20.3,6.3)--(15.4,0)--(20.3,-6.3);
	\end{braid}
	}
	\\
	&=-
\resizebox{115mm}{22mm}{
\begin{braid}\tikzset{baseline=-.3em}
        \draw[dots] (17.5,7)--(19.2,7);
	\draw[dots] (17.5,-7)--(19.2,-7);
	\draw[dots] (24.2,7)--(25.7,7);
	\draw[dots] (24.2,-7)--(25.7,-7);
	\draw[dots] (29.9,7)--(31.6,7);
	\draw[dots] (29.9,-7)--(31.6,-7);
	\draw[dots] (36.8,7)--(38.5,7);
	\draw[dots] (36.8,-7)--(38.5,-7);
	\draw[dots] (42,7)--(43.7,7);
	\draw[dots] (42,-7)--(43.7,-7);
	\draw[dots] (48.9,7)--(50.6,7);
	\draw[dots] (48.9,-7)--(50.6,-7);
	\draw (14.7,7) node{$j$};
	\draw (14.7,-7) node{$j$};
	\draw (16.2,7) node{$j+1$};
	\draw (16.2,-7) node{$j+1$};
	\draw (19.5,7) node{$i$};
	\draw (19.5,-7) node{$i$};
	\draw (20.4,7) node{\color{blue}$\ell$\color{black}};
	\draw (20.4,-7) node{\color{blue}$\ell$\color{black}};
	\braidbox{21.2}{23.6}{6.3}{7.6}{};
	\draw (22.4,7) node{$(\ell-1)^2$};
	\braidbox{21.2}{23.6}{-6.3}{-7.6}{};
	\draw (22.4,-7) node{$(\ell-1)^2$};
	\braidbox{26}{28.3}{6.3}{7.6}{};
	\draw (27.2,7) node{$(i+1)^2$};
	\braidbox{26}{28.3}{-6.3}{-7.6}{};
	\draw (27.2,-7) node{$(i+1)^2$};
	\draw (29.2,7) node{$i$};
	\draw (29.2,-7) node{$i$};
	\draw (32.5,7) node{$j+1$};
	\draw (32.5,-7) node{$j+1$};
	\draw (34,7) node{$j$};
	\draw (34,-7) node{$j$};
	\draw (35.5,7) node{$j-1$};
	\draw (35.5,-7) node{$j-1$};
	\draw (38.7,7) node{$1$};
	\draw (38.7,-7) node{$1$};
	\redbraidbox{39.5}{40.5}{6.3}{7.6}{};
	\draw (40,7) node{$\color{red}0\,\,0\color{black}$};
	\redbraidbox{39.5}{40.5}{-6.3}{-7.6}{};
	\draw (40,-7) node{$\color{red}0\,\,0\color{black}$};
	\draw (41.3,7) node{$1$};
	\draw (41.3,-7) node{$1$};
	\draw (44.7,7) node{$j-1$};
	\draw (44.7,-7) node{$j-1$};
	\draw (46.2,7) node{$j$};
	\draw (46.2,-7) node{$j$};
	\draw (47.7,7) node{$j+1$};
	\draw (47.7,-7) node{$j+1$};
	\draw (50.9,7) node{$i$};
	\draw (50.9,-7) node{$i$};
	\draw (14.7,6.3)--(29,0)--(34,-6.3);
	\draw (14.7,-6.3)--(29,0)--(34,6.3);
	\draw (46.2,-6.3)--(40.4,0)--(46.2,6.3);
	\draw (16.2,6.3)--(41.4,0)--(47.7,-6.3);
	\draw (16.2,-6.3)--(41.4,0)--(47.7,6.3);
	\draw(19.5,6.3)--(44.7,0)--(50.9,-6.3);
	\draw(19.5,-6.3)--(44.7,0)--(50.9,6.3);
	\draw (44.7,6.3)--(39.5,0)--(44.7,-6.3);
	\draw (41.3,6.3)--(35.6,0)--(41.3,-6.3);
	\draw[red](40.4,6.3)--(34.7,0)--(40.4,-6.3);
	\draw[red](39.6,6.3)--(33.9,0)--(39.6,-6.3);
	\draw (38.7,6.3)--(32.9,0)--(38.7,-6.3);
	\draw (35.5,6.3)--(30.1,0)--(35.5,-6.3);
	\draw (32.5,6.3)--(27.4,0)--(32.5,-6.3);
	\draw (29.2,6.3)--(24.3,0)--(29.2,-6.3);
	\draw (28,6.3)--(23,0)--(28,-6.3);
	\draw (26.5,6.3)--(21.5,0)--(26.5,-6.3);
	\draw (23.2,6.3)--(18.2,0)--(23.2,-6.3);
	\draw (21.7,6.3)--(16.7,0)--(21.7,-6.3);
	\draw[blue] (20.3,6.3)--(15.4,0)--(20.3,-6.3);
	\end{braid}.
	}
\end{align*}
We apply the braid relation (\ref{R7}) for 
$\begin{braid}\tikzset{baseline=.4em}
	\draw (0,0) node{$j$};
        \draw (1.1,0) node{$j+1$};
        \draw (2.2,0) node{$j$};
		\draw(0,0.5)--(2.2,1.5);
\draw(2.2,0.5)--(0,1.5);
\draw(1.1,0.5)--(0,1)--(1.1,1.5);
	\end{braid}$. 
The term corresponding to 
$\begin{braid}\tikzset{baseline=.4em}
	\draw (0,0) node{$j$};
        \draw (1.1,0) node{$j+1$};
        \draw (2.2,0) node{$j$};
		\draw(0,0.5)--(2.2,1.5);
\draw(2.2,0.5)--(0,1.5);
\draw(1.1,0.5)--(2.2,1)--(1.1,1.5);
	\end{braid}$
is $0$ since the word beginning with $\ell(\ell-1)^2\cdots(i+2)^2i\cdots (j-2)j$ is not cuspidal by Lemma~\ref{LCuspExpl}. So 
\begin{align*}
\Theta^i_j&=
\resizebox{117mm}{22mm}{
\begin{braid}\tikzset{baseline=-.3em}
        \draw[dots] (17.5,7)--(19.2,7);
	\draw[dots] (17.5,-7)--(19.2,-7);
	\draw[dots] (24.2,7)--(25.7,7);
	\draw[dots] (24.2,-7)--(25.7,-7);
	\draw[dots] (29.9,7)--(31.6,7);
	\draw[dots] (29.9,-7)--(31.6,-7);
	\draw[dots] (36.8,7)--(38.5,7);
	\draw[dots] (36.8,-7)--(38.5,-7);
	\draw[dots] (42,7)--(43.7,7);
	\draw[dots] (42,-7)--(43.7,-7);
	\draw[dots] (48.9,7)--(50.6,7);
	\draw[dots] (48.9,-7)--(50.6,-7);
	\draw (14.7,7) node{$j$};
	\draw (14.7,-7) node{$j$};
	\draw (16.2,7) node{$j+1$};
	\draw (16.2,-7) node{$j+1$};
	\draw (19.5,7) node{$i$};
	\draw (19.5,-7) node{$i$};
	\draw (20.4,7) node{\color{blue}$\ell$\color{black}};
	\draw (20.4,-7) node{\color{blue}$\ell$\color{black}};
	\braidbox{21.2}{23.6}{6.3}{7.6}{};
	\draw (22.4,7) node{$(\ell-1)^2$};
	\braidbox{21.2}{23.6}{-6.3}{-7.6}{};
	\draw (22.4,-7) node{$(\ell-1)^2$};
	\braidbox{26}{28.3}{6.3}{7.6}{};
	\draw (27.2,7) node{$(i+1)^2$};
	\braidbox{26}{28.3}{-6.3}{-7.6}{};
	\draw (27.2,-7) node{$(i+1)^2$};
	\draw (29.2,7) node{$i$};
	\draw (29.2,-7) node{$i$};
	\draw (32.5,7) node{$j+1$};
	\draw (32.5,-7) node{$j+1$};
	\draw (34,7) node{$j$};
	\draw (34,-7) node{$j$};
	\draw (35.5,7) node{$j-1$};
	\draw (35.5,-7) node{$j-1$};
	\draw (38.7,7) node{$1$};
	\draw (38.7,-7) node{$1$};
	\redbraidbox{39.5}{40.5}{6.3}{7.6}{};
	\draw (40,7) node{$\color{red}0\,\,0\color{black}$};
	\redbraidbox{39.5}{40.5}{-6.3}{-7.6}{};
	\draw (40,-7) node{$\color{red}0\,\,0\color{black}$};
	\draw (41.3,7) node{$1$};
	\draw (41.3,-7) node{$1$};
	\draw (44.7,7) node{$j-1$};
	\draw (44.7,-7) node{$j-1$};
	\draw (46.2,7) node{$j$};
	\draw (46.2,-7) node{$j$};
	\draw (47.7,7) node{$j+1$};
	\draw (47.7,-7) node{$j+1$};
	\draw (50.9,7) node{$i$};
	\draw (50.9,-7) node{$i$};
	\draw (14.7,6.3)--(28,0)--(14.7,-6.3);
	\draw (34,-6.3)--(29,0)--(34,6.3);
	\draw (46.2,-6.3)--(40.4,0)--(46.2,6.3);
	\draw (16.2,6.3)--(41.4,0)--(47.7,-6.3);
	\draw (16.2,-6.3)--(41.4,0)--(47.7,6.3);
	\draw(19.5,6.3)--(44.7,0)--(50.9,-6.3);
	\draw(19.5,-6.3)--(44.7,0)--(50.9,6.3);
	\draw (44.7,6.3)--(39.5,0)--(44.7,-6.3);
	\draw (41.3,6.3)--(35.6,0)--(41.3,-6.3);
	\draw[red](40.4,6.3)--(34.7,0)--(40.4,-6.3);
	\draw[red](39.6,6.3)--(33.9,0)--(39.6,-6.3);
	\draw (38.7,6.3)--(32.9,0)--(38.7,-6.3);
	\draw (35.5,6.3)--(30.1,0)--(35.5,-6.3);
	\draw (32.5,6.3)--(28.1,0)--(32.5,-6.3);
	\draw (29.2,6.3)--(24.3,0)--(29.2,-6.3);
	\draw (28,6.3)--(23,0)--(28,-6.3);
	\draw (26.5,6.3)--(21.5,0)--(26.5,-6.3);
	\draw (23.2,6.3)--(18.2,0)--(23.2,-6.3);
	\draw (21.7,6.3)--(16.7,0)--(21.7,-6.3);
	\draw[blue] (20.3,6.3)--(15.4,0)--(20.3,-6.3);
	\end{braid}
	}
	\\
	&=
\resizebox{117mm}{22mm}{
\begin{braid}\tikzset{baseline=-.3em}
        \draw[dots] (17.5,7)--(19.2,7);
	\draw[dots] (17.5,-7)--(19.2,-7);
	\draw[dots] (24.2,7)--(25.7,7);
	\draw[dots] (24.2,-7)--(25.7,-7);
	\draw[dots] (29.9,7)--(31.6,7);
	\draw[dots] (29.9,-7)--(31.6,-7);
	\draw[dots] (36.8,7)--(38.5,7);
	\draw[dots] (36.8,-7)--(38.5,-7);
	\draw[dots] (42,7)--(43.7,7);
	\draw[dots] (42,-7)--(43.7,-7);
	\draw[dots] (48.9,7)--(50.6,7);
	\draw[dots] (48.9,-7)--(50.6,-7);
	\draw (14.7,7) node{$j$};
	\draw (14.7,-7) node{$j$};
	\draw (16.2,7) node{$j+1$};
	\draw (16.2,-7) node{$j+1$};
	\draw (19.5,7) node{$i$};
	\draw (19.5,-7) node{$i$};
	\draw (20.4,7) node{\color{blue}$\ell$\color{black}};
	\draw (20.4,-7) node{\color{blue}$\ell$\color{black}};
	\braidbox{21.2}{23.6}{6.3}{7.6}{};
	\draw (22.4,7) node{$(\ell-1)^2$};
	\braidbox{21.2}{23.6}{-6.3}{-7.6}{};
	\draw (22.4,-7) node{$(\ell-1)^2$};
	\braidbox{26}{28.3}{6.3}{7.6}{};
	\draw (27.2,7) node{$(i+1)^2$};
	\braidbox{26}{28.3}{-6.3}{-7.6}{};
	\draw (27.2,-7) node{$(i+1)^2$};
	\draw (29.2,7) node{$i$};
	\draw (29.2,-7) node{$i$};
	\draw (32.5,7) node{$j+1$};
	\draw (32.5,-7) node{$j+1$};
	\draw (34,7) node{$j$};
	\draw (34,-7) node{$j$};
	\draw (35.5,7) node{$j-1$};
	\draw (35.5,-7) node{$j-1$};
	\draw (38.7,7) node{$1$};
	\draw (38.7,-7) node{$1$};
	\redbraidbox{39.5}{40.5}{6.3}{7.6}{};
	\draw (40,7) node{$\color{red}0\,\,0\color{black}$};
	\redbraidbox{39.5}{40.5}{-6.3}{-7.6}{};
	\draw (40,-7) node{$\color{red}0\,\,0\color{black}$};
	\draw (41.3,7) node{$1$};
	\draw (41.3,-7) node{$1$};
	\draw (44.7,7) node{$j-1$};
	\draw (44.7,-7) node{$j-1$};
	\draw (46.2,7) node{$j$};
	\draw (46.2,-7) node{$j$};
	\draw (47.7,7) node{$j+1$};
	\draw (47.7,-7) node{$j+1$};
	\draw (50.9,7) node{$i$};
	\draw (50.9,-7) node{$i$};
	\draw (14.7,6.3)--(14.7,-6.3);
	\draw (34,-6.3)--(29,0)--(34,6.3);
	\draw (46.2,-6.3)--(40.4,0)--(46.2,6.3);
	\draw (16.2,6.3)--(41.4,0)--(47.7,-6.3);
	\draw (16.2,-6.3)--(41.4,0)--(47.7,6.3);
	\draw(19.5,6.3)--(44.7,0)--(50.9,-6.3);
	\draw(19.5,-6.3)--(44.7,0)--(50.9,6.3);
	\draw (44.7,6.3)--(39.5,0)--(44.7,-6.3);
	\draw (41.3,6.3)--(35.6,0)--(41.3,-6.3);
	\draw[red](40.4,6.3)--(34.7,0)--(40.4,-6.3);
	\draw[red](39.6,6.3)--(33.9,0)--(39.6,-6.3);
	\draw (38.7,6.3)--(32.9,0)--(38.7,-6.3);
	\draw (35.5,6.3)--(30.1,0)--(35.5,-6.3);
	\draw (32.5,6.3)--(28,0)--(32.5,-6.3);
	\draw (29.2,6.3)--(24.4,0)--(29.2,-6.3);
	\draw (28,6.3)--(23.1,0)--(28,-6.3);
	\draw (26.5,6.3)--(21.6,0)--(26.5,-6.3);
	\draw (23.2,6.3)--(18.3,0)--(23.2,-6.3);
	\draw (21.7,6.3)--(16.8,0)--(21.7,-6.3);
	\draw[blue] (20.3,6.3)--(15.5,0)--(20.3,-6.3);
	\end{braid},
	}
\end{align*}
which is $\Theta^i_{j+1}$ as required. 
\end{proof}

\section{The element \texorpdfstring{$\Upsilon$}{}}\label{sec:upsilon}

Define the element $\Upsilon\in R_{2\de}$ via 
\begin{equation}\label{EUpsilon}
\Upsilon=(\psi_p\psi_{p+1}\cdots \psi_{2p-1})(\psi_{p-1}\psi_p\cdots \psi_{2p-2})\cdots(\psi_{2}\psi_{3}\cdots \psi_{p+1})(\psi_{2}\cdots \psi_p).
\index{$\Upsilon$}
\end{equation}
Thus, in terms of Khovanov-Lauda diagrams, we have 
$$
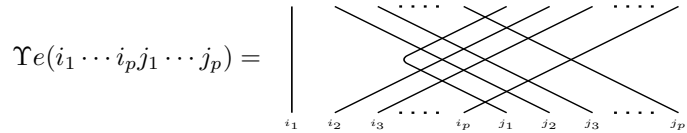
\Upsilon e(i_1\cdots i_pj_1\cdots j_p)=
\resizebox{60mm}{8mm}{
\begin{braid}\tikzset{baseline=-.3em}
	\draw(0,-3)node[below]{$i_1$}--(0,3);
	\draw(2,-3)node[below]{$i_2$}--(12,3);
	\draw(4,-3)node[below]{$i_3$}--(14,3);
	\draw[dots] (5,-3)--(7,-3);
	\draw[dots] (5,3)--(7,3);
	\draw(8,-3)node[below]{$i_p$}--(18,3);
	\draw(10,-3)node[below]{$j_1$}--(5,0)--(10,3);
	\draw(12,-3)node[below]{$j_2$}--(2,3);
	\draw(14,-3)node[below]{$j_3$}--(4,3);
	\draw[dots] (15,-3)--(17,-3);
	\draw[dots] (15,3)--(17,3);
	\draw(18,-3)node[below]{$j_p$}--(8,3);
	\end{braid}.
	}
$$
For example, we have the element 
$\ga^{j,i}\Upsilon\ga^{i,j}$ of $B_2$
$$
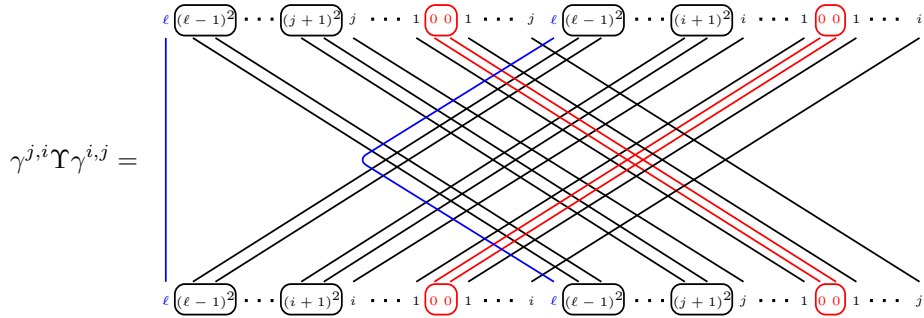
\ga^{j,i}\Upsilon\ga^{i,j}=
\resizebox{108mm}{22mm}{
\begin{braid}\tikzset{baseline=-.3em}
	\draw (0,6) node{\color{blue}$\ell$\color{black}};
	\braidbox{0.6}{2.9}{5.3}{6.6}{};
	\draw (1.8,6) node{$(\ell-1)^2$};
	\draw[dots] (3.5,6)--(5,6);
	\braidbox{5.3}{7.6}{5.3}{6.6}{};
	\draw (6.5,6) node{$(j+1)^2$};
	\draw (8.3,6) node{$j$};
	\draw[dots] (9.1,6)--(10.4,6);
	\draw (11.1,6) node{$1$};
	\redbraidbox{11.7}{12.7}{5.3}{6.6}{};
	\draw (12.2,6) node{$\color{red}0\,\,0\color{black}$};
	\draw (13.4,6) node{$1$};
	\draw[dots] (14.1,6)--(15.5,6);
	\draw (16.2,6) node{$j$};
	\draw (0,-6) node{\color{blue}$\ell$\color{black}};
	\braidbox{0.6}{2.9}{-5.3}{-6.6}{};
	\draw (1.8,-6) node{$(\ell-1)^2$};
	\draw[dots] (3.5,-6)--(5,-6);
	\braidbox{5.3}{7.6}{-5.3}{-6.6}{};
	\draw (6.5,-6) node{$(i+1)^2$};
	\draw (8.3,-6) node{$i$};
	\draw[dots] (9.1,-6)--(10.4,-6);
	\draw (11.1,-6) node{$1$};
	\redbraidbox{11.7}{12.7}{-5.3}{-6.6}{};
	\draw (12.2,-6) node{$\color{red}0\,\,0\color{black}$};
	\draw (13.4,-6) node{$1$};
	\draw[dots] (14.1,-6)--(15.5,-6);
	\draw (16.2,-6) node{$i$};
\draw (17.2,6) node{\color{blue}$\ell$\color{black}};
	\braidbox{17.8}{20.1}{5.3}{6.6}{};
	\draw (19,6) node{$(\ell-1)^2$};
	\draw[dots] (20.8,6)--(22.3,6);
	\braidbox{22.6}{24.9}{5.3}{6.6}{};
	\draw (23.8,6) node{$(i+1)^2$};
	\draw (25.6,6) node{$i$};
	\draw[dots] (26.3,6)--(28,6);
	\draw (28.3,6) node{$1$};
	\redbraidbox{29}{29.9}{5.3}{6.6}{};
	\draw (29.4,6) node{$\color{red}0\,\,0\color{black}$};
	\draw (30.6,6) node{$1$};
	\draw[dots] (31.2,6)--(32.6,6);
	\draw (33.4,6) node{$i$};	
	\draw (17.2,-6) node{\color{blue}$\ell$\color{black}};
	\braidbox{17.8}{20.1}{-5.3}{-6.6}{};
	\draw (19,-6) node{$(\ell-1)^2$};
	\draw[dots] (20.8,-6)--(22.3,-6);
	\braidbox{22.6}{24.9}{-5.3}{-6.6}{};
	\draw (23.8,-6) node{$(j+1)^2$};
	\draw (25.6,-6) node{$j$};
	\draw[dots] (26.3,-6)--(28,-6);
	\draw (28.3,-6) node{$1$};
	\redbraidbox{29}{29.9}{-5.3}{-6.6}{};
	\draw (29.4,-6) node{$\color{red}0\,\,0\color{black}$};
	\draw (30.6,-6) node{$1$};
	\draw[dots] (31.2,-6)--(32.6,-6);
	\draw (33.4,-6) node{$j$};	
	\draw[blue](0,-5.2)--(0,5.2);
	\draw(1.2,-5.2)--(18.3,5.2);
	\draw(2.2,-5.2)--(19.3,5.2);
	\draw(1.2,5.2)--(18.3,-5.2);
	\draw(2.2,5.2)--(19.3,-5.2);
	
	\draw(6,-5.2)--(23,5.2);
	\draw(7,-5.2)--(24,5.2);
	\draw(6,5.2)--(23,-5.2);
	\draw(7,5.2)--(24,-5.2);
	
	\draw(8.3,-5.2)--(25.6,5.2);
	\draw(8.3,5.2)--(25.6,-5.2);
	
	\draw(11.1,-5.2)--(28.3,5.2);
	\draw(11.1,5.2)--(28.3,-5.2);
	
	\draw[red](11.9,-5.2)--(29,5.2);
	\draw[red](12.6,-5.2)--(29.7,5.2);
	\draw[red](11.9,5.2)--(29,-5.2);
	\draw[red](12.6,5.2)--(29.7,-5.2);
	
	\draw(13.4,-5.2)--(30.6,5.2);
	\draw(16.2,-5.2)--(33.4,5.2);
	\draw(13.4,5.2)--(30.6,-5.2);
	\draw(16.2,5.2)--(33.4,-5.2);
	
	\draw[blue](17.2,5.2)--(8.5,0)--(17.2,-5.2);
	\end{braid}.
	}
$$
Here, and throughout this section, we follow our usual convention (\ref{EDoubleRed}) for the order of odd crossings.

For the following lemma, recall the definition of the $P^{(t)}$'s from $\S$\ref{sec:uptau}.

\begin{Lemma}\label{lem:poly_P}
Let $i,j \in J$ and $2 \leq t \leq p$. Then, in ${\bar R}_{2\de}$,
$$P^{(t)}\ga^{i,j} =
\begin{cases}
(y_{r^j_g}-y_{t+p})(y_{s^j_g}-y_{t+p})P^{(t-1)} \ga^{i,j} &\text{if }g \neq 0,\\
(y_{r^j_0}^2-y_{t+p}^2)(y_{s^j_0}^2-y_{t+p}^2)P^{(t-1)} \ga^{i,j} &\text{if }g = 0,
\end{cases}
$$
where $g:= g^i_t$. Furthermore, in ${\bar R}_{2\de}$,
$$P^{(1)}\ga^{i,j} = \Upsilon \ga^{i,j} + (y_1 - y_{p+1})P^{(0)}\ga^{i,j}.
$$
\end{Lemma}

\begin{proof}
For the first statement, we prove a stronger statement. Namely that, in ${\bar R}_{2\de}$,
$$P^{(t)}{\hat \ga}^{i,j} =
\begin{cases}
(y_{r^j_g}-y_{t+p})(y_{s^j_g}-y_{t+p})P^{(t-1)}{\hat \ga^{i,j}} &\text{if }g \neq 0,\\
(y_{r^j_0}^2-y_{t+p}^2)(y_{s^j_0}^2-y_{t+p}^2)P^{(t-1)}{\hat \ga^{i,j}} &\text{if }g = 0,
\end{cases}
$$
where ${\hat \ga}^{i,j}:=e({\hat \bg}^i {\hat \bg}^j)$. Certainly this implies the result. We first assume that $g \neq 0$. Then
$$
P^{(t)}{\hat \ga}^{i,j}=
\resizebox{104mm}{23mm}{
\begin{braid}\tikzset{baseline=-.3em}
	\draw[blue] (0,6) node{$\ell$};
	\draw (1.5,6) node{$\ell-1$};
	\draw[dots] (3,6)--(4.5,6);
	\draw[darkgreen] (5.7,6) node{$g$};
	\draw[dots] (6.9,6)--(8.4,6);
	\draw[darkgreen] (9.6,6) node{$g$};
	\draw[dots] (10.8,6)--(12.3,6);
	\draw (13.5,6) node{$j$};
	\draw[blue] (15,6) node{$\ell$};
    \draw (16.5,6) node{$\ell-1$};
	\draw[dots] (18,6)--(19.5,6);
	\draw (21,6) node{$ $};
	\draw[darkgreen] (22.5,6) node{$g$};
	\draw (24,6) node{$ $};
	\draw[dots] (25.5,6)--(27,6);
	\draw (28.5,6) node{$i$};

	\draw[blue] (0,-6) node{$\ell$};
	\draw (1.5,-6) node{$\ell-1$};
	\draw[dots] (3,-6)--(4.5,-6);
	\draw (6,-6) node{$ $};
	\draw[darkgreen] (7.5,-6) node{$g$};
	\draw (9,-6) node{$ $};
	\draw[dots] (10.5,-6)--(12,-6);
	\draw (13.5,-6) node{$i$};
	\draw[blue] (15,-6) node{$\ell$};
	\draw (16.5,-6) node{$\ell-1$};
	\draw[dots] (18,-6)--(19.5,-6);
	\draw[darkgreen] (20.7,-6) node{$g$};
	\draw[dots] (21.9,-6)--(23.4,-6);
	\draw[darkgreen] (24.6,-6) node{$g$};
	\draw[dots] (25.8,-6)--(27.3,-6);
	\draw (28.5,-6) node{$j$};

	\draw[blue](0,-5.2)--(15,5.2);
	\draw(1.5,-5.2)--(16.5,5.2);
	\draw(6,-5.2)--(21,5.2);
	\draw[darkgreen](7.5,-5.2)--(22.5,5.2);
	\draw(9,-5.2)--(24,5.2);
	\draw(13.5,-5.2)--(28.5,5.2);

	\draw[blue](0,5.2)--(15,-5.2);
	\draw(1.5,5.2)--(16.5,-5.2);
	\draw[darkgreen](5.7,5.2)--(20.7,-5.2);
	\draw[darkgreen](9.6,5.2)--(24.6,-5.2);
	\draw(13.5,5.2)--(28.5,-5.2);

	\filldraw[color=blue!60, fill=black!5, thin] (7.5,0) circle (0.15);
	\filldraw[color=black!60, fill=black!5, thin] (8.25,0.52) circle (0.15);
	\filldraw[color=black!60, fill=black!5, thin] (10.35,1.976) circle (0.15);
	\filldraw[color=black!60, fill=black!5, thin] (12.3,3.328) circle (0.15);
	\filldraw[color=black!60, fill=black!5, thin] (14.25,4.68) circle (0.15);

	\filldraw[color=black!60, fill=black!5, thin] (8.25,-0.52) circle (0.15);
	\filldraw[color=black!60, fill=black!5, thin] (9,0) circle (0.15);
	\filldraw[color=black!60, fill=black!5, thin] (11.1,1.456) circle (0.15);
	\filldraw[color=black!60, fill=black!5, thin] (13.05,2.808) circle (0.15);
	\filldraw[color=black!60, fill=black!5, thin] (15,4.16) circle (0.15);

	\filldraw[color=black!60, fill=black!5, thin] (10.5,-2.08) circle (0.15);
	\filldraw[color=black!60, fill=black!5, thin] (11.25,-1.56) circle (0.15);
	\filldraw[color=black!60, fill=black!5, thin] (13.35,-0.104) circle (0.15);
	\filldraw[color=black!60, fill=black!5, thin] (15.3,1.248) circle (0.15);
	\filldraw[color=black!60, fill=black!5, thin] (17.25,2.6) circle (0.15);

	\filldraw[color=black!60, fill=black!5, thin] (11.25,-2.6) circle (0.15);
	\filldraw[color=black!60, fill=black!5, thin] (12,-2.08) circle (0.15);
	\filldraw[color=darkgreen!60, fill=black!5, thin] (14.1,-0.624) circle (0.15);
	\filldraw[color=darkgreen!60, fill=black!5, thin] (16.05,0.728) circle (0.15);
	\filldraw[color=black!60, fill=black!5, thin] (18,2.08) circle (0.15);
	\end{braid}.
	}
$$
Using (\ref{EOpening}), we have
$$
P^{(t)}{\hat \ga}^{i,j}=
\resizebox{104mm}{23mm}
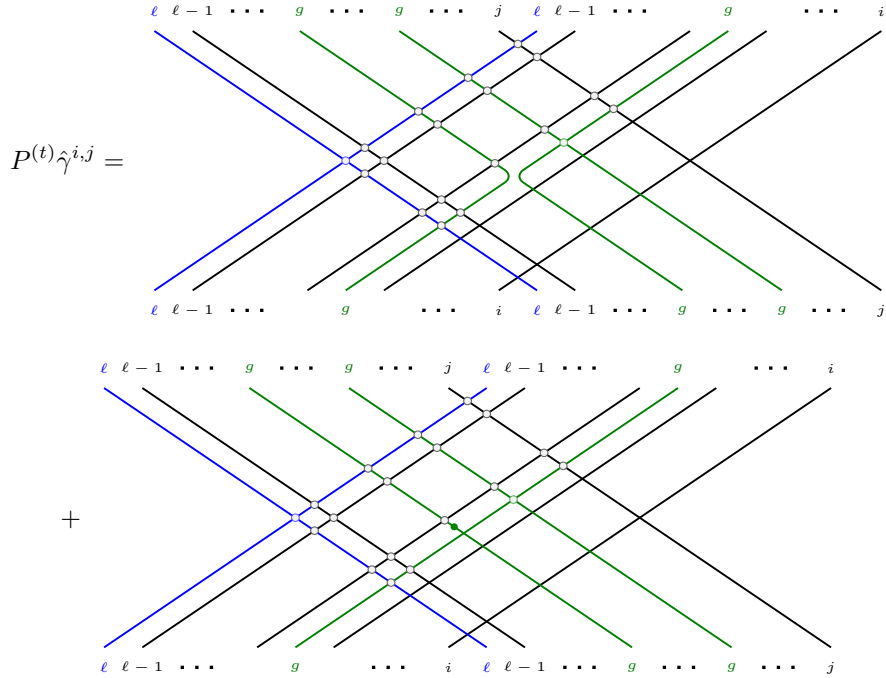
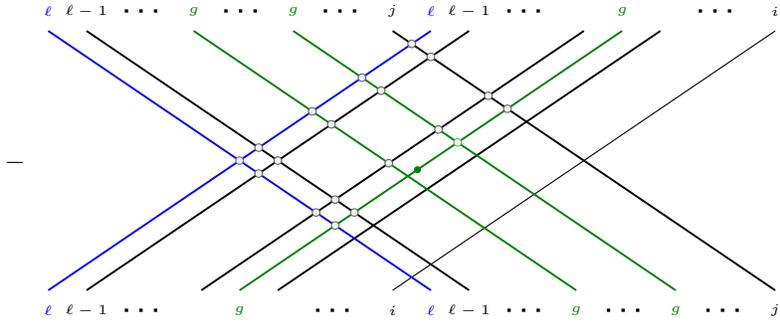
{
\begin{braid}\tikzset{baseline=-.3em}
	\draw[blue] (0,6) node{$\ell$};
	\draw (1.5,6) node{$\ell-1$};
	\draw[dots] (3,6)--(4.5,6);
	\draw[darkgreen] (5.7,6) node{$g$};
	\draw[dots] (6.9,6)--(8.4,6);
	\draw[darkgreen] (9.6,6) node{$g$};
	\draw[dots] (10.8,6)--(12.3,6);
	\draw (13.5,6) node{$j$};
	\draw[blue] (15,6) node{$\ell$};
    \draw (16.5,6) node{$\ell-1$};
	\draw[dots] (18,6)--(19.5,6);
	\draw (21,6) node{$ $};
	\draw[darkgreen] (22.5,6) node{$g$};
	\draw (24,6) node{$ $};
	\draw[dots] (25.5,6)--(27,6);
	\draw (28.5,6) node{$i$};

	\draw[blue] (0,-6) node{$\ell$};
	\draw (1.5,-6) node{$\ell-1$};
	\draw[dots] (3,-6)--(4.5,-6);
	\draw (6,-6) node{$ $};
	\draw[darkgreen] (7.5,-6) node{$g$};
	\draw (9,-6) node{$ $};
	\draw[dots] (10.5,-6)--(12,-6);
	\draw (13.5,-6) node{$i$};
	\draw[blue] (15,-6) node{$\ell$};
	\draw (16.5,-6) node{$\ell-1$};
	\draw[dots] (18,-6)--(19.5,-6);
	\draw[darkgreen] (20.7,-6) node{$g$};
	\draw[dots] (21.9,-6)--(23.4,-6);
	\draw[darkgreen] (24.6,-6) node{$g$};
	\draw[dots] (25.8,-6)--(27.3,-6);
	\draw (28.5,-6) node{$j$};

	\draw[blue](0,-5.2)--(15,5.2);
	\draw(1.5,-5.2)--(16.5,5.2);
	\draw(6,-5.2)--(21,5.2);
	\draw(9,-5.2)--(24,5.2);
	\draw(13.5,-5.2)--(28.5,5.2);

	\draw[blue](0,5.2)--(15,-5.2);
	\draw(1.5,5.2)--(16.5,-5.2);
	\draw[darkgreen](9.6,5.2)--(24.6,-5.2);
	\draw(13.5,5.2)--(28.5,-5.2);
	
	\draw[darkgreen](7.5,-5.2)--(14.1,-0.624)--(5.7,5.2);
    \draw[darkgreen](20.7,-5.2)--(14.1,-0.624)--(22.5,5.2);

	\filldraw[color=blue!60, fill=black!5, thin] (7.5,0) circle (0.15);
	\filldraw[color=black!60, fill=black!5, thin] (8.25,0.52) circle (0.15);
	\filldraw[color=black!60, fill=black!5, thin] (10.35,1.976) circle (0.15);
	\filldraw[color=black!60, fill=black!5, thin] (12.3,3.328) circle (0.15);
	\filldraw[color=black!60, fill=black!5, thin] (14.25,4.68) circle (0.15);

	\filldraw[color=black!60, fill=black!5, thin] (8.25,-0.52) circle (0.15);
	\filldraw[color=black!60, fill=black!5, thin] (9,0) circle (0.15);
	\filldraw[color=black!60, fill=black!5, thin] (11.1,1.456) circle (0.15);
	\filldraw[color=black!60, fill=black!5, thin] (13.05,2.808) circle (0.15);
	\filldraw[color=black!60, fill=black!5, thin] (15,4.16) circle (0.15);

	\filldraw[color=black!60, fill=black!5, thin] (10.5,-2.08) circle (0.15);
	\filldraw[color=black!60, fill=black!5, thin] (11.25,-1.56) circle (0.15);
	\filldraw[color=black!60, fill=black!5, thin] (13.35,-0.104) circle (0.15);
	\filldraw[color=black!60, fill=black!5, thin] (15.3,1.248) circle (0.15);
	\filldraw[color=black!60, fill=black!5, thin] (17.25,2.6) circle (0.15);

	\filldraw[color=black!60, fill=black!5, thin] (11.25,-2.6) circle (0.15);
	\filldraw[color=black!60, fill=black!5, thin] (12,-2.08) circle (0.15);
	\filldraw[color=darkgreen!60, fill=black!5, thin] (16.05,0.728) circle (0.15);
	\filldraw[color=black!60, fill=black!5, thin] (18,2.08) circle (0.15);
	\end{braid}
	}
$$

$$
+\resizebox{104mm}{23mm}{
\begin{braid}\tikzset{baseline=-.3em}
	\draw[blue] (0,6) node{$\ell$};
	\draw (1.5,6) node{$\ell-1$};
	\draw[dots] (3,6)--(4.5,6);
	\draw[darkgreen] (5.7,6) node{$g$};
	\draw[dots] (6.9,6)--(8.4,6);
	\draw[darkgreen] (9.6,6) node{$g$};
	\draw[dots] (10.8,6)--(12.3,6);
	\draw (13.5,6) node{$j$};
	\draw[blue] (15,6) node{$\ell$};
    \draw (16.5,6) node{$\ell-1$};
	\draw[dots] (18,6)--(19.5,6);
	\draw (21,6) node{$ $};
	\draw[darkgreen] (22.5,6) node{$g$};
	\draw (24,6) node{$ $};
	\draw[dots] (25.5,6)--(27,6);
	\draw (28.5,6) node{$i$};

	\draw[blue] (0,-6) node{$\ell$};
	\draw (1.5,-6) node{$\ell-1$};
	\draw[dots] (3,-6)--(4.5,-6);
	\draw (6,-6) node{$ $};
	\draw[darkgreen] (7.5,-6) node{$g$};
	\draw (9,-6) node{$ $};
	\draw[dots] (10.5,-6)--(12,-6);
	\draw (13.5,-6) node{$i$};
	\draw[blue] (15,-6) node{$\ell$};
	\draw (16.5,-6) node{$\ell-1$};
	\draw[dots] (18,-6)--(19.5,-6);
	\draw[darkgreen] (20.7,-6) node{$g$};
	\draw[dots] (21.9,-6)--(23.4,-6);
	\draw[darkgreen] (24.6,-6) node{$g$};
	\draw[dots] (25.8,-6)--(27.3,-6);
	\draw (28.5,-6) node{$j$};

	\draw[blue](0,-5.2)--(15,5.2);
	\draw(1.5,-5.2)--(16.5,5.2);
	\draw(6,-5.2)--(21,5.2);
	\draw[darkgreen](7.5,-5.2)--(22.5,5.2);
	\draw(9,-5.2)--(24,5.2);
	\draw(13.5,-5.2)--(28.5,5.2);

	\draw[blue](0,5.2)--(15,-5.2);
	\draw(1.5,5.2)--(16.5,-5.2);
	\draw[darkgreen](5.7,5.2)--(20.7,-5.2);
	\draw[darkgreen](9.6,5.2)--(24.6,-5.2);
	\draw(13.5,5.2)--(28.5,-5.2);

	\filldraw[color=blue!60, fill=black!5, thin] (7.5,0) circle (0.15);
	\filldraw[color=black!60, fill=black!5, thin] (8.25,0.52) circle (0.15);
	\filldraw[color=black!60, fill=black!5, thin] (10.35,1.976) circle (0.15);
	\filldraw[color=black!60, fill=black!5, thin] (12.3,3.328) circle (0.15);
	\filldraw[color=black!60, fill=black!5, thin] (14.25,4.68) circle (0.15);

	\filldraw[color=black!60, fill=black!5, thin] (8.25,-0.52) circle (0.15);
	\filldraw[color=black!60, fill=black!5, thin] (9,0) circle (0.15);
	\filldraw[color=black!60, fill=black!5, thin] (11.1,1.456) circle (0.15);
	\filldraw[color=black!60, fill=black!5, thin] (13.05,2.808) circle (0.15);
	\filldraw[color=black!60, fill=black!5, thin] (15,4.16) circle (0.15);

	\filldraw[color=black!60, fill=black!5, thin] (10.5,-2.08) circle (0.15);
	\filldraw[color=black!60, fill=black!5, thin] (11.25,-1.56) circle (0.15);
	\filldraw[color=black!60, fill=black!5, thin] (13.35,-0.104) circle (0.15);
	\filldraw[color=black!60, fill=black!5, thin] (15.3,1.248) circle (0.15);
	\filldraw[color=black!60, fill=black!5, thin] (17.25,2.6) circle (0.15);

	\filldraw[color=black!60, fill=black!5, thin] (11.25,-2.6) circle (0.15);
	\filldraw[color=black!60, fill=black!5, thin] (12,-2.08) circle (0.15);
	\filldraw[color=darkgreen!60, fill=black!5, thin] (16.05,0.728) circle (0.15);
	\filldraw[color=black!60, fill=black!5, thin] (18,2.08) circle (0.15);

	\smallgreendot (13.725,-0.364);
	
	\end{braid}
	}
$$

$$
-\resizebox{104mm}{23mm}{
\begin{braid}\tikzset{baseline=-.3em}
	\draw[blue] (0,6) node{$\ell$};
	\draw (1.5,6) node{$\ell-1$};
	\draw[dots] (3,6)--(4.5,6);
	\draw[darkgreen] (5.7,6) node{$g$};
	\draw[dots] (6.9,6)--(8.4,6);
	\draw[darkgreen] (9.6,6) node{$g$};
	\draw[dots] (10.8,6)--(12.3,6);
	\draw (13.5,6) node{$j$};
	\draw[blue] (15,6) node{$\ell$};
    \draw (16.5,6) node{$\ell-1$};
	\draw[dots] (18,6)--(19.5,6);
	\draw (21,6) node{$ $};
	\draw[darkgreen] (22.5,6) node{$g$};
	\draw (24,6) node{$ $};
	\draw[dots] (25.5,6)--(27,6);
	\draw (28.5,6) node{$i$};

	\draw[blue] (0,-6) node{$\ell$};
	\draw (1.5,-6) node{$\ell-1$};
	\draw[dots] (3,-6)--(4.5,-6);
	\draw (6,-6) node{$ $};
	\draw[darkgreen] (7.5,-6) node{$g$};
	\draw (9,-6) node{$ $};
	\draw[dots] (10.5,-6)--(12,-6);
	\draw (13.5,-6) node{$i$};
	\draw[blue] (15,-6) node{$\ell$};
	\draw (16.5,-6) node{$\ell-1$};
	\draw[dots] (18,-6)--(19.5,-6);
	\draw[darkgreen] (20.7,-6) node{$g$};
	\draw[dots] (21.9,-6)--(23.4,-6);
	\draw[darkgreen] (24.6,-6) node{$g$};
	\draw[dots] (25.8,-6)--(27.3,-6);
	\draw (28.5,-6) node{$j$};

	\draw[blue](0,-5.2)--(15,5.2);
	\draw(1.5,-5.2)--(16.5,5.2);
	\draw(6,-5.2)--(21,5.2);
	\draw[darkgreen](7.5,-5.2)--(22.5,5.2);
	\draw(9,-5.2)--(24,5.2);
	\draw[thin](13.5,-5.2)--(28.5,5.2);

	\draw[blue](0,5.2)--(15,-5.2);
	\draw(1.5,5.2)--(16.5,-5.2);
	\draw[darkgreen](5.7,5.2)--(20.7,-5.2);
	\draw[darkgreen](9.6,5.2)--(24.6,-5.2);
	\draw(13.5,5.2)--(28.5,-5.2);

	\filldraw[color=blue!60, fill=black!5, thin] (7.5,0) circle (0.15);
	\filldraw[color=black!60, fill=black!5, thin] (8.25,0.52) circle (0.15);
	\filldraw[color=black!60, fill=black!5, thin] (10.35,1.976) circle (0.15);
	\filldraw[color=black!60, fill=black!5, thin] (12.3,3.328) circle (0.15);
	\filldraw[color=black!60, fill=black!5, thin] (14.25,4.68) circle (0.15);

	\filldraw[color=black!60, fill=black!5, thin] (8.25,-0.52) circle (0.15);
	\filldraw[color=black!60, fill=black!5, thin] (9,0) circle (0.15);
	\filldraw[color=black!60, fill=black!5, thin] (11.1,1.456) circle (0.15);
	\filldraw[color=black!60, fill=black!5, thin] (13.05,2.808) circle (0.15);
	\filldraw[color=black!60, fill=black!5, thin] (15,4.16) circle (0.15);

	\filldraw[color=black!60, fill=black!5, thin] (10.5,-2.08) circle (0.15);
	\filldraw[color=black!60, fill=black!5, thin] (11.25,-1.56) circle (0.15);
	\filldraw[color=black!60, fill=black!5, thin] (13.35,-0.104) circle (0.15);
	\filldraw[color=black!60, fill=black!5, thin] (15.3,1.248) circle (0.15);
	\filldraw[color=black!60, fill=black!5, thin] (17.25,2.6) circle (0.15);

	\filldraw[color=black!60, fill=black!5, thin] (11.25,-2.6) circle (0.15);
	\filldraw[color=black!60, fill=black!5, thin] (12,-2.08) circle (0.15);
	\filldraw[color=darkgreen!60, fill=black!5, thin] (16.05,0.728) circle (0.15);
	\filldraw[color=black!60, fill=black!5, thin] (18,2.08) circle (0.15);

	\smallgreendot (14.475,-0.364);
	
	\end{braid}
	}.
$$
Repeated application of (\ref{ECircleBraid}) gives that the first summand above is equal to
$$
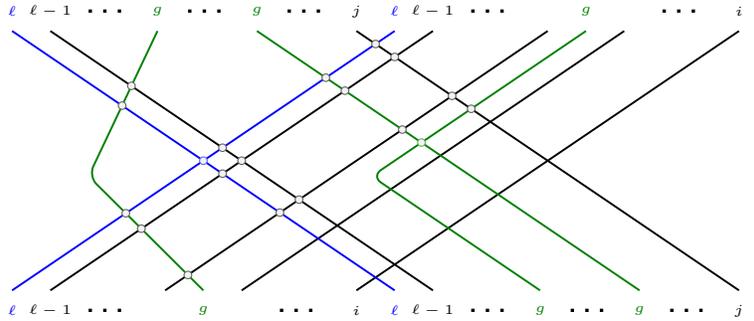
\resizebox{104mm}{23mm}{
\begin{braid}\tikzset{baseline=-.3em}
	\draw[blue] (0,6) node{$\ell$};
	\draw (1.5,6) node{$\ell-1$};
	\draw[dots] (3,6)--(4.5,6);
	\draw[darkgreen] (5.7,6) node{$g$};
	\draw[dots] (6.9,6)--(8.4,6);
	\draw[darkgreen] (9.6,6) node{$g$};
	\draw[dots] (10.8,6)--(12.3,6);
	\draw (13.5,6) node{$j$};
	\draw[blue] (15,6) node{$\ell$};
    \draw (16.5,6) node{$\ell-1$};
	\draw[dots] (18,6)--(19.5,6);
	\draw (21,6) node{$ $};
	\draw[darkgreen] (22.5,6) node{$g$};
	\draw (24,6) node{$ $};
	\draw[dots] (25.5,6)--(27,6);
	\draw (28.5,6) node{$i$};

	\draw[blue] (0,-6) node{$\ell$};
	\draw (1.5,-6) node{$\ell-1$};
	\draw[dots] (3,-6)--(4.5,-6);
	\draw (6,-6) node{$ $};
	\draw[darkgreen] (7.5,-6) node{$g$};
	\draw (9,-6) node{$ $};
	\draw[dots] (10.5,-6)--(12,-6);
	\draw (13.5,-6) node{$i$};
	\draw[blue] (15,-6) node{$\ell$};
	\draw (16.5,-6) node{$\ell-1$};
	\draw[dots] (18,-6)--(19.5,-6);
	\draw[darkgreen] (20.7,-6) node{$g$};
	\draw[dots] (21.9,-6)--(23.4,-6);
	\draw[darkgreen] (24.6,-6) node{$g$};
	\draw[dots] (25.8,-6)--(27.3,-6);
	\draw (28.5,-6) node{$j$};

	\draw[blue](0,-5.2)--(15,5.2);
	\draw(1.5,-5.2)--(16.5,5.2);
	\draw(6,-5.2)--(21,5.2);
	\draw(9,-5.2)--(24,5.2);
	\draw(13.5,-5.2)--(28.5,5.2);

	\draw[blue](0,5.2)--(15,-5.2);
	\draw(1.5,5.2)--(16.5,-5.2);
	\draw[darkgreen](9.6,5.2)--(24.6,-5.2);
	\draw(13.5,5.2)--(28.5,-5.2);
	
	\draw[darkgreen](7.5,-5.2)--(3,-0.624)--(5.7,5.2);
    \draw[darkgreen](20.7,-5.2)--(14.1,-0.624)--(22.5,5.2);

	\filldraw[color=blue!60, fill=black!5, thin] (7.5,0) circle (0.15);
	\filldraw[color=black!60, fill=black!5, thin] (8.25,0.52) circle (0.15);
	\filldraw[color=black!60, fill=black!5, thin] (12.3,3.328) circle (0.15);
	\filldraw[color=black!60, fill=black!5, thin] (14.25,4.68) circle (0.15);

	\filldraw[color=black!60, fill=black!5, thin] (8.25,-0.52) circle (0.15);
	\filldraw[color=black!60, fill=black!5, thin] (9,0) circle (0.15);
	\filldraw[color=black!60, fill=black!5, thin] (13.05,2.808) circle (0.15);
	\filldraw[color=black!60, fill=black!5, thin] (15,4.16) circle (0.15);

	\filldraw[color=black!60, fill=black!5, thin] (10.5,-2.08) circle (0.15);
	\filldraw[color=black!60, fill=black!5, thin] (11.25,-1.56) circle (0.15);
	\filldraw[color=black!60, fill=black!5, thin] (15.3,1.248) circle (0.15);
	\filldraw[color=black!60, fill=black!5, thin] (17.25,2.6) circle (0.15);

	\filldraw[color=darkgreen!60, fill=black!5, thin] (16.05,0.728) circle (0.15);
	\filldraw[color=black!60, fill=black!5, thin] (18,2.08) circle (0.15);

	\filldraw[color=black!60, fill=black!5, thin] (4.459,-2.108) circle (0.15);
	\filldraw[color=black!60, fill=black!5, thin] (5.068,-2.726) circle (0.15);
	\filldraw[color=black!60, fill=black!5, thin] (6.892,-4.582) circle (0.15);

	\filldraw[color=black!60, fill=black!5, thin] (4.314,2.209) circle (0.15);
	\filldraw[color=black!60, fill=black!5, thin] (4.678,3) circle (0.15);

	\end{braid}
	},
$$
which, by Lemma \ref{LCuspExpl}, is equal to zero. Now, by (\ref{EDotPastCircle}), the second summand is equal to
$$
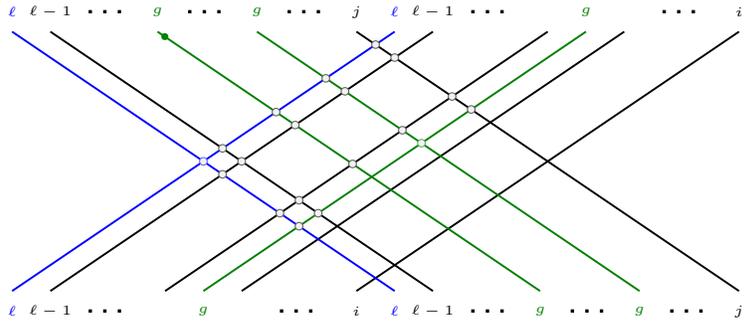
\resizebox{104mm}{23mm}{
\begin{braid}\tikzset{baseline=-.3em}
	\draw[blue] (0,6) node{$\ell$};
	\draw (1.5,6) node{$\ell-1$};
	\draw[dots] (3,6)--(4.5,6);
	\draw[darkgreen] (5.7,6) node{$g$};
	\draw[dots] (6.9,6)--(8.4,6);
	\draw[darkgreen] (9.6,6) node{$g$};
	\draw[dots] (10.8,6)--(12.3,6);
	\draw (13.5,6) node{$j$};
	\draw[blue] (15,6) node{$\ell$};
    \draw (16.5,6) node{$\ell-1$};
	\draw[dots] (18,6)--(19.5,6);
	\draw (21,6) node{$ $};
	\draw[darkgreen] (22.5,6) node{$g$};
	\draw (24,6) node{$ $};
	\draw[dots] (25.5,6)--(27,6);
	\draw (28.5,6) node{$i$};

	\draw[blue] (0,-6) node{$\ell$};
	\draw (1.5,-6) node{$\ell-1$};
	\draw[dots] (3,-6)--(4.5,-6);
	\draw (6,-6) node{$ $};
	\draw[darkgreen] (7.5,-6) node{$g$};
	\draw (9,-6) node{$ $};
	\draw[dots] (10.5,-6)--(12,-6);
	\draw (13.5,-6) node{$i$};
	\draw[blue] (15,-6) node{$\ell$};
	\draw (16.5,-6) node{$\ell-1$};
	\draw[dots] (18,-6)--(19.5,-6);
	\draw[darkgreen] (20.7,-6) node{$g$};
	\draw[dots] (21.9,-6)--(23.4,-6);
	\draw[darkgreen] (24.6,-6) node{$g$};
	\draw[dots] (25.8,-6)--(27.3,-6);
	\draw (28.5,-6) node{$j$};

	\draw[blue](0,-5.2)--(15,5.2);
	\draw(1.5,-5.2)--(16.5,5.2);
	\draw(6,-5.2)--(21,5.2);
	\draw[darkgreen](7.5,-5.2)--(22.5,5.2);
	\draw(9,-5.2)--(24,5.2);
	\draw(13.5,-5.2)--(28.5,5.2);

	\draw[blue](0,5.2)--(15,-5.2);
	\draw(1.5,5.2)--(16.5,-5.2);
	\draw[darkgreen](5.7,5.2)--(20.7,-5.2);
	\draw[darkgreen](9.6,5.2)--(24.6,-5.2);
	\draw(13.5,5.2)--(28.5,-5.2);

	\filldraw[color=blue!60, fill=black!5, thin] (7.5,0) circle (0.15);
	\filldraw[color=black!60, fill=black!5, thin] (8.25,0.52) circle (0.15);
	\filldraw[color=black!60, fill=black!5, thin] (10.35,1.976) circle (0.15);
	\filldraw[color=black!60, fill=black!5, thin] (12.3,3.328) circle (0.15);
	\filldraw[color=black!60, fill=black!5, thin] (14.25,4.68) circle (0.15);

	\filldraw[color=black!60, fill=black!5, thin] (8.25,-0.52) circle (0.15);
	\filldraw[color=black!60, fill=black!5, thin] (9,0) circle (0.15);
	\filldraw[color=black!60, fill=black!5, thin] (11.1,1.456) circle (0.15);
	\filldraw[color=black!60, fill=black!5, thin] (13.05,2.808) circle (0.15);
	\filldraw[color=black!60, fill=black!5, thin] (15,4.16) circle (0.15);

	\filldraw[color=black!60, fill=black!5, thin] (10.5,-2.08) circle (0.15);
	\filldraw[color=black!60, fill=black!5, thin] (11.25,-1.56) circle (0.15);
	\filldraw[color=black!60, fill=black!5, thin] (13.35,-0.104) circle (0.15);
	\filldraw[color=black!60, fill=black!5, thin] (15.3,1.248) circle (0.15);
	\filldraw[color=black!60, fill=black!5, thin] (17.25,2.6) circle (0.15);

	\filldraw[color=black!60, fill=black!5, thin] (11.25,-2.6) circle (0.15);
	\filldraw[color=black!60, fill=black!5, thin] (12,-2.08) circle (0.15);
	\filldraw[color=darkgreen!60, fill=black!5, thin] (16.05,0.728) circle (0.15);
	\filldraw[color=black!60, fill=black!5, thin] (18,2.08) circle (0.15);

	\smallgreendot (5.988,5);
	
	\end{braid}
	}.
$$
Once again by (\ref{EOpening}), this is equal to
$$
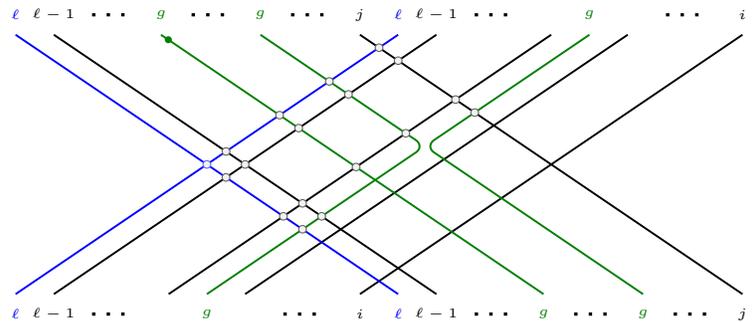
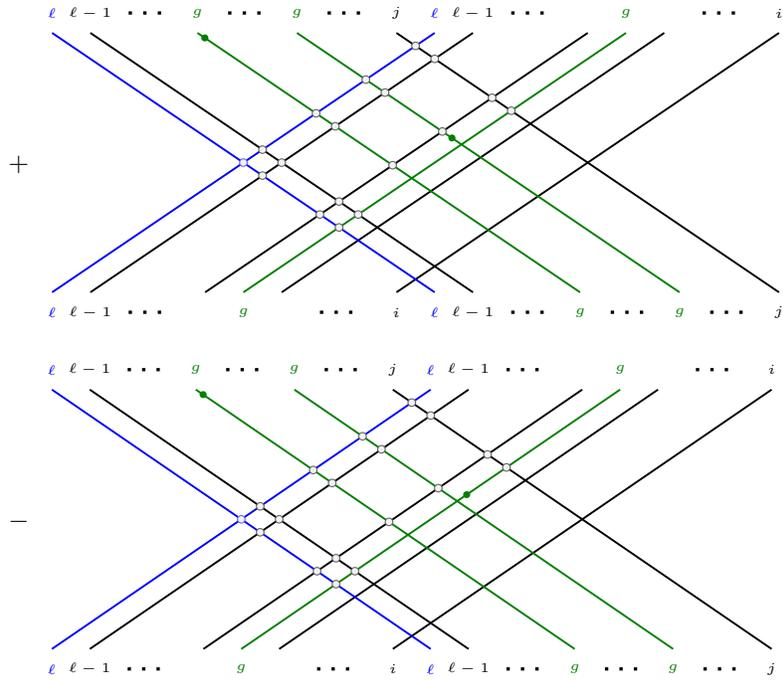
\resizebox{104mm}{23mm}{
\begin{braid}\tikzset{baseline=-.3em}
	\draw[blue] (0,6) node{$\ell$};
	\draw (1.5,6) node{$\ell-1$};
	\draw[dots] (3,6)--(4.5,6);
	\draw[darkgreen] (5.7,6) node{$g$};
	\draw[dots] (6.9,6)--(8.4,6);
	\draw[darkgreen] (9.6,6) node{$g$};
	\draw[dots] (10.8,6)--(12.3,6);
	\draw (13.5,6) node{$j$};
	\draw[blue] (15,6) node{$\ell$};
    \draw (16.5,6) node{$\ell-1$};
	\draw[dots] (18,6)--(19.5,6);
	\draw (21,6) node{$ $};
	\draw[darkgreen] (22.5,6) node{$g$};
	\draw (24,6) node{$ $};
	\draw[dots] (25.5,6)--(27,6);
	\draw (28.5,6) node{$i$};

	\draw[blue] (0,-6) node{$\ell$};
	\draw (1.5,-6) node{$\ell-1$};
	\draw[dots] (3,-6)--(4.5,-6);
	\draw (6,-6) node{$ $};
	\draw[darkgreen] (7.5,-6) node{$g$};
	\draw (9,-6) node{$ $};
	\draw[dots] (10.5,-6)--(12,-6);
	\draw (13.5,-6) node{$i$};
	\draw[blue] (15,-6) node{$\ell$};
	\draw (16.5,-6) node{$\ell-1$};
	\draw[dots] (18,-6)--(19.5,-6);
	\draw[darkgreen] (20.7,-6) node{$g$};
	\draw[dots] (21.9,-6)--(23.4,-6);
	\draw[darkgreen] (24.6,-6) node{$g$};
	\draw[dots] (25.8,-6)--(27.3,-6);
	\draw (28.5,-6) node{$j$};

	\draw[blue](0,-5.2)--(15,5.2);
	\draw(1.5,-5.2)--(16.5,5.2);
	\draw(6,-5.2)--(21,5.2);
	\draw(9,-5.2)--(24,5.2);
	\draw(13.5,-5.2)--(28.5,5.2);

	\draw[blue](0,5.2)--(15,-5.2);
	\draw(1.5,5.2)--(16.5,-5.2);
	\draw[darkgreen](5.7,5.2)--(20.7,-5.2);
	\draw(13.5,5.2)--(28.5,-5.2);

	\draw[darkgreen](7.5,-5.2)--(16.05,0.728)--(9.6,5.2);
	\draw[darkgreen](24.6,-5.2)--(16.05,0.728)--(22.5,5.2);

	\filldraw[color=blue!60, fill=black!5, thin] (7.5,0) circle (0.15);
	\filldraw[color=black!60, fill=black!5, thin] (8.25,0.52) circle (0.15);
	\filldraw[color=black!60, fill=black!5, thin] (10.35,1.976) circle (0.15);
	\filldraw[color=black!60, fill=black!5, thin] (12.3,3.328) circle (0.15);
	\filldraw[color=black!60, fill=black!5, thin] (14.25,4.68) circle (0.15);

	\filldraw[color=black!60, fill=black!5, thin] (8.25,-0.52) circle (0.15);
	\filldraw[color=black!60, fill=black!5, thin] (9,0) circle (0.15);
	\filldraw[color=black!60, fill=black!5, thin] (11.1,1.456) circle (0.15);
	\filldraw[color=black!60, fill=black!5, thin] (13.05,2.808) circle (0.15);
	\filldraw[color=black!60, fill=black!5, thin] (15,4.16) circle (0.15);

	\filldraw[color=black!60, fill=black!5, thin] (10.5,-2.08) circle (0.15);
	\filldraw[color=black!60, fill=black!5, thin] (11.25,-1.56) circle (0.15);
	\filldraw[color=black!60, fill=black!5, thin] (13.35,-0.104) circle (0.15);
	\filldraw[color=black!60, fill=black!5, thin] (15.3,1.248) circle (0.15);
	\filldraw[color=black!60, fill=black!5, thin] (17.25,2.6) circle (0.15);

	\filldraw[color=black!60, fill=black!5, thin] (11.25,-2.6) circle (0.15);
	\filldraw[color=black!60, fill=black!5, thin] (12,-2.08) circle (0.15);
	\filldraw[color=black!60, fill=black!5, thin] (18,2.08) circle (0.15);

	\smallgreendot (5.988,5);
	
	\end{braid}
	}
$$

$$
+\resizebox{104mm}{23mm}{
\begin{braid}\tikzset{baseline=-.3em}
	\draw[blue] (0,6) node{$\ell$};
	\draw (1.5,6) node{$\ell-1$};
	\draw[dots] (3,6)--(4.5,6);
	\draw[darkgreen] (5.7,6) node{$g$};
	\draw[dots] (6.9,6)--(8.4,6);
	\draw[darkgreen] (9.6,6) node{$g$};
	\draw[dots] (10.8,6)--(12.3,6);
	\draw (13.5,6) node{$j$};
	\draw[blue] (15,6) node{$\ell$};
    \draw (16.5,6) node{$\ell-1$};
	\draw[dots] (18,6)--(19.5,6);
	\draw (21,6) node{$ $};
	\draw[darkgreen] (22.5,6) node{$g$};
	\draw (24,6) node{$ $};
	\draw[dots] (25.5,6)--(27,6);
	\draw (28.5,6) node{$i$};

	\draw[blue] (0,-6) node{$\ell$};
	\draw (1.5,-6) node{$\ell-1$};
	\draw[dots] (3,-6)--(4.5,-6);
	\draw (6,-6) node{$ $};
	\draw[darkgreen] (7.5,-6) node{$g$};
	\draw (9,-6) node{$ $};
	\draw[dots] (10.5,-6)--(12,-6);
	\draw (13.5,-6) node{$i$};
	\draw[blue] (15,-6) node{$\ell$};
	\draw (16.5,-6) node{$\ell-1$};
	\draw[dots] (18,-6)--(19.5,-6);
	\draw[darkgreen] (20.7,-6) node{$g$};
	\draw[dots] (21.9,-6)--(23.4,-6);
	\draw[darkgreen] (24.6,-6) node{$g$};
	\draw[dots] (25.8,-6)--(27.3,-6);
	\draw (28.5,-6) node{$j$};

	\draw[blue](0,-5.2)--(15,5.2);
	\draw(1.5,-5.2)--(16.5,5.2);
	\draw(6,-5.2)--(21,5.2);
	\draw[darkgreen](7.5,-5.2)--(22.5,5.2);
	\draw(9,-5.2)--(24,5.2);
	\draw(13.5,-5.2)--(28.5,5.2);

	\draw[blue](0,5.2)--(15,-5.2);
	\draw(1.5,5.2)--(16.5,-5.2);
	\draw[darkgreen](5.7,5.2)--(20.7,-5.2);
	\draw[darkgreen](9.6,5.2)--(24.6,-5.2);
	\draw(13.5,5.2)--(28.5,-5.2);

	\filldraw[color=blue!60, fill=black!5, thin] (7.5,0) circle (0.15);
	\filldraw[color=black!60, fill=black!5, thin] (8.25,0.52) circle (0.15);
	\filldraw[color=black!60, fill=black!5, thin] (10.35,1.976) circle (0.15);
	\filldraw[color=black!60, fill=black!5, thin] (12.3,3.328) circle (0.15);
	\filldraw[color=black!60, fill=black!5, thin] (14.25,4.68) circle (0.15);

	\filldraw[color=black!60, fill=black!5, thin] (8.25,-0.52) circle (0.15);
	\filldraw[color=black!60, fill=black!5, thin] (9,0) circle (0.15);
	\filldraw[color=black!60, fill=black!5, thin] (11.1,1.456) circle (0.15);
	\filldraw[color=black!60, fill=black!5, thin] (13.05,2.808) circle (0.15);
	\filldraw[color=black!60, fill=black!5, thin] (15,4.16) circle (0.15);

	\filldraw[color=black!60, fill=black!5, thin] (10.5,-2.08) circle (0.15);
	\filldraw[color=black!60, fill=black!5, thin] (11.25,-1.56) circle (0.15);
	\filldraw[color=black!60, fill=black!5, thin] (13.35,-0.104) circle (0.15);
	\filldraw[color=black!60, fill=black!5, thin] (15.3,1.248) circle (0.15);
	\filldraw[color=black!60, fill=black!5, thin] (17.25,2.6) circle (0.15);

	\filldraw[color=black!60, fill=black!5, thin] (11.25,-2.6) circle (0.15);
	\filldraw[color=black!60, fill=black!5, thin] (12,-2.08) circle (0.15);
	\filldraw[color=black!60, fill=black!5, thin] (18,2.08) circle (0.15);

	\smallgreendot (5.988,5);
	\smallgreendot (15.675,0.988);
	\end{braid}
	}
$$

$$
-\resizebox{104mm}{23mm}{
\begin{braid}\tikzset{baseline=-.3em}
	\draw[blue] (0,6) node{$\ell$};
	\draw (1.5,6) node{$\ell-1$};
	\draw[dots] (3,6)--(4.5,6);
	\draw[darkgreen] (5.7,6) node{$g$};
	\draw[dots] (6.9,6)--(8.4,6);
	\draw[darkgreen] (9.6,6) node{$g$};
	\draw[dots] (10.8,6)--(12.3,6);
	\draw (13.5,6) node{$j$};
	\draw[blue] (15,6) node{$\ell$};
    \draw (16.5,6) node{$\ell-1$};
	\draw[dots] (18,6)--(19.5,6);
	\draw (21,6) node{$ $};
	\draw[darkgreen] (22.5,6) node{$g$};
	\draw (24,6) node{$ $};
	\draw[dots] (25.5,6)--(27,6);
	\draw (28.5,6) node{$i$};

	\draw[blue] (0,-6) node{$\ell$};
	\draw (1.5,-6) node{$\ell-1$};
	\draw[dots] (3,-6)--(4.5,-6);
	\draw (6,-6) node{$ $};
	\draw[darkgreen] (7.5,-6) node{$g$};
	\draw (9,-6) node{$ $};
	\draw[dots] (10.5,-6)--(12,-6);
	\draw (13.5,-6) node{$i$};
	\draw[blue] (15,-6) node{$\ell$};
	\draw (16.5,-6) node{$\ell-1$};
	\draw[dots] (18,-6)--(19.5,-6);
	\draw[darkgreen] (20.7,-6) node{$g$};
	\draw[dots] (21.9,-6)--(23.4,-6);
	\draw[darkgreen] (24.6,-6) node{$g$};
	\draw[dots] (25.8,-6)--(27.3,-6);
	\draw (28.5,-6) node{$j$};

	\draw[blue](0,-5.2)--(15,5.2);
	\draw(1.5,-5.2)--(16.5,5.2);
	\draw(6,-5.2)--(21,5.2);
	\draw[darkgreen](7.5,-5.2)--(22.5,5.2);
	\draw(9,-5.2)--(24,5.2);
	\draw(13.5,-5.2)--(28.5,5.2);

	\draw[blue](0,5.2)--(15,-5.2);
	\draw(1.5,5.2)--(16.5,-5.2);
	\draw[darkgreen](5.7,5.2)--(20.7,-5.2);
	\draw[darkgreen](9.6,5.2)--(24.6,-5.2);
	\draw(13.5,5.2)--(28.5,-5.2);

	\filldraw[color=blue!60, fill=black!5, thin] (7.5,0) circle (0.15);
	\filldraw[color=black!60, fill=black!5, thin] (8.25,0.52) circle (0.15);
	\filldraw[color=black!60, fill=black!5, thin] (10.35,1.976) circle (0.15);
	\filldraw[color=black!60, fill=black!5, thin] (12.3,3.328) circle (0.15);
	\filldraw[color=black!60, fill=black!5, thin] (14.25,4.68) circle (0.15);

	\filldraw[color=black!60, fill=black!5, thin] (8.25,-0.52) circle (0.15);
	\filldraw[color=black!60, fill=black!5, thin] (9,0) circle (0.15);
	\filldraw[color=black!60, fill=black!5, thin] (11.1,1.456) circle (0.15);
	\filldraw[color=black!60, fill=black!5, thin] (13.05,2.808) circle (0.15);
	\filldraw[color=black!60, fill=black!5, thin] (15,4.16) circle (0.15);

	\filldraw[color=black!60, fill=black!5, thin] (10.5,-2.08) circle (0.15);
	\filldraw[color=black!60, fill=black!5, thin] (11.25,-1.56) circle (0.15);
	\filldraw[color=black!60, fill=black!5, thin] (13.35,-0.104) circle (0.15);
	\filldraw[color=black!60, fill=black!5, thin] (15.3,1.248) circle (0.15);
	\filldraw[color=black!60, fill=black!5, thin] (17.25,2.6) circle (0.15);

	\filldraw[color=black!60, fill=black!5, thin] (11.25,-2.6) circle (0.15);
	\filldraw[color=black!60, fill=black!5, thin] (12,-2.08) circle (0.15);
	\filldraw[color=black!60, fill=black!5, thin] (18,2.08) circle (0.15);

	\smallgreendot (5.988,5);
	\smallgreendot (16.425,0.988);
	\end{braid},
	}
$$
This time we require repeated applications of both (\ref{EPartialCircleBraid}) and (\ref{ECircleBraid}) then Lemma \ref{LCuspExpl} to give that this first summand is zero. Applying (\ref{EDotPastCircle}) to the last two summands we get
$$
\resizebox{104mm}{23mm}{
\begin{braid}\tikzset{baseline=-.3em}
	\draw[blue] (0,6) node{$\ell$};
	\draw (1.5,6) node{$\ell-1$};
	\draw[dots] (3,6)--(4.5,6);
	\draw[darkgreen] (5.7,6) node{$g$};
	\draw[dots] (6.9,6)--(8.4,6);
	\draw[darkgreen] (9.6,6) node{$g$};
	\draw[dots] (10.8,6)--(12.3,6);
	\draw (13.5,6) node{$j$};
	\draw[blue] (15,6) node{$\ell$};
    \draw (16.5,6) node{$\ell-1$};
	\draw[dots] (18,6)--(19.5,6);
	\draw (21,6) node{$ $};
	\draw[darkgreen] (22.5,6) node{$g$};
	\draw (24,6) node{$ $};
	\draw[dots] (25.5,6)--(27,6);
	\draw (28.5,6) node{$i$};

	\draw[blue] (0,-6) node{$\ell$};
	\draw (1.5,-6) node{$\ell-1$};
	\draw[dots] (3,-6)--(4.5,-6);
	\draw (6,-6) node{$ $};
	\draw[darkgreen] (7.5,-6) node{$g$};
	\draw (9,-6) node{$ $};
	\draw[dots] (10.5,-6)--(12,-6);
	\draw (13.5,-6) node{$i$};
	\draw (15,-6) node{$\ell$};
	\draw (16.5,-6) node{$\ell-1$};
	\draw[dots] (18,-6)--(19.5,-6);
	\draw[darkgreen] (20.7,-6) node{$g$};
	\draw[dots] (21.9,-6)--(23.4,-6);
	\draw[darkgreen] (24.6,-6) node{$g$};
	\draw[dots] (25.8,-6)--(27.3,-6);
	\draw (28.5,-6) node{$j$};

	\draw[blue](0,-5.2)--(15,5.2);
	\draw(1.5,-5.2)--(16.5,5.2);
	\draw(6,-5.2)--(21,5.2);
	\draw[darkgreen](7.5,-5.2)--(22.5,5.2);
	\draw(9,-5.2)--(24,5.2);
	\draw(13.5,-5.2)--(28.5,5.2);

	\draw[blue](0,5.2)--(15,-5.2);
	\draw(1.5,5.2)--(16.5,-5.2);
	\draw[darkgreen](5.7,5.2)--(20.7,-5.2);
	\draw[darkgreen](9.6,5.2)--(24.6,-5.2);
	\draw(13.5,5.2)--(28.5,-5.2);

	\filldraw[color=blue!60, fill=black!5, thin] (7.5,0) circle (0.15);
	\filldraw[color=black!60, fill=black!5, thin] (8.25,0.52) circle (0.15);
	\filldraw[color=black!60, fill=black!5, thin] (10.35,1.976) circle (0.15);
	\filldraw[color=black!60, fill=black!5, thin] (12.3,3.328) circle (0.15);
	\filldraw[color=black!60, fill=black!5, thin] (14.25,4.68) circle (0.15);

	\filldraw[color=black!60, fill=black!5, thin] (8.25,-0.52) circle (0.15);
	\filldraw[color=black!60, fill=black!5, thin] (9,0) circle (0.15);
	\filldraw[color=black!60, fill=black!5, thin] (11.1,1.456) circle (0.15);
	\filldraw[color=black!60, fill=black!5, thin] (13.05,2.808) circle (0.15);
	\filldraw[color=black!60, fill=black!5, thin] (15,4.16) circle (0.15);

	\filldraw[color=black!60, fill=black!5, thin] (10.5,-2.08) circle (0.15);
	\filldraw[color=black!60, fill=black!5, thin] (11.25,-1.56) circle (0.15);
	\filldraw[color=black!60, fill=black!5, thin] (13.35,-0.104) circle (0.15);
	\filldraw[color=black!60, fill=black!5, thin] (15.3,1.248) circle (0.15);
	\filldraw[color=black!60, fill=black!5, thin] (17.25,2.6) circle (0.15);

	\filldraw[color=black!60, fill=black!5, thin] (11.25,-2.6) circle (0.15);
	\filldraw[color=black!60, fill=black!5, thin] (12,-2.08) circle (0.15);
	\filldraw[color=black!60, fill=black!5, thin] (18,2.08) circle (0.15);

	\smallgreendot (5.988,5);
	\smallgreendot (9.888,5);
	\end{braid}
	}
$$

$$
-\resizebox{104mm}{23mm}{
\begin{braid}\tikzset{baseline=-.3em}
	\draw[blue] (0,6) node{$\ell$};
	\draw (1.5,6) node{$\ell-1$};
	\draw[dots] (3,6)--(4.5,6);
	\draw[darkgreen] (5.7,6) node{$g$};
	\draw[dots] (6.9,6)--(8.4,6);
	\draw[darkgreen] (9.6,6) node{$g$};
	\draw[dots] (10.8,6)--(12.3,6);
	\draw (13.5,6) node{$j$};
	\draw[blue] (15,6) node{$\ell$};
    \draw (16.5,6) node{$\ell-1$};
	\draw[dots] (18,6)--(19.5,6);
	\draw (21,6) node{$ $};
	\draw[darkgreen] (22.5,6) node{$g$};
	\draw (24,6) node{$ $};
	\draw[dots] (25.5,6)--(27,6);
	\draw (28.5,6) node{$i$};

	\draw[blue] (0,-6) node{$\ell$};
	\draw (1.5,-6) node{$\ell-1$};
	\draw[dots] (3,-6)--(4.5,-6);
	\draw (6,-6) node{$ $};
	\draw[darkgreen] (7.5,-6) node{$g$};
	\draw (9,-6) node{$ $};
	\draw[dots] (10.5,-6)--(12,-6);
	\draw (13.5,-6) node{$i$};
	\draw[blue] (15,-6) node{$\ell$};
	\draw (16.5,-6) node{$\ell-1$};
	\draw[dots] (18,-6)--(19.5,-6);
	\draw[darkgreen] (20.7,-6) node{$g$};
	\draw[dots] (21.9,-6)--(23.4,-6);
	\draw[darkgreen] (24.6,-6) node{$g$};
	\draw[dots] (25.8,-6)--(27.3,-6);
	\draw (28.5,-6) node{$j$};

	\draw[blue](0,-5.2)--(15,5.2);
	\draw(1.5,-5.2)--(16.5,5.2);
	\draw(6,-5.2)--(21,5.2);
	\draw[darkgreen](7.5,-5.2)--(22.5,5.2);
	\draw(9,-5.2)--(24,5.2);
	\draw(13.5,-5.2)--(28.5,5.2);

	\draw[blue](0,5.2)--(15,-5.2);
	\draw(1.5,5.2)--(16.5,-5.2);
	\draw[darkgreen](5.7,5.2)--(20.7,-5.2);
	\draw[darkgreen](9.6,5.2)--(24.6,-5.2);
	\draw(13.5,5.2)--(28.5,-5.2);

	\filldraw[color=blue!60, fill=black!5, thin] (7.5,0) circle (0.15);
	\filldraw[color=black!60, fill=black!5, thin] (8.25,0.52) circle (0.15);
	\filldraw[color=black!60, fill=black!5, thin] (10.35,1.976) circle (0.15);
	\filldraw[color=black!60, fill=black!5, thin] (12.3,3.328) circle (0.15);
	\filldraw[color=black!60, fill=black!5, thin] (14.25,4.68) circle (0.15);

	\filldraw[color=black!60, fill=black!5, thin] (8.25,-0.52) circle (0.15);
	\filldraw[color=black!60, fill=black!5, thin] (9,0) circle (0.15);
	\filldraw[color=black!60, fill=black!5, thin] (11.1,1.456) circle (0.15);
	\filldraw[color=black!60, fill=black!5, thin] (13.05,2.808) circle (0.15);
	\filldraw[color=black!60, fill=black!5, thin] (15,4.16) circle (0.15);

	\filldraw[color=black!60, fill=black!5, thin] (10.5,-2.08) circle (0.15);
	\filldraw[color=black!60, fill=black!5, thin] (11.25,-1.56) circle (0.15);
	\filldraw[color=black!60, fill=black!5, thin] (13.35,-0.104) circle (0.15);
	\filldraw[color=black!60, fill=black!5, thin] (15.3,1.248) circle (0.15);
	\filldraw[color=black!60, fill=black!5, thin] (17.25,2.6) circle (0.15);

	\filldraw[color=black!60, fill=black!5, thin] (11.25,-2.6) circle (0.15);
	\filldraw[color=black!60, fill=black!5, thin] (12,-2.08) circle (0.15);
	\filldraw[color=black!60, fill=black!5, thin] (18,2.08) circle (0.15);

	\smallgreendot (5.988,5);
	\smallgreendot (22.212,5);
	\end{braid}
	}.
$$
Using (\ref{EPhiDifCol}) this becomes
$$
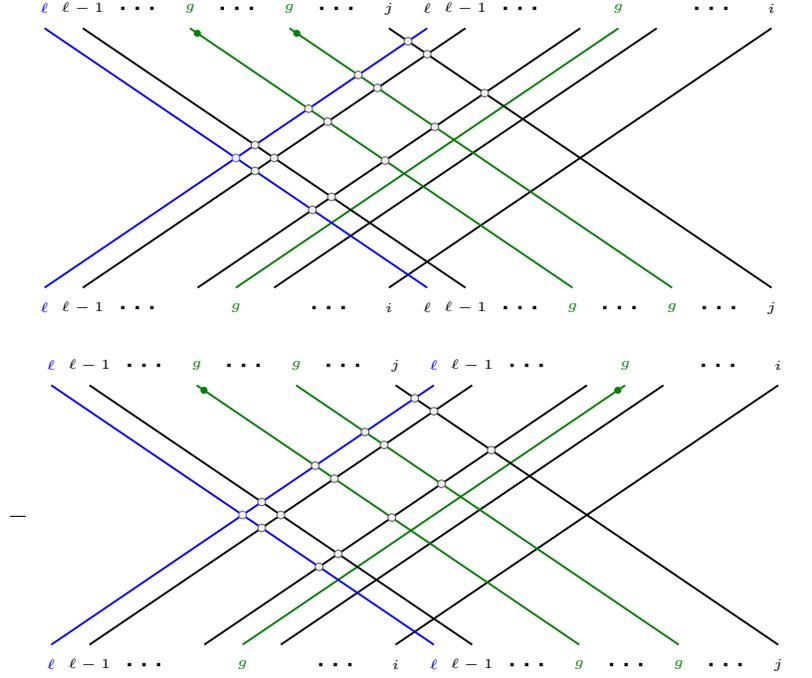
\resizebox{104mm}{23mm}{
\begin{braid}\tikzset{baseline=-.3em}
	\draw[blue] (0,6) node{$\ell$};
	\draw (1.5,6) node{$\ell-1$};
	\draw[dots] (3,6)--(4.5,6);
	\draw[darkgreen] (5.7,6) node{$g$};
	\draw[dots] (6.9,6)--(8.4,6);
	\draw[darkgreen] (9.6,6) node{$g$};
	\draw[dots] (10.8,6)--(12.3,6);
	\draw (13.5,6) node{$j$};
	\draw (15,6) node{$\ell$};
    \draw (16.5,6) node{$\ell-1$};
	\draw[dots] (18,6)--(19.5,6);
	\draw (21,6) node{$ $};
	\draw[darkgreen] (22.5,6) node{$g$};
	\draw (24,6) node{$ $};
	\draw[dots] (25.5,6)--(27,6);
	\draw (28.5,6) node{$i$};

	\draw[blue] (0,-6) node{$\ell$};
	\draw (1.5,-6) node{$\ell-1$};
	\draw[dots] (3,-6)--(4.5,-6);
	\draw (6,-6) node{$ $};
	\draw[darkgreen] (7.5,-6) node{$g$};
	\draw (9,-6) node{$ $};
	\draw[dots] (10.5,-6)--(12,-6);
	\draw (13.5,-6) node{$i$};
	\draw (15,-6) node{$\ell$};
	\draw (16.5,-6) node{$\ell-1$};
	\draw[dots] (18,-6)--(19.5,-6);
	\draw[darkgreen] (20.7,-6) node{$g$};
	\draw[dots] (21.9,-6)--(23.4,-6);
	\draw[darkgreen] (24.6,-6) node{$g$};
	\draw[dots] (25.8,-6)--(27.3,-6);
	\draw (28.5,-6) node{$j$};

	\draw[blue](0,-5.2)--(15,5.2);
	\draw(1.5,-5.2)--(16.5,5.2);
	\draw(6,-5.2)--(21,5.2);
	\draw[darkgreen](7.5,-5.2)--(22.5,5.2);
	\draw(9,-5.2)--(24,5.2);
	\draw(13.5,-5.2)--(28.5,5.2);

	\draw[blue](0,5.2)--(15,-5.2);
	\draw(1.5,5.2)--(16.5,-5.2);
	\draw[darkgreen](5.7,5.2)--(20.7,-5.2);
	\draw[darkgreen](9.6,5.2)--(24.6,-5.2);
	\draw(13.5,5.2)--(28.5,-5.2);

	\filldraw[color=blue!60, fill=black!5, thin] (7.5,0) circle (0.15);
	\filldraw[color=black!60, fill=black!5, thin] (8.25,0.52) circle (0.15);
	\filldraw[color=black!60, fill=black!5, thin] (10.35,1.976) circle (0.15);
	\filldraw[color=black!60, fill=black!5, thin] (12.3,3.328) circle (0.15);
	\filldraw[color=black!60, fill=black!5, thin] (14.25,4.68) circle (0.15);

	\filldraw[color=black!60, fill=black!5, thin] (8.25,-0.52) circle (0.15);
	\filldraw[color=black!60, fill=black!5, thin] (9,0) circle (0.15);
	\filldraw[color=black!60, fill=black!5, thin] (11.1,1.456) circle (0.15);
	\filldraw[color=black!60, fill=black!5, thin] (13.05,2.808) circle (0.15);
	\filldraw[color=black!60, fill=black!5, thin] (15,4.16) circle (0.15);

	\filldraw[color=black!60, fill=black!5, thin] (10.5,-2.08) circle (0.15);
	\filldraw[color=black!60, fill=black!5, thin] (11.25,-1.56) circle (0.15);
	\filldraw[color=black!60, fill=black!5, thin] (13.35,-0.104) circle (0.15);
	\filldraw[color=black!60, fill=black!5, thin] (15.3,1.248) circle (0.15);
	\filldraw[color=black!60, fill=black!5, thin] (17.25,2.6) circle (0.15);

	\smallgreendot (5.988,5);
	\smallgreendot (9.888,5);
	\end{braid}
	}
$$

$$
-\resizebox{104mm}{23mm}{
\begin{braid}\tikzset{baseline=-.3em}
	\draw[blue] (0,6) node{$\ell$};
	\draw (1.5,6) node{$\ell-1$};
	\draw[dots] (3,6)--(4.5,6);
	\draw[darkgreen] (5.7,6) node{$g$};
	\draw[dots] (6.9,6)--(8.4,6);
	\draw[darkgreen] (9.6,6) node{$g$};
	\draw[dots] (10.8,6)--(12.3,6);
	\draw (13.5,6) node{$j$};
	\draw[blue] (15,6) node{$\ell$};
    \draw (16.5,6) node{$\ell-1$};
	\draw[dots] (18,6)--(19.5,6);
	\draw (21,6) node{$ $};
	\draw[darkgreen] (22.5,6) node{$g$};
	\draw (24,6) node{$ $};
	\draw[dots] (25.5,6)--(27,6);
	\draw (28.5,6) node{$i$};

	\draw[blue] (0,-6) node{$\ell$};
	\draw (1.5,-6) node{$\ell-1$};
	\draw[dots] (3,-6)--(4.5,-6);
	\draw (6,-6) node{$ $};
	\draw[darkgreen] (7.5,-6) node{$g$};
	\draw (9,-6) node{$ $};
	\draw[dots] (10.5,-6)--(12,-6);
	\draw (13.5,-6) node{$i$};
	\draw[blue] (15,-6) node{$\ell$};
	\draw (16.5,-6) node{$\ell-1$};
	\draw[dots] (18,-6)--(19.5,-6);
	\draw[darkgreen] (20.7,-6) node{$g$};
	\draw[dots] (21.9,-6)--(23.4,-6);
	\draw[darkgreen] (24.6,-6) node{$g$};
	\draw[dots] (25.8,-6)--(27.3,-6);
	\draw (28.5,-6) node{$j$};

	\draw[blue](0,-5.2)--(15,5.2);
	\draw(1.5,-5.2)--(16.5,5.2);
	\draw(6,-5.2)--(21,5.2);
	\draw[darkgreen](7.5,-5.2)--(22.5,5.2);
	\draw(9,-5.2)--(24,5.2);
	\draw(13.5,-5.2)--(28.5,5.2);

	\draw[blue](0,5.2)--(15,-5.2);
	\draw(1.5,5.2)--(16.5,-5.2);
	\draw[darkgreen](5.7,5.2)--(20.7,-5.2);
	\draw[darkgreen](9.6,5.2)--(24.6,-5.2);
	\draw(13.5,5.2)--(28.5,-5.2);

	\filldraw[color=blue!60, fill=black!5, thin] (7.5,0) circle (0.15);
	\filldraw[color=black!60, fill=black!5, thin] (8.25,0.52) circle (0.15);
	\filldraw[color=black!60, fill=black!5, thin] (10.35,1.976) circle (0.15);
	\filldraw[color=black!60, fill=black!5, thin] (12.3,3.328) circle (0.15);
	\filldraw[color=black!60, fill=black!5, thin] (14.25,4.68) circle (0.15);

	\filldraw[color=black!60, fill=black!5, thin] (8.25,-0.52) circle (0.15);
	\filldraw[color=black!60, fill=black!5, thin] (9,0) circle (0.15);
	\filldraw[color=black!60, fill=black!5, thin] (11.1,1.456) circle (0.15);
	\filldraw[color=black!60, fill=black!5, thin] (13.05,2.808) circle (0.15);
	\filldraw[color=black!60, fill=black!5, thin] (15,4.16) circle (0.15);

	\filldraw[color=black!60, fill=black!5, thin] (10.5,-2.08) circle (0.15);
	\filldraw[color=black!60, fill=black!5, thin] (11.25,-1.56) circle (0.15);
	\filldraw[color=black!60, fill=black!5, thin] (13.35,-0.104) circle (0.15);
	\filldraw[color=black!60, fill=black!5, thin] (15.3,1.248) circle (0.15);
	\filldraw[color=black!60, fill=black!5, thin] (17.25,2.6) circle (0.15);

	\smallgreendot (5.988,5);
	\smallgreendot (22.212,5);
	\end{braid}
	},
$$
which, in non-diagrammatic form, is $y_{r^j_g}(y_{s^j_g}-y_{t+p})P^{(t-1)}{\hat \gamma^{i,j}}$. Similarly the second summand in our original sum is equal to $y_{t+p}(y_{s^j_g}-y_{t+p})P^{(t-1)}{\hat \gamma^{i,j}}$. Therefore,
$$
P^{(t)}{\hat \gamma^{i,j}}=(y_{r^j_g}-y_{t+p})(y_{s^j_g}-y_{t+p})P^{(t-1)}{\hat \gamma^{i,j}},
$$
as desired.

We now assume $g=0$. This case is similar but we outline the main differences.
$$
P^{(t)}{\hat \ga}^{i,j}=
\resizebox{104mm}{23mm}
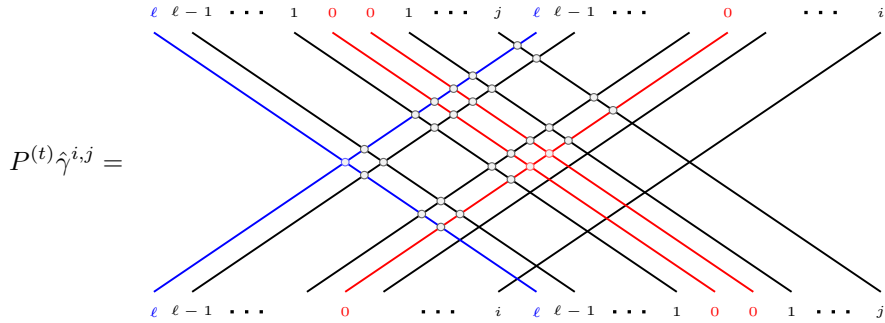
{
\begin{braid}\tikzset{baseline=-.3em}
	\draw[blue] (0,6) node{$\ell$};
	\draw (1.5,6) node{$\ell-1$};
	\draw[dots] (3,6)--(4.5,6);
	\draw (5.5,6) node{$1$};
	\draw[red] (7,6) node{$0$};
	\draw[red] (8.5,6) node{$0$};
	\draw (10,6) node{$1$};
	\draw[dots] (11,6)--(12.5,6);
	\draw (13.5,6) node{$j$};
	\draw[blue] (15,6) node{$\ell$};
    \draw (16.5,6) node{$\ell-1$};
	\draw[dots] (18,6)--(19.5,6);
	\draw (21,6) node{$ $};
	\draw[red] (22.5,6) node{$0$};
	\draw (24,6) node{$ $};
	\draw[dots] (25.5,6)--(27,6);
	\draw (28.5,6) node{$i$};

	\draw[blue] (0,-6) node{$\ell$};
    \draw (1.5,-6) node{$\ell-1$};
	\draw[dots] (3,-6)--(4.5,-6);
	\draw (6,-6) node{$ $};
	\draw[red] (7.5,-6) node{$0$};
	\draw (9,-6) node{$ $};
	\draw[dots] (10.5,-6)--(12,-6);
	\draw (13.5,-6) node{$i$};
	\draw[blue] (15,-6) node{$\ell$};
	\draw (16.5,-6) node{$\ell-1$};
	\draw[dots] (18,-6)--(19.5,-6);
	\draw (20.5,-6) node{$1$};
	\draw[red] (22,-6) node{$0$};
	\draw[red] (23.5,-6) node{$0$};
	\draw (25,-6) node{$1$};
	\draw[dots] (26,-6)--(27.5,-6);
	\draw (28.5,-6) node{$j$};

	\draw[blue](0,-5.2)--(15,5.2);
	\draw(1.5,-5.2)--(16.5,5.2);
	\draw(6,-5.2)--(21,5.2);
	\draw[red](7.5,-5.2)--(22.5,5.2);
	\draw(9,-5.2)--(24,5.2);
	\draw(13.5,-5.2)--(28.5,5.2);

	\draw[blue](0,5.2)--(15,-5.2);
	\draw(1.5,5.2)--(16.5,-5.2);
	\draw(5.5,5.2)--(20.5,-5.2);
	\draw[red](7,5.2)--(22,-5.2);
	\draw[red](8.5,5.2)--(23.5,-5.2);
	\draw(10,5.2)--(25,-5.2);
	\draw(13.5,5.2)--(28.5,-5.2);

	\filldraw[color=blue!60, fill=black!5, thin] (7.5,0) circle (0.15);
	\filldraw[color=black!60, fill=black!5, thin] (8.25,0.52) circle (0.15);
	\filldraw[color=black!60, fill=black!5, thin] (10.25,1.907) circle (0.15);
	\filldraw[color=black!60, fill=black!5, thin] (11,2.427) circle (0.15);
	\filldraw[color=black!60, fill=black!5, thin] (11.75,2.947) circle (0.15);
	\filldraw[color=black!60, fill=black!5, thin] (12.5,3.467) circle (0.15);
	\filldraw[color=black!60, fill=black!5, thin] (14.25,4.68) circle (0.15);

	\filldraw[color=black!60, fill=black!5, thin] (8.25,-0.52) circle (0.15);
	\filldraw[color=black!60, fill=black!5, thin] (9,0) circle (0.15);
	\filldraw[color=black!60, fill=black!5, thin] (11,1.387) circle (0.15);
	\filldraw[color=black!60, fill=black!5, thin] (11.75,1.907) circle (0.15);
	\filldraw[color=black!60, fill=black!5, thin] (12.5,2.427) circle (0.15);
	\filldraw[color=black!60, fill=black!5, thin] (13.25,2.947) circle (0.15);
	\filldraw[color=black!60, fill=black!5, thin] (15,4.16) circle (0.15);

	\filldraw[color=black!60, fill=black!5, thin] (10.5,-2.08) circle (0.15);
	\filldraw[color=black!60, fill=black!5, thin] (11.25,-1.56) circle (0.15);
	\filldraw[color=black!60, fill=black!5, thin] (13.25,-0.173) circle (0.15);
	\filldraw[color=black!60, fill=black!5, thin] (14,0.347) circle (0.15);
    \filldraw[color=black!60, fill=black!5, thin] (14.75,0.867) circle (0.15);
    \filldraw[color=black!60, fill=black!5, thin] (15.5,1.387) circle (0.15);
	\filldraw[color=black!60, fill=black!5, thin] (17.25,2.6) circle (0.15);

	\filldraw[color=black!60, fill=black!5, thin] (11.25,-2.6) circle (0.15);
	\filldraw[color=black!60, fill=black!5, thin] (12,-2.08) circle (0.15);
	\filldraw[color=black!60, fill=black!5, thin] (14,-0.693) circle (0.15);
	\filldraw[color=red!60, fill=black!5, thin] (14.75,-0.173) circle (0.15);
	\filldraw[color=red!60, fill=black!5, thin] (15.5,0.347) circle (0.15);
	\filldraw[color=black!60, fill=black!5, thin] (16.25,0.867) circle (0.15);
	\filldraw[color=black!60, fill=black!5, thin] (18,2.08) circle (0.15);
	\end{braid}
	}.
$$
Using (\ref{EOpening0}), this is equal to
$$
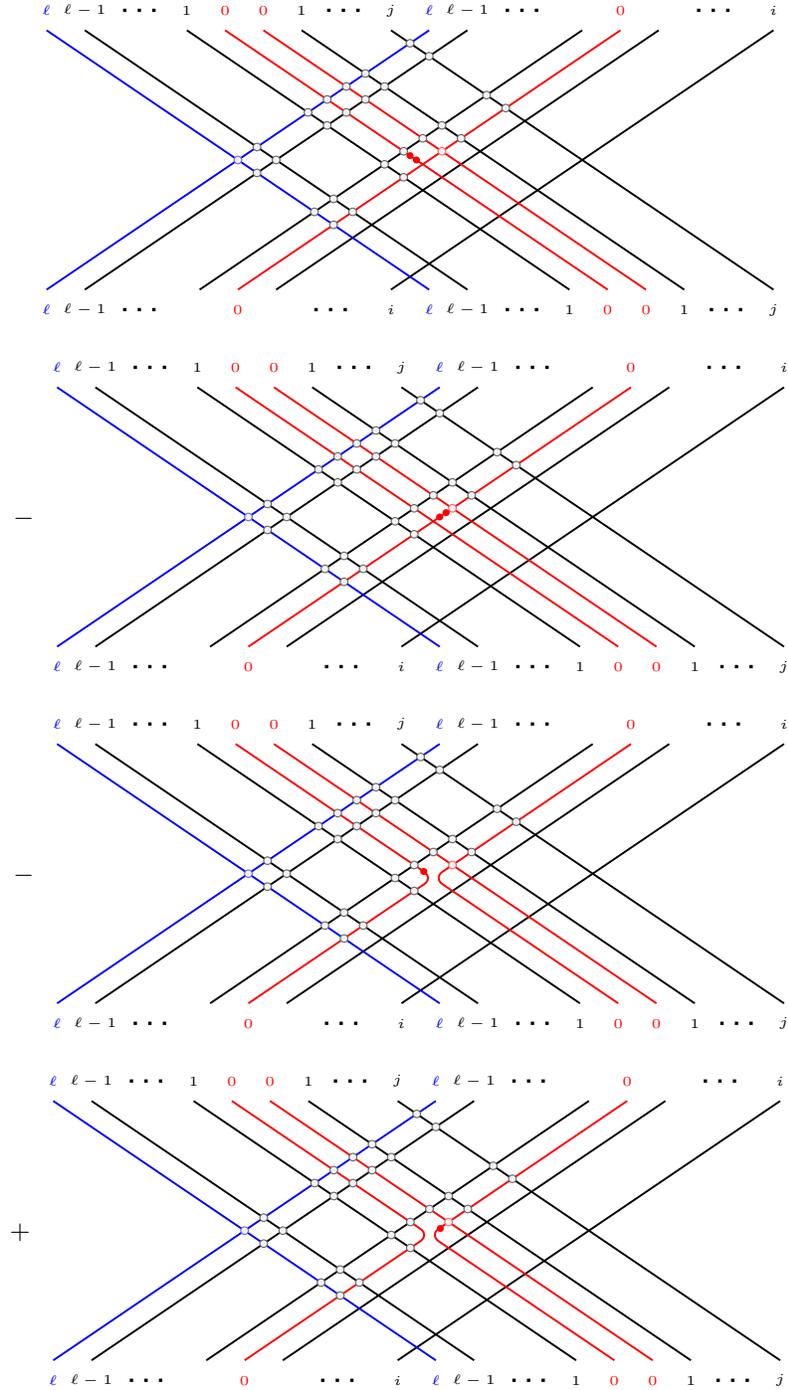
\resizebox{104mm}{23mm}{
\begin{braid}\tikzset{baseline=-.3em}
	\draw[blue] (0,6) node{$\ell$};
	\draw (1.5,6) node{$\ell-1$};
	\draw[dots] (3,6)--(4.5,6);
	\draw (5.5,6) node{$1$};
	\draw[red] (7,6) node{$0$};
	\draw[red] (8.5,6) node{$0$};
	\draw (10,6) node{$1$};
	\draw[dots] (11,6)--(12.5,6);
	\draw (13.5,6) node{$j$};
	\draw[blue] (15,6) node{$\ell$};
    \draw (16.5,6) node{$\ell-1$};
	\draw[dots] (18,6)--(19.5,6);
	\draw (21,6) node{$ $};
	\draw[red] (22.5,6) node{$0$};
	\draw (24,6) node{$ $};
	\draw[dots] (25.5,6)--(27,6);
	\draw (28.5,6) node{$i$};

	\draw[blue] (0,-6) node{$\ell$};
    \draw (1.5,-6) node{$\ell-1$};
	\draw[dots] (3,-6)--(4.5,-6);
	\draw (6,-6) node{$ $};
	\draw[red] (7.5,-6) node{$0$};
	\draw (9,-6) node{$ $};
	\draw[dots] (10.5,-6)--(12,-6);
	\draw (13.5,-6) node{$i$};
	\draw[blue] (15,-6) node{$\ell$};
	\draw (16.5,-6) node{$\ell-1$};
	\draw[dots] (18,-6)--(19.5,-6);
	\draw (20.5,-6) node{$1$};
	\draw[red] (22,-6) node{$0$};
	\draw[red] (23.5,-6) node{$0$};
	\draw (25,-6) node{$1$};
	\draw[dots] (26,-6)--(27.5,-6);
	\draw (28.5,-6) node{$j$};

	\draw[blue](0,-5.2)--(15,5.2);
	\draw(1.5,-5.2)--(16.5,5.2);
	\draw(6,-5.2)--(21,5.2);
	\draw[red](7.5,-5.2)--(22.5,5.2);
	\draw(9,-5.2)--(24,5.2);
	\draw(13.5,-5.2)--(28.5,5.2);

	\draw[blue](0,5.2)--(15,-5.2);
	\draw(1.5,5.2)--(16.5,-5.2);
	\draw(5.5,5.2)--(20.5,-5.2);
	\draw[red](7,5.2)--(22,-5.2);
	\draw[red](8.5,5.2)--(23.5,-5.2);
	\draw(10,5.2)--(25,-5.2);
	\draw(13.5,5.2)--(28.5,-5.2);

	\filldraw[color=blue!60, fill=black!5, thin] (7.5,0) circle (0.15);
	\filldraw[color=black!60, fill=black!5, thin] (8.25,0.52) circle (0.15);
	\filldraw[color=black!60, fill=black!5, thin] (10.25,1.907) circle (0.15);
	\filldraw[color=black!60, fill=black!5, thin] (11,2.427) circle (0.15);
	\filldraw[color=black!60, fill=black!5, thin] (11.75,2.947) circle (0.15);
	\filldraw[color=black!60, fill=black!5, thin] (12.5,3.467) circle (0.15);
	\filldraw[color=black!60, fill=black!5, thin] (14.25,4.68) circle (0.15);

	\filldraw[color=black!60, fill=black!5, thin] (8.25,-0.52) circle (0.15);
	\filldraw[color=black!60, fill=black!5, thin] (9,0) circle (0.15);
	\filldraw[color=black!60, fill=black!5, thin] (11,1.387) circle (0.15);
	\filldraw[color=black!60, fill=black!5, thin] (11.75,1.907) circle (0.15);
	\filldraw[color=black!60, fill=black!5, thin] (12.5,2.427) circle (0.15);
	\filldraw[color=black!60, fill=black!5, thin] (13.25,2.947) circle (0.15);
	\filldraw[color=black!60, fill=black!5, thin] (15,4.16) circle (0.15);

	\filldraw[color=black!60, fill=black!5, thin] (10.5,-2.08) circle (0.15);
	\filldraw[color=black!60, fill=black!5, thin] (11.25,-1.56) circle (0.15);
	\filldraw[color=black!60, fill=black!5, thin] (13.25,-0.173) circle (0.15);
	\filldraw[color=black!60, fill=black!5, thin] (14,0.347) circle (0.15);
    \filldraw[color=black!60, fill=black!5, thin] (14.75,0.867) circle (0.15);
    \filldraw[color=black!60, fill=black!5, thin] (15.5,1.387) circle (0.15);
	\filldraw[color=black!60, fill=black!5, thin] (17.25,2.6) circle (0.15);

	\filldraw[color=black!60, fill=black!5, thin] (11.25,-2.6) circle (0.15);
	\filldraw[color=black!60, fill=black!5, thin] (12,-2.08) circle (0.15);
	\filldraw[color=black!60, fill=black!5, thin] (14,-0.693) circle (0.15);
	\filldraw[color=red!60, fill=black!5, thin] (15.5,0.347) circle (0.15);
	\filldraw[color=black!60, fill=black!5, thin] (16.25,0.867) circle (0.15);
	\filldraw[color=black!60, fill=black!5, thin] (18,2.08) circle (0.15);

	\smallreddot (14.5,0);
	\smallreddot (14.25,0.174);

	\end{braid}
	}
$$

$$
-\resizebox{104mm}{23mm}{
\begin{braid}\tikzset{baseline=-.3em}
	\draw[blue] (0,6) node{$\ell$};
	\draw (1.5,6) node{$\ell-1$};
	\draw[dots] (3,6)--(4.5,6);
	\draw (5.5,6) node{$1$};
	\draw[red] (7,6) node{$0$};
	\draw[red] (8.5,6) node{$0$};
	\draw (10,6) node{$1$};
	\draw[dots] (11,6)--(12.5,6);
	\draw (13.5,6) node{$j$};
	\draw[blue] (15,6) node{$\ell$};
    \draw (16.5,6) node{$\ell-1$};
	\draw[dots] (18,6)--(19.5,6);
	\draw (21,6) node{$ $};
	\draw[red] (22.5,6) node{$0$};
	\draw (24,6) node{$ $};
	\draw[dots] (25.5,6)--(27,6);
	\draw (28.5,6) node{$i$};

	\draw[blue] (0,-6) node{$\ell$};
    \draw (1.5,-6) node{$\ell-1$};
	\draw[dots] (3,-6)--(4.5,-6);
	\draw (6,-6) node{$ $};
	\draw[red] (7.5,-6) node{$0$};
	\draw (9,-6) node{$ $};
	\draw[dots] (10.5,-6)--(12,-6);
	\draw (13.5,-6) node{$i$};
	\draw[blue] (15,-6) node{$\ell$};
	\draw (16.5,-6) node{$\ell-1$};
	\draw[dots] (18,-6)--(19.5,-6);
	\draw (20.5,-6) node{$1$};
	\draw[red] (22,-6) node{$0$};
	\draw[red] (23.5,-6) node{$0$};
	\draw (25,-6) node{$1$};
	\draw[dots] (26,-6)--(27.5,-6);
	\draw (28.5,-6) node{$j$};

	\draw[blue](0,-5.2)--(15,5.2);
	\draw(1.5,-5.2)--(16.5,5.2);
	\draw(6,-5.2)--(21,5.2);
	\draw[red](7.5,-5.2)--(22.5,5.2);
	\draw(9,-5.2)--(24,5.2);
	\draw(13.5,-5.2)--(28.5,5.2);

	\draw[blue](0,5.2)--(15,-5.2);
	\draw(1.5,5.2)--(16.5,-5.2);
	\draw(5.5,5.2)--(20.5,-5.2);
	\draw[red](7,5.2)--(22,-5.2);
	\draw[red](8.5,5.2)--(23.5,-5.2);
	\draw(10,5.2)--(25,-5.2);
	\draw(13.5,5.2)--(28.5,-5.2);

	\filldraw[color=blue!60, fill=black!5, thin] (7.5,0) circle (0.15);
	\filldraw[color=black!60, fill=black!5, thin] (8.25,0.52) circle (0.15);
	\filldraw[color=black!60, fill=black!5, thin] (10.25,1.907) circle (0.15);
	\filldraw[color=black!60, fill=black!5, thin] (11,2.427) circle (0.15);
	\filldraw[color=black!60, fill=black!5, thin] (11.75,2.947) circle (0.15);
	\filldraw[color=black!60, fill=black!5, thin] (12.5,3.467) circle (0.15);
	\filldraw[color=black!60, fill=black!5, thin] (14.25,4.68) circle (0.15);

	\filldraw[color=black!60, fill=black!5, thin] (8.25,-0.52) circle (0.15);
	\filldraw[color=black!60, fill=black!5, thin] (9,0) circle (0.15);
	\filldraw[color=black!60, fill=black!5, thin] (11,1.387) circle (0.15);
	\filldraw[color=black!60, fill=black!5, thin] (11.75,1.907) circle (0.15);
	\filldraw[color=black!60, fill=black!5, thin] (12.5,2.427) circle (0.15);
	\filldraw[color=black!60, fill=black!5, thin] (13.25,2.947) circle (0.15);
	\filldraw[color=black!60, fill=black!5, thin] (15,4.16) circle (0.15);

	\filldraw[color=black!60, fill=black!5, thin] (10.5,-2.08) circle (0.15);
	\filldraw[color=black!60, fill=black!5, thin] (11.25,-1.56) circle (0.15);
	\filldraw[color=black!60, fill=black!5, thin] (13.25,-0.173) circle (0.15);
	\filldraw[color=black!60, fill=black!5, thin] (14,0.347) circle (0.15);
    \filldraw[color=black!60, fill=black!5, thin] (14.75,0.867) circle (0.15);
    \filldraw[color=black!60, fill=black!5, thin] (15.5,1.387) circle (0.15);
	\filldraw[color=black!60, fill=black!5, thin] (17.25,2.6) circle (0.15);

	\filldraw[color=black!60, fill=black!5, thin] (11.25,-2.6) circle (0.15);
	\filldraw[color=black!60, fill=black!5, thin] (12,-2.08) circle (0.15);
	\filldraw[color=black!60, fill=black!5, thin] (14,-0.693) circle (0.15);
	\filldraw[color=red!60, fill=black!5, thin] (15.5,0.347) circle (0.15);
	\filldraw[color=black!60, fill=black!5, thin] (16.25,0.867) circle (0.15);
	\filldraw[color=black!60, fill=black!5, thin] (18,2.08) circle (0.15);

	\smallreddot (15,0);
	\smallreddot (15.25,0.174);

	\end{braid}
	}
$$

$$
-\resizebox{104mm}{23mm}{
\begin{braid}\tikzset{baseline=-.3em}
	\draw[blue] (0,6) node{$\ell$};
	\draw (1.5,6) node{$\ell-1$};
	\draw[dots] (3,6)--(4.5,6);
	\draw (5.5,6) node{$1$};
	\draw[red] (7,6) node{$0$};
	\draw[red] (8.5,6) node{$0$};
	\draw (10,6) node{$1$};
	\draw[dots] (11,6)--(12.5,6);
	\draw (13.5,6) node{$j$};
	\draw[blue] (15,6) node{$\ell$};
    \draw (16.5,6) node{$\ell-1$};
	\draw[dots] (18,6)--(19.5,6);
	\draw (21,6) node{$ $};
	\draw[red] (22.5,6) node{$0$};
	\draw (24,6) node{$ $};
	\draw[dots] (25.5,6)--(27,6);
	\draw (28.5,6) node{$i$};

	\draw[blue] (0,-6) node{$\ell$};
    \draw (1.5,-6) node{$\ell-1$};
	\draw[dots] (3,-6)--(4.5,-6);
	\draw (6,-6) node{$ $};
	\draw[red] (7.5,-6) node{$0$};
	\draw (9,-6) node{$ $};
	\draw[dots] (10.5,-6)--(12,-6);
	\draw (13.5,-6) node{$i$};
	\draw[blue] (15,-6) node{$\ell$};
	\draw (16.5,-6) node{$\ell-1$};
	\draw[dots] (18,-6)--(19.5,-6);
	\draw (20.5,-6) node{$1$};
	\draw[red] (22,-6) node{$0$};
	\draw[red] (23.5,-6) node{$0$};
	\draw (25,-6) node{$1$};
	\draw[dots] (26,-6)--(27.5,-6);
	\draw (28.5,-6) node{$j$};

	\draw[blue](0,-5.2)--(15,5.2);
	\draw(1.5,-5.2)--(16.5,5.2);
	\draw(6,-5.2)--(21,5.2);
	\draw(9,-5.2)--(24,5.2);
	\draw(13.5,-5.2)--(28.5,5.2);

	\draw[blue](0,5.2)--(15,-5.2);
	\draw(1.5,5.2)--(16.5,-5.2);
	\draw(5.5,5.2)--(20.5,-5.2);
	\draw[red](8.5,5.2)--(23.5,-5.2);
	\draw(10,5.2)--(25,-5.2);
	\draw(13.5,5.2)--(28.5,-5.2);

	\draw[red](7.5,-5.2)--(14.75,-0.173)--(7,5.2);
	\draw[red](22,-5.2)--(14.75,-0.173)--(22.5,5.2);

	\filldraw[color=blue!60, fill=black!5, thin] (7.5,0) circle (0.15);
	\filldraw[color=black!60, fill=black!5, thin] (8.25,0.52) circle (0.15);
	\filldraw[color=black!60, fill=black!5, thin] (10.25,1.907) circle (0.15);
	\filldraw[color=black!60, fill=black!5, thin] (11,2.427) circle (0.15);
	\filldraw[color=black!60, fill=black!5, thin] (11.75,2.947) circle (0.15);
	\filldraw[color=black!60, fill=black!5, thin] (12.5,3.467) circle (0.15);
	\filldraw[color=black!60, fill=black!5, thin] (14.25,4.68) circle (0.15);

	\filldraw[color=black!60, fill=black!5, thin] (8.25,-0.52) circle (0.15);
	\filldraw[color=black!60, fill=black!5, thin] (9,0) circle (0.15);
	\filldraw[color=black!60, fill=black!5, thin] (11,1.387) circle (0.15);
	\filldraw[color=black!60, fill=black!5, thin] (11.75,1.907) circle (0.15);
	\filldraw[color=black!60, fill=black!5, thin] (12.5,2.427) circle (0.15);
	\filldraw[color=black!60, fill=black!5, thin] (13.25,2.947) circle (0.15);
	\filldraw[color=black!60, fill=black!5, thin] (15,4.16) circle (0.15);

	\filldraw[color=black!60, fill=black!5, thin] (10.5,-2.08) circle (0.15);
	\filldraw[color=black!60, fill=black!5, thin] (11.25,-1.56) circle (0.15);
	\filldraw[color=black!60, fill=black!5, thin] (13.25,-0.173) circle (0.15);
	\filldraw[color=black!60, fill=black!5, thin] (14,0.347) circle (0.15);
    \filldraw[color=black!60, fill=black!5, thin] (14.75,0.867) circle (0.15);
    \filldraw[color=black!60, fill=black!5, thin] (15.5,1.387) circle (0.15);
	\filldraw[color=black!60, fill=black!5, thin] (17.25,2.6) circle (0.15);

	\filldraw[color=black!60, fill=black!5, thin] (11.25,-2.6) circle (0.15);
	\filldraw[color=black!60, fill=black!5, thin] (12,-2.08) circle (0.15);
	\filldraw[color=black!60, fill=black!5, thin] (14,-0.693) circle (0.15);
	\filldraw[color=red!60, fill=black!5, thin] (15.5,0.347) circle (0.15);
	\filldraw[color=black!60, fill=black!5, thin] (16.25,0.867) circle (0.15);
	\filldraw[color=black!60, fill=black!5, thin] (18,2.08) circle (0.15);

	\smallreddot (14.375,0.087);

	\end{braid}
	}
$$

$$
+\resizebox{104mm}{23mm}{
\begin{braid}\tikzset{baseline=-.3em}
	\draw[blue] (0,6) node{$\ell$};
	\draw (1.5,6) node{$\ell-1$};
	\draw[dots] (3,6)--(4.5,6);
	\draw (5.5,6) node{$1$};
	\draw[red] (7,6) node{$0$};
	\draw[red] (8.5,6) node{$0$};
	\draw (10,6) node{$1$};
	\draw[dots] (11,6)--(12.5,6);
	\draw (13.5,6) node{$j$};
	\draw[blue] (15,6) node{$\ell$};
    \draw (16.5,6) node{$\ell-1$};
	\draw[dots] (18,6)--(19.5,6);
	\draw (21,6) node{$ $};
	\draw[red] (22.5,6) node{$0$};
	\draw (24,6) node{$ $};
	\draw[dots] (25.5,6)--(27,6);
	\draw (28.5,6) node{$i$};

	\draw[blue] (0,-6) node{$\ell$};
    \draw (1.5,-6) node{$\ell-1$};
	\draw[dots] (3,-6)--(4.5,-6);
	\draw (6,-6) node{$ $};
	\draw[red] (7.5,-6) node{$0$};
	\draw (9,-6) node{$ $};
	\draw[dots] (10.5,-6)--(12,-6);
	\draw (13.5,-6) node{$i$};
	\draw[blue] (15,-6) node{$\ell$};
	\draw (16.5,-6) node{$\ell-1$};
	\draw[dots] (18,-6)--(19.5,-6);
	\draw (20.5,-6) node{$1$};
	\draw[red] (22,-6) node{$0$};
	\draw[red] (23.5,-6) node{$0$};
	\draw (25,-6) node{$1$};
	\draw[dots] (26,-6)--(27.5,-6);
	\draw (28.5,-6) node{$j$};

	\draw[blue](0,-5.2)--(15,5.2);
	\draw(1.5,-5.2)--(16.5,5.2);
	\draw(6,-5.2)--(21,5.2);
	\draw(9,-5.2)--(24,5.2);
	\draw(13.5,-5.2)--(28.5,5.2);

	\draw[blue](0,5.2)--(15,-5.2);
	\draw(1.5,5.2)--(16.5,-5.2);
	\draw(5.5,5.2)--(20.5,-5.2);
	\draw[red](8.5,5.2)--(23.5,-5.2);
	\draw(10,5.2)--(25,-5.2);
	\draw(13.5,5.2)--(28.5,-5.2);

	\draw[red](7.5,-5.2)--(14.75,-0.173)--(7,5.2);
	\draw[red](22,-5.2)--(14.75,-0.173)--(22.5,5.2);

	\filldraw[color=blue!60, fill=black!5, thin] (7.5,0) circle (0.15);
	\filldraw[color=black!60, fill=black!5, thin] (8.25,0.52) circle (0.15);
	\filldraw[color=black!60, fill=black!5, thin] (10.25,1.907) circle (0.15);
	\filldraw[color=black!60, fill=black!5, thin] (11,2.427) circle (0.15);
	\filldraw[color=black!60, fill=black!5, thin] (11.75,2.947) circle (0.15);
	\filldraw[color=black!60, fill=black!5, thin] (12.5,3.467) circle (0.15);
	\filldraw[color=black!60, fill=black!5, thin] (14.25,4.68) circle (0.15);

	\filldraw[color=black!60, fill=black!5, thin] (8.25,-0.52) circle (0.15);
	\filldraw[color=black!60, fill=black!5, thin] (9,0) circle (0.15);
	\filldraw[color=black!60, fill=black!5, thin] (11,1.387) circle (0.15);
	\filldraw[color=black!60, fill=black!5, thin] (11.75,1.907) circle (0.15);
	\filldraw[color=black!60, fill=black!5, thin] (12.5,2.427) circle (0.15);
	\filldraw[color=black!60, fill=black!5, thin] (13.25,2.947) circle (0.15);
	\filldraw[color=black!60, fill=black!5, thin] (15,4.16) circle (0.15);

	\filldraw[color=black!60, fill=black!5, thin] (10.5,-2.08) circle (0.15);
	\filldraw[color=black!60, fill=black!5, thin] (11.25,-1.56) circle (0.15);
	\filldraw[color=black!60, fill=black!5, thin] (13.25,-0.173) circle (0.15);
	\filldraw[color=black!60, fill=black!5, thin] (14,0.347) circle (0.15);
    \filldraw[color=black!60, fill=black!5, thin] (14.75,0.867) circle (0.15);
    \filldraw[color=black!60, fill=black!5, thin] (15.5,1.387) circle (0.15);
	\filldraw[color=black!60, fill=black!5, thin] (17.25,2.6) circle (0.15);

	\filldraw[color=black!60, fill=black!5, thin] (11.25,-2.6) circle (0.15);
	\filldraw[color=black!60, fill=black!5, thin] (12,-2.08) circle (0.15);
	\filldraw[color=black!60, fill=black!5, thin] (14,-0.693) circle (0.15);
	\filldraw[color=red!60, fill=black!5, thin] (15.5,0.347) circle (0.15);
	\filldraw[color=black!60, fill=black!5, thin] (16.25,0.867) circle (0.15);
	\filldraw[color=black!60, fill=black!5, thin] (18,2.08) circle (0.15);

	\smallreddot (15.175,0.087);

	\end{braid}
	}.
$$
Similar to the $g \neq 0$ case, we apply (\ref{EDotPastCircle}) and (\ref{ECircleBraid}) several times and then Lemma \ref{LCuspExpl} to obtain that the last two summands are zero. Following similar steps to the $g \neq 0$ case, the first summand is equal to
$$
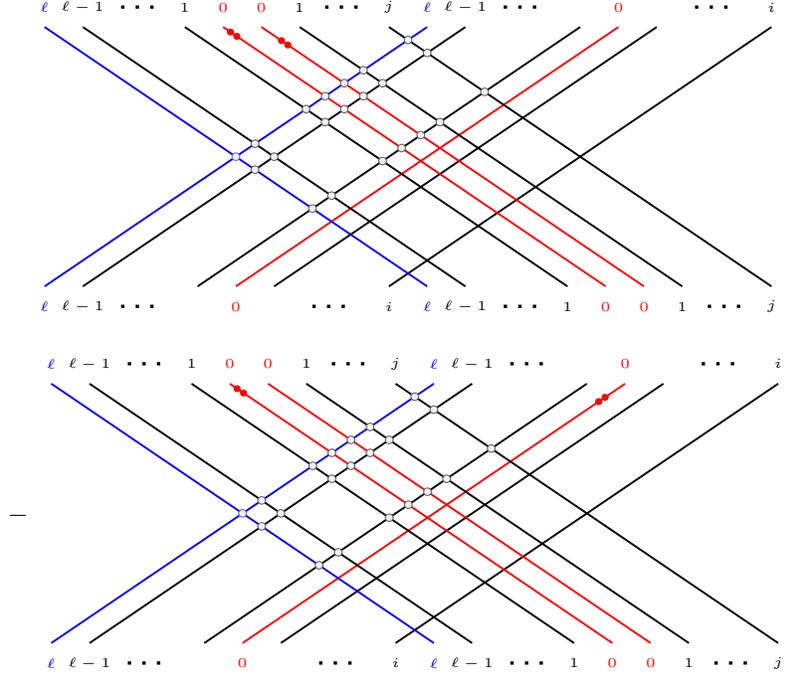
\resizebox{104mm}{23mm}{
\begin{braid}\tikzset{baseline=-.3em}
	\draw[blue] (0,6) node{$\ell$};
	\draw (1.5,6) node{$\ell-1$};
	\draw[dots] (3,6)--(4.5,6);
	\draw (5.5,6) node{$1$};
	\draw[red] (7,6) node{$0$};
	\draw[red] (8.5,6) node{$0$};
	\draw (10,6) node{$1$};
	\draw[dots] (11,6)--(12.5,6);
	\draw (13.5,6) node{$j$};
	\draw[blue] (15,6) node{$\ell$};
    \draw (16.5,6) node{$\ell-1$};
	\draw[dots] (18,6)--(19.5,6);
	\draw (21,6) node{$ $};
	\draw[red] (22.5,6) node{$0$};
	\draw (24,6) node{$ $};
	\draw[dots] (25.5,6)--(27,6);
	\draw (28.5,6) node{$i$};

	\draw[blue] (0,-6) node{$\ell$};
    \draw (1.5,-6) node{$\ell-1$};
	\draw[dots] (3,-6)--(4.5,-6);
	\draw (6,-6) node{$ $};
	\draw[red] (7.5,-6) node{$0$};
	\draw (9,-6) node{$ $};
	\draw[dots] (10.5,-6)--(12,-6);
	\draw (13.5,-6) node{$i$};
	\draw[blue] (15,-6) node{$\ell$};
	\draw (16.5,-6) node{$\ell-1$};
	\draw[dots] (18,-6)--(19.5,-6);
	\draw (20.5,-6) node{$1$};
	\draw[red] (22,-6) node{$0$};
	\draw[red] (23.5,-6) node{$0$};
	\draw (25,-6) node{$1$};
	\draw[dots] (26,-6)--(27.5,-6);
	\draw (28.5,-6) node{$j$};

	\draw[blue](0,-5.2)--(15,5.2);
	\draw(1.5,-5.2)--(16.5,5.2);
	\draw(6,-5.2)--(21,5.2);
	\draw[red](7.5,-5.2)--(22.5,5.2);
	\draw(9,-5.2)--(24,5.2);
	\draw(13.5,-5.2)--(28.5,5.2);

	\draw[blue](0,5.2)--(15,-5.2);
	\draw(1.5,5.2)--(16.5,-5.2);
	\draw(5.5,5.2)--(20.5,-5.2);
	\draw[red](7,5.2)--(22,-5.2);
	\draw[red](8.5,5.2)--(23.5,-5.2);
	\draw(10,5.2)--(25,-5.2);
	\draw(13.5,5.2)--(28.5,-5.2);

	\filldraw[color=blue!60, fill=black!5, thin] (7.5,0) circle (0.15);
	\filldraw[color=black!60, fill=black!5, thin] (8.25,0.52) circle (0.15);
	\filldraw[color=black!60, fill=black!5, thin] (10.25,1.907) circle (0.15);
	\filldraw[color=black!60, fill=black!5, thin] (11,2.427) circle (0.15);
	\filldraw[color=black!60, fill=black!5, thin] (11.75,2.947) circle (0.15);
	\filldraw[color=black!60, fill=black!5, thin] (12.5,3.467) circle (0.15);
	\filldraw[color=black!60, fill=black!5, thin] (14.25,4.68) circle (0.15);

	\filldraw[color=black!60, fill=black!5, thin] (8.25,-0.52) circle (0.15);
	\filldraw[color=black!60, fill=black!5, thin] (9,0) circle (0.15);
	\filldraw[color=black!60, fill=black!5, thin] (11,1.387) circle (0.15);
	\filldraw[color=black!60, fill=black!5, thin] (11.75,1.907) circle (0.15);
	\filldraw[color=black!60, fill=black!5, thin] (12.5,2.427) circle (0.15);
	\filldraw[color=black!60, fill=black!5, thin] (13.25,2.947) circle (0.15);
	\filldraw[color=black!60, fill=black!5, thin] (15,4.16) circle (0.15);

	\filldraw[color=black!60, fill=black!5, thin] (10.5,-2.08) circle (0.15);
	\filldraw[color=black!60, fill=black!5, thin] (11.25,-1.56) circle (0.15);
	\filldraw[color=black!60, fill=black!5, thin] (13.25,-0.173) circle (0.15);
	\filldraw[color=black!60, fill=black!5, thin] (14,0.347) circle (0.15);
    \filldraw[color=black!60, fill=black!5, thin] (14.75,0.867) circle (0.15);
    \filldraw[color=black!60, fill=black!5, thin] (15.5,1.387) circle (0.15);
	\filldraw[color=black!60, fill=black!5, thin] (17.25,2.6) circle (0.15);

	\smallreddot (7.288,5);
	\smallreddot (7.538,4.827);

	\smallreddot (9.288,4.654);
	\smallreddot (9.537,4.481);

	\end{braid}
	}
$$

$$
-\resizebox{104mm}{23mm}{
\begin{braid}\tikzset{baseline=-.3em}
	\draw[blue] (0,6) node{$\ell$};
	\draw (1.5,6) node{$\ell-1$};
	\draw[dots] (3,6)--(4.5,6);
	\draw (5.5,6) node{$1$};
	\draw[red] (7,6) node{$0$};
	\draw[red] (8.5,6) node{$0$};
	\draw (10,6) node{$1$};
	\draw[dots] (11,6)--(12.5,6);
	\draw (13.5,6) node{$j$};
	\draw[blue] (15,6) node{$\ell$};
    \draw (16.5,6) node{$\ell-1$};
	\draw[dots] (18,6)--(19.5,6);
	\draw (21,6) node{$ $};
	\draw[red] (22.5,6) node{$0$};
	\draw (24,6) node{$ $};
	\draw[dots] (25.5,6)--(27,6);
	\draw (28.5,6) node{$i$};

	\draw[blue] (0,-6) node{$\ell$};
    \draw (1.5,-6) node{$\ell-1$};
	\draw[dots] (3,-6)--(4.5,-6);
	\draw (6,-6) node{$ $};
	\draw[red] (7.5,-6) node{$0$};
	\draw (9,-6) node{$ $};
	\draw[dots] (10.5,-6)--(12,-6);
	\draw (13.5,-6) node{$i$};
	\draw[blue] (15,-6) node{$\ell$};
	\draw (16.5,-6) node{$\ell-1$};
	\draw[dots] (18,-6)--(19.5,-6);
	\draw (20.5,-6) node{$1$};
	\draw[red] (22,-6) node{$0$};
	\draw[red] (23.5,-6) node{$0$};
	\draw (25,-6) node{$1$};
	\draw[dots] (26,-6)--(27.5,-6);
	\draw (28.5,-6) node{$j$};

	\draw[blue](0,-5.2)--(15,5.2);
	\draw(1.5,-5.2)--(16.5,5.2);
	\draw(6,-5.2)--(21,5.2);
	\draw[red](7.5,-5.2)--(22.5,5.2);
	\draw(9,-5.2)--(24,5.2);
	\draw(13.5,-5.2)--(28.5,5.2);

	\draw[blue](0,5.2)--(15,-5.2);
	\draw(1.5,5.2)--(16.5,-5.2);
	\draw(5.5,5.2)--(20.5,-5.2);
	\draw[red](7,5.2)--(22,-5.2);
	\draw[red](8.5,5.2)--(23.5,-5.2);
	\draw(10,5.2)--(25,-5.2);
	\draw(13.5,5.2)--(28.5,-5.2);

	\filldraw[color=blue!60, fill=black!5, thin] (7.5,0) circle (0.15);
	\filldraw[color=black!60, fill=black!5, thin] (8.25,0.52) circle (0.15);
	\filldraw[color=black!60, fill=black!5, thin] (10.25,1.907) circle (0.15);
	\filldraw[color=black!60, fill=black!5, thin] (11,2.427) circle (0.15);
	\filldraw[color=black!60, fill=black!5, thin] (11.75,2.947) circle (0.15);
	\filldraw[color=black!60, fill=black!5, thin] (12.5,3.467) circle (0.15);
	\filldraw[color=black!60, fill=black!5, thin] (14.25,4.68) circle (0.15);

	\filldraw[color=black!60, fill=black!5, thin] (8.25,-0.52) circle (0.15);
	\filldraw[color=black!60, fill=black!5, thin] (9,0) circle (0.15);
	\filldraw[color=black!60, fill=black!5, thin] (11,1.387) circle (0.15);
	\filldraw[color=black!60, fill=black!5, thin] (11.75,1.907) circle (0.15);
	\filldraw[color=black!60, fill=black!5, thin] (12.5,2.427) circle (0.15);
	\filldraw[color=black!60, fill=black!5, thin] (13.25,2.947) circle (0.15);
	\filldraw[color=black!60, fill=black!5, thin] (15,4.16) circle (0.15);

	\filldraw[color=black!60, fill=black!5, thin] (10.5,-2.08) circle (0.15);
	\filldraw[color=black!60, fill=black!5, thin] (11.25,-1.56) circle (0.15);
	\filldraw[color=black!60, fill=black!5, thin] (13.25,-0.173) circle (0.15);
	\filldraw[color=black!60, fill=black!5, thin] (14,0.347) circle (0.15);
    \filldraw[color=black!60, fill=black!5, thin] (14.75,0.867) circle (0.15);
    \filldraw[color=black!60, fill=black!5, thin] (15.5,1.387) circle (0.15);
	\filldraw[color=black!60, fill=black!5, thin] (17.25,2.6) circle (0.15);

	\smallreddot (7.288,5);
	\smallreddot (7.538,4.827);

	\smallreddot (21.713,4.654);
	\smallreddot (21.463,4.481);
	
	\end{braid}
	},
$$
which, in non-diagrammatic form, is $y_{r^j_0}^2(y_{s^j_0}^2-y_{t+p}^2)P^{(t-1)}{\hat \gamma^{i,j}}$. Similarly the second summand is equal to $y_{t+p}^2(y_{s^j_0}^2-y_{t+p}^2)P^{(t-1)}{\hat \gamma^{i,j}}$. Therefore,
$$
P^{(t)}{\hat \gamma^{j,i}}=(y_{r^j_g}^2-y_{t+p}^2)(y_{s^j_g}^2-y_{t+p}^2)P^{(t-1)}{\hat \gamma^{i,j}}
$$
and the claim is proved.

For the second claim, we again prove a stronger statement. Namely that, in ${\bar R}_{2\de}$,
$$P^{(1)}{\hat \ga^{i,j}} = \Upsilon {\hat \ga^{i,j}} + (y_1 - y_{p+1})P^{(0)}{\hat \ga^{i,j}}.
$$
Now, by (\ref{EOpening}),
$$
P^{(1)}{\hat \ga}^{i,j}=
\begin{braid}\tikzset{baseline=-.3em}
	\draw[blue] (0,3) node{\scriptsize{$\ell$}};
	\draw (1.5,3) node{\scriptsize{$\ell-1$}};
	\draw[dots] (3,3)--(4.5,3);
	\draw (6,3) node{\scriptsize{$j$}};
	\draw (7.5,3) node{\scriptsize{$\ell$}};
	\draw (9,3) node{\scriptsize{$\ell-1$}};
	\draw[dots] (10.5,3)--(12,3);
	\draw (13.5,3) node{\scriptsize{$i$}};

	\draw[blue] (0,-3) node{\scriptsize{$\ell$}};
	\draw (1.5,-3) node{\scriptsize{$\ell-1$}};
	\draw[dots] (3,-3)--(4.5,-3);
	\draw (6,-3) node{\scriptsize{$i$}};
	\draw (7.5,-3) node{\scriptsize{$\ell$}};
	\draw (9,-3) node{\scriptsize{$\ell-1$}};
	\draw[dots] (10.5,3)--(12,3);
	\draw (13.5,-3) node{\scriptsize{$j$}};

	\draw[blue](0,-2.6)--(7.5,2.6);
	\draw(1.5,-2.6)--(9,2.6);
	\draw(6,-2.6)--(13.5,2.6);

	\draw[blue](0,2.6)--(7.5,-2.6);
	\draw(1.5,2.6)--(9,-2.6);
	\draw(6,2.6)--(13.5,-2.6);

	\filldraw[color=blue!60, fill=black!5, thin] (3.75,0) circle (0.17);
	\filldraw[color=black!60, fill=black!5, thin] (4.5,0.52) circle (0.17);
	\filldraw[color=black!60, fill=black!5, thin] (6.75,2.08) circle (0.17);
	\end{braid}
=
\begin{braid}\tikzset{baseline=-.3em}
	\draw[blue] (0,3) node{\scriptsize{$\ell$}};
	\draw (1.5,3) node{\scriptsize{$\ell-1$}};
	\draw[dots] (3,3)--(4.5,3);
	\draw (6,3) node{\scriptsize{$j$}};
	\draw (7.5,3) node{\scriptsize{$\ell$}};
	\draw (9,3) node{\scriptsize{$\ell-1$}};
	\draw[dots] (10.5,3)--(12,3);
	\draw (13.5,3) node{\scriptsize{$i$}};

	\draw[blue] (0,-3) node{\scriptsize{$\ell$}};
	\draw (1.5,-3) node{\scriptsize{$\ell-1$}};
	\draw[dots] (3,-3)--(4.5,-3);
	\draw (6,-3) node{\scriptsize{$i$}};
	\draw (7.5,-3) node{\scriptsize{$\ell$}};
	\draw (9,-3) node{\scriptsize{$\ell-1$}};
	\draw[dots] (10.5,3)--(12,3);
	\draw (13.5,-3) node{\scriptsize{$j$}};

	\draw(1.5,-2.6)--(9,2.6);
	\draw(6,-2.6)--(13.5,2.6);

	\draw(1.5,2.6)--(9,-2.6);
	\draw(6,2.6)--(13.5,-2.6);

	\draw[blue](0,-2.6)--(3.75,0)--(0,2.6);
    \draw[blue](7.5,-2.6)--(3.75,0)--(7.5,2.6);

	\filldraw[color=black!60, fill=black!5, thin] (4.5,0.52) circle (0.17);
	\filldraw[color=black!60, fill=black!5, thin] (6.75,2.08) circle (0.17);

	\end{braid}
$$

$$
+\begin{braid}\tikzset{baseline=-.3em}
	\draw[blue] (0,3) node{\scriptsize{$\ell$}};
	\draw (1.5,3) node{\scriptsize{$\ell-1$}};
	\draw[dots] (3,3)--(4.5,3);
	\draw (6,3) node{\scriptsize{$j$}};
	\draw (7.5,3) node{\scriptsize{$\ell$}};
	\draw (9,3) node{\scriptsize{$\ell-1$}};
	\draw[dots] (10.5,3)--(12,3);
	\draw (13.5,3) node{\scriptsize{$i$}};

	\draw[blue] (0,-3) node{\scriptsize{$\ell$}};
	\draw (1.5,-3) node{\scriptsize{$\ell-1$}};
	\draw[dots] (3,-3)--(4.5,-3);
	\draw (6,-3) node{\scriptsize{$i$}};
	\draw (7.5,-3) node{\scriptsize{$\ell$}};
	\draw (9,-3) node{\scriptsize{$\ell-1$}};
	\draw[dots] (10.5,3)--(12,3);
	\draw (13.5,-3) node{\scriptsize{$j$}};

	\draw[blue](0,-2.6)--(7.5,2.6);
	\draw(1.5,-2.6)--(9,2.6);
	\draw(6,-2.6)--(13.5,2.6);

	\draw[blue](0,2.6)--(7.5,-2.6);
	\draw(1.5,2.6)--(9,-2.6);
	\draw(6,2.6)--(13.5,-2.6);

	\filldraw[color=black!60, fill=black!5, thin] (4.5,0.52) circle (0.17);
	\filldraw[color=black!60, fill=black!5, thin] (6.75,2.08) circle (0.17);

	\bluedot (3.375,0.26);
	\end{braid}
-
\begin{braid}\tikzset{baseline=-.3em}
	\draw[blue] (0,3) node{\scriptsize{$\ell$}};
	\draw (1.5,3) node{\scriptsize{$\ell-1$}};
	\draw[dots] (3,3)--(4.5,3);
	\draw (6,3) node{\scriptsize{$j$}};
	\draw (7.5,3) node{\scriptsize{$\ell$}};
	\draw (9,3) node{\scriptsize{$\ell-1$}};
	\draw[dots] (10.5,3)--(12,3);
	\draw (13.5,3) node{\scriptsize{$i$}};

	\draw[blue] (0,-3) node{\scriptsize{$\ell$}};
	\draw (1.5,-3) node{\scriptsize{$\ell-1$}};
	\draw[dots] (3,-3)--(4.5,-3);
	\draw (6,-3) node{\scriptsize{$i$}};
	\draw (7.5,-3) node{\scriptsize{$\ell$}};
	\draw (9,-3) node{\scriptsize{$\ell-1$}};
	\draw[dots] (10.5,-3)--(12,-3);
	\draw (13.5,-3) node{\scriptsize{$j$}};

	\draw[blue](0,-2.6)--(7.5,2.6);
	\draw(1.5,-2.6)--(9,2.6);
	\draw(6,-2.6)--(13.5,2.6);

	\draw[blue](0,2.6)--(7.5,-2.6);
	\draw(1.5,2.6)--(9,-2.6);
	\draw(6,2.6)--(13.5,-2.6);

	\filldraw[color=black!60, fill=black!5, thin] (4.5,0.52) circle (0.17);
	\filldraw[color=black!60, fill=black!5, thin] (6.75,2.08) circle (0.17);

	\bluedot (4.125,0.26);
	\end{braid}.
$$
By (\ref{EPhiDifCol}), the first summand is $\Upsilon {\hat \ga^{i,j}}$ and, by (\ref{EDotPastCircle}) and (\ref{EPhiDifCol}), the last two summands are $(y_1 - y_{p+1})P^{(0)}{\hat \ga^{i,j}}$, proving the claim.
\end{proof}

Recall the notation $B_{1^d}^{(n)}$ from (\ref{EBNPar}).

\begin{Lemma} \label{LKeySameColors}
Let $i\in J$. Then 
$$
\ga^{i,i}\Upsilon\ga^{i,i}\equiv(-1)^i(z_1^2-z_2^2)\ga^{i,i}\pmod{B_{1^2}^{(1)}}.
$$ 
\end{Lemma}
\begin{proof}
We provide details for the generic case $1<i<\ell-1$, the special cases $i=0$, $i=1$ and $i=\ell-1$, being similar, are left to the reader. Applying Lemma~\ref{LA} and using the fact that a word starting with $\ell-1$ is not cuspidal by Lemma~\ref{LCuspExpl}, $\ga^{i,i}\Upsilon\ga^{i,i}$ equals  
$$
\resizebox{126mm}{23mm}{
\begin{braid}\tikzset{baseline=-.3em}
	\draw (0,6) node{\color{blue}$\ell$\color{black}};
	\braidbox{0.6}{2.9}{5.3}{6.6}{};
	\draw (1.8,6) node{$(\ell-1)^2$};
	\draw[dots] (3.5,6)--(5,6);
	\braidbox{5.3}{7.6}{5.3}{6.6}{};
	\draw (6.5,6) node{$(i+1)^2$};
	\draw (8.3,6) node{$i$};
	\draw (9.5,6) node{$i-1$};
	\draw[dots] (10.8,6)--(12.5,6);
	\draw (12.7,6) node{$1$};
	\redbraidbox{13.4}{14.3}{5.3}{6.6}{};
	\draw (13.8,6) node{$\color{red}0\,\,0\color{black}$};
	\draw (15,6) node{$1$};
	\draw[dots] (15.6,6)--(17,6);
	\draw (18.2,6) node{$i-1$};
	\draw (19.5,6) node{$i$};
	\draw (0,-6) node{\color{blue}$\ell$\color{black}};
	\braidbox{0.6}{2.9}{-6.7}{-5.4}{};
	\draw (1.8,-6) node{$(\ell-1)^2$};
	\draw[dots] (3.5,-6)--(5,-6);
	\braidbox{5.3}{7.6}{-6.7}{-5.4}{};
	\draw (6.5,-6) node{$(i+1)^2$};
	\draw (8.3,-6) node{$i$};
	\draw (9.5,-6) node{$i-1$};
	\draw[dots] (10.8,-6)--(12.5,-6);
	\draw (12.7,-6) node{$1$};
	\redbraidbox{13.4}{14.3}{-6.7}{-5.4}{};
	\draw (13.8,-6) node{$\color{red}0\,\,0\color{black}$};
	\draw (15,-6) node{$1$};
	\draw[dots] (15.6,-6)--(17,-6);
	\draw (18.2,-6) node{$i-1$};
	\draw (19.5,-6) node{$i$};

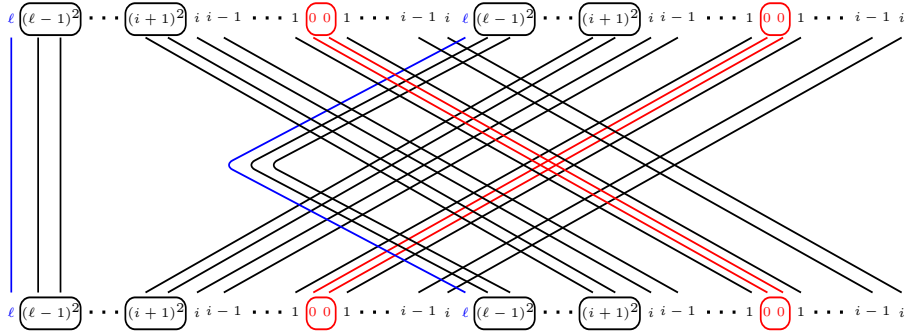
	\draw (20.3,6) node{\color{blue}$\ell$\color{black}};
	\braidbox{20.9}{23.2}{5.3}{6.6}{};
	\draw (22.1,6) node{$(\ell-1)^2$};
	\draw[dots] (23.8,6)--(25.3,6);
	\braidbox{25.6}{27.9}{5.3}{6.6}{};
	\draw (26.8,6) node{$(i+1)^2$};
	\draw (28.6,6) node{$i$};
	\draw (29.8,6) node{$i-1$};
	\draw[dots] (31.1,6)--(32.8,6);
	\draw (33,6) node{$1$};
	\redbraidbox{33.7}{34.6}{5.3}{6.6}{};
	\draw (34.1,6) node{$\color{red}0\,\,0\color{black}$};
	\draw (35.3,6) node{$1$};
	\draw[dots] (35.9,6)--(37.3,6);
	\draw (38.5,6) node{$i-1$};
	\draw (39.8,6) node{$i$};
	\draw (20.3,-6) node{\color{blue}$\ell$\color{black}};
	\braidbox{20.9}{23.2}{-6.7}{-5.4}{};
	\draw (22.1,-6) node{$(\ell-1)^2$};
	\draw[dots] (23.8,-6)--(25.1,-6);
	\braidbox{25.6}{27.9}{-6.7}{-5.4}{};
	\draw (26.8,-6) node{$(i+1)^2$};
	\draw (28.6,-6) node{$i$};
	\draw (29.8,-6) node{$i-1$};
	\draw[dots] (31.1,-6)--(32.8,-6);
	\draw (33,-6) node{$1$};
	\redbraidbox{33.7}{34.6}{-6.7}{-5.4}{};
	\draw (34.1,-6) node{$\color{red}0\,\,0\color{black}$};
	\draw (35.3,-6) node{$1$};
	\draw[dots] (35.9,-6)--(37.3,-6);
	\draw (38.5,-6) node{$i-1$};
	\draw (39.8,-6) node{$i$};
	\draw[blue](0,-5.2)--(0,5.2);
	\draw(1.2,-5.2)--(1.2,5.2);
	\draw(2.2,-5.2)--(2.2,5.2);
	\draw(6,-5.2)--(26,5.2);
	\draw(7,-5.2)--(27,5.2);
	\draw(8.3,-5.2)--(28.6,5.2);
	\draw(9.5,-5.2)--(29.8,5.2);
	\draw(12.7,-5.2)--(33.1,5.2);
	\draw[red](13.5,-5.2)--(33.8,5.2);
	\draw[red](14.2,-5.2)--(34.5,5.2);
	\draw(15,-5.2)--(35.3,5.2);
	\draw(18.2,-5.2)--(38.5,5.2);
	\draw(39.8,-5.2)--(19.5,5.2);
	\draw(39.8,5.2)--(19.5,-5.2);
	\draw[blue](20.3,-5.2)--(9.5,0)--(20.3,5.2);
	\draw(21.3,-5.2)--(10.5,0)--(21.3,5.2);
	\draw(22.3,-5.2)--(11.5,0)--(22.3,5.2);
	\draw(26,-5.2)--(6,5.2);
	\draw(27,-5.2)--(7,5.2);
	\draw(28.6,-5.2)--(8.3,5.2);
	\draw(29.8,-5.2)--(9.5,5.2);
	\draw(38.5,-5.2)--(18.2,5.2);
	\draw(35.3,-5.2)--(15,5.2);
	\draw[red](34.5,-5.2)--(14.2,5.2);
	\draw[red](33.8,-5.2)--(13.5,5.2);
	\draw(33.1,-5.2)--(12.7,5.2);
	\end{braid}.
	}
$$
Repeatedly applying Lemma~\ref{LA1/2} and using the fact that a word starting with $\ell(\ell-1)^2\dots(j+2)^2j$ is not cuspidal by Lemma~\ref{LCuspExpl}, the above diagram equals 
$$
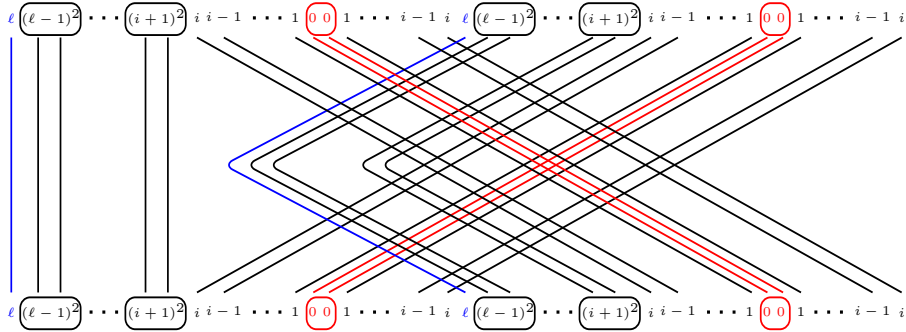
\resizebox{126mm}{23mm}{
\begin{braid}\tikzset{baseline=-.3em}
	\draw (0,6) node{\color{blue}$\ell$\color{black}};
	\braidbox{0.6}{2.9}{5.3}{6.6}{};
	\draw (1.8,6) node{$(\ell-1)^2$};
	\draw[dots] (3.5,6)--(5,6);
	\braidbox{5.3}{7.6}{5.3}{6.6}{};
	\draw (6.5,6) node{$(i+1)^2$};
	\draw (8.3,6) node{$i$};
	\draw (9.5,6) node{$i-1$};
	\draw[dots] (10.8,6)--(12.5,6);
	\draw (12.7,6) node{$1$};
	\redbraidbox{13.4}{14.3}{5.3}{6.6}{};
	\draw (13.8,6) node{$\color{red}0\,\,0\color{black}$};
	\draw (15,6) node{$1$};
	\draw[dots] (15.6,6)--(17,6);
	\draw (18.2,6) node{$i-1$};
	\draw (19.5,6) node{$i$};
	\draw (0,-6) node{\color{blue}$\ell$\color{black}};
	\braidbox{0.6}{2.9}{-6.7}{-5.4}{};
	\draw (1.8,-6) node{$(\ell-1)^2$};
	\draw[dots] (3.5,-6)--(5,-6);
	\braidbox{5.3}{7.6}{-6.7}{-5.4}{};
	\draw (6.5,-6) node{$(i+1)^2$};
	\draw (8.3,-6) node{$i$};
	\draw (9.5,-6) node{$i-1$};
	\draw[dots] (10.8,-6)--(12.5,-6);
	\draw (12.7,-6) node{$1$};
	\redbraidbox{13.4}{14.3}{-6.7}{-5.4}{};
	\draw (13.8,-6) node{$\color{red}0\,\,0\color{black}$};
	\draw (15,-6) node{$1$};
	\draw[dots] (15.6,-6)--(17,-6);
	\draw (18.2,-6) node{$i-1$};
	\draw (19.5,-6) node{$i$};

	\draw (20.3,6) node{\color{blue}$\ell$\color{black}};
	\braidbox{20.9}{23.2}{5.3}{6.6}{};
	\draw (22.1,6) node{$(\ell-1)^2$};
	\draw[dots] (23.8,6)--(25.3,6);
	\braidbox{25.6}{27.9}{5.3}{6.6}{};
	\draw (26.8,6) node{$(i+1)^2$};
	\draw (28.6,6) node{$i$};
	\draw (29.8,6) node{$i-1$};
	\draw[dots] (31.1,6)--(32.8,6);
	\draw (33,6) node{$1$};
	\redbraidbox{33.7}{34.6}{5.3}{6.6}{};
	\draw (34.1,6) node{$\color{red}0\,\,0\color{black}$};
	\draw (35.3,6) node{$1$};
	\draw[dots] (35.9,6)--(37.3,6);
	\draw (38.5,6) node{$i-1$};
	\draw (39.8,6) node{$i$};
	\draw (20.3,-6) node{\color{blue}$\ell$\color{black}};
	\braidbox{20.9}{23.2}{-6.7}{-5.4}{};
	\draw (22.1,-6) node{$(\ell-1)^2$};
	\draw[dots] (23.8,-6)--(25.1,-6);
	\braidbox{25.6}{27.9}{-6.7}{-5.4}{};
	\draw (26.8,-6) node{$(i+1)^2$};
	\draw (28.6,-6) node{$i$};
	\draw (29.8,-6) node{$i-1$};
	\draw[dots] (31.1,-6)--(32.8,-6);
	\draw (33,-6) node{$1$};
	\redbraidbox{33.7}{34.6}{-6.7}{-5.4}{};
	\draw (34.1,-6) node{$\color{red}0\,\,0\color{black}$};
	\draw (35.3,-6) node{$1$};
	\draw[dots] (35.9,-6)--(37.3,-6);
	\draw (38.5,-6) node{$i-1$};
	\draw (39.8,-6) node{$i$};
	\draw[blue](0,-5.2)--(0,5.2);
	\draw(1.2,-5.2)--(1.2,5.2);
	\draw(2.2,-5.2)--(2.2,5.2);
	\draw(6,-5.2)--(6,5.2);
	\draw(7,-5.2)--(7,5.2);
	\draw(8.3,-5.2)--(28.6,5.2);
	\draw(9.5,-5.2)--(29.8,5.2);
	\draw(12.7,-5.2)--(33.1,5.2);
	\draw[red](13.5,-5.2)--(33.8,5.2);
	\draw[red](14.2,-5.2)--(34.5,5.2);
	\draw(15,-5.2)--(35.3,5.2);
	\draw(18.2,-5.2)--(38.5,5.2);
	\draw(39.8,-5.2)--(19.5,5.2);
	\draw(39.8,5.2)--(19.5,-5.2);
	\draw[blue](20.3,-5.2)--(9.5,0)--(20.3,5.2);
	\draw(21.3,-5.2)--(10.5,0)--(21.3,5.2);
	\draw(22.3,-5.2)--(11.5,0)--(22.3,5.2);
	\draw(26,-5.2)--(15.5,0)--(26,5.2);
	\draw(27,-5.2)--(16.5,0)--(27,5.2);
	\draw(28.6,-5.2)--(8.3,5.2);
	\draw(29.8,-5.2)--(9.5,5.2);
	\draw(38.5,-5.2)--(18.2,5.2);
	\draw(35.3,-5.2)--(15,5.2);
	\draw[red](34.5,-5.2)--(14.2,5.2);
	\draw[red](33.8,-5.2)--(13.5,5.2);
	\draw(33.1,-5.2)--(12.7,5.2);
	\end{braid}.
	}
$$
We apply Lemma~\ref{LA3/4} to write this diagram as a sum of five terms $I+II+III+IV+V$ corresponding to the five summands in the right hand side of the lemma. 

First we check that $I=0$. Indeed, 
$$
I=
\resizebox{119mm}{23mm}
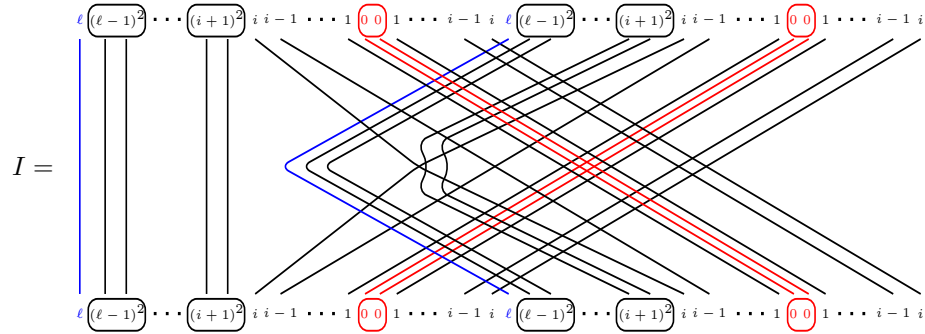
{
\begin{braid}\tikzset{baseline=-.3em}
	\draw (0,6) node{\color{blue}$\ell$\color{black}};
	\braidbox{0.6}{2.9}{5.3}{6.6}{};
	\draw (1.8,6) node{$(\ell-1)^2$};
	\draw[dots] (3.5,6)--(5,6);
	\braidbox{5.3}{7.6}{5.3}{6.6}{};
	\draw (6.5,6) node{$(i+1)^2$};
	\draw (8.3,6) node{$i$};
	\draw (9.5,6) node{$i-1$};
	\draw[dots] (10.8,6)--(12.5,6);
	\draw (12.7,6) node{$1$};
	\redbraidbox{13.4}{14.3}{5.3}{6.6}{};
	\draw (13.8,6) node{$\color{red}0\,\,0\color{black}$};
	\draw (15,6) node{$1$};
	\draw[dots] (15.6,6)--(17,6);
	\draw (18.2,6) node{$i-1$};
	\draw (19.5,6) node{$i$};
	\draw (0,-6) node{\color{blue}$\ell$\color{black}};
	\braidbox{0.6}{2.9}{-6.7}{-5.4}{};
	\draw (1.8,-6) node{$(\ell-1)^2$};
	\draw[dots] (3.5,-6)--(5,-6);
	\braidbox{5.3}{7.6}{-6.7}{-5.4}{};
	\draw (6.5,-6) node{$(i+1)^2$};
	\draw (8.3,-6) node{$i$};
	\draw (9.5,-6) node{$i-1$};
	\draw[dots] (10.8,-6)--(12.5,-6);
	\draw (12.7,-6) node{$1$};
	\redbraidbox{13.4}{14.3}{-6.7}{-5.4}{};
	\draw (13.8,-6) node{$\color{red}0\,\,0\color{black}$};
	\draw (15,-6) node{$1$};
	\draw[dots] (15.6,-6)--(17,-6);
	\draw (18.2,-6) node{$i-1$};
	\draw (19.5,-6) node{$i$};

	\draw (20.3,6) node{\color{blue}$\ell$\color{black}};
	\braidbox{20.9}{23.2}{5.3}{6.6}{};
	\draw (22.1,6) node{$(\ell-1)^2$};
	\draw[dots] (23.8,6)--(25.3,6);
	\braidbox{25.6}{27.9}{5.3}{6.6}{};
	\draw (26.8,6) node{$(i+1)^2$};
	\draw (28.6,6) node{$i$};
	\draw (29.8,6) node{$i-1$};
	\draw[dots] (31.1,6)--(32.8,6);
	\draw (33,6) node{$1$};
	\redbraidbox{33.7}{34.6}{5.3}{6.6}{};
	\draw (34.1,6) node{$\color{red}0\,\,0\color{black}$};
	\draw (35.3,6) node{$1$};
	\draw[dots] (35.9,6)--(37.3,6);
	\draw (38.5,6) node{$i-1$};
	\draw (39.8,6) node{$i$};
	\draw (20.3,-6) node{\color{blue}$\ell$\color{black}};
	\braidbox{20.9}{23.2}{-6.7}{-5.4}{};
	\draw (22.1,-6) node{$(\ell-1)^2$};
	\draw[dots] (23.8,-6)--(25.1,-6);
	\braidbox{25.6}{27.9}{-6.7}{-5.4}{};
	\draw (26.8,-6) node{$(i+1)^2$};
	\draw (28.6,-6) node{$i$};
	\draw (29.8,-6) node{$i-1$};
	\draw[dots] (31.1,-6)--(32.8,-6);
	\draw (33,-6) node{$1$};
	\redbraidbox{33.7}{34.6}{-6.7}{-5.4}{};
	\draw (34.1,-6) node{$\color{red}0\,\,0\color{black}$};
	\draw (35.3,-6) node{$1$};
	\draw[dots] (35.9,-6)--(37.3,-6);
	\draw (38.5,-6) node{$i-1$};
	\draw (39.8,-6) node{$i$};
	\draw[blue](0,-5.2)--(0,5.2);
	\draw(1.2,-5.2)--(1.2,5.2);
	\draw(2.2,-5.2)--(2.2,5.2);
	\draw(6,-5.2)--(6,5.2);
	\draw(7,-5.2)--(7,5.2);
	\draw(8.3,-5.2)--(16,0)--(28.6,5.2);
	\draw(9.5,-5.2)--(29.8,5.2);
	\draw(12.7,-5.2)--(33.1,5.2);
	\draw[red](13.5,-5.2)--(33.8,5.2);
	\draw[red](14.2,-5.2)--(34.5,5.2);
	\draw(15,-5.2)--(35.3,5.2);
	\draw(18.2,-5.2)--(38.5,5.2);
	\draw(39.8,-5.2)--(19.5,5.2);
	\draw(39.8,5.2)--(19.5,-5.2);
	\draw[blue](20.3,-5.2)--(9.5,0)--(20.3,5.2);
	\draw(21.3,-5.2)--(10.5,0)--(21.3,5.2);
	\draw(22.3,-5.2)--(11.5,0)--(22.3,5.2);
	\draw(26,-5.2)--(16,-1)--(16.5,0)--(16,1)--(26,5.2);
	\draw(27,-5.2)--(17,-1)--(17.5,0)--(17,1)--(27,5.2);
	\draw(28.6,-5.2)--(16,0)--(8.3,5.2);
	\draw(29.8,-5.2)--(9.5,5.2);
	\draw(38.5,-5.2)--(18.2,5.2);
	\draw(35.3,-5.2)--(15,5.2);
	\draw[red](34.5,-5.2)--(14.2,5.2);
	\draw[red](33.8,-5.2)--(13.5,5.2);
	\draw(33.1,-5.2)--(12.7,5.2);
	\end{braid}.
	}
$$
Using braid relations, we get
$$
I=
\resizebox{119mm}{20mm}
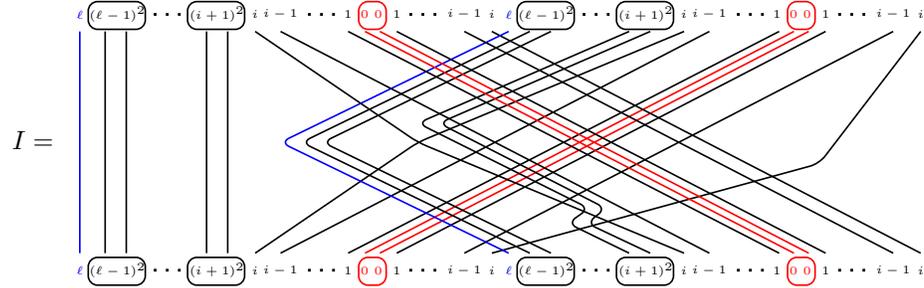
{
\begin{braid}\tikzset{baseline=-.3em}
	\draw (0,6) node{\color{blue}$\ell$\color{black}};
	\braidbox{0.6}{2.9}{5.3}{6.6}{};
	\draw (1.8,6) node{$(\ell-1)^2$};
	\draw[dots] (3.5,6)--(5,6);
	\braidbox{5.3}{7.6}{5.3}{6.6}{};
	\draw (6.5,6) node{$(i+1)^2$};
	\draw (8.3,6) node{$i$};
	\draw (9.5,6) node{$i-1$};
	\draw[dots] (10.8,6)--(12.5,6);
	\draw (12.7,6) node{$1$};
	\redbraidbox{13.4}{14.3}{5.3}{6.6}{};
	\draw (13.8,6) node{$\color{red}0\,\,0\color{black}$};
	\draw (15,6) node{$1$};
	\draw[dots] (15.6,6)--(17,6);
	\draw (18.2,6) node{$i-1$};
	\draw (19.5,6) node{$i$};
	\draw (0,-6) node{\color{blue}$\ell$\color{black}};
	\braidbox{0.6}{2.9}{-6.7}{-5.4}{};
	\draw (1.8,-6) node{$(\ell-1)^2$};
	\draw[dots] (3.5,-6)--(5,-6);
	\braidbox{5.3}{7.6}{-6.7}{-5.4}{};
	\draw (6.5,-6) node{$(i+1)^2$};
	\draw (8.3,-6) node{$i$};
	\draw (9.5,-6) node{$i-1$};
	\draw[dots] (10.8,-6)--(12.5,-6);
	\draw (12.7,-6) node{$1$};
	\redbraidbox{13.4}{14.3}{-6.7}{-5.4}{};
	\draw (13.8,-6) node{$\color{red}0\,\,0\color{black}$};
	\draw (15,-6) node{$1$};
	\draw[dots] (15.6,-6)--(17,-6);
	\draw (18.2,-6) node{$i-1$};
	\draw (19.5,-6) node{$i$};

	\draw (20.3,6) node{\color{blue}$\ell$\color{black}};
	\braidbox{20.9}{23.2}{5.3}{6.6}{};
	\draw (22.1,6) node{$(\ell-1)^2$};
	\draw[dots] (23.8,6)--(25.3,6);
	\braidbox{25.6}{27.9}{5.3}{6.6}{};
	\draw (26.8,6) node{$(i+1)^2$};
	\draw (28.6,6) node{$i$};
	\draw (29.8,6) node{$i-1$};
	\draw[dots] (31.1,6)--(32.8,6);
	\draw (33,6) node{$1$};
	\redbraidbox{33.7}{34.6}{5.3}{6.6}{};
	\draw (34.1,6) node{$\color{red}0\,\,0\color{black}$};
	\draw (35.3,6) node{$1$};
	\draw[dots] (35.9,6)--(37.3,6);
	\draw (38.5,6) node{$i-1$};
	\draw (39.8,6) node{$i$};
	\draw (20.3,-6) node{\color{blue}$\ell$\color{black}};
	\braidbox{20.9}{23.2}{-6.7}{-5.4}{};
	\draw (22.1,-6) node{$(\ell-1)^2$};
	\draw[dots] (23.8,-6)--(25.1,-6);
	\braidbox{25.6}{27.9}{-6.7}{-5.4}{};
	\draw (26.8,-6) node{$(i+1)^2$};
	\draw (28.6,-6) node{$i$};
	\draw (29.8,-6) node{$i-1$};
	\draw[dots] (31.1,-6)--(32.8,-6);
	\draw (33,-6) node{$1$};
	\redbraidbox{33.7}{34.6}{-6.7}{-5.4}{};
	\draw (34.1,-6) node{$\color{red}0\,\,0\color{black}$};
	\draw (35.3,-6) node{$1$};
	\draw[dots] (35.9,-6)--(37.3,-6);
	\draw (38.5,-6) node{$i-1$};
	\draw (39.8,-6) node{$i$};
	\draw[blue](0,-5.2)--(0,5.2);
	\draw(1.2,-5.2)--(1.2,5.2);
	\draw(2.2,-5.2)--(2.2,5.2);
	\draw(6,-5.2)--(6,5.2);
	\draw(7,-5.2)--(7,5.2);
	\draw(8.3,-5.2)--(16,0)--(28.6,5.2);
	\draw(9.5,-5.2)--(29.8,5.2);
	\draw(12.7,-5.2)--(33.1,5.2);
	\draw[red](13.5,-5.2)--(33.8,5.2);
	\draw[red](14.2,-5.2)--(34.5,5.2);
	\draw(15,-5.2)--(35.3,5.2);
	\draw(18.2,-5.2)--(38.5,5.2);
	\draw(39.8,-5.2)--(19.5,5.2);
	\draw(39.8,5.2)--(35,-1)--(19.5,-5.2);
	\draw[blue](20.3,-5.2)--(9.5,0)--(20.3,5.2);
	\draw(21.3,-5.2)--(10.5,0)--(21.3,5.2);
	\draw(22.3,-5.2)--(11.5,0)--(22.3,5.2);
	\draw(26,-5.2)--(23.1,-3.8)--(24.1,-3.1)--(16,0.9)--(26,5.2);
	\draw(27,-5.2)--(24,-3.8)--(24.9,-3.1)--(17,0.9)--(27,5.2);
	\draw(28.6,-5.2)--(16,0)--(8.3,5.2);
	\draw(29.8,-5.2)--(9.5,5.2);
	\draw(38.5,-5.2)--(18.2,5.2);
	\draw(35.3,-5.2)--(15,5.2);
	\draw[red](34.5,-5.2)--(14.2,5.2);
	\draw[red](33.8,-5.2)--(13.5,5.2);
	\draw(33.1,-5.2)--(12.7,5.2);
	\end{braid}.
	}
$$
Applying Lemma~\ref{LA3/4} again, we write $I$ as a sum of five terms, all of which are zero since words starting with $\ell(\ell-1)^2\cdots (i+2)^2i$ and $\ell(\ell-1)^2\cdots (i+1)^2(i-1)$ are not cuspidal by Lemma~\ref{LCuspExpl}. 

We next consider the fifth summand 
$$
V=-
\resizebox{116mm}{20mm}
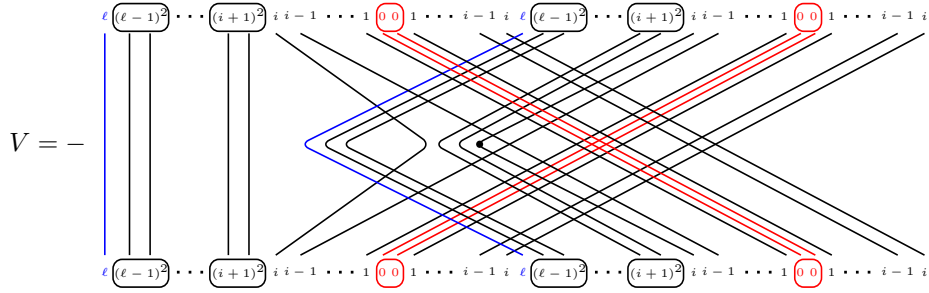
{
\begin{braid}\tikzset{baseline=-.3em}
	\draw (0,6) node{\color{blue}$\ell$\color{black}};
	\braidbox{0.6}{2.9}{5.3}{6.6}{};
	\draw (1.8,6) node{$(\ell-1)^2$};
	\draw[dots] (3.5,6)--(5,6);
	\braidbox{5.3}{7.6}{5.3}{6.6}{};
	\draw (6.5,6) node{$(i+1)^2$};
	\draw (8.3,6) node{$i$};
	\draw (9.5,6) node{$i-1$};
	\draw[dots] (10.8,6)--(12.5,6);
	\draw (12.7,6) node{$1$};
	\redbraidbox{13.4}{14.3}{5.3}{6.6}{};
	\draw (13.8,6) node{$\color{red}0\,\,0\color{black}$};
	\draw (15,6) node{$1$};
	\draw[dots] (15.6,6)--(17,6);
	\draw (18.2,6) node{$i-1$};
	\draw (19.5,6) node{$i$};
	\draw (0,-6) node{\color{blue}$\ell$\color{black}};
	\braidbox{0.6}{2.9}{-6.7}{-5.4}{};
	\draw (1.8,-6) node{$(\ell-1)^2$};
	\draw[dots] (3.5,-6)--(5,-6);
	\braidbox{5.3}{7.6}{-6.7}{-5.4}{};
	\draw (6.5,-6) node{$(i+1)^2$};
	\draw (8.3,-6) node{$i$};
	\draw (9.5,-6) node{$i-1$};
	\draw[dots] (10.8,-6)--(12.5,-6);
	\draw (12.7,-6) node{$1$};
	\redbraidbox{13.4}{14.3}{-6.7}{-5.4}{};
	\draw (13.8,-6) node{$\color{red}0\,\,0\color{black}$};
	\draw (15,-6) node{$1$};
	\draw[dots] (15.6,-6)--(17,-6);
	\draw (18.2,-6) node{$i-1$};
	\draw (19.5,-6) node{$i$};
	\draw (20.3,6) node{\color{blue}$\ell$\color{black}};
	\braidbox{20.9}{23.2}{5.3}{6.6}{};
	\draw (22.1,6) node{$(\ell-1)^2$};
	\draw[dots] (23.8,6)--(25.3,6);
	\braidbox{25.6}{27.9}{5.3}{6.6}{};
	\draw (26.8,6) node{$(i+1)^2$};
	\draw (28.6,6) node{$i$};
	\draw (29.8,6) node{$i-1$};
	\draw[dots] (31.1,6)--(32.8,6);
	\draw (33,6) node{$1$};
	\redbraidbox{33.7}{34.6}{5.3}{6.6}{};
	\draw (34.1,6) node{$\color{red}0\,\,0\color{black}$};
	\draw (35.3,6) node{$1$};
	\draw[dots] (35.9,6)--(37.3,6);
	\draw (38.5,6) node{$i-1$};
	\draw (39.8,6) node{$i$};
	\draw (20.3,-6) node{\color{blue}$\ell$\color{black}};
	\braidbox{20.9}{23.2}{-6.7}{-5.4}{};
	\draw (22.1,-6) node{$(\ell-1)^2$};
	\draw[dots] (23.8,-6)--(25.1,-6);
	\braidbox{25.6}{27.9}{-6.7}{-5.4}{};
	\draw (26.8,-6) node{$(i+1)^2$};
	\draw (28.6,-6) node{$i$};
	\draw (29.8,-6) node{$i-1$};
	\draw[dots] (31.1,-6)--(32.8,-6);
	\draw (33,-6) node{$1$};
	\redbraidbox{33.7}{34.6}{-6.7}{-5.4}{};
	\draw (34.1,-6) node{$\color{red}0\,\,0\color{black}$};
	\draw (35.3,-6) node{$1$};
	\draw[dots] (35.9,-6)--(37.3,-6);
	\draw (38.5,-6) node{$i-1$};
	\draw (39.8,-6) node{$i$};
	\draw[blue](0,-5.2)--(0,5.2);
	\draw(1.2,-5.2)--(1.2,5.2);
	\draw(2.2,-5.2)--(2.2,5.2);
	\draw(6,-5.2)--(6,5.2);
	\draw(7,-5.2)--(7,5.2);
	\draw(8.3,-5.2)--(15.8,0)--(8.3,5.2);
	\draw(9.5,-5.2)--(29.8,5.2);
	\draw(12.7,-5.2)--(33.1,5.2);
	\draw[red](13.5,-5.2)--(33.8,5.2);
	\draw[red](14.2,-5.2)--(34.5,5.2);
	\draw(15,-5.2)--(35.3,5.2);
	\draw(18.2,-5.2)--(38.5,5.2);
	\draw(39.8,-5.2)--(19.5,5.2);
	\draw(39.8,5.2)--(19.5,-5.2);
	\draw[blue](20.3,-5.2)--(9.5,0)--(20.3,5.2);
	\draw(21.3,-5.2)--(10.5,0)--(21.3,5.2);
	\draw(22.3,-5.2)--(11.5,0)--(22.3,5.2);
	\draw(26,-5.2)--(16,0)--(26,5.2);
	\draw(27,-5.2)--(17,0)--(27,5.2);
	\draw(28.6,-5.2)--(18,0)--(28.6,5.2);
	\draw(29.8,-5.2)--(9.5,5.2);
	\draw(38.5,-5.2)--(18.2,5.2);
	\draw(35.3,-5.2)--(15,5.2);
	\draw[red](34.5,-5.2)--(14.2,5.2);
	\draw[red](33.8,-5.2)--(13.5,5.2);
	\draw(33.1,-5.2)--(12.7,5.2);
	\blackdot (18.2,0);
	\end{braid}.
	}
$$
We move the dot down its string using the dot-crossing relations (\ref{R5}). The only error term arising from the $(i,i)$-crossing is zero since a word starting with $\ell(\ell-1)^2\cdots(i+1)^2(i-1)$ is not cuspidal by Lemma~\ref{LCuspExpl}. Moreover, in view of (\ref{ERelB3}), when the dot reaches the bottom it becomes equal to $\ga^{i,i}z_2\ga^{i,i}$. Thus $V=-X\ga^{i,i}z_2\ga^{i,i}$, where 
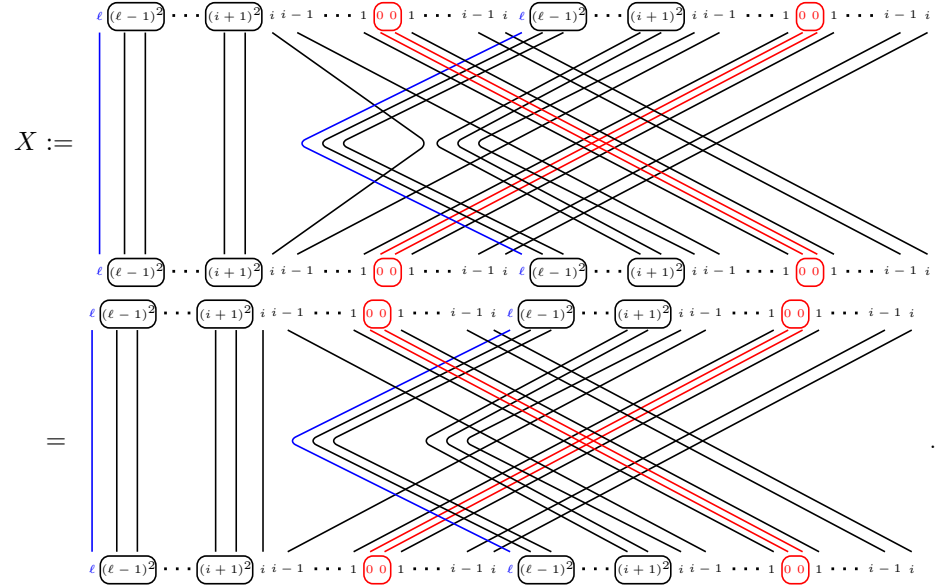
\begin{align*}
X&:=
\resizebox{116mm}{20mm}{
\begin{braid}\tikzset{baseline=-.3em}
	\draw (0,6) node{\color{blue}$\ell$\color{black}};
	\braidbox{0.6}{2.9}{5.3}{6.6}{};
	\draw (1.8,6) node{$(\ell-1)^2$};
	\draw[dots] (3.5,6)--(5,6);
	\braidbox{5.3}{7.6}{5.3}{6.6}{};
	\draw (6.5,6) node{$(i+1)^2$};
	\draw (8.3,6) node{$i$};
	\draw (9.5,6) node{$i-1$};
	\draw[dots] (10.8,6)--(12.5,6);
	\draw (12.7,6) node{$1$};
	\redbraidbox{13.4}{14.3}{5.3}{6.6}{};
	\draw (13.8,6) node{$\color{red}0\,\,0\color{black}$};
	\draw (15,6) node{$1$};
	\draw[dots] (15.6,6)--(17,6);
	\draw (18.2,6) node{$i-1$};
	\draw (19.5,6) node{$i$};
	\draw (0,-6) node{\color{blue}$\ell$\color{black}};
	\braidbox{0.6}{2.9}{-6.7}{-5.4}{};
	\draw (1.8,-6) node{$(\ell-1)^2$};
	\draw[dots] (3.5,-6)--(5,-6);
	\braidbox{5.3}{7.6}{-6.7}{-5.4}{};
	\draw (6.5,-6) node{$(i+1)^2$};
	\draw (8.3,-6) node{$i$};
	\draw (9.5,-6) node{$i-1$};
	\draw[dots] (10.8,-6)--(12.5,-6);
	\draw (12.7,-6) node{$1$};
	\redbraidbox{13.4}{14.3}{-6.7}{-5.4}{};
	\draw (13.8,-6) node{$\color{red}0\,\,0\color{black}$};
	\draw (15,-6) node{$1$};
	\draw[dots] (15.6,-6)--(17,-6);
	\draw (18.2,-6) node{$i-1$};
	\draw (19.5,-6) node{$i$};
	\draw (20.3,6) node{\color{blue}$\ell$\color{black}};
	\braidbox{20.9}{23.2}{5.3}{6.6}{};
	\draw (22.1,6) node{$(\ell-1)^2$};
	\draw[dots] (23.8,6)--(25.3,6);
	\braidbox{25.6}{27.9}{5.3}{6.6}{};
	\draw (26.8,6) node{$(i+1)^2$};
	\draw (28.6,6) node{$i$};
	\draw (29.8,6) node{$i-1$};
	\draw[dots] (31.1,6)--(32.8,6);
	\draw (33,6) node{$1$};
	\redbraidbox{33.7}{34.6}{5.3}{6.6}{};
	\draw (34.1,6) node{$\color{red}0\,\,0\color{black}$};
	\draw (35.3,6) node{$1$};
	\draw[dots] (35.9,6)--(37.3,6);
	\draw (38.5,6) node{$i-1$};
	\draw (39.8,6) node{$i$};
	\draw (20.3,-6) node{\color{blue}$\ell$\color{black}};
	\braidbox{20.9}{23.2}{-6.7}{-5.4}{};
	\draw (22.1,-6) node{$(\ell-1)^2$};
	\draw[dots] (23.8,-6)--(25.1,-6);
	\braidbox{25.6}{27.9}{-6.7}{-5.4}{};
	\draw (26.8,-6) node{$(i+1)^2$};
	\draw (28.6,-6) node{$i$};
	\draw (29.8,-6) node{$i-1$};
	\draw[dots] (31.1,-6)--(32.8,-6);
	\draw (33,-6) node{$1$};
	\redbraidbox{33.7}{34.6}{-6.7}{-5.4}{};
	\draw (34.1,-6) node{$\color{red}0\,\,0\color{black}$};
	\draw (35.3,-6) node{$1$};
	\draw[dots] (35.9,-6)--(37.3,-6);
	\draw (38.5,-6) node{$i-1$};
	\draw (39.8,-6) node{$i$};
	\draw[blue](0,-5.2)--(0,5.2);
	\draw(1.2,-5.2)--(1.2,5.2);
	\draw(2.2,-5.2)--(2.2,5.2);
	\draw(6,-5.2)--(6,5.2);
	\draw(7,-5.2)--(7,5.2);
	\draw(8.3,-5.2)--(15.8,0)--(8.3,5.2);
	\draw(9.5,-5.2)--(29.8,5.2);
	\draw(12.7,-5.2)--(33.1,5.2);
	\draw[red](13.5,-5.2)--(33.8,5.2);
	\draw[red](14.2,-5.2)--(34.5,5.2);
	\draw(15,-5.2)--(35.3,5.2);
	\draw(18.2,-5.2)--(38.5,5.2);
	\draw(39.8,-5.2)--(19.5,5.2);
	\draw(39.8,5.2)--(19.5,-5.2);
	\draw[blue](20.3,-5.2)--(9.5,0)--(20.3,5.2);
	\draw(21.3,-5.2)--(10.5,0)--(21.3,5.2);
	\draw(22.3,-5.2)--(11.5,0)--(22.3,5.2);
	\draw(26,-5.2)--(16,0)--(26,5.2);
	\draw(27,-5.2)--(17,0)--(27,5.2);
	\draw(28.6,-5.2)--(18,0)--(28.6,5.2);
	\draw(29.8,-5.2)--(9.5,5.2);
	\draw(38.5,-5.2)--(18.2,5.2);
	\draw(35.3,-5.2)--(15,5.2);
	\draw[red](34.5,-5.2)--(14.2,5.2);
	\draw[red](33.8,-5.2)--(13.5,5.2);
	\draw(33.1,-5.2)--(12.7,5.2);
	\end{braid}
	}
	\\
	&=
	\resizebox{116mm}{20mm}{
	\begin{braid}\tikzset{baseline=-.3em}
	\draw (0,6) node{\color{blue}$\ell$\color{black}};
	\braidbox{0.6}{2.9}{5.3}{6.6}{};
	\draw (1.8,6) node{$(\ell-1)^2$};
	\draw[dots] (3.5,6)--(5,6);
	\braidbox{5.3}{7.6}{5.3}{6.6}{};
	\draw (6.5,6) node{$(i+1)^2$};
	\draw (8.3,6) node{$i$};
	\draw (9.5,6) node{$i-1$};
	\draw[dots] (10.8,6)--(12.5,6);
	\draw (12.7,6) node{$1$};
	\redbraidbox{13.4}{14.3}{5.3}{6.6}{};
	\draw (13.8,6) node{$\color{red}0\,\,0\color{black}$};
	\draw (15,6) node{$1$};
	\draw[dots] (15.6,6)--(17,6);
	\draw (18.2,6) node{$i-1$};
	\draw (19.5,6) node{$i$};
	\draw (0,-6) node{\color{blue}$\ell$\color{black}};
	\braidbox{0.6}{2.9}{-6.7}{-5.4}{};
	\draw (1.8,-6) node{$(\ell-1)^2$};
	\draw[dots] (3.5,-6)--(5,-6);
	\braidbox{5.3}{7.6}{-6.7}{-5.4}{};
	\draw (6.5,-6) node{$(i+1)^2$};
	\draw (8.3,-6) node{$i$};
	\draw (9.5,-6) node{$i-1$};
	\draw[dots] (10.8,-6)--(12.5,-6);
	\draw (12.7,-6) node{$1$};
	\redbraidbox{13.4}{14.3}{-6.7}{-5.4}{};
	\draw (13.8,-6) node{$\color{red}0\,\,0\color{black}$};
	\draw (15,-6) node{$1$};
	\draw[dots] (15.6,-6)--(17,-6);
	\draw (18.2,-6) node{$i-1$};
	\draw (19.5,-6) node{$i$};
	\draw (20.3,6) node{\color{blue}$\ell$\color{black}};
	\braidbox{20.9}{23.2}{5.3}{6.6}{};
	\draw (22.1,6) node{$(\ell-1)^2$};
	\draw[dots] (23.8,6)--(25.3,6);
	\braidbox{25.6}{27.9}{5.3}{6.6}{};
	\draw (26.8,6) node{$(i+1)^2$};
	\draw (28.6,6) node{$i$};
	\draw (29.8,6) node{$i-1$};
	\draw[dots] (31.1,6)--(32.8,6);
	\draw (33,6) node{$1$};
	\redbraidbox{33.7}{34.6}{5.3}{6.6}{};
	\draw (34.1,6) node{$\color{red}0\,\,0\color{black}$};
	\draw (35.3,6) node{$1$};
	\draw[dots] (35.9,6)--(37.3,6);
	\draw (38.5,6) node{$i-1$};
	\draw (39.8,6) node{$i$};
	\draw (20.3,-6) node{\color{blue}$\ell$\color{black}};
	\braidbox{20.9}{23.2}{-6.7}{-5.4}{};
	\draw (22.1,-6) node{$(\ell-1)^2$};
	\draw[dots] (23.8,-6)--(25.1,-6);
	\braidbox{25.6}{27.9}{-6.7}{-5.4}{};
	\draw (26.8,-6) node{$(i+1)^2$};
	\draw (28.6,-6) node{$i$};
	\draw (29.8,-6) node{$i-1$};
	\draw[dots] (31.1,-6)--(32.8,-6);
	\draw (33,-6) node{$1$};
	\redbraidbox{33.7}{34.6}{-6.7}{-5.4}{};
	\draw (34.1,-6) node{$\color{red}0\,\,0\color{black}$};
	\draw (35.3,-6) node{$1$};
	\draw[dots] (35.9,-6)--(37.3,-6);
	\draw (38.5,-6) node{$i-1$};
	\draw (39.8,-6) node{$i$};
	\draw[blue](0,-5.2)--(0,5.2);
	\draw(1.2,-5.2)--(1.2,5.2);
	\draw(2.2,-5.2)--(2.2,5.2);
	\draw(6,-5.2)--(6,5.2);
	\draw(7,-5.2)--(7,5.2);
	\draw(8.3,-5.2)--(8.3,5.2);
	\draw(9.5,-5.2)--(29.8,5.2);
	\draw(12.7,-5.2)--(33.1,5.2);
	\draw[red](13.5,-5.2)--(33.8,5.2);
	\draw[red](14.2,-5.2)--(34.5,5.2);
	\draw(15,-5.2)--(35.3,5.2);
	\draw(18.2,-5.2)--(38.5,5.2);
	\draw(39.8,-5.2)--(19.5,5.2);
	\draw(39.8,5.2)--(19.5,-5.2);
	\draw[blue](20.3,-5.2)--(9.5,0)--(20.3,5.2);
	\draw(21.3,-5.2)--(10.5,0)--(21.3,5.2);
	\draw(22.3,-5.2)--(11.5,0)--(22.3,5.2);
	\draw(26,-5.2)--(16,0)--(26,5.2);
	\draw(27,-5.2)--(17,0)--(27,5.2);
	\draw(28.6,-5.2)--(18,0)--(28.6,5.2);
	\draw(29.8,-5.2)--(9.5,5.2);
	\draw(38.5,-5.2)--(18.2,5.2);
	\draw(35.3,-5.2)--(15,5.2);
	\draw[red](34.5,-5.2)--(14.2,5.2);
	\draw[red](33.8,-5.2)--(13.5,5.2);
	\draw(33.1,-5.2)--(12.7,5.2);
	\end{braid}.
	}
\end{align*}
Similarly, $III+IV=0$ thanks to (\ref{ERelB2}), and $II=-X\ga^{i,i}z_1\ga^{i,i}$. Thus,
\begin{equation}\label{E160721}
\ga^{i,i}\Upsilon\ga^{i,i}=I+II+III+IV+V=X\ga^{i,i}(-z_1-z_2)\ga^{i,i}.
\end{equation}
Using the braid relation for 
$\begin{braid}\tikzset{baseline=.4em}
	\draw (0,0) node{$i-1$};
        \draw (1.1,0) node{$i$};
        \draw (2.2,0) node{$i-1$};
		\draw(0,0.5)--(2.2,1.5);
\draw(2.2,0.5)--(0,1.5);
\draw(1.1,0.5)--(0,1)--(1.1,1.5);
	\end{braid}$ 
 and the fact that the word beginning with $\ell(\ell-1)^2\cdots(i+2)^2(i-1)$ is not cuspidal by Lemma~\ref{LCuspExpl}, we get 
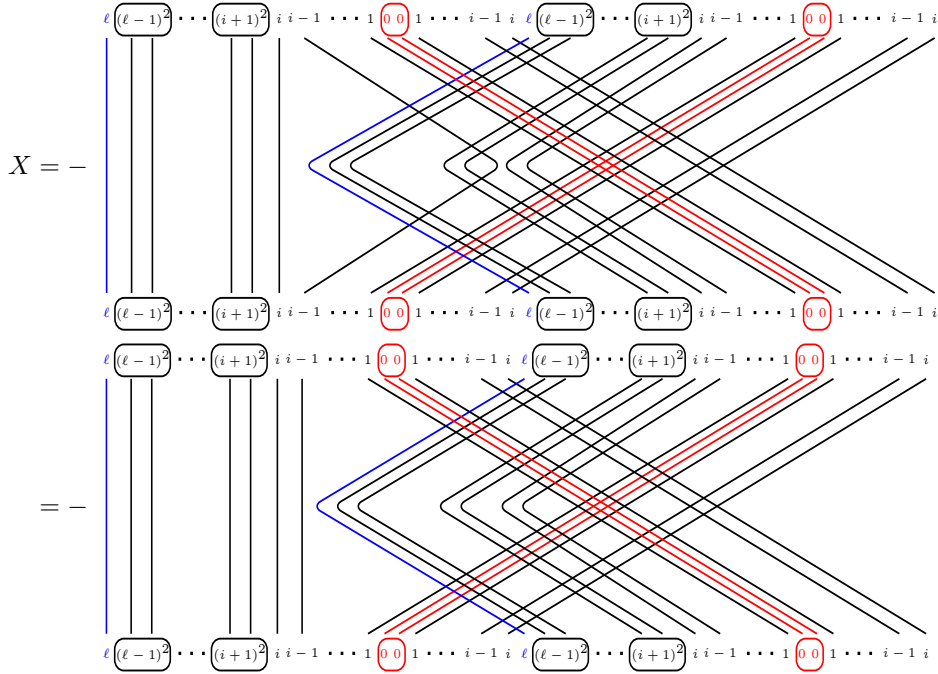
\begin{align*}
X
&=-
\resizebox{116mm}{23mm}{
\begin{braid}\tikzset{baseline=-.3em}
	\draw (0,6) node{\color{blue}$\ell$\color{black}};
	\braidbox{0.6}{2.9}{5.3}{6.6}{};
	\draw (1.8,6) node{$(\ell-1)^2$};
	\draw[dots] (3.5,6)--(5,6);
	\braidbox{5.3}{7.6}{5.3}{6.6}{};
	\draw (6.5,6) node{$(i+1)^2$};
	\draw (8.3,6) node{$i$};
	\draw (9.5,6) node{$i-1$};
	\draw[dots] (10.8,6)--(12.5,6);
	\draw (12.7,6) node{$1$};
	\redbraidbox{13.4}{14.3}{5.3}{6.6}{};
	\draw (13.8,6) node{$\color{red}0\,\,0\color{black}$};
	\draw (15,6) node{$1$};
	\draw[dots] (15.6,6)--(17,6);
	\draw (18.2,6) node{$i-1$};
	\draw (19.5,6) node{$i$};
	\draw (0,-6) node{\color{blue}$\ell$\color{black}};
	\braidbox{0.6}{2.9}{-6.7}{-5.4}{};
	\draw (1.8,-6) node{$(\ell-1)^2$};
	\draw[dots] (3.5,-6)--(5,-6);
	\braidbox{5.3}{7.6}{-6.7}{-5.4}{};
	\draw (6.5,-6) node{$(i+1)^2$};
	\draw (8.3,-6) node{$i$};
	\draw (9.5,-6) node{$i-1$};
	\draw[dots] (10.8,-6)--(12.5,-6);
	\draw (12.7,-6) node{$1$};
	\redbraidbox{13.4}{14.3}{-6.7}{-5.4}{};
	\draw (13.8,-6) node{$\color{red}0\,\,0\color{black}$};
	\draw (15,-6) node{$1$};
	\draw[dots] (15.6,-6)--(17,-6);
	\draw (18.2,-6) node{$i-1$};
	\draw (19.5,-6) node{$i$};
	\draw (20.3,6) node{\color{blue}$\ell$\color{black}};
	\braidbox{20.9}{23.2}{5.3}{6.6}{};
	\draw (22.1,6) node{$(\ell-1)^2$};
	\draw[dots] (23.8,6)--(25.3,6);
	\braidbox{25.6}{27.9}{5.3}{6.6}{};
	\draw (26.8,6) node{$(i+1)^2$};
	\draw (28.6,6) node{$i$};
	\draw (29.8,6) node{$i-1$};
	\draw[dots] (31.1,6)--(32.8,6);
	\draw (33,6) node{$1$};
	\redbraidbox{33.7}{34.6}{5.3}{6.6}{};
	\draw (34.1,6) node{$\color{red}0\,\,0\color{black}$};
	\draw (35.3,6) node{$1$};
	\draw[dots] (35.9,6)--(37.3,6);
	\draw (38.5,6) node{$i-1$};
	\draw (39.8,6) node{$i$};
	\draw (20.3,-6) node{\color{blue}$\ell$\color{black}};
	\braidbox{20.9}{23.2}{-6.7}{-5.4}{};
	\draw (22.1,-6) node{$(\ell-1)^2$};
	\draw[dots] (23.8,-6)--(25.1,-6);
	\braidbox{25.6}{27.9}{-6.7}{-5.4}{};
	\draw (26.8,-6) node{$(i+1)^2$};
	\draw (28.6,-6) node{$i$};
	\draw (29.8,-6) node{$i-1$};
	\draw[dots] (31.1,-6)--(32.8,-6);
	\draw (33,-6) node{$1$};
	\redbraidbox{33.7}{34.6}{-6.7}{-5.4}{};
	\draw (34.1,-6) node{$\color{red}0\,\,0\color{black}$};
	\draw (35.3,-6) node{$1$};
	\draw[dots] (35.9,-6)--(37.3,-6);
	\draw (38.5,-6) node{$i-1$};
	\draw (39.8,-6) node{$i$};
	\draw[blue](0,-5.2)--(0,5.2);
	\draw(1.2,-5.2)--(1.2,5.2);
	\draw(2.2,-5.2)--(2.2,5.2);
	\draw(6,-5.2)--(6,5.2);
	\draw(7,-5.2)--(7,5.2);
	\draw(8.3,-5.2)--(8.3,5.2);
	\draw(9.5,-5.2)--(19,0)--(9.5,5.2);
	\draw(12.7,-5.2)--(33.1,5.2);
	\draw[red](13.5,-5.2)--(33.8,5.2);
	\draw[red](14.2,-5.2)--(34.5,5.2);
	\draw(15,-5.2)--(35.3,5.2);
	\draw(18.2,-5.2)--(38.5,5.2);
	\draw(39.8,-5.2)--(19.5,5.2);
	\draw(39.8,5.2)--(19.5,-5.2);
	\draw[blue](20.3,-5.2)--(9.5,0)--(20.3,5.2);
	\draw(21.3,-5.2)--(10.5,0)--(21.3,5.2);
	\draw(22.3,-5.2)--(11.5,0)--(22.3,5.2);
	\draw(26,-5.2)--(16,0)--(26,5.2);
	\draw(27,-5.2)--(17,0)--(27,5.2);
	\draw(28.6,-5.2)--(19,0)--(28.6,5.2);
	\draw(29.8,-5.2)--(20,0)--(29.8,5.2);
	\draw(38.5,-5.2)--(18.2,5.2);
	\draw(35.3,-5.2)--(15,5.2);
	\draw[red](34.5,-5.2)--(14.2,5.2);
	\draw[red](33.8,-5.2)--(13.5,5.2);
	\draw(33.1,-5.2)--(12.7,5.2);
	\end{braid}
	}
	\\
	&=-
	\resizebox{116mm}{23mm}{
	\begin{braid}\tikzset{baseline=-.3em}
	\draw (0,6) node{\color{blue}$\ell$\color{black}};
	\braidbox{0.6}{2.9}{5.3}{6.6}{};
	\draw (1.8,6) node{$(\ell-1)^2$};
	\draw[dots] (3.5,6)--(5,6);
	\braidbox{5.3}{7.6}{5.3}{6.6}{};
	\draw (6.5,6) node{$(i+1)^2$};
	\draw (8.3,6) node{$i$};
	\draw (9.5,6) node{$i-1$};
	\draw[dots] (10.8,6)--(12.5,6);
	\draw (12.7,6) node{$1$};
	\redbraidbox{13.4}{14.3}{5.3}{6.6}{};
	\draw (13.8,6) node{$\color{red}0\,\,0\color{black}$};
	\draw (15,6) node{$1$};
	\draw[dots] (15.6,6)--(17,6);
	\draw (18.2,6) node{$i-1$};
	\draw (19.5,6) node{$i$};
	\draw (0,-6) node{\color{blue}$\ell$\color{black}};
	\braidbox{0.6}{2.9}{-6.7}{-5.4}{};
	\draw (1.8,-6) node{$(\ell-1)^2$};
	\draw[dots] (3.5,-6)--(5,-6);
	\braidbox{5.3}{7.6}{-6.7}{-5.4}{};
	\draw (6.5,-6) node{$(i+1)^2$};
	\draw (8.3,-6) node{$i$};
	\draw (9.5,-6) node{$i-1$};
	\draw[dots] (10.8,-6)--(12.5,-6);
	\draw (12.7,-6) node{$1$};
	\redbraidbox{13.4}{14.3}{-6.7}{-5.4}{};
	\draw (13.8,-6) node{$\color{red}0\,\,0\color{black}$};
	\draw (15,-6) node{$1$};
	\draw[dots] (15.6,-6)--(17,-6);
	\draw (18.2,-6) node{$i-1$};
	\draw (19.5,-6) node{$i$};
	\draw (20.3,6) node{\color{blue}$\ell$\color{black}};
	\braidbox{20.9}{23.2}{5.3}{6.6}{};
	\draw (22.1,6) node{$(\ell-1)^2$};
	\draw[dots] (23.8,6)--(25.3,6);
	\braidbox{25.6}{27.9}{5.3}{6.6}{};
	\draw (26.8,6) node{$(i+1)^2$};
	\draw (28.6,6) node{$i$};
	\draw (29.8,6) node{$i-1$};
	\draw[dots] (31.1,6)--(32.8,6);
	\draw (33,6) node{$1$};
	\redbraidbox{33.7}{34.6}{5.3}{6.6}{};
	\draw (34.1,6) node{$\color{red}0\,\,0\color{black}$};
	\draw (35.3,6) node{$1$};
	\draw[dots] (35.9,6)--(37.3,6);
	\draw (38.5,6) node{$i-1$};
	\draw (39.8,6) node{$i$};
	\draw (20.3,-6) node{\color{blue}$\ell$\color{black}};
	\braidbox{20.9}{23.2}{-6.7}{-5.4}{};
	\draw (22.1,-6) node{$(\ell-1)^2$};
	\draw[dots] (23.8,-6)--(25.1,-6);
	\braidbox{25.6}{27.9}{-6.7}{-5.4}{};
	\draw (26.8,-6) node{$(i+1)^2$};
	\draw (28.6,-6) node{$i$};
	\draw (29.8,-6) node{$i-1$};
	\draw[dots] (31.1,-6)--(32.8,-6);
	\draw (33,-6) node{$1$};
	\redbraidbox{33.7}{34.6}{-6.7}{-5.4}{};
	\draw (34.1,-6) node{$\color{red}0\,\,0\color{black}$};
	\draw (35.3,-6) node{$1$};
	\draw[dots] (35.9,-6)--(37.3,-6);
	\draw (38.5,-6) node{$i-1$};
	\draw (39.8,-6) node{$i$};
	\draw[blue](0,-5.2)--(0,5.2);
	\draw(1.2,-5.2)--(1.2,5.2);
	\draw(2.2,-5.2)--(2.2,5.2);
	\draw(6,-5.2)--(6,5.2);
	\draw(7,-5.2)--(7,5.2);
	\draw(8.3,-5.2)--(8.3,5.2);
	\draw(9.5,-5.2)--(9.5,5.2);
	\draw(12.7,-5.2)--(33.1,5.2);
	\draw[red](13.5,-5.2)--(33.8,5.2);
	\draw[red](14.2,-5.2)--(34.5,5.2);
	\draw(15,-5.2)--(35.3,5.2);
	\draw(18.2,-5.2)--(38.5,5.2);
	\draw(39.8,-5.2)--(19.5,5.2);
	\draw(39.8,5.2)--(19.5,-5.2);
	\draw[blue](20.3,-5.2)--(10,0)--(20.3,5.2);
	\draw(21.3,-5.2)--(11,0)--(21.3,5.2);
	\draw(22.3,-5.2)--(12,0)--(22.3,5.2);
	\draw(26,-5.2)--(16,0)--(26,5.2);
	\draw(27,-5.2)--(17,0)--(27,5.2);
	\draw(28.6,-5.2)--(19,0)--(28.6,5.2);
	\draw(29.8,-5.2)--(20,0)--(29.8,5.2);
	\draw(38.5,-5.2)--(18.2,5.2);
	\draw(35.3,-5.2)--(15,5.2);
	\draw[red](34.5,-5.2)--(14.2,5.2);
	\draw[red](33.8,-5.2)--(13.5,5.2);
	\draw(33.1,-5.2)--(12.7,5.2);
	\end{braid}.
	}
\end{align*}
Continuing like this and using braid relations $i-2$ more times, we get 
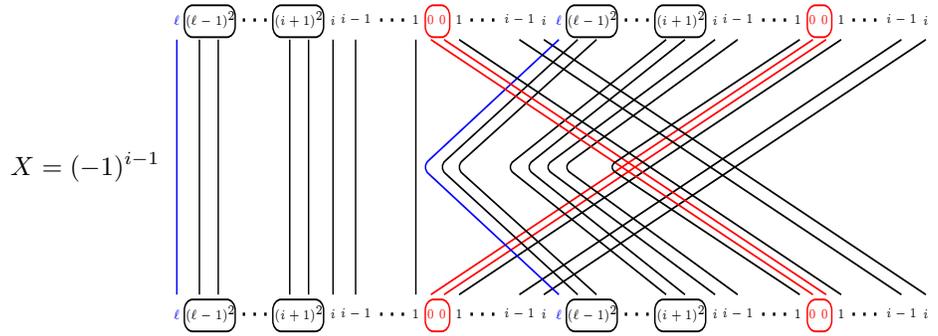
\begin{align*}
X
&=(-1)^{i-1}
\resizebox{106mm}{23mm}{
	\begin{braid}\tikzset{baseline=-.3em}
	\draw (0,6) node{\color{blue}$\ell$\color{black}};
	\braidbox{0.6}{2.9}{5.3}{6.6}{};
	\draw (1.8,6) node{$(\ell-1)^2$};
	\draw[dots] (3.5,6)--(5,6);
	\braidbox{5.3}{7.6}{5.3}{6.6}{};
	\draw (6.5,6) node{$(i+1)^2$};
	\draw (8.3,6) node{$i$};
	\draw (9.5,6) node{$i-1$};
	\draw[dots] (10.8,6)--(12.5,6);
	\draw (12.7,6) node{$1$};
	\redbraidbox{13.4}{14.3}{5.3}{6.6}{};
	\draw (13.8,6) node{$\color{red}0\,\,0\color{black}$};
	\draw (15,6) node{$1$};
	\draw[dots] (15.6,6)--(17,6);
	\draw (18.2,6) node{$i-1$};
	\draw (19.5,6) node{$i$};
	\draw (0,-6) node{\color{blue}$\ell$\color{black}};
	\braidbox{0.6}{2.9}{-6.7}{-5.4}{};
	\draw (1.8,-6) node{$(\ell-1)^2$};
	\draw[dots] (3.5,-6)--(5,-6);
	\braidbox{5.3}{7.6}{-6.7}{-5.4}{};
	\draw (6.5,-6) node{$(i+1)^2$};
	\draw (8.3,-6) node{$i$};
	\draw (9.5,-6) node{$i-1$};
	\draw[dots] (10.8,-6)--(12.5,-6);
	\draw (12.7,-6) node{$1$};
	\redbraidbox{13.4}{14.3}{-6.7}{-5.4}{};
	\draw (13.8,-6) node{$\color{red}0\,\,0\color{black}$};
	\draw (15,-6) node{$1$};
	\draw[dots] (15.6,-6)--(17,-6);
	\draw (18.2,-6) node{$i-1$};
	\draw (19.5,-6) node{$i$};
	\draw (20.3,6) node{\color{blue}$\ell$\color{black}};
	\braidbox{20.9}{23.2}{5.3}{6.6}{};
	\draw (22.1,6) node{$(\ell-1)^2$};
	\draw[dots] (23.8,6)--(25.3,6);
	\braidbox{25.6}{27.9}{5.3}{6.6}{};
	\draw (26.8,6) node{$(i+1)^2$};
	\draw (28.6,6) node{$i$};
	\draw (29.8,6) node{$i-1$};
	\draw[dots] (31.1,6)--(32.8,6);
	\draw (33,6) node{$1$};
	\redbraidbox{33.7}{34.6}{5.3}{6.6}{};
	\draw (34.1,6) node{$\color{red}0\,\,0\color{black}$};
	\draw (35.3,6) node{$1$};
	\draw[dots] (35.9,6)--(37.3,6);
	\draw (38.5,6) node{$i-1$};
	\draw (39.8,6) node{$i$};
	\draw (20.3,-6) node{\color{blue}$\ell$\color{black}};
	\braidbox{20.9}{23.2}{-6.7}{-5.4}{};
	\draw (22.1,-6) node{$(\ell-1)^2$};
	\draw[dots] (23.8,-6)--(25.1,-6);
	\braidbox{25.6}{27.9}{-6.7}{-5.4}{};
	\draw (26.8,-6) node{$(i+1)^2$};
	\draw (28.6,-6) node{$i$};
	\draw (29.8,-6) node{$i-1$};
	\draw[dots] (31.1,-6)--(32.8,-6);
	\draw (33,-6) node{$1$};
	\redbraidbox{33.7}{34.6}{-6.7}{-5.4}{};
	\draw (34.1,-6) node{$\color{red}0\,\,0\color{black}$};
	\draw (35.3,-6) node{$1$};
	\draw[dots] (35.9,-6)--(37.3,-6);
	\draw (38.5,-6) node{$i-1$};
	\draw (39.8,-6) node{$i$};
	\draw[blue](0,-5.2)--(0,5.2);
	\draw(1.2,-5.2)--(1.2,5.2);
	\draw(2.2,-5.2)--(2.2,5.2);
	\draw(6,-5.2)--(6,5.2);
	\draw(7,-5.2)--(7,5.2);
	\draw(8.3,-5.2)--(8.3,5.2);
	\draw(9.5,-5.2)--(9.5,5.2);
	\draw(12.7,-5.2)--(12.7,5.2);
	\draw[red](13.5,-5.2)--(33.8,5.2);
	\draw[red](14.2,-5.2)--(34.5,5.2);
	\draw(15,-5.2)--(35.3,5.2);
	\draw(18.2,-5.2)--(38.5,5.2);
	\draw(39.8,-5.2)--(19.5,5.2);
	\draw(39.8,5.2)--(19.5,-5.2);
	\draw[blue](20.3,-5.2)--(13,0)--(20.3,5.2);
	\draw(21.3,-5.2)--(13.9,0)--(21.3,5.2);
	\draw(22.3,-5.2)--(14.8,0)--(22.3,5.2);
	\draw(26,-5.2)--(17.5,0)--(26,5.2);
	\draw(27,-5.2)--(18.5,0)--(27,5.2);
	\draw(28.6,-5.2)--(19.5,0)--(28.6,5.2);
	\draw(29.8,-5.2)--(20.5,0)--(29.8,5.2);
	\draw(38.5,-5.2)--(18.2,5.2);
	\draw(35.3,-5.2)--(15,5.2);
	\draw[red](34.5,-5.2)--(14.2,5.2);
	\draw[red](33.8,-5.2)--(13.5,5.2);
	\draw(33.1,-5.2)--(22.9,0)--(33.1,5.2);
	\end{braid}.
	}
\end{align*}
Applying Lemma~\ref{LB} and using non-cuspidality of a word starting with $\ell(\ell-1)^2\cdots(i+1)^2i\cdots 20$, cf.  Lemma~\ref{LCuspExpl}, we get
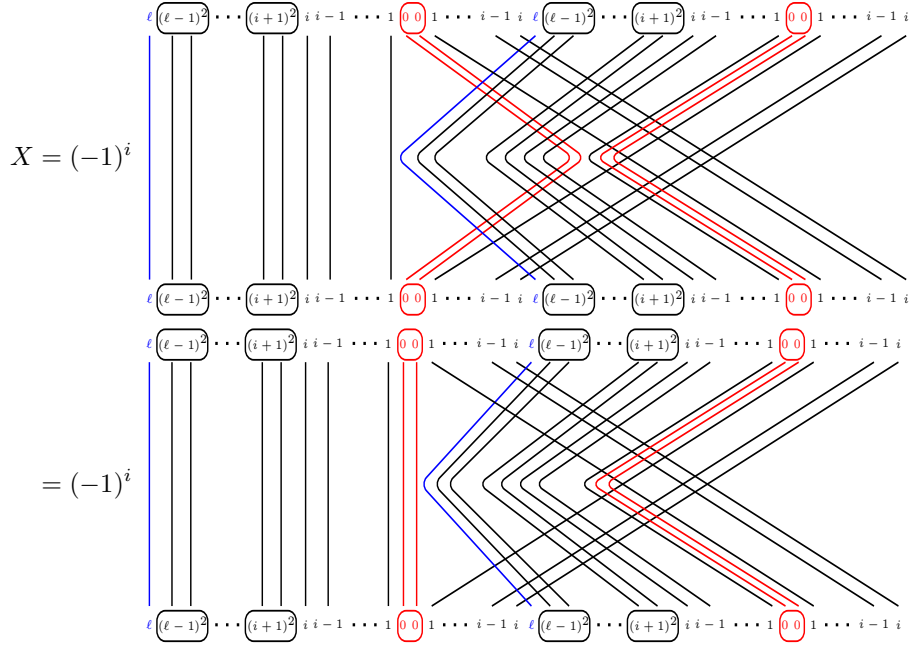
\begin{align*}
X
&=(-1)^{i}
\resizebox{106mm}{22mm}{
	\begin{braid}\tikzset{baseline=-.3em}
	\draw (0,6) node{\color{blue}$\ell$\color{black}};
	\braidbox{0.6}{2.9}{5.3}{6.6}{};
	\draw (1.8,6) node{$(\ell-1)^2$};
	\draw[dots] (3.5,6)--(5,6);
	\braidbox{5.3}{7.6}{5.3}{6.6}{};
	\draw (6.5,6) node{$(i+1)^2$};
	\draw (8.3,6) node{$i$};
	\draw (9.5,6) node{$i-1$};
	\draw[dots] (10.8,6)--(12.5,6);
	\draw (12.7,6) node{$1$};
	\redbraidbox{13.4}{14.3}{5.3}{6.6}{};
	\draw (13.8,6) node{$\color{red}0\,\,0\color{black}$};
	\draw (15,6) node{$1$};
	\draw[dots] (15.6,6)--(17,6);
	\draw (18.2,6) node{$i-1$};
	\draw (19.5,6) node{$i$};
	\draw (0,-6) node{\color{blue}$\ell$\color{black}};
	\braidbox{0.6}{2.9}{-6.7}{-5.4}{};
	\draw (1.8,-6) node{$(\ell-1)^2$};
	\draw[dots] (3.5,-6)--(5,-6);
	\braidbox{5.3}{7.6}{-6.7}{-5.4}{};
	\draw (6.5,-6) node{$(i+1)^2$};
	\draw (8.3,-6) node{$i$};
	\draw (9.5,-6) node{$i-1$};
	\draw[dots] (10.8,-6)--(12.5,-6);
	\draw (12.7,-6) node{$1$};
	\redbraidbox{13.4}{14.3}{-6.7}{-5.4}{};
	\draw (13.8,-6) node{$\color{red}0\,\,0\color{black}$};
	\draw (15,-6) node{$1$};
	\draw[dots] (15.6,-6)--(17,-6);
	\draw (18.2,-6) node{$i-1$};
	\draw (19.5,-6) node{$i$};
	\draw (20.3,6) node{\color{blue}$\ell$\color{black}};
	\braidbox{20.9}{23.2}{5.3}{6.6}{};
	\draw (22.1,6) node{$(\ell-1)^2$};
	\draw[dots] (23.8,6)--(25.3,6);
	\braidbox{25.6}{27.9}{5.3}{6.6}{};
	\draw (26.8,6) node{$(i+1)^2$};
	\draw (28.6,6) node{$i$};
	\draw (29.8,6) node{$i-1$};
	\draw[dots] (31.1,6)--(32.8,6);
	\draw (33,6) node{$1$};
	\redbraidbox{33.7}{34.6}{5.3}{6.6}{};
	\draw (34.1,6) node{$\color{red}0\,\,0\color{black}$};
	\draw (35.3,6) node{$1$};
	\draw[dots] (35.9,6)--(37.3,6);
	\draw (38.5,6) node{$i-1$};
	\draw (39.8,6) node{$i$};
	\draw (20.3,-6) node{\color{blue}$\ell$\color{black}};
	\braidbox{20.9}{23.2}{-6.7}{-5.4}{};
	\draw (22.1,-6) node{$(\ell-1)^2$};
	\draw[dots] (23.8,-6)--(25.1,-6);
	\braidbox{25.6}{27.9}{-6.7}{-5.4}{};
	\draw (26.8,-6) node{$(i+1)^2$};
	\draw (28.6,-6) node{$i$};
	\draw (29.8,-6) node{$i-1$};
	\draw[dots] (31.1,-6)--(32.8,-6);
	\draw (33,-6) node{$1$};
	\redbraidbox{33.7}{34.6}{-6.7}{-5.4}{};
	\draw (34.1,-6) node{$\color{red}0\,\,0\color{black}$};
	\draw (35.3,-6) node{$1$};
	\draw[dots] (35.9,-6)--(37.3,-6);
	\draw (38.5,-6) node{$i-1$};
	\draw (39.8,-6) node{$i$};
	\draw[blue](0,-5.2)--(0,5.2);
	\draw(1.2,-5.2)--(1.2,5.2);
	\draw(2.2,-5.2)--(2.2,5.2);
	\draw(6,-5.2)--(6,5.2);
	\draw(7,-5.2)--(7,5.2);
	\draw(8.3,-5.2)--(8.3,5.2);
	\draw(9.5,-5.2)--(9.5,5.2);
	\draw(12.7,-5.2)--(12.7,5.2);
	\draw[red](13.5,-5.2)--(22.3,0)--(13.5,5.2);
	\draw[red](14.2,-5.2)--(22.9,0)--(14.2,5.2);
	\draw(15,-5.2)--(35.3,5.2);
	\draw(18.2,-5.2)--(38.5,5.2);
	\draw(39.8,-5.2)--(19.5,5.2);
	\draw(39.8,5.2)--(19.5,-5.2);
	\draw[blue](20.3,-5.2)--(13,0)--(20.3,5.2);
	\draw(21.3,-5.2)--(13.9,0)--(21.3,5.2);
	\draw(22.3,-5.2)--(14.8,0)--(22.3,5.2);
	\draw(26,-5.2)--(17.5,0)--(26,5.2);
	\draw(27,-5.2)--(18.5,0)--(27,5.2);
	\draw(28.6,-5.2)--(19.5,0)--(28.6,5.2);
	\draw(29.8,-5.2)--(20.5,0)--(29.8,5.2);
	\draw(38.5,-5.2)--(18.2,5.2);
	\draw(35.3,-5.2)--(15,5.2);
	\draw[red](34.5,-5.2)--(24.2,0)--(34.5,5.2);
	\draw[red](33.8,-5.2)--(23.5,0)--(33.8,5.2);
	\draw(33.1,-5.2)--(22.9,0)--(33.1,5.2);
	\end{braid}
	}
	\\
	&
	=(-1)^{i}
	\resizebox{106mm}{22mm}{
	\begin{braid}\tikzset{baseline=-.3em}
	\draw (0,6) node{\color{blue}$\ell$\color{black}};
	\braidbox{0.6}{2.9}{5.3}{6.6}{};
	\draw (1.8,6) node{$(\ell-1)^2$};
	\draw[dots] (3.5,6)--(5,6);
	\braidbox{5.3}{7.6}{5.3}{6.6}{};
	\draw (6.5,6) node{$(i+1)^2$};
	\draw (8.3,6) node{$i$};
	\draw (9.5,6) node{$i-1$};
	\draw[dots] (10.8,6)--(12.5,6);
	\draw (12.7,6) node{$1$};
	\redbraidbox{13.4}{14.3}{5.3}{6.6}{};
	\draw (13.8,6) node{$\color{red}0\,\,0\color{black}$};
	\draw (15,6) node{$1$};
	\draw[dots] (15.6,6)--(17,6);
	\draw (18.2,6) node{$i-1$};
	\draw (19.5,6) node{$i$};
	\draw (0,-6) node{\color{blue}$\ell$\color{black}};
	\braidbox{0.6}{2.9}{-6.7}{-5.4}{};
	\draw (1.8,-6) node{$(\ell-1)^2$};
	\draw[dots] (3.5,-6)--(5,-6);
	\braidbox{5.3}{7.6}{-6.7}{-5.4}{};
	\draw (6.5,-6) node{$(i+1)^2$};
	\draw (8.3,-6) node{$i$};
	\draw (9.5,-6) node{$i-1$};
	\draw[dots] (10.8,-6)--(12.5,-6);
	\draw (12.7,-6) node{$1$};
	\redbraidbox{13.4}{14.3}{-6.7}{-5.4}{};
	\draw (13.8,-6) node{$\color{red}0\,\,0\color{black}$};
	\draw (15,-6) node{$1$};
	\draw[dots] (15.6,-6)--(17,-6);
	\draw (18.2,-6) node{$i-1$};
	\draw (19.5,-6) node{$i$};
	\draw (20.3,6) node{\color{blue}$\ell$\color{black}};
	\braidbox{20.9}{23.2}{5.3}{6.6}{};
	\draw (22.1,6) node{$(\ell-1)^2$};
	\draw[dots] (23.8,6)--(25.3,6);
	\braidbox{25.6}{27.9}{5.3}{6.6}{};
	\draw (26.8,6) node{$(i+1)^2$};
	\draw (28.6,6) node{$i$};
	\draw (29.8,6) node{$i-1$};
	\draw[dots] (31.1,6)--(32.8,6);
	\draw (33,6) node{$1$};
	\redbraidbox{33.7}{34.6}{5.3}{6.6}{};
	\draw (34.1,6) node{$\color{red}0\,\,0\color{black}$};
	\draw (35.3,6) node{$1$};
	\draw[dots] (35.9,6)--(37.3,6);
	\draw (38.5,6) node{$i-1$};
	\draw (39.8,6) node{$i$};
	\draw (20.3,-6) node{\color{blue}$\ell$\color{black}};
	\braidbox{20.9}{23.2}{-6.7}{-5.4}{};
	\draw (22.1,-6) node{$(\ell-1)^2$};
	\draw[dots] (23.8,-6)--(25.1,-6);
	\braidbox{25.6}{27.9}{-6.7}{-5.4}{};
	\draw (26.8,-6) node{$(i+1)^2$};
	\draw (28.6,-6) node{$i$};
	\draw (29.8,-6) node{$i-1$};
	\draw[dots] (31.1,-6)--(32.8,-6);
	\draw (33,-6) node{$1$};
	\redbraidbox{33.7}{34.6}{-6.7}{-5.4}{};
	\draw (34.1,-6) node{$\color{red}0\,\,0\color{black}$};
	\draw (35.3,-6) node{$1$};
	\draw[dots] (35.9,-6)--(37.3,-6);
	\draw (38.5,-6) node{$i-1$};
	\draw (39.8,-6) node{$i$};
	\draw[blue](0,-5.2)--(0,5.2);
	\draw(1.2,-5.2)--(1.2,5.2);
	\draw(2.2,-5.2)--(2.2,5.2);
	\draw(6,-5.2)--(6,5.2);
	\draw(7,-5.2)--(7,5.2);
	\draw(8.3,-5.2)--(8.3,5.2);
	\draw(9.5,-5.2)--(9.5,5.2);
	\draw(12.7,-5.2)--(12.7,5.2);
	\draw[red](13.5,-5.2)--(13.5,5.2);
	\draw[red](14.2,-5.2)--(14.2,5.2);
	\draw(15,-5.2)--(35.3,5.2);
	\draw(18.2,-5.2)--(38.5,5.2);
	\draw(39.8,-5.2)--(19.5,5.2);
	\draw(39.8,5.2)--(19.5,-5.2);
	\draw[blue](20.3,-5.2)--(14.4,0)--(20.3,5.2);
	\draw(21.3,-5.2)--(15.1,0)--(21.3,5.2);
	\draw(22.3,-5.2)--(15.8,0)--(22.3,5.2);
	\draw(26,-5.2)--(17.5,0)--(26,5.2);
	\draw(27,-5.2)--(18.5,0)--(27,5.2);
	\draw(28.6,-5.2)--(19.5,0)--(28.6,5.2);
	\draw(29.8,-5.2)--(20.5,0)--(29.8,5.2);
	\draw(38.5,-5.2)--(18.2,5.2);
	\draw(35.3,-5.2)--(15,5.2);
	\draw[red](34.5,-5.2)--(24.2,0)--(34.5,5.2);
	\draw[red](33.8,-5.2)--(23.5,0)--(33.8,5.2);
	\draw(33.1,-5.2)--(22.9,0)--(33.1,5.2);
	\end{braid}.
	}
\end{align*}
Thus that $X=(-1)^i\Theta^i_1$. In view of (\ref{E160721}) and  Lemma~\ref{L110621Gen}, we get 
\begin{align*}
\ga^{i,i}\Upsilon\ga^{i,i}&
=(-1)^i\ga^{i,i}(-z_1+z_2+(-1)^i(c_1+c_2))\ga^{i,i}(-z_1-z_2)\ga^{i,i}
\\&=
(-1)^i\ga^{i,i}(z_1^2-z_2^2)\ga^{i,i}-\ga^{i,i}(c_1+c_2)(z_1+z_2)\ga^{i,i},
\end{align*}
which implies the lemma. 
\end{proof}

\begin{Lemma} \label{LKeyDifferentColors} 
Let $i,j\in J$ with $i\neq j$. Then  
$\ga^{j,i}\Upsilon\ga^{i,j}\in B_{1^2}^{(1)}$.
\end{Lemma}
\begin {proof}
In view of Corollary~\ref{CBBasis}, the relations in $B_2$ allow us to write $\ga^{j,i}\Upsilon\ga^{i,j}$ as an element of  $\bigoplus_{c_1,c_2\in\Z_{\geq 0}}\ga^{j,i}z_1^{c_1} z_2^{c_2} \,\Zig_\ell^{\otimes 2}\ga^{i,j}.$ The result now follows by degrees. 
\end{proof}

\printindex

\end{document}